\documentclass[reqno,10pt]{amsart}
\usepackage{amscd}
\usepackage{amsfonts}
\usepackage{amssymb}
\usepackage{latexsym}
\usepackage{color}
\usepackage{esint}

 \usepackage{cancel}

\newtheorem{theorem}{Theorem}

\newtheorem{definition}{Definition}
\newtheorem{lemma}{Lemma}
\newtheorem{proposition}[theorem]{Proposition}
\newtheorem{remark}{Remark}

\let\e=\varepsilon

\let\pt=\partial
\let\p=\partial

\let\O=\Omega

\let\o=\omega
\let\g=\gamma

\let\b=\beta

\numberwithin{equation}{section}

\let\hide\iffalse
\let\unhide\fi

\newcommand{\R}{\mathbb{R}}

\renewcommand{\S}{\mathbb{S}}

\newcommand{\be}{\begin{equation}}
\newcommand{\bm}{\begin{multline}}
\newcommand{\ee}{\end{equation}}
\newcommand{\dd}{\mathrm{d}}

\newcommand{\xb}{x_{\mathbf{b}}}
\newcommand{\tb}{t_{\mathbf{b}}}
\newcommand{\vb}{v_{\mathbf{b}}}
\newcommand{\xbt}{\tilde{x}_{\mathbf{b}}}

\newcommand{\xft}{\tilde{x}_{\mathbf{f}}}

\newcommand{\xf}{x_{\mathbf{f}}}
\newcommand{\tf}{t_{\mathbf{f}}}
\newcommand{\vf}{v_{\mathbf{f}}}

\newcommand{\xba}{x^1_{\mathbf{b}}}
\newcommand{\xbb}{x^2_{\mathbf{b}}}

\newcommand{\Bes}{\begin{eqnarray*}}
\newcommand{\Ees}{\end{eqnarray*}}
\newcommand{\Be}{\begin{equation} }
\newcommand{\Ee}{\end{equation}}

\def\p{\partial}

\def\O{\Omega}
\def\R{\mathbb{R}}

\def\B{\begin{equation}}
\def\E{\end{equation}}
\def\BN{\begin{eqnarray*}}
\def\EN{\end{eqnarray*}}

\begin{document}
 \date{}

\title{Global strong solutions of the Vlasov-Poisson-Boltzmann system in bounded domains}

 \author{Yunbai Cao}
   
 \author{Chanwoo Kim}
 \email{chanwoo.kim@wisc.edu, ckim.pde@gmail.com}
  \address{Department of Mathematics, University of Wisconsin, Madison, WI 53706 USA}

 \author{Donghyun Lee}

 \begin{abstract}
 When dilute charged particles are confined in a bounded domain, boundary effects are crucial in the global dynamics. We construct a unique global-in-time solution to the Vlasov-Poisson-Boltzmann system in convex domains with the diffuse boundary condition. The construction is based on an $L^{2}$-$L^{\infty}$ framework with a novel \textit{nonlinear-normed energy estimate} of a distribution function in weighted $W^{1,p}$-spaces and a $C^{2,\delta}$-estimate of the self-consistent electric potential. Moreover we prove an exponential convergence of the distribution function toward the global Maxwellian. 
 \end{abstract}
 
\maketitle
 
 \tableofcontents

\section{Introduction}
\label{intro}

The object of kinetic theory is the modeling of many particles whose behavior is encoded in a distribution function on the phase space, which is denoted by $F(t,x,v)$ for $(t,x,v) \in [0, \infty) \times  {\O} \times \R^{3}$. Here $\O$ denotes an open subset of $\R^{3}$. Dynamics and collision processes of dilute charged particles subjected to an field $E$ can be modeled by the Vlasov-Boltzmann equation:
\begin{equation}\label{Boltzmann_E}
\partial_{t} F + v\cdot \nabla_{x} F + E\cdot \nabla_{v} F = Q(F,F).
\end{equation}
The collision operator measures the change rate in binary hard sphere collisions and takes the form of
\begin{equation}\begin{split}\label{Q}
Q(F_{1},F_{2}) (v)&: = Q_\mathrm{gain}(F_1,F_2)-Q_\mathrm{loss}(F_1,F_2)\\
&: =\int_{\R^3} \int_{\S^2} 
|(v-u) \cdot \omega| [F_1 (u^\prime) F_2 (v^\prime) - F_1 (u) F_2 (v)]
 \dd \omega \dd u,
\end{split}\end{equation}  
where $u^\prime = u - [(u-v) \cdot \omega] \omega$ and $v^\prime = v + [(u-v) \cdot \omega] \omega$. The collision operator enjoys a collision invariance: for any measurable $G$,  
\Be\label{collison_invariance}
\int_{\R^{3}} \begin{bmatrix}1 & v & \frac{|v|^{2}-3}{2}\end{bmatrix} Q(G,G) \dd v = \begin{bmatrix}0 & 0 & 0 \end{bmatrix} .
\Ee 
It is well-known that a global Maxwellian $\mu$ 
satisfies $Q(\mu,\mu)=0$, where
\Be\label{Maxwellian}
\mu(v):= \frac{1}{(2\pi)^{3/2}} \exp\bigg(
 - \frac{|v |^{2}}{2 }
 \bigg).
\Ee

Due to its importance of the Boltzmann equation in theories and their applications,
there have been explosive research results in analytic study of the equation. The nonlinear energy method has led to solutions of many open problems (\cite{Guo_P,Guo_M}) including global classical solutions and their asymptotic stability of Boltzmann equation coupled with either the Poisson equation or the Maxwell system for electromagnetism when the initial data are close to the
Maxwellian $\mu$. For large-amplitude solutions, renormalized solutions have been studied extensively from the end of 80's (see \cite{DL,Laure} and references therein). An asymptotic stability of some large-amplitude solutions is established in \cite{DV,GMM}, provided certain a priori strong Sobolev estimates can be verified. Such high regularity insures an $L^\infty$-control of solutions which is crucial to handle the quadratic nonlinearity. It should be noted that all of these results deal with idealized periodic domains or
whole space, in which the solutions can remain bounded in high Sobolev spaces.

In many important physical applications, e.g. semiconductor and tokamak, the charged dilute
gas is confined within a container, and its interaction with the boundary often plays a crucial role in global dynamics. The interaction of the gas with a boundary is described by a suitable boundary conditions \cite{Maxwell}. In this paper we consider one of the physical conditions, a so-called diffuse boundary condition:
\Be\label{diffuse_BC}
F(t,x,v)=c_\mu \mu(v) \int_{n(x ) \cdot u>0} F(t,x,u) \{n(x) \cdot u\}\dd u \ \  {for} \ (x,v) \in \gamma_-.
\Ee
Here, an outgoing set is defined as 
\Be\label{gamma_-}
\gamma_-: = \{(x,v) \in \p\O \times \R^3: n(x) \cdot v<0\},
\Ee
and $n(x)$ is the outward unit normal at a boundary point $x$. A number $c_\mu$ is chosen to be $\sqrt{2\pi}$ so that $c_\mu \int_{n(x) \cdot u > 0} \mu(u) \{n(x) \cdot u\} \dd u=1$ with $\mu$ in (\ref{Maxwellian}). 
Due to this normalization, the distribution of (\ref{diffuse_BC}) enjoys a null flux condition at the boundary: 
\Be\label{null_flux}
\int_{\R^{3}} F(t,x,v) \{n(x) \cdot v\} \dd v =0 \ \ {for}  \ x \in\p\O.
\Ee

In general, high regularity may not be expected for solutions of the Boltzmann equation in physical bounded domains. Such a drastic boundary effect has been
demonstrated recently by the second author and his collaborators as the formation and propagation of discontinuity in non-convex domains in \cite{Kim11,EGKM}, and a non-existence of some second order derivative at the boundary even in convex domains \cite{GKTT1}. Evidently these results show that smooth solutions are generally unavailable to the boundary problems (cf. \cite{DV}).

In order to overcome such critical difficulty, an $L^2$-$L^\infty$ framework has been developed in \cite{Guo10} to study the Boltzmann equation with various boundary conditions. The core of the method lays in a direct approach (without taking derivatives) to achieve a pointwise bound using trajectory of the transport operator, which leads substantial development in various directions (\cite{EGKM2,EGKM,GKTT1,GKTT2,KL1,KL2,KL3,K2,strain,GJ,KY,KL3,CKL}). Among others we briefly point out some relevant works to the Vlasov-Boltzmann equation (\ref{Boltzmann_E}) in bounded domains. In \cite{KL1}, the second and the third authors established the global well-posedness and asymptotic stability of the Boltzmann equation near the global Maxwellian for the specular
reflection boundary condition with or without \textit{given} small external fields. In \cite{EGKM2}, the second author and other collaborators studied diffusive hydrodynamic limit of Boltzmann to Navier-Stokes-Fourier system in both steady and unsteady cases in the presence of \textit{given} small external fields.

However, in the study of charged particles generating \textit{self-consistent fields}, many fundamental boundary problems of the kinetic models are still widely open. This category of problems is important in the application to plasma or galaxy. One of the major difficulties is that trajectories are curved and behave in a very complicated way when they hit the boundary. 

The field $E$, that we are interested in, is associated with an electrostatic potential $\phi$ as
\Be\label{Field}
E(t,x): =  - \nabla  \phi(t,x),
\Ee
where the potential is determined by the Poisson equation: 
\Be
- \Delta \phi(t,x)    =  \int_{\mathbb{R}^{3}} F(t,x,v) \dd v - \rho_0  \ \  {in} \ \O , \label{Poisson}\\
\Ee
with the \textit{zero} Neumann boundary condition 
\Be
 \frac{\p \phi}{\p n}     =  0  \ \   {on} \ \p\O.  \label{phi_BC}
\Ee
Here a background density $\rho_{0}$ in (\ref{Poisson}) is a constant number. 
 %
 
The coupled system of (\ref{Boltzmann_E}) with (\ref{Field}) and (\ref{Poisson}) is a so-called \textit{Vlasov-Poisson-Boltzmann system} (VPB), which describes the dynamics of electrons in the absence of a magnetic field. This model has been considered as a fundamental collisional plasma model and has attracted much attention (see \cite{DD,Guo_P,GJ,GlasseyStrauss,Michler} and the references therein).

From (\ref{collison_invariance}) and (\ref{null_flux}), a strong solution of VPB with the diffuse BC (\ref{diffuse_BC}) preserves total mass:
\Be\label{conservation_mass}
\iint_{\O\times\R^{3}} F(t,x,v) \dd v \dd x \equiv  \iint_{\O\times\R^{3}} F(0,x,v) \dd v \dd x
\ \  {for \ all } \  t\geq 0.
\Ee
We set, without loss of generality, 
\Be\label{neutral_condition}
\rho_0  = \frac{1}{|\Omega|}\iint_{\O \times \mathbb{R}^{3}} F(0,x,v )  \dd v \dd x
. \ \ \   {(a \ neutral \  condition)}
\Ee
Then $\int_\O \left\{\int_{\mathbb{R}^{3}} F(t,x,v) \dd v - \rho_0\right\} \dd x =0$ for all $t>0$ from (\ref{conservation_mass}). This zero-mean condition guarantees a solvability of the Poisson equation (\ref{Poisson}) with the Neumann boundary condition (\ref{phi_BC}).

Without loss of generality, by some rescaling, we assume that
\Be\label{same_mass_mu}
\iint_{\O \times \R^3} \mu \dd v \dd x = \iint_{\O \times \R^3} F(0,x,v) \dd v \dd x.
\Ee
We consider a perturbation around $\mu$:
 \Be
F(t,x,v)= \mu  + \sqrt{
\mu} f(t,x,v ).\label{perturbation}  
\Ee 
The corresponding problem of (\ref{perturbation}) is given by 
\Be
 \p_t f + v\cdot \nabla_x f - \nabla_x  \phi _f \cdot \nabla_v f + \frac{v}{2} \cdot \nabla_x \phi_f  f +   Lf  
 =   \Gamma(f,f) - v\cdot \nabla_x \phi_f  \sqrt{\mu}  ,\label{eqtn_f} 
\Ee
 with $ f(0,x,v) = f_0 (x,v)$ and 
 \Be
 - \Delta \phi _f (t,x)= \int_{\R^{3}} f (t,x,v)\sqrt{ \mu(v)} \dd v  \ \  {in} \ \O , 
\ \ 
 \ \frac{\p \phi_f}{\p n}=0 \ \  {on} \ \p\O,\label{phi_f}\Ee
 and 
 \Be
 f(t,x,v)
= c_\mu \sqrt{\mu} \int_{n(x) \cdot u>0} 
f(t,x,u) \sqrt{\mu(u)}\{n(x) \cdot u\} \dd u \ \  {for} \ (x,v) \in\gamma_-.\label{BC_f}
\Ee
We note that from (\ref{conservation_mass}) and (\ref{same_mass_mu})
\Be\label{mass_conservation_f}
\int_{\O}\int_{\R^3} f(t,x,v) \sqrt{\mu(v) } \dd v   \dd x =0 \ \    {for \  all } \ t \geq 0.
\Ee
The standard notions of $L, \nu, K,$ and $\Gamma$ are defined in (\ref{L_decomposition})-(\ref{Gamma_def}).

\hide
$f(0,x,v) = f_0 (x,v)$ and 
\begin{equation}\label{eqtn_f} 
 \p_t f + v\cdot \nabla_x f - \nabla_x  \phi _f \cdot \nabla_v f + \frac{v}{2} \cdot \nabla_x \phi_f  f +   Lf 
 =   \Gamma(f,f) - v\cdot \nabla_x \phi_f  \sqrt{\mu}  ,
\end{equation} 
and
\begin{equation}\label{phi_f}
- \Delta \phi _f (t,x)= \int_{\R^{3}} f (t,x,v)\sqrt{ \mu(v)} \dd v , 
\ \ 
 \ \frac{\p \phi_f}{\p n}=0 \ \ \text{on} \ \p\O,
\end{equation} 
with
\Be\label{BC_f}
f(t,x,v)
= c_\mu \sqrt{\mu} \int_{n(x) \cdot u>0} 
f(t,x,u) \sqrt{\mu(u)}\{n(x) \cdot u\} \dd u \ \ \text{for} \ (x,v) \in\gamma_-.
\Ee \unhide

\hide

The equation of characteristics is 
\begin{equation}
\frac{dX(s;t,x,v)}{ds} = V(s;t,x,v),  \ \ \ 
\frac{dV(s;t,x,v)}{ds} = - \nabla_x  \phi(s, X(s;t,x,v))   .
\end{equation}
\unhide

\subsection{A New Kinetic Weight} 
 An intrinsic feature of the transport equation in bounded domains is the singular behavior of its derivatives. Let us consider a solution of the Vlasov-Poisson equation 
\Be\label{Vlasov_eqtn}
\p_t f + v\cdot \nabla_x f - \nabla_x \phi_f \cdot \nabla_vf =0, 
\Ee
where $\phi_f(t,x)$ satisfies (\ref{phi_f}). Assume that the boundary condition is determined by some given function $g$ as
\Be\label{inflow_BC}
f(t,x,v) = g(t,x,v) \ \  { for } \  (x,v) \in\gamma_-.
\Ee 
Throughout this paper we extend $\phi_f$ for a \textit{negative time} as
\Be\label{negative_t_extension}
\phi_f(t,x) :=   e^{-|t|}\phi_{f_0}(x) \ \  {in } \ (\ref{phi_f}) \ \  {for} \ - \infty<t<0.
\Ee
Note that, for all $t \in \R$, $\phi_f(t,x)$  only depends on $f(t,x,v)$ for the non-negative time $t\geq 0$. 
 
The \textit{characteristics (trajectories)} are determined by the Hamilton ODEs
\Be\label{hamilton_ODE}
\frac{d}{ds} \left[ \begin{matrix}X^f(s;t,x,v)\\ V^f(s;t,x,v)\end{matrix} \right] = \left[ \begin{matrix}V^f(s;t,x,v)\\ 
- \nabla_x \phi_f
(s, X^f(s;t,x,v))\end{matrix} \right]     ,
\Ee
for $- \infty< s ,  t < \infty$ with $(X^f(t;t,x,v), V^f(t;t,x,v)) =  (x,v)$. 

For $(t,x,v) \in \R  \times  \O \times \R^3$, we define \textit{the backward exit time} $\tb^f(t,x,v)$ as   
\Be\label{tb}
\tb^f (t,x,v) := \sup \{s \geq 0 : X^f(\tau;t,x,v) \in \O \ \  for \ all   \ \tau \in (t-s,t) \}.
\Ee
Furthermore, we define $\xb^f (t,x,v) := X^f(t-\tb^f(t,x,v);t,x,v)$ and $\vb^f (t,x,v) := V^f(t-\tb^f(t,x,v);t,x,v)$. We note that all the definitions of $X^f(s;t,x,v)$ and $V^f(s;t,x,v)$ only depend on $f(t,x,v)$ for a non-negative time and our extension (\ref{negative_t_extension}). Hence our definition $\tb^f(t,x,v)$ does not depend on $f(t,x,v)$ for $t<0$.

If we compel $f$ to solve (\ref{Vlasov_eqtn}) and (\ref{inflow_BC}) even for a negative time with $\phi_f$ defined for all $t \in \R$ as (\ref{negative_t_extension}) then 
\Be\notag
f(t,x,v)= g(t-\tb^f(t,x,v), \xb^f(t,x,v), \vb^f(t,x,v)) \ \ {for } \  t\geq \tb^f(t,x,v).
\Ee
From direct computations (see (\ref{nabla_tb}) and (\ref{xbvb_xv})), derivatives have a singularity in general as
\Be\label{deriv_singular}
\nabla_{x} f(t,x,v)   \sim  \nabla_x \xb^f(t,x,v) 
\sim \frac{1}{n(\xb^f(t,x,v)) \cdot \vb^f(t,x,v)}. 
\Ee
Inspired by this observation, we introduce the following notion.
\begin{definition}[Kinetic Weight] Suppose $f(t,x,v)$ is given for $t\geq 0$, $(x,v) \in \bar{\O} \times \R^3$ and satisfies (\ref{mass_conservation_f}). We define $\phi_f (t,x)$ as a unique solution of (\ref{phi_f}) for $t\geq0$ and extend $\phi_f (t,x)$ for a negative time as well by (\ref{negative_t_extension}). Suppose $(X^f(s;t,x,v), V^f(s;t,x,v))$ solves (\ref{hamilton_ODE}) and $\tb^f$ is defined as (\ref{tb}). For $\e>0$, we define a \textit{kinetic weight} as 
\Be\label{alphaweight}\begin{split}
\alpha_{f, \e }(t,x,v) : =& \  
\chi \Big(\frac{t-\tb^f(t,x,v)+\e}{\e}\Big)
|n(\xb^f(t,x,v)) \cdot \vb^f(t,x,v)| \\
&+ \Big[1- \chi \Big(\frac{t-\tb^f(t,x,v) +\e}{\e}\Big)\Big].
\end{split}\Ee
Here we use a smooth function $\chi: \R \rightarrow [0,1]$ satisfying
\Be\label{chi}
\begin{split}
\chi(\tau)  =0,  \     \tau\leq 0, \  {and} \  \ 
\chi(\tau)  = 1    ,  \  \tau\geq 1, \\ 
 \frac{d}{d\tau}\chi(\tau)  \in [0,4] \ \    for \ all  \    \tau \in \R.
\end{split}
\Ee
\end{definition}\hide
\Be\label{weight}
\alpha(t,x,v) = 
\begin{cases} \ 
\mathbf{1}_{\tb(t,x,v) \leq t+1}
|n(\xb(t,x,v)) \cdot \vb(t,x,v)| 
 \ \ \ \text{for}  \ \tb(t,x,v)< \infty
 ,\\
 \ \ \  \ \ \ 
  \ \ \ \ \  \ \ \   \ |\vb(t,x,v)|  \  \ \  \ \ \  \ \ \ \ \ \  \ \ \  \ \ \  \ \ \  \text{for}  \ \tb(t,x,v)= \infty
.
\end{cases}
\Ee\unhide
Note that $\alpha_{f,\e}(0,x,v)$ is solely determined by $\phi_f(t,x)= e^{-|t|} \phi_{f_0}(x)$ for $t\leq 0$, which only depends on $f_0$. Moreover, we have 
\Be\label{alpha_match_nv}
\alpha_{f,\e}(t,x,v) = |n(x) \cdot v| \ \ on \ \gamma_-.
\Ee
For the sake of simplicity, we could drop the superscription $^f$ in $X^f, V^f, \tb^f, \xb^f, \vb^f$ unless this could cause any confusion.

One of the crucial properties of the new kinetic weight in (\ref{alphaweight}) is an invariance under the action of the Vlasov operator: 
\Be\label{alpha_invariant}
\big[\p_t + v\cdot \nabla_x - \nabla_x \phi_f \cdot \nabla_v \big] \alpha_{f,\e}(t,x,v) =0.
\Ee
This is due to the fact that the characteristics solve a deterministic system (\ref{hamilton_ODE}) (See the proof in the appendix). This crucial invariant property under the Vlasov operator is one of the key points in our approach. 

Importantly we note that several different versions of kinetic weight have been used, for convex domains, in \cite{GKTT1} for $E \equiv 0$, and in \cite{Guo_V,Hwang,CK} for $E \neq 0$ but when an extra \textit{favorable sign condition} 
\Be\notag
E\cdot n>0 \ \ on \ \p\O
\Ee
is imposed (recall that $n$ is the outward normal at the boundary). When $E\equiv0$, this weight takes a form of $\tilde{\alpha}$ in (\ref{beta}) which is basically equivalent to $\{\text{dist} (x, \p\O)^2 + |n(x) \cdot v|^2 + |v|^2\text{dist} (x, \p\O)\}^{1/2}$ (also see \cite{GKTT1}). For $E \neq 0$, when the Vlasov operator hits $(\tilde{\alpha})^2$ the outcome contains $|v| \times \text{dist} (x, \p\O)$. This term can be bounded as  
\Be\notag
\frac{1}{|v|} \times (\tilde{\alpha})^2  \ \ with \  a  \ harmful  \ factor \  \frac{1}{|v|}.
\Ee
Clearly we lose control for small velocity. Such \textit{loss of control for small velocity} is indeed a generic difficulty in the study of characteristics with external fields $E$: as the velocity gets smaller, the curvature effect $ |v|^2\text{dist} (x, \p\O)$ vanishes quadratically and therefore even small field can change the characteristics drastically when the velocity is small. While, the favorable sign condition prevents a possible bad behavior of trajectories of small velocity (an interior point can reach the boundary tangentially via characteristics). In fact, in this case, an extra term $(E\cdot n) \times \text{dist}(x,\p\O)$ can be added to $(\tilde{\alpha})^2$ (for example see the definition above Lemma 7 in \cite{CK}), and it leads a control of $|v| \times \text{dist} (x, \p\O)$ by $\frac{1}{E\cdot n } \times |v| (\tilde{\alpha})^2$. Unfortunately such approach fails in our case $E\cdot n|_{\p\O}=0$ as (\ref{phi_BC}). Another problem of this weight is that, $\tilde{\alpha}$ is not exactly invariant under the transport operator (or Vlasov operator) but an extra $\langle v\rangle$-multiplier appears as 
\Be\notag
|\big[\p_t + v\cdot \nabla_x - \nabla_x \phi_f \cdot \nabla_v \big]  \tilde{\alpha}|\lesssim \langle v\rangle \tilde{\alpha}.
\Ee
Conceivably this causes a \textit{super-exponential growth} in $\tilde{\alpha}$-weighted $W^{1,p}$ estimate, which seems quite hard to lead an asymptotical stability of $\mu$.

\hide

After certain modification, similar forms of weight (for example see the definition above Lemma 7 in \cite{CK}) can be successfully adopted in $E\neq 0$ cases only when a favorable \textit{sign condition} of $E\cdot n>0$ ($n$ is the outward normal) is imposed at the boundary.  
We note that $(\tilde{\alpha})^2$ is basically equivalent to $\{\text{dist} (x, \p\O)^2 + |n(x) \cdot v|^2 + |v|^2\text{dist} (x, \p\O)\}$.

Unfortunately there are two major reasons that any of such weighs in \cite{Guo_V, Hwang, GKTT1,CK} is not applicable in our case.

In this \textit{interesting} work the author extends some result of [5] to the case with external fields $E$ for the diffuse reflection boundary condition. A sign condition (1.8) of $n \cdot E$ at the boundary is imposed and it seems crucial in the proof. The proof of Lemma 7 heavily relies on this sign condition: When the Vlasov operator hits the distance function it contains a term $\sim |v||\xi(x)|$ which is $\frac{1}{|v|} \times |v|^2 |\xi(x)|$ or $\frac{|v|}{|\xi(x)|} \times |\xi(x)|^2$. Even though $|v|^2 |\xi(x)|$ and $|\xi(x)|^2$ can be controlled, there is harmful factors $\frac{1}{|v|}$ or $\frac{|v|}{|\xi(x)|}$. Therefore in general we lose the control for small velocity or near the boundary. In other words if the velocity is small then the curvature effect vanishes quadratically and therefore even very small external field can be dominant. The author adds new term in $\beta$ in page 14, the last term of $\beta$, which gives an extra control of $|\xi(x)|$ as long as the sign condition (1.8) holds.

 in order to derive invariance up to some multiplier under the action of Vlasov operator in \cite{Guo_V, Hwang,CK}. Even small field can change the characteristics drastically when the velocity is small. This sign condition prevents a possible bad behavior of trajectories (an interior point can reach the boundary tangentially via characteritics) of \textit{small velocity}. Secondly, $\tilde{\alpha}$ is not exactly invariant under the transport operator (or Vlasov operator) but an extra $\langle v\rangle$-multiplier appears as $|[transport \ operator ] \tilde{\alpha}|\lesssim \langle v\rangle \tilde{\alpha}$. This causes a \textit{super-exponential growth} in $\tilde{\alpha}$-weighted $W^{1,p}$ estimate, which seems quite hard to lead an asymptotical stability of $\mu$. \unhide

In this paper we establish an exponential asymptotical stability of $\mu$ when the potential satisfies zero Neumann boundary condition (\ref{phi_BC}) (therefore no \textit{favorable sign condition} of 
$E\cdot n$ at the boundary). 
For two species problem such a sign condition would not be helpful anymore since a favorable sign of electrons would be bad sign for ions (and vise versa). The new method of this paper, which does not require the favorable sign condition, would serve an effective machinery of two species problems in global well-posedness (\cite{CK2}) and the hydrodynamic limit
of various two species Vlasov-Boltzmann equation (\cite{JM,AS}).

 \subsection{Main Results}

 Construction of a unique global solution and proving its asymptotic stability of VPB in general domains has been a challenging open problem for any boundary condition. 
 The main goal of this paper is to provide the \textit{first} construction of a unique global \textit{strong} solution of VPB system with the diffuse boundary condition when the domain is $C^3$ and \textit{convex.} Moreover an asymptotic stability of the global Maxwellian $\mu$ is studied. 

Here a $C^{3}$ domain means that for any ${p} \in \partial{\Omega}$, there exist sufficiently small $\delta_{1}>0, \delta_{2}>0$, and an one-to-one and onto $C^{3}$-map
	\begin{equation}\label{eta}
	\begin{split}
	\eta_{{p}}:  \{ x_{{ \parallel}} \in \mathbb{R}^{2}: |x_{ \parallel}| < \delta_1  \}  \ &\rightarrow  \ \p\Omega \cap B({p}, \delta_{2}),\\
	x_{{ \parallel}}=(x_{ \parallel,1},x_{\parallel,2} )	 \ &\mapsto \    \eta_{{p}}  (x_{ \parallel,1},x_{\parallel,2} ).
	\end{split}
	\end{equation} 
	A \textit{convex} domain means that there exists $C_\O>0$ such that for all $p \in \p\O$ and  $\eta_p$ and for all $x_\parallel$ in (\ref{eta}) 
\begin{equation}\label{convexity_eta}
\begin{split}
\sum_{i,j=1}^{2} \zeta_{i} \zeta_{j}\p_{i} \p_{j} \eta _{{p}}   ( x_{\parallel }  )\cdot  
n ( x_{\parallel } )
  \leq    - C_{\Omega} |\zeta|^{2}  \ \ 
 { for \ all}   \ \zeta \in \mathbb{R}^{2}.
\end{split}
\end{equation}

We denote  
\Be
\label{weight}
w_\vartheta(v) =  e^{\vartheta|v|^2}.
\Ee

\begin{theorem} 
\label{main_existence}
Assume a bounded open $C^3$ domain $\O \subset\R^3$ is convex as (\ref{convexity_eta}). Let $0< \tilde{\vartheta}< \vartheta\ll1$. Assume a compatibility condition 
\Be\label{compatibility_condition}
f_0 (x,v) = c_\mu \sqrt{\mu(v)} \int_{n(x) \cdot u>0} f_0 (x,u)\sqrt{\mu(u)} \{n(x) \cdot u\} \dd u   \ \ \text{on} \ \gamma_-.
\Ee

There exists a small constant $0< \e_0 \ll 1$ such that for all $0< \e \leq \e_0$ if an initial datum $F_0 = \mu+ \sqrt{\mu}f_0\geq 0$ satisfies (\ref{same_mass_mu}) and 
\Be\label{small_initial_stronger}
 \|w_\vartheta f_0 \|_{L^\infty(\bar{\O} \times \R^3)}< \e,\Ee
and 
\Be\label{W1p_initial}
\begin{split}
 \| w_{\tilde{\vartheta}} \alpha_{f_0, \e }^\beta \nabla_{x,v } f_0 \|_{ {L}^{p } ( {\O} \times \R^3)}
 <\e
\ \
\text{for}  \  \ 3< p < 6, \ \ 
1-\frac{2}{p }
 < \beta<
\frac{2}{3}
,\end{split}
\Ee
and 
\Be\label{nabla_v_f_0_bounded}
 \|   w_{\tilde{\vartheta}}   \nabla_{v } f_0 \|_{ {L}^{3 } ( {\O} \times \R^3)}< \infty,
\Ee
\hide
\Be\label{W1p_initial}
 \| w_{\tilde{\vartheta}} \alpha_{f, \e }^\beta \nabla_{x,v} f_0 \|_{ {L}^{p} ( {\O} \times \R^3)}<\e
\ \ \ \text{for} \ \ 3<p<6, \ 
1-\frac{2}{p}
 < \beta<
\frac{3}{2}
  ,
\Ee\unhide
%
%
then there exists a unique global-in-time solution $(f, \phi_f)$ to (\ref{eqtn_f}), (\ref{phi_f}), (\ref{BC_f}) such that $F(t)= \mu+ \sqrt{\mu} f(t) \geq 0$. Moreover there exists $\lambda_{\infty} > 0$ such that 
\Be\begin{split}\label{main_Linfty}
 \sup_{ t \geq0}e^{\lambda_{\infty} t} \| w_\vartheta f(t)\|_{L^\infty(\bar{\O} \times \R^3)}+ 
 \sup_{ t \geq0}e^{\lambda_{\infty} t} \| \phi_f(t)  \|_{C^{2}(\O)}  \lesssim 1,
\end{split}\Ee
and, for some $C>0$,
\Be\label{W1p_main}
 \| w_{\tilde{\vartheta}} \alpha_{f, \e }^{\beta } \nabla_{x,v} f(t)  \|_{L^{ p} ( {\O} \times \R^3)} 
 \lesssim e^{Ct} \ \ \text{for all } t \geq 0
,
\Ee
and, for $0< \delta= \delta(p,\beta) $,
\Be\label{nabla_v f_31}
\| \nabla_v f (t) \|_{L^3_x (\O) L^{1+\delta }_v (\R^3)} \lesssim_t 1  \ \ \text{for all } \  t\geq 0.
\Ee

Furthermore, if $(f, \phi_f)$ and $(g, \phi_g)$ are both solutions to (\ref{eqtn_f}), (\ref{phi_f}), (\ref{BC_f})
then 
\Be\label{stability_1+}
\| f(t) - g(t) \|_{L^{1+\delta} (\O \times \R^3)} \lesssim_t \| f(0) - g(0) \|_{L^{1+\delta} (\O \times \R^3)} \ \ \text{for all } \  t\geq 0.
\Ee 

\end{theorem}

\hide
\begin{remark} A large class of functional spaces satisfy the condition (\ref{W1p_initial}). Any initial datum $f_0$, whose $1^{\text{st}}$ order weak derivative belongs to $L^p$ for some $p>3$ and small and $\nabla_{x,v} f_0$ decays fast enough as $|v| \rightarrow\infty$, satisfies (\ref{W1p_initial}). \end{remark}\unhide

\begin{remark}A large class of functional spaces satisfy the condition (\ref{W1p_initial}). For example, any initial
datum $f_0$, whose weak derivatives $\nabla_{x,v}f_0$ are small in $L^p(\O \times \R^3)$ for some $p > 3$, and $|\nabla_{x,v} f_0|$ decays fast enough as $|v| \rightarrow \infty$, satisfies (\ref{W1p_initial}). \end{remark}

\hide
Actually we can replace the condition (\ref{W1p_initial}) by a weaker condition: There is $3<p_0<6$ such that 
\Be\label{W1p_initial}
 \|  \langle v\rangle w_{\tilde{\vartheta}}   \nabla_{x,v} f_0 \|_{ {L}^{p_0} ( {\O} \times \R^3)}<\e
,\Ee

Here the kinetic weight $\alpha_{f,\e}(t,x,v)$ is defined in (\ref{weight}). The condition (\ref{W1p_initial}) implies that $\nabla_{x} f_0$ can be singular. 
\unhide

\begin{remark}We note that an exponential growth bound of (\ref{W1p_main}) is a stronger result than the result of \cite{GKTT1}. In \cite{GKTT1}, the upper bound has a super-exponential growth $e^{C t^{2}}$ even in the absence of an external potential. The exponential growth of (\ref{W1p_main}) is crucially used (in an interpolation with an exponential decay of $\phi_f(t)$ in $C^{1, 1-\e}$ for $0< \e \ll1$, Lemma \ref{lemma_interpolation}) to obtain a decay of $\phi_f(t)$ in $C^2$. And this decay-in-time of $\phi_f$ in $C^2$ is one of the crucial ingredient to construct a global-in-time solution.\end{remark}

\hide
\begin{remark}
A natural question is a regularity estimate (\ref{W1p_main}) for $p\geq 6$. If (\ref{W1p_initial}) holds for $6 \leq p < \infty$ then there exists a small time $T>0$ such that 
\Be\label{W1p_main_p_big}
\sup_{0 \leq t \leq T}\big\{ \|  \alpha^\beta \nabla_{x,v} f(t)  \|_{W^{1,p} ( {\O} \times \R^3)} + \big\| \phi(t) \big\|_{W^{3, p}( 
\O )}\big\} \lesssim 1 \ \ \text{for} \ \ 6 \leq p < \infty.
\Ee 
We show this higher regularity result in the appendix.\end{remark}\unhide

\hide
\begin{remark}If (\ref{W1p_initial}) holds for $6 \leq p < \infty$, in addition to the assumptions of Theorem \ref{main_existence}, then for $6 \leq p < \infty$
\Be\label{W1p_main_p_big}
 \|  \alpha^\beta \nabla_{x,v} f(t)  \|_{W^{1,p} ( {\O} \times \R^3)} + \big\| \phi(t) \big\|_{W^{3, p}( 
\O )} \lesssim_t 1 \ \ \text{for all } t \geq 0.
\Ee 
\end{remark}\unhide

\hide\begin{theorem}[Uniqueness]
In addition to (\ref{small_initial}) and (\ref{W1p_initial}), we further assume that for $0< \tilde{\vartheta}< \vartheta$
\Be
\label{extra_con_uniqueness}
\| e^{\tilde{\vartheta}|v|^2} \nabla_v f_0 \|_{L^3 (\O \times \R^3)}< \infty, \ \ 
 \| e^{\tilde{\vartheta}|v|^2} \alpha^\beta \nabla_{x,v} f_0 \|_{ {L}^{p} ( {\O} \times \R^3)}< \e.
\Ee
If $f_0$ and $g_0$ satisfy all above assumptions and if $f$ and $g$ are corresponding solutions then for $0< \delta \ll 1$

Hence the solution, which is constructed in Theorem \ref{main_existence} with an extra condition (\ref{extra_con_uniqueness}) on the initial datum, is unique.
\end{theorem}\unhide
\begin{remark}
As far as the authors know, Theorem 1 provides the \textit{first} unique global-in-time solution to the Vlasov-Poisson-Boltzmann system in bounded domains with physical boundary condition. Moreover, the result of (\ref{main_Linfty}) is the \text{first} proof of asymptotic stability toward the global Maxwellian for  Vlasov-Poisson-Boltzmann system with physical boundary. Also the result of (\ref{W1p_main}) and (\ref{nabla_v f_31}) are the \textit{first} regularity results of Vlasov-Poisson-Boltzmann system in bounded domains with a physical boundary. 
\end{remark}

\begin{remark}We note that the norm in (\ref{W1p_main}) is \textit{nonlinear} in $f$. In the construction of solutions via a sequence argument this nonlinearity causes subtle issue on the convergence (see the proof of Theorem \ref{local_existence}.) 
\end{remark}

 \subsection{Nonlinear-Normed Energy Estimates and an Interpolation}
 
 In the energy-type estimate of $\nabla_{x,v}f$ in $\alpha_{f,\e}^\beta$-weighted $L^p$-norm, the operator $v\cdot \nabla_x$ causes a boundary term to be controlled:
$$
 \int^t_0 \int_{\gamma_-}  
|\alpha_{f,\e}^\beta \nabla_{x,v} f| ^p
\dd s.
$$
Considering the singularity of (\ref{deriv_singular}) and (\ref{alpha_match_nv}), this integrand is integrable if
\Be\label{beta_lower_intro}
\beta> \frac{p-2}{p}    \ \  {so \ that} \ \   |n \cdot v|^{p \beta - p + 1} \in L^1_{loc}(\R^3).
\Ee

On the other hand, in the bulk, we have two terms to be controlled:
\Be
\int^t_0 \iint_{\O \times \R^3}  \nabla_x^2 \phi_f \nabla_v f \alpha_{f,\e}^{p \beta} |\nabla_{x,v} f|^{p-1},\label{We_need_phi_C2} 
\Ee
and
 \Be\label{We_need_phi_K}
\int^t_0\iint_{\O \times \R^3}K  \nabla_{x,v}f  \alpha_{f,\e}^{p \beta} |\nabla_{x,v} f|^{p-1}.
 \Ee 
To handle (\ref{We_need_phi_C2}) we need \textit{a bound of} $\phi_f (t)$ in $C^2_x.$ Unfortunately such estimate is a boarder line case of the well-known Schauder elliptic regularity theory in (\ref{phi_f}) when $\int_{\R^3} f \sqrt{\mu} \dd v$ is merely continuous or bounded. A key observation is that, for $\frac{1}{p}+ \frac{1}{p^*}=1$, we have
%
\Be\label{bound_wp}
\left\|\int_{\R^3}\nabla_x f \sqrt{\mu} \dd v \right\|_{L^p_x(\O)}\lesssim  \sup_{x} \left\|  \frac{ \sqrt{\mu}}{\alpha_{f,\e}^{ \beta}} \right\|_{L^{p^*} (\R^3)}    \left\|  \alpha_{f,\e}^\beta  \nabla_x f \right\|_{L^p(\O \times \R^3)} ,
\Ee
which leads $C^{2,0+}$-bound of $\phi_f$ by the Morrey inequality for $p>3$ as long as 
\Be\label{alpha_integrable}
\alpha_{f,\e}^{- \beta p^*} \in L^1_{loc}(\{ v \in \R^3\}) \ \   for  \ some \ \ \beta p^*> \frac{p-2}{p-1}.
\Ee

We note that $\phi_f(t)$ has an exponential decay in weaker H\"older spaces $C^{1,1-}$ if $f$ decays exponentially in $L^\infty$. As long as the bound (\ref{bound_wp}) grows at most \textit{exponentially in time},  \Be\label{unif_phi}
\textit{ an exponential decay of }  \phi_f(t)   \   in  \ C^2_x
\Ee  can be verified through the following interpolation in H\"older spaces.
\begin{lemma}\label{lemma_interpolation}Assume $\O \subset \R^3$ with a smooth boundary $\p\O$. For $0< D_1<1$, $0< D_2<1$, and $\Lambda_0>0$,
		\Be\begin{split}\label{phi_interpolation}
			\|\nabla^2_x \phi(t )\|_{L^\infty (\O)}
			\lesssim_{\O, D_1, D_2}&
\ 		\	e^{D_1 \Lambda_0t}\|  \phi(t)\|_{C^{1,1-D_1}(\O)}\\
			&
			+ e^{- D_2 \Lambda_0t}\|  \phi(t)\|_{C^{2, D_2}(\O)} \ \ for \ all \  t \geq 0.
		\end{split}\Ee
	\end{lemma}The proof of this lemma is given at Section \ref{sec_global}.

 \subsection{Desingularization via Mixing in Velocity}

To prove (\ref{alpha_integrable}), a major difficulty arises from the \textit{non-local} feature of $\alpha_{f,\e}$, which is determined on the characteristics 
at $t-\tb^f(t,x,v)$. We exam the integrability of $\alpha_{f,\e}^{- \beta p^*}$ by employing a change of variables 
\Be\label{COV_intro}
v  \mapsto (\xb^f(t,x,v) , \tb^f (t,x,v)).
\Ee
  By the direct computation, the Jacobian is equivalent to (see (\ref{dxdtdv}))
\Be\label{Jacobian_intro}\begin{split}
 &
 \frac{|\tb^f|^3}{\alpha_{f,\e} (t,x,v)} \\
  &
   \times  \det \bigg[
\text{Id}_{3 \times 3} + \frac{1}{\tb^f} \int^{t-\tb^f}_{t} \int^s_t \nabla_v X^f(\tau ) \nabla^2_x \phi_f(\tau;X^f(\tau )) \dd \tau \dd s
\bigg],\end{split}
\Ee
where $X^f(\tau)= X^f(\tau;t,x,v)$.

Importantly we note that, for having a \textit{uniform-in-time} positive lower bound of (\ref{Jacobian_intro}), it is necessary to have  
\Be\label{linear_growth_intro}
|\nabla_v X^f(\tau;t,x,v)|\lesssim |t-\tau| \  \  for \ all     \ \tau \leq t.
\Ee
 In the presence of a time-dependent potential we have a \textit{non-autonomous system} from (\ref{hamilton_ODE}):
\Be\begin{split}\label{non-auto}
 \frac{d}{ds} 
\begin{bmatrix}
\nabla_vX^f(s;t,x,v)\\
\nabla_vV^f(s;t,x,v)
\end{bmatrix} 
=    
\begin{bmatrix}
0_{3 \times 3} & \mathrm{Id}_{3 \times 3}\\
- \nabla_x^2 \phi_f(s,X^f(s)) & 0_{3 \times 3} 
\end{bmatrix}\begin{bmatrix}
\nabla_vX^f(s;t,x,v)\\
\nabla_vV^f(s;t,x,v)
\end{bmatrix}.
\end{split}
\Ee
Using an exponential decay from (\ref{unif_phi}), we prove (\ref{linear_growth_intro}) in Lemma \ref{est_X_v} and therefore conclude (\ref{COV_intro}) as long as (\ref{unif_phi}) can be verified.

\hide
be an exponential growth factor in $|\nabla_v X^f(\tau;t,x,v)|$ so that the change of variables (\ref{COV_intro}) is only valid for a small time interval. Using (\ref{unif_phi}) crucially 
which leads to the linear growth bound  and 
\unhide

\hide

Another serious difficulty is that the characteristic ODE is non-autonomous since the field is time-dependent. In order to produce global-in-time uniform bound they are force to use the decay of the field in this ODE. They discretize the time and use induction estimate for the multiplication of matrices. By very delicate estimate they are able to control derivatives of the characteristic ODE in a uniform fashion.

However, such an approach has not been successful in the study of Boltzmann equation due to the non-local
nature of the Boltzmann collision operator, which mixes up different velocities so that their distance towards
?0 can not be controlled.

To overcome this, the PI and collaborators employ a new weight function 
\begin{equation}\label{weight_steady}
\alpha(t,x,v) \sim |v_{\mathbf{b}}(t,x,v) \cdot n(x_{\mathbf{b}}(t,x,v))|,
\end{equation}
where $(\xb,\vb)$ is the characteristics, which starts at $(t,x,v)$, evaluated at the time when it hits the boundary $t-\tb(t,x,v)$. We note that $\tb(t,x,v)$ is the traveling time that the characteristics spends to hit the boundary starting from $(t,x,v)$. Obviously the new weight function is invariant along the characteristics. However since the new weight function is not evaluated ``locally" but somewhere on the characteristics it is very hard to compute the first term in (\ref{bound_wp})
\Be\label{int_alpha}
\int_{|v|\lesssim 1} \frac{\dd v}{\alpha(t,x,v)^{\text{some number}}}.
\Ee 
In order to compute this type of integrations they use a geometric change of variables 
\Be
v \mapsto (\xb, \tb), \ \ \ \text{with the Jacobian} \ \ \det\left(\frac{\p (\xb, \tb) }{\p v}\right) = \frac{|\tb|^3}{|n(\xb) \cdot \vb|}.
\Ee 
\unhide

Applying the change of variables (\ref{COV_intro}), 
$\alpha_{f,\e}$-factor in the Jacobian (\ref{Jacobian_intro}) cancels out the singularity in (\ref{alpha_integrable}) and leads ${\alpha_{f,\e}^{1- \beta p^*}}/{(\tb^f)^3}$ instead. Then we carefully use a lower bound of $\tb^f\gtrsim \frac{|\xb^f-x|}{\max |V^f|}$ and a bound $\alpha_{f,\e}\lesssim \frac{|(x-\xb^f) \cdot n(\xb^f)|}{\tb^f}$ (see (\ref{n_bv_b})) to have
\Be\label{alpha_bounded_intro}
\begin{split}
&\int_{|v| \lesssim 1} \alpha_{f,\e}^{- \beta p^*} \dd v\\
 \lesssim  & \ \int_{\text{boundary}} \frac{|(x- \xb^f) \cdot n(\xb^f)|^{1- \beta p^*}}{|x-\xb^f|^{3- \beta p^*}} \dd \xb^f
+ \text{good terms}< \infty, 
\end{split}
\Ee
which turns to be bounded as long as $\beta p^*<1$. 

This estimate is good to control (\ref{We_need_phi_C2}) but we have to restrict a range of $p$ due to $\frac{1}{|v-u|}$-singularity of $K$ in (\ref{We_need_phi_K}) (see (\ref{k_estimate})). We bound $
K(\frac{1}{\alpha_{f,\e}^\beta } \alpha_{f,\e}^\beta \nabla_{x,v} f)$ by $\big\| \frac{1}{|v-u|} \frac{1}{\alpha_{f,\e}^\beta (u)} \big\|_{L^{p^*}_u} \times  \| \alpha_{f,\e}^\beta \nabla_{x,v} f \|_{L^p_u}$. Viewing $\big\| \frac{1}{|v-u|} \frac{1}{\alpha_{f,\e}^\beta (u)} \big\|_{L^{p^*}_u}$ as $
\big|\frac{1}{|v-\cdot|^{p^*}} * \frac{1}{\alpha_{f,\e}(\cdot)^{\beta p^*}} \big|^{1/{p^*}}$ we apply the Hardy-Littlewood-Sobolev inequality to have
\Be\label{HLS_intro}
  \Big\|\frac{1}{|v- \cdot|^{p^*}} * \frac{1}{\alpha_{f,\e}(\cdot)^{\beta p^*}} \Big\| _{L^{p/{p^*}}} 
\lesssim \Big\| \frac{1}{\alpha_{f,\e}^{\beta}} \Big\|_{L_{loc}^{3/2}} + \text{good terms} .
\Ee
This causes another restriction $\beta < \frac{2}{3}$ from (\ref{alpha_bounded_intro}) and then $p<6$ from (\ref{beta_lower_intro}).

\hide

Clearly this integration (over 2D hypersurface) is bounded as long as the power of weight in (\ref{int_alpha}) is less than $1$. Then they apply Gronwall-type estimate to derive an exponential growth estimate (\ref{bound_wp}) which is good enough to construct a global solution and prove asymptotic stability.

We prove $\int_{\R^3} f \sqrt{\mu} \dd v$ belongs to some H\"older space 

they establish that 
\Be\label{W1p}
\text{a weight} \times\nabla f \in L^p \ \ \  \text{for}   \ \ p>3
\Ee

In order to use $L^p$--$L^\infty$ bootstrap-type estimate in (\ref{p_infty}) this condition is a minimal requirement. Moreover $\phi$ has to be $C^2$ to admit a classical solution of VPB.

By the Morrey's inequality then they obtain a bound of $\int_{\R^3} f \sqrt{\mu} \dd v$ in ${C_x^{0,1-3/p}}$ as
\Be
 \bigg\| \frac{1}{\text{a weight}}\bigg\|_{L^{\frac{p}{p-3}}_v} \times 
 \bigg\|\text{a weight} \times\nabla f (t)\bigg\|_{L^p}
\Ee

\subsection{}

 One of the major difficulty rises up in the regularity estimate of $f$. The phase boundary $\partial\Omega \times \mathbb{R}^{3}$ is \textit{characteristic} but not \textit{uniformly characteristic} at the grazing set $\{(x,v) \in \partial\Omega \times \mathbb{R}^{3} :  n(x) \cdot v=0\}$. Moreover $Lf$ and $\Gamma(f,f)$ is nonlocal. Therefore solutions cannot be localized around the boundary but may be influenced globally in $x$ and $v$. Due to these facts we have an inevitable singularity for the spatial normal derivative of $F$ at the boundary $x\in \partial\Omega$ as 
\begin{equation}\label{boundary_singularity}
\partial_{n}f(t,x,v) \ \sim \ \frac{\Gamma(f,f)}{n(x) \cdot v}  \ \notin \ L^{1}_{loc}.
\end{equation}
In \cite{GKTT1} the PI and collaborators established regularity estimates of $f$ in $W^{1,p}$ for $1\leq p \leq \infty$ with the aid of an weight function 
\Be
 \sqrt{ |n(x) \cdot v| + \text{dist} (x,\O) \times |v|^2 }.
\Ee
Here, $n(x)$ is some extension of outward normal direction of $\p\O$ and $\text{dist} (x,\O)$ stands the distance from $x$ to the boundary. Unfortunately the results of \cite{GKTT1} has \textit{super-exponential growth in time} even without an external field $E\equiv 0$ as  
\Be
\big\| \text{a weight}^{\text{some number}} \nabla_x f(t) \big\|_{L^p} \lesssim e^{t^2}.
\Ee 
  Such a severe growth originates from an unbounded $|v|$-multiplier in a crucial invariant $v\cdot \nabla_x \alpha \sim |v| \times \alpha$. In addition to this ``large velocity issue", an external field causes another singularity at ``small velocity": Without an extra sign condition on the normal derivative of $\phi$ in (\ref{phi}) such as 
  \Be\label{sign_phi}
  \frac{\p \phi}{\p n}\Big|_{\p\O}\geq \delta_0 >0,
  \Ee
  one has $[v\cdot\nabla_x + E \cdot \nabla_v] \alpha \sim \{|v|+ \frac{1}{|v|}\} \alpha$ for properly modified $\alpha$. This small velocity issue is \textit{intrinsic} since even extremely small fields can change the characteristics drastically when the velocity is small.


The domain is bounded open 
{\color{red}[Write the equation of $f$.]}

There might be FOUR main issues in this problem: (1) $\rho \in C^{0,\gamma}$; (2) severe time growth in the velocity lemma, i.e. $ e^{-C|v||t-s|} \alpha(x,v)\leq\alpha(X(s;t,x,v),V(s;t,x,v)) \leq e^{C|v||t-s| } \alpha(x,v)$ for $C\gg1$; (3) the uniqueness; and (4) Uniform-in-time estimate. Recall that 
\begin{equation}\label{alpha_old}
\alpha(x,v) = \sqrt{ |\nabla \xi(x) \cdot v| ^{2 } -2 \xi(x) \{ v\cdot \nabla_{x}^{2} \xi(x) \cdot v\}}.
\end{equation}
(Note that this $\alpha$ is $\sqrt{\alpha \  \text{in} \ \text{\cite{GKTT1}} }$)

The Morrey's inequality said that for $n < p \leq \infty$ and $\O \subset \R^{n}$ with a smooth boundary $\p\O$
\Be\label{morrey}
\| u \|_{C^{0,1-n/p} (\bar{\O})} \lesssim_{n,p,\O} \| u \|_{W^{1,p} (\O)}.
\Ee

From (\ref{morrey}),
\Be\begin{split}\notag
 \Big\| \int_{\R^{3}} f(t,x,v) \sqrt{\mu(v)}\dd v \Big\|_{C^{0,1-n/p} (\bar{\O})}  
 \lesssim_{n,p,\O}  \Big\| \int_{\R^{3}} \nabla_{x} f(t,x,v) \sqrt{\mu(v)}\dd v \Big\|_{L^{p} (\O)}.
\end{split}
\Ee
By Holder inequality 
\Be\begin{split}\notag
 \Big|\int_{\R^{3}} \nabla_{x} f(t,x,v) \sqrt{\mu(v)}\dd v\Big|
 \leq  \Big\|\frac{\sqrt{\mu(\cdot)}}{ \alpha(t,x,\cdot) ^{\beta}}  \Big\|_{L^{\frac{p}{p-1} }(\R^{3})}
 \Big\| \alpha(t,x,\cdot )^{\beta} \nabla_{x} f(t,x,\cdot ) \Big\|_{L^{p} (\R^{3})}
\end{split}\Ee

Assume that we can establish a bound for $\frac{p-2}{p}<\beta< \frac{p-1}{p}$ 
\Be\label{W1p}
 \Big\| \alpha(t )^{\beta} \nabla_{x} f(t ) \Big\|_{L^{p} (\O \times \R^{3}) }
 \lesssim e^{C t } 
 .
\Ee
Note that 
\Be\notag
\frac{p-2}{p-1}<\beta \times \frac{p}{p-1}< 1.
\Ee
Then we may be able to conclude that 
\Be\label{2_growth}
\| \phi_f(t) \|_{C^{2,1-3/p}} \lesssim e^{Ct}.
\Ee
This could be a reasonable estimate from \cite{GKTT1}. However the estimate there is 
\[
\Big\| e^{- L \langle v\rangle t} \alpha(t,x,v)^\beta \nabla _x f(t,x,v)\Big\|_{p} \lesssim_t 1.
\]
The factor of $e^{- L \langle v\rangle t}$ would destroy an exponential growth of 
\[
\Big\|\frac{\sqrt{\mu(\cdot)}}{ e^{- \beta L \langle \cdot\rangle t}   \alpha(t,x,\cdot) ^{\beta}}  \Big\|_{L^{\frac{p}{p-1} }(\R^{3})} \sim e^{t^2},
\]
since by $\tilde{r} = |v| - t/2$
\[
e^{|v|(t - |v|)} = e^{ (\frac{t}{2} + \tilde{r})(\frac{t}{2} - \tilde{r})  } = e^{\frac{t^2}{4}} e^{- (\tilde{r})^2}
.\]

If we can prove following two lemmas then I think there is a chance to remove such $e^{-L \langle v\rangle t}$ in the weight function.\unhide

 \subsection{$L^\infty$-Estimate}
Finally we use an $L^2$-$L^\infty$ bootstrap argument to derive an exponential decay of $f$ in $L^\infty$. The key of this process is to derive a positive lower bound of 
\Be\begin{split}\label{det>1}
\det\Bigg( - (t-s) \text{Id}_{3\times 3}   - \int^s_t \int^\tau_t \frac{\p X^f(\tau^\prime )}{\p v} \nabla_x^2 \phi(\tau^\prime,X^f(\tau^\prime )) \dd \tau^\prime \dd \tau\Bigg),
\end{split}\Ee
except for a small set of $s$. Here $X^f(\tau^\prime)= X^f(\tau^\prime;t,x,v)$. Again as (\ref{Jacobian_intro}) it is crucial to verify (\ref{unif_phi}) and (\ref{linear_growth_intro}) for obtaining a uniform-in-time positive lower bound of (\ref{det>1}). Finally we can close the estimates by proving an exponential growth bound of $\| \alpha_{f,\e}^\beta \nabla_{x,v} f \|_{L^p(\O \times \R^3)}$ from the Gronwall inequality and an exponential decay of $f$ in $L^\infty$ and therefore achieve (\ref{unif_phi}) by Lemma \ref{lemma_interpolation}.

  \subsection{A Priori $L^3_xL^{1+}_v$-Estimate of $\nabla_v f$ and $L^{1+}$-Stability}
For constructing a solution and proving its uniqueness, we need some \textit{stability} estimate of the difference of the solutions $f-g$. The difficulty comes from the term of $\nabla_x \phi_f \cdot \nabla_ v f$. To prove $L^{q}$-stability for $q=1+$ we have, by Sobolev embedding $\nabla_x \phi_{f-g} \in W^{1, q } (\O) \subset L(\O)^{ \frac{3q}{3-q}}$, 
\Be
\begin{split}
\notag
&\iint |\nabla_x \phi_{f-g} \cdot \nabla_v f| |f-g|^{q-1}\\
\lesssim & \   \|\nabla_x \phi_{f-g}\|_{L_x^{ \frac{3q}{3-q}}}  \left\| \| \nabla_v f \|_{L^q_v} \right\|_{L^3_x} \left\| |f-g|^{q-1}\right\|_{L^{\frac{q}{q-1}}_{x,v}  }. 
\end{split}
\Ee
We note that $\nabla_v f$ is bounded at $\gamma_-$, from the boundary condition (\ref{BC_f}). However the equation of $\nabla_v f$ has $\nabla_x f$ as a forcing term. Therefore the key term to bound $\left\| \| \nabla_v f \|_{L^{q}_v} \right\|_{L^3_x}$ for $q=1+  $ is
\Be\notag\begin{split}
&\left\|\left\|\int^t_0 \nabla_x f(s,X(s;t,x,v), V(s;t,x,v)) \dd s \right\|_{L^{1+}_v}\right\|_{L^{3}_x}\\
\lesssim & \ \int^t_0
\| \alpha_{f,\e}^{- \beta} \|_{L^\infty_x L^{\frac{p}{p-1}+}_v}
\|\alpha_{f,\e}^{\beta} \nabla_x f \|_{L^{p}_{x,v}}
 \ \  {for} \ p>3.
\end{split}\Ee
From (\ref{beta_lower_intro}), $ \beta \left(\frac{p}{p-1} +\right) <1$ and therefore $\| \alpha_{f,\e}^{- \beta} \|_{L^\infty_x L^{\frac{p}{p-1}+}_v}$ is bounded from (\ref{alpha_integrable}).

\section{Preliminary and In-Flow Problems}
We define standard notions in (\ref{eqtn_f}). For the hard sphere cross section (\ref{Q}) and the global Maxwellian (\ref{Maxwellian}),
 \begin{equation}\label{L_decomposition}
Lf:= \frac{1}{\sqrt{\mu}}[Q(\mu, \sqrt{\mu}f) + Q(\sqrt{\mu}f,\mu)] : =  \nu (v)f-Kf.
\end{equation}%
Here the collision frequency is defined as
\Be\label{collision_frequency}
\nu (v):= \int_{\mathbb{S}^{2}}\int_{%
\mathbb{R}^{3}}|(v-u) \cdot \omega| \mu(u)\mathrm{d}%
u\mathrm{d}\omega \sim \langle v\rangle
:=\sqrt{1+ |v|^2}
,
\Ee
 and
\begin{equation}\label{def_Kf}
\begin{split}
Kf  &:=K_{2}f -K_{1}f\\
&:= \frac{1}{\sqrt{\mu}} \big[Q_{\mathrm{gain}}(\mu, \sqrt{\mu}f) + Q_{\mathrm{gain}}(  \sqrt{\mu}f,\mu)\big]
 - \frac{1}{\sqrt{\mu}} Q_{\mathrm{loss}} (\sqrt{\mu}f, \mu)\\
&:=\int_{\mathbb{R}^{3}}\mathbf{k}_{2}(v,u)f(u)\mathrm{d}u - \int_{\mathbb{R}^{3}}\mathbf{k}_{1}(v,u)f(u)\mathrm{d}u.
\end{split}
\end{equation}%
It is well-known (See \cite{gl}) that, for some constants $C_{\mathbf{k}_{1}}>0$ and $C_{\mathbf{k}_{2}}>0$,
\begin{equation}\label{k_estimate}
\begin{split}
 \mathbf{k}_{1}(v,u)=  & C_{\mathbf{k}_{1}}  |v-u|   e^{-\frac{|v|^{2} +|u|^{2}}{4}} ,
\\
\mathbf{k}_{2}(v,u) =& C_{\mathbf{k}_{2}} 
|v-u|^{-1} e^{- \frac{|v-u|^2}{8}} 
e^{-   \frac{  | |v|^2- |u|^2   |^2}{8|v-u|^2}}.
\end{split}
\end{equation}%
%
%
\hide

\begin{equation}
Kf(v)=\int_{\mathbb{R}^{3}}\mathbf{k}(u,v)f(u)\mathrm{d}u, \ \ \ \mathbf{k}(u,v)%
=\mathbf{k}_{2}(u,v)-\mathbf{k}_{1}(u,v),  \label{def_Kf}
\end{equation}%

\unhide
The nonlinear operator is defined as 
\begin{equation}\label{Gamma_def}
 \begin{split}
  \Gamma (f_{1},f_{2}) 
 :=& 
  \  \Gamma_{\mathrm{gain}} (f_{1},f_{2})   
 -  \Gamma_{\mathrm{loss}} (f_{1},f_{2})   
 \\
 : =  &
   \ 
 \frac{1}{\sqrt{\mu}}Q_{\mathrm{gain}}(\sqrt{\mu} f_1, \sqrt{\mu} f_2) 
 -  \frac{1}{\sqrt{\mu}}Q_{\mathrm{loss}}(\sqrt{\mu} f_1, \sqrt{\mu} f_2) .
\end{split}
\end{equation}

From now on, in this section, we prove basic estimates of initial-boundary problems of the transport equation in the presence of a time-dependent field $E(t,x)$
\Be\label{transport_E}
\p_t f + v\cdot \nabla_x f + E \cdot \nabla_v f + \psi f = H,
\Ee
where $H=H(t,x,v)$ and $\psi= \psi(t,x,v)\geq 0$. We assume that $E$ is defined for all $t \in \R$. Throughout this section $(X(s;t,x,v),V(s;t,x,v))$ denotes the characteristic which is determined by (\ref{hamilton_ODE}) with replacing $- \nabla_x \phi_f$ by $E$. \begin{lemma}\label{cannot_graze}Assume that $\O$ is convex (\ref{convexity_eta}). Suppose that $\sup_t\| E(t) \|_{C^1_x} < \infty$ and 
	\Be\label{nE=0}
	n(x) \cdot E(t,x) =0 \ \ \text{for } x \in \p\O \ \text{and for all  } t.
	\Ee
	Assume $(t,x,v) \in \R_+ \times \bar{\O} \times \R^3$ and $t+1 \geq \tb(t,x,v)$. If $x \in \p\O$ then we further assume that $n(x) \cdot v > 0$. Then we have 
	\Be
	n(\xb(t,x,v)) \cdot \vb(t,x,v) <0.\label{no_graze}
	\Ee
\end{lemma}
\begin{proof}\hide  \textit{Step 1.}  We claim that for all $(t,x) \in [ 0, \infty) \times \bar{\O}$ as $N \rightarrow \infty$
	\Be\label{conv_NN}
	\mathbf{1}_{\tb(t,x,u)< N}
	\mathbf{1}_{n(\xb(t,x,u)) \cdot \vb(t,x,u) < -\frac{1}{N}}
	\nearrow \mathbf{1}_{\tb(t,x,u)< \infty} \ \ \text{almost every } u \in \R^3.
	\Ee
	First we prove that, for fixed $N \in \mathbb{N}$, as $M \rightarrow \infty$
	\Be\label{conv_NM}
	\mathbf{1}_{\tb(t,x,u)< N}
	\mathbf{1}_{n(\xb(t,x,u)) \cdot \vb(t,x,u) < -\frac{1}{M}}
	\nearrow \mathbf{1}_{\tb(t,x,u)< N} \ \ \text{almost every } u \in \R^3.
	\Ee
	Since $\mathbf{1}_{\tb(t,x,u)< N}$ converges to $\mathbf{1}_{\tb(t,x,u)< \infty}$ as $N \rightarrow \infty$ we can apply Cantor's diagonal argument to conclude (\ref{conv_NN}) from (\ref{conv_NM}).\unhide
	

	\textit{Step 1.} Note that we can locally parametrize the trajectory (see Lemma 15 in \cite{GKTT1} or \cite{KL2} for details). We consider local parametrization (\ref{eta}). We drop the subscript $p$ for the sake of simplicity. If $X(s;t,x,v)$ is near the boundary then we can define $(X_n, X_\parallel)$ to satisfy 
	\Be\label{X_local}
	X(s;t,x,v)  =   \eta (X_\parallel (s;t,x,v)) + X_n(s;t,x,v) [- n(X_\parallel(s;t,x,v))].
	\Ee
	\hide $\eta$ near $\p\O$ such that 
	\Be
	\eta : (x,y,z)\in B(0,r_{1})\cap \R^{3}_{+} \mapsto B(p,r_{2})\cap\O,
	\Ee
	where $\eta(0)=p\in\p\O$ and $\eta(x,y,0) \in \p\O$ for some $r_{1}, r_{2} > 0$. Then for $X(s;t,x,v)$, there exist a unique $x_{*}\in B(p,r_{2})\cap \p\O$ such that
	\Be \label{X def}
	|X(s;t,x,v)-x_{*}| \leq \sup_{x\in B(p,r_{2})\cap \p\O} |X(s;t,x,v)-x|,
	\Ee 
	since $\eta$ is bijective. We define 
	\[
	(X_{\parallel},0) = \eta^{-1}(x_{*}) \quad \text{and} \quad  X_{n} := |X(s;t,x,v)-x_{*}|.
	\]\unhide
	For the normal velocity 
	we define
	\Be\label{def_V_n}
	V_{n}(s;t,x,v) := V(s;t,x,v)\cdot [-n(X_{\parallel}(s;t,x,v))].
	\Ee
	We define $V_{\parallel}$ tangential to the level set $  \big( \eta(X_{\parallel}) + X_{n}(-n(X_{\parallel})) \big)$ for fixed $X_{n}$. Note that 
	\[
	\frac{\p   \big( \eta(x_{\parallel} ) + x_{n}(-n(x_{\parallel})) \big) }{\p {x_{\parallel, i}}}  \perp n(x_{\parallel}) \ \  \text{for} \ i=1,2.
	\]
	We define $(V_{\parallel,1}, V_{\parallel,2})$ as 
	\Be\label{def_V_parallel}
	V_{\parallel, i} :=\Big( V - V_{n}[-n(X_{\parallel})]\Big) \cdot   \Big( 
	\frac{\p   \eta(X_{\parallel} )}{\p {x_{\parallel, i}} } + X_{n}\Big[- \frac{\p n(X_{\parallel})}{\p x_{\parallel, i}}\Big] \Big)  .
	\Ee
	\hide \Be\notag
	\begin{split} 
		&\sum_{i=1,2} V_{\parallel,i} \p_{i} \big( \eta(X_{\parallel},0) + X_{n}(-n(X_{\parallel})) \big)  \\
		&= \nabla_{\parallel} \big( \eta(X_{\parallel},0) + X_{n}(-n(X_{\parallel})) \big) V_{\parallel} = V - V_{n}(-n(X_{\parallel}))  .
	\end{split}
	\Ee\unhide
	Therefore we obtain
	\Be \label{V_local}
	V(s;t,x,u)   = V_n [- n(X_\parallel)] + V_\parallel \cdot \nabla_{x_\parallel} \eta (X_\parallel ) 
	- X_n V_\parallel \cdot \nabla_{x_\parallel} n (X_\parallel).
	\Ee  
	
	Directly we have 
	\Be \begin{split}\notag
		\dot{X}(s;t,x,u) &=\dot{X}_{\parallel} \cdot \nabla_{x_\parallel}\eta (X_\parallel) + \dot{X}_n [- n(X_\parallel)] - X_{n}\dot{X}_{\parallel}  \cdot \nabla_{x_\parallel} n(X_{\parallel}) .
	\end{split}\Ee
	Comparing coefficients of normal and tangential components, we obtain that 
	\Be\label{dot_Xn_Vn}
	\dot{X}_{n}(s;t,x,v) = V_{n}(s;t,x,v) , \ \  \dot{X}_{\parallel}(s;t,x,v) = V_{\parallel}(s;t,x,v).
	\Ee
	
	On the other hand, from (\ref{V_local}),
	\Be \begin{split} \label{Vdotn}
		\dot{V} (s) &=  \dot{V}_{n} [-n(X_{\parallel})] - V_{n} \nabla_{x_\parallel} n(X_{\parallel})\dot{X}_{\parallel} + V_{\parallel}\cdot\nabla^{2}_{x_\parallel}\eta(X_{\parallel}) \dot{X}_{\parallel}   \\
		&\quad+ \dot{V}_{\parallel} \cdot \nabla_{x_\parallel}\eta(X_{\parallel})  - \dot{X}_{n}\nabla_{x_\parallel} n(X_{\parallel})V_{\parallel} - X_{n}\nabla_{x_\parallel} n(X_{\parallel})\dot{V}_{\parallel} \\
		&\quad- X_{n} V_{\parallel}\cdot\nabla_{x_\parallel}^{2}n(X_{\parallel})\dot{X}_{\parallel}. 
	\end{split}\Ee
	From $(\ref{Vdotn})\cdot [-n(X_{\parallel})]$, (\ref{dot_Xn_Vn}), and $\dot{V}=E$, we obtain that 
	\Be \begin{split}\label{hamilton_ODE_perp}
		\dot{V}_n (s) 
		&=  [V_\parallel (s)\cdot \nabla^2 \eta (X_\parallel(s)) \cdot V_\parallel(s) ] \cdot n(X_\parallel(s))  
		\\
		&\quad +  E (s , X (s ) ) \cdot [-n(X_\parallel(s)) ]  \\
		&\quad- X_n (s) [V_\parallel(s) \cdot \nabla^2 n (X_\parallel(s)) \cdot V_\parallel(s)]  \cdot n(X_\parallel(s)) .
	\end{split}\Ee
	
	\vspace{4pt}
	
	\textit{Step 2.} We prove (\ref{no_graze}) by the contradiction argument. Assume we choose $(t,x,v)$ satisfying the assumptions of Lemma \ref{cannot_graze}. Let us assume
	\Be\label{initial_00}
	X_n (t-\tb;t,x,v)  +V_n (t-\tb;t,x,v) =0.
	\Ee
	
	First we choose $0<\e \ll 1$ such that $X_n(s;t,x,v) \ll 1$ and 
	\Be\label{Vn_positive}
	V_n (s;t,x,v) \geq0 \ \ \text{for} \  t- \tb(t,x,v)<s<t-\tb(t,x,v) + \e.
	\Ee
	The sole case that we cannot choose such $\e>0$ is when there exists $0< \delta\ll1$ such that $V_n(s;t,x,v)<0$ for all $s \in ( t-\tb(t,x,v), t-\tb(t,x,v) + \delta)$. But from (\ref{dot_Xn_Vn}) for $s \in ( t-\tb(t,x,v), t-\tb(t,x,v) + \delta)$
	\Be\begin{split}\notag
	0 \leq & \  X_n(s;t,x,v) \\
	  =   & \ X_n(t-\tb(t,x,v);t,x,v)  +  \int^s_{t-\tb(t,x,v)} V_n (\tau; t,x,v) \dd \tau <  0.
	\end{split}\Ee
	
	Now with $\e>0$ in (\ref{Vn_positive}), temporarily we define that $t_* := t-\tb(t,x,v) + \e$, $x_* = X(t-\tb(t,x,v) + \e; t,x,v),$ and $v_* = V(t-\tb(t,x,v) + \e; t,x,v)$. Then $(X_n(s;t,x,v), X_\parallel (s;t,x,v)) = (X_n(s; t_*, x_*, v_*), X_\parallel (s; t_*, x_*, v_*))$ and \\$(V_n(s;t,x,v), V_\parallel (s;t,x,v)) = (V_n(s; t_*, x_*, v_*), V_\parallel (s; t_*, x_*, v_*))$.
	\hide
	From (\ref{hamilton_ODE}), for $(X_n(s), X_\parallel (s)) = (X_n(s; t_*, x_*, u_*), X_\parallel (s; t_*, x_*, u_*))$ and $(V_n(s), V_\parallel (s)) = (V_n(s; t_*, x_*, u_*), V_\parallel (s; t_*, x_*, u_*))$, 
	\Be \begin{split}\label{hamilton_ODE_perp}
		\dot{X}_n (s)  =& V_n (s),\\
		\dot{V}_n (s)  =& 
		[V_\parallel (s)\cdot \nabla^2 \eta (X_\parallel(s)) \cdot V_\parallel(s) ] \cdot n(X_\parallel(s))  
		+   \nabla\phi (s , X (s ) ) \cdot n(X_\parallel(s))  \\
		&+   X_n (s) [V_\parallel(s) \cdot \nabla^2 n (X_\parallel(s)) \cdot V_\parallel(s)]  \cdot n(X_\parallel(s)) .
	\end{split}\Ee\unhide
	
	Now we consider the RHS of (\ref{hamilton_ODE_perp}). From (\ref{convexity_eta}), the first term $[V_\parallel(s) \cdot \nabla^2 \eta (X_\parallel(s)) \cdot V_\parallel(s) ] \cdot n(X_\parallel(s))\leq 0$. By an expansion and (\ref{nE=0}) we can bound the second term 
	\Be\begin{split}\label{expansion_E}
		&E (s , X(s )) \cdot n(X_\parallel(s ) )\\
		=&   \ E (s , X_n(s ), X_\parallel(s ) ) \cdot n(X_\parallel (s )) \\
		=& \  E (s , 0, X_\parallel(s ) ) \cdot n(X_\parallel (s )) 
		+ \| E (s) \|_{C_x^1}  O( |X_n(s )| )\\
		=&  \ \| E (s) \|_{C_x^1}  O( |X_n(s )| ).
	\end{split}\Ee
	From (\ref{hamilton_ODE}) and assumptions of Lemma \ref{cannot_graze},
	$$|V_\parallel (s;t,x,v )|\leq |v| + \tb(t,x,v) \| E \|_\infty  \leq |v| +  (1+t) \| E \|_\infty.$$ 
	Combining the above results with (\ref{hamilton_ODE_perp}), we conclude that 
	\Be\notag
	\dot{V}_n(s;t_*,x_*,v_*) \lesssim  ( |v| + (1+ t) \| E\|_\infty  )^2X_n(s;t_*,x_*,v_*) ,
	\Ee
	and hence from (\ref{dot_Xn_Vn}) for $t-\tb(t,x,v)\leq s \leq t_*$
	\Be\label{ODE_X+V}
	\begin{split}
		&\frac{d}{ds} [X_n (s;t_*,x_*,v_* )  +V_n (s ;t_*,x_*,v_*) ]\\
		\lesssim & \ ( |v| + (1+ t) \| E\|_\infty  )^2  [X_n (s;t_*,x_*,v_* )  +V_n (s;t_*,x_*,v_* ) ].\end{split}
	\Ee
	By the Gronwall inequality and (\ref{initial_00}), for $t-\tb(t,x,v)\leq s \leq t_*$
	\Be \begin{split}\notag
		& [X_n (s;t_*,x_*,v_*)  +V_n (s;t_*,x_*,v_*) ]  \\
		\lesssim & \   [X_n (t-\tb(t,x,u))  +V_n (t-\tb(t,x,u)) ] e^{C \e ( |v| + (1+ t) \| E\|_\infty  )^2) }\\
		=&  \ 0.
	\end{split}\Ee 
	
	From (\ref{Vn_positive}) we conclude that $X_n (s;t,x,v) \equiv 0$ and $V_n (s;t,x,v) \equiv 0$ for all $s \in [t-\tb(t,x,u), t-\tb(t,x,u) + \e]$. We can continue this argument successively to deduce that $X_n (s;t,x,v) \equiv 0$ and $V_n (s;t,x,v) \equiv 0$ for all $s \in [t-\tb(t,x,v), t]$. Therefore $x_n =0 = v_n$ which implies $x \in \p\O$ and $n(x) \cdot v =0$. This is a contradiction since we chose $n(x) \cdot v>0$ if $x \in \p\O$.\end{proof}

\hide\begin{lemma} If $n(\xb(t,x,v)) \cdot \vb(t,x,v) \neq 0$ then $(\tb,\xb,\vb)$ is differentiable and 
	\Be\begin{split}\label{computation_tb_x}
		\frac{\p\tb}{\p x_i}  =& \  \frac{1}{n(\xb) \cdot \vb}n(\xb) \cdot  \left[
		e_i + \int^{t-\tb}_t \int^s_t \Big(\frac{\p X(\tau )}{\p x_i} \cdot \nabla\Big) E(\tau, X(\tau )) \dd \tau \dd s 
		\right]  ,\\
		\frac{ \p\xb}{\p x_i} = & \  e_i - \frac{\p \tb}{\p x_i} \vb + \int^{t-\tb}_{t} \int^s_t    \Big(\frac{\p X(\tau )}{\p x_i} \cdot \nabla\Big) E(\tau, X(\tau ))  \dd \tau \dd s,\\
		\frac{\p \vb}{\p x_i} = & \ - \frac{\p \tb}{\p x_i} E(t-\tb, \xb) + \int^{t-\tb}_t   \Big(\frac{\p X(\tau )}{\p x_i} \cdot \nabla\Big) E(\tau, X(\tau )) \dd \tau,\\
		\frac{\p\tb}{\p v_i}  =& \  \frac{1}{n(\xb) \cdot \vb}n(\xb) \cdot  \left[
		e_i + \int^{t-\tb}_t \int^s_t \Big(\frac{\p X(\tau )}{\p v_i} \cdot \nabla\Big) E(\tau, X(\tau )) \dd \tau \dd s 
		\right]  ,\\ 
		\frac{\p \xb}{\p v_i} = & \ - \tb e_i - \frac{
			\p \tb}{\p v_i} \vb + \int^{t-\tb}_{t} \int^s_t    \Big(\frac{\p X(\tau )}{\p v_i} \cdot \nabla\Big) E(\tau, X(\tau ))  \dd \tau \dd s ,\\
		\frac{\p \vb}{\p v_i} = & \ e_i - \frac{\p \tb}{\p v_i} E(t-\tb, \xb) + \int^{t-\tb}_t
		\Big(\frac{\p X(\tau )}{\p v_i} \cdot \nabla\Big) E(\tau, X(\tau ))  \dd \tau.
	\end{split} \Ee
	where $(X(\tau),V(\tau))= (X(\tau;t,x,v), V(\tau;t,x,v))$.
\end{lemma}
\begin{proof}
	These equalities can be derived from direct computations and the implicit function theorem. See \cite{KL1} for details. 
\end{proof}\unhide

\begin{lemma} \label{lem COV}
	For fixed $t$, a map 
	\Be\label{map_xv_boundary}
	(x,v) \in \O \times \R^3 \mapsto (t-\tb(t,x,v), \xb(t,x,v),\vb(t,x,v)) \in \R \times  \gamma_-
	\Ee
	is one-to-one and 
	\Be\label{Jac_xv_boundary}
	\begin{split}
	&\left|\det\left( \frac{\p (t-\tb(t,x,v),  x_{\mathbf{b} } (t,x,v), \vb(t,x,v))}{\p (x,v)} \right)
	\right|
	\\
	=  & \ \frac{1}{|\vb(t,x,v) \cdot n(\xb(t,x,v))|}.
	\end{split}\Ee
Also, a map
	\Be\label{map_boundary_boundary}
	(t,x,v) \in \R \times \gamma_+ \mapsto (t-\tb(t,x,v), \xb(t,x,v), \vb(t,x,v))  \in \R \times \gamma_-
	\Ee
	is one-to-one and 
	\Be\begin{split}
	&\left|	\det \left( \frac{\p (t-\tb(t,x,v), \xb(t,x,v), \vb(t,x,v))}{\p (t,x ,v)}\right)\right|\\
	= & \ \frac{|n(x) \cdot v|}{|n(\xb(t,x,v)) \cdot \vb(t,x,v)|}.\label{Jac_boundary_boundary}
	\end{split}
	\Ee\hide

	For fixed $s>0$ so that $t-\tb (t,x, v ) < s < t$,
	\Be \label{mapC+}
	\begin{split}
		( t, x , v )  \in   [0,T]\times\gamma_{+} &\mapsto (X(s;t, x, v ) , V(s; t, x, v ) ) \in  \Omega\times\R^{3}  
	\end{split}
	\Ee 
	is injective. 
	For fixed $s>0$ so that $t < s < t + \tf (t,x, v ) $,  
	\Be \label{mapC-}
	\begin{split} 
		( t, {x} , v )  \in [0,T]\times\gamma_{-} \mapsto  (X(s;t,x , v ) , V(s; t, x, v ) )\in\Omega\times\R^{3} 
	\end{split}
	\Ee
	is also injective. For both maps (\ref{mapC+}) and (\ref{mapC-}),
	\Be\label{Jacobian_C+}
	\det\left(\frac{\p (X(s;t, x, v ) , V(s; t, x, v ) ) }{\p  ( t, x , v ) }\right) = |
	n(x)	 \cdot v
	|.
	\Ee
	
	Moreover 
	\Be \label{map+}
	\begin{split}
		& ( t, s, x , v )
		\in [0,T] \times \{ - \min \{t, \tb(t,x,v)\}< s<0 \} \times \gamma_+\\
		\mapsto 	&	
		\big( t+s, X(t+s;t, x, v ) , V(t+s; t, x, v ) \big)  
		\in  [0,T]\times\Omega\times\R^{3},
	\end{split}
	\Ee 
	and 
	\Be \label{map-} 
	\begin{split}
		& ( t, s,x , v )  
		\in 	    [0,T] \times\{ 0<s < \min \{\tf(t,x,v),T-t\} \} \times\gamma_{-}\\
		\mapsto  & \big(t+s, X(t+s;t,x, v ), V(t+s; t,x, v ) \big)
		\in [0,T] \times \O \times \R^3
	\end{split}\Ee
	are both injective, and for both maps (\ref{map+}) and (\ref{map-})
	\Be\label{Jac+}
	\det \left(\frac{\p (t+s, X(t+s;t,x,v), V(t+s;t,x,v ))}{\p (t,s,x,v)} \right)=   |n(x) \cdot 
	v
	|.
	\Ee
	\hide

	For $t\in[0,T]$ and $(\eta(\bar{x}), v ) \in \gamma_+$, let us fix $s>0$ so that $0 < s < \tb(t,\eta(\bar{x}),v)$. Then map
	\Be \label{map+}
	\begin{split}
		\mathcal{A}_{+} &: [0,T]^{2}\times\gamma_{+} \mapsto [0,T]\times\Omega\times\R^{3},  \\	
		\mathcal{A}_{+}( t, s, \bar{x} , v ) &= \big( t-s, X(t-s;t, \eta(\bar{x}), v ) , V(t-s; t, \eta(\bar{x}), v ) \big)  
	\end{split}
	\Ee 
	is well-defined. If we denote Jacobian of $\mathcal{A}_{+}$ as $J_{\mathcal{A}_{+}}$, then
	\Be \label{Jac+}
	\det J_{\mathcal{A}_{+}} = - v \cdot (\p_1 \eta(\bar{x}) \times \p_2 \eta(\bar{x}) ).
	\Ee
	Similarly, for $t\in[0,T]$ and $(\eta(\bar{x}), v) \in \gamma_{-}$, we choose $s>0$ so that $0 < s < \tf (t,\eta(\bar{x}), v ) $. For well-defined map
	
	we obtain
	\Be \label{Jac-}
	\det J_{\mathcal{A}_{-}} = v \cdot (\p_1 \eta(\bar{x}) \times \p_2 \eta(\bar{x}) ).  \\
	\Ee\unhide\unhide
\end{lemma}
\begin{proof}Both maps (\ref{map_xv_boundary}) and (\ref{map_boundary_boundary}) are clearly one-to-one since the characteristics solve (\ref{hamilton_ODE}), and hence are deterministic. 
	From (\ref{eta}) we denote $\eta (x_{\mathbf{b},\parallel}(t,x,v))=\eta(x_{\mathbf{b},1} (t,x,v),x_{\mathbf{b},2} (t,x,v)  )= x_{\mathbf{b}} (t,x,v)$. We use the notations $\p_s X = \frac{\p X(s;t,x,v)}{\p s}$, $\p_t X = \frac{\p X(s;t,x,v)}{\p t}$, $\p_{x_i} X = \frac{\p X(s;t,x,v)}{\p x_i}$, and $\p_{v_i} X = \frac{\p X(s;t,x,v)}{\p v_i}$. Then 
	\Be\label{nabla_tb}
	\begin{split}
		\p_{x_i} \tb(t,x,v) &= \frac{\p_{x_i} X(t-\tb(t,x,v) ;t,x,v) \cdot n(\xb(t,x,v) )
		}{ \vb(t,x,v)  \cdot n(\xb(t,x,v) )},\\
		\p_{v_i} \tb(t,x,v) &= \frac{\p_{v_i} X(t-\tb(t,x,v) ;t,x,v) \cdot n(\xb(t,x,v) ) }{ \vb(t,x,v)  \cdot n(\xb(t,x,v) )}.
	\end{split}
	\Ee
	From (\ref{hamilton_ODE})
	\Be
	\begin{split}\label{xbvb_xv}
		\nabla_{x,v} \xb(t,x,v) 
		=&     \nabla_{x,v} [X(t-\tb;t,x,v)] \\
		=&  \nabla_{x,v} \tb \vb + \nabla_{x,v} X(t-\tb;t,x,v),\\
		\nabla_{x,v} \vb(t,x,v) 
		=&  \nabla_{x,v} [V(t-\tb;t,x,v)] \\
		=& \ \nabla_{x,v} \tb E(t-\tb,\xb) + \nabla_{x,v} V(t-\tb;t,x,v).
	\end{split}
	\Ee

	From (\ref{xbvb_xv})
	\Be
	\begin{split}\label{tbxbvb_xv}
		& \frac{\p (t-\tb(t,x,v),   x_{\mathbf{b},\parallel} (t,x,v), \vb(t,x,v))}{\p (x,v)}\\
		=& \left[\begin{matrix}
			- \nabla_x \tb  & -\nabla_v \tb \\
			( - \nabla_x \tb\vb 
			- \nabla_x X)  \cdot \nabla (\eta^{-1})_1& 
			(- \nabla_v \tb\vb 
			- \nabla_v X) \cdot \nabla (\eta^{-1})_1\\
			(- \nabla_x \tb\vb 
			- \nabla_x X) \cdot \nabla (\eta^{-1})_2& 
			(- \nabla_v \tb\vb 
			- \nabla_v X) \cdot \nabla (\eta^{-1})_2\\
			- \nabla_x \tb E   - \nabla_x V 
			& - \nabla_v \tb E  - \nabla_v V 
		\end{matrix}\right],
	\end{split}
	\Ee
	where $(\tb,X,V)= (\tb(t,x,v),X(t-\tb(t,x,v);t,x,v),V(t-\tb(t,x,v);t,x,v))$ and $E= E(t-\tb,\xb)$. 
	
	We multiply $-\vb \cdot \nabla(\eta^{-1})_1$ to the first row and add it to the second row. Then we multiply $-\vb \cdot \nabla(\eta^{-1})_2$ to the first row and add it to third second row. Next we multiply $E_i$ to the first row and add it to the $(i+3)$-row (elementary row operations). Then we obtain a matrix, with the same determinant as (\ref{tbxbvb_xv}),
	\Be\notag\label{tbxbvb_xv_1}
	\left[
	\begin{matrix}
		- \nabla_x \tb & - \nabla_v \tb\\
		- \nabla_x X  \cdot \nabla (\eta^{-1})_1& - \nabla_v X  \cdot \nabla (\eta^{-1})_1\\
		- \nabla_x X  \cdot \nabla (\eta^{-1})_2& - \nabla_v X  \cdot \nabla (\eta^{-1})_2\\
		- \nabla_x V & - \nabla_v V
	\end{matrix}
	\right].
	\Ee
	Then using (\ref{nabla_tb}) we get another matrix with the same determinant 
	\Be\notag\label{tbxbvb_xv_1}\begin{split}
		&\det\left[
		\begin{matrix}
			-  \frac{n(\xb)}{\vb \cdot n(\xb)} \cdot \nabla_x X  &- \frac{n(\xb)}{\vb \cdot n(\xb)}\cdot   \nabla_v X 
			\\
			- \nabla_x X  \cdot \nabla (\eta^{-1})_1& - \nabla_v X  \cdot \nabla (\eta^{-1})_1\\
			- \nabla_x X  \cdot \nabla (\eta^{-1})_2& - \nabla_v X  \cdot \nabla (\eta^{-1})_2\\
			- \nabla_x V & - \nabla_v V
		\end{matrix}
		\right]\\
		= &	
		- 
		\det\left[\begin{matrix}
			\nabla_x X & \nabla_x V\\
			\nabla_v X & \nabla_v V
		\end{matrix}\right]\times \det
		\left[
		\begin{matrix}
			\frac{n(\xb)}{\vb\cdot n(\xb)} & \nabla \eta^{-1} & 0
			\\
			0 & 0  & \text{Id}_{3\times 3}
		\end{matrix}
		\right]
		.
	\end{split}\Ee
	By the Liouville theorem $\det\left[\begin{matrix}
	\nabla_x X & \nabla_x V\\
	\nabla_v X & \nabla_v V
	\end{matrix}\right]=1$. Note that we can always assume $\p_{{x_\parallel,1}} \eta \cdot \p_{{x_\parallel,2}} \eta =0$ by reparametrization in (\ref{eta}). Then 
	\Be\begin{split}\notag
		& \Big| \det
		\left[
		\begin{matrix}
			\frac{n(\xb)}{\vb\cdot n(\xb)} & \nabla \eta^{-1} & 0
			\\
			0 & 0  & \text{Id}_{3\times 3}
		\end{matrix}
		\right] \Big|
		=  \Big| \det
		\left[
		\begin{matrix}\frac{n(\xb)}{\vb\cdot n(\xb)} & \nabla \eta^{-1}\end{matrix}
		\right] \Big| \\
		= & \frac{1}{ |\vb\cdot n(\xb)| |\p_1 \eta(\xb)  \times \p_2 \eta(\xb)  | }
		.
	\end{split}	\Ee
	Therefore we conclude (\ref{Jac_xv_boundary}).
	
	Next we consider (\ref{Jac_boundary_boundary}). By direct computations and (\ref{xbvb_xv}), (\ref{transport_tb}) 
		\Be\begin{aligned}\notag
			&\frac{\p (t-\tb(t,\eta(x_\parallel),v), x_{\mathbf{b},\parallel}(t,\eta(x_\parallel),v), \vb(t,\eta(x_\parallel),v))}{\p (t,x_\parallel,v)}\\
			=&
			\left[\begin{matrix} 
				1- \frac{\p \tb}{\p t} & - \nabla_{x_\parallel} \eta (x_\parallel) \cdot \nabla_x \tb&
				\\
				\left((1- \frac{\p \tb}{\p t}) \vb + \p_t X  \right)\cdot \nabla (\eta^{-1})_1&  \left(\nabla_{x_\parallel} \eta  \cdot 
				(\nabla_{x_\parallel} \tb \vb + \nabla_{x_\parallel}X )
				\right) \cdot \nabla (\eta^{-1})_1 &\\
				\left((1- \frac{\p \tb}{\p t}) \vb + \p_t X\right) \cdot \nabla (\eta^{-1})_2 &\left(  \nabla_{x_\parallel} \eta   \cdot(\nabla_{x_\parallel} \tb \vb + \nabla_{x_\parallel}X ) \right)  \cdot \nabla (\eta^{-1})_2&  \\
				(1- \frac{\p \tb}{\p t}) E(t-\tb, \xb)
				+ \p_t V
				& \nabla_{x_\parallel} \eta (x_\parallel) \cdot
				(\nabla_x \tb  E+ \nabla_x V)
				&
			\end{matrix}\right. \\
			& \ \ \ \   \ \
			\left. \begin{matrix} 
		 - \nabla_{v}   \tb \\
		  (\nabla_{v} \tb \vb + \nabla_{v}X ) \cdot \nabla (\eta^{-1})_1\\
		  (\nabla_{v} \tb \vb + \nabla_{v}X )   \cdot \nabla (\eta^{-1})_2\\
		  \nabla_v \tb E + \nabla_v V
			\end{matrix}\right].
		\end{aligned}\Ee

	By elementary row operations, we obtain 
	\Be\begin{split}\label{Jac_bdry_bdry_1}
		\left[\begin{matrix} 
			1- \frac{\p \tb}{\p t} & - \nabla_{x_\parallel} \eta (x_\parallel) \cdot \nabla_x \tb& - \nabla_{v}   \tb
			\\
			\p_t X   \cdot \nabla (\eta^{-1})_1&  \left(\nabla_{x_\parallel} \eta  \cdot 
			\nabla_{x_\parallel}X  
			\right) \cdot \nabla (\eta^{-1})_1 &  \nabla_{v}X   \cdot \nabla (\eta^{-1})_1\\
			\p_t X  \cdot \nabla (\eta^{-1})_2 &\left(  \nabla_{x_\parallel} \eta   \cdot  \nabla_{x_\parallel}X   \right)  \cdot \nabla (\eta^{-1})_2&    \nabla_{v}X     \cdot \nabla (\eta^{-1})_2\\
			\p_t V  & \nabla_{x_\parallel} \eta (x_\parallel) \cdot
			\nabla_x V 
			&  \nabla_v V
		\end{matrix}\right].
	\end{split}\Ee
	\hide

	The maps are injective since the characteristics are a solution of (\ref{hamilton_ODE}).

	For $x \in \p\O$ we have $x= \eta(x_\parallel)$ locally from (\ref{eta}). Then we compute (\ref{Jac+}),
	\Be
	\begin{split}\label{XV_txv}
		& \frac{\p ( X( s;t, \eta(x_\parallel), v),  V( s;t, \eta(x_\parallel), v) )}{\p (t, x_\parallel, v)}\\
		= &  \begin{bmatrix}
			\p_t X( s;t, \eta(x_\parallel), v) & \nabla_{x_\parallel} X( s;t, \eta(x_\parallel), v) & \nabla_v X( s;t, \eta(x_\parallel), v)\\
			\p_t V( s;t, \eta(x_\parallel), v)& \nabla_{x_\parallel} V( s;t, \eta(x_\parallel), v) & \nabla_v V( s;t, \eta(x_\parallel), v)
		\end{bmatrix}\\
		= &  
		\begin{bmatrix}
			\p_t X( s;t, \eta(x_\parallel), v)& \nabla_{x_\parallel} \eta (x_\parallel) \cdot \nabla_x X( s;t, \eta(x_\parallel), v)  & \nabla_v X( s;t, \eta(x_\parallel), v)\\
			\p_t V( s;t, \eta(x_\parallel), v) &  \nabla_{x_\parallel} \eta (x_\parallel) \cdot \nabla_x V( s;t, \eta(x_\parallel), v)  & \nabla_v V( s;t, \eta(x_\parallel), v)
		\end{bmatrix}.
	\end{split}
	\Ee\unhide
	
	Since the characteristics is deterministic, for either $0< \e \ll 1$ or $0< - \e \ll 1$ with $X(t+ \e;t,\eta(x_\parallel), v) \in \bar{\O}$, 
	\Be \label{identi}
	\begin{split}
		\tb(  t+ \varepsilon, X(t+ \varepsilon;t, \eta(x_\parallel),v), V(t+ \varepsilon;t, \eta(x_\parallel),v))& =  \tb(t, \eta(x_\parallel), v) + \e.
	\end{split}
	\Ee\hide
	for both $\mathcal{C}_{+}$ case,
	\Be \label{g+}
	(s,t,(\eta(\bar{x}),v)) \in (t-\tb(t,\eta(\bar{x}),v), t) \times [0,T]\times \gamma_{+},\quad s-t \leq \varepsilon \leq 0,
	\Ee
	and $\mathcal{C}_{-}$ case,
	\Be \label{g-}
	(s,t,(\eta(\bar{x}),v)) \in (t, t + \tf(t,\eta(\bar{x}),v)) \times [0,T]\times \gamma_{-},\quad 0\leq \varepsilon \leq s-t.
	\Ee\unhide
	Differentiate (\ref{identi}) with respect to $\e$ and use (\ref{hamilton_ODE}), and then set $\e=0$ to have
	\Be\label{transport_tb} 
	[\p_t + v\cdot \nabla_x + E(t,\eta(x_\parallel)) \cdot \nabla_v]\tb(t,\eta(x_\parallel), v) 
	= 1.
	\Ee
	Also we have 
	\Be \label{identi}
	\begin{split}
		X(s; t+ \varepsilon, X(t+ \varepsilon;t, \eta(x_\parallel),v), V(t+ \varepsilon;t, \eta(x_\parallel),v))& = X(s;t, \eta(x_\parallel), v),  \\
		V(s; t+ \varepsilon, X(t+ \varepsilon;t, \eta(x_\parallel),v), V(t+ \varepsilon;t, \eta(x_\parallel),v))& = V(s;t, \eta(x_\parallel), v).
	\end{split}
	\Ee\hide
	for both $\mathcal{C}_{+}$ case,
	\Be \label{g+}
	(s,t,(\eta(\bar{x}),v)) \in (t-\tb(t,\eta(\bar{x}),v), t) \times [0,T]\times \gamma_{+},\quad s-t \leq \varepsilon \leq 0,
	\Ee
	and $\mathcal{C}_{-}$ case,
	\Be \label{g-}
	(s,t,(\eta(\bar{x}),v)) \in (t, t + \tf(t,\eta(\bar{x}),v)) \times [0,T]\times \gamma_{-},\quad 0\leq \varepsilon \leq s-t.
	\Ee\unhide
	We differentiate the above identity with respect to $\e$ and use (\ref{hamilton_ODE}), and then set $\e=0$ to have
	\Be \begin{split}\label{pt_XV}
		[\p_t + v\cdot \nabla_x + E(s,X(s;t,\eta(x_\parallel), v)) \cdot \nabla_v]X(s;t,\eta(x_\parallel), v) & = 0,\\
		[\p_t + v\cdot \nabla_x+ E(s,X(s;t,\eta(x_\parallel), v)) \cdot \nabla_v]V(s;t,\eta(x_\parallel), v) & = 0.
	\end{split}
	\Ee\hide
	or equivalently     
	\Be\label{pt_XV}
	\begin{split}
		\begin{bmatrix}
			\p_t X( s;t, \eta(x_\parallel), v)\\
			\p_t V( s;t, \eta(x_\parallel), v)
		\end{bmatrix}
		=&
		\begin{bmatrix}
			\nabla_x X( s;t, \eta(x_\parallel), v)  & \nabla_v X( s;t, \eta(x_\parallel ), v)\\
			\nabla_x V( s;t, \eta(x_\parallel), v)  & \nabla_v V(s;t, \eta(x_\parallel), v)
		\end{bmatrix}\\
		&\times 
		\begin{bmatrix}
			-v  \\
			-E (s,X(s;t,\eta(x_\parallel), v))
		\end{bmatrix}.
	\end{split}
	\Ee\unhide
	From (\ref{transport_tb}) and (\ref{pt_XV}) 
	\Be\begin{split}
		 (\ref{Jac_bdry_bdry_1})
		= &  		
		\left[\begin{matrix}
		v \cdot \nabla_x \tb + E \cdot \nabla_v \tb\\
		(- v\cdot \nabla_x X - E \cdot \nabla_v X) \cdot \nabla(\eta^{-1})_1\\
		(- v\cdot \nabla_x X - E \cdot \nabla_v X) \cdot \nabla(\eta^{-1})_2\\
- v\cdot \nabla_x V - E \cdot \nabla_v V
		\end{matrix}
		\right.
\\				&  \  \ 
		\left.\begin{matrix}
			- \nabla_{x_\parallel} \eta \cdot \nabla_x \tb & - \nabla_v \tb\\
			\left(\nabla_{x_\parallel} \eta  \cdot 
			\nabla_{x_\parallel}X  
			\right) \cdot \nabla (\eta^{-1})_1 &  \nabla_{v}X   \cdot \nabla (\eta^{-1})_1
			\\
						  \left(\nabla_{x_\parallel} \eta  \cdot 
			\nabla_{x_\parallel}X  
			\right) \cdot \nabla (\eta^{-1})_2 &  \nabla_{v}X   \cdot \nabla (\eta^{-1})_2
			\\
			    \nabla_{x_\parallel} \eta(x_\parallel) \cdot \nabla_x V & \nabla_v V
		\end{matrix}
		\right].\notag
	\end{split}\Ee
	Now we multiply $v \cdot \p_{x_{\parallel,i}} \eta$ to the $(i+1)$-column for $i=1,2$ and add these to the first column. Then we multiply $E_i$ to the $(i+3)$-column for $i=1,2,3$ and add these to the first column. And then we use (\ref{nabla_tb}) to get
	\Be
	\left[\begin{matrix}
		(n(x_\parallel) \cdot v)  \frac{n(x_\parallel) \cdot \nabla X(t-\tb) \cdot n(\xb)}{n(\xb) \cdot \vb}
		& - \nabla_{x_\parallel} \eta \cdot \frac{\nabla_x X\cdot n(\xb)}{n(\xb) \cdot\vb} 
		& - \frac{ \nabla_v X \cdot n(\xb)}{n(\xb) \cdot \vb}\\
		- v\cdot n  \p_n X    \cdot \nabla(\eta^{-1})_1&
		\left(\nabla_{x_\parallel} \eta  \cdot 
		\nabla_{x_\parallel}X  
		\right) \cdot \nabla (\eta^{-1})_1 &  \nabla_{v}X   \cdot \nabla (\eta^{-1})_1
		\\
		- v\cdot n \p_n X  \cdot \nabla(\eta^{-1})_2
		&  \left(\nabla_{x_\parallel} \eta  \cdot 
		\nabla_{x_\parallel}X  
		\right) \cdot \nabla (\eta^{-1})_2 &  \nabla_{v}X   \cdot \nabla (\eta^{-1})_2
		\\
		- v\cdot n \p_n V    & \nabla_{x_\parallel} \eta(x_\parallel) \cdot \nabla_x V & \nabla_v V
	\end{matrix}
	\right].\notag
	\Ee
	Therefore the determinant of (\ref{Jac_bdry_bdry_1}) equals
	\Be
	\begin{split}
		\frac{n(x_\parallel) \cdot v}{n(\xb) \cdot \vb}
		&\times
		\det
		\left[
		\begin{matrix}
			n(x_\parallel) & \nabla_{x_\parallel} \eta  & 0 \\
			0&0 & \text{Id}_{3 \times 3}
		\end{matrix}\right]
		\det\left[
		\begin{matrix}
			\nabla_x X & \nabla_x V\\
			\nabla_v X & \nabla_v V
		\end{matrix}\right]\\
		&\times 
		\det \left[
		\begin{matrix}
			n(\xb) & \nabla(\eta^{-1}) & 0\\
			0 & 0 & \text{Id}_{3 \times 3}
		\end{matrix}\right],\notag
	\end{split}	\Ee
	which implies (\ref{Jac_boundary_boundary}).\hide
	From (\ref{XV_txv}) and (\ref{pt_XV}),
	\Be
	\begin{split}\label{XV_txv_final}
		(\ref{XV_txv})
		= &\begin{bmatrix}
			\nabla_x X( s;t, \eta(x_\parallel), v)  & \nabla_v X( s;t, \eta(x_\parallel), v)\\
			\nabla_x V( s;t, \eta(x_\parallel), v)  & \nabla_v V(s;t, \eta(x_\parallel), v)
		\end{bmatrix}  \\
		&\times 
		\begin{bmatrix}
			-v & 
			\p_{x_\parallel} \eta (x_\parallel)  & 0_{3 \times 3} \\
			-E(s, X(s;t,\eta(x_\parallel),v))
			& 
			0_{3 \times 2} & \mathrm{Id}_{3\times 3}
		\end{bmatrix}.  
	\end{split}
	\Ee
	Since (from Liouville theorem)
	$$\det \begin{bmatrix}
	\nabla_x X( s;t, \eta(x_\parallel), v)  & \nabla_v X( s;t, \eta(x_\parallel), v)\\
	\nabla_x V( s;t, \eta(x_\parallel), v)  & \nabla_v V(s;t, \eta(x_\parallel), v)
	\end{bmatrix} =1,
	$$
	we conclude that
	\Be \label{cov2}
	\begin{split} 
		(\ref{XV_txv})	 &=\det \begin{bmatrix}
			-v & 
			\p_{x_\parallel} \eta (x_\parallel)  & 0_{3 \times 3} \\
			-E(s, X(s;t,\eta(x_\parallel),v))
			& 
			0_{3 \times 2} & \mathrm{Id}_{3\times 3}
		\end{bmatrix}\\
		&=   - v \cdot (\p_1 \eta(x_\parallel) \times \p_2 \eta(x_\parallel) ).  
	\end{split}\Ee
	Since the surface measure equals $\dd S_x = |\p_1 \eta (x_\parallel) \times \p_2 \eta (x_\parallel)| \dd x_\parallel$ we conclude (\ref{Jacobian_C+}).
	
	For (\ref{Jac+}) we compute 
	\Be\begin{split}\notag
		&\det\left(\frac{\p (t+s, X(t+s;t,\eta(x_\parallel),v), V(t+s;t,\eta(x_\parallel),v ))}{\p (t,s,x_\parallel,v)}\right)\\
		=&\det  \left[\begin{matrix}
			1 &1 & 0 & 0 \\
			\p_s X(t+s)  + \p_t X(t+s) & \p_s X(t+s)  & \nabla_{x_\parallel} X(t+s) &  \nabla_v X(t+s)\\
			\p_s V(t+s)  + \p_t V(t+s) & \p_s V(t+s)  & \nabla_{x_\parallel}  V(t+s) &  \nabla_v V(t+s)
		\end{matrix}
		\right]\\
		=& \det  \left[\begin{matrix}
			0&1  & 0 & 0 \\
			\p_t X(t+s) &  \p_s X(t+s)  & \nabla_{x_\parallel}  X(t+s) &  \nabla_v X(t+s)\\
			\p_t V(t+s) & \p_s V(t+s)  & \nabla_{x_\parallel}  V(t+s) &  \nabla_v V(t+s)
		\end{matrix}
		\right]\\
		=& \det  \left[\begin{matrix}
			\p_t X(t+s)    & \nabla_{x_\parallel}  X(t+s) &  \nabla_v X(t+s)\\
			\p_t V(t+s)    & \nabla_{x_\parallel}  V(t+s) &  \nabla_v V(t+s)
		\end{matrix}
		\right],
	\end{split}	\Ee
	which equals (\ref{XV_txv}).\hide
	
	Now we prove (\ref{Jac+}) and (\ref{Jac-}). In this proof, we use the notation $ \partial_t X(t-s;t,x,v ) :=  \frac{ \partial X(s;t,x,v ) }{\partial t } |_{(t-s;t,x,v )} $, and $\partial_s X(t-s;t,x,v ) :=  \frac{ \partial X(s;t,x,v ) }{\partial s } |_{(t-s;t,x,v )} $ to denote the partial derivatives of $X$ for the first and second coordinates, and similar for $V$. We note that these are different to the derivative with respect to $s$ or $t$ which we denote $\frac{d X(t-s; t,x,v ) }{ ds }$ and $\frac{d X(t-s; t,x,v ) }{ dt }$.  
	
	Let us compute $\mathcal{A}_{+}$ case first. We have 
	\Be
	\begin{split}
		\frac{d X(t-s; t,x,v ) }{ ds } &= - \partial_s X(t-s; t,x,v ),  \\
		\frac{d X(t-s; t,x,v ) }{ dt } &= \partial_s X(t-s; t,x,v ) +  \partial_t X(t-s; t,x,v ),  \\
	\end{split}
	\Ee
	and from direct computation and elementary column operation,
	\tiny
	\Be \label{77} 
	\begin{split}
		&\det J_{\mathcal{A}_{+} }\\
		&= 
		\det
		\begin{bmatrix}
			1 & -1 & 0_{1\times 2 } &  0_{1 \times 3 } \\
			\partial_s X(t-s; t,x,v ) +  \partial_t X(t-s; t,x,v ) & -\partial_s X(t-s; t,x,v ) & \partial_{\bar x } X(t-s; t,x,v ) & \partial_v X(t-s; t,x,v ) \\
			\partial_s V(t-s; t,x,v ) +  \partial_t V(t-s; t,x,v ) & -\partial_s V(t-s; t,x,v ) & \partial_{\bar x } V(t-s; t,x,v ) & \partial_v V(t-s; t,x,v ) \\
		\end{bmatrix}  \\
		&= \det \begin{bmatrix}
			1 & 0 & 0_{1\times 2 } &  0_{1 \times 3 } \\
			\partial_s X(t-s; t,x,v ) +  \partial_t X(t-s; t,x,v ) & \partial_t X(t-s; t,x,v ) & \partial_{\bar x } X(t-s; t,x,v ) & \partial_v X(t-s; t,x,v ) \\
			\partial_s V(t-s; t,x,v ) +  \partial_t V(t-s; t,x,v ) & \partial_t V(t-s; t,x,v ) & \partial_{\bar x } V(t-s; t,x,v ) & \partial_v V(t-s; t,x,v ) \\
		\end{bmatrix}  \\
		&= \det \begin{bmatrix}
			\partial_t X(t-s; t,x,v ) & \partial_{\bar x } X(t-s; t,x,v ) & \partial_v X(t-s; t,x,v ) \\
			\partial_t V(t-s; t,x,v ) & \partial_{\bar x } V(t-s; t,x,v ) & \partial_v V(t-s; t,x,v ) \\
		\end{bmatrix} .
	\end{split} 
	\Ee
	\normalsize
	
	Note that $0\leq s \leq \tb(t,x,v)$ implies $t-\tb(t,x,v) \leq t-s \leq t$, so condition (\ref{g+}) is satisfied. From (\ref{cov2}), we obtain (\ref{Jac+}).  \\
	
	For $\mathcal{A}_{-}$ case, we apply similar computation as (\ref{77}):
	\tiny
	\Be \label{77-} 
	\begin{split}
		&\det J_{\mathcal{A}_{-} } \\
		&= 
		\det
		\begin{bmatrix}
			1 & 1 & 0_{1\times 2 } &  0_{1 \times 3 } \\
			\partial_s X(t+s; t,x,v ) +  \partial_t X(t+s; t,x,v ) & \partial_s X(t+s; t,x,v ) & \partial_{\bar x } X(t+s; t,x,v ) & \partial_v X(t+s; t,x,v ) \\
			\partial_s V(t+s; t,x,v ) +  \partial_t V(t+s; t,x,v ) & \partial_s V(t+s; t,x,v ) & \partial_{\bar x } V(t+s; t,x,v ) & \partial_v V(t+s; t,x,v ) \\
		\end{bmatrix}  \\
		&= \det \begin{bmatrix}
			1 & 0 & 0_{1\times 2 } &  0_{1 \times 3 } \\
			\partial_s X(t+s; t,x,v ) +  \partial_t X(t+s; t,x,v ) & -\partial_t X(t+s; t,x,v ) & \partial_{\bar x } X(t+s; t,x,v ) & \partial_v X(t+s; t,x,v ) \\
			\partial_s V(t+s; t,x,v ) +  \partial_t V(t+s; t,x,v ) & -\partial_t V(t+s; t,x,v ) & \partial_{\bar x } V(t+s; t,x,v ) & \partial_v V(t+s; t,x,v ) \\
		\end{bmatrix}  \\
		&= - \det \begin{bmatrix}
			\partial_t X(t+s; t,x,v ) & \partial_{\bar x } X(t+s; t,x,v ) & \partial_v X(t+s; t,x,v ) \\
			\partial_t V(t+s; t,x,v ) & \partial_{\bar x } V(t+s; t,x,v ) & \partial_v V(t+s; t,x,v ) \\
		\end{bmatrix} .
	\end{split} 
	\Ee
	\normalsize
	Note that $0\leq s \leq \tf(t,x,v)$ implies $t \leq t+s \leq t+\tf(t,x,v)$, so condition (\ref{g-}) holds. From (\ref{cov2}), we obtain (\ref{Jac-}).  \\
	\unhide\unhide
\end{proof}\hide

\begin{lemma}
	\Be
	\frac{\p (t-\tb(t,x,v), \xb(t,x,v), \vb(t,x,v))}{\p (t,x,v)}
	\Ee
\end{lemma}
\begin{proof}
	\Be\begin{split}
		&\frac{\p (t-\tb(t,\eta(x_\parallel),v), \xb(t,\eta(x_\parallel),v), \vb(t,\eta(x_\parallel),v))}{\p (t,x_\parallel,v)}\\
		=&
		\left[\begin{matrix} 
			1- \frac{\p \tb}{\p t} & - \nabla_{x_\parallel} \eta (x_\parallel) \cdot \nabla_x \tb& - \nabla_{v}   \tb\\
			(1- \frac{\p \tb}{\p t}) \vb + \p_t X(t-\tb;t,x,v) &  \nabla_{x_\parallel} \eta (x_\parallel ) \cdot \nabla_x \xb & \nabla_v \xb\\
			(1- \frac{\p \tb}{\p t}) E(t-\tb, \xb)  & \nabla_{x_\parallel} \eta (x_\parallel) \cdot \nabla_x \vb & \nabla_v \vb
		\end{matrix}\right]
	\end{split}\Ee
\end{proof}\unhide

\begin{lemma}\label{lem COV2}
	For fixed $s>0$ so that $t-\tb (t,x, v ) < s < t$,
	\Be \label{mapC+}
	\begin{split}
		( t, x , v )  \in   [0,T]\times\gamma_{+} &\mapsto (X(s;t, x, v ) , V(s; t, x, v ) ) \in  \Omega\times\R^{3}  
	\end{split}
	\Ee 
	is injective. 
	For fixed $s>0$ so that $t < s < t + \tf (t,x, v ) $,  
	\Be \label{mapC-}
	\begin{split} 
		( t, {x} , v )  \in [0,T]\times\gamma_{-} \mapsto  (X(s;t,x , v ) , V(s; t, x, v ) )\in\Omega\times\R^{3} 
	\end{split}
	\Ee
	is also injective. For both maps (\ref{mapC+}) and (\ref{mapC-}),
	\Be\label{Jacobian_C+}
	\left|\det\left(\frac{\p (X(s;t, x, v ) , V(s; t, x, v ) ) }{\p  ( t, x , v ) }\right)\right| = |
	n(x)	 \cdot v
	|.
	\Ee
	
Moreover 
	\Be \label{map+}
	\begin{split}
		& ( t, s, x , v )
		\in [0,T] \times \{ - \min \{t, \tb(t,x,v)\}< s<0 \} \times \gamma_+\\
		\mapsto 	&	
		\big( t+s, X(t+s;t, x, v ) , V(t+s; t, x, v ) \big)  
		\in  [0,T]\times\Omega\times\R^{3},
	\end{split}
	\Ee 
	and 
	\Be \label{map-} 
	\begin{split}
		& ( t, s,x , v )  
		\in 	    [0,T] \times\{ 0<s < \min \{\tf(t,x,v),T-t\} \} \times\gamma_{-}\\
		\mapsto  & \big(t+s, X(t+s;t,x, v ), V(t+s; t,x, v ) \big)
		\in [0,T] \times \O \times \R^3
	\end{split}\Ee
	are both injective, and for both maps (\ref{map+}) and (\ref{map-})
	\Be\label{Jac+}
	\left|\det \left(\frac{\p (t+s, X(t+s;t,x,v), V(t+s;t,x,v ))}{\p (t,s,x,v)} \right)\right|=   |n(x) \cdot 
	v
	|.
	\Ee
\end{lemma}
\begin{proof}
	
	The maps are injective since the characteristics are a solution of (\ref{hamilton_ODE}).

	For $x \in \p\O$ we have $x= \eta(x_\parallel)$ locally from (\ref{eta}). Then we compute (\ref{Jac+}),
	\Be
	\begin{split}\label{XV_txv}
		& \frac{\p ( X( s;t, \eta(x_\parallel), v),  V( s;t, \eta(x_\parallel), v) )}{\p (t, x_\parallel, v)}\\
		= &  \begin{bmatrix}
			\p_t X( s;t, \eta(x_\parallel), v) & \nabla_{x_\parallel} X( s;t, \eta(x_\parallel), v) & \nabla_v X( s;t, \eta(x_\parallel), v)\\
			\p_t V( s;t, \eta(x_\parallel), v)& \nabla_{x_\parallel} V( s;t, \eta(x_\parallel), v) & \nabla_v V( s;t, \eta(x_\parallel), v)
		\end{bmatrix}\\
		= &  
		\left[\begin{matrix}
			\p_t X( s;t, \eta(x_\parallel), v)& \nabla_{x_\parallel} \eta (x_\parallel) \cdot \nabla_x X( s;t, \eta(x_\parallel), v)  \\
			\p_t V( s;t, \eta(x_\parallel), v) &  \nabla_{x_\parallel} \eta (x_\parallel) \cdot \nabla_x V( s;t, \eta(x_\parallel), v)    \end{matrix}\right.\\
			& \ \    \left. \begin{matrix}
			\nabla_v X( s;t, \eta(x_\parallel), v)
			\\
			\nabla_v V( s;t, \eta(x_\parallel), v)\
		\end{matrix}\right].
	\end{split}
	\Ee
	From (\ref{XV_txv}) and (\ref{pt_XV}),
	\Be
	\begin{split}\label{XV_txv_final}
		(\ref{XV_txv})
		= &\begin{bmatrix}
			\nabla_x X( s;t, \eta(x_\parallel), v)  & \nabla_v X( s;t, \eta(x_\parallel), v)\\
			\nabla_x V( s;t, \eta(x_\parallel), v)  & \nabla_v V(s;t, \eta(x_\parallel), v)
		\end{bmatrix}  \\
		&\times 
		\begin{bmatrix}
			-v & 
			\p_{x_\parallel} \eta (x_\parallel)  & 0_{3 \times 3} \\
			-E(s, X(s;t,\eta(x_\parallel),v))
			& 
			0_{3 \times 2} & \mathrm{Id}_{3\times 3}
		\end{bmatrix}.  
	\end{split}
	\Ee
	From Liouville theorem, 
	we conclude that
	\Be \label{cov2}
	\begin{split} 
		(\ref{XV_txv})	 &=\det \begin{bmatrix}
			-v & 
			\p_{x_\parallel} \eta (x_\parallel)  & 0_{3 \times 3} \\
			-E(s, X(s;t,\eta(x_\parallel),v))
			& 
			0_{3 \times 2} & \mathrm{Id}_{3\times 3}
		\end{bmatrix}\\
		&=   - v \cdot (\p_1 \eta(x_\parallel) \times \p_2 \eta(x_\parallel) ).  
	\end{split}\Ee
	Since the surface measure equals $\dd S_x = |\p_1 \eta (x_\parallel) \times \p_2 \eta (x_\parallel)| \dd x_\parallel$ we conclude (\ref{Jacobian_C+}).
	
	For (\ref{Jac+}) we compute 
	\Be\begin{split}\notag
		&\det\left(\frac{\p (t+s, X(t+s;t,\eta(x_\parallel),v), V(t+s;t,\eta(x_\parallel),v ))}{\p (t,s,x_\parallel,v)}\right)\\
		=&\det  \left[\begin{matrix}
			1 &1 & 0 & 0 \\
			\p_s X(t+s)  + \p_t X(t+s) & \p_s X(t+s)  & \nabla_{x_\parallel} X(t+s) &  \nabla_v X(t+s)\\
			\p_s V(t+s)  + \p_t V(t+s) & \p_s V(t+s)  & \nabla_{x_\parallel}  V(t+s) &  \nabla_v V(t+s)
		\end{matrix}
		\right]\\
		=& \det  \left[\begin{matrix}
			0&1  & 0 & 0 \\
			\p_t X(t+s) &  \p_s X(t+s)  & \nabla_{x_\parallel}  X(t+s) &  \nabla_v X(t+s)\\
			\p_t V(t+s) & \p_s V(t+s)  & \nabla_{x_\parallel}  V(t+s) &  \nabla_v V(t+s)
		\end{matrix}
		\right]\\
		=& \det  \left[\begin{matrix}
			\p_t X(t+s)    & \nabla_{x_\parallel}  X(t+s) &  \nabla_v X(t+s)\\
			\p_t V(t+s)    & \nabla_{x_\parallel}  V(t+s) &  \nabla_v V(t+s)
		\end{matrix}
		\right],
	\end{split}	\Ee
	which equals (\ref{XV_txv}).
\end{proof}

We define \textit{the forward exit time}
\Be
\tf (t,x,v) : = \sup \{s\geq 0: X(\tau;t,x,v) \in \O \ \ \text{for all }   \tau \in (t,t+s)   \}.\label{forward_exit}
\Ee 
\begin{lemma} \label{int_id}
	Suppose $h(t,x,v) \in L^1 ( [0,T] \times \Omega \times \mathbb R^3) $. Then 
	\[
	\begin{split}
	&	\int_0 ^T  \iint _{\Omega \times \mathbb R^3} h (t,x,v) \dd v \dd x \dd t \\ &=  \iint_{\Omega \times \mathbb R^3 } \int_{- \min \{ T,\tb(T,x,v) \}} ^0\\
	& \ \ \ \ \ \ \ \   \ \ \ \ \ \ \ \   \ \ \ \ \ \ \ \ \times  h (T +s, X(T+s; T,x,v), V(T+s;T,x,v) ) \dd s \dd v \dd x \\ & + \int_0 ^T \int_{\gamma_+ } \int_{-\min \{ t, \tb(t,x,v)\} } ^0 h(t +s , X(t+s;t,x,v), V(t +s;t,x,v) )\dd s \dd\gamma \dd t.
	\end{split}
	\]
\end{lemma}
\begin{proof}
	The region $\{ (t,x,v) \in [0,T] \times \Omega \times \mathbb R^3 \}$ is a disjoint union of 
	\Be
	\begin{split}\notag
		A:= \{ (t,x,v) \in [0,T] \times \Omega \times \mathbb R^3: \tf (t,x,v) + t\le T \} , \\
		B:= \{ (t,x,v) \in [0,T] \times \Omega \times \mathbb R^3: \tf (t,x,v)+ t > T \} .\end{split}
	\Ee 
	We also define
	\Be
	\begin{split}\notag
		A^\prime &:=  \{ (t,s,x,v) \in [0,T]\times \R \times \gamma_+ :  -\min\{ \tb(t,x,v), t\} \leq s\leq 0 \},  \\
		B^\prime &:= \{ (s,x,v ) \in [0,T] \times \Omega \times \mathbb R^3 : s \leq \tb(T,x,v) \}.
	\end{split}
	\Ee 
	
	Let us denote (\ref{map+}) by $\mathcal{A}_+: A^\prime \to A$. Since $\tf(t+s, X(t+s ; t,x,v) , V(t+s;t,x,v) ) + (t+s) = -s + (t+s) = t \le T $ if $(x,v) \in \gamma_+$ and $- \min \{\tb(t,x,v),t\}<s<0$, this map is well-defined. For any $(t,x,v) \in A$, we have
	\[
	(t + \tf(t,x,v), -\tf(t,x,v), X(t + \tf(t,x,v);t,x,v), V(t + \tf(t,x,v);t,x,v) ) \in A',
	\] 
	since $t + \tf \leq T$ and $\tb ( t + \tf, X(t + \tf;t,x,v ) , V(t + \tf;t,x,v ) ) > \tf $. Moreover 
	\[
	\mathcal A_+ (t + \tf, -\tf, X(t + \tf;t,x,v), V(t + \tf;t,x,v) ) = (t,x,v) 
	\] 
	implies that $\mathcal A_+$ is surjective. From Lemma \ref{lem COV}, $\mathcal A_+$ is bijective. 
	Applying the change of variable $\mathcal A_+$ with the Jacobian (\ref{Jac+}),  
	\Be \notag\begin{split}
		& \iiint_A  h(t,x,v ) \dd t \dd x \dd v 
		\\ 
		\hide= & \int_0^T \int_{\gamma_+ } \int_0 ^{\min\{ \tb(t, x,v) , t \} } h(t-s, X(t-s;t,x,v ) ,V(t-s; t,x,v )  \\
		& \ \ \ \ \quad\quad \times \big| v \cdot  \frac{ \p_1 \eta({x}_\parallel) \times \p_2 \eta({x}_\parallel) }{   | \p_1 \eta({x}_\parallel) \times \p_2 \eta({x}_\parallel) |} \big| | (\p_1 \eta({x}_\parallel) \times \p_2 \eta({x}_\parallel) )| \dd s \dd {x}_\parallel \dd t
		\\\unhide = & \int_0^T \int_{\gamma_+ } \int^0 _{-\min\{ \tb(t,x,v) , t \} } h(t+s, X(t+s;t,x,v ) ,V(t+s; t,x,v )    \dd s \dd \gamma \dd t.
	\end{split} \Ee   
	
	Next, we consider a map 
	\Be\label{Map_B}
	(s,x,v) \in B^\prime\mapsto (T -s , X(T-s; T,x,v ), V(T-s; T,x,v ) ) \in B.
	\Ee
	%
	From $ \tf(T-s, X(T-s; T,x,v ),V(T-s; T,x,v ) ) + (T-s) > s + (T-s) = T$, the map is well-defined. Since the characteristic is deterministic, the map is injective. 
	Moreover it is also surjective since for $(s,x,v) \in B$ we have $(T-s, X(T; s, x,v ) , V(T;s,x,v) ) \in B^\prime$ and $(T-s, X(T; s, x,v ) , V(T;s,x,v) )\mapsto (s,x,v)$ by this map (\ref{Map_B}). It is well-known that this map has a unit Jacobian (Liouville theorem).
	By the change of variable of (\ref{Map_B}) and $s \mapsto -s$
	\Be\notag
	\begin{split}
		&\iiint_B h(t,x,v ) \dd t \dd x\dd v   \\
		= &   \iint_{\Omega \times \mathbb R^3 } \int^0_{- \min (T,\tb(T,x,v) ) }  h(T+s, X(T+s;T,x,v ) ,V(T+s; T,x,v ) \dd s \dd x \dd v.
	\end{split}
	\Ee\end{proof}

The first result is an ``energy estimate" to the transport operator. 

\begin{lemma} [Green's identity]  \label{lem_Green}
	For $p \in [1, \infty)$, we assume $f \in L^{p}_{loc} (\R_+ \times \Omega \times \mathbb R^3 )$ satisfies
	\Be\notag
	\partial_t f + v \cdot \nabla_x f + E \cdot \nabla_v f  \in L^p_{loc} (\R_+; L^p (\Omega \times \mathbb R^3 ) ),  \ \
	f   \in L^p_{loc} (\R_+; L^p (\gamma_+ ) ).	
	\Ee
	Then $f \in C^0_{loc}( \R_+ ; L^p (\Omega \times \mathbb R^3 ) )$ and $f  \in L^p_{loc} (\R_+; L^p (\gamma_-) )$. 
	
	Moreover
	\Be \label{Greedid}
	\begin{split}
		\| f(T)\|_p ^p + \int_0^{T} |f|_{ p,+ } ^p &= \| f(0) \|_p^p + \int_0^{T} |f|_{ p,- }^p  \\
		&  +p \int_0^{T} \iint _{\Omega \times \mathbb R^3 }   \{ \partial_t + v \cdot \nabla_x f + E \cdot \nabla_v f \} |f|^{p-2} f.
	\end{split} 
	\Ee
\end{lemma}
\begin{proof}\hide
	Using H\"older's inequality we have 
	\[
	\begin{split}
	&\| (\partial_t f + v \cdot \nabla_x f + E \cdot \nabla_v f) |f|^{p-2} f  \|_{L^1 ( [0,T] \times \Omega \times \mathbb R^3 )}  \\
	&\quad \leq \|(\partial_t f + v \cdot \nabla_x f + E \cdot \nabla_v f) \|_{L^p ([0,T] \times \Omega \times \mathbb R^3 ) } \| f \|^{p-1}_{L^{p} ([0,T] \times \Omega \times \mathbb R^3 )} < \infty.
	\end{split}
	\]\unhide
	By Lemma \ref{int_id}, 
	\Be
	\begin{split}\label{Green1}
		& p \int_0^{T}  \iint_{\Omega \times \mathbb R^3 } \{\partial_t f + v \cdot \nabla_x f + E \cdot \nabla_v f\} |f|^{p-2} f  \dd x\dd v\dd t  
		\\ 
		= & \ p \iint_{\Omega \times \mathbb R^3 } \int_{- \min \{ T,\tb(T,x,v) \}} ^0\{\partial_t f + v \cdot \nabla_x f + E \cdot \nabla_v f\} \\
		&  \ \ \ \ \ \ \ \ \ \ \ \     \times |f|^{p-2} f(T  +s, X(T +s; T ,x,v), V(T +s;T ,x,v) ) \dd s \dd v \dd x  \\ 
		+	&  p  \int_0 ^{T} \int_{\gamma_+ } \int_{-\min \{ t, \tb(t,x,v)\} } ^0 \{\partial_t f + v \cdot \nabla_x f + E \cdot \nabla_v f\} \\
		&  \ \ \ \ \ \ \ \ \ \ \ \  \times  |f|^{p-2} f (t +s , X(t+s;t,x,v), V(t +s;t,x,v) ) \dd s \dd\gamma \dd t  .
	\end{split}\Ee
	Note that 
	\Be
	\begin{split}\notag
		&\frac{d}{ds}  |f( t+s, X(t+s; t,x,v ), V(t+s; t,x,v ))|^p \\
		= & \ \{\partial_t f + v \cdot \nabla_x f + E \cdot \nabla_v f\} \\
		& \times  |f|^{p-2} f (t +s , X(t+s;t,x,v), V(t +s;t,x,v) ) .
	\end{split}
	\Ee	
	We plug this into (\ref{Green1}) and apply the integration by parts in $s$ to obtain 
	\Be\begin{split} \label{green}
		\frac{1}{p} \times (\ref{Green1}) 
		&= \iint_{\Omega \times \mathbb R^3 } |f(T ,x,v )|^p  \dd x \dd v  + \int_0^{T} \int_{\gamma ^+ } |f (t,x,v )|^p \dd\gamma \dd t  \\
		&  - \underbrace{ \iint_{\Omega \times \mathbb R^3 } \textbf{1}_{  T  \ge \tb(T ,x,v)  } |f  (T -\tb,\xb,\vb ) |^p\dd x \dd v }_{(\ref{green})_1 } \\
		&  - \underbrace{ \iint_{\Omega \times \mathbb R^3 } \textbf{1}_{  T < \tb(T ,x,v) } |f (0,X(0;T ,x,v) ,V(0;T ,x,v ) )|^p \dd x \dd v }_{(\ref{green})_2}  \\
		&  - \underbrace{ \int_0^{T } \int_{\gamma_+ } \textbf{1}_{ t \ge \tb (t,x,v )  }|f (t -\tb, \xb,\vb )|^p \dd\gamma \dd t }_{(\ref{green})_3}  \\
		&   - \underbrace{ \int_0^{T} \int_{\gamma_+ }  \textbf{1}_{  t < \tb (t,x,v )  }|f (0, X(0;t,x,v ),V(0;t,x,v) ) |^p \dd\gamma \dd t }_{(\ref{green})_4} .
	\end{split} 
	\Ee

	We claim that \Be \label{O12}
	(\ref{green})_2 + (\ref{green})_4= \iint_{\O\times\R^{3}} |f_{0}|^{p} \dd x \dd v.
	\Ee 
	If $T< \tb(T,x,v)$ then $\tf (0,X(0;T,x,v), V(0;T,x,v)) >T$. On the other hand, if $t< \tb(t,x,v)$ then $\tf(0,X(0;t,x,v),V(0;t,x,v))<T$. Now we apply the change of variables $(x,v) \mapsto (X(0;t,x,v),V(0;t,x,v))$ to $(\ref{green})_2$ and to $(\ref{green})_4$ to conclude (\ref{O12}).
	
	Next we claim that 
	\Be\label{O34}
	(\ref{green})_1+ (\ref{green})_3 = \int^T_0 |f|^p_{p,-}.
	\Ee
	We split 
	\Be\label{split_tf_T}
	[0,T] \times \gamma_- = \{ \tf (s,x,v) + s>T \} \cup \{ \tf(s,x,v)+s \leq T\}.
	\Ee
	We consider the map (\ref{map_xv_boundary}) for fixed $T>0$ from $\O \times \R^3$ to $\{ (s,x,v) \in [0,T] \times \gamma_-: \tf (s,x,v) + s>T \}$. This map is onto since $\tb(T,X(T;s,x,v),V(T;s,x,v)) = T-s$ and therefore $(X(T;s,x,v),V(T;s,x,v)) \mapsto (s,x,v )$. By the change of variables with (\ref{Jac_xv_boundary})
	\[
	(\ref{green})_1 = \int^T_0  \iint_{\gamma_-} \mathbf{1}_{\tf (s,x,v) + s>T} |f(s,x,v)|^p \dd \gamma \dd s 	.\]

	For $(\ref{green})_3$ we consider (\ref{map_boundary_boundary}) from $[0,T]\times \gamma_+$ to $\{ (s,x,v) \in [0,T] \times \gamma_-: \tf (s,x,v) + s\leq T \}$. This map is onto since $\tb(s+ \tf(s,x,v), X(s+ \tf;s,x,v), V(s+ \tf;s,x,v))= \tf(s,x,v)$ and therefore $(s+ \tf(s,x,v), X(s+ \tf;s,x,v), V(s+ \tf;s,x,v)) \mapsto (s,x,v)$. Applying the change of variables of (\ref{Jac_boundary_boundary}), we obtain that 
	\[
	(\ref{green})_2 = \int^T_0  \iint_{\gamma_-} \mathbf{1}_{\tf (s,x,v) + s\leq T} |f(s,x,v)|^p \dd \gamma \dd s 	.
	\]	\end{proof}

	\hide We consider the map $(x,v) \mapsto  (T-\tb(T,x,v), \xb(T,x,v),\vb(T,x,v))$ for $T\geq \tb(T,x,v)$. Clearly this the inverse map of (\ref{mapC-}). Also we have $(T-\tb(T,x,v))+ \tf (T-\tb(T,x,v), \xb(T,x,v),\vb(T,x,v))
	> T$. Therefore
	\Be\notag
	\begin{split}
		(\ref{green})_1 = \int_0^T \int_{\gamma_-} 
		\mathbf{1}_{\tf(t,x,v)+ t>T}
		|f(t,x,v)|^p
		\dd \gamma \dd t.
	\end{split}
	\Ee

	Note that if $T \geq \tb(T,x,v)$ then $\tf(T-\tb(T,x,v), \xb (T,x,v), \vb(T,x,v))>\tb(T,x,v)$. On the other hand if $t>\tb(t,x,v)$ and $(x,v) \in \gamma_+$ then $\tf (t-\tb(t,x,v),\xb(t,x,v), \vb(t,x,v))< \tb(t,x,v)$.

	through change of variables. For $O_{1}$, we consider the map 
	\Be \label{A1}
	\begin{split}
		&\mathcal{O}_{1} : \{ (x,v) \in \Omega \times \mathbb R^3 : T' < \tb(T',x,v) \} \mapsto \{ (x,v) \in \Omega \times \mathbb R^3 : \tf(0,x,v) >  T' \},  \\
		&\mathcal{O}_{1}(x,v) = (X(0; T',x,v), V(0; T',x,v) ) .  
	\end{split}
	\Ee	
	$\mathcal{O}_{1}$ is well-defined as $\tf(0, X(0; T',x,v), V(0; T',x,v)) > T'$ since $x \in \Omega$.  $\mathcal{O}_1$ is injective as the characteristic flow is unique. And for any $ (x,v) \in \Omega \times \mathbb R^3$ such that $\tf(0,x,v) > T'$, $\mathcal{O}_1 ( X(T';0,x,v ), V(T';0,x,v ) ) = (x,v )$ implies that $\mathcal{O}_1$ is surjective. Therefore $\mathcal{O}_1$ is a bijection. Since $\mathcal{O}_{1}$ is measure preserving, we have
	\begin{equation} \label{green1}
	\begin{split}
	O_{1} &:= \iint_{\Omega \times \mathbb R^3 } |f | ^p (X(0; T',x,v), V(0;T',x,v)) \textbf{1}_{\{ T' < \tb(T',x,v) \}}  dx dv  \\
	&= \iint_{\Omega \times \mathbb R^3 } |f_0|^p \textbf{1}_{\{ \tf(0,x,v) > T' \} } dxdv .  \\
	\end{split}
	\end{equation}
	For $O_{2}$,  
	we recall $\mathcal{C}_{+}$ in (\ref{map+}) and pick $s=0$, 
	\Be \label{A2}
	\begin{split}
		&\mathcal{C}_{+} : \{ (t,x,v) \in  (0, T'] \times \gamma_+ : t < \tb(t,x,v) \} \mapsto \{ (x,v) \in \Omega \times \mathbb R^3 : \tf(0,x,v) \le T' \},  \\
		&\mathcal{C}_{+}(t,x,v) = (X(0;t,x,v), V(0;t,x,v)) .  
	\end{split}
	\Ee	
	$\mathcal{C}_{+}$ is well-defined as $\tf( 0, X(0;t,x,v), V(0;t,x,v) ) = t \le T'$. $\mathcal{C}_{+}$ is injective as the characteristic flow is unique and for any $(x,v) \in \Omega \times \mathbb R^3$ such that $\tf(0,x,v) \le T'$, 
	\[
	\mathcal{C}_{+} ( \tf, X(\tf; 0,x,v) , V(\tf; 0,x,v) ) = (x,v),
	\]
	implies that $\mathcal{C}_{+}$ in (\ref{A2}) is surjective. Therefore, $\mathcal{C}_{+}$ is a bijection and we apply (\ref{cov2}) to peform change of variable :
	\begin{equation} \label{green2}
	\begin{split}
	O_{2} &:= \int_0^{T' } \int_{ \gamma_+ } |f_0|^p (X(0;t,x,v), V(0;t,x,v) ) \textbf{1}_{ \{ t < \tb(t,x,v) \} } d\gamma d t  \\
	&= \iint_{\Omega \times \mathbb R^3 } |f_0 |^p \textbf{1}_{ \{ \tf(0,x,v) \le T' \} }   dx dv. 
	\end{split}
	\end{equation}
	From (\ref{green1}) and (\ref{green2}), we get (\ref{O12}). \\

	\textit{Step 2.} For underbraced $M_{1}$ and $M_{2}$ in (\ref{green}), we show that 
	\Be \label{M12}
	M_{1} + M_{2} = \int_{0}^{T^{\prime}}|f|^{p}_{ p,- } .
	\Ee
	For $M_{1}$, we recall $\mathcal{C}_{+}$, $\mathcal{C}_{-}$ from (\ref{map+}), (\ref{map-}) and consider their composition 
	\Be
	\begin{split}
		&\mathcal{C}_{-}^{-1} \circ \mathcal{C}_{+} : \{ (t,x,v ) \in [0, T' ] \times \gamma_+ : t \ge \tb(t,x,v) \} \mapsto \{ (s,x,v) \in [0, T' ) \times \gamma_- : T'  \ge s + \tf(s,x,v) \},  \\
		&\mathcal{C}_{-}^{-1} \circ \mathcal{C}_{+}(t,x,v ) = (t -\tb(t,x,v) , \xb, \vb ). 	   
	\end{split}
	\Ee
	$\mathcal{C}_{-}^{-1} \circ \mathcal{C}_{+}$ is well defined as $\tf( t -\tb ,\xb ,\vb ) + (t - \tb) = \tb + t - \tb = t \le T'$ and
	\[
	( X(s;t,x,v), V(s;t,x,v) ) = ( X(s; t-\tb, \xb, \vb), V(s; t-\b, \xb, \vb) ).
	\]
	$\mathcal{C}_{-}^{-1} \circ \mathcal{C}_{+}$ is injective as the characteristic flow is unique and for any $(s,x,v ) \in [0, T') \times \gamma_-$ such that $T'  \ge s + \tf(s,x,v)$,
	\[
	\mathcal{C}_3 (s + \tf, X(s + \tf; s,x,v ) , V(s + \tf; s,x,v )) = (s,x,v ),
	\]
	implies that $\mathcal{C}_{-}^{-1} \circ \mathcal{C}_{+}$ is surjective. Therefore $\mathcal{C}_{-}^{-1} \circ \mathcal{C}_{+}$ is a bijection.  \\
	
	If we use $\eta_{2}$ to denote local parametrization near $\xb(t,x,v)$, we can apply (\ref{cov2}) to compute Jacobian of $\mathcal{C}_{-}^{-1} \circ \mathcal{C}_{+}$,
	\Be \notag
	\begin{split}
		\det J_{\mathcal{C}_{-}^{-1} \circ \mathcal{C}_{+}} &:= \det \frac{\p(t -\tb(t,x,v) , \xb, \vb )}{\p(t,x,v )}  \\
		&= \det \frac{\p( X(s;t,x,v), V(s;t,x,v) )}{\p(t,x,v )} 
		\det \frac{\p(t -\tb(t,x,v) , \xb, \vb )}{\p( X(s; t-\tb, \xb, \vb), V(s; t-\b, \xb, \vb) )}  \\
		&= v \cdot (\p_1 \eta(\bar{x}) \times \p_2 \eta(\bar{x}) ) \times \frac{1}{ \vb \cdot (\p_1 \eta_{2}(\overline{\xb}) \times \p_2 \eta_{2}(\overline{\xb}) ) }.
	\end{split}
	\Ee
	Thus we get
	\begin{equation} \label{green3}
	\begin{split}
	M_{1} &:= \int_0^{T' } \int_{\gamma_+ } |f|^p ( t- \tb(t,x,v), \xb,\vb) \textbf{1}_{ \{ t  \ge \tb(t,x,v) \} } d\gamma dt  \\
	&= \int_0^{T' } \int_{\gamma_- } |f|^p (t,x,v) \textbf{1}_{ \{T'   \ge s +  \tf(s,x,v) \} } d\gamma ds.  
	\end{split}
	\end{equation}

	\hide
	
	Suppose without lose of generality that locally we have $x = \eta_1 (x_1, x_2)$ and $\xb = \eta_2 ( \xba, \xbb$). Then the change of variable is:
	\[
	(t,x_1,x_2, v ) \mapsto ( t -\tb (t,\eta_1 (x_1,x_2) , v ) , \eta_2 ^{-1} (X(t -\tb; t ,\eta_1(x_1,x_2) ,v)), V(t-\tb; t, \eta_1(x_1,x_2), v) )
	\]

	We compute the Jacobian matrix $J$ of this change of variable:
	\Be
	\begin{split}\label{XVtb_txv}
		J = &\frac{\p (  t -\tb, \eta_2 ^{-1} (X(t -\tb; t ,\eta_1(\bar x ) ,v)), V(t-\tb; t, \eta_1(\bar x), v)  )}{\p (t, \bar x, v)}\\
		= & \ \begin{bmatrix}
			1 - \partial_t \tb & - \partial_{\bar x} \tb  & - \nabla_v \tb \\
			\nabla_x (\eta_2 ^{-1} ) \cdot ( \partial_s X + \partial_t X) & \nabla_x (\eta_2 ^{-1} ) \cdot ( \nabla_x X \cdot \partial_{\bar x} \eta_1 - \partial_s X \cdot \partial_{\bar x } \tb )   & \nabla_x (\eta_2 ^{-1} ) \cdot ( \nabla_v X  - \partial_s X \cdot \nabla_v \tb ) \\
			\partial_s V + \partial_t V & \nabla_x V \cdot \partial_{\bar x } \eta_1 - \partial_s V \cdot \partial_{\bar x } \tb & \nabla_v V - \partial_s V \cdot \nabla_v \tb 
		\end{bmatrix}\\
	\end{split}
	\Ee
	
	Let
	\[ 
	A =
	\begin{bmatrix}
	-\partial_s X & \partial_{\bar x } \eta_2 & 0_{3 \times 3} \\
	-\partial_s V &  0_{3 \times 2} & \mathrm{Id}_{3 \times 3}
	\end{bmatrix}
	\]
	
	Then we have
	\tiny
	\[ \begin{split}
	& A \cdot  J \\  =  & 
	\begin{bmatrix}
	-\partial_s X & \partial_{\bar x } \eta & 0_{3 \times 3} \\
	-\partial_s V &  0_{3 \times 2} & \mathrm{Id}_{3 \times 3}
	\end{bmatrix}
	\cdot
	\begin{bmatrix}
	1 & - \partial_{\bar x} \tb  & - \nabla_v \tb \\
	\nabla_x (\eta_2 ^{-1} ) \cdot ( \partial_s X + \partial_t X) & \nabla_x (\eta_2 ^{-1} ) \cdot ( \nabla_x X \cdot \partial_{\bar x} \eta_1 - \partial_s X \cdot \partial_{\bar x } \tb )   & \nabla_x (\eta_2 ^{-1} ) \cdot ( \nabla_v X  - \partial_s X \cdot \nabla_v \tb ) \\
	\partial_s V + \partial_t V  & \nabla_x V \cdot \partial_{\bar x } \eta_1 - \partial_s V \cdot \partial_{\bar x } \tb & \nabla_v V - \partial_s V \cdot \nabla_v \tb
	\end{bmatrix}
	\\ 
	= & \begin{bmatrix}
	- \partial_s X + \partial_{\bar x } \eta_2 \cdot \nabla_x( \eta_2 ^{-1} ) ( \partial_s X + \partial_t X) &  \partial_s X \cdot \partial_{\bar x } \tb + \partial_{\bar x } \eta_2 \cdot  \nabla_x (\eta_2 ^{-1} ) \cdot ( \nabla_x X \cdot \partial_{\bar x} \eta_1 - \partial_s X \cdot \partial_{\bar x } \tb )   & \partial_s X \cdot \nabla_v \tb + \partial_{\bar x} \eta_2 \cdot \nabla_x (\eta_2 ^{-1} ) \cdot ( \nabla_v X  - \partial_s X \cdot \nabla_v \tb )   \\
	-\partial_s V + \partial_s V + \partial_t V & \partial_s V \cdot \partial_{\bar x} \tb + \nabla_x V \cdot \partial_{\bar x } \eta_1 - \partial_s V \cdot \partial_{\bar x } \tb & \partial_s V \cdot \nabla_v \tb + \nabla_v V - \partial_s V \cdot \nabla_v \tb
	\end{bmatrix}
	\\
	= & \begin{bmatrix}
	- \partial_t X & \nabla_x X \cdot \partial_{\bar x } \eta_1 & \nabla_v X \\
	-\partial_t V & \nabla_x V \cdot \partial_{ \bar x } \eta_1 & \nabla_v V 
	\end{bmatrix}
	\end{split} \]
	\normalsize
	
	Since
	\[
	\partial_{\bar x }  \eta_2 \cdot \nabla_x (\eta_2 ^{-1} ) = \nabla_x ( \eta_2 \circ \eta_2 ^{-1} )  = \mathrm{Id}_{3\times 3 }.
	\]
	
	And by the computation of previous change of variable (\ref{cov2}) we have
	\[ \begin{split}
	\det  (A \cdot  J )   = & \det
	\begin{bmatrix}
	- \partial_s X & \nabla_x X \cdot \partial_{\bar x } \eta_1 & \nabla_v X \\
	-\partial_s V & \nabla_x V \cdot \partial_{ \bar x } \eta_1 & \nabla_v V 
	\end{bmatrix}
	\\ = &  \ \det \begin{bmatrix}
	-v & 
	\p_{\bar{x}} \eta_1   & 0_{3 \times 3} \\
	-E
	& 
	0_{3 \times 2} & \mathrm{Id}_{3\times 3}
	\end{bmatrix}\\
	\\ =&  - v \cdot (\partial_1 \eta_1 (\bar x ) \times \partial_2 \eta_1 (\bar x ) 
	\end{split} \]
	
	Since
	\[ \begin{split}
	& \det  (A  ) 
	\\  = & \det
	\begin{bmatrix}
	-\partial_s X & \partial_{\bar x } \eta & 0_{3 \times 3} \\
	-\partial_s V &  0_{3 \times 2} & \mathrm{Id}_{3 \times 3}
	\end{bmatrix}
	\\ =&  - \vb \cdot (\partial_1 \eta_2 (\overline {\xb} ) \times \partial_2 \eta_2 (\overline {\xb} ) 
	\end{split} \]
	
	Therefore
	\[
	\det(J ) = \frac{v \cdot (\partial_1 \eta_1 (\bar x ) \times \partial_2 \eta_1 (\bar x ) }{\vb \cdot (\partial_1 \eta_2 (\overline {\xb} ) \times \partial_2 \eta_2 (\overline {\xb} ) }
	\]
	So
	\begin{equation} \label{green3}
	\int_0^{T' } \int_{\gamma_- } |f|^p (t,x,v) \textbf{1}_{ \{T'   \ge s +  \tf(s,x,v) \} } d\gamma ds = \int_0^{T' } \int_{\gamma_+ } |f|^p ( t- \tb(t,x,v), \xb,\vb) \textbf{1}_{ \{ t  \ge \tb(t,x,v) \} } d\gamma dt
	\end{equation}
	
	\unhide
	\noindent For $M_{2}$, we recall $\mathcal{C}_{-}$ in (\ref{map-}) and choose $s=T' - \tb(T',x,v)$, 
	\Be
	\begin{split}
		&\mathcal{C}_{-} : \{ (x,v) \in \Omega \times \mathbb R^3 : T' \ge \tb(T',x,v) \} \mapsto \{ (s,x,v ) \in [0, T') \times \gamma_- : T'  < s + \tf(s,x,v) \},  \\
		&\mathcal{C}_{-}(x,v) = (T' - \tb(T',x,v) ,\xb ,\vb ).
	\end{split}
	\Ee
	$\mathcal{C}_{-}$ is well defined as $\tf(T' - \tb ,\xb ,\vb ) + (T' - \tb )  > \tb +(T' - \tb) = T$ as $x \in \Omega$ is in the interior. $\mathcal{C}_{-}$ is injective as the characteristic flow is unique and for any $\{ (s,x,v ) \in [0, T') \times \gamma_-$ such that $T'  < s + \tf(s,x,v) $, 
	\[
	\mathcal{C}_{-} ( X( T'; s, x,v ) , V(T';s,x,v ) ) = (s,x,v ),
	\] implies that $\mathcal{C}_{-}$ is surjective. Therefore $\mathcal{C}_{-}$ is a bijection and using (\ref{cov2}),
	\begin{equation} \label{green4}
	\begin{split}
	M_{2} &:= \iint_{ \Omega \times \mathbb R^3 } |f |^p (T' - \tb(T',x,v) , \xb, \vb ) \textbf{1}_{ \{ T' \ge \tb(T',x,v) \}} dxdv  \\
	&= \int_0^{T' } \int_{\gamma_- } |f|^p (t,x,v) \textbf{1}_{ \{T'   < s + \tf(s,x,v) \} }   d\gamma dt .
	\end{split}
	\end{equation}
	From (\ref{green3}) and (\ref{green4}), we prove (\ref{M12}). Finally, we apply (\ref{O12}) and (\ref{M12}) to (\ref{green}) conclude (\ref{Greedid}) and $f \in C^0( [0, T] ; L^p (\Omega \times \mathbb R^3 ) $, $f_{\gamma_+ } \in L^p ([0,T]; L^p (\gamma) )$.  \unhide

We define 
\Be\label{nongrazing_e}
\gamma_+^\e : = \{(x,v) \in \gamma_+: |n(x) \cdot v| \leq \e   \text{ or } |v|\geq1/ \e   \}.
\Ee

\begin{lemma} \label{le:ukai} 
	Assume 
	that, for $\Lambda_1>0$, $\delta_1>0$, 
	\begin{equation}
	\label{decay_E}
	\sup_{t \geq 0} e^{\Lambda_1 t} \| E(t) \|_{\infty} \leq \delta_1 \ll1 .
	\end{equation}
	We also assume $\frac{1}{C}\langle v\rangle \leq \psi(t,x,v)\leq C \langle v\rangle$ for some $C>0$. For $\varepsilon$ satisfying
	\Be\label{lower_bound_e}
	\e> \frac{2\delta_1}{  \Lambda_1}>0,
	\Ee
	there exists a
	constant $C_{\delta_1, \Lambda_1 ,\Omega }>0$ 
	for all $t\geq 0$, such that 
	\begin{equation} \label{case:decay}
	\begin{split}
	&\int_{0}^{t}\int_{\gamma _{+}\setminus \gamma _{+}^{\varepsilon }}|h|\mathrm{%
		d}\gamma \mathrm{d}s\\
	&\leq C_{\delta_1, \Lambda_1 ,\Omega }\bigg\{  
	||h_{0}||_{1}+\int_{0}^{t}  \| h(s)\|_{1} \dd s\\
	& \ \ \ \   \ \  \ \ \ \ \ \ \ \ \ +
	\int_{0}^{t} 
	\big{\Vert} 
	[
	\partial
	_{t}+v\cdot \nabla _{x}+E \cdot \nabla_v + \psi ]h(s)\big{\Vert} _{1} \mathrm{d}s  \bigg\}.\end{split}
	\end{equation}
	If $E \in L^\infty$ does not decay but
	\begin{equation} \label{nondecay}
	\| E (t) \|_{\infty} \leq \delta,
	\end{equation}
	then for $\varepsilon > 0$,
	\begin{equation} \label{case:nondecay} 
	\begin{split}
	&\int_{0}^{t}\int_{\gamma _{+}\setminus \gamma _{+}^{\varepsilon }}|h|\mathrm{%
		d}\gamma \mathrm{d}s\\
	&\leq C_{\delta, t, \varepsilon, \Omega }
	\bigg\{  
	||h_{0}||_{1}+\int_{0}^{t}  \| h(s)\|_{1} \dd s\\
	& \ \ \ \   \ \  \ \ \ \ \ \ \ \ \ +
	\int_{0}^{t} 
	\big{\Vert} 
	[
	\partial
	_{t}+v\cdot \nabla _{x}+E \cdot \nabla_v + \psi ]h(s)\big{\Vert} _{1} \mathrm{d}s  \bigg\},
	\end{split}
	\end{equation}
	where we have time-dependent constant $C_{\delta,t, \varepsilon, \O}>0$. \hide Furthermore, for any $(t,x,v)$ in $[0, \infty)\times \Omega \times \mathbb{R}^{3}$
	the function $
	h
	(t+s^{\prime },
	X(t+s^\prime;t,x,v)
	,
	V(t+s^\prime;t,x,v))
	$ is absolutely continuous in
	$s^{\prime } \in (-\min \{t_{\mathbf{b}}(t,x,v),t\}, t_{%
		\mathbf{f}}(t,x, v) )$.\unhide
\end{lemma}
\begin{proof}
	For $t-\tb(t,x,v) \leq s \leq t$, from (\ref{decay_E}),
	\Be \label{diffV}
	\begin{split} 
		|V(s;t,x,v)-v| &\leq \int_{s}^{t} |  E (\tau, X(\tau;t,x,v))| \dd \tau  
		\leq   \frac{\delta_1}{\Lambda_1} .
	\end{split}
	\Ee
	Therefore, for $|v| \geq \e$ with the condition (\ref{lower_bound_e}),
	\Be
	\begin{split}\notag
		|X(t-\tb;t,x,v)-x| &\geq  \left|\int^t_{t-\tb} v\right| -  \int^t_{t-\tb}| V(\tau;t,x,v)-v | \dd \tau  \\
		&\geq |v| \tb(t,x,v) -\frac{\delta_1}{\Lambda_1} \tb(t,x,v) \geq \frac{\delta_1}{\Lambda_1}\tb(t,x,v).
	\end{split}
	\Ee
	This gives an upper bound as
	\Be\label{upper_tb}
	\tb(t,x,v) \leq \frac{\Lambda_1}{\delta_1} \times \mathrm{diam} (\O)  \ \ \text{for } |v|\geq \e .
	\Ee
	
	Now we consider a lower bound of $\tb(t,x,v)$ for $(t,x,v) \in [0, \infty) \times \gamma_{+} \backslash \gamma_{+}^{\varepsilon} $. Since
	\Be \label{limit}
	\lim_{ \substack{
			y \to x_1 \\ y\in \partial \Omega }} \frac{ | (x_1 - y ) \cdot n(x_1 ) | }{ |x_1 - y | } = 0,\quad \text{for}\quad x_{1} \in \p\O,
	\Ee
	we have $| (x_1 - y ) \cdot n(x_1) | \le C_{\Omega} |x_1 - y|^{2} $ for all $x_1, y \in \partial \Omega$. Thus from (\ref{diffV}),
	\Be \label{distb}
	\begin{split}
		C_{\Omega}|x - \xb |^2 &\geq | (\xb - x ) \cdot n(x) | = \left| \int_{t - \tb} ^t  V(\tau;t,x,v )  \dd\tau \cdot n(x) \right| 
		\\ 
		&\geq \left| \int_{t -\tb}^t v  \cdot n(x)  \dd \tau\right| -   \int_{t-\tb} ^t |v - V(\tau;t,x,v) | \dd\tau     \\ 
		&\geq\Big\{ | v \cdot n(x) | - \frac{\delta_1}{\Lambda_1} \Big\}  \tb>  \frac{\varepsilon}{2} \tb .  \\ 
	\end{split} 
	\Ee
	Also we have
	\[ \begin{split}
	|x - \xb |^2 &= \left| \tb v - \int_{t -\tb} ^t ( v - V(s) ) \dd s \right| ^2  \\ 
	& \le 2 \tb ^2 |v|^2 + 2 \left( \int_{t -\tb} ^t |v - V(s) | \dd s \right) ^2  \\ 
	& \le 2 \tb^2 |v|^2 + 2 (\tb \frac{ \delta_1}{\Lambda_1})^2
	< 2 \tb^2 |v| ^2 + \frac{1}{2} \tb^{2} \varepsilon ^2.
	\end{split} \]
	Combining above two estimates, we get 
	$
	C_\Omega (  2 \tb^2 |v| ^2 + \frac{1}{2} \tb^2 \varepsilon ^2) > \frac{\varepsilon}{2} \tb,
	$
	and dividing $|v|^2 \tb$ on both sides we get
	\Be \label{lower-tb}
	\tb(t,x,v) > \frac{1}{5 C_\Omega} \varepsilon^{3}.
	\Ee

	If $h$ solves (\ref{transport_E}), then for $(t,x,v) \in [0,T] \times \gamma_+ $ and $ - \min \{ t ,\tb (t,x,v) \}  \leq s\leq 0$,  
	\begin{equation} \label{h_expand}
	\begin{split}
	h(t,x,v) &= h(t+s, X(t+s),V(t+s) ) e^{ -\int_{t-s}^{t} \psi (\tau',X(\tau'),V(\tau')) \dd\tau'}  \\ 
	& +  \int_{t+s}^{t} e^{-\int_{\tau}^{t} \psi(\tau^{\prime},X(\tau^{\prime}),V(\tau^{\prime})) \dd\tau' }
	H(\tau,X(\tau;t,x,v), V(\tau;t,x,v))  \dd \tau,
	\end{split}
	\end{equation}
	where $X(s) = X(s; t ,x,v ),V(s) = V(s; t,x,v)$.
	
	\hide
	
	\begin{equation} \label{h_expand}
	\begin{split}
	h(t,x,v) &= h(t-s, X(t-s),V(t-s) ) e^{\int_0^{-s} \psi (V(t + \tau' )) d\tau'}  \\ 
	&\quad +  \int_{-s}^0 e^{\int_0^\tau \phi(V(t + \tau' )) d\tau' } \big[ \partial_t h + V(t + \tau ) \cdot \nabla_x h  \\
	&\quad + E(X(t + \tau )) \cdot \nabla_v h + \phi (V(t + \tau)) h \big]\Big\vert_{( t + \tau , X(t + \tau ) , V(t + \tau ))} d \tau,
	\end{split}
	\end{equation}
	
	{\color{red}
		\begin{equation} \label{h_expand}
		\begin{split}
		h(t,x,v) &= h(t-s, X(t-s),V(t-s) ) e^{\int_0^{-s} \phi (V(t + \tau' )) d\tau'}  \\ 
		&\quad +  \int_{-s}^0 e^{\int_0^\tau \phi(V(t + \tau' )) d\tau' } \big[ \partial_t h + V(t + \tau ) \cdot \nabla_x h  \\
		&\quad + E(X(t + \tau )) \cdot \nabla_v h + \phi (V(t + \tau)) h \big]\Big\vert_{( t + \tau , X(t + \tau ) , V(t + \tau ))} d \tau,
		\end{split}
		\end{equation}
		where we abbreviated $X(t + \tau ) = X(t + \tau; t ,x,v ) $, and $V(t + \tau ) = V( t + \tau; t,x,v  )$.  \\
	} 
	
	\unhide
	
	Then
	\Be \label{expan_int}
	\begin{split}
		& \min \{ t, \tb(t,x,v) \}  \times |h (t,x,v)|\\
		=  & \int_{ - \min \{ t, \tb(t,x,v) \} } ^0 |h(t,x,v)| \dd s  \\
		\leq &  \int_{ - \min \{ t, \tb(t,x,v) \} } ^0 |h(t + s, X(t+s), V(t+s) ) | \dd s   \\
		&  +  \int_{ - \min \{ t, \tb(t,x,v) \} } ^0 \int_{t+s}^{t} |H(\tau,X(\tau;t,x,v), V(\tau;t,x,v)) |  \dd \tau \dd s  \\ 
		\leq& \int_{ - \min \{ t, \tb(t,x,v) \} }^{0} |h(t+s, X(t+s),V(t+s) ) | \dd s  \\ 
		&  +  \tb(t,x,v) \int_{ t - \min \{ t, \tb(t,x,v) \}}^{t} |H(\tau,X(\tau;t,x,v), V(\tau;t,x,v)) |  \dd \tau .  
	\end{split} 
	\Ee

	\hide
	
	\[ \begin{split}
	\min & \{ t, \tb(t,x,v) \}  \times |h (t,x,v)| = \int_{ - \min \{ t, \tb(t,x,v) \} } ^0 |h(t,x,v)| ds
	\\   \le & \int_{ - \min \{ t, \tb(t,x,v) \} } ^0 |h(t + s, X(t+s),V(t+s) ) | ds 
	\\ & +  \int_{ - \min \{ t, \tb(t,x,v) \} } ^0 \int_{s}^0 | \left[ \partial_t h + V(t + \tau ) \cdot \nabla_x h + E(X(t + \tau )) \cdot \nabla_v h + \phi (V(t + \tau)) h \right ] 
	\\ & ( t + \tau , X(t + \tau ) , V(t + \tau )) | d \tau ds
	\\ \le & \int_{ - \min \{ t, \tb(t,x,v) \} } ^0 |h(t+s, X(t+s),V(t+s) ) | ds 
	\\ & +  \tb(t,x,v) \times \int_{ - \min \{ t, \tb(t,x,v) \} } ^0  | \left[ \partial_t h + V(t + \tau ) \cdot \nabla_x h + E(X(t + \tau )) \cdot \nabla_v h + \phi (V(t + \tau)) h \right ] 
	\\ & ( t + \tau , X(t + \tau ) , V(t + \tau )) | d \tau 
	\end{split} \]
	
	\unhide

	Let $\epsilon_1 := \frac{1}{5 C_\Omega} \varepsilon ^3 $. Then from (\ref{lower-tb}), for $(t,x,v) \in [\epsilon_1 ,T] \times \gamma_+ \backslash \gamma_+^\varepsilon $, we integrate over $\int_{\epsilon_1} ^T \int_{\gamma_+ \backslash \gamma_+^\varepsilon} $ to get 
	\Be \label{T away 0}
	\begin{split}
		& \epsilon_1  \int_{\epsilon_1} ^T \int_{\gamma_+ \backslash \gamma_+^\varepsilon} | h(t,x,v) | \dd\gamma \dd t  \\ 
		\leq 	&  \int_{\epsilon_1} ^T \int_{\gamma_+ \backslash \gamma_+^\varepsilon}   \displaystyle {\min_{ [\epsilon_1 , T ] \times [ \gamma_+ \backslash \gamma_+^\varepsilon ]} } \{ t, \tb (t,x,v) \} \times | h(t,x,v) |   \dd\gamma \dd t \\ 
		\leq&  \int_0 ^T \int_{\gamma_+ \backslash \gamma_+^\varepsilon } \int_{-\min\{t, \tb (t,x,v) \} } ^0 |h(t+s, X(t+s) ,V(t+s)) | \dd s \dd\gamma \dd t  \\ 
		&	 +\sup_{\gamma_+ \backslash \gamma_+^\e} \tb(t,x,v)    \int_0^T \int_{\gamma_+ \backslash \gamma_+ ^\varepsilon } \int^{t }_{t - \min \{ t, \tb(t,x,v) \}} \big| H (\tau,X(\tau ), V(\tau ))\big| 
		 \\
		\lesssim&   \int_0^T \| h  \|_1  + \int_0^T \big\| [\partial_t + v \cdot \nabla_x - \nabla\phi \cdot \nabla_v + \psi ] h  \big\| _1   ,
	\end{split}
	\Ee
	where we used Lemma \ref{int_id} and (\ref{upper_tb}) in the last inequality.  
	
	On the other hand, from our choice $\varepsilon$ and $\epsilon_1$, 
	\Be\label{tb>t}
	\tb(t,x,v) > t \  \text{ for all } (t,x,v) \in [0, \epsilon_1] \times \gamma_+ \backslash \gamma_+^\varepsilon . 
	\Ee
	Integrating  
	\[
	\begin{split}
	|h(t,x,v)|   &\leq | h_0( X(0),V(0) )|   +  \int_{-t}^0 \big|H ( t + \tau , X(t + \tau ) , V(t + \tau )) \big|\dd \tau
	\end{split}
	\] 
	over $\int_{0}^{\epsilon_1} \int_{\gamma_+ \backslash \gamma_+^\varepsilon} $, we get 
	\Be \label{T near}
	\begin{split}
		&\int_{0}^{\epsilon_1} \int_{\gamma_+ \backslash \gamma_+^\varepsilon} |h(t,x,v)| \\
		=&  \int_{0}^{\epsilon_1} \int_{\gamma_+ \backslash \gamma_+^\varepsilon} | h_0( X(0),V(0) )| \dd\gamma \dd t  \\ 
		&+ \int_{0}^{\epsilon_1} \int_{\gamma_+ \backslash \gamma_+^\varepsilon} \int_{-t}^0 \big|  H( t + \tau , X(t + \tau ) , V(t + \tau )) \big|\dd \tau \dd\gamma \dd t.
	\end{split}
	\Ee 
	Applying the change of variables of (\ref{map+}), the first term in RHS of (\ref{T near}) is bounded by $\| h_0 \|_1$. From (\ref{tb>t}) and (\ref{map+}), the second term in RHS is bounded by 
	%
	\[
	\begin{split}
	\int_0^{\epsilon_1} \| (\partial_t + v \cdot \nabla_x - \nabla\phi \cdot \nabla_v + \psi ) h \| _1 \dd t.
	\end{split}
	\]

	\hide
	for any $(t,x,v) \in [0, \epsilon_1] \times \gamma_+ \backslash \gamma_+^\epsilon$ with $\tb(t,x,v) > \epsilon_{1} \geq t$, 
	\[
	\begin{split}
	&\int_0^{\epsilon_1} \int_{\gamma_+ \backslash \gamma_+^\varepsilon} |h_0(X(0;t,x,v), V(0;t,x,v)) | d\gamma dt \\ 
	&\leq \int_0^{T } \int_{ \gamma_+ } |h_0| (X(0;t,x,v), V(0;t,x,v) ) \textbf{1}_{ \{ t < \tb(t,x,v) \} }    d\gamma d t   \\
	&= \iint_{\Omega \times \mathbb R^3 } |h_0 | \textbf{1}_{ \{ \tf(0,x,v) \le T \} }   dx dv  
	\leq \| h_0 \|_1.
	\end{split}
	\]
	Applying above two estimates to (\ref{T near}), 
	\Be \label{T near 0}
	\int_0 ^{\epsilon_1} \int_{\gamma_+ \backslash \gamma_+^\varepsilon } | h(t,x,v) | d\gamma dt \lesssim_{\Omega, \varepsilon} \| h_0 \|_1 + \int_0^{\epsilon_1} \| [\partial_t + v \cdot \nabla_x - \nabla\phi \cdot \nabla_v + \psi ] h \| _1 dt.
	\Ee\unhide
	Finally we combine (\ref{T near}) and (\ref{T away 0}) to obtain (\ref{case:decay}).

	The proof of (\ref{case:nondecay}) is similar. In the last line of (\ref{expan_int}), we have $t$ multiplier instead of $\tb(t,x,v)$. The rest of proof is same.\hide
	
	Next we prove (\ref{case:nondecay}) in the case of (\ref{nondecay}).
	First we compute upper bound of $\tb(t,x,v)$. When backward trajectory $V(s;t,x,v)$ is well-defined for $t-\tb(t,x,v) \leq s \leq t$, from (\ref{nondecay}), $|V(s;t,x,v)-v| \leq \delta(t-s)$. \\
	
	Now we give an upper bound of $\tb(t,x,v)$ for $(t,x,v) \in [0, \infty) \times \gamma_{+} \backslash \gamma_{+}^{\varepsilon} $. From (\ref{limit}),
	we have $| (x_1 - y ) \cdot n(x_1) | \le C_{\Omega} |x_1 - y|^{2} $ for all $ y \in \partial \Omega$. Similar as (\ref{distb}),
	\[ \begin{split}
	C_{\Omega}|x - \xb |^2 &\geq | (\xb - x ) \cdot n(x) | = | \int_{t - \tb} ^t  V(\tau;t,x,v )  d\tau \cdot n(x) | 
	\\ 
	&\geq \tb | v \cdot n(x) | - \int_{t -\tb}^t |v - V(\tau;t,x,v) | d\tau   
	\geq \varepsilon \tb - \frac{\delta}{2}\tb^{2} ,
	\end{split} \]
	and we also have
	\[ \begin{split}
	|x - \xb |^2 &= | \tb v - \int_{t -\tb} ^t ( v - V(s) ) ds | ^2 
	\le 2 \tb^2 |v|^2 + 2 ( \frac{\delta}{2}\tb^{2} )^2
	< \frac{2\tb^2}{\varepsilon^{2}} + \frac{\delta^{2}}{2} \tb^{4}.
	\end{split} \]
	Combining above two estimates, we get 
	\[
	\varepsilon  \leq \big( \frac{\delta}{2} + \frac{2}{\varepsilon^{2}} \big) \tb + \frac{\delta^{2}}{2} \tb^{3},\quad \tb > 0,
	\]
	and hence by monotone property, there exist $L_{\varepsilon,\delta}$ such that above inequality implies $\tb \geq L_{\varepsilon,\delta}$. Now we directly apply (\ref{expan_int}) with $\epsilon_{1} = L_{\varepsilon,\delta}$. Then we follow (\ref{T away 0}) to get
	\Be \label{N T away 0}
	\begin{split}
		& \epsilon_1  \int_{\epsilon_1} ^T \int_{\gamma_+ \backslash \gamma_+^\varepsilon} | h(t,x,v) | d\gamma dt  \leq  \int_{\epsilon_1} ^T \int_{\gamma_+ \backslash \gamma_+^\varepsilon} \Big( \displaystyle {\min_{ [\epsilon_1 , T ] \times [ \gamma_+ \backslash \gamma_+^\varepsilon ]} } \{ t, \tb (t,x,v) \} \times | h(t,x,v) | \Big) d\gamma dt \\ 
		& \leq \int_0 ^T \int_{\gamma_+ \backslash \gamma_+^\varepsilon } \int_{-\min\{t, \tb (t,x,v) \} } ^0 |h(t+s, X(t+s) ,V(t+s)) | ds d\gamma dt  \\ 
		& + t \int_0^T \int_{\gamma_+ \backslash \gamma_+ ^\varepsilon } \int_{- \min \{ t, \tb (t,x,v) \}  }^0 \big| \big( \p_{t} + v\cdot\nabla_{x} - \nabla\phi(t,x)\nabla_{v} + \psi(t,x,v) \big) \big|_{ (t+s,X(t+s;t,x,v), V(t+s;t,x,v)) } ds d\gamma dt  \\
		& \leq C_{t,\delta,\varepsilon,\O} \Big[ \int_0^T \| h(t) \|_1 dt + \int_0^T \big\| [\partial_t + v \cdot \nabla_x - \nabla\phi \cdot \nabla_v + \psi ] h(t) \big\| _1 dt \Big],
	\end{split}
	\Ee
	similar as (\ref{T away 0}). Note that we have $t$-dependent $C_{t,\delta,\varepsilon,\O}$ in above estimate since we do not have uniform upper bound of $\tb$ unlike to (\ref{T away 0}).  
	
	For $\int^{\epsilon_1}_{0} \int_{\gamma_+ \backslash \gamma_+^\varepsilon} | h(t,x,v) | d\gamma dt$, we can use exactly same estimate as decaying potential case with $\epsilon_{1} = L_{\varepsilon,\delta}$ to obtain (\ref{T near 0}). With (\ref{N T away 0}), we obtain (\ref{case:nondecay}). \unhide
\end{proof}

\begin{lemma}If $n(\xb(t,x,v)) \cdot \vb(t,x,v) \neq 0$ then $(\tb,\xb,\vb)$ is differentiable and 
	\Be\begin{split}\label{computation_tb_x}
		\frac{\p\tb}{\p x_i}  =& \  \frac{1}{n(\xb) \cdot \vb}n(\xb) \cdot  \left[
		e_i + \int^{t-\tb}_t \int^s_t \Big(\frac{\p X(\tau )}{\p x_i} \cdot \nabla\Big) E(\tau, X(\tau )) \dd \tau \dd s 
		\right]  ,\\
		\frac{ \p\xb}{\p x_i} = & \  e_i - \frac{\p \tb}{\p x_i} \vb + \int^{t-\tb}_{t} \int^s_t    \Big(\frac{\p X(\tau )}{\p x_i} \cdot \nabla\Big) E(\tau, X(\tau ))  \dd \tau \dd s,\\
		\frac{\p \vb}{\p x_i} = & \ - \frac{\p \tb}{\p x_i} E(t-\tb, \xb) + \int^{t-\tb}_t   \Big(\frac{\p X(\tau )}{\p x_i} \cdot \nabla\Big) E(\tau, X(\tau )) \dd \tau,\\
		\frac{\p\tb}{\p v_i}  =& \  \frac{1}{n(\xb) \cdot \vb}n(\xb) \cdot  \left[
		e_i + \int^{t-\tb}_t \int^s_t \Big(\frac{\p X(\tau )}{\p v_i} \cdot \nabla\Big) E(\tau, X(\tau )) \dd \tau \dd s 
		\right]  ,\\ 
		\frac{\p \xb}{\p v_i} = & \ - \tb e_i - \frac{
			\p \tb}{\p v_i} \vb + \int^{t-\tb}_{t} \int^s_t    \Big(\frac{\p X(\tau )}{\p v_i} \cdot \nabla\Big) E(\tau, X(\tau ))  \dd \tau \dd s ,\\
		\frac{\p \vb}{\p v_i} = & \ e_i - \frac{\p \tb}{\p v_i} E(t-\tb, \xb) + \int^{t-\tb}_t
		\Big(\frac{\p X(\tau )}{\p v_i} \cdot \nabla\Big) E(\tau, X(\tau ))  \dd \tau.
	\end{split} \Ee
\end{lemma}
\begin{proof}
	The equalities are derived from direct computations and an implicit function theorem. For details see \cite{KL1}.
\end{proof}

\begin{proposition}
	Assume the compatibility condition
	\Be\label{compatibility_inflow}
	f_0(x,v) = g(0,x,v)\quad  \text{for} \quad (x,v) \in \gamma_- .
	\Ee
	Let $p \in [1, \infty )$ and $0 < \vartheta < 1/4$. \hide
	
	Define $D^\epsilon : = \{ (x,v) \in \Omega \times \mathbb R^3 : (\xb, \vb) \in \gamma_- \setminus \gamma_-^\epsilon \} .$ \unhide
	Assume
	\Be \begin{split}\label{assumption_inflow}
		\nabla_x f_0 , \nabla_v f_0    \in L^p (\Omega \times \mathbb R^3 ),
		\\ \nabla_{x,v} \tb  \partial_t g, \nabla_{x,v} \vb \nabla_v g, \nabla_{x,v} \xb \partial_{\xb } g , \nabla_{x,v} \tb \psi g \in L^p ( [0, T] \times \gamma_- ) ,
		\\ \hide\frac{ n }{ n \cdot \vb } \Big \{ \partial_t g + \sum_{i =1}^2 ( \vb \cdot \tau_i ) \partial_{\tau_i } g  + \nu g - H +  E \cdot \nabla_v g \Big \} 
		\\ +  \frac{ n \cdot \iint \partial_x E }{ n \cdot \vb } \Big \{ \partial_t g + \sum_{i =1}^2 ( \vb \cdot \tau_i) \partial_{ \tau_i } g   - \nu g + H \Big \} \in L^p ([0, T ] \times \gamma_- ),
		\\ \unhide  \nabla_x H ,\nabla_v H \in L^p  ([0, T ] \times \Omega \times \mathbb R^3 ),
		\\ 
		e^{- \vartheta |v|^2 } \nabla_x \psi, e^{-\vartheta |v|^2 } \nabla_v \psi \in L^p ([ 0, T ] \times \Omega \times \mathbb R^3 ),
		\\ e^{\vartheta |v|^2 } f_0 \in L^\infty ( \Omega \times \mathbb R^3 ) , e^{\vartheta |v|^2 } g \in L^\infty ( [0, T] \times \gamma_- ),
		\\ e^{\vartheta |v|^2 } H \in L^\infty ([0, T] \times \Omega \times \mathbb R^3 ).
	\end{split} \Ee
	Then for any $T > 0$, there exists a unique solution $f$ to (\ref{transport_E}) such that $\nabla_{x,v} f\in  C^0 ([0,T] ; L^p (\O \times \R^3)) \cap L^1((0,T); L^p (\gamma))$.
	
	\hide$f$ to 
	\[ \{ \partial_t + v \cdot \nabla_x + E \cdot \nabla_v + \nu \} f = H \]
	such that $f, \partial_t, \nabla_x f ,\nabla_v f \in C^0( [ 0, T ] ; L^p (\O \times \R^3) ) $ and their traces satisfy
	\[ \begin{split}
	\partial_t f|_{\gamma_-
	}  = \partial_t g, \nabla_v f |_{\gamma
	} = \nabla_v g, \nabla_x f|_{\gamma_- 
	} = \nabla_x g, \quad \text{on} \quad \gamma_- 
	,
	\\ \nabla_x f(0,x,v) = \nabla_x f_0, \nabla_v f(0,x,v) = \nabla_v f_0, \quad \text{in} \quad\O \times \R^3,
	\\ \partial_t f(0,x,v) = \partial_t f_0, \quad \text{in} \quad \O \times \R^3.
	\end{split} \]
	
	Moreover
	\Be\label{inflow_energy}
	\begin{split}
		&\| \nabla_{x,v} f(t) \|_p^p + \int^t_0 | \nabla_{x,v} f |_{+,p}^p \\
		= & \  \| \nabla_{x,v} f _0\|_p^p + \int^t_0 | \nabla_{x,v} g |_{-,p}^p\\
		&+ p \int^t_0 \iint_{\O \times \R^3} \{ \nabla_{ x,v} H - [\nabla_{x,v} v] \nabla_x f - [\nabla_{x,v} \psi] f \}
		|\nabla_{ x,v} f|^{p-2} \nabla_{ x,v}f.
	\end{split}
	\Ee

	\unhide\end{proposition}

\begin{proof}Along the characteristics
	\Be \begin{split}\label{H_trajectory}
		f(t,x,v)   
		&  =   \mathbf{1}_{t< \tb} e^{- \int_0^t \psi (\tau, X(\tau), V(\tau)) \dd \tau} f_0 (X(0), V(0) )  \\
		&+\mathbf{1}_{t< \tb}  \int_0^t e^{- \int _s^t \psi(\tau, X(\tau), V(\tau)) \dd \tau} H (s, X(s), V(s ) ) \dd s \\
		&+  \mathbf{1}_{t> \tb} e^{- \int_{t-\tb}^t \psi (\tau, X(\tau), V(\tau)) \dd \tau} g (t-\tb,X(t-\tb), V(t-\tb) )\\
		&+  \mathbf{1}_{t> \tb}\int_{t-\tb}^t e^{- \int _s^t \psi(\tau, X(\tau), V(\tau)) \dd \tau} H (s, X(s), V(s ) ) \dd s  
		,
	\end{split} 
	\Ee
	where $(X(s),V(s))= (X(s;t,x,v), V(s;t,x,v))$.
	%
	%
	%
	%
	\hide\[ \begin{split}
	f(t,x,v) = & f_0(X(0) , V(0) ) e^{-\int_0^t \nu }- \int_0^t \frac{d}{ds} \left( f(t-s, X(t-s) , V(t-s) ) e^{\int_0^s \nu } \right) ds
	\\ =  & f_0 (X(0), V(0) ) e^{- \int_0^t \nu } + \int_0^t e^{- \int _0^s \nu} H (t-s, X(t-s), V(t-s ) ) ds 
	\end{split} \]
	
	where $H = \{ \partial_t + v \cdot \nabla_x  + E \cdot \nabla_v  + \nu \} f $, $\nu = \nu(\tau') = \nu ( t - \tau' , X(t- \tau'), V(t- \tau' ) )$.
	
	And for $t > \tb$:
	\[ \begin{split}
	f(t,x,v) =  e^{-\int_0^ {\tb}  \nu} g(t -\tb, \xb,\vb ) + \int_0 ^ {\tb} e^{-\int_0^s \nu } H(t-s ,X(t-s), V(t-s) ) ds 
	\end{split} \]
	
	We can rewrite it as:
	
	\[ \begin{split}
	f(t,x,v) = & \textbf{1}_{ \{ t \le \tb \} } e^{-\int_0^t \nu } f_0 (X(0) ,V(0)) + \textbf{1}_{ \{ t > \tb \} } e^{-\int_0^ {\tb}  \nu} g(t -\tb, \xb,\vb ) \\ & + \int_0^{\min \{ \tb, t\} } e^{-\int_0^s \nu } H(t-s ,X(t-s), V(t-s) ) ds 
	\end{split} .\]\unhide


	By direct computations, we have\hide
	\Be \begin{split}
		\partial_t & f(t,x,v ) \textbf{1}_{\{ t \neq \tb \}} 
		\\ = & - \textbf{1}_{ \{ t < \tb \} } e^{- \int_0^t \psi } \left[ \nu f_0 + v \cdot \nabla_x f_0 + E \cdot \nabla_x f_0 - H_{| t= 0 }  + \int_0^t \partial_t \nu \times f_0 \right]  
		\\ &  + \textbf{1}_{\{ t > \tb \} } e^{- \int_0 ^ {\tb} \nu } \left[ \partial_t g - \int_0^{\tb} \partial_t \nu \times g \right] ( t - \tb, \xb,\vb )
		\\ & + \int_0^{ \min \{ t, \tb \} } e^{-\int_0 ^s  \nu } \left[ \partial_t H - \int_0^s \partial_t \nu \times H \right] ( t-s, X(t-s), V(t-s ) ) ds.
	\end{split} \]\unhide
	 \Be \begin{split}\notag
			&\nabla_{x,v}   f(t,x,v ) \textbf{1}_{ \{ t\neq \tb \}} 
			\\
			= &    \   \textbf{1}_{\{ t < \tb \}} e^{-\int_0^t \psi } \Big\{ \nabla_{x,v} X(0) \cdot \nabla_x f_0  +\nabla_{x,v} V (0)\cdot  \nabla_v f_0 \\
			&   \ \ \ \ \  \ \ \ \ \ \ \ \ \    - f_0  \int_0^t ( \nabla_{x,v} X (0)\cdot \nabla_{x} \psi  +\nabla_{x,v} V(0) \cdot \nabla_v \psi )   \Big\} 
			\\ +&       \textbf{1}_{\{ t > \tb \} } e^{- \int^t_{t-\tb} \psi } 
			\Big\{
			- \nabla_{x,v} \tb  \p_t g(t-\tb) + \nabla_{x,v} \xb \cdot \nabla_x g(t-\tb) \\
			& \ \ \ \ \ \ \ \ \ \      + \nabla_{x,v} \vb \cdot \nabla_v g(t-\tb)
			- \nabla_{x,v} \tb \psi(t-\tb)  g (t-\tb)\\
			& \ \ \ \ \ \ \ \ \ \       -g(t-\tb) \int^t_{t-\tb} (\nabla_{x,v} X(\tau) \cdot \nabla_x  \psi(\tau) + \nabla_{x,v} V(\tau) \cdot \nabla_v \psi(\tau) ) \dd \tau
			\Big\}
			\\ +&   \int^t_{ \max \{ 0, t-\tb \}} e^{- \int_s ^t \psi }  \Big\{\nabla_{x,v} X(s)\cdot   \nabla_x H(s)   +\nabla_{x,v}  V(s)\cdot  \nabla_v H (s)\\
			&\qquad \qquad   \ \ - H(s)\int^t_s (  \nabla_{x,v}  X(\tau) \cdot \nabla_x \psi(\tau)  + \nabla_{x,v}  V (\tau)\cdot \nabla_v \psi(\tau)   ) \dd \tau \Big\}  \dd s,
		\end{split} \Ee 
	where $\nabla_{x,v} \tb, \nabla_{x,v} \xb, \nabla_{x,v} \vb$ in (\ref{computation_tb_x}).
	
	From (\ref{hamilton_ODE}) with replacing $- \nabla_x \phi_f$ by $E$,  
	\Be\label{nabla_Hamilton}
	\begin{split}
	&\frac{d}{ds} \left[ \begin{matrix}\nabla_{x,v} X(s;t,x,v) \\ \nabla_{x,v} V(s;t,x,v)\end{matrix} \right]
	\\= & \  \left[\begin{matrix} 
	0_{3\times 3} & \text{Id}_{3\times 3}\\
	\nabla_x E(s,X(s;t,x,v)) & 0_{3\times 3} 
	\end{matrix}\right]
	\left[ \begin{matrix}\nabla_{x,v} X(s;t,x,v) \\ \nabla_{x,v} V(s;t,x,v)\end{matrix} \right].
	\end{split}\Ee
	Then easily we have for $C=C(\| \nabla_x E \|_\infty^{1/2})$
	\Be
	|\nabla_{x,v} X(s;t,x,v)|+ 
	|\nabla_{x,v} V(s;t,x,v)| \lesssim e^{ C |t-s|}.\notag
	\Ee
	By the change of variables in Lemma \ref{lem COV} and Lemma \ref{int_id}, we have 
	\Be\label{uniform_nabla_f_l}
	\begin{split}
		&\| \nabla_{x,v} f (t) \mathbf{1}_{t \neq \tb} \|_{p}\\
		\lesssim_t & \ \|  \nabla_{x,v} f_0 \|_p  + \| e^{\vartheta |v|^2} f_0 \|_\infty \| e^{- \vartheta |v|^2} \nabla_{x,v} \psi \|_p\\
	& 	+ \bigg[\int^t_0 \| \nabla_{x,v} \tb \p_t g \|_{p}^p  + \| \nabla_{x,v} \xb \p_x g \|_{p}^p\\
		& \ \ \  \ 
	  + \| \nabla_{x,v} \vb \p_v g \|_{p}^p 	+ \| \nabla_{x,v} \tb \psi g \|_p^p\bigg]^{1/p}\\
		&+  \left[\int^t_0  \| \nabla_{x,v} H \|_p^p + \| e^{\vartheta|v|^{2}}H \|_\infty \| e^{-\vartheta|v|^{2}}\nabla_{x,v} \psi   \|_p^p\right]^{1/p}.
	\end{split}
	\Ee

	\hide

	\[ \begin{split}
	\nabla_v & f(t,x,v ) \textbf{1}_{\{ t\neq \tb \}}
	\\ = & \textbf{1}_{\{ t < \tb \} } e^{-\int_0^t \nu } \left[ \nabla_x f_0 \cdot \partial_v X + \nabla_v f_0 \cdot \partial_v V - \int_0^t (\nabla_x \nu \cdot \partial_v X + \nabla_v \nu \cdot \partial_v V) \times f_0 \right] (X(0), V(0))
	\\ & + \textbf{1}_{\{ t > \tb \}} e^{-\int_0^{\tb } \nu } \left\{ - \nu \nabla_v \tb g - \int_0^{\tb} (\nabla_x \nu \cdot \nabla_v X + \nabla_v \nu \cdot \partial_v V) \times g + \nabla_v g + \nabla_v \tb H \right\} (t - \tb, \xb,\vb )
	\\ & + \int^t_{\max \{ 0, t-\tb \}} e^{-\int_0^s \nu} \bigg [ \nabla_x H \cdot \partial_v X + \partial_v H \cdot \partial_v V
	\\& \qquad \qquad \qquad \qquad - \int_0^s (\nabla_x \nu \cdot \partial_v X + \nabla_v \nu \cdot \partial_v V) \times H \bigg ] (t-s , X(t-s),V(t -s)) ds
	\end{split} \]
	
	Now we compute $\nabla_x [g(t-\tb,\xb,\vb ) ]$ and $\nabla_v [g(t-\tb,\xb,\vb) ]$:
	
	\[ \begin{split}
	\nabla_x  & [ g(t -\tb, \xb, \vb ) ]  
	\\= &  - \nabla_x \tb \partial_t g + \nabla_x \xb \nabla_{\tau } g + \nabla_v g \nabla_x \vb
	\\  = & -\frac{n(\xb) }{ n(\xb) \cdot \vb } \partial_t g - \frac{ n(\xb) \cdot \iint \partial_x E }{ n(\xb ) \cdot \vb } \partial_t g 
	\\ & + \tau_1 \partial_{\tau_1} g + \tau_2 \partial_{\tau_2} g - \frac{ n(\xb ) } { n(\xb ) \cdot \vb } \left( \vb \cdot \tau_1 \partial_{\tau_1} g + \vb \cdot \tau_2 \partial_{\tau_2} g \right) - \frac{ n(\xb) \cdot \iint \partial_x E } {n (\xb) \cdot \vb } \left( \vb \cdot \tau_1 \partial_{\tau_1} g + \vb \cdot \tau_2 \partial_{\tau_2 } g \right)
	\\ & - \frac{n(\xb)}{n (\xb) \cdot \vb } ( E \cdot \nabla_v g ) - \nabla_v g \cdot \int \partial_x E 
	\end{split} \]
	
	\[ \begin{split}
	\nabla_v  & [ g(t -\tb, \xb, \vb ) ]  
	\\= &  - \nabla_v \tb \partial_t g + \nabla_v \xb \nabla_{\tau } g + \nabla_v g \nabla_v \vb
	\\  = & -\frac{\tb n(\xb) }{ n(\xb) \cdot \vb } \partial_t g - \frac{ n(\xb) \cdot \iint \partial_v E }{ n(\xb ) \cdot \vb } \partial_t g 
	\\ & - \tb ( \tau_1 \partial_{\tau_1} g +  \tau_2 \partial_{\tau_2} g) - \tb \frac{ n(\xb ) } { n(\xb ) \cdot \vb } \left( \vb \cdot \tau_1 \partial_{\tau_1} g + \vb \cdot \tau_2 \partial_{\tau_2} g \right) - \frac{ n \cdot \iint \partial_v E } {n \cdot \vb } \left( \vb \cdot \tau_1 \partial_{\tau_1} g + \vb \cdot \tau_2 \partial_{\tau_2 } g \right)
	\\ & - \frac{\tb n(\xb)}{n (\xb) \cdot \vb } ( E \cdot \nabla_v g ) - \nabla_v g \cdot \int \partial_v E + \nabla_v g
	\end{split} \]
	
	where
	\[ \begin{split}
	& \iint \partial_x E = \int _{ t -\tb}^t \int_s^t \partial_x E(X(\tau ) ) d\tau ds, \, \iint \partial_v E =  \int _{ t -\tb}^t \int_s^t \partial_v E(X(\tau ) ) d\tau ds
	\\ & \int \partial_x E = \int_{t-\tb}^t \partial_x E(X(s))ds,  \, \int \partial_v E = \int _{ t -\tb}^t \partial_v E(X(s))ds 
	\end{split} \]
	
	Plug into the previous equation we eventually have:
	
	\[ \begin{split}
	\nabla_x & f(t,x,v ) \textbf{1}_{ \{ t\neq \tb \}} 
	\\ = & \textbf{1}_{\{ t < \tb \}} e^{-\int_0^t \nu } \left[  \nabla_x f_0 \cdot \partial_x X + \nabla_v f_0 \cdot \partial_x V - \int_0^t ( \nabla_x \nu \cdot \partial_x X + \nabla_v \nu \cdot \partial_x V) \times f_0  \right] (X(0), V(0)) 
	\\ & + \textbf{1}_{ t > \tb } e^{- \int_0 ^{\tb } } \Bigg \{ \sum_{i =1 }^2  \tau_i \partial_{\tau_i } g  - \nabla_v g \cdot \int \partial_x E-  \int_0^{\tb} (\nabla_x \nu \cdot \partial_x X + \nabla_v \nu \cdot \partial_x V) \times g 
	\\ & \qquad \qquad \qquad \qquad- \frac{ n }{ n \cdot \vb } \Big \{ \partial_t g + \sum_{i =1}^2 ( \vb \cdot \tau_i ) \partial_{\tau_i } g  + \nu g - H +  E \cdot \nabla_v g \Big \}
	\\  & \qquad\qquad\qquad\qquad- \frac{ n \cdot \iint \partial_x E }{ n \cdot \vb } \Big \{ \partial_t g + \sum_{i =1}^2 ( \vb \cdot \tau_i) \partial_{ \tau_i } g  - \nu g + H \Big \}
	\Bigg \} ( t -\tb, \xb, \vb )
	\\ & + \int_0^{ \min \{ t, \tb \}} e^{- \int_0 ^s \nu }   \bigg [  \nabla_x H \cdot \partial_x X  + \partial_v H \cdot \partial_x V
	\\  & \qquad \qquad \qquad \qquad -\int_0^s ( \nabla_x \nu \cdot \partial_x X +  \nabla_v \nu \cdot \partial_x V) \times H \bigg ]   (t-s, X(t-s ) , V(t-s ) ) ds.
	\end{split} \]
	
	\[ \begin{split}
	\nabla_v & f(t,x,v ) \textbf{1}_{\{ t\neq \tb \}}
	\\ = & \textbf{1}_{\{ t < \tb \} } e^{-\int_0^t \nu } \left[ \nabla_x f_0 \cdot \partial_v X + \nabla_v f_0 \cdot \partial_v V - \int_0^t (\nabla_x \nu \cdot \partial_v X + \nabla_v \nu \cdot \partial_v V) \times f_0 \right] (X(0), V(0))
	\\ & - \textbf{1}_{t > \tb } \tb e^{-\int _0 ^ {\tb} \nu } \Bigg \{ \sum_{i =1}^2 \tau_i \partial_{\tau_i } g 
	- \frac{ n}{ n \cdot \vb } \Big \{ \partial_t g + \sum_{i=1}^2 ( \vb \cdot \tau_i ) \partial_{\tau_i } g + \nu g - H   + E \cdot \nabla_v g \Big \}  \Bigg \}
	\\ & + \textbf{1}_{t > \tb} e^{-\int_0^{\tb}  \nu } \Bigg \{ \nabla_v g - \int_0 ^ {\tb}  ( \nabla_x \nu \cdot \partial_v X + \nabla_v \nu \cdot \partial_v V ) \times g - \nabla_v g \cdot \int \partial_v E \Bigg \}
	\\ & - \textbf{1}_{t > \tb } e^{-\int_0^{\tb} \nu } \Bigg \{ - \frac{ n \cdot \iint \partial_v E }{ n \cdot \vb } \Big \{ \partial_t g + \sum_{i =1}^2 ( \vb \cdot \tau_i) \partial_{ \tau_i } g  - \nu g + H \Big \} \Bigg \} ( t -\tb, \xb, \vb )
	\\ & + \int_0^{\min \{ t, \tb \}} e^{-\int_0^s \nu} \bigg [ \nabla_x H \cdot \partial_v X + \partial_v H \cdot \partial_v V
	\\& \qquad \qquad \qquad \qquad - \int_0^s (\nabla_x \nu \cdot \partial_v X + \nabla_v \nu \cdot \partial_v V) \times H \bigg ] (t-s , X(t-s),V(t -s)) ds
	\end{split} \]
	
	Therefore by the change of variable in (\ref{cov2}) we have:
	\[
	\| f(t) \textbf{1}_{ \{t \neq \tb \} } \| _p
	\lesssim \| f_0 \|_p + \left [ \int_0^t \int_{\gamma_- } | g|^p d\gamma ds \right ] ^{1/p} + \left [ \int_0^t \| H \|_p^p ds \right ]^{1/p},
	\]
	
	\[ \begin{split}
	\| \partial_t & f(t) \textbf{1}_{ \{ t \neq \tb \} } \|_p
	\\ & \lesssim_t \| v \cdot \nabla_x f_0 + \nu f_0 + E \cdot \nabla_x f_0 - H(0,\cdot,\cdot ) \|_p
	\\ & + \left [ \int_0^t \int_{\gamma_- } | \partial_t g | ^p d \gamma ds \right]^{1/p} + \left[ \int_0^t \| \partial_t H \|_p^p \right] ^{1/p}
	\\ & + \{ \| e^ {\theta |v|^2 } f_0 \|_\infty + \| e^{\theta |v|^2 } H\|_\infty + | e^{\theta |v| ^2} g |_\infty \} \left[ \int_0^t \| e^{-\theta |v|^2 } \partial_t \nu \|_p^p \right]^{1/p},
	\end{split} \]
	
	And since we cut off $ \{ (x,v) \in \Omega \times \mathbb R^3: n(\xb ) \cdot \vb < \epsilon \}$, we have $\tb $ is always bounded therefore $|\partial_x X| + | \partial_v X | + |\partial_x V | + |\partial_v V | $ is bounded. Thus by the change of variable (\ref{cov2}) we have
	
	\[ \begin{split}
	\| \nabla_x & f(t) \textbf{1}_{ \{ t \neq \tb \} } \|_p
	\\ & \lesssim_t \| \nabla_x f_0 \|_p + \| \nabla_v f_0 \|_p + \left[ \int_0^t \| \nabla_x H \|_p^p +  \| \nabla_v H \|_p^p \right]^{1/p } + \left[ \int_0^t \int_{\gamma _- } | \nabla_v g |^p d\gamma ds \right] ^{1/p}
	\\ & + \{ \| e^ {\theta |v|^2 } f_0 \|_\infty + \| e^{\theta |v|^2 } H\|_\infty + | e^{\theta |v| ^2} g |_\infty \} \left[ \int_0^t \| e^{-\theta |v|^2 } \partial_t \nu \|_p^p  + \| e^{-\theta |v|^2 } \nabla_v \nu \|_p^p \right]^{1/p}
	\\ & + \Bigg [ \int_0^t \int_{\gamma _- } d \gamma ds | \{ \sum_{i =1 }^2  \tau_i \partial_{\tau_i } g - \frac{ n }{ n \cdot \vb } \Big \{ \partial_t g + \sum_{i =1}^2 ( \vb \cdot \tau_i ) \partial_{\tau_i } g  + \nu g - H +  E \cdot \nabla_v g \Big \}
	\\  & \qquad\qquad\qquad\qquad- \frac{ n \cdot \iint \partial_x E }{ n \cdot \vb } \Big \{ \partial_t g + \sum_{i =1}^2 ( \vb \cdot \tau_i) \partial_{ \tau_i } g  - \nu g + H \Big \} | ^p \Bigg ] ^{1/p}.
	\end{split} \]
	
	and
	
	\[ \begin{split}
	\| \nabla_v & f(t) \textbf{1}_{ \{ t \neq \tb \} } \|_p
	\\ & \lesssim_t \| \nabla_x f_0 \|_p + \| \nabla_v f_0 \|_p + \left[ \int_0^t \| \nabla_x H \|_p^p +  \| \nabla_v H \|_p^p \right]^{1/p } + \left[ \int_0^t \int_{\gamma _- } | \nabla_v g |^p d\gamma ds \right] ^{1/p}
	\\ & + \{ \| e^ {\theta |v|^2 } f_0 \|_\infty + \| e^{\theta |v|^2 } H\|_\infty + | e^{\theta |v| ^2} g |_\infty \} \left[ \int_0^t \| e^{-\theta |v|^2 } \partial_t \nu \|_p^p  + \| e^{-\theta |v|^2 } \nabla_v \nu \|_p^p \right]^{1/p}
	\\ & + \Bigg [ \int_0^t \int_{\gamma _- } d \gamma ds | \{ \sum_{i =1 }^2  \tau_i \partial_{\tau_i } g - \frac{ n }{ n \cdot \vb } \Big \{ \partial_t g + \sum_{i =1}^2 ( \vb \cdot \tau_i ) \partial_{\tau_i } g  + \nu g - H +  E \cdot \nabla_v g \Big \}
	\\  & \qquad\qquad\qquad\qquad- \frac{ n \cdot \iint \partial_v E }{ n \cdot \vb } \Big \{ \partial_t g + \sum_{i =1}^2 ( \vb \cdot \tau_i) \partial_{ \tau_i } g  - \nu g + H \Big \} | ^p \Bigg ] ^{1/p}.
	\end{split} \]\unhide
	From our hypothesis (\ref{assumption_inflow}), these terms are bounded and therefore,
	\[
	\nabla _x f \textbf{1}_{ \{ t \neq \tb \} }, \nabla_v f \textbf{1}_{ \{ t \neq \tb \} }   \in L^\infty ([0,T]; L^p(\Omega \times \mathbb R^3 )).
	\]
	
	Now we prove that $\nabla_{x,v} f \mathbf{1}_{t \neq \tb}$ is weak derivatives of $f$. We claim the following: Let $\Phi (t,x,v) \in C_c ^\infty ((0,T) \times {\Omega} \times \mathbb R^3 )$ then we have
	\Be\label{weak_derivative}
	\int_0^T \iint_{\Omega \times \mathbb R^3 } f \nabla_{x,v} \Phi = - \int_0 ^T \iint_{ \Omega \times \mathbb R^3 } \nabla_{x,v} f \textbf{1}_{ \{ t \neq \tb \} } \Phi.
	\Ee

	Fix the test function $\Phi (t,x,v ) $. From Lemma \ref{cannot_graze}, we have (\ref{no_graze}) for all $(t,x,v) \in \R+_ \times  \O \times \R^3$. Since $\text{supp} (\Phi)$ is compact and $n(\xb(t,x,v))\cdot \vb(t,x,v)$ is continuous if (\ref{no_graze}) we have
	\Be
	\sup_{(t,x,v) \in \text{supp} (\Phi)}|n(\xb(t,x,v)) \cdot \vb(t,x,v)| > \delta_{\Phi} >0.\label{nv_lower_bound_Phi}
	\Ee
	Therefore, from (\ref{computation_tb_x}), $\tb(t,x,v)$ is differentiable in $\text{supp} (\Phi)$ and hence 
	\Be
	\mathcal{M} : = \{ (\tb(t,x,v), x,v) : (t,x,v) \in\text{supp} (\Phi)\}
	\Ee
	is a $C^1$-manifold in $\R_+ \times \O \times \R^3$ with the normal direction 
	\Be\notag
	n_{\mathcal{M}}: = \frac{1}{\sqrt{1+ |\nabla_x \tb|^2 + |\nabla_v \tb|^2}} (1, - \nabla_x \tb, - \nabla_v \tb) \in \R^7.
	\Ee
	
	Now we take $C^1$-approximation $(f_0^l, H^l, g^l)$ of $(f_0, H, g)$ such that, as $l \uparrow \infty$ 
	\Be\label{approximation_f_g_H}
	\begin{split}
		&\| f^l_0 - f_0 \|_{W^{1,p}(\O \times B_{R}(0))} + \| g^l - g \|_{W^{1,p}([0,T] \times \{ \gamma_{-} \cap \ \p\O\times B_{R}(0) \} )}  \\
		&\quad + \| H^l- H \|_{W^{1,p}([0,T] \times \O \times B_{R}(0))} \rightarrow 0,
	\end{split}
	\Ee
	where $B_{R}(0)$ includes projection (onto velocity phase) of $\text{supp}(\Phi)$.
	Then by the trace theorem, for $0< \delta \ll 1$ 
	\[
	f_0^l (x,v) \rightarrow f_0 (x,v) \ \text{and} \ 
	g^l (0,x,v) \rightarrow g(0,x,v) \ \ \text{in} \ L^{1}(\gamma_- \backslash \gamma_-^\delta).
	\]
	
	We define 
	\Be\begin{split}\label{f^l}
		f^l (t,x,v) = & \ \mathbf{1}_{t< \tb} e^{- \int^t_0 \psi} f^l_0 (X(0), V(0)) \\
		& + \mathbf{1}_{t> \tb} e^{- \int^t_{t-\tb} \psi} g^l (t-\tb, \xb, \vb)\\
		&+ \int^t_{\max\{t-\tb, 0\}} e^{- \int^t_s \psi} H^l (s, X(s),V(s)) \dd s,
	\end{split}\Ee
	and $f^l_{+} (t,x,v) : = \mathbf{1}_{t>\tb} f^l(t,x,v) $ and $f^l_{-} (t,x,v) : = \mathbf{1}_{t<\tb} f^l(t,x,v) $. We can derive the same estimate as (\ref{uniform_nabla_f_l}). This implies that up to subsequence 
	\Be\label{weak_convergent_g_l}
	\nabla f^l_{\pm}  \rightharpoonup \nabla f_{\pm} \ \ \text{weakly in } L^p.
	\Ee
	
	By the Gauss theorem, for a standard basis $e_i \in \R_t \times \R_x^3 \times \R^3_v$ with $i=2,3,\cdots, 7$,
	\Be\begin{split}\label{difference_+_-}
		&\iiint e_i \cdot \nabla_{x,v} \Phi  f^l \dd x \dd v  \dd t  \\
		= &   \iint  [f^l_+ - f^l_-] \Phi  \mathbf{e} \cdot n_\mathcal{M}  \dd \mathcal{M}\\
		-&\Big\{ \iiint_{t>\tb} \Phi   e_i \cdot \nabla_{x,v} f^l_+  +  \iiint_{t<\tb} \Phi   e_i \cdot \nabla_{x,v} f^l_-
		\Big\}.
	\end{split} \Ee
	
	From (\ref{f^l}), 
	\Be\begin{split}\notag
		\lim_{t \downarrow \tb}f_+^l (t,x,v) 
		=    e^{-\int^t_0 \psi} g^l (0, \xb, \vb) + \int^t_0 e^{- \int^t_s \psi} H^l (s,X(s),V(s)) \dd s .
	\end{split}\Ee
	Note that as $t \uparrow \tb$, $(X(0;t,x,v), V(0;t,x,v))\rightarrow (\xb, \vb) \in \gamma_-$. Therefore 
	\Be\begin{split}\notag
		\lim_{t \uparrow \tb}f_-^l (t,x,v) 
		=     e^{-\int^t_0 \psi} f^l_0 (\xb, \vb) + \int^t_0 e^{- \int^t_s \psi} H^l (s,X(s),V(s)) \dd s.
	\end{split}\Ee
	Finally (\ref{compatibility_inflow}), (\ref{approximation_f_g_H}), and above two limits imply that, as $l \rightarrow \infty$, 
	\Be
	\begin{split}\notag
		&\|f^l_+ (t,x,v) - f^l_- (t,x,v)\|_{L^1}\\
		= & \ \| e^{- \int^t_{t-\tb} \psi} g^l (t-\tb, \xb,\vb) - e^{- \int^t_0 \psi} f_0^l (X(0), V(0))\|_{L^1} + O(t-\tb) \\
		\lesssim & \ \|  e^{-\int^t_0 \psi}[g(t-\tb, \xb, \vb) - f_0  (\xb,\vb)] \|_{L^1}+ O(t-\tb) \\
		\rightarrow & \ 0. 
	\end{split}
	\Ee 
	Therefore the first term of RHS in (\ref{difference_+_-}) converges to zero as $l \rightarrow \infty$. From (\ref{weak_convergent_g_l}) the second term of RHS in (\ref{difference_+_-}) converges to the RHS of (\ref{weak_derivative}). This proves (\ref{weak_derivative}).
	\hide
	There exists $\delta = \delta_{\phi } > 0$ such that $\phi \equiv 0$ for $t \ge \frac{1}{\delta}$, or dist$(x,\partial \Omega) < \delta$, or $|v| \ge \frac{1}{\delta}$. Let $\phi (t,x,v) \neq 0 $ and $(t,x,v) \in \mathcal M$, so $t = \tb(x,v)$. We have:
	\[
	\delta \le \text{dist}(x,\Omega) \le | x - \xb | = \int_0 ^t | V(s) | ds
	\]
	
	(It seems that I don't have a lower bound for $|V(s)|$ here as the external field $|E|$ could be relatively big comparing 
	to $\delta$, so it's possible for $|V(s)|$ to be arbitrarily small even if $t$ is bounded from above...)\unhide
	%
\end{proof}

\hide

\begin{proposition}
	Lef $f$ be a solution of $\{ \partial_t + v \cdot \nabla_x + E \cdot \nabla_v + \nu \} f = H$. Assume $f_0 (x,v) = g(0,x,v)$ for all $(x,v) \in \gamma_- $.   Define $D^\epsilon : = \{ (x,v) \in \Omega \times \mathbb R^3 : (\xb, \vb) \in \gamma_- \setminus \gamma_-^\epsilon \} .$
	
	For any fixed $p \in [ 2, \infty]$, $0 < \theta 1/4$ and $\beta > 0$, assume
	\[ \begin{split}
	\alpha^\beta \nabla_x f_0, \alpha^\beta \nabla_v f_0, \alpha^\beta[-v \cdot \nabla_x f_0 - \nu ( 0,\cdot,\cdot )f_0  - E \cdot \nabla_v f_0+ H(0,\cdot,\cdot ) ] & \in L^p (\Omega \times \mathbb R^3 ),
	\\ e^{ - \overline \omega \int_0^t (|V|+1) } \alpha^\beta \partial_t g, e^{ -\overline \omega \int_0^t (|V|+1) } \alpha^\beta \nabla_v g, e^{-\overline \omega \int_0^t (|V|+1) } \alpha^\beta \partial_{\tau_i } g &  \in L^p ( [0,T] \times \gamma_- ),
	\\ \frac{ e^{ -\overline \omega \int_0 ^t (|V|+1) } \alpha ^\beta } { n \cdot \vb } \Big \{ \partial_t g + \sum ( v \cdot \tau_i ) \partial_{\tau_i } g + \nu g - H + E\cdot \nabla_v g \Big \}  
	\\   +  \frac{ e^{-\overline \omega \int_0 ^t (|V|+1) } \alpha ^ \beta n \cdot \iint \partial_x E  }{ n \cdot \vb } \Big \{ \partial_t g + \sum ( \vb \cdot \tau_i ) \partial_ {\tau_i } g - \nu g + H \Big \} & \in L^p ([0, T] \times \gamma_- ),
	\\ e^{ - \overline \omega \int_0^t (|V|+1) } \alpha^\beta \partial_t H, e^{ -\overline \omega \int_0^t (|V|+1) } \alpha^\beta \nabla_v H, e^{-\overline \omega \int_0^t (|V|+1) } \alpha^\beta \nabla_x H &  \in L^p ( [0,T] \times \Omega \times \mathbb R^3 ),
	\\  e^{ -\theta |v| ^2 } e^{ - \overline \omega \int_0^t (|V|+1) } \alpha^\beta \partial_t \nu , e^{-\theta |v| ^2 }  e^{ -\overline \omega \int_0^t (|V|+1) } \alpha^\beta \nabla_v \nu, e^{-\theta |v|^2} e^{-\overline \omega \int_0^t (|V|+1) } \alpha^\beta \nabla_x \nu &  \in L^p ( [0,T] \times \Omega \times \mathbb R^3 ),
	\\ e^{\theta |v|^2} f_0 \in L^\infty( \Omega \times \mathbb R^3 ),  e^{\theta |v|^2 } g \in L^\infty ( [0, T ] \times \gamma_-), e^{\theta |v| ^2 } H  & \in L^\infty ([ 0, T] \times \Omega \times \mathbb R^3 ).
	\end{split} \]

\end{proposition}

\begin{proof}
	
	First we assume $f_0$, $g$ and $H$ have compact supports in $\{ v \in \mathbb R^3: |v| < m \}$.
	
	By the previous lemma we have for $\overline \omega\gg 1$ depending only on $E$ and $\xi$, that for any $0 \le s_1, s_2 \le t $ and any $(x,v) \in \Omega \times B(0,m)$ with characteristic flow $X(s),V(s)$ such that $X(t) = x, V(t) =v$:
	
	\[
	e^{ - \overline \omega \int _{s_1} ^{s_2 }( |V(\tau )  | +1) d \tau} \alpha (X(s_1),V(s_1) ) \le \alpha (X(s_2),V(s_2) ) \le e^{ \overline\omega  \int_{s_1}^{s_2} ( |V(\tau)| +1)d\tau} \alpha (X(s_1),V(s_1)).
	\]
	
	Thus we have:
	
	\begin{equation} \label{exp}
	\sup_{ t \le \tb } \frac{ e^{- \overline \omega \int_0^t (|V(\tau) |+1)d \tau } \alpha^\beta (x,v ) }{\alpha ^\beta (X(0),V(0)) } \le e^{ C_{m, \beta } t }, \,
	\sup_{ t \ge \tb } \frac{ e^{- \overline \omega \int_0^t (|V(\tau)| +1)d \tau } \alpha^\beta (x,v ) }{e^{- \overline \omega \int_0 ^{t -\tb}  (| V(\tau ) | +1)d \tau} \alpha ^\beta (\xb,\vb) } \le e^{ C_{m, \beta } \tb } ,\end{equation}
	\[
	\sup_{ \max \{ t-\tb, 0 \} \le s \le t } \frac{ e^{- \overline \omega \int_0^t (|V(\tau)|+1) d \tau } \alpha^\beta (x,v ) }{e^{- \overline \omega \int_0 ^{t -s} ( | V(\tau ) |+1) d \tau} \alpha ^\beta (X(t-s),V(t-s)) } \le e^{ C_{m, \beta } s },
	\]

	Multiplying $e^{ - \overline \omega \int_0^t  (|V(\tau ) | +1) d \tau } \alpha ^\beta (x,v)$ to the previous computations and then using the change of variables from (\ref{det}) and the bound from (\ref{exp}), we get
	\[ \begin{split}
	\| e^{ -\overline \omega \int_0 ^t (|V|+1) }  \alpha ^\beta \partial_t f(t) \|_{L^p (D^\epsilon)}  \lesssim_{t ,m, \beta} & \| \alpha^\beta [ v \cdot \nabla_x f_0 + \nu f_0 + E \cdot \nabla_v f_0 - H(0, \cdot, \cdot ) ] \|_p
	\\ &+ \left [ \int_0^t | e^{- \overline \omega  \int_0 ^s (|V|+1) } \alpha ^\beta \partial_t g(s) | _{\gamma, p } ^ p + \| e^{ -\overline \omega \int_0 ^ s ( | V|+1) }   \alpha ^ \beta \partial_t H(s) \|_p ^ p  \right ] ^{1/p}
	\\ & + C' \left[ \int_0 ^t \| e^{- \theta |v | ^2 } e^{ -\overline \omega  \int_ 0 ^ s (|V| +1)} \alpha ^ \beta \partial_t \nu \|_p ^ p \right ] ^{1/p },
	\end{split} \]
	
	\[ \begin{split}
	\| & e^{ - \overline \omega \int_0 ^t( |V|+1)} \alpha ^\beta \nabla_x f(t) \|_{L^p (D^\epsilon)} \lesssim_{ t,m,\beta } \| \alpha ^\beta \nabla_x f_0 \|_p + \| \alpha^\beta \nabla_v f_0 \|_p
	\\ & + \Bigg [ \int_0 ^t \bigg | \frac{ e^{ -\overline \omega \int_0 ^s (|V|+1) } \alpha ^\beta } { n \cdot \vb } \Big \{ \partial_t g + \sum ( v \cdot \tau_i ) \partial_{\tau_i } g + \nu g - H + E\cdot \nabla_v g \Big \} 
	\\ & \qquad \qquad \qquad \qquad  +  \frac{ e^{-\overline \omega \int_0 ^s (|V|+1) } \alpha ^ \beta n \cdot \iint \partial_x E  }{ n \cdot \vb } \Big \{ \partial_t g + \sum ( \vb \cdot \tau_i ) \partial_ {\tau_i } g - \nu g + H \Big \} \bigg | _{\gamma , p } ^ p  ds \Bigg ]^{1/p }
	\\ & + \Bigg [ \int_0 ^t \sum_{ i =1}^2 \Big | e^{- \overline \omega \int_0 ^s (|V|+1) } \alpha ^\beta \partial_{\tau_i } g (s) \Big | _{\gamma, p } ^p + \Big | e^{- \overline \omega \int_0 ^s (|V|+1) } \alpha ^ \beta \nabla_v g (s)\Big | _{ \gamma, p } ^p 
	\\ & \qquad \qquad \qquad \qquad + \| e^{ -\overline \omega \int_0 ^s (|V|+1) } \alpha ^ \beta \nabla_x H(s) \|_p ^ p + \| e^{ -\overline \omega \int_0 ^s (|V|+1) } \alpha ^ \beta \nabla_v H(s) \|_p ^p ds \Bigg ] ^{ 1/ p}
	\\ & + C' \Bigg [ \int_0 ^t \| e^{ -\theta |v| ^ 2} e^{ - \overline \omega \int_0 ^s (|V|+1) } \alpha ^\beta \nabla_x v \|_p ^ p +  \| e^{ -\theta |v| ^ 2} e^{ - \overline \omega \int_0 ^s (|V|+1) } \alpha ^\beta \nabla_v v \|_p ^ p ds \Bigg ] ^{1/p },
	\end{split} \]
	
	\[ \begin{split}
	\| & e^{ - \overline \omega \int_0 ^t (|V|+1)} \alpha ^\beta \nabla_v f(t) \|_{ L^p (D^\epsilon ) } \lesssim_{ t,m,\beta } \| \alpha ^\beta \nabla_x f_0 \|_p + \| \alpha^\beta \nabla_v f_0 \|_p
	\\ & + \Bigg [ \int_0 ^t \bigg | \frac{ e^{ -\overline \omega \int_0 ^s (|V|+1) } \alpha ^\beta } { n \cdot \vb } \Big \{ \partial_t g + \sum ( v \cdot \tau_i ) \partial_{\tau_i } g + \nu g - H + E\cdot \nabla_v g \Big \} 
	\\ & \qquad \qquad \qquad \qquad  +  \frac{ e^{-\overline \omega \int_0 ^s (|V|+1) } \alpha ^ \beta n \cdot \iint \partial_v E  }{ n \cdot \vb } \Big \{ \partial_t g + \sum ( \vb \cdot \tau_i ) \partial_ {\tau_i } g - \nu g + H \Big \} \bigg | _{\gamma , p } ^ p  ds \Bigg ]^{1/p }
	\\ & + \Bigg [ \int_0 ^t \sum_{ i =1}^2 \Big | e^{- \overline \omega \int_0 ^s (|V|+1) } \alpha ^\beta \partial_{\tau_i } g (s) \Big | _{\gamma, p } ^p + \Big | e^{- \overline \omega \int_0 ^s (|V|+1) } \alpha ^ \beta \nabla_v g (s)\Big | _{ \gamma, p } ^p 
	\\ & \qquad \qquad \qquad \qquad + \| e^{ -\overline \omega \int_0 ^s (|V|+1) } \alpha ^ \beta \nabla_x H(s) \|_p ^ p + \| e^{ -\overline \omega \int_0 ^s (|V|+1) } \alpha ^ \beta \nabla_v H(s) \|_p ^p ds \Bigg ] ^{ 1/ p}
	\\ & + C' \Bigg [ \int_0 ^ t \| e^{ -\theta |v| ^ 2} e^{ - \overline \omega \int_0 ^s (|V|+1) } \alpha ^\beta \nabla_x v \|_p ^ p +  \| e^{ -\theta |v| ^ 2} e^{ - \overline \omega \int_0 ^s (|V|+1) } \alpha ^\beta \nabla_v v \|_p ^ p ds \Bigg ] ^{1/p },
	\end{split} \]
	
\end{proof}

\hide

\begin{lemma} \label{lem COV}
	Let us use $\bar{x}$ to denote horizontal local coordinate for parametrization $\eta$, i.e. $\eta(\bar{x}) \in \p\O$. For fixed $s>0$ and any $(\eta(\bar{x}), v ) \in \gamma_+$ such that  $t-\tb (t,\eta(\bar{x}), v ) < s < t$, a map
	\Be \label{map+}
	\mathcal{M}_{+} : ( t, \bar{x} , v ) \mapsto (X(s;t, \eta(\bar{x}), v ) , V(s; t, \eta(\bar{x}), v ) ) \in \Omega\times\R^{3}
	\Ee
	is well-defined. If we denote Jacobian of $\mathcal{M}_{+}$ as $J_{\mathcal{M}_{+}}$, then
	\Be
	\det J_{\mathcal{M}_{+}} = - v \cdot (\p_1 \eta(\bar{x}) \times \p_2 \eta(\bar{x}) ).
	\Ee
	Similary, for fixed $s>0$ and any $(\eta(\bar{x}), v) \in \gamma_{-}$ such that  $t < s < t + \tf (t,\eta(\bar{x}), v ) $, the following map
	\Be \label{map-}
	\mathcal{M}_{-} : ( t, \bar{x} , v ) \mapsto ( X(s;t, \eta(\bar{x}), v ) , V(s; t, \eta(\bar{x}), v ) ) \in \Omega\times\R^{3}
	\Ee
	is well-defined and also
	\Be
	\det J_{\mathcal{M}_{-}} = - v \cdot (\p_1 \eta(\bar{x}) \times \p_2 \eta(\bar{x}) ),
	\Ee
	holds.
\end{lemma}
\begin{proof}
	We compute $\mathcal{M}_{+}$ case first:

	\Be
	\begin{split}\label{XV_txv}
		&\frac{\p ( X( s;t, \eta(\bar{x}), v),  V( s;t, \eta(\bar{x}), v) )}{\p (t, \bar{x}, v)}\\
		= & \ \begin{bmatrix}
			\p_t X( s;t, \eta(\bar{x}), v) & \nabla_{\bar{x}} X( s;t, \eta(\bar{x}), v) & \nabla_v X( s;t, \eta(\bar{x}), v)\\
			\p_t V( s;t, \eta(\bar{x}), v)& \nabla_{\bar{x}} V( s;t, \eta(\bar{x}), v) & \nabla_v V( s;t, \eta(\bar{x}), v)
		\end{bmatrix}\\
		= & \ 
		\begin{bmatrix}
			\p_t X( s;t, \eta(\bar{x}), v)& \nabla_{\bar{x}} \eta (\bar{x}) \cdot \nabla_x X( s;t, \eta(\bar{x}), v)  & \nabla_v X( s;t, \eta(\bar{x}), v)\\
			\p_t V( s;t, \eta(\bar{x}), v) &  \nabla_{\bar{x}} \eta (\bar{x}) \cdot \nabla_x V( s;t, \eta(\bar{x}), v)  & \nabla_v V( s;t, \eta(\bar{x}), v)
		\end{bmatrix}
	\end{split}
	\Ee
	

	From the deterministic property of trajectory phase yields the following identity,
	\Be
	\begin{split}
		X(s; t+ \delta, X(t+ \delta;t, \eta(\bar{x}),v), V(t+ \delta;t, \eta(\bar{x}),v))& = X(s;t, \eta(\bar{x}), v),  \\
		V(s; t+ \delta, X(t+ \delta;t, \eta(\bar{x}),v), V(t+ \delta;t, \eta(\bar{x}),v))& = V(s;t, \eta(\bar{x}), v),  \\
	\end{split}
	\Ee
	for $s-t \leq \delta \leq 0$. Therefore 
	\Be \begin{split}\notag
		[\p_t + v\cdot \nabla_x -\nabla_x \phi  (t, \eta(\bar{x})) \cdot \nabla_v]X(s;t, \eta(\bar{x}), v) & = 0,\\
		[\p_t + v\cdot \nabla_x -\nabla_x \phi  (t, \eta(\bar{x})) \cdot \nabla_v]V(s;t, \eta(\bar{x}), v) & = 0.
	\end{split}
	\Ee
	Equivalently     
	\Be\label{pt_XV}
	\begin{split}
		& \begin{bmatrix}
			\p_t X( s;t, \eta(\bar{x}), v)\\
			\p_t V( s;t, \eta(\bar{x}), v)
		\end{bmatrix}\\
		&=
		\begin{bmatrix}
			\nabla_x X( s;t, \eta(\bar{x}), v)  & \nabla_v X( s;t, \eta(\bar{x}), v)\\
			\nabla_x V( s;t, \eta(\bar{x}), v)  & \nabla_v V(s;t, \eta(\bar{x}), v)
		\end{bmatrix}
		\begin{bmatrix}
			-v  \\
			\nabla \phi (t, \eta(\bar{x})) 
		\end{bmatrix}
	\end{split}
	\Ee
	From (\ref{XV_txv}) and (\ref{pt_XV}) we conclude that 
	\Be
	\begin{split}\label{XV_txv_final}
		&\frac{\p ( X( s;t, \eta(\bar{x}), v),  V( s;t, \eta(\bar{x}), v) )}{\p (t, \bar{x}, v)}\\
		= &\begin{bmatrix}
			\nabla_x X( s;t, \eta(\bar{x}), v)  & \nabla_v X( s;t, \eta(\bar{x}), v)\\
			\nabla_x V( s;t, \eta(\bar{x}), v)  & \nabla_v V(s;t, \eta(\bar{x}), v)
		\end{bmatrix}  
		\begin{bmatrix}
			-v & 
			\p_{\bar{x}} \eta   & 0_{3 \times 3} \\
			\nabla \phi(t, \eta(\bar{x}))
			& 
			0_{3 \times 2} & \mathrm{Id}_{3\times 3}
		\end{bmatrix}.  
	\end{split}
	\Ee
	Since 
	$$\det \begin{bmatrix}
	\nabla_x X( s;t, \eta(\bar{x}), v)  & \nabla_v X( s;t, \eta(\bar{x}), v)\\
	\nabla_x V( s;t, \eta(\bar{x}), v)  & \nabla_v V(s;t, \eta(\bar{x}), v)
	\end{bmatrix} =1,
	$$
	we conclude that
	\Be \label{cov2}
	\begin{split} 
		&\det\left(\frac{\p ( X( s;t, \eta(\bar{x}), v),  V( s;t, \eta(\bar{x}), v) )}{\p (t, \bar{x}, v)}\right) \\
		= &  \ \det \begin{bmatrix}
			-v & 
			\p_{\bar{x}} \eta   & 0_{3 \times 3} \\
			\nabla \phi(t, \eta(\bar{x}))
			& 
			0_{3 \times 2} & \mathrm{Id}_{3\times 3}
		\end{bmatrix}\\
		=& \ - v \cdot (\p_1 \eta(\bar{x}) \times \p_2 \eta(\bar{x}) ).
	\end{split}\Ee
	
	For $\mathcal{M}_{-}$ case, we note that (\ref{identi}) holds for $0 \leq \delta \leq s-t$. Hence, the computation for Jacobian of $\mathcal{M}_{-}$ is exactly same as $\mathcal{M}_{+}$ case. \\
\end{proof}

\begin{lemma}
	Suppose $h(t,x,v) \in L^1 ( [0,T] \times \Omega \times \mathbb R^3) $ then:
	\[
	\begin{split}
	\int_0 ^T & \iint _{\Omega \times \mathbb R^3} h (t,x,v) dv dx dt \\ = & \iint_{\Omega \times \mathbb R^3 } \int_{- \min ( T,\tb(T,x,v) )} ^0 h (T +s, X(T+s; T,x,v), V(T+s;T,x,v) ) ds dv dx \\ & + \int_0 ^T \int_{\gamma^+ } \int_{-\min ( t, \tb(t,x,v)) } ^0 h(t +s , X(t+s;t,x,v), V(t +s;t,x,v) )ds d\gamma dt.
	\end{split}
	\]
\end{lemma}

\begin{proof}
	The region $\{ (t,x,v) \in [0,T] \times \Omega \times \mathbb R^3 \}$ is the disjoint union of 
	\[
	A:= \{ (t,x,v) \in [0,T] \times \Omega \times \mathbb R^3: \tf (t,x,v) + t\le T \} 
	\] and 
	\[ B:= \{ (t,x,v) \in [0,T] \times \Omega \times \mathbb R^3: \tf (t,x,v)+ t > T \} .\]
	
	Now let:
	\[
	A': =  \{ (t,s,x,v) \in [0,T]^2 \times \gamma^+ : s < \tb(t,x,v), s \le t \},
	\] and
	\[
	B':= \{ (s,x,v ) \in [0,T] \times \Omega \times \mathbb R^3 : s < \tb(T,x,v) \}
	\] 
	
	Consider the map $\mathcal A: A' \to A$ with
	\[
	\mathcal A(t,s,x,v) = (t-s, X(t-s; t,x,v ) , V(t-s;t,x,v) ).
	\]
	Since $\tf(t-s, X(t-s ; t,x,v) , V(t-s;t,x,v) ) + (t-s) = s + (t-s) = t \le T $, $\mathcal A$ is well-defined. And since the characteristic flow is deterministic, $\alpha$ is injective. And for any $(t,x,v) \in A$, we have
	\[
	(t + \tf(t,x,v), \tf(t,x,v), X(t + \tf(t,x,v);t,x,v), V(t + \tf(t,x,v);t,x,v) ) \in A'
	\] 
	since $\tf \le t + \tf$ and $\tb ( t + \tf, X(t + \tf;t,x,v ) , V(t + \tf;t,x,v ) ) > \tf $ as $x \in \Omega$ is in the interior.
	Moreover 
	\[
	\mathcal A (t + \tf, \tf, X(t + \tf;t,x,v), V(t + \tf;t,x,v) ) = (t,x,v) ,
	\] 
	so $\mathcal A$ is surjective. Therefore $\mathcal A$ is bijective with inverse $\mathcal A^{-1} (t,x,v ) =  (t + \tf, \tf, X(t + \tf;t,x,v), V(t + \tf;t,x,v) )$ \\
	
	Suppose locally at $x \in \gamma^+$ we have $x = \eta(\bar{x} ) $, and let 
	\[
	J_{\mathcal A} = \frac{ (t-s, X(t-s; t,x,v ) , V(t-s;t,x,v) )}{\partial(t,s,\bar x, v )}
	\]
	be the Jacobian matrix of $\mathcal A$. Let us compute the determinant of $J_{\mathcal A }$. Here we use the notation $ \partial_t X(t-s;t,x,v ) :=  \frac{ \partial X(s;t,x,v ) }{\partial t } |_{(t-s;t,x,v )} $, and $\partial_s X(t-s;t,x,v ) :=  \frac{ \partial X(s;t,x,v ) }{\partial s } |_{(t-s;t,x,v )} $ to denote the partial derivatives of $X$ for the first and second coordinates, and similar for $V$. Not to confuse with the derivative with respect to $s$ or $t$ which we denote $\frac{d X(t-s; t,x,v ) }{ ds }$ and $\frac{d X(t-s; t,x,v ) }{ dt }$. So we actually have: $\frac{d X(t-s; t,x,v ) }{ ds } = - \partial_s X(t-s; t,x,v ) $, and $\frac{d X(t-s; t,x,v ) }{ dt } = \partial_s X(t-s; t,x,v ) +  \partial_t X(t-s; t,x,v )$.
	Now we have:
	\tiny
	\[ \begin{split}
	J_{\mathcal A } = 
	\begin{bmatrix}
	1 & -1 & 0_{1\times 2 } &  0_{1 \times 3 } \\
	\partial_s X(t-s; t,x,v ) +  \partial_t X(t-s; t,x,v ) & -\partial_s X(t-s; t,x,v ) & \partial_{\bar x } X(t-s; t,x,v ) & \partial_v X(t-s; t,x,v ) \\
	\partial_s V(t-s; t,x,v ) +  \partial_t V(t-s; t,x,v ) & -\partial_s V(t-s; t,x,v ) & \partial_{\bar x } V(t-s; t,x,v ) & \partial_v V(t-s; t,x,v ) \\
	\end{bmatrix}
	\end{split} \]
	\normalsize
	
	From elementary column operation,
	\tiny
	\begin{equation} \label{77}
	\begin{split}
	\det & (J_{\mathcal A }) = \det(J_{\mathcal A }') \\
	= & \det \left \{ \begin{bmatrix}
	1 & 0 & 0_{1\times 2 } &  0_{1 \times 3 } \\
	\partial_s X(t-s; t,x,v ) +  \partial_t X(t-s; t,x,v ) & \partial_t X(t-s; t,x,v ) & \partial_{\bar x } X(t-s; t,x,v ) & \partial_v X(t-s; t,x,v ) \\
	\partial_s V(t-s; t,x,v ) +  \partial_t V(t-s; t,x,v ) & \partial_t V(t-s; t,x,v ) & \partial_{\bar x } V(t-s; t,x,v ) & \partial_v V(t-s; t,x,v ) \\
	\end{bmatrix} \right \} \\
	= & \det \left \{ \begin{bmatrix}
	\partial_t X(t-s; t,x,v ) & \partial_{\bar x } X(t-s; t,x,v ) & \partial_v X(t-s; t,x,v ) \\
	\partial_t V(t-s; t,x,v ) & \partial_{\bar x } V(t-s; t,x,v ) & \partial_v V(t-s; t,x,v ) \\
	\end{bmatrix} \right \} .
	\end{split} \end{equation}
	\normalsize
	
	Therefore, from Lemma \ref{lem COV},
	\[
	\begin{split}
	\det  (J_{\mathcal A })  &=  \ - v \cdot (\p_1 \eta(\bar{x}) \times \p_2 \eta(\bar{x}) ).  \\
	\end{split}
	\]
	Hence,
	\small
	\Be \label{intA}\begin{split}
		& \iiint_A  h(t,x,v ) dt dxdv 
		\\ = & \int_0^T \int_{\gamma^+ } \int_0 ^{\min( \tb(t, x,v) , t ) } h(t-s, X(t-s;t,x,v ) ,V(t-s; t,x,v )  \\
		&\quad\quad \times \big| v \cdot  \frac{ \p_1 \eta(\bar{x}) \times \p_2 \eta(\bar{x}) }{   | \p_1 \eta(\bar{x}) \times \p_2 \eta(\bar{x}) |} \big| | (\p_1 \eta(\bar{x}) \times \p_2 \eta(\bar{x}) )| ds d \bar x dt
		\\ = & \int_0^T \int_{\gamma^+ } \int_0 ^{\min( \tb(t,x,v) , t ) } h(t-s, X(t-s;t,x,v ) ,V(t-s; t,x,v )    ds d \gamma dt
	\end{split} \Ee  \\
	\normalsize
	
	Now consider the map $\mathcal B: B' \to B$ with
	\[
	\mathcal B(s,x,v) = (T -s , X(T-s, T,x,v ), V(T-s, T,x,v ) ).
	\]
	Since $ \tf(T-s, X(T-s, T,x,v ),V(T-s, T,x,v ) ) + (T-s) > s + (T-s) = T$, $\mathcal B$ is well-defined. And since the characteristic flow is deterministic, $\beta$ is injective. And for any $(t,x,v) \in B$, we have
	\[
	(T-t, X(T; t, x,v ) , V(T;t,x,v) ) \in B'
	\]
	since $\tb(T, X(T; t,x,v ) , V(T; t,x,v ) ) > T -t$ as $x \in \Omega$ is in the interior. Moreover 
	\[
	\mathcal B (T-t, X(T; t, x,v ) , V(T;t,x,v) ) = (t,x,v),
	\]
	so $\mathcal B$ is surjective. Therefore $\mathcal B$ is bijective with inverse $\mathcal B ^{-1} (t,x,v) =  (T-t, X(T; t, x,v ) , V(T;t,x,v) )$
	
	And since $\mathcal{B}$ is a measure preserving mapping, we obtain
	\Be \label{intB}
	\begin{split}
		&\iiint_B h(t,x,v ) dtdxdv   \\
		&\quad =  \iint_{\Omega \times \mathbb R^3 } \int_0^{ \min (T,\tb(T,x,v) ) }  h(T-s, X(T-s;T,x,v ) ,V(T-s; T,x,v ) ds dx dv.
	\end{split}
	\Ee
	We combining (\ref{intA}) and (\ref{intB}), and apply change of variable $s\rightarrow -s$ to finish the proof.
\end{proof}

\hide
\begin{lemma}
	Suppose $h(t,x,v) \in L^1 ( [0,T] \times \Omega \times \mathbb R^3 $ then:
	\[
	\begin{split}
	\int_0 ^T & \iint _{\Omega \times \mathbb R^3} h (t,x,v) dv dx dt \\ = & \iint_{\Omega \times \mathbb R^3 } \int_{- \min ( T,\tb(T,x,v) )} ^0 h (T +s, X(T+s; T,x,v), V(T+s;T,x,v) ) ds dv dx \\ & + \int_0 ^T \int_{\gamma^+ } \int_{-\min ( t, \tb(t,x,v)) } ^0 h(t +s , X(t+s;t,x,v), V(t +s;t,x,v) )ds d\gamma dt.
	\end{split}
	\]
\end{lemma}

\begin{proof}
	The region $\{ (t,x,v) \in [0,T] \times \Omega \times \mathbb R^3 \}$ is the disjoint union of 
	\[
	A:= \{ (t,x,v) \in [0,T] \times \Omega \times \mathbb R^3: \tf (t,x,v) + t\le T \} 
	\] and 
	\[ B:= \{ (t,x,v) \in [0,T] \times \Omega \times \mathbb R^3: \tf (t,x,v)+ t > T \} .\]
	
	Now let:
	\[
	A': =  \{ (t,s,x,v) \in [0,T]^2 \times \gamma^+ : s < \tb(t,x,v), s \le t \},
	\] and
	\[
	B':= \{ (s,x,v ) \in [0,T] \times \Omega \times \mathbb R^3 : s < \tb(T,x,v) \}
	\]
	
	Consider the map $\mathcal A: A' \to A$ with
	\[
	\mathcal A(t,s,x,v) = (t-s, X(t-s; t,x,v ) , V(t-s;t,x,v) ).
	\]
	Since $\tf(t-s, X(t-s ; t,x,v) , V(t-s;t,x,v) ) + (t-s) = s + (t-s) = t \le T $, $\mathcal A$ is well-defined. And since the characteristic flow is deterministic, $\alpha$ is injective. And for any $(t,x,v) \in A$, we have
	\[
	(t + \tf(t,x,v), \tf(t,x,v), X(t + \tf(t,x,v);t,x,v), V(t + \tf(t,x,v);t,x,v) ) \in A'
	\] 
	since $\tf \le t + \tf$ and $\tb ( t + \tf, X(t + \tf;t,x,v ) , V(t + \tf;t,x,v ) ) > \tf $ as $x \in \Omega$ is in the interior.
	Moreover 
	\[
	\mathcal A (t + \tf, \tf, X(t + \tf;t,x,v), V(t + \tf;t,x,v) ) = (t,x,v) ,
	\] 
	so $\mathcal A$ is surjective. Therefore $\mathcal A$ is bijective with inverse $\mathcal A^{-1} (t,x,v ) =  (t + \tf, \tf, X(t + \tf;t,x,v), V(t + \tf;t,x,v) )$ 
	
	Suppose locally at $x \in \gamma^+$ we have $x = \eta(\bar x ) $, and let $J_{\mathcal A} = \frac{ (t-s, X(t-s; t,x,v ) , V(t-s;t,x,v) )}{\partial(t,s,\bar x, v )}$ be the Jacobian matrix of $\mathcal A$. Let's compute the determinant of $J_{\mathcal A }$. Here let's use the notation $ \partial_t X(t-s;t,x,v ) :=  \frac{ \partial X(s;t,x,v ) }{\partial t } |_{(t-s;t,x,v )} $, and $\partial_s X(t-s;t,x,v ) :=  \frac{ \partial X(s;t,x,v ) }{\partial s } |_{(t-s;t,x,v )} $ to denote the partial derivatives of $X$ for the first and second coordinates, and similar for $V$. Not to confuse with the derivative with respect to $s$ or $t$ which we denote $\frac{d X(t-s; t,x,v ) }{ ds }$ and $\frac{d X(t-s; t,x,v ) }{ dt }$. So we actually have: $\frac{d X(t-s; t,x,v ) }{ ds } = - \partial_s X(t-s; t,x,v ) $, and $\frac{d X(t-s; t,x,v ) }{ dt } = \partial_s X(t-s; t,x,v ) +  \partial_t X(t-s; t,x,v )$.
	Now we have:
	\small
	\[ \begin{split}
	J_{\mathcal A } = 
	\begin{bmatrix}
	1 & -1 & 0_{1\times 2 } &  0_{1 \times 3 } \\
	\partial_s X(t-s; t,x,v ) +  \partial_t X(t-s; t,x,v ) & -\partial_s X(t-s; t,x,v ) & \partial_{\bar x } X(t-s; t,x,v ) & \partial_v X(t-s; t,x,v ) \\
	\partial_s V(t-s; t,x,v ) +  \partial_t V(t-s; t,x,v ) & -\partial_s V(t-s; t,x,v ) & \partial_{\bar x } V(t-s; t,x,v ) & \partial_v V(t-s; t,x,v ) \\
	\end{bmatrix}
	\end{split} \]
	\normalsize
	
	Let $J_{\mathcal A} '$ be the matrix obtained by ddding the first column of $J_{\mathcal A}$ to its second column, so
	\small
	\begin{equation} \label{77}
	\begin{split}
	\det & (J_{\mathcal A }) = \det(J_{\mathcal A }') \\
	= & \det \left \{ \begin{bmatrix}
	1 & 0 & 0_{1\times 2 } &  0_{1 \times 3 } \\
	\partial_s X(t-s; t,x,v ) +  \partial_t X(t-s; t,x,v ) & \partial_t X(t-s; t,x,v ) & \partial_{\bar x } X(t-s; t,x,v ) & \partial_v X(t-s; t,x,v ) \\
	\partial_s V(t-s; t,x,v ) +  \partial_t V(t-s; t,x,v ) & \partial_t V(t-s; t,x,v ) & \partial_{\bar x } V(t-s; t,x,v ) & \partial_v V(t-s; t,x,v ) \\
	\end{bmatrix} \right \} \\
	= & \det \left \{ \begin{bmatrix}
	\partial_t X(t-s; t,x,v ) & \partial_{\bar x } X(t-s; t,x,v ) & \partial_v X(t-s; t,x,v ) \\
	\partial_t V(t-s; t,x,v ) & \partial_{\bar x } V(t-s; t,x,v ) & \partial_v V(t-s; t,x,v ) \\
	\end{bmatrix} \right \} .
	\end{split} \end{equation}
	\normalsize

	To calculate this determinant let's now consider the change of variable: for fixed $s>0$, for any $(\eta(x_1,x_2), v ) \in \gamma_+$ such that  $\tb (\eta(x_1, x_2), v ) > s$, we have $( t, x_1,x_2 , v ) \in [0,s) \times \gamma_+ \mapsto (X(s;t, \eta(x_1,x_2), v ) , V(s; t, \eta(x_1,x_2), v ) ) \in \Omega$. We can compute the determinant of this change of variable:

	\Be
	\begin{split}\label{XV_txv}
		&\frac{\p ( X( s;t, \eta(\bar{x}), v),  V( s;t, \eta(\bar{x}), v) )}{\p (t, \bar{x}, v)}\\
		= & \ \begin{bmatrix}
			\p_t X( s;t, \eta(\bar{x}), v) & \nabla_{\bar{x}} X( s;t, \eta(\bar{x}), v) & \nabla_v X( s;t, \eta(\bar{x}), v)\\
			\p_t V( s;t, \eta(\bar{x}), v)& \nabla_{\bar{x}} V( s;t, \eta(\bar{x}), v) & \nabla_v V( s;t, \eta(\bar{x}), v)
		\end{bmatrix}\\
		= & \ 
		\begin{bmatrix}
			\p_t X( s;t, \eta(\bar{x}), v)& \nabla_{\bar{x}} \eta (\bar{x}) \cdot \nabla_x X( s;t, \eta(\bar{x}), v)  & \nabla_v X( s;t, \eta(\bar{x}), v)\\
			\p_t V( s;t, \eta(\bar{x}), v) &  \nabla_{\bar{x}} \eta (\bar{x}) \cdot \nabla_x V( s;t, \eta(\bar{x}), v)  & \nabla_v V( s;t, \eta(\bar{x}), v)
		\end{bmatrix}
	\end{split}
	\Ee
	

	Note that 
	\Be\notag
	\begin{split}
		X(s; t+ \Delta, X(t+ \Delta;t, \eta(\bar{x}),v), V(t+ \Delta;t, \eta(\bar{x}),v))& = X(s;t, \eta(\bar{x}), v),\\
		V(s; t+ \Delta, X(t+ \Delta;t, \eta(\bar{x}),v), V(t+ \Delta;t, \eta(\bar{x}),v))& = V(s;t, \eta(\bar{x}), v).\\
	\end{split}
	\Ee
	Therefore 
	\Be \begin{split}\notag
		[\p_t + v\cdot \nabla_x -\nabla_x \phi  (t, \eta(\bar{x})) \cdot \nabla_v]X(s;t, \eta(\bar{x}), v) & = 0,\\
		[\p_t + v\cdot \nabla_x -\nabla_x \phi  (t, \eta(\bar{x})) \cdot \nabla_v]V(s;t, \eta(\bar{x}), v) & = 0.
	\end{split}
	\Ee
	Equivalently     
	\Be\label{pt_XV}
	\begin{split}
		& \begin{bmatrix}
			\p_t X( s;t, \eta(\bar{x}), v)\\
			\p_t V( s;t, \eta(\bar{x}), v)
		\end{bmatrix}\\
		&=
		\begin{bmatrix}
			\nabla_x X( s;t, \eta(\bar{x}), v)  & \nabla_v X( s;t, \eta(\bar{x}), v)\\
			\nabla_x V( s;t, \eta(\bar{x}), v)  & \nabla_v V(s;t, \eta(\bar{x}), v)
		\end{bmatrix}
		\begin{bmatrix}
			-v  \\
			\nabla \phi (t, \eta(\bar{x})) 
		\end{bmatrix}
	\end{split}
	\Ee
	From (\ref{XV_txv}) and (\ref{pt_XV}) we conclude that 
	\Be
	\begin{split}\label{XV_txv_final}
		&\frac{\p ( X( s;t, \eta(\bar{x}), v),  V( s;t, \eta(\bar{x}), v) )}{\p (t, \bar{x}, v)}\\
		= &\begin{bmatrix}
			\nabla_x X( s;t, \eta(\bar{x}), v)  & \nabla_v X( s;t, \eta(\bar{x}), v)\\
			\nabla_x V( s;t, \eta(\bar{x}), v)  & \nabla_v V(s;t, \eta(\bar{x}), v)
		\end{bmatrix}  
		\begin{bmatrix}
			-v & 
			\p_{\bar{x}} \eta   & 0_{3 \times 3} \\
			\nabla \phi(t, \eta(\bar{x}))
			& 
			0_{3 \times 2} & \mathrm{Id}_{3\times 3}
		\end{bmatrix}.  
	\end{split}
	\Ee
	Since $$\det \begin{bmatrix}
	\nabla_x X( s;t, \eta(\bar{x}), v)  & \nabla_v X( s;t, \eta(\bar{x}), v)\\
	\nabla_x V( s;t, \eta(\bar{x}), v)  & \nabla_v V(s;t, \eta(\bar{x}), v)
	\end{bmatrix} =1,$$
	we conclude that
	\Be \label{cov2}
	\begin{split} 
		&\det\left(\frac{\p ( X( s;t, \eta(\bar{x}), v),  V( s;t, \eta(\bar{x}), v) )}{\p (t, \bar{x}, v)}\right) \\
		= &  \ \det \begin{bmatrix}
			-v & 
			\p_{\bar{x}} \eta   & 0_{3 \times 3} \\
			\nabla \phi(t, \eta(\bar{x}))
			& 
			0_{3 \times 2} & \mathrm{Id}_{3\times 3}
		\end{bmatrix}\\
		=& \ - v \cdot (\p_1 \eta(\bar{x}) \times \p_2 \eta(\bar{x}) ).
	\end{split}\Ee
	
	Therefore from (\ref{77}) and (\ref{XV_txv_final}) we have
	\[
	\begin{split}
	\det  (J_{\mathcal A })  = & \det(J_{\mathcal A }') \\
	=&  \det \left \{ \begin{bmatrix}
	\partial_t X(t-s; t,x,v ) & \partial_{\bar x } X(t-s; t,x,v ) & \partial_v X(t-s; t,x,v ) \\
	\partial_t V(t-s; t,x,v ) & \partial_{\bar x } V(t-s; t,x,v ) & \partial_v V(t-s; t,x,v ) \\
	\end{bmatrix} \right \} \\
	= &\det \left \{ \begin{bmatrix}
	\nabla_x X( t-s;t, \eta(\bar{x}), v)  & \nabla_v X( t-s;t, \eta(\bar{x}), v)\\
	\nabla_x V( t-s;t, \eta(\bar{x}), v)  & \nabla_v V(t-s;t, \eta(\bar{x}), v)
	\end{bmatrix}  
	\begin{bmatrix}
	-v & 
	\p_{\bar{x}} \eta   & 0_{3 \times 3} \\
	\nabla \phi(t, \eta(\bar{x}))
	& 
	0_{3 \times 2} & \mathrm{Id}_{3\times 3}
	\end{bmatrix} \right \} \\
	= &  \ \det \begin{bmatrix}
	-v & 
	\p_{\bar{x}} \eta   & 0_{3 \times 3} \\
	\nabla \phi(t, \eta(\bar{x}))
	& 
	0_{3 \times 2} & \mathrm{Id}_{3\times 3}
	\end{bmatrix}\\
	=& \ - v \cdot (\p_1 \eta(\bar{x}) \times \p_2 \eta(\bar{x}) ).
	\end{split}
	\]

	Therefore
	\small
	\[ \begin{split}
	& \iiint_A  h(t,x,v ) dt dxdv 
	\\ = & \int_0^T \int_{\gamma^+ } \int_0 ^{\min( \tb(t, x,v) , t ) } h(t-s, X(t-s;t,x,v ) ,V(t-s; t,x,v )  ( v \cdot  \frac{ \p_1 \eta(\bar{x}) \times \p_2 \eta(\bar{x}) }{   | \p_1 \eta(\bar{x}) \times \p_2 \eta(\bar{x}) |}) | (\p_1 \eta(\bar{x}) \times \p_2 \eta(\bar{x}) )| ds d \bar x dt
	\\ = & \int_0^T \int_{\gamma^+ } \int_0 ^{\min( \tb(t,x,v) , t ) } h(t-s, X(t-s;t,x,v ) ,V(t-s; t,x,v )    ds d \gamma dt
	\end{split} \]
	\normalsize
	Now consider the map $\mathcal B: B' \to B$ with
	\[
	\mathcal B(s,x,v) = (T -s , X(T-s, T,x,v ), V(T-s, T,x,v ) ).
	\]
	Since $ \tf(T-s, X(T-s, T,x,v ),V(T-s, T,x,v ) ) + (T-s) > s + (T-s) = T$, $\mathcal B$ is well-defined. And since the characteristic flow is deterministic, $\beta$ is injective. And for any $(t,x,v) \in B$, we have
	\[
	(T-t, X(T; t, x,v ) , V(T;t,x,v) ) \in B'
	\]
	since $\tb(T, X(T; t,x,v ) , V(T; t,x,v ) ) > T -t$ as $x \in \Omega$ is in the interior. Moreover 
	\[
	\mathcal B (T-t, X(T; t, x,v ) , V(T;t,x,v) ) = (t,x,v),
	\]
	so $\mathcal B$ is surjective. Therefore $\mathcal B$ is bijective with inverse $\mathcal B ^{-1} (t,x,v) =  (T-t, X(T; t, x,v ) , V(T;t,x,v) )$
	
	And since $\beta$ is a measure preserving change of variable we have:
	\[
	\iiint_B h(t,x,v ) dtdxdv =  \iint_{\Omega \times \mathbb R^3 } \int_0^{ \min (T,\tb(T,x,v) ) }  h(T-s, X(T-s;T,x,v ) ,V(T-s; T,x,v ) ds dx dv
	\]
	
	Thus:
	\[
	\begin{split}
	\int_0 ^T  \iint _{\Omega \times \mathbb R^3} h (t,x,v) dv dx dt  & = \iiint_A h(t,x,v ) dt dxdv + \iiint_B h(t,x,v ) dt dxdv 
	\\ = & \iint_{\Omega \times \mathbb R^3 } \int_{- \min ( T,\tb(T, x,v) )} ^0 h (T +s, X(T+s; T,x,v), V(T+s;T,x,v) ) ds dv dx 
	\\ & + \int_0 ^T \int_{\gamma^+ } \int_{-\min ( t, \tb(t,x,v)) } ^0 h(t +s , X(t+s;t,x,v), V(t +s;t,x,v) )ds d\gamma dt.
	\end{split}
	\]
	as wanted.
\end{proof}\unhide

The first result is an ``energy estimate" to the transport operator. 

\begin{lemma} [Green's identity] \label{lem_Green}
	For $p \in [1, \infty)$ assume $f$, $\partial_t f + v \cdot \nabla_x f + E \cdot \nabla_v f \in L^p ([0,T]; L^p (\Omega \times \mathbb R^3 ) )$ and $f_{\gamma_- } \in L^p ( [0,T]; L^p (\gamma ) )$. Then $f \in C^0( [0, T] ; L^p (\Omega \times \mathbb R^3 ) $ and $f_{\gamma_+ } \in L^p ([0,T]; L^p (\gamma) )$ and for almost every $T' \in [0, T]$:
	\[ \begin{split}
	\| f(T')\|_p ^p + \int_0^{T'} |f|_{ \gamma_+, p } ^p = \| f(0) \|_p^p + \int_0^{T'} |f|_{ \gamma_-, p }^p + \int_0^{T'} \iint _{\Omega \times \mathbb R^3 }  p \{ \partial_t + v \cdot \nabla_x f + E \cdot \nabla_v f \} |f|^{p-2} f.
	\end{split} \]
	
\end{lemma}

\begin{proof}
	First let's assume $f \in C^1$.
	
	For almost every $T' \in [0,T]$, By Holder's inequality we have $ \| (\partial_t f + v \cdot \nabla_x f + E \cdot \nabla_v f) |f|^{p-2} f  \|_{L^1 ( [0,T] \times \Omega \times \mathbb R^3 )}  \le \|(\partial_t f + v \cdot \nabla_x f + E \cdot \nabla_v f) \|_{L^p ([0,T] \times \Omega \times \mathbb R^3 ) } \| |f|^{p-1}\|_{L^{p/(p-1)} ([0,T] \times \Omega \times \mathbb R^3 )} < \infty$. Thus by Lemma (\ref{int_id}) we have:
	\small
	\[ \begin{split}
	& \int_0^{T'}  \int_{\Omega \times \mathbb R^3 }  p(\partial_t f + v \cdot \nabla_x f + E \cdot \nabla_v f) |f|^{p-2} f  dxdvdt  
	\\ = & \iint_{\Omega \times \mathbb R^3 } \int_{- \min ( T',\tb(T',x,v) )} ^0 p(\partial_t f + v \cdot \nabla_x f + E \cdot \nabla_v f) |f|^{p-2} f (T' +s, X(T'+s; T',x,v), V(T'+s;T',x,v) ) ds dv dx 
	\\ & + \int_0 ^{T'} \int_{\gamma^+ } \int_{-\min ( t, \tb(t,x,v)) } ^0 p(\partial_t f + v \cdot \nabla_x f + E \cdot \nabla_v f) |f|^{p-2} f(t +s , X(t+s;t,x,v), V(t +s;t,x,v) )ds d\gamma dt.
	\end{split} \]
	\normalsize
	
	Since
	\[ \begin{split}
	\frac{d}{ds } |f|^p & ( T'+s, X(T'+s; T',x,v ), V(T'+s; T',x,v )) 
	\\ =  & p(\partial_t f + v \cdot \nabla_x f + E \cdot \nabla_v f) |f|^{p-2} f (T' +s, X(T'+s; T',x,v), V(T'+s;T',x,v) ),
	\end{split} \]
	
	and
	\[ \begin{split}
	\frac{d}{ds } |f|^p & ( t+s, X(t+s; t,x,v ), V(t+s; t,x,v )) 
	\\ =  & p(\partial_t f + v \cdot \nabla_x f + E \cdot \nabla_v f) |f|^{p-2} f (t +s, X(t+s; t,x,v), V(t+s;t,x,v) ).
	\end{split} \]
	
	We have
	\begin{equation} \label{green} \begin{split}
	\int_0^{T'} &  \int_{\Omega \times \mathbb R^3 }  p(\partial_t f + v \cdot \nabla_x f + E \cdot \nabla_v f) |f|^{p-2} f dxdvdt  
	\\ = & \iint_{\Omega \times \mathbb R^3 } \int_{- \min ( T',\tb(T',x,v) )} ^0 \frac{d}{ds } |f|^p ( T'+s, X(T'+s; T',x,v ), V(T'+s; T',x,v ))ds dv dx 
	\\ & + \int_0 ^{T'} \int_{\gamma^+ } \int_{-\min ( t, \tb(t,x,v)) } ^0 \frac{d}{ds } |f|^p ( t+s, X(t+s; t,x,v ), V(t+s; t,x,v )) ds d\gamma dt.
	\\ = & \iint_{\Omega \times \mathbb R^3 } |f|^p (T',x,v ) dxdv - \iint_{\Omega \times \mathbb R^3 } \textbf{1}_{\{ T' \ge \tb(T',x,v) \}} |f|^p (T'-\tb,\xb,\vb ) dx dv 
	\\&  - \iint_{\Omega \times \mathbb R^3 } \textbf{1}_{\{ T' < \tb(T',x,v) \}} |f|^p (0,X(0;T',x,v) ,V(0;T',x,v ) ) dx dv
	\\ & + \int_0^{T'} \int_{\gamma ^+ } |f|^p (t,x,v ) d\gamma dt - \int_0^{T'} \int_{\gamma^+ } \textbf{1}_{\{ t \ge \tb (t,x,v ) \} }|f|^p (t -\tb, \xb,\vb ) d\gamma dt
	\\& - \int_0^{T'} \int_{\gamma^+ }  \textbf{1}_{\{ t < \tb (t,x,v ) \} }|f|^p (0, X(0;t,x,v ),V(0;t,x,v) ) d\gamma dt
	\end{split} \end{equation}

	Now consider the map $ \mathcal A_1 : (x,v) \mapsto (X(0; T',x,v), V(0; T',x,v) ) $ from $\{ (x,v) \in \Omega \times \mathbb R^3 : T' < \tb(T',x,v) \}$ to $\{ (x,v) \in \Omega \times \mathbb R^3 : \tf(0,x,v) >  T' \}$. This map is well defined as $\tf(0, X(0; T',x,v), V(0; T',x,v)) > T'$ since $x \in \Omega$ is in the interior. $\mathcal A_1$ is injective as the characteristic flow is unique. And for any $ (x,v) \in \Omega \times \mathbb R^3$ such that $\tf(0,x,v) > T'$, we have $X(T';0,x,v ) \in \Omega$ and $\mathcal A_1 ( X(T';0,x,v ), V(T';0,x,v ) ) = (x,v )$, so $\mathcal A_1$ is surjective. Therefore $\mathcal A_1$ is a bijection.
	And since the trajectory of this change of variable is measure preserving, we have
	\begin{equation} \label{green1}
	\iint_{\Omega \times \mathbb R^3 } |f_0|^p \textbf{1}_{\{ \tf(0,x,v) > T' \} } dxdv = \iint_{\Omega \times \mathbb R^3 } |f_0 | ^p (X(0; T',x,v), V(0;T',x,v)) \textbf{1}_{\{ T' < \tb(T',x,v) \}}  dx dv
	\end{equation}
	
	Consider the map $\mathcal A_2: (t,x,v ) \mapsto (X(0;t,x,v), V(0;t,x,v))$ from $ \{ (t,x,v) \in  (0, T'] \times \gamma_+ : t < \tb(t,x,v) \}$ to $\{ (x,v) \in \Omega \times \mathbb R^3 : \tf(0,x,v) \le T' \}$. This map is well defined as $\tf( 0, X(0;t,x,v), V(0;t,x,v) ) = t \le T'$. $\mathcal A_2$ is injective as the characteristic flow is unique. And for any $(x,v) \in \Omega \times \mathbb R^3$ such that $\tf(0,x,v) \le T'$, we have $(\tf, X(\tf; 0,x,v) , V(\tf; 0,x,v) ) \in (0, T'] \times \gamma_+$ and $\tb(\tf, X(\tf; 0,x,v) , V(\tf; 0,x,v)) > \tf$ as $x \in \Omega$ is in the interior; moreover, $\mathcal A_2 ( \tf, X(\tf; 0,x,v) , V(\tf; 0,x,v) ) = (x,v)$, so $\mathcal A_2$ is surjective. Therefore $\mathcal A_2$ is a bijection.

	So by our change of variable computation (\ref{cov2}) we have:
	\begin{equation} \label{green2}
	\iint_{\Omega \times \mathbb R^3 } |f_0 |^p \textbf{1}_{ \{ \tf(0,x,v) \le T' \} }   dx dv = \int_0^{T' } \int_{ \gamma_+ } |f_0|^p (X(0;t,x,v), V(0;t,x,v) ) \textbf{1}_{ \{ t < \tb(t,x,v) \} }     d\gamma d t
	\end{equation}
	
	Therefore we have
	\[ \begin{split}
	\iint_{\Omega \times \mathbb R^3 }  |f_0|^p dx dv =  & \iint_{\Omega \times \mathbb R^3 } |f_0 |^p \textbf{1}_{ \{ \tf (0,x,v) > T' \} } dx dv + \iint_{\Omega \times \mathbb R^3 } |f_0 |^p \textbf{1}_{ \{ \tf(0,x,v) \le T' \} } dx dv
	\\ = &  \iint_{\Omega \times \mathbb R^3 } |f_0 | ^p (X(0;T',x,v), V(0;T',x,v)) \textbf{1}_{\{ T' < \tb(T',x,v) \}}  dx dv 
	\\ & + \int_0^{T' } \int_{ \gamma_+ } |f_0|^p (X(0;t,x,v), V(0;t,x,v) ) \textbf{1}_{ \{ t < \tb(t,x,v) \} }     d\gamma d t
	\end{split} \]
	
	Now consider the map $\mathcal A_3 : (t,x,v ) \mapsto (t -\tb(t,x,v) , \xb, \vb ) $ from $\{ (t,x,v ) \in [0, T' ] \times \gamma_+ : t \ge \tb(t,x,v) \} $ to $\{ (s,x,v) \in [0, T' ) \times \gamma_- : T'  \ge s + \tf(s,x,v) \}$. This map is well defined as $\tf( t -\tb ,\xb ,\vb ) + (t - \tb) = \tb + t - \tb = t \le T'$. $\mathcal A_3$ is injective as the characteristic flow is unique. And for any $(s,x,v ) \in [0, T') \times \gamma_-$ such that $s + \tf(s,x,v) \le T'$, we have $(s + \tf, X(s + \tf; s,x,v ) , V(s + \tf; x,v ) ) \in [ 0, T'] \times \gamma_+$ and $\tb ( s + \tf, X(s + \tf; s,x,v ) , V(s + \tf; s,x,v ) )  = \tf \le s + \tf$; moreover, $\mathcal A_3 (s + \tf, X(s + \tf; s,x,v ) , V(s + \tf; s,x,v )) = (s,x,v ) $, so $\mathcal A_3$ is surjective. Therefore $\mathcal A_3$ is a bijection.
	
	Suppose without lose of generality that locally we have $x = \eta_1 (x_1, x_2)$ and $\xb = \eta_2 ( \xba, \xbb$). Then the change of variable is:
	\[
	(t,x_1,x_2, v ) \mapsto ( t -\tb (t,\eta_1 (x_1,x_2) , v ) , \eta_2 ^{-1} (X(t -\tb; t ,\eta_1(x_1,x_2) ,v)), V(t-\tb; t, \eta_1(x_1,x_2), v) )
	\]

	We compute the Jacobian matrix $J$ of this change of variable:
	\Be
	\begin{split}\label{XVtb_txv}
		J = &\frac{\p (  t -\tb, \eta_2 ^{-1} (X(t -\tb; t ,\eta_1(\bar x ) ,v)), V(t-\tb; t, \eta_1(\bar x), v)  )}{\p (t, \bar x, v)}\\
		= & \ \begin{bmatrix}
			1 - \partial_t \tb & - \partial_{\bar x} \tb  & - \nabla_v \tb \\
			\nabla_x (\eta_2 ^{-1} ) \cdot ( \partial_s X + \partial_t X) & \nabla_x (\eta_2 ^{-1} ) \cdot ( \nabla_x X \cdot \partial_{\bar x} \eta_1 - \partial_s X \cdot \partial_{\bar x } \tb )   & \nabla_x (\eta_2 ^{-1} ) \cdot ( \nabla_v X  - \partial_s X \cdot \nabla_v \tb ) \\
			\partial_s V + \partial_t V & \nabla_x V \cdot \partial_{\bar x } \eta_1 - \partial_s V \cdot \partial_{\bar x } \tb & \nabla_v V - \partial_s V \cdot \nabla_v \tb 
		\end{bmatrix}\\
	\end{split}
	\Ee
	
	Let
	\[ 
	A =
	\begin{bmatrix}
	-\partial_s X & \partial_{\bar x } \eta_2 & 0_{3 \times 3} \\
	-\partial_s V &  0_{3 \times 2} & \mathrm{Id}_{3 \times 3}
	\end{bmatrix}
	\]
	
	Then we have
	\tiny
	\[ \begin{split}
	& A \cdot  J \\  =  & 
	\begin{bmatrix}
	-\partial_s X & \partial_{\bar x } \eta & 0_{3 \times 3} \\
	-\partial_s V &  0_{3 \times 2} & \mathrm{Id}_{3 \times 3}
	\end{bmatrix}
	\cdot
	\begin{bmatrix}
	1 & - \partial_{\bar x} \tb  & - \nabla_v \tb \\
	\nabla_x (\eta_2 ^{-1} ) \cdot ( \partial_s X + \partial_t X) & \nabla_x (\eta_2 ^{-1} ) \cdot ( \nabla_x X \cdot \partial_{\bar x} \eta_1 - \partial_s X \cdot \partial_{\bar x } \tb )   & \nabla_x (\eta_2 ^{-1} ) \cdot ( \nabla_v X  - \partial_s X \cdot \nabla_v \tb ) \\
	\partial_s V + \partial_t V  & \nabla_x V \cdot \partial_{\bar x } \eta_1 - \partial_s V \cdot \partial_{\bar x } \tb & \nabla_v V - \partial_s V \cdot \nabla_v \tb
	\end{bmatrix}
	\\ 
	= & \begin{bmatrix}
	- \partial_s X + \partial_{\bar x } \eta_2 \cdot \nabla_x( \eta_2 ^{-1} ) ( \partial_s X + \partial_t X) &  \partial_s X \cdot \partial_{\bar x } \tb + \partial_{\bar x } \eta_2 \cdot  \nabla_x (\eta_2 ^{-1} ) \cdot ( \nabla_x X \cdot \partial_{\bar x} \eta_1 - \partial_s X \cdot \partial_{\bar x } \tb )   & \partial_s X \cdot \nabla_v \tb + \partial_{\bar x} \eta_2 \cdot \nabla_x (\eta_2 ^{-1} ) \cdot ( \nabla_v X  - \partial_s X \cdot \nabla_v \tb )   \\
	-\partial_s V + \partial_s V + \partial_t V & \partial_s V \cdot \partial_{\bar x} \tb + \nabla_x V \cdot \partial_{\bar x } \eta_1 - \partial_s V \cdot \partial_{\bar x } \tb & \partial_s V \cdot \nabla_v \tb + \nabla_v V - \partial_s V \cdot \nabla_v \tb
	\end{bmatrix}
	\\
	= & \begin{bmatrix}
	- \partial_t X & \nabla_x X \cdot \partial_{\bar x } \eta_1 & \nabla_v X \\
	-\partial_t V & \nabla_x V \cdot \partial_{ \bar x } \eta_1 & \nabla_v V 
	\end{bmatrix}
	\end{split} \]
	\normalsize
	
	Since
	\[
	\partial_{\bar x }  \eta_2 \cdot \nabla_x (\eta_2 ^{-1} ) = \nabla_x ( \eta_2 \circ \eta_2 ^{-1} )  = \mathrm{Id}_{3\times 3 }.
	\]
	
	And by the computation of previous change of variable (\ref{cov2}) we have
	\[ \begin{split}
	\det  (A \cdot  J )   = & \det
	\begin{bmatrix}
	- \partial_s X & \nabla_x X \cdot \partial_{\bar x } \eta_1 & \nabla_v X \\
	-\partial_s V & \nabla_x V \cdot \partial_{ \bar x } \eta_1 & \nabla_v V 
	\end{bmatrix}
	\\ = &  \ \det \begin{bmatrix}
	-v & 
	\p_{\bar{x}} \eta_1   & 0_{3 \times 3} \\
	-E
	& 
	0_{3 \times 2} & \mathrm{Id}_{3\times 3}
	\end{bmatrix}\\
	\\ =&  - v \cdot (\partial_1 \eta_1 (\bar x ) \times \partial_2 \eta_1 (\bar x ) 
	\end{split} \]
	
	Since
	\[ \begin{split}
	& \det  (A  ) 
	\\  = & \det
	\begin{bmatrix}
	-\partial_s X & \partial_{\bar x } \eta & 0_{3 \times 3} \\
	-\partial_s V &  0_{3 \times 2} & \mathrm{Id}_{3 \times 3}
	\end{bmatrix}
	\\ =&  - \vb \cdot (\partial_1 \eta_2 (\overline {\xb} ) \times \partial_2 \eta_2 (\overline {\xb} ) 
	\end{split} \]
	
	Therefore
	\[
	\det(J ) = \frac{v \cdot (\partial_1 \eta_1 (\bar x ) \times \partial_2 \eta_1 (\bar x ) }{\vb \cdot (\partial_1 \eta_2 (\overline {\xb} ) \times \partial_2 \eta_2 (\overline {\xb} ) }
	\]

	So
	\begin{equation} \label{green3}
	\int_0^{T' } \int_{\gamma_- } |f|^p (t,x,v) \textbf{1}_{ \{T'   \ge s +  \tf(s,x,v) \} } d\gamma ds = \int_0^{T' } \int_{\gamma_+ } |f|^p ( t- \tb(t,x,v), \xb,\vb) \textbf{1}_{ \{ t  \ge \tb(t,x,v) \} } d\gamma dt
	\end{equation}
	Finally consider the map $\mathcal A_4: ( x,v ) \mapsto (T' - \tb(T',x,v) ,\xb ,\vb )$ from $\{ (x,v) \in \Omega \times \mathbb R^3 : T' \ge \tb(T',x,v) \} $ to $\{ (s,x,v ) \in [0, T') \times \gamma_- : T'  < s + \tf(s,x,v) \} $. This map is well defined as $\tf(T' - \tb ,\xb ,\vb ) + (T' - \tb )  > \tb +(T' - \tb) = T$ as $x \in \Omega$ is in the interior. $\mathcal A_4$ is injective as the characteristic flow is unique. And for any $\{ (s,x,v ) \in [0, T') \times \gamma_-$ such that $T'  < s + \tf(s,x,v) $, we have $ (X( T'; s, x,v ) , V(T';s,x,v ) ) \in \Omega \times \mathbb R^3$ and $ \tb ( T', X( T'; s, x,v ) , V(T';s,x,v )) = T'-s \le T'$; moreover, $\mathcal A_4 ( X( T'; s, x,v ) , V(T';s,x,v ) ) = (s,x,v )$, so $\mathcal A_4$ is surjective. Therefore $\mathcal A_4$ is a bijection. 
	
	And by the computation in the change of variable (\ref{XVtb_txv}) we have:
	\begin{equation} \label{green4}
	\int_0^{T' } \int_{\gamma_- } |f|^p (t,x,v) \textbf{1}_{ \{T'   < s + \tf(s,x,v) \} }   d\gamma dt =   \iint_{ \Omega \times \mathbb R^3 } |f |^p (T' - \tb(T',x,v) , \xb, \vb ) \textbf{1}_{ \{ T' \ge \tb(T',x,v) \}} dxdv.
	\end{equation}
	
	Now substitute all these identities (\ref{green1}), (\ref{green2}), (\ref{green3}),(\ref{green4}) into equation (\ref{green}) we finally get:
	
	\[ \begin{split}
	\int_0^{T'} &  \int_{\Omega \times \mathbb R^3 }  p(\partial_t f + v \cdot \nabla_x f + E \cdot \nabla_v f) |f|^{p-2} f dxdvdt  
	\\ = & \iint_{\Omega \times \mathbb R^3 } |f|^p (T',x,v ) dxdv - \iint_{\Omega \times \mathbb R^3 } \textbf{1}_{\{ T' \ge \tb(T',x,v) \}} |f|^p (T'-\tb,\xb,\vb ) dx dv 
	\\&  - \iint_{\Omega \times \mathbb R^3 } \textbf{1}_{\{ T' < \tb(T',x,v) \}} |f|^p (0,X(0;T',x,v) ,V(0;T',x,v ) ) dx dv
	\\ & + \int_0^{T'} \int_{\gamma ^+ } |f|^p (t,x,v ) d\gamma dt - \int_0^{T'} \int_{\gamma^+ } \textbf{1}_{\{ t \ge \tb (t,x,v ) \} }|f|^p (t -\tb, \xb,\vb ) d\gamma dt
	\\& - \int_0^{T'} \int_{\gamma^+ }  \textbf{1}_{\{ t < \tb (t,x,v ) \} }|f|^p (0, X(0;t,x,v ),V(0;t,x,v) ) d\gamma dt
	\\ = & \iint_{\Omega \times \mathbb R^3 } |f|^p (T',x,v ) dxdv -\int_0^{T' } \int_{\gamma_- } |f|^p (t,x,v) \textbf{1}_{ \{T'   < s + \tf(s,x,v) \} }   d\gamma dt 
	\\&  - \iint_{\Omega \times \mathbb R^3 } |f_0|^p \textbf{1}_{\{ \tf(0,x,v) > T' \} } dxdv
	\\ & + \int_0^{T'} \int_{\gamma ^+ } |f|^p (t,x,v ) d\gamma dt - \int_0^{T' } \int_{\gamma_- } |f|^p (t,x,v) \textbf{1}_{ \{T'   \ge s +  \tf(s,x,v) \} } d\gamma ds
	\\& - \iint_{\Omega \times \mathbb R^3 } |f_0 |^p \textbf{1}_{ \{ \tf(0,x,v) \le T' \} }   dx dv 
	\\ = & \iint_{\Omega \times \mathbb R^3 } |f|^p (T',x,v ) dxdv + \int_0^{T'} \int_{\gamma ^+ } |f|^p (t,x,v ) d\gamma dt
	\\ & - \iint_{\Omega \times \mathbb R^3 } |f_0 |^p dxdv - \int_0^{T' } \int_{\gamma_- } |f|^p ( t,x,v) d\gamma dt
	\end{split} \]
	
	as wanted.
	
	Note since the left hand side of the above equality is finite, and by our assumption all the terms on the right hand side except $\int_0^{T'} \int_{\gamma _+ } |f|^p (t,x,v ) d\gamma dt$ is finite, thus $f \in L^p ( [0,T] ; L^p (\gamma_+ ) ) $.\end{proof}
\hide

?????? Finally for $f \in L^p ([0,T] ; L^p ( \Omega \times \mathbb R^3 )  ) $ such that $\partial_t f + v \cdot \nabla_x f + E \cdot \nabla_v f  \in L^p ([0,T] ; L^p ( \Omega \times \mathbb R^3 )  ) $. We can find a sequence of $C^1$ functions $\{ f_n \}$ such that $ f_n \to f$ in $L^p ([0,T] ; L^p ( \Omega \times \mathbb R^3 )  ),  $and $\partial_t f_n + v \cdot \nabla_x f_n + E \cdot \nabla_v f_n \to \partial_t f + v \cdot \nabla_x f + E \cdot \nabla_v f$ in $L^p ([0,T] ; L^p ( \Omega \times \mathbb R^3 )  )  $. Then we get the desired equality for $f$.

\unhide

\hide

\begin{lemma}[Green's Identity]
	\label{Green} For $p\in[1,\infty)$ assume that $h,\partial_t h + v\cdot
	\nabla_x h-\nabla\phi \cdot \nabla_v h \in L^p ([0,T];L^p(\Omega\times\mathbb{R}^3))$ and $%
	h_{\gamma_-} \in L^p ([0,T];L^p(\gamma))$. Then $h  \in C^0( [0,T];
	L^p(\Omega\times\mathbb{R}^3))$ and $h_{\gamma_+} \in
	L^p ([0,T];L^p(\gamma))$ and for almost every $t\in [0,T]$ :
	\Be\begin{split}
		&\| h(t) \|_{p}^p + \int_0^t |h |_{\gamma_+,p}^p \\
		= & \  \| h(0)\|_p^p + \int_0^t
		|h |_{\gamma_-,p}^p + \int_0^t \iint_{\Omega\times\mathbb{R}^3 }
		\{\partial_t h + v\cdot \nabla_x h- \nabla \phi \cdot \nabla_v h\} |h|^{p-2}h.
	\end{split}\Ee
\end{lemma}
\begin{proof}Note that 
	\Be\begin{split}\notag
		&\frac{d}{ds} |h(s,X(s; \cdot ,x,v), V(s;\cdot, x,v))|^p \\
		= & \ p[\p_t + v\cdot \nabla_x-\nabla \phi \cdot \nabla_v] h(s,X(s; \cdot ,x,v), V(s;\cdot, x,v))\\
		& \times  h(s,X(s; \cdot ,x,v), V(s;\cdot, x,v)) |  h(s,X(s; \cdot ,x,v), V(s;\cdot, x,v))|^{p-2}.
	\end{split}
	\Ee
	Therefore 
	\Be
	\begin{split}
		&\int^{s_1}_{s_2} [\p_t + v\cdot \nabla_x - \nabla_x \phi \cdot \nabla_v] |h(s,X(s;\cdot, x,v) ,V(s;\cdot, x,v))|^p \dd s\\
		=& \ \int^{s_1}_{s_2} \frac{d}{ds}|h(s,X(s;\cdot, x,v) ,V(s;\cdot, x,v))|^p \dd s\\
		=& \ |h(s_1,X(s_1;\cdot, x,v) ,V(s_1;\cdot, x,v))|^p - |h(s_2,X(s_2;\cdot, x,v) ,V(s_2;\cdot, x,v))|^p.
	\end{split}
	\Ee
	From (\ref{integration}) we conclude the lemma.

\end{proof}

\unhide

The second result if a trace theorem in the presence of an external potential. Note that the result is local in $(x,v)$, but the estimate is uniform in time provided that a potential decays exponentially.  
\begin{lemma}
	\label{le:ukai} Assume we have $\phi(t,x) \in C^1_x$ such that, for $\lambda>0$, $\delta>0$, 
	\begin{equation}
	\label{decay_phi}
	\sup_{t \geq 0} e^{\lambda t} \| \nabla \phi (t) \|_{\infty} \leq \delta.
	\end{equation}
	We also assume $\frac{1}{C}\langle v\rangle \leq \psi(t,x,v)\leq C \langle v\rangle$ for some $C>0$. For any small parameter 
	\Be\label{lower_bound_e}
	\e> \frac{2\delta}{  \lambda}>0,
	\Ee
	there exists a
	constant $C_{\delta, \e ,\Omega }>0$, which does not depend on $t$, such that for any $h$ in $%
	L^{1}_{loc}
	([0,\infty);L^{1}(\Omega \times \mathbb{R}^{3}))$ with $\partial _{t}h+v\cdot
	\nabla _{x}h
	- \nabla \phi \cdot \nabla_v h + \psi h$ in $L^{1}_{loc}( [0, \infty);L^{1}(\Omega \times \mathbb{R}%
	^{3}))$, we have, for all $t\geq 0$,
	\begin{equation*}
	\begin{split}
	&\int_{0}^{t}\int_{\gamma _{+}\setminus \gamma _{+}^{\varepsilon }}|h|\mathrm{%
		d}\gamma \mathrm{d}s\\
	&\leq C_{\delta, \varepsilon ,\Omega }\left\{ \
	||h_{0}||_{1}+\int_{0}^{t}  \| h(s)\|_{1}+\big{\Vert} 
	[
	\partial
	_{t}+v\cdot \nabla _{x}- \nabla_x\phi \cdot \nabla_v + \psi ]h(s)\big{\Vert} _{1} \mathrm{d}s\ \right\}.\end{split}
	\end{equation*}%
	Furthermore, for any $(t,x,v)$ in $[0, \infty)\times \Omega \times \mathbb{R}^{3}$
	the function $
	h
	(t+s^{\prime },
	X(t+s^\prime;t,x,v)
	,
	V(t+s^\prime;t,x,v))
	$ is absolutely continuous in
	$s^{\prime } \in (-\min \{t_{\mathbf{b}}(t,x,v),t\}, t_{%
		\mathbf{f}}(t,x, v) )$.
\end{lemma}
\begin{proof}From (\ref{decay_phi}), for $0 \leq s \leq t$
	\Be
	\begin{split}\notag
		&|V(s;t,x,v)-v|\\
		\leq 
		&  \ 
		\int_s^t |  \nabla\phi (\tau, X(\tau;t,x,v))| \dd \tau  \\
		\leq   &
		\ 
		\frac{\delta}{\lambda} e^{- \lambda s} \{1- e^{-\lambda (t-s)}\} \\
		\leq
		& \ 
		\frac{\delta}{\lambda} .
	\end{split}
	\Ee
	Hence if $|v| \geq \e$ with the condition (\ref{lower_bound_e}) then 
	\Be
	\begin{split}\notag
		&|X(s;t,x,v)-x|\\
		\geq  & \ \left|\int^t_s v\right| -  \int^t_s| V(\tau;t,x,v)-v | \dd \tau \\
		\geq & \  |v| (t-s) -\frac{\delta}{\lambda} (t-s) \\
		\geq & \  \frac{\delta}{\lambda}(t-s) .
	\end{split}
	\Ee
	Therefore 
	\Be\label{upper_tb}
	\tb(t,x,v) \leq \frac{\lambda}{\delta} \times \mathrm{diam} (\O) .
	\Ee
	
	Now let's give an upper bound for $\tb(t,x,v)$ for any $(t,x,v) \in [0, \infty) \times (\gamma_+ \setminus \gamma_+^\epsilon )$:
	
	Note that since for $x_1 \in \partial \Omega$,
	\[
	\lim_{ y \to x_1, y\in \partial \Omega } \frac{ | (x_1 - y ) \cdot n(x_1 ) | }{ |x_1 - y | } = 0.
	\] 
	Hence we have $| (x_1 - y ) \cdot n(x_1) | \le C_{\Omega} |x_1 - y| $ for all $ y \in \partial \Omega$.
	
	Thus
	\[ \begin{split}
	C_{\Omega}|x - \xb |^2 \ge & | (\xb - x ) \cdot n(x) | = | \int_{t - \tb} ^t  V(\tau;t,x,v )  d\tau \cdot n(x) | 
	\\ \ge &  \left| \int_{t -\tb}^t v  d \tau \cdot n(x) \right| - \left| \int_{t-\tb} ^t |v - V(\tau;t,x,v) | d\tau \cdot n(x) \right|
	\\ \ge & \tb | v \cdot n(x) | - \int_{t -\tb}^t |v - V(\tau;t,x,v) | d\tau
	\\ > & \tb \epsilon - \tb \frac{\delta}{\lambda}
	\\ > & \frac{\epsilon}{2} \tb
	\end{split} \]
	
	On the other hand
	\[ \begin{split}
	|x - \xb |^2 & = | \tb v - \int_{t -\tb} ^t ( v - V(s) ) ds | ^2
	\\ & \le 2 \tb ^2 |v|^2 + 2 \left( \int_{t -\tb} ^t |v - V(s) | ds \right) ^2
	\\ & \le 2 \tb^2 |v|^2 + 2 (\tb \frac{ \delta}{\lambda})^2
	\\ & < 2 \tb^2 |v| ^2 + \frac{1}{2} \tb \epsilon ^2.
	\end{split} \]
	Thus
	\[
	C_\Omega (  2 \tb^2 |v| ^2 + \frac{1}{2} \tb^2 \epsilon ^2) > \frac{\epsilon}{2} \tb.
	\]
	Dividing $|v|^2 \tb$ on both sides we get
	\[
	C_\Omega \frac{5}{2} \tb > C_\Omega (  2 \tb  + \frac{1}{2} \tb \frac{\epsilon ^2}{|v|^2}) > \frac{\epsilon}{2 |v|^2} > \frac{\epsilon^3}{2} .
	\]
	Therefore $\tb(t,x,v) > \frac{1}{5 C_\Omega} \epsilon ^3$.
	
	Now by lemma (\ref{int_id}) we have 
	\begin{equation} \label{h_id}
	\begin{split}
	\int_0 ^T & \iint _{\Omega \times \mathbb R^3} h (t,x,v) dv dx dt \\ = & \iint_{\Omega \times \mathbb R^3 } \int_{- \min ( T,\tb(T,x,v) )} ^0 h (T +s, X(T+s; T,x,v), V(T+s;T,x,v) ) ds dv dx \\ & + \int_0 ^T \int_{\gamma^+ } \int_{-\min ( t, \tb(t,x,v)) } ^0 h(t +s , X(t+s;t,x,v), V(t +s;t,x,v) )ds d\gamma dt.
	\end{split}
	\end{equation}

	For $(t,x,v) \in [0,T] \times \gamma_+ $ and $0 \le s \le \min \{ t ,\tb (t,x,v) \} $, we have
	
	\begin{equation} \label{h_expand}
	\begin{split}
	h(t,x,v) =  &h(t-s, X(t-s),V(t-s) ) e^{\int_0^{-s} \phi (V(t + \tau' )) d\tau'} \\ &+  \int_{-s}^0 e^{\int_0^\tau \phi(V(t + \tau' )) d\tau' } \left[ \partial_t h + V(t + \tau ) \cdot \nabla_x h + E(X(t + \tau )) \cdot \nabla_v h + \phi (V(t + \tau)) h \right ] \\ & ( t + \tau , X(t + \tau ) , V(t + \tau )) d \tau.
	\end{split}
	\end{equation}
	Here $X(t + \tau ) = X(t + \tau; t ,x,v ) $, and $V(t + \tau ) = V( t + \tau; t,x,v  )$.
	
	Thus
	\[ \begin{split}
	\min & \{ t, \tb(t,x,v) \}  \times |h (t,x,v)| = \int_{ - \min \{ t, \tb(t,x,v) \} } ^0 |h(t,x,v)| ds
	\\   \le & \int_{ - \min \{ t, \tb(t,x,v) \} } ^0 |h(t + s, X(t+s),V(t+s) ) | ds 
	\\ & +  \int_{ - \min \{ t, \tb(t,x,v) \} } ^0 \int_{s}^0 | \left[ \partial_t h + V(t + \tau ) \cdot \nabla_x h + E(X(t + \tau )) \cdot \nabla_v h + \phi (V(t + \tau)) h \right ] 
	\\ & ( t + \tau , X(t + \tau ) , V(t + \tau )) | d \tau ds
	\\ \le & \int_{ - \min \{ t, \tb(t,x,v) \} } ^0 |h(t+s, X(t+s),V(t+s) ) | ds 
	\\ & +  \tb(t,x,v) \times \int_{ - \min \{ t, \tb(t,x,v) \} } ^0  | \left[ \partial_t h + V(t + \tau ) \cdot \nabla_x h + E(X(t + \tau )) \cdot \nabla_v h + \phi (V(t + \tau)) h \right ] 
	\\ & ( t + \tau , X(t + \tau ) , V(t + \tau )) | d \tau 
	\end{split} \]

	Now let $\epsilon_1 = \frac{1}{5 C_\Omega} \epsilon ^3 $, and for $(t,x,v) \in [\epsilon_1 ,T] \times \gamma_+ \setminus \gamma_+^\epsilon $, we integrate over $\int_{\epsilon_1} ^T \int_{\gamma_+ \setminus \gamma_+^\epsilon} \int_{-\min \{ t, \tb (t,x,v) \} } ^0 $ to get 
	\[
	\begin{split}
	& \epsilon_1  \int_{\epsilon_1} ^T \int_{\gamma_+ \setminus \gamma_+^\epsilon} | h(t,x,v) | d\gamma dt 
	\\ & \lesssim  \int_{\epsilon_1} ^T \int_{\gamma_+ \setminus \gamma_+^\epsilon} \left( \displaystyle {\min_{ [\epsilon_1 , T ] \times [ \gamma_+ \setminus \gamma_+^\epsilon ]} } \{ t, \tb (t,x,v) \} \times | h(t,x,v) | \right) d\gamma dt 
	\\ & \lesssim \int_0 ^T \int_{\gamma_+ \setminus \gamma_+^\epsilon } \int_{-\min\{t, \tb (t,x,v) \} } ^0 |h(t+s, X(t+s) ,V(t+s)) | ds d\gamma dt 
	\\ & + \tb(t,x,v) \int_0^T \int_{\gamma_+ \setminus \gamma_+ ^\epsilon } \int_{- \min \{ t, \tb (t,x,v) \}  }^0 |\partial_t h + V(t + \tau ) \cdot \nabla_x h + E(X(t + \tau )) \cdot \nabla_v h + \phi (V(t + \tau)) h|  \\ & ( t + \tau , X(t + \tau ) , V(t + \tau )) d\tau d\gamma dt
	\\ & \lesssim \int_0^T \| h(t) \|_1 dt + \int_0^T \| [\partial_t + v \cdot \nabla_x + E \cdot \nabla_v + \phi ] h(t) \| _1 dt
	\end{split}
	\]
	where in the last inequality we have used the identity (\ref{h_id}), and that $\tb(t,x,v) \le \frac{\lambda}{\delta} \times \mathrm{diam} (\O) $
	
	We now only need to show 
	\[
	\int_0 ^{\epsilon_1} \int_{\gamma_+ \setminus \gamma_+^\epsilon } | h(t,x,v) | d\gamma dt \lesssim_{\Omega, \epsilon, \epsilon_1} \| h_0 \|_1 + \int_0^{\epsilon_1} \| [\partial_t + v \cdot \nabla_x + E \cdot \nabla_v + \phi ] h(t) \| _1 dt.
	\]
	Because of our choice $\epsilon$ and $\epsilon_1$, $\tb(t,x,v) > t $ for all $(t,x,v) \in [0, \epsilon_1] \times \gamma_+ \setminus \gamma_+^\epsilon $. Then 
	\[
	\begin{split}
	|h(t,x,v)| =  &| h_0( X(0),V(0) )|  \\ &+  \int_{-t}^0 | \left[ \partial_t h + V(t + \tau ) \cdot \nabla_x h + E(X(t + \tau )) \cdot \nabla_v h + \phi (V(t + \tau)) h \right ] \\ & ( t + \tau , X(t + \tau ) , V(t + \tau )) |d \tau,
	\end{split}
	\] 
	where the second contribution is bounded, again from (\ref{h_id}), by 
	\[
	\begin{split}
	\int_0^{\epsilon_1}  & \int_{\gamma_+ \setminus \gamma_+^\epsilon} \int_{-t}^0   | \left[ \partial_t h + V(t + \tau ) \cdot \nabla_x h + E(X(t + \tau )) \cdot \nabla_v h + \phi (V(t + \tau)) h \right ]  \\ & ( t + \tau , X(t + \tau ) , V(t + \tau )) |d \tau
	\\ & \lesssim \int_0^{\epsilon_1} \| [\partial_t + v \cdot \nabla_x + E \cdot \nabla_v + \phi ] h(t) \| _1 dt.
	\end{split}
	\]
	And the contribution of the initial data is estimated by the same change of variable: $ (t,x,v) \mapsto (X(0;t,x,v), V(0;t,x,v)) \in \Omega \times \mathbb R^3$, as $\tb(t,x,v) > t$ for any $(t,x,v) \in [0, \epsilon_1] \times \gamma_+ \setminus \gamma_+^\epsilon$. So by (\ref{green2}) we have:
	\[
	\begin{split}
	\int_0^{\epsilon_1} & \int_{\gamma_+ \setminus \gamma_+^\epsilon} |h_0(X(0;t,x,v), V(0;t,x,v)) | d\gamma dt 
	\\ & \le \int_0^{T } \int_{ \gamma_+ } |h_0| (X(0;t,x,v), V(0;t,x,v) ) \textbf{1}_{ \{ t < \tb(t,x,v) \} }     d\gamma d t
	\\& = \iint_{\Omega \times \mathbb R^3 } |h_0 | \textbf{1}_{ \{ \tf(0,x,v) \le T' \} }   dx dv  
	\\ & \le \| h_0 \|_1.
	\end{split}
	\]\end{proof}

\unhide

\hide
Even the proof is a modification of the proof of in \cite{GKTT1}, we heavily rely on the decay assumption (\ref{decay_phi}) to obtain an uniform-in-time estimate. Therefore here we write the details of the proof.

Now for $(t,x,v) \in [\varepsilon_1,T]\times \gamma_+\setminus \gamma_+^\varepsilon,$
we integrate as $\int_{\varepsilon_1}^T \int_{\gamma_+\setminus\gamma_+^\varepsilon}
\int^0_{\min\{t,t_{\mathbf{b}}(x,v)\}}
(\ref{h_eq})
\dd s
\dd \gamma
\dd t$ to get
\begin{equation}\notag
\begin{split}
&\min\{\varepsilon_1 , \varepsilon^3\} \times \int^T_{\varepsilon_1} \int_{\gamma_+\setminus\gamma^\varepsilon_+} | h(t,x,v) | \mathrm{d} \gamma \mathrm{d} t \\
& \lesssim  \ \min_{[\varepsilon_1,T]\times [\gamma_+\setminus\gamma_+^\varepsilon]}\{t,t_{\mathbf{b}}(x,v)\} \times \int^T_{\varepsilon_1} \int_{\gamma_+\setminus\gamma^\varepsilon_+} |h(t,x,v) |\mathrm{d} \gamma \mathrm{d} t\\
& \lesssim \int_0^T  \int_{\gamma_+} \int^0_{- \min \{t,t_{\mathbf{b}}(x,v)\}}|h(t+s,x+sv,v)|   \mathrm{d} s  \mathrm{d} \gamma \mathrm{d} t \\
&  \ \ +   \int_0^T \int_{\gamma_+} \int^0_{-\min\{t,t_{\mathbf{b}}(x,v)\}}  |\partial_t h + v\cdot \nabla_x h + \varphi h|(t+ \tau, x+\tau v,v)\mathrm{d} \tau \mathrm{d} \gamma\mathrm{d} t\\
& \lesssim  \int_0^T ||h(t)||_1\mathrm{d} t + \int^T_0 || [\partial_t  + v\cdot \nabla_x   + \varphi]h(t)  ||_1\mathrm{d} t,
\end{split}
\end{equation}
where we have used the integration identity (\ref{integration}), and (40) of \cite{Guo10} to obtain
$t_{\mathbf{b}}(x,v) \geq C_{\Omega}  {|n(x)\cdot v|}/ {|v|^2} \geq C_{\Omega}  \varepsilon^3$ for $(x,v) \in \gamma_+\backslash \gamma_+^\varepsilon$. Now we choose $\varepsilon_1 = \varepsilon_1(\Omega,\varepsilon)$ as
\[
\varepsilon_1 \leq C_{\Omega} \varepsilon^3 \leq \inf_{(x,v)\in \gamma_+\setminus \gamma_+^\varepsilon} t_{\mathbf{b}}(x,v).
\]
We only need to show, for $\varepsilon_1 \leq  C_{\Omega} \varepsilon^3,$
\begin{equation}\notag
\begin{split}
\int^{\varepsilon_1}_0 \int_{\gamma_+\setminus \gamma_+^\varepsilon} |h(t,x,v)| \mathrm{d} \gamma \mathrm{d} t \lesssim_{\Omega, \varepsilon,\varepsilon_1} || h_0||_1 + \int_0^{\varepsilon_1} ||[\partial_t + v\cdot \nabla_x +\varphi] h(t) ||_1\mathrm{d} t.
\end{split}
\end{equation}
Because of our choice $\varepsilon$ and $\varepsilon_1$,  $t_{\mathbf{b}}(x,v) >t$ for all $(t,x,v) \in [0,\varepsilon_1]\times \gamma_+\setminus\gamma_+^\varepsilon.$ Then
\[
|h(t,x,v)| \lesssim |h_0(x-tv,v)|+ \int_0^t \Big{|}[\partial_t + v\cdot \nabla_x + \varphi(v)]h(s,x-(t-s)v,v)\Big{|}\mathrm{d} s,
\]
where the second contribution is bounded, from (\ref{integration}), by
\begin{equation}\notag
\begin{split}
&\int^{\varepsilon_1}_0 \int_{\gamma_+\setminus \gamma_+^\varepsilon} \int_0^t \Big{|}[\partial_t + v\cdot \nabla_x + \varphi(v)]h(s,x-(t-s)v,v)\Big{|}\mathrm{d} s \mathrm{d} \gamma \mathrm{d} t \\
\lesssim& \int_0^{\varepsilon_1} || [\partial_t + v\cdot \nabla_x + \varphi(v)]h(t)||_1 \mathrm{d} t.
\end{split}
\end{equation}

Consider the initial datum contribution of $|h_0(x-tv,v)|$: We may assume $\partial_{x_3}\xi(x_0)\neq 0$. By the implicit function theorem $\partial\Omega$ can be represented locally by the graph $\eta=\eta(x_1,x_2)$ satisfying $\xi(x_1,x_2,\eta(x_1,x_2))=0$ and $(\partial_{x_1}\eta(x_1,x_2),\partial_{x_2}\eta(x_1,x_2))= (-\partial_{x_1}\xi / \partial_{x_3}\xi,-\partial_{x_2}\xi / \partial_{x_3}\xi )$ at $(x_1,x_2,\eta(x_1,x_2)).$ We define the change of variables
\begin{equation}\notag
\begin{split}
(x,t)\in \partial\Omega \cap\{x\sim  x_0\} \times [0,\varepsilon_1] \mapsto y = x-tv \in \bar{\Omega},
\end{split}
\end{equation}
where $
\left|\frac{\partial y}{\partial(x,t)}\right|  = -v_1\frac{\partial_{x_1}\xi}{\partial_{x_3}\xi} - v_2 \frac{\partial_{x_2}\xi}{\partial_{x_3}\xi}- v_3.
$

Therefore
\begin{equation*}
\begin{split}
|n(x)\cdot v|\mathrm{d} S_x \mathrm{d} t &= (n(x)\cdot v) \Big[1+ \big(\frac{\partial_{x_1}\xi}{\partial_{x_3}\xi}\big)^2
+ \big(\frac{\partial_{x_2}\xi}{\partial_{x_3}\xi}\big)^2 \Big]^{1/2} \mathrm{d} x_1 \mathrm{d} x_2 \mathrm{d} t\\
& = \left[ -v_1\frac{\partial_{x_1}\xi}{\partial_{x_3}\xi} - v_2 \frac{\partial_{x_2}\xi}{\partial_{x_3}\xi}- v_3 \right]\mathrm{d} x_1 \mathrm{d} x_2 \mathrm{d} t = \mathrm{d} y,
\end{split}
\end{equation*}
and
$\int_0^{\varepsilon_1} \int_{\gamma_+\setminus \gamma_+^\varepsilon \cap \{x\sim x_0\}}
|h_0(x-tv,v)| \mathrm{d} \gamma \mathrm{d} t \ \lesssim_{\varepsilon,\varepsilon_1, x_0}
\ \iint_{\Omega\times\mathbb{R}^3} |h_0(y,v)| \mathrm{d} y \mathrm{d} v.$ Since $\partial\Omega$
is compact we can choose a finite covers of $\partial\Omega$ and repeat the same argument for each piece to conclude \begin{equation}\notag
\int_0^{\varepsilon_1} \int_{\gamma_+\setminus \gamma_+^\varepsilon } |h_0(x-tv,v)| \mathrm{d} \gamma \mathrm{d} t \ \lesssim_{\Omega, \varepsilon,\varepsilon_1} \ \iint_{\Omega\times\mathbb{R}^3} |h_0(y,v)| \mathrm{d} y \mathrm{d} v.
\end{equation}\unhide

\hide

\begin{lemma}
	\label{le:ukai} Assume we have $\phi(t,x) \in C^1_x$ such that, for $\lambda>0$, $\delta>0$, 
	\begin{equation}
	\label{decay_phi}
	\sup_{t \geq 0} e^{\lambda t} \| \nabla \phi (t) \|_{\infty} \leq \delta.
	\end{equation}
	We also assume $\frac{1}{C}\langle v\rangle \leq \psi(t,x,v)\leq C \langle v\rangle$ for some $C>0$. For any small parameter 
	\Be\label{lower_bound_e}
	\e> \frac{2\delta}{  \lambda}>0,
	\Ee
	there exists a
	constant $C_{\delta, \e ,\Omega }>0$, which does not depend on $t$, such that for any $h$ in $%
	L^{1}_{loc}
	([0,\infty);L^{1}(\Omega \times \mathbb{R}^{3}))$ with $\partial _{t}h+v\cdot
	\nabla _{x}h
	- \nabla \phi \cdot \nabla_v h + \psi h$ in $L^{1}_{loc}( [0, \infty);L^{1}(\Omega \times \mathbb{R}%
	^{3}))$, we have, for all $t\geq 0$,
	\begin{equation*}
	\begin{split}
	&\int_{0}^{t}\int_{\gamma _{+}\setminus \gamma _{+}^{\varepsilon }}|h|\mathrm{%
		d}\gamma \mathrm{d}s\\
	&\leq C_{\delta, \varepsilon ,\Omega }\left\{ \
	||h_{0}||_{1}+\int_{0}^{t}  \| h(s)\|_{1}+\big{\Vert} 
	[
	\partial
	_{t}+v\cdot \nabla _{x}- \nabla_x\phi \cdot \nabla_v + \psi ]h(s)\big{\Vert} _{1} \mathrm{d}s\ \right\}.\end{split}
	\end{equation*}%
	Furthermore, for any $(t,x,v)$ in $[0, \infty)\times \Omega \times \mathbb{R}^{3}$
	the function $
	h
	(t+s^{\prime },
	X(t+s^\prime;t,x,v)
	,
	V(t+s^\prime;t,x,v))
	$ is absolutely continuous in
	$s^{\prime } \in (-\min \{t_{\mathbf{b}}(t,x,v),t\}, t_{%
		\mathbf{f}}(t,x, v) )$.
\end{lemma}
\begin{proof}From (\ref{decay_phi}), for $0 \leq s \leq t$
	\Be
	\begin{split}\notag
		&|V(s;t,x,v)-v|\\
		\leq 
		&  \ 
		\int_s^t |  \nabla\phi (\tau, X(\tau;t,x,v))| \dd \tau  \\
		\leq   &
		\ 
		\frac{\delta}{\lambda} e^{- \lambda s} \{1- e^{-\lambda (t-s)}\} \\
		\leq
		& \ 
		\frac{\delta}{\lambda} .
	\end{split}
	\Ee
	Hence if $|v| \geq \e$ with the condition (\ref{lower_bound_e}) then 
	\Be
	\begin{split}\notag
		&|X(s;t,x,v)-x|\\
		\geq  & \ \left|\int^t_s v\right| -  \int^t_s| V(\tau;t,x,v)-v | \dd \tau \\
		\geq & \  |v| (t-s) -\frac{\delta}{\lambda} (t-s) \\
		\geq & \  \frac{\delta}{\lambda}(t-s) .
	\end{split}
	\Ee
	Therefore 
	\Be\label{upper_tb}
	\tb(t,x,v) \leq \frac{\lambda}{\delta} \times \mathrm{diam} (\O).
	\Ee

	\hide
	Even the proof is a modification of the proof of in \cite{GKTT1}, we heavily rely on the decay assumption (\ref{decay_phi}) to obtain an uniform-in-time estimate. Therefore here we write the details of the proof.

	Now for $(t,x,v) \in [\varepsilon_1,T]\times \gamma_+\setminus \gamma_+^\varepsilon,$
	we integrate as $\int_{\varepsilon_1}^T \int_{\gamma_+\setminus\gamma_+^\varepsilon}
	\int^0_{\min\{t,t_{\mathbf{b}}(x,v)\}}
	(\ref{h_eq})
	\dd s
	\dd \gamma
	\dd t$ to get
	\begin{equation}\notag
	\begin{split}
	&\min\{\varepsilon_1 , \varepsilon^3\} \times \int^T_{\varepsilon_1} \int_{\gamma_+\setminus\gamma^\varepsilon_+} | h(t,x,v) | \mathrm{d} \gamma \mathrm{d} t \\
	& \lesssim  \ \min_{[\varepsilon_1,T]\times [\gamma_+\setminus\gamma_+^\varepsilon]}\{t,t_{\mathbf{b}}(x,v)\} \times \int^T_{\varepsilon_1} \int_{\gamma_+\setminus\gamma^\varepsilon_+} |h(t,x,v) |\mathrm{d} \gamma \mathrm{d} t\\
	& \lesssim \int_0^T  \int_{\gamma_+} \int^0_{- \min \{t,t_{\mathbf{b}}(x,v)\}}|h(t+s,x+sv,v)|   \mathrm{d} s  \mathrm{d} \gamma \mathrm{d} t \\
	&  \ \ +   \int_0^T \int_{\gamma_+} \int^0_{-\min\{t,t_{\mathbf{b}}(x,v)\}}  |\partial_t h + v\cdot \nabla_x h + \varphi h|(t+ \tau, x+\tau v,v)\mathrm{d} \tau \mathrm{d} \gamma\mathrm{d} t\\
	& \lesssim  \int_0^T ||h(t)||_1\mathrm{d} t + \int^T_0 || [\partial_t  + v\cdot \nabla_x   + \varphi]h(t)  ||_1\mathrm{d} t,
	\end{split}
	\end{equation}
	where we have used the integration identity (\ref{integration}), and (40) of \cite{Guo10} to obtain
	$t_{\mathbf{b}}(x,v) \geq C_{\Omega}  {|n(x)\cdot v|}/ {|v|^2} \geq C_{\Omega}  \varepsilon^3$ for $(x,v) \in \gamma_+\backslash \gamma_+^\varepsilon$. Now we choose $\varepsilon_1 = \varepsilon_1(\Omega,\varepsilon)$ as
	\[
	\varepsilon_1 \leq C_{\Omega} \varepsilon^3 \leq \inf_{(x,v)\in \gamma_+\setminus \gamma_+^\varepsilon} t_{\mathbf{b}}(x,v).
	\]
	We only need to show, for $\varepsilon_1 \leq  C_{\Omega} \varepsilon^3,$
	\begin{equation}\notag
	\begin{split}
	\int^{\varepsilon_1}_0 \int_{\gamma_+\setminus \gamma_+^\varepsilon} |h(t,x,v)| \mathrm{d} \gamma \mathrm{d} t \lesssim_{\Omega, \varepsilon,\varepsilon_1} || h_0||_1 + \int_0^{\varepsilon_1} ||[\partial_t + v\cdot \nabla_x +\varphi] h(t) ||_1\mathrm{d} t.
	\end{split}
	\end{equation}
	Because of our choice $\varepsilon$ and $\varepsilon_1$,  $t_{\mathbf{b}}(x,v) >t$ for all $(t,x,v) \in [0,\varepsilon_1]\times \gamma_+\setminus\gamma_+^\varepsilon.$ Then
	\[
	|h(t,x,v)| \lesssim |h_0(x-tv,v)|+ \int_0^t \Big{|}[\partial_t + v\cdot \nabla_x + \varphi(v)]h(s,x-(t-s)v,v)\Big{|}\mathrm{d} s,
	\]
	where the second contribution is bounded, from (\ref{integration}), by
	\begin{equation}\notag
	\begin{split}
	&\int^{\varepsilon_1}_0 \int_{\gamma_+\setminus \gamma_+^\varepsilon} \int_0^t \Big{|}[\partial_t + v\cdot \nabla_x + \varphi(v)]h(s,x-(t-s)v,v)\Big{|}\mathrm{d} s \mathrm{d} \gamma \mathrm{d} t \\
	\lesssim& \int_0^{\varepsilon_1} || [\partial_t + v\cdot \nabla_x + \varphi(v)]h(t)||_1 \mathrm{d} t.
	\end{split}
	\end{equation}
	
	Consider the initial datum contribution of $|h_0(x-tv,v)|$: We may assume $\partial_{x_3}\xi(x_0)\neq 0$. By the implicit function theorem $\partial\Omega$ can be represented locally by the graph $\eta=\eta(x_1,x_2)$ satisfying $\xi(x_1,x_2,\eta(x_1,x_2))=0$ and $(\partial_{x_1}\eta(x_1,x_2),\partial_{x_2}\eta(x_1,x_2))= (-\partial_{x_1}\xi / \partial_{x_3}\xi,-\partial_{x_2}\xi / \partial_{x_3}\xi )$ at $(x_1,x_2,\eta(x_1,x_2)).$ We define the change of variables
	\begin{equation}\notag
	\begin{split}
	(x,t)\in \partial\Omega \cap\{x\sim  x_0\} \times [0,\varepsilon_1] \mapsto y = x-tv \in \bar{\Omega},
	\end{split}
	\end{equation}
	where $
	\left|\frac{\partial y}{\partial(x,t)}\right|  = -v_1\frac{\partial_{x_1}\xi}{\partial_{x_3}\xi} - v_2 \frac{\partial_{x_2}\xi}{\partial_{x_3}\xi}- v_3.
	$
	
	Therefore
	\begin{equation*}
	\begin{split}
	|n(x)\cdot v|\mathrm{d} S_x \mathrm{d} t &= (n(x)\cdot v) \Big[1+ \big(\frac{\partial_{x_1}\xi}{\partial_{x_3}\xi}\big)^2
	+ \big(\frac{\partial_{x_2}\xi}{\partial_{x_3}\xi}\big)^2 \Big]^{1/2} \mathrm{d} x_1 \mathrm{d} x_2 \mathrm{d} t\\
	& = \left[ -v_1\frac{\partial_{x_1}\xi}{\partial_{x_3}\xi} - v_2 \frac{\partial_{x_2}\xi}{\partial_{x_3}\xi}- v_3 \right]\mathrm{d} x_1 \mathrm{d} x_2 \mathrm{d} t = \mathrm{d} y,
	\end{split}
	\end{equation*}
	and
	$\int_0^{\varepsilon_1} \int_{\gamma_+\setminus \gamma_+^\varepsilon \cap \{x\sim x_0\}}
	|h_0(x-tv,v)| \mathrm{d} \gamma \mathrm{d} t \ \lesssim_{\varepsilon,\varepsilon_1, x_0}
	\ \iint_{\Omega\times\mathbb{R}^3} |h_0(y,v)| \mathrm{d} y \mathrm{d} v.$ Since $\partial\Omega$
	is compact we can choose a finite covers of $\partial\Omega$ and repeat the same argument for each piece to conclude \begin{equation}\notag
	\int_0^{\varepsilon_1} \int_{\gamma_+\setminus \gamma_+^\varepsilon } |h_0(x-tv,v)| \mathrm{d} \gamma \mathrm{d} t \ \lesssim_{\Omega, \varepsilon,\varepsilon_1} \ \iint_{\Omega\times\mathbb{R}^3} |h_0(y,v)| \mathrm{d} y \mathrm{d} v.
	\end{equation}\unhide\end{proof}\unhide
\unhide


\section{Desingularization via Mixing in the Velocity space}
The main purpose of this section is proving Proposition \ref{prop_int_alpha}.

\begin{proposition}\label{prop_int_alpha}
	Assume $E(t,x) \in C^1_x$ is given, and both (\ref{nE=0}) and (\ref{decay_E}) hold, and 
	\Be\label{decay_phi_2}
	\sup_{t\geq 0}e^{\Lambda_2 t} \| \nabla_x E(t)\|_{\infty} \leq \delta_2 \ll1 .
	\Ee
	Then for all $0<\sigma<1$ and $N>1$ and for all $s\geq 0,$ $x \in \bar{\O}$,
	\Be\label{NLL_split2}
	\int_{|u| \leq N } 
	\frac{   \dd u
	}{\alpha_{f,\e}(s,x,u)^{ \sigma}}
	\lesssim_{\sigma, \O, \Lambda_1, \delta_1, \Lambda_2, \delta_2,N}  1, 
	\Ee
	and, for any $0< \kappa\leq2$,  
	\Be\label{NLL_split3} 
	\int_{  |u|\geq N} \frac{e^{-C|v-u|^2}}{|v-u|^{2-\kappa}} \frac{1}{\alpha_{f,\e}(s,x,u)^\sigma} \dd u
	\lesssim_{\sigma, \O, \Lambda_1, \delta_1, \Lambda_2, \delta_2,N,\kappa}  1.
	\Ee
\end{proposition}

The key element of the proof is the next change of variable formula.  
\begin{lemma}\label{COV_boundary}Assume (\ref{decay_E}) and (\ref{decay_phi_2}) and $
\delta_2 (\Lambda_2)^{-2} \ll1 $. For $x \in \bar{\O}$ and $t\geq 0$, we define a map
	\Be
	u\in\{ u \in \R^3: \tb(t,x,u) \leq t+1 \} \mapsto (\tb(t,x,u), \xb(t,x,u)) \in \R \times \p\O,
	\Ee
	with $\det \left( \frac{\p (\tb(t,x,u), \xb(t,x,u)) }{\p u}\right) \gtrsim \frac{\tb(t,x,u)^3}{|n(\xb(t,x,u)) \cdot \vb(t,x,u)| }$.
	
	For $g\geq 0$ and all $x \in \bar{\O}$ and $t\geq 0$
	\Be \label{COV_xbtb}\begin{split}
		&\int_{\R^{3}} 
		\mathbf{1}_{\tb(t,x,u)< t+1}g(t,x,u)  \dd u\\
		\lesssim & \  \int _{\p\O} \int \frac{  |n(\xb )  \cdot \vb (t,x,u)|}{|\tb |^{3}} g(t,x, u(\xb, \tb))\dd \tb \dd S_{\xb},
	\end{split}\Ee 
	where the lower end of $\tb$-integration has a lower bound as 
	\Be\label{tb_lower}
	\tb(t,x,u) \geq \frac{|\xb(t,x,u) -x|}{\max_{t-\tb(t,x,u) \leq s \leq t} |V(s;t,x,u)|} .
	\Ee 
\end{lemma}

We need the following lemma to prove Lemma \ref{COV_boundary}.
\begin{lemma}\label{est_X_v}
	Assume (\ref{decay_phi_2}) with $\Lambda_2+ \delta_2+ \e \leq 1$.
	\hide
	\Be\begin{split}\notag
		\frac{d}{ds} 
		\begin{bmatrix}
			X(s;t,x,v)\\
			V(s;t,x,v)
		\end{bmatrix} = 
		\begin{bmatrix}
			V(s;t,x,v)\\
			- \nabla_x \phi (s,X(s;t,x,v))
		\end{bmatrix}.
	\end{split}\Ee\unhide
	Then there exists $C>0$ such that 
	\Be \label{result_X_v}\begin{split}
		|\nabla_vX(s;t,x,v)|
	 \leq& \  Ce^{C \delta_2 (\Lambda_{2})^{-2}   } 
		|t-s|
	,\\ &  \ \ \text{for all}   \ \max(t - \tb(t,x,v), - \e )\leq s \leq t.
\end{split}	\Ee  
\end{lemma}
\begin{proof}Recall that $\nabla_v X
(t;t,x,v)= \nabla_v x=0$ and $\nabla_v V(t;t,x,v)=\nabla_v v=\text{Id}_{3,3}$. From (\ref{nabla_Hamilton}) we have 
\Be\begin{split}\label{|Hamilton|}
\frac{d}{ds} |\nabla_v X(s;t,x,v)| \lesssim & \ |\nabla_v V(s;t,x,v)|,\\
\frac{d}{ds} |\nabla_v V(s;t,x,v)| \lesssim&  \ |\nabla_x E(s, X(s;t,x,v))| |\nabla_v X(s;t,x,v)| \\
\lesssim  & \   \delta_2 e^{-\Lambda_2 |s|} |\nabla_v X(s;t,x,v)|.
\end{split}\Ee

Then
\Be\begin{split}\notag
|\nabla_v X(s;t,x,v)| &\lesssim \int^t_s|\nabla_v V(s^\prime;t,x,v)| \dd s^\prime\\
&\lesssim
  |t-s| + \int^t_s  \int_{s^\prime}^t \delta_2 e^{-\Lambda_2  |s^{\prime\prime} |} |\nabla_v X(s^{\prime\prime};t,x,v)|  \dd s^{\prime\prime}
 \dd s^\prime\\
 &
\lesssim  |t-s| 
 + \int_s^t  \int_{s }^{s^{\prime\prime}} \delta_2 e^{-\Lambda_2 |s^{\prime\prime} | } |\nabla_v X(s^{\prime\prime};t,x,v)|  \dd s^\prime \dd s^{\prime\prime}\\
 &\lesssim   |t-s|  + \int_s^t  |s^{\prime\prime}-s|  \delta_2 e^{-\Lambda_2  |s^{\prime\prime} | } |\nabla_v X(s^{\prime\prime};t,x,v)|    \dd s^{\prime\prime}.
\end{split}\Ee
By the Gronwall's inequality (set $\tilde{s}:=t-s$ if necessary), there exists some $C \gg1$ such that 
\Be\notag
\begin{split}
|\nabla_v X(s;t,x,v)| & \lesssim |t-s| \exp\bigg( \int^t_0 |s^{\prime\prime}-s|  \delta_2 e^{-\Lambda_2 |s^{\prime\prime}|}      \dd s^{\prime\prime} \\
&  \ \ \ \ \ \ \ \ \   \ \ \ \ \ \ \ \ \ \ 
+ \int_{\min(s, 0)}^0 |s^{\prime\prime}-s|  \delta_2 e^{-\Lambda_2 |s^{\prime\prime}|}      \dd s^{\prime\prime} 
\bigg)\\
&\lesssim |t-s| \exp\left( \frac{\delta_2}{(\Lambda_2)^2} \int^\infty_0
(\Lambda_2 s^{\prime\prime})
   e^{-(\Lambda_2 s^{\prime\prime})}     \dd (\Lambda_2 s^{\prime\prime}) 
   + (\e)^2 \delta_2 
   \right)\\
   &\lesssim|t-s| e^{(\e)^2 \delta_2} e^{C  {\delta_2}{ (\Lambda_2)^{-2} 
   }}\\
&   \lesssim|t-s|  e^{C  {\delta_2}{ (\Lambda_2)^{-2} 
   }}, \  \ 
 for  \ all \   \max(t - \tb(t,x,v), - \e )\leq s \leq t .
\end{split}
\Ee
\hide

From the same argument,  
 estimate and the second line of (\ref{|Hamilton|}), for $C \gg C_1$,
\Be\begin{split}\notag
|\nabla_v V(s;t,x,v)| & \lesssim  1 +
\int^t_s
\int^t_{s^\prime}
\delta_2 e^{-\Lambda_2 |s^\prime|} |\nabla_v V(s^{\prime\prime};t,x,v)|
\dd s^{\prime\prime}
\dd s^\prime
\\
&\lesssim 1+ 
\int^t_s|s^{\prime\prime} - s |  
\delta_2 e^{-\Lambda_2 |s^\prime|} |\nabla_v V(s^{\prime\prime};t,x,v)|
\dd s^{\prime\prime} 
\\
& \lesssim 1+ 
\int^t_s  
\\
&\leq   C  e^{C   {\delta_2}{(\Lambda_2)^{-2}}}, \  \ 
 for  \ all \   \max(t - \tb(t,x,v), - \e )\leq s \leq t .
\end{split}\Ee \unhide
By choosing $C\gg 1$ large enough we conclude (\ref{result_X_v}).
\hide
Now following the similar estimate, we obtain that 
\Be\notag
\begin{split}
|\nabla_v V(s;t,x,v)| \lesssim & \ 1 +
\int^t_s 
\delta_2 e^{-\Lambda_2 s^\prime}\int^t_{s^\prime} 
|\nabla_v V(s^{\prime\prime};t,x,v)|
\dd s^{\prime\prime}
\dd s^\prime
\\
\lesssim & \ 1+ 
\int^t_s 
\int^{s^{\prime\prime}}_s \delta_2 e^{-\Lambda_2 s^\prime} \dd s^\prime
|\nabla_v V(s^{\prime\prime};t,x,v)|
\dd s^{\prime\prime}\\
\lesssim & \ 1+ 
\int^t_s 
\delta_2 \frac{
e^{-\Lambda_2 s } - e^{-\Lambda_2 s^{\prime\prime} } 
}{\Lambda_2}
|\nabla_v V(s^{\prime\prime};t,x,v)|
\dd s^{\prime\prime}
\end{split}
\Ee\unhide
\end{proof}

\hide

In the proof of this lemma we utilize an exponential decay of $\phi$ crucially to exclude a possible exponential growth of a non-autonomous system, as the following way.  
\begin{lemma}\label{lemma_nonauto}
	For some $\delta^\prime>0$, $\lambda^\prime>0$ we assume that $A(s;t)\geq 0$ and $B(s;t)\geq 0$ exist for all $0 \leq s \leq  t$ and satisfy 
	\Be\label{2_diff_ineq}
	\left| \frac{d}{ds} \begin{bmatrix} A(s;t)\\
		B(s;t)
	\end{bmatrix}\right|
	\leq \begin{bmatrix}
		0 & 1 \\
		\delta^\prime e^{-\lambda^\prime s}  & 0
	\end{bmatrix}
	\begin{bmatrix} A(s;t)\\
		B(s;t)
	\end{bmatrix}, \ \ \  \begin{bmatrix}
		A(t;t) \\  B(t;t)
	\end{bmatrix} =
	\begin{bmatrix}
		0 \\ 
		1
	\end{bmatrix}.
	\Ee
	Then for all $0 \leq s \leq t$
	\Be\label{2_diff_bound}
	\begin{bmatrix} A(s;t)\\
		B(s;t)
	\end{bmatrix}
	\leq Ce^{C \lambda^\prime \sqrt{\delta^\prime}}
	\begin{bmatrix}t-s \\ 1
	\end{bmatrix}.
	\Ee
\end{lemma}

\begin{proof}[\textbf{Proof of Lemma \ref{est_X_v}}]
	Recall (\ref{non-auto}).

	Now we set $A(s;t) : = |\nabla_v X(s;t,x,v)|$ and $B(s;t) : = |\nabla_v V(s;t,x,v)|$. From (\ref{decay_phi_2}) we have a bound of (\ref{2_diff_ineq}) with $\delta^\prime = \delta_2$ and $\lambda^\prime = \lambda_2$. Now we apply Lemma \ref{lemma_nonauto} and conclude (\ref{result_X_v}).\end{proof} \begin{proof}[\textbf{Proof of Lemma \ref{lemma_nonauto}}] It is convenient to adopt a backward temporal variable
	\Be\label{tilde_s}
	\tilde{s} : = t-s \in [0, t].
	\Ee
	We define 
	\Be\label{tilde_AB}
	\tilde{A}(\tilde{s};t): = {A}(t- \tilde{s};t) \ \  \text{and} \ \  {B}(\tilde{s};t): = B(t- \tilde{s};t).  
	\Ee
	
	\hide
	\textit{Step 1. }\textit{Comparison principle}:  For a non-negative matrix $M(\tilde{s};t)   \in \mathbb{R}^{n \times n}$,
	let us assume that $H(\tilde{s};t) \in \mathbb{R}^n$ and $J(\tilde{s};t)\in \mathbb{R}^n$ satisfy 
	\Be \label{eqnts_H_J}
	\frac{dH(\tilde{s};t)}{d\tilde{s}}  =  M (\tilde{s};t) H(\tilde{s};t),\ \
	\frac{dJ(\tilde{s};t)}{d\tilde{s}}  \leq  M (\tilde{s};t) J(\tilde{s};t) \ \ \text{for all }  \tilde{s}  \in [0,t].
	\Ee 
	If $H(0;t) \geq J(0;t) \geq 0$ and bounded then 
	\Be\label{HJ}
	H(\tilde{s};t) \geq J(\tilde{s};t) \geq 0 \text{ for all }  \tilde{s}  \in [0,t].
	\Ee
	
	Let us introduce $H^\e$ solving the same equation 
	with an $\e$-perturbed initial datum $H^\e (0;t) = H(0;t)+ \e (1, \cdots, 1)^T$. By the Gronwall's inequality, $|H^\e (\tilde{s};t) - H(\tilde{s};t)|\lesssim_n \e \exp\big(\int_0^{\tilde{s}} |M(\tilde{\tau};t)|   \dd \tau\big)$. Therefore for fixed $t$ and $\tilde{s}$, $H^\e(\tilde{s};t) \rightarrow H(\tilde{s};t)$ as $\e \downarrow 0$. Moreover from a temporal integration over the difference of equations in (\ref{eqnts_H_J}),
	\Be \begin{split}\label{H-J}
		H^\e(\tilde{s};t) - J(\tilde{s};t)   &\geq     H (0;t) - J(0;t)  + \e    (1, \cdots, 1)^T \\
		& +  \int_0^{\tilde{s}} M(\tau ;t)[H^\e(\tau ;t) - J(\tau ;t)]  \dd \tau 
		.\end{split}
	\Ee 
	
	We prove the estimate (\ref{HJ}) with exchanging $H(\tilde{s};t)$ by $H^\e (\tilde{s};t)$ first. Assume the statement is false. Then there exists $\tilde{s}_* \in [0,t)$ such that $H^\e(\tilde{s};t)\geq J(\tilde{s};t)$ for $\tilde{s} \in [0,\tilde{s}_*]$ but for some $i=1,2,\cdots, n$ we have $H^\e_i(\tilde{\tau}; t)< J_i(\tilde{\tau}; t)$ for $\tilde{\tau}> \tilde{s}_*$ with $|\tilde{\tau} - \tilde{s}_*| \ll 1$. 
	
	On the other hand, from (\ref{H-J}) and other conditions, for all $i$,  
	\Bes
	H ^\e(\tilde{\tau};t) - J (\tilde{\tau};t) &\geq& \e + \int^{\tilde{\tau}}_{\tilde{s}_*} M(\tau;t) [H ^\e(\tau;t) - J (\tau;t)] \dd \tau\\
	&\geq& \e - O_{M,H,J}(| {\tilde{\tau}} -{\tilde{s}_*}  |).
	\Ees
	If $| {\tilde{\tau}} -{\tilde{s}_*}  |$ is small enough then $H^\e (\tilde{\tau};t) > J(\tilde{\tau};t)$. This is a contradiction. 
	
	Then by passing a limit $\e \downarrow 0$ we prove (\ref{HJ}).

	\vspace{4pt}
	
	\unhide
	
	\textit{Step 1. }  
	We consider a temporally discretized problem: for $k \in \mathbb{N}$ with $1 \leq k \leq N$,
	\Be\label{auton}
	\begin{split}
		\frac{d}{d\tilde{s}} \begin{bmatrix}
			\tilde{A}(\tilde{s};N\tilde{T})\\
			\tilde{B}(\tilde{s};N\tilde{T})
		\end{bmatrix}
		=
		\begin{bmatrix}
			0 & 1 \\
			\delta^\prime e^{- \lambda^\prime (N-k )\tilde{T} } & 0
		\end{bmatrix}\begin{bmatrix}
			\tilde{A}(\tilde{s};N\tilde{T})\\
			\tilde{B}(\tilde{s};N\tilde{T})
		\end{bmatrix}  \\ \text{for} \  (k-1) \tilde{T}\leq \tilde{s} \leq k\tilde{T} ,\end{split}
	\Ee
	with an initial datum 
	\Be
	\begin{bmatrix}
		\tilde{A} (0;N \tilde{T}) \\  \tilde{B}(0;N \tilde{T})
	\end{bmatrix} =
	\begin{bmatrix}
		0 \\ 
		1
	\end{bmatrix}.\label{initial_AB}
	\Ee
	Since (\ref{auton}) is an autonomous ODE system we have a formula of a solution.  
	Note that we have the following diagonalization for $a \geq 0$ 
	\Be\notag\label{diagonalization_a}
	\begin{bmatrix}
		\frac{1}{2} & \frac{1}{2 \sqrt{a}}\\
		\frac{1}{2} & \frac{-1}{2 \sqrt{a}}
	\end{bmatrix}
	\begin{bmatrix}
		0 & 1 \\
		a & 0
	\end{bmatrix}
	\begin{bmatrix}
		1 & 1\\
		\sqrt{a}   &-\sqrt{a} 
	\end{bmatrix}
	= \begin{bmatrix}
		\sqrt{a}   & 0 \\
		0 &  -\sqrt{a} 
	\end{bmatrix}.
	\Ee
	Applying this diagonalization to (\ref{auton}) with $a = \delta^\prime e^{- \lambda^\prime  (N-k)\tilde{T}}$, 
	we deduce that 
	for $k \in \mathbb{N}$ with $1 \leq k \leq N$,
	\Be\begin{split}
		& \begin{bmatrix}
			\tilde{A}( k\tilde{T};N\tilde{T})\\
			\tilde{B}( k\tilde{T};N\tilde{T})
		\end{bmatrix}\\
		= &
		\begin{bmatrix}
			\cosh\big( \tilde{T}\sqrt{\delta^\prime e^{- \lambda ^\prime \tilde{T} (N-k)   }} \big)
			&
			\frac{\sinh\big( \tilde{T}\sqrt{\delta^\prime e^{- \lambda^\prime \tilde{T} (N-k)  }} \big)}{\sqrt{\delta^\prime e^{- \lambda^\prime \tilde{T} (N-k)  }} }
			\\
			\sqrt{\delta^\prime e^{- \lambda^\prime \tilde{T} (N-k)  }}  \sinh\big( \tilde{T}\sqrt{\delta^\prime e^{- \lambda^\prime \tilde{T} (N-k)  }} \big)
			& 
			\cosh\big(\tilde{T}\sqrt{\delta^\prime e^{-\lambda^\prime \tilde{T} (N-k)  }} \big)
		\end{bmatrix}\\
		& \times 
		\begin{bmatrix}
			\tilde{A}( (k-1)\tilde{T};N\tilde{T})\\
			\tilde{B}( (k-1)\tilde{T};N\tilde{T})
		\end{bmatrix}.
	\end{split}\label{k_k-1}
	\Ee 
	
	Throughout \textit{Step 2}- \textit{Step 4}, we claim that for $1\leq k \leq N$
	\Be\label{tilde_A_k}\begin{split}  
		\tilde{A}(k\tilde{T}; N \tilde{T}) 
		&=
		\sum_{i=1}^k \sum_{\substack{ \ell_1 + \cdots \ell_i =k \\ \ell_1, \cdots \ell_i \geq 1 }} 
		(1- e^{ -\frac{\lambda^\prime \tilde{T}}{2}})^{i-1}
		e^{- \frac{\lambda^\prime\tilde{T}}{2}(k-\ell_1 - (i-1)) }\\
		& \ \ \ \ \  \times   \frac{\sinh\big( \sum_{j=1}^{\ell_1} \tilde{T}\sqrt{\delta^\prime} e^{- \frac{\lambda^\prime \tilde{T}}{2}(N-j) }\big)}{
			\sqrt{\delta^\prime} e^{- \frac{ \lambda^\prime \tilde{T}}{2} (N-\ell_1)} 
		}\\
		&  \ \ \ \ \   \times \prod_{m=2}^i \cosh \big(
		\sum_{j=1}^{\ell_m} \tilde{T} \sqrt{\delta^\prime} e^{- \frac{\lambda^\prime \tilde{T}}{2} (N- (\ell_1 +  \cdots + \ell_{m-1} +j) ) }
		\big), 
	\end{split}\Ee
	and 
	\Be\label{tilde_B_k}\begin{split}
		\begin{array}{rl}
			&\tilde{B}(k\tilde{T}; N \tilde{T}) \\
			=&  \cosh \big( \sum_{j=1}^{k}  \tilde{T} \sqrt{\delta^\prime} e^{- \frac{\lambda^\prime \tilde{T} }{2} (N-j)}  \big)
			\\
			& \  + 
			\sum_{i=2}^k \sum_{\substack{ \ell_1 + \cdots + \ell_i =k\\  \ell_1 , \cdots , \ell_i \geq 1}}  (e^{\frac{\lambda^\prime}{2} \tilde{T} }-1)^{i-1} \\
			& \ \ \   \times 
			\sinh \big(
			\sum_{j=1}^{\ell_i} \tilde{T} \sqrt{\delta^\prime} e^{- \frac{\lambda^\prime \tilde{T}}{2} (N- (\ell_1 + \cdots + \ell_{i-1} + j))}
			\big)\\
			& \ \ \  \times \prod_{m=2}^{i-1} \cosh \big( 
			\sum_{j=1}^{\ell_m} \tilde{T} \sqrt{\delta^\prime} e^{- \frac{\lambda^\prime \tilde{T}}{2} (N- (\ell_1 + \cdots + \ell_{m-1} + j))}
			\big)\\
			&  \ \ \   \times \sinh\big(
			\sum_{j=1}^{\ell_1} \tilde{T} \sqrt{\delta^\prime} e^{- \frac{\lambda^\prime \tilde{T}}{2} (N-j)}
			\big)
			. \end{array}
	\end{split}\Ee

	\vspace{4pt} \textit{Step 2. } Let us check $k=1$ first. From (\ref{k_k-1}) and (\ref{initial_AB}),
	\Be\begin{split}\label{AB_1}
		\begin{bmatrix}
			\tilde{A}( \tilde{T};N\tilde{T})\\
			\tilde{B}( \tilde{T};N\tilde{T})
		\end{bmatrix} 
		= \begin{bmatrix}
			\frac{\sinh\big( \tilde{T}\sqrt{\delta^\prime e^{- \lambda^\prime (N-1)\tilde{T} }} \big)}{\sqrt{\delta^\prime e^{- \lambda^\prime (N-1)\tilde{T} }} }
			\\
			\cosh\big(\tilde{T}\sqrt{\delta^\prime e^{-\lambda^\prime (N-1)\tilde{T} }} \big)
		\end{bmatrix}.
	\end{split}\Ee
	
	For $k=2$, from (\ref{k_k-1}) and (\ref{AB_1}),
	\Be\begin{split}\notag
		\begin{array}{rl}
			&\begin{bmatrix}
				\tilde{A}( 2\tilde{T};N\tilde{T})\\
				\tilde{B}( 2\tilde{T};N\tilde{T})
			\end{bmatrix} \\
			= &
			\begin{bmatrix}
				\cosh\big( \tilde{T}\sqrt{\delta^\prime e^{- \lambda ^\prime (N-2)\tilde{T} }} \big)
				&
				\frac{\sinh\big( \tilde{T}\sqrt{\delta^\prime e^{- \lambda^\prime (N-2)\tilde{T} }} \big)}{\sqrt{\delta^\prime e^{- \lambda^\prime (N-2)\tilde{T} }} }
				\\
				\sqrt{\delta^\prime e^{- \lambda^\prime (N-2)\tilde{T} }}  \sinh\big( \tilde{T}\sqrt{\delta^\prime e^{- \lambda^\prime (N-2)\tilde{T} }} \big)
				& 
				\cosh\big(\tilde{T}\sqrt{\delta^\prime e^{-\lambda^\prime (N-2)\tilde{T} }} \big)
			\end{bmatrix}\\
			& \times 
			\begin{bmatrix}
				\frac{\sinh\big( \tilde{T}\sqrt{\delta^\prime e^{- \lambda^\prime (N-1)\tilde{T} }} \big)}{\sqrt{\delta^\prime e^{- \lambda^\prime (N-1)\tilde{T} }} }
				\\
				\cosh\big(\tilde{T}\sqrt{\delta^\prime e^{-\lambda^\prime (N-1)\tilde{T} }} \big)
			\end{bmatrix}.\end{array}
	\end{split}\Ee
	Now we use the identity 
	\Be\begin{split}\label{hyp_identity}
		\sinh (a+b) &= \sinh a \cosh b + \cosh a \sinh b,\\
		\cosh (a+b) &=  \cosh a \cosh b + \sinh a \sinh b.
	\end{split}\Ee
	Then we can check (\ref{tilde_A_k}) and (\ref{tilde_B_k}) for $k=2$ as
	\Be\begin{split} \label{A_2}
		\begin{array}{rl}
			&\tilde{A}( 2\tilde{T};N\tilde{T}) \\
			=& \cosh\big( \tilde{T}\sqrt{\delta^\prime e^{- \lambda ^\prime (N-2)\tilde{T} }} \big)\frac{\sinh\big( \tilde{T}\sqrt{\delta^\prime e^{- \lambda^\prime (N-1)\tilde{T} }} \big)}{\sqrt{\delta^\prime e^{- \lambda^\prime (N-1)\tilde{T} }} }\\
			&  + \frac{\sinh\big( \tilde{T}\sqrt{\delta^\prime e^{- \lambda^\prime (N-2)\tilde{T} }} \big)}{\sqrt{\delta^\prime e^{- \lambda^\prime (N-2)\tilde{T} }} }\cosh\big(\tilde{T}\sqrt{\delta^\prime e^{-\lambda^\prime (N-1)\tilde{T} }} \big)\\
			=& \frac{  \sinh\big( \sum_{j=1}^2
				\tilde{T}\sqrt{\delta^\prime} e^{- \frac{\lambda^\prime \tilde{T} }{2} (N-j) }
				\big)}{\sqrt{\delta^\prime} e^{- \frac{\lambda^\prime \tilde{T} }{2} (N-2) } } \\
			& +(1- e^{-\frac{\lambda^\prime \tilde{T}}{2} })
			\frac{ \sinh\big( \tilde{T}\sqrt{\delta^\prime} e^{- \frac{\lambda^\prime \tilde{T} }{2} (N-1) } \big)}{\sqrt{\delta^\prime} e^{- \frac{\lambda^\prime  \tilde{T} }{2}(N-1)}}
			\cosh\big( \tilde{T}\sqrt{\delta^\prime} e^{- \frac{\lambda ^\prime \tilde{T} }{2} (N-2) } \big),
		\end{array}
	\end{split}\Ee
	and 
	\Be
	\begin{split}\label{B_2}
		\begin{array}{rl}
			&\tilde{B}( 2\tilde{T};N\tilde{T}) \\
			=&
			\sqrt{\delta^\prime e^{- \lambda^\prime (N-2)\tilde{T} }}  \sinh\big( \tilde{T}\sqrt{\delta^\prime e^{- \lambda^\prime (N-2)\tilde{T} }} \big)
			\frac{\sinh\big( \tilde{T}\sqrt{\delta^\prime e^{- \lambda^\prime (N-1)\tilde{T} }} \big)}{\sqrt{\delta^\prime e^{- \lambda^\prime (N-1)\tilde{T} }} }\\&
			+
			\cosh\big(\tilde{T}\sqrt{\delta^\prime e^{-\lambda^\prime (N-2)\tilde{T} }} \big)
			\cosh\big(\tilde{T}\sqrt{\delta^\prime e^{-\lambda^\prime (N-1)\tilde{T} }} \big)\\
			=&
			\cosh\big( \sum_{j=1}^2\tilde{T}\sqrt{\delta^\prime} e^{-\frac{\lambda^\prime \tilde{T} }{2} (N-j) }  \big)
			\\
			&+(e^{\frac{\lambda^\prime \tilde{T}}{2}}-1) \sinh\big( \tilde{T}\sqrt{\delta^\prime} e^{- \frac{\lambda^\prime \tilde{T} }{2} (N-2) } \big)\sinh\big( \tilde{T}\sqrt{\delta^\prime} e^{- \frac{\lambda^\prime  \tilde{T}}{2} (N-1) } \big).\end{array}
	\end{split}
	\Ee

	\vspace{4pt} \textit{Step 3. }  Now we assume that (\ref{tilde_A_k}) and (\ref{tilde_B_k}) hold for $1 \leq k \leq N-1$. From (\ref{k_k-1}), (\ref{tilde_A_k}), and (\ref{tilde_B_k}), we obtain that $\tilde{A}((k+1)\tilde{T}; N \tilde{T})$ equals 
	\begin{eqnarray}
	&&  \left\{  \begin{array}{rl}
	&\cosh \big( \tilde{T} \sqrt{\delta^\prime} e^{- \frac{\lambda^\prime \tilde{T}}{2} (N- (k+1))  }\big)\\
	& \times \sum_{i=1}^k \sum_{\substack{ \ell_1 + \cdots \ell_i =k \\ \ell_1, \cdots \ell_i \geq 1 }} 
	(1- e^{ -\frac{\lambda^\prime \tilde{T}}{2}})^{i-1}
	e^{- \frac{\lambda^\prime\tilde{T}}{2}(k-\ell_1 - (i-1)) }\\
	&  \ \ \ \ \   \times \prod_{m=2}^i \cosh \big(
	\sum_{j=1}^{\ell_m} \tilde{T} \sqrt{\delta^\prime} e^{- \frac{\lambda^\prime \tilde{T}}{2} (N- (\ell_1 +  \cdots + \ell_{m-1} +j) ) }
	\big)\\
	& \ \ \ \ \  \times   \frac{\sinh\big( \sum_{j=1}^{\ell_1} \tilde{T}\sqrt{\delta^\prime} e^{- \frac{\lambda^\prime \tilde{T}}{2}(N-j) }\big)}{
		\sqrt{\delta^\prime} e^{- \frac{ \lambda^\prime \tilde{T}}{2} (N-\ell_1)} 
	}
	\end{array} \right .\label{A_k+1_1} \\
	&&  \left\{ \begin{array}{rl}&+\frac{\sinh\big( \tilde{T} \sqrt{\delta^\prime} e^{- \frac{\lambda^\prime \tilde{T}}{2} (N- (k+1))  }\big)}{
		\sqrt{\delta^\prime} e^{- \frac{\lambda^\prime \tilde{T}}{2} (N- (k+1))  } 
	}
	\cosh \big( \sum_{j=1}^{k}  \tilde{T} \sqrt{\delta^\prime} e^{- \frac{\lambda^\prime \tilde{T} }{2} (N-j)}  \big)\end{array}
	\right . \label{A_k+1_2}\\
	&&   \left\{ \begin{array}{rl}
	&+ \frac{\sinh\big( \tilde{T} \sqrt{\delta^\prime} e^{- \frac{\lambda^\prime \tilde{T}}{2} (N- (k+1))  }\big)}{
		\sqrt{\delta^\prime} e^{- \frac{\lambda^\prime \tilde{T}}{2} (N- (k+1))  } 
	}\\
	& \times 
	\sum_{i=2}^k \sum_{\substack{ \ell_1 + \cdots + \ell_i =k\\  \ell_1 , \cdots , \ell_i \geq 1}}  (e^{\frac{\lambda^\prime}{2} \tilde{T} }-1)^{i-1}  \\
	& \ \ \ \ \ \times 
	\sinh \big(
	\sum_{j=1}^{\ell_i} \tilde{T} \sqrt{\delta^\prime} e^{- \frac{\lambda^\prime \tilde{T}}{2} (N- (\ell_1 + \cdots + \ell_{i-1} + j))}
	\big)
	\\
	& \ \ \ \ \ \times \prod_{m=2}^{i-1} \cosh \big( 
	\sum_{j=1}^{\ell_m} \tilde{T} \sqrt{\delta^\prime} e^{- \frac{\lambda^\prime \tilde{T}}{2} (N- (\ell_1 + \cdots + \ell_{m-1} + j))}
	\big)\\
	&  \ \ \ \ \   \times \sinh\big(
	\sum_{j=1}^{\ell_1} \tilde{T} \sqrt{\delta^\prime} e^{- \frac{\lambda^\prime \tilde{T}}{2} (N-j)}
	\big)
	.\end{array} \right .\label{A_k+1_3}
	\end{eqnarray}
	
	First we split 
	\Be\begin{split}\notag
		(\ref{A_k+1_1})  =    (\ref{A_k+1_1}) _{i=1} +  (\ref{A_k+1_1}) _{i>1} 
		:=  
		\Big\{   \cdots \times \sum_{i=1}^1 \cdots \Big\}   
		+  \Big\{   \cdots \times \sum_{i=2}^{k }\cdots \Big\} .
	\end{split}\Ee
	Then by $\frac{1}{e^{- \frac{ \lambda^\prime \tilde{T}}{2} (N-k)}}
	= \frac{1}{e^{- \frac{ \lambda^\prime \tilde{T}}{2} (N-(k+1))}} + \frac{1
		- e^{-\frac{\lambda^\prime \tilde{T}}{2}}
	}{e^{- \frac{ \lambda^\prime \tilde{T}}{2} (N-k)}},
	$ 
	\Be\begin{split}\label{A_k+1_1+2}
		&(\ref{A_k+1_1})_{i=1} + (\ref{A_k+1_2}) \\
		= & \frac{\sinh\big(\sum_{j=1}^{k+1}  \tilde{T} \sqrt{\delta^\prime} e^{- \frac{\lambda^\prime \tilde{T} }{2} (N-j)} \big)}{
			\sqrt{\delta^\prime} e^{- \frac{\lambda^\prime \tilde{T}}{2} (N- (k+1))  } 
		}\\   +&  
		(1
		- e^{-\frac{\lambda^\prime \tilde{T}}{2}})
		\cosh \big( \tilde{T} \sqrt{\delta^\prime} e^{- \frac{\lambda^\prime \tilde{T}}{2} (N- (k+1))  }\big)
		\frac{\sinh\big( \sum_{j=1}^{k} \tilde{T}\sqrt{\delta^\prime} e^{- \frac{\lambda^\prime \tilde{T}}{2}(N-j) }\big)}{
			\sqrt{\delta^\prime} e^{- \frac{ \lambda^\prime \tilde{T}}{2} (N-k)} 
		}.
	\end{split}\Ee

	Now we rewrite $(\ref{A_k+1_1})_{i>1}$ as 
	\Be\begin{split}\notag
		\begin{array}{rl}
			(\ref{A_k+1_1})_{i>1} 
			= &
			\cosh \big( \tilde{T} \sqrt{\delta^\prime} e^{- \frac{\lambda^\prime \tilde{T}}{2} (N- (k+1))  }\big)\\
			& \times \sum_{i=2}^k \sum_{\substack{ \ell_1 + \cdots \ell_i =k \\ \ell_1, \cdots \ell_i \geq 1 }} 
			(1- e^{ -\frac{\lambda^\prime \tilde{T}}{2}})^{i-1}
			e^{- \frac{\lambda^\prime\tilde{T}}{2}(k-\ell_1 - (i-1)) }\\
			& \ \ \ \ \ \times \cosh \big(
			\sum_{j=1}^{\ell_i} \tilde{T} \sqrt{\delta^\prime} e^{- \frac{\lambda^\prime \tilde{T}}{2} (N- (\ell_1 +  \cdots + \ell_{i-1} +j) ) }
			\big)
			\\ 
			&  \ \ \ \ \   \times \prod_{m=2}^{i-1} \cosh \big(
			\sum_{j=1}^{\ell_m} \tilde{T} \sqrt{\delta^\prime} e^{- \frac{\lambda^\prime \tilde{T}}{2} (N- (\ell_1 +  \cdots + \ell_{m-1} +j) ) }
			\big)\\
			& \ \ \ \ \  \times   \frac{\sinh\big( \sum_{j=1}^{\ell_1} \tilde{T}\sqrt{\delta^\prime} e^{- \frac{\lambda^\prime \tilde{T}}{2}(N-j) }\big)}{
				\sqrt{\delta^\prime} e^{- \frac{ \lambda^\prime \tilde{T}}{2} (N-\ell_1)} 
			}.\end{array}
	\end{split}
	\Ee
	Let us consider (\ref{A_k+1_3}) as polynomials of hyperbolic sine and hyperbolic cosine functions. Then the coefficient of (\ref{A_k+1_3}) equals $\frac{ (1-e^{-\frac{\lambda^\prime\tilde{T}}{2}  } )^{i-1} e^{\frac{\lambda^\prime \tilde{T}}{2} (i-1)}}{
		\sqrt{\delta^\prime} e^{- \frac{\lambda^\prime \tilde{T}}{2} (N- (k+1))  }}$ which equals
	\[
	\frac{(1-e^{-\frac{\lambda^\prime\tilde{T}}{2}  } )^{i-1}e^{- \frac{\lambda^\prime \tilde{T}}{2} (k - \ell_1 - (i-1))} }{\sqrt{\delta^\prime} e^{-\frac{\lambda^\prime \tilde{T}}{2} (N-\ell_1)}} e^{- \frac{\lambda^\prime \tilde{T}}{2}} .
	\]
	Note that the coefficient of (\ref{A_k+1_3}) equals the coefficient of $(\ref{A_k+1_1})_{i>1}$ multiplied by $e^{- \frac{\lambda^\prime \tilde{T}}{2}} $. We split $(\ref{A_k+1_1})_{i>1}$ using $1=  e^{- \frac{\lambda^\prime \tilde{T}}{2}}+(1-e^{- \frac{\lambda^\prime \tilde{T}}{2}}  ) $. We combine the first piece of this splitting of $(\ref{A_k+1_1})_{i>1}$ and (\ref{A_k+1_3}) using (\ref{hyp_identity}). Then we can conclude the following identity:
	\begin{eqnarray}
	&& (\ref{A_k+1_1})_{i>1}+ (\ref{A_k+1_3})\nonumber\\ 
	& & \left .\begin{array}{rl}
	& =
	\sum_{i=2}^k \sum_{\substack{ \ell_1 + \cdots \ell_i =k \\ \ell_1, \cdots \ell_i \geq 1 }} 
	(1- e^{ -\frac{\lambda^\prime \tilde{T}}{2}})^{i-1}
	e^{- \frac{\lambda^\prime\tilde{T}}{2}(k + 1-\ell_1 - ( 
		i-1) 
		) }\\
	&  \ \ \ \ \ \times \cosh \big(
	\sum_{j=1}^{\ell_i + 1} \tilde{T} \sqrt{\delta^\prime}
	e^{-\frac{\lambda^\prime \tilde{T}}{2} (N - (\ell_1 + \cdots \ell_{i-1} + j)) }
	\big)\\
	&  \ \ \ \ \   \times \prod_{m=2}^{i-1} \cosh \big(
	\sum_{j=1}^{\ell_m} \tilde{T} \sqrt{\delta^\prime} e^{- \frac{\lambda^\prime \tilde{T}}{2} (N- (\ell_1 +  \cdots + \ell_{m-1} +j) ) }
	\big)\\
	& \ \ \ \ \  \times   \frac{\sinh\big( \sum_{j=1}^{\ell_1} \tilde{T}\sqrt{\delta^\prime} e^{- \frac{\lambda^\prime \tilde{T}}{2}(N-j) }\big)}{
		\sqrt{\delta^\prime} e^{- \frac{ \lambda^\prime \tilde{T}}{2} (N-\ell_1)} 
	} \end{array}\right\}\label{72_74_A} \\
	&&
	\left . \begin{array}{rl}   &+
	\cosh \big( \tilde{T} \sqrt{\delta^\prime} e^{- \frac{\lambda^\prime \tilde{T}}{2} (N- (k+1))  }\big)\\
	&  \ \ \ \ \ \times \sum_{i=2}^k \sum_{\substack{ \ell_1 + \cdots \ell_i =k \\ \ell_1, \cdots \ell_i \geq 1 }} 
	(1- e^{ -\frac{\lambda^\prime \tilde{T}}{2}})^{i }
	e^{- \frac{\lambda^\prime\tilde{T}}{2}(k-\ell_1 - (i-1)) }\\
	&  \ \ \ \ \   \times \prod_{m=2}^{i
	} \cosh \big(
	\sum_{j=1}^{\ell_m} \tilde{T} \sqrt{\delta^\prime} e^{- \frac{\lambda^\prime \tilde{T}}{2} (N- (\ell_1 +  \cdots + \ell_{m-1} +j) ) }
	\big)\\
	& \ \ \ \ \  \times   \frac{\sinh\big( \sum_{j=1}^{\ell_1} \tilde{T}\sqrt{\delta^\prime} e^{- \frac{\lambda^\prime \tilde{T}}{2}(N-j) }\big)}{
		\sqrt{\delta^\prime} e^{- \frac{ \lambda^\prime \tilde{T}}{2} (N-\ell_1)} 
	}.
	\end{array} \right\}\label{72_74_B}
	\end{eqnarray}

	%
	%
	
	Now we consider (\ref{72_74_A}). For $i=2,\cdots, k$, we change a summation index $\ell_i+1$ to $ \ell_i$. Since the original $\ell_i$ is running from $1$ to $k - (i-1)$ the newly defined $\ell \in \{2, \cdots, k+1 - (i-1)  \}$. 
	Then we can rewrite (\ref{72_74_A}) as 
	\Be\label{72_74_A_dash}
	\begin{array}{rl}
		(\ref{72_74_A})  & =
		\sum_{i=2}^k \sum_{\substack{ \ell_1 + \cdots +  \ell_{i-1 } + {\ell}_i  =k+1 \\ \ell_1, \cdots,  \ell_{i-1} \geq 1,  {\ell}_i \geq 2 }} 
		(1- e^{ -\frac{\lambda^\prime \tilde{T}}{2}})^{i-1}
		e^{- \frac{\lambda^\prime\tilde{T}}{2}(k + 1-\ell_1 - ( 
			i-1) 
			) }\\
		&  \ \ \     \ \  \times \prod_{m=2}^{i } \cosh \big(
		\sum_{j=1}^{\ell_m} \tilde{T} \sqrt{\delta^\prime} e^{- \frac{\lambda^\prime \tilde{T}}{2} (N- (\ell_1 +  \cdots + \ell_{m-1} +j) ) }
		\big)\\
		& \ \ \   \ \   \times   \frac{\sinh\big( \sum_{j=1}^{\ell_1} \tilde{T}\sqrt{\delta^\prime} e^{- \frac{\lambda^\prime \tilde{T}}{2}(N-j) }\big)}{
			\sqrt{\delta^\prime} e^{- \frac{ \lambda^\prime \tilde{T}}{2} (N-\ell_1)} 
		} . \end{array} 
	\Ee

	On the other hand, by changing an index $i +1$ to $i  \in \{ 3,4, \cdots, k+1 \}$, we can rewrite (\ref{72_74_B}) as 
	\Be\notag
	\begin{array}{rl}    
		&    \sum_{ {i}=3}^{k+1} \sum_{\substack{ \ell_1 + \cdots \ell_{ {i}-1} =k \\ \ell_1, \cdots \ell_{ {i}-1} \geq 1 }} 
		(1- e^{ -\frac{\lambda^\prime \tilde{T}}{2}})^{ {i}-1}
		e^{- \frac{\lambda^\prime\tilde{T}}{2}((k+1)-\ell_1 - ( {i}-1)) }\\
		& \ \  \times  \cosh \big( \tilde{T} \sqrt{\delta^\prime} e^{- \frac{\lambda^\prime \tilde{T}}{2} (N- (k+1))  }\big) \\
		&  \ \  \times \prod_{m=2}^{{i}-1} \cosh \big(
		\sum_{j=1}^{\ell_m} \tilde{T} \sqrt{\delta^\prime} e^{- \frac{\lambda^\prime \tilde{T}}{2} (N- (\ell_1 +  \cdots + \ell_{m-1} +j) ) }
		\big)\\
		&  \ \  \times   \frac{\sinh\big( \sum_{j=1}^{\ell_1} \tilde{T}\sqrt{\delta^\prime} e^{- \frac{\lambda^\prime \tilde{T}}{2}(N-j) }\big)}{
			\sqrt{\delta^\prime} e^{- \frac{ \lambda^\prime \tilde{T}}{2} (N-\ell_1)} 
		}. 
	\end{array} 
	\Ee
	Note that if we set $\ell_{i}=1$ then for $\ell_1 + \cdots + \ell_{i-1} = k$ and $m=i$
	\Be \notag
	\sum_{j=1}^{\ell_m} \tilde{T} \sqrt{\delta^\prime} e^{- \frac{\lambda^\prime \tilde{T}}{2} (N- (\ell_1 +  \cdots + \ell_{m-1} +j) ) } 
	=  \tilde{T} \sqrt{\delta^\prime} e^{- \frac{\lambda^\prime \tilde{T}}{2} (N- (k +1) ) }.
	\Ee
	Hence we deduce 
	\Be\label{72_74_B_dash}
	\begin{array}{rl}    
		(\ref{72_74_B})=&   \sum_{i=3}^{k+1} \sum_{\substack{ \ell_1 + \cdots \ell_{i} =k+1 \\      \ell_1, \cdots  , \ell_{{i}-1} \geq 1, \ \ell_{i}=1 }} 
		(1- e^{ -\frac{\lambda^\prime \tilde{T}}{2}})^{{i}-1}
		e^{- \frac{\lambda^\prime\tilde{T}}{2}((k+1)-\ell_1 - ({i}-1)) }\\
		& \ \   \times \prod_{m=2}^{{i} } \cosh \big(
		\sum_{j=1}^{\ell_m} \tilde{T} \sqrt{\delta^\prime} e^{- \frac{\lambda^\prime \tilde{T}}{2} (N- (\ell_1 +  \cdots + \ell_{m-1} +j) ) }
		\big)\\
		& \ \    \times   \frac{\sinh\big( \sum_{j=1}^{\ell_1} \tilde{T}\sqrt{\delta^\prime} e^{- \frac{\lambda^\prime \tilde{T}}{2}(N-j) }\big)}{
			\sqrt{\delta^\prime} e^{- \frac{ \lambda^\prime \tilde{T}}{2} (N-\ell_1)} 
		} .
	\end{array} 
	\Ee
	
	Combining (\ref{72_74_A_dash}) and (\ref{72_74_B_dash}), we conclude that 
	\Be\label{72_74_A+B}
	\begin{split}
		&(\ref{72_74_A})+ 
		(\ref{72_74_B})\\
		& \begin{array}{rl}    
			=&   \sum_{i=3}^{k+1} \sum_{\substack{ \ell_1 + \cdots \ell_{i} =k+1 \\      \ell_1, \cdots  , \ell_{{i} } \geq 1 }} 
			(1- e^{ -\frac{\lambda^\prime \tilde{T}}{2}})^{{i}-1}
			e^{- \frac{\lambda^\prime\tilde{T}}{2}((k+1)-\ell_1 - ({i}-1)) }\\
			& \ \   \times \prod_{m=2}^{{i} } \cosh \big(
			\sum_{j=1}^{\ell_m} \tilde{T} \sqrt{\delta^\prime} e^{- \frac{\lambda^\prime \tilde{T}}{2} (N- (\ell_1 +  \cdots + \ell_{m-1} +j) ) }
			\big)\\
			& \ \    \times   \frac{\sinh\big( \sum_{j=1}^{\ell_1} \tilde{T}\sqrt{\delta^\prime} e^{- \frac{\lambda^\prime \tilde{T}}{2}(N-j) }\big)}{
				\sqrt{\delta^\prime} e^{- \frac{ \lambda^\prime \tilde{T}}{2} (N-\ell_1)} 
		}\end{array} \\
		& \begin{array}{rl}
			+ &
			\sum_{\substack{ \ell_1  + {\ell}_2 =k+1 \\ \ell_1  \geq 1,  {\ell}_2 \geq 2 }} 
			(1- e^{ -\frac{\lambda^\prime \tilde{T}}{2}}) 
			e^{- \frac{\lambda^\prime\tilde{T}}{2}( (k + 1)-\ell_1 - 1 
				) }\\
			&  \ \ \     \ \  \times  \cosh \big(
			\sum_{j=1}^{\ell_2} \tilde{T} \sqrt{\delta^\prime} e^{- \frac{\lambda^\prime \tilde{T}}{2} (N- (\ell_1  +j) ) }
			\big)\\
			& \ \ \   \ \   \times   \frac{\sinh\big( \sum_{j=1}^{\ell_1} \tilde{T}\sqrt{\delta^\prime} e^{- \frac{\lambda^\prime \tilde{T}}{2}(N-j) }\big)}{
				\sqrt{\delta^\prime} e^{- \frac{ \lambda^\prime \tilde{T}}{2} (N-\ell_1)} 
			} . \end{array}  \end{split}
	\Ee
	Finally from (\ref{A_k+1_1+2}) and (\ref{72_74_A+B}) we conclude the identity (\ref{tilde_A_k}) for $A((k+1) \tilde{T}; N \tilde{T})$. Note that the first term in the RHS of (\ref{A_k+1_1+2}) corresponds to $i=1$ case of $A((k+1) \tilde{T}; N \tilde{T})$. And $i=2$ case of $A((k+1) \tilde{T}; N \tilde{T})$ is obtained by combining second term of the RHS of (\ref{A_k+1_1+2}) and the last term of (\ref{72_74_A+B}). Finally $i \geq 3$ case of $A((k+1) \tilde{T}; N \tilde{T})$ is exactly same as the first term in the RHS of (\ref{72_74_A+B}).

	\vspace{4pt} \textit{Step 4. } Again we assume (\ref{tilde_A_k}) and (\ref{tilde_B_k}) for $1 \leq k \leq N-1$. From (\ref{k_k-1}), (\ref{tilde_A_k}), and (\ref{tilde_B_k}), we derive that $B((k+1) \tilde{T}; N \tilde{T})$ equals 
	\begin{eqnarray}
	&&  \left\{  \begin{array}{rl}
	& \sum_{i=1}^k \sum_{\substack{ \ell_1 + \cdots \ell_i =k \\ \ell_1, \cdots \ell_i \geq 1 }} 
	(  e^{ \frac{\lambda^\prime \tilde{T}}{2}}-1)^{i-1} e^{\frac{\lambda^\prime \tilde{T}}{2}}
	\sinh \big( \tilde{T} \sqrt{\delta^\prime} e^{- \frac{\lambda^\prime \tilde{T}}{2} (N- (k+1))  }\big)\\
	&  \ \ \ \ \   \times \prod_{m=2}^i \cosh \big(
	\sum_{j=1}^{\ell_m} \tilde{T} \sqrt{\delta^\prime} e^{- \frac{\lambda^\prime \tilde{T}}{2} (N- (\ell_1 +  \cdots + \ell_{m-1} +j) ) }
	\big)\\
	& \ \ \ \ \  \times  
	\sinh\big( \sum_{j=1}^{\ell_1} \tilde{T}\sqrt{\delta^\prime} e^{- \frac{\lambda^\prime \tilde{T}}{2}(N-j) }\big)
	\end{array} \right .\label{B_k+1_1} \\
	&&  \left\{ \begin{array}{rl}&+
	\cosh\big( \tilde{T} \sqrt{\delta^\prime} e^{- \frac{\lambda^\prime \tilde{T}}{2} (N- (k+1))  }\big)
	\cosh \big( \sum_{j=1}^{k}  \tilde{T} \sqrt{\delta^\prime} e^{- \frac{\lambda^\prime \tilde{T} }{2} (N-j)}  \big)\end{array}
	\right . \label{B_k+1_2}\\
	&&   \left\{ \begin{array}{rl}
	&+ 
	\sum_{i=2}^k \sum_{\substack{ \ell_1 + \cdots + \ell_i =k\\  \ell_1 , \cdots , \ell_i \geq 1}}  (e^{\frac{\lambda^\prime \tilde{T}  }{2}}-1)^{i-1} 
	\cosh\big( \tilde{T} \sqrt{\delta^\prime} e^{- \frac{\lambda^\prime \tilde{T}}{2} (N- (k+1))  }\big)
	\\
	& \ \ \ \ \ \times 
	\sinh \big(
	\sum_{j=1}^{\ell_i} \tilde{T} \sqrt{\delta^\prime} e^{- \frac{\lambda^\prime \tilde{T}}{2} (N- (\ell_1 + \cdots + \ell_{i-1} + j))}
	\big)
	\\
	& \ \ \ \ \ \times \prod_{m=2}^{i-1} \cosh \big( 
	\sum_{j=1}^{\ell_m} \tilde{T} \sqrt{\delta^\prime} e^{- \frac{\lambda^\prime \tilde{T}}{2} (N- (\ell_1 + \cdots + \ell_{m-1} + j))}
	\big)\\
	&  \ \ \ \ \   \times \sinh\big(
	\sum_{j=1}^{\ell_1} \tilde{T} \sqrt{\delta^\prime} e^{- \frac{\lambda^\prime \tilde{T}}{2} (N-j)}
	\big)
	.\end{array} \right .\label{B_k+1_3}
	\end{eqnarray}
	Using $e^{\frac{\lambda^\prime \tilde{T}}{2}}= 1+(e^{\frac{\lambda^\prime \tilde{T}}{2}}-1)$, we split (\ref{B_k+1_1}) into several pieces:
	\Be\begin{split}\label{B_split}
		(\ref{B_k+1_1})  =  & (\ref{B_k+1_1})_{1,i=1} + (\ref{B_k+1_1})_{1,i>1} +(\ref{B_k+1_1})_2
		\\
		: =& \sum_{i=1}^1 \sum (  e^{ \frac{\lambda^\prime \tilde{T}}{2}}-1)^{i-1} \sinh  \cdots +  \sum_{i>1} \sum (  e^{ \frac{\lambda^\prime \tilde{T}}{2}}-1)^{i-1}\sinh \cdots \\& + 
		\sum_{i=1}^k \sum (  e^{ \frac{\lambda^\prime \tilde{T}}{2}}-1)^{i }\sinh \cdots.
	\end{split}\Ee
	
	Note that $(\ref{B_k+1_1})_{1,i=1}$ equals $\sinh\big( \tilde{T} \sqrt{\delta^\prime} e^{- \frac{\lambda^\prime \tilde{T}}{2} (N- (k+1))  }\big)
	\sinh \big( \sum_{j=1}^{k}  \tilde{T} \sqrt{\delta^\prime} e^{- \frac{\lambda^\prime \tilde{T} }{2} (N-j)}  \big)$. We combine $(\ref{B_k+1_1})_{1,i=1}$ with (\ref{B_k+1_2}) and use (\ref{hyp_identity}) to derive 
	\Be\label{B_1,1+2}
	(\ref{B_k+1_1})_{1,i=1}  + (\ref{B_k+1_2})= \cosh \big( \sum_{j=1}^{k+1}  \tilde{T} \sqrt{\delta^\prime} e^{- \frac{\lambda^\prime \tilde{T} }{2} (N-j)}  \big).
	\Ee
	
	Now we combine $(\ref{B_k+1_1})_{1, i>1}$ with (\ref{B_k+1_3}). Note that $(\ref{B_k+1_1})_{1, i>1}$ equals 
	\[
	\begin{array}{rl}
	& \sum_{i=2}^k \sum_{\substack{ \ell_1 + \cdots \ell_i =k \\ \ell_1, \cdots \ell_i \geq 1 }} 
	(  e^{ \frac{\lambda^\prime \tilde{T}}{2}}-1)^{i-1}
	\sinh \big( \tilde{T} \sqrt{\delta^\prime} e^{- \frac{\lambda^\prime \tilde{T}}{2} (N- (k+1))  }\big)\\
	& \ \ \ \ \ \times 
	\cosh \big(
	\sum_{j=1}^{\ell_i} \tilde{T} \sqrt{\delta^\prime} e^{- \frac{\lambda^\prime \tilde{T}}{2} (N- (\ell_1 +  \cdots + \ell_{i-1} +j) ) }
	\big)
	\\
	&  \ \ \ \ \   \times \prod_{m=2}^{i-1} \cosh \big(
	\sum_{j=1}^{\ell_m} \tilde{T} \sqrt{\delta^\prime} e^{- \frac{\lambda^\prime \tilde{T}}{2} (N- (\ell_1 +  \cdots + \ell_{m-1} +j) ) }
	\big)\\
	& \ \ \ \ \  \times  
	\sinh\big( \sum_{j=1}^{\ell_1} \tilde{T}\sqrt{\delta^\prime} e^{- \frac{\lambda^\prime \tilde{T}}{2}(N-j) }\big).
	\end{array}
	\]
	Using (\ref{hyp_identity}), we deduce that $(\ref{B_k+1_1})_{1, i>1}   +  (\ref{B_k+1_3})$ equals 
	\Be\notag
	\begin{array}{rl}
		& \sum_{i=2}^k \sum_{\substack{ \ell_1 + \cdots \ell_i =k \\ \ell_1, \cdots \ell_i \geq 1 }} 
		(  e^{ \frac{\lambda^\prime \tilde{T}}{2}}-1)^{i-1} 
		\\
		& \ \ \ \ \ \times 
		\sinh \big(
		\sum_{j=1}^{\ell_i} \tilde{T} \sqrt{\delta^\prime} e^{- \frac{\lambda^\prime \tilde{T}}{2} (N- (\ell_1 +  \cdots + \ell_{i-1} +j) ) }
		+  \tilde{T} \sqrt{\delta^\prime} e^{- \frac{\lambda^\prime \tilde{T}}{2} (N- (k+1)) }
		\big)
		\\
		&  \ \ \ \ \   \times \prod_{m=2}^{i-1} \cosh \big(
		\sum_{j=1}^{\ell_m} \tilde{T} \sqrt{\delta^\prime} e^{- \frac{\lambda^\prime \tilde{T}}{2} (N- (\ell_1 +  \cdots + \ell_{m-1} +j) ) }
		\big)\\
		& \ \ \ \ \  \times  
		\sinh\big( \sum_{j=1}^{\ell_1} \tilde{T}\sqrt{\delta^\prime} e^{- \frac{\lambda^\prime \tilde{T}}{2}(N-j) }\big).
	\end{array}
	\Ee
	Now we change an index $\ell_i + 1 \in \{2, \cdots, k+1\}$ to $\tilde{\ell}_i$ and $\ell_j\equiv \tilde{\ell}_j$ for all $j<\ell_i$. Note that $\tilde{\ell}_1 + \cdots \tilde{\ell}_i= k+1$. With this new $\tilde{\ell}_i$, we have 
	\Be\begin{split}\notag
		\begin{array}{rl}
			&\sum_{j=1}^{\tilde{\ell}_i} \tilde{T} \sqrt{\delta^\prime} e^{- \frac{\lambda^\prime \tilde{T}}{2} (N- (\tilde{\ell}_1 +  \cdots +\tilde{\ell}_{i-1} +j) ) }
			\\
			&=\sum_{j=1}^{ {\ell}_{i-1}} \tilde{T} \sqrt{\delta^\prime} e^{- \frac{\lambda^\prime \tilde{T}}{2} (N- ( {\ell}_1 +  \cdots + {\ell}_{i-2} +j) ) }
			+ \sum_{j=\ell_i}^{ {\ell}_i+1} \tilde{T} \sqrt{\delta^\prime} e^{- \frac{\lambda^\prime \tilde{T}}{2} (N- ( {\ell}_1 +  \cdots + {\ell}_{i-1} +j) ) }\\
			&=\sum_{j=1}^{ {\ell}_{i }} \tilde{T} \sqrt{\delta^\prime} e^{- \frac{\lambda^\prime \tilde{T}}{2} (N- ( {\ell}_1 +  \cdots + {\ell}_{i-2} +j) ) }  + \tilde{T} \sqrt{\delta^\prime} e^{- \frac{\lambda^\prime \tilde{T}}{2} (N - (k+1))  },
	\end{array}\end{split}
	\Ee 
	where we have used $\ell_1 + \cdots \ell_i =k$. Removing the tilde on $\ell$ we derive 
	\Be\label{B_1,1+3}
	\begin{split}
		&(\ref{B_k+1_1})_{1, i>1}   +  (\ref{B_k+1_3}) \\
		& \begin{array}{rl}
			&= \sum_{i=2}^{k+1} \sum_{\substack{ \ell_1 + \cdots \ell_i =k+1 \\ \ell_1, \cdots \ell_{i-1} \geq 1, \ell_i\geq 2 }} 
			(  e^{ \frac{\lambda^\prime \tilde{T}}{2}}-1)^{i-1} 
			\\
			& \ \ \  \ \ \times 
			\sinh \big(
			\sum_{j=1}^{\ell_i} \tilde{T} \sqrt{\delta^\prime} e^{- \frac{\lambda^\prime \tilde{T}}{2} (N- (\ell_1 +  \cdots + \ell_{i-1} +j) ) }
			\big)
			\\
			&  \ \ \  \ \  \times \prod_{m=2}^{i-1} \cosh \big(
			\sum_{j=1}^{\ell_m} \tilde{T} \sqrt{\delta^\prime} e^{- \frac{\lambda^\prime \tilde{T}}{2} (N- (\ell_1 +  \cdots + \ell_{m-1} +j) ) }
			\big)\\
			& \ \ \  \ \   \times  
			\sinh\big( \sum_{j=1}^{\ell_1} \tilde{T}\sqrt{\delta^\prime} e^{- \frac{\lambda^\prime \tilde{T}}{2}(N-j) }\big).
		\end{array}\end
		{split}
		\Ee
		Note that $i=k+1$ is null since if $\ell_{k+1}\geq 2$ and $\ell_1, \cdots \ell_k \geq 1$ then $\ell_1  + \cdots \ell_i = k+1$ cannot hold. 
		
		Now we consider $(\ref{B_k+1_1})_2$. We change a summation index $\tilde{i} = i+1 \in \{2, \cdots, k+1\}$. Then $\ell_1 + \cdots + \ell_{\tilde{i}-1} =k$. Let us set $\ell_{\tilde{i}}=1$. Then $\ell_1 +\cdot + \ell_{\tilde{i}} = k+1$. Removing the tilde on $i$ we derive that
		\Be\label{B_1,2_new}
		\begin{split} 
			\begin{array}{rl}
				(\ref{B_k+1_1})_2 &= \sum_{i=2}^{k+1} \sum_{\substack{ \ell_1 + \cdots \ell_i =k+1 \\ \ell_1, \cdots \ell_{i-1} \geq 1, \ell_i=1 }} 
				(  e^{ \frac{\lambda^\prime \tilde{T}}{2}}-1)^{i-1} 
				\\
				& \ \ \  \ \ \times 
				\sinh \big(
				\sum_{j=1}^{\ell_i} \tilde{T} \sqrt{\delta^\prime} e^{- \frac{\lambda^\prime \tilde{T}}{2} (N- (\ell_1 +  \cdots + \ell_{i-1} +j) ) }
				\big)
				\\
				&  \ \ \  \ \  \times \prod_{m=2}^{i-1} \cosh \big(
				\sum_{j=1}^{\ell_m} \tilde{T} \sqrt{\delta^\prime} e^{- \frac{\lambda^\prime \tilde{T}}{2} (N- (\ell_1 +  \cdots + \ell_{m-1} +j) ) }
				\big)\\
				& \ \ \  \ \   \times  
				\sinh\big( \sum_{j=1}^{\ell_1} \tilde{T}\sqrt{\delta^\prime} e^{- \frac{\lambda^\prime \tilde{T}}{2}(N-j) }\big).
			\end{array}\end
			{split}
			\Ee
			From (\ref{B_1,1+2}), (\ref{B_1,1+3}), (\ref{B_1,2_new}), we conclude (\ref{tilde_B_k}) for $B((k+1)\tilde{T}; N \tilde{T})$. Note that (\ref{B_1,1+2}) is the first term in (\ref{tilde_B_k}) and sum of (\ref{B_1,1+3}), (\ref{B_1,2_new}) yields $i\geq 2$ terms of (\ref{tilde_B_k}). Now by the mathematical induction, from \textit{Step 3}, \textit{Step 4}, and \textit{Step 5}, we prove (\ref{tilde_A_k}) and (\ref{tilde_B_k}) for all $1 \leq k \leq N$.
			
			\vspace{4pt}
			\textit{Step 5. } With (\ref{tilde_A_k}) and (\ref{tilde_B_k}), we are ready to prove (\ref{2_diff_bound}). It suffices to prove 
			\Be\label{claim_AB}
			\tilde{A}(k \tilde{T}; N\tilde{T}) \leq C e^{C \lambda^\prime \sqrt{\delta^\prime}} k\tilde{T}, \ \ \  \tilde{B}(k \tilde{T}; N\tilde{T}) \lesssim1.
			\Ee
			Using geometric series, we have 
			\Be\label{sinh_gs}
			\begin{split}
				\begin{array}{rl}
					&\sinh\big( \sum_{j=1}^{\ell_1} \tilde{T}\sqrt{\delta^\prime} e^{- \frac{\lambda^\prime \tilde{T}}{2}(N-j) }\big)\\
					\leq
					& \ \sinh \Big( \tilde{T} \sqrt{\delta^\prime} e^{- \frac{\lambda^\prime \tilde{T}}{2} (N- \ell_1)} \frac{1- e^{-\frac{\lambda^\prime \tilde{T}}{2} \ell_1}}{1- e^{-\frac{\lambda^\prime \tilde{T}}{2}}}
					\Big),\\
					&\cosh \big(
					\sum_{j=1}^{\ell_m} \tilde{T} \sqrt{\delta^\prime} e^{- \frac{\lambda^\prime \tilde{T}}{2} (N- (\ell_1 +  \cdots + \ell_{m-1} +j) ) }
					\big)\\
					\leq & \ \cosh \Big(  
					\tilde{T} \sqrt{\delta^\prime} e^{- \frac{\lambda^\prime \tilde{T}}{2} \big(N- (\ell_1 + \cdots + \ell_m)\big)} \frac{1- e^{-\frac{\lambda^\prime \tilde{T}}{2} \ell_m}}{1- e^{-\frac{\lambda^\prime \tilde{T}}{2}}}
					\Big).\end{array}
			\end{split}
			\Ee   
			From 
			\Be\label{hyp_sin}
			\sinh x \leq \frac{3}{2} x \ \ \text{for} \ \ 0\leq x \ll 1,
			\Ee
			we derive, for $N-k>1$ and $\tilde{T}\gg 1$,
			\Be\notag
			\frac{\sinh\big( \sum_{j=1}^{\ell_1} \tilde{T}\sqrt{\delta^\prime} e^{- \frac{\lambda^\prime \tilde{T}}{2}(N-j) }\big)}{\sqrt{\delta^\prime} e^{- \frac{ \lambda^\prime \tilde{T}}{2} (N-\ell_1)} } \leq 2 \tilde{T}.
			\Ee
			From $\frac{1}{1- e^{- \frac{\lambda^\prime \tilde{T}} {2}}}< \frac{3}{2}$ for $\tilde{T}\gg_{\lambda^\prime}1$ and 
			\[
			\sum_{m=2}^ie^{- \frac{\lambda^\prime \tilde{T}}{2}  (N- (\ell_1 + \cdots + \ell_m) )} 
			\leq \sum_{\ell=1}^{k}e^{- \frac{\lambda^\prime \tilde{T}}{2}  (N-\ell )}  \leq \frac{3}{2} e^{- \frac{\lambda^\prime \tilde{T}}{2}  (N-k )},
			\]
			and from $\cosh x \leq e^x$, we deduce that
			\Be
			\begin{split}\notag
				\begin{array}{rl}
					& \prod_{m=2}^i \cosh \big( 
					\sum_{j=1}^{\ell_m} \tilde{T} \sqrt{\delta^\prime} e^{- \frac{\lambda^\prime \tilde{T}}{2} (N- (\ell_1 +  \cdots + \ell_{m-1} +j) ) }
					\big)\\
					\leq & \  \prod_{m=2}^i \exp  \Big(  \frac{3}{2}
					\tilde{T} \sqrt{\delta^\prime} e^{- \frac{\lambda^\prime \tilde{T}}{2} \big(N- (\ell_1 + \cdots + \ell_m)\big)}  
					\Big)\\
					\leq & \ \exp  \Big(\frac{3}{2}
					\tilde{T} \sqrt{\delta^\prime} \sum_{\ell=1}^k  
					e^{- \frac{\lambda^\prime \tilde{T}}{2} \big(N- \ell\big)}  
					\Big)\\
					\leq & 
					\exp
					\Big(
					2 \tilde{T} \sqrt{\delta^\prime} e^{- \frac{\lambda^\prime \tilde{T}}{2} (N-k )}
					\Big)
					.\end{array}
			\end{split}
			\Ee 
			Hence we have a bound  
			\Be\label{tilde_A_k}\begin{split}\notag
				\begin{array}{rl}
					& \tilde{A}(k\tilde{T}; N \tilde{T}) \\
					\leq & \
					\sum_{i=1}^k \sum_{\substack{ \ell_1 + \cdots \ell_i =k \\ \ell_1, \cdots \ell_i \geq 1 }} 
					(1- e^{ -\frac{\lambda^\prime \tilde{T}}{2}})^{i-1}
					e^{- \frac{\lambda^\prime\tilde{T}}{2}(k-\ell_1 - (i-1)) }\\
					& \times 
					2\tilde{T}
					\exp
					\Big(
					2 \tilde{T} \sqrt{\delta^\prime} e^{- \frac{\lambda^\prime \tilde{T}}{2}  (N-k )}
					\Big)\\
					\leq & \ 2\tilde{T} e^{2 \tilde{T} \sqrt{\delta^\prime}} \sum_{i=1}^k \sum_{\substack{ \ell_1 + \cdots \ell_i =k \\ \ell_1, \cdots \ell_i \geq 1 }} 
					(1- e^{ -\frac{\lambda^\prime \tilde{T}}{2}})^{i-1}
					e^{- \frac{\lambda^\prime\tilde{T}}{2}(k-\ell_1 - (i-1)) }.\end{array}
			\end{split}\Ee
			
			We change this summation of $\sum_{i=1}^k \sum_{\substack{ \ell_1 + \cdots +\ell_i =k \\ \ell_1, \cdots \ell_i \geq 1 }} \cdots$ to 
			\Be\notag
			2\tilde{T} e^{2 \tilde{T} \sqrt{\delta^\prime}} \sum_{\ell_1=1}^k \sum_{i=1}^{k- (\ell_1-1)} 
			(1- e^{ -\frac{\lambda^\prime \tilde{T}}{2}})^{i-1}
			e^{- \frac{\lambda^\prime\tilde{T}}{2}(k-\ell_1 - (i-1)) }
			\sum_{\substack{ \ell_2 + \cdots + \ell_i = k- \ell_1 \\
					\ell_2, \cdots, \ell_i \geq 1}}1.
			\Ee
			Then by defining $\bar{\ell}_2+1= \ell_2, \cdots, \bar{\ell}_i+1= \bar{\ell}_i$, 
			\Be\notag
			\sum_{\substack{ \ell_2 + \cdots + \ell_i = k- \ell_1 \\
					\ell_2, \cdots, \ell_i \geq 1}} 1=  \sum_{\substack{\bar{\ell}_2 + \cdots +\bar{\ell}_i = k- \ell_1- (i-1) \\
					\bar{\ell}_2, \cdots, \bar{\ell}_i \geq 0}} 1.
			\Ee
			We have a formula for a $(i-1)$-combination with repetitions from $k- \ell_1- (i-1)$ elements:
			\Be\notag
			\binom{ k- \ell_1- (i-1) + (i-1) -1  }{ i-1} = \binom{k-\ell_1-1}{i-1}.
			\Ee
			
			Now we need to consider
			\Be\notag
			2\tilde{T} e^{2 \tilde{T} \sqrt{\delta^\prime}} \sum_{\ell_1=1}^k e^{- \frac{\lambda^\prime \tilde{T}}{2} (k-\ell_1)}\sum_{i=1}^{k- (\ell_1-1)}\binom{k-\ell_1-1}{i-1} 
			(e^{ \frac{\lambda^\prime \tilde{T}}{2}}-1)^{i-1}.
			\Ee
			From $(1+x)^n= \sum_{\ell=0}^n \binom{n}{\ell} x^\ell$, we can bound the above term by 
			\Be\notag
			\begin{split}
				&2\tilde{T} e^{2 \tilde{T} \sqrt{\delta^\prime}} \sum_{\ell_1=1}^k e^{- \frac{\lambda^\prime \tilde{T}}{2} (k-\ell_1)}
				e^{ \frac{\lambda^\prime \tilde{T}}{2} ({k- (\ell_1-1)})}   
				\leq 2\tilde{T} e^{2 \tilde{T} \sqrt{\delta^\prime}} \sum_{\ell_1=1}^k e^{\frac{\lambda^\prime \tilde{T}}{2}}\\
				\leq &\   2e^{2 \tilde{T} \sqrt{\delta^\prime}} e^{\frac{\lambda^\prime \tilde{T}}{2}} \times  k\tilde{T} .
			\end{split}
			\Ee
			This proves the first estimate of (\ref{claim_AB}).

			Now we estimate $\tilde{B}(k\tilde{T}; N\tilde{T})$. From (\ref{sinh_gs}) and (\ref{hyp_sin}),
			\Be \notag
			\begin{split}
				\cosh \Big( \tilde{T}\sqrt{\delta^{\prime}} \sum_{j=1}^{\ell_{1}} e^{ -\frac{\lambda^{\prime}\tilde{T}}{2}(N-j) } \Big) 
				&\leq \exp \left(2 \tilde{T} \sqrt{\delta^\prime} e^{- \frac{\lambda^\prime \tilde{T}}{2} (N- \ell_1)} \right) ,  \\
				\sinh \Big( \tilde{T}\sqrt{\delta^{\prime}} \sum_{j=1}^{\ell_{i}} e^{ -\frac{\lambda^{\prime}\tilde{T}}{2}(N-(\ell_{1}+\cdots+\ell_{i-1}+j)) } \Big) 
				&\leq 2 \tilde{T}\sqrt{\delta^{\prime}} e^{ -\frac{\lambda^{\prime}\tilde{T}}{2}(N- (\ell_{1}+\cdots+\ell_{i})) } \\
				&= 2 \tilde{T}\sqrt{\delta^{\prime}} e^{ -\frac{\lambda^{\prime}\tilde{T}}{2}(N- k) } ,
			\end{split}
			\Ee
			for $\tilde{T} \gg 1$ and $\ell_1+ \cdots + \ell_i=k$. Therefore from (\ref{tilde_B_k})
			\Be
			\begin{split}\notag
				&\tilde{B}(k\tilde{T};N\tilde{T}) \\
				\leq& \ 2e^{2 \tilde{T} \sqrt{\delta^\prime}} e^{\frac{\lambda^\prime \tilde{T}}{2}} \left\{
				1 +  
				\sum_{\ell_{1}=1}^{k-1}
				\sum_{i=2}^{k-(\ell_{1}-1)} 
				\binom{k-\ell_1-1}{i-1}
				(e^{\frac{\lambda^{\prime}\tilde{T}}{2}} - 1)^{i-1} 
				e^{-\frac{\lambda^{\prime}\tilde{T}}{2}(N-\ell_{1})} 
				e^{-\frac{\lambda^{\prime}\tilde{T}}{2}(N-k)}\right\} \\  
				\leq& \ 2e^{2 \tilde{T} \sqrt{\delta^\prime}} e^{\frac{\lambda^\prime \tilde{T}}{2}} \left\{
				1 +  
				\sum_{\ell_{1}=1}^{k-1}
				e^{\frac{\lambda^{\prime}\tilde{T}}{2}(k-\ell_{1}-1)}   
				e^{-\frac{\lambda^{\prime}\tilde{T}}{2}(N-\ell_{1})} 
				e^{-\frac{\lambda^{\prime}\tilde{T}}{2}(N-k)} \right\}\\  
				\leq& \ 2e^{2 \tilde{T} \sqrt{\delta^\prime}} e^{\frac{\lambda^\prime \tilde{T}}{2}}  \left\{
				1 +  
				(k-1) e^{-\frac{\lambda^{\prime}\tilde{T}}{2}(2N-2k+1)}\right\}\\
				\leq& \  4e^{2 \tilde{T} \sqrt{\delta^\prime}} e^{\frac{\lambda^\prime \tilde{T}}{2}} .  
			\end{split}
			\Ee

			\hide
			\vspace{4pt}
			\textit{Step 7. } Now we claim (\ref{2_diff_bound}). Since all the entries in (\ref{2_diff_ineq}) are non-negative, we suffice to assume $A(s;t)$ and $B(s;t)$ in (\ref{2_diff_ineq}) are monotone increasing for maximal case. Moreover, $(2,1)$ entry in (\ref{auton}) overestimates $\delta^{\prime}e^{-\lambda^{\prime}s}$ in (\ref{2_diff_ineq}). Therefore $\tilde{A}$ and $\tilde{B}$ overestimate $A$ and $B$. Considering the following step function
			$
			\sum_{i=1}^{N} \tilde{A}(k\tilde{T}; N\tilde{T}) \mathbf{1}_{(k-1)\tilde{T} < s \leq k\tilde{T}}(s),
			$ we obtain linear growth when $(k-1)\tilde{T} < \tilde{s} := t-s \leq k\tilde{T}$,
			\Be
			\begin{split}
				A(s;t) &\leq \sum_{j=1}^{N} \tilde{A}(j\tilde{T}; N\tilde{T}) \mathbf{1}_{(j-1)\tilde{T} < t-s \leq j\tilde{T}}(s)  \\
				&\leq C_{\tilde{T},\delta^{\prime}} k \ \lesssim_{\tilde{T}, \delta^{\prime}} \ (t-s).
			\end{split}
			\Ee
			Uniform bound of $B(s,t)$ is also obtained by 
			
			\Be
			\begin{split}
				B(s;t) &\leq \sum_{j=1}^{N} \tilde{B}(j\tilde{T}; N\tilde{T}) \mathbf{1}_{(j-1)\tilde{T} < \tilde{s} \leq j\tilde{T}}(\tilde{s})  \leq C_{\tilde{T}, \delta^{\prime}}.
			\end{split}
			\Ee

			Now we claim (\ref{2_diff_bound}). Since all the entries in (\ref{2_diff_ineq}) is non-negative, we suffice to assume $A(s;t)$ and $B(s;t)$ are monotone increasing for maximal case.. Therefore, by considering the following step function
			$
			\sum_{k=1}^{N} A(k\tilde{T}; N\tilde{T}) \mathbf{1}_{(k-1)\tilde{T} < s \leq k\tilde{T}}(s),
			$, we obtain linear growth
			\Be
			\begin{split}\notag
				A(\tilde{s};t) &\leq \sum_{k=1}^{N} A(k\tilde{T}; N\tilde{T}) \mathbf{1}_{(k-1)\tilde{T} < \tilde{s} \leq k\tilde{T}}(\tilde{s})  \\
				&\leq C_{\tilde{T},\delta^{\prime}} \tilde{s} =  C_{\tilde{T},\delta^{\prime}} (t-s).
			\end{split}
			\Ee
			Uniform bound of $B(\tilde{s},t)$ is also obtained by 
			\Be
			\begin{split}\notag
				B(\tilde{s};t) &\leq \sum_{k=1}^{N} B(k\tilde{T}; N\tilde{T}) \mathbf{1}_{(k-1)\tilde{T} < \tilde{s} \leq k\tilde{T}}(\tilde{s})  \leq C_{\tilde{T}, \delta^{\prime}}.
			\end{split}
			\Ee\unhide

			\vspace{4pt}
			
			\textit{Step 6. }\textit{Comparison principle}:  For a non-negative matrix $M(\tilde{s};t)   \in \mathbb{R}^{n \times n}$,
			let us assume that $H(\tilde{s};t) \in \mathbb{R}^n$ and $J(\tilde{s};t)\in \mathbb{R}^n$ satisfy 
			\Be \label{eqnts_H_J}
			\frac{dH(\tilde{s};t)}{d\tilde{s}}  =  M (\tilde{s};t) H(\tilde{s};t),\ \
			\frac{dJ(\tilde{s};t)}{d\tilde{s}}  \leq  M (\tilde{s};t) J(\tilde{s};t) \ \ \text{for all }  \tilde{s}  \in [0,t].
			\Ee 
			If $H(0;t) \geq J(0;t) \geq 0$ and bounded then 
			\Be\label{HJ}
			H(\tilde{s};t) \geq J(\tilde{s};t) \geq 0 \text{ for all }  \tilde{s}  \in [0,t].
			\Ee
			
			Let us introduce $H^\e$ solving the same equation 
			with an $\e$-perturbed initial datum $H^\e (0;t) = H(0;t)+ \e (1, \cdots, 1)^T$. By the Gronwall's inequality, $|H^\e (\tilde{s};t) - H(\tilde{s};t)|\lesssim_n \e \exp\big(\int_0^{\tilde{s}} |M(\tilde{\tau};t)|   \dd \tau\big)$. Therefore for fixed $t$ and $\tilde{s}$, $H^\e(\tilde{s};t) \rightarrow H(\tilde{s};t)$ as $\e \downarrow 0$. Moreover from a temporal integration over the difference of equations in (\ref{eqnts_H_J}),
			\Be \begin{split}\label{H-J}
				H^\e(\tilde{s};t) - J(\tilde{s};t)   &\geq     H (0;t) - J(0;t)  + \e    (1, \cdots, 1)^T \\
				& +  \int_0^{\tilde{s}} M(\tau ;t)[H^\e(\tau ;t) - J(\tau ;t)]  \dd \tau 
				.\end{split}
			\Ee 
			
			We prove the estimate (\ref{HJ}) with exchanging $H(\tilde{s};t)$ by $H^\e (\tilde{s};t)$ first. Assume the statement is false. Then there exists $\tilde{s}_* \in [0,t)$ such that $H^\e(\tilde{s};t)\geq J(\tilde{s};t)$ for $\tilde{s} \in [0,\tilde{s}_*]$ but for some $i=1,2,\cdots, n$ we have $H^\e_i(\tilde{\tau}; t)< J_i(\tilde{\tau}; t)$ for $\tilde{\tau}> \tilde{s}_*$ with $|\tilde{\tau} - \tilde{s}_*| \ll 1$. 
			
			On the other hand, from (\ref{H-J}) and other conditions, for all $i$,  
			\Bes
			H ^\e(\tilde{\tau};t) - J (\tilde{\tau};t) &\geq& \e + \int^{\tilde{\tau}}_{\tilde{s}_*} M(\tau;t) [H ^\e(\tau;t) - J (\tau;t)] \dd \tau\\
			&\geq& \e - O_{M,H,J}(| {\tilde{\tau}} -{\tilde{s}_*}  |).
			\Ees
			If $| {\tilde{\tau}} -{\tilde{s}_*}  |$ is small enough then $H^\e (\tilde{\tau};t) > J(\tilde{\tau};t)$. This is a contradiction. By passing a limit $\e \downarrow 0$ we prove (\ref{HJ}). 
			
			Now we apply above comparison principle to (\ref{2_diff_ineq}) and (\ref{auton}) as
			\Be
			\begin{split}\notag
				M(\tilde{s};t) = 
				\begin{bmatrix}
					0 & 1 \\
					\delta^\prime e^{- \lambda^\prime (N-k )\tilde{T} } & 0
				\end{bmatrix}
				,\quad 
				J(\tilde{s};t) = 
				\begin{bmatrix}
					A(s;t) \\  B(s;t)
				\end{bmatrix} 
				,\quad	
				H(\tilde{s};t) = 
				\begin{bmatrix}
					\tilde{A}(\tilde{s}; N\tilde{T}) \\  \tilde{B}(\tilde{s}; N\tilde{T})
				\end{bmatrix} .
			\end{split}
			\Ee
			We use (\ref{HJ}) and (\ref{claim_AB}) to prove (\ref{2_diff_bound}),
			\Be\notag
			\begin{bmatrix}
				A(s;t) \\  B(s;t)
			\end{bmatrix} 
			\leq
			\begin{bmatrix}
				\tilde{A}(\tilde{s}; N\tilde{T}) \\  \tilde{B}(\tilde{s}; N\tilde{T})
			\end{bmatrix}
			\leq 
			Ce^{C \lambda^\prime \sqrt{\delta^\prime}}
			\begin{bmatrix}
				k \\ 1
			\end{bmatrix}
			\leq Ce^{C \lambda^\prime \sqrt{\delta^\prime}}
			\begin{bmatrix}
				t-s \\ 1
			\end{bmatrix},
			\Ee
			where $C$ depends on $\tilde{T}$ which is large but fixed. 
		\end{proof}
		
		\unhide
		
		Now we are ready to prove Lemma \ref{COV_boundary}.
		\begin{proof}[\textbf{Proof of Lemma \ref{COV_boundary}}] The lower bound (\ref{tb_lower}) is a direct consequence of the identity
			$$\xb(t,x,u) = x + \int^{t-\tb(t,x,u)}_t V(s;t,x,u) \dd s.$$ 
			
			\hide\vspace{4pt}
			
			\textit{Step 1.} We claim that for all $(t,x) \in [ 0, \infty) \times \bar{\O}$ as $N \rightarrow \infty$
			\Be\label{conv_NN}
			\mathbf{1}_{\tb(t,x,u)< N}
			\mathbf{1}_{n(\xb(t,x,u)) \cdot \vb(t,x,u) < -\frac{1}{N}}
			\nearrow \mathbf{1}_{\tb(t,x,u)< \infty} \ \ \text{almost every } u \in \R^3.
			\Ee
			First we prove that, for fixed $N \in \mathbb{N}$, as $M \rightarrow \infty$
			\Be\label{conv_NM}
			\mathbf{1}_{\tb(t,x,u)< N}
			\mathbf{1}_{n(\xb(t,x,u)) \cdot \vb(t,x,u) < -\frac{1}{M}}
			\nearrow \mathbf{1}_{\tb(t,x,u)< N} \ \ \text{almost every } u \in \R^3.
			\Ee
			Since $\mathbf{1}_{\tb(t,x,u)< N}$ converges to $\mathbf{1}_{\tb(t,x,u)< \infty}$ as $N \rightarrow \infty$ we can apply Cantor's diagonal argument to conclude (\ref{conv_NN}) from (\ref{conv_NM}).
			
			Assume that $\tb(t,x,u)<N$ and $u \neq 0$. If $x \in \p\O$ then we further assume that $n(x) \cdot u > 0$. If $x \in \p\O$ and $n(x) \cdot u<0$ then $(\xb(t,x,u), \vb(t,x,u)) = (x,u)$ and hence $n(\xb(t,x,u)) \cdot \vb(t,x,u)= n(x) \cdot u <0$, which implies (\ref{conv_NM}). 
			
			We show that if $\tb(t,x,u)<N, u \neq 0$, and $n(x) \cdot u > 0$ if $x \in \p\O$, then we have
			\Be\label{no_graze}
			n(\xb(t,x,u)) \cdot \vb(t,x,u)<0,
			\Ee
			so that (\ref{conv_NM}) can be proven. Note that locally we can parametrize the trajectory (see Lemma 15 in \cite{GKTT1} or \cite{KL2} for details). We consider local parametrization (\ref{eta}).  (We drop the subscript $p$ for the sake of simplicity.) We define $(X_n, X_\parallel)$ to satisfy 
			\Be\label{X_local}
			X(s;t,x,u)  =   \eta (X_\parallel,0) + X_n [- n(X_\parallel)].
			\Ee
			
			\hide $\eta$ near $\p\O$ such that 
			\Be
			\eta : (x,y,z)\in B(0,r_{1})\cap \R^{3}_{+} \mapsto B(p,r_{2})\cap\O,
			\Ee
			where $\eta(0)=p\in\p\O$ and $\eta(x,y,0) \in \p\O$ for some $r_{1}, r_{2} > 0$. Then for $X(s;t,x,v)$, there exist a unique $x_{*}\in B(p,r_{2})\cap \p\O$ such that
			\Be \label{X def}
			|X(s;t,x,v)-x_{*}| \leq \sup_{x\in B(p,r_{2})\cap \p\O} |X(s;t,x,v)-x|,
			\Ee 
			since $\eta$ is bijective. We define 
			\[
			(X_{\parallel},0) = \eta^{-1}(x_{*}) \quad \text{and} \quad  X_{n} := |X(s;t,x,v)-x_{*}|.
			\]\unhide
			For velocity 
			we define
			\[
			V_{n} := V(s;t,x,v)\cdot (-n(X_{\parallel})).
			\]
			To define $V_{\parallel}$, we consider level set $  \big( \eta(X_{\parallel},0) + X_{n}(-n(X_{\parallel})) \big)$ for fixed $X_{n}$. Then we have  
			\[
			\p_{i} \big( \eta(X_{\parallel},0) + X_{n}(-n(X_{\parallel})) \big) \perp n(X_{\parallel}),\quad i=1,2, 
			\]
			and above two independent vectors span tangential plane at $X(s;t,x,v)$. We define $(V_{\parallel,1}, V_{\parallel,2})$ as each coefficient, i.e.
			\Be\notag
			\begin{split} 
				&\sum_{i=1,2} V_{\parallel,i} \p_{i} \big( \eta(X_{\parallel},0) + X_{n}(-n(X_{\parallel})) \big)  \\
				&= \nabla_{\parallel} \big( \eta(X_{\parallel},0) + X_{n}(-n(X_{\parallel})) \big) V_{\parallel} = V - V_{n}(-n(X_{\parallel}))  .
			\end{split}
			\Ee
			Therefore we obtain
			\Be \label{V_local}
			V(s;t,x,u)   = V_n [- n(X_\parallel)] + V_\parallel \cdot \nabla_{x_\parallel} \eta (X_\parallel,0) 
			- X_n V_\parallel \cdot \nabla_{x_\parallel} n (X_\parallel).
			\Ee  
			For the proof of (\ref{no_graze}) we use the contradiction argument: let us assume
			\Be\label{initial_00}
			X_n (t-\tb;t,x,u)  +V_n (t-\tb;t,x,u) =0.
			\Ee
			
			First we choose $0<\e \ll 1$ such that $X_n(s;t,x,u) \ll 1$ and 
			\Be\label{Vn_positive}
			V_n (s;t,x,u) \geq0 \ \ \text{for} \  t- \tb(t,x,u)<s<t-\tb(t,x,u) + \e.
			\Ee
			The sole case that we cannot choose such $\e>0$ is when $V_n(s;t,x,u)<0$ for $s \in ( t-\tb(t,x,u), t-\tb(t,x,u) + \delta)$ with $0< \delta\ll1$. But this cannot be true: from (\ref{hamilton_ODE}) for $s \in ( t-\tb(t,x,u), t-\tb(t,x,u) + \delta)$
			$$
			0 \leq X_n(s;t,x,u)   =  X_n(t-\tb(t,x,u);t,x,u)  +  \int^s_{t-\tb(t,x,u)} V_n (\tau; t,x,u) \dd \tau <  0.$$
			
			Temporarily we define that $t_* := t-\tb(t,x,u) + \e$, $x_* = X(t-\tb(t,x,u) + \e; t,x,u),$ and $u_* = V(t-\tb(t,x,u) + \e; t,x,u)$. 
			
			From (\ref{hamilton_ODE}) , for $(X_n(s), X_\parallel (s)) = (X_n(s; t_*, x_*, u_*), X_\parallel (s; t_*, x_*, u_*))$ and $(V_n(s), V_\parallel (s)) = (V_n(s; t_*, x_*, u_*), V_\parallel (s; t_*, x_*, u_*))$,  
			\Be \begin{split}
				\dot{X}(s;t,x,u) &= \nabla_{\parallel}\eta (X_\parallel)\dot{X}_{\parallel} + \dot{X}_n [- n(X_\parallel)] - X_{n}\nabla_{\parallel} n(X_{\parallel})\dot{X}_{\parallel}  \\
				&= V_n [- n(X_\parallel)] + \nabla_{x_\parallel} \eta (X_\parallel)V_\parallel 
				- X_n \nabla_{x_\parallel} n (X_\parallel)  V_\parallel .
			\end{split}\Ee
			Comparing coefficients of normal and tangential components, $\dot{X}_{n} = V_{n}$ and $\dot{X}_{\parallel} = V_{\parallel}$. From (\ref{V_local}),
			\Be \begin{split} \label{Vdotn}
				\dot{V} (s) &= -\dot{V}_{n} n(X_{\parallel}) - V_{n}\nabla_{\parallel} n(X_{\parallel})\dot{X}_{\parallel} + V_{\parallel}\cdot\nabla^{2}_{\parallel}\eta(X_{\parallel})\dot{X}_{\parallel} + \nabla_{\parallel}\eta(X_{\parallel})\dot{V}_{\parallel}  \\
				&\quad - \dot{X}_{n}\nabla_{\parallel} n(X_{\parallel})V_{\parallel} - X_{n}\nabla_{\parallel} n(X_{\parallel})\dot{V}_{\parallel} - X_{n} V_{\parallel}\cdot\nabla_{\parallel}^{2}n(X_{\parallel})\dot{X}_{\parallel}. 
			\end{split}\Ee
			From $(\ref{Vdotn})\cdot n(X_{\parallel})$ and (\ref{hamilton_ODE}) we obtain
			\Be \begin{split}\label{hamilton_ODE_perp}
				\dot{X}_n (s) &= V_{n},  \\
				\dot{V}_n (s) 
				&=  [V_\parallel (s)\cdot \nabla^2 \eta (X_\parallel(s)) \cdot V_\parallel(s) ] \cdot n(X_\parallel(s))  
				+   \nabla\phi (s , X (s ) ) \cdot n(X_\parallel(s))  \\
				&\quad - X_n (s) [V_\parallel(s) \cdot \nabla^2 n (X_\parallel(s)) \cdot V_\parallel(s)]  \cdot n(X_\parallel(s)) .
			\end{split}\Ee
			
			\hide
			From (\ref{hamilton_ODE}), for $(X_n(s), X_\parallel (s)) = (X_n(s; t_*, x_*, u_*), X_\parallel (s; t_*, x_*, u_*))$ and $(V_n(s), V_\parallel (s)) = (V_n(s; t_*, x_*, u_*), V_\parallel (s; t_*, x_*, u_*))$, 
			\Be \begin{split}\label{hamilton_ODE_perp}
				\dot{X}_n (s)  =& V_n (s),\\
				\dot{V}_n (s)  =& 
				[V_\parallel (s)\cdot \nabla^2 \eta (X_\parallel(s)) \cdot V_\parallel(s) ] \cdot n(X_\parallel(s))  
				+   \nabla\phi (s , X (s ) ) \cdot n(X_\parallel(s))  \\
				&+   X_n (s) [V_\parallel(s) \cdot \nabla^2 n (X_\parallel(s)) \cdot V_\parallel(s)]  \cdot n(X_\parallel(s)) .
			\end{split}\Ee\unhide
			From (\ref{convexity_eta}), $[V_\parallel(s) \cdot \nabla^2 \eta (X_\parallel(s)) \cdot V_\parallel(s) ] \cdot n(X_\parallel(s))\leq 0$. By an expansion and (\ref{phi_BC}) we deduce that
			\Bes
			&&\nabla\phi (s , X(s )) \cdot n(X_\parallel(s ) )\\
			& =&  \nabla\phi (s , X_n(s ), X_\parallel(s ) ) \cdot n(X_\parallel (s )) \\
			&=& \nabla\phi (s , 0, X_\parallel(s ) ) \cdot n(X_\parallel (s )) 
			+ \| \phi (s) \|_{C^2}  O(|X_n(s )|)\\
			&=&  \| \phi (s) \|_{C^2}  O(|X_n(s )|).
			\Ees
			Again from (\ref{hamilton_ODE}), $$|V_\parallel (s )|\leq |u| + \tb(t,x,u) \sup_{t- \tb \leq s \leq t} |\nabla \phi(s,X(s;t,x,u)) |  \leq |u| +  N \| \phi \|_{C^1}.$$ 
			Therefore we conclude that 
			\Be\notag
			\dot{V}_n(s;t_*,x_*,u_*) \leq C_{\eta, N , \phi }(1+ |u|^2  )X_n(s;t_*,x_*,u_*) ,
			\Ee
			and hence 
			\Be\label{ODE_X+V}
			\frac{d}{ds} [X_n (s )  +V_n (s ) ]
			\lesssim (1+ |u|^2)  [X_n (s )  +V_n (s ) ].
			\Ee
			By Gronwall's inequality and (\ref{initial_00})
			\Be \begin{split}\notag
				& [X_n (s)  +V_n (s) ]  \\
				\lesssim & \   [X_n (t-\tb(t,x,u))  +V_n (t-\tb(t,x,u)) ] e^{C (1+ |u|^2) |(t-\tb (t,x,u) )- s|}\\
				=&  \ 0.
			\end{split}\Ee 
			From (\ref{Vn_positive}) we conclude that $X_n (s) \equiv 0$ and $V_n (s) \equiv 0$ for $s \in [t-\tb(t,x,u), t-\tb(t,x,u) + \e]$. We can continue this argument successively to deduce that $X_n (s) \equiv 0$ and $V_n (s) \equiv 0$ for all $s \in [t-\tb(t,x,u), t]$. Then $x \in \p\O$ and $n(x) \cdot u =0$. This is a contradiction since we chose $n(x) \cdot u>0$.


			\unhide
			
			
			Now we prove the inequality of (\ref{COV_xbtb}). Except measure zero set of $u$ (when $x \in \p\O$ we only need to exclude $u \cdot n(x)=0$) we have $n(\xb(t,x,u)) \cdot \vb(t,x,u)\neq 0$ from Lemma \ref{cannot_graze}. \hide

			Note that we only need to prove that $$\int_{\R^3} \mathbf{1}_{\tb(t,x,u)< N}
			\mathbf{1}_{n(\xb(t,x,u)) \cdot \vb(t,x,u) < -\frac{1}{N}}
			g(x,u) \dd u$$ is bounded by the RHS of (\ref{COV_xbtb}) uniformly for all $N>0$. Once it is done we can pass $N \rightarrow \infty$ and then apply the monotone convergence theorem with (\ref{conv_NN}) to conclude (\ref{COV_xbtb}).

			Choose $(x,u) \in \bar{\O} \times \R^3$ such that $\tb(t,x,u)< N$ and $n(\xb(t,x,u)) \cdot \vb(t,x,u) < -\frac{1}{N}$.\unhide Then by the implicit function theorem $\tb$ and $\xb$ are differentiable locally.

			Let us choose $p \in \p\O$ which is close to $\xb(t,x,u)$. Recall (\ref{eta}) and the notations $v_{\parallel,i}$ and $v_n$ in (\ref{def_V_n}) and (\ref{def_V_parallel}). Here, we temporarily define three-dimensional local parametrization $\eta_{p}$ near $\xb(t,x,u)$ so that it is consistent with boundary parametrization (\ref{eta}). We use $x_{\parallel,i}$ and $x_{n}$ for (\ref{X_local}). Let us pick $p\in\p\O$ very close to $\xb(t,x,u)\in\p\O$ and we define
			\begin{equation}\label{etain}
			\begin{split}
			{\eta}_{{p}}:  \{ (x_{\parallel,1}, x_{\parallel,2}, x_{n}) \in \mathbb{R}^{3}:  x_{n} > 0  \} &\cap B(0; \delta^{\prime}_{1}) \ \rightarrow  \ \Omega \cap B({p}; \delta^{\prime}_{2}),\\
			(x_{\parallel,1}, x_{\parallel,2}, x_{n})	 \ &\mapsto \   {\eta}_{p} (x_{\parallel,1}, x_{\parallel,2}, x_{n}),  \\
			{\eta}_{p}(x_{\parallel,1},x_{\parallel,2},x_{n}) =  {\eta}_{p}(x_{\parallel,1},x_{\parallel,2},0) &+ x_{n}\big( -n(x_{\parallel,1},x_{\parallel,2}) \big),
			\end{split}
			\end{equation}
			for sufficiently small $\delta^{\prime}_{1}, \delta^{\prime}_{2} \ll 1$ and we assume $|p-\xb(t,x,u)| < \delta^{\prime}_{2}$. We used notation $n(x_{\parallel,1},x_{\parallel,2})$ to denote $n(\eta_{p}(x_{\parallel,1},x_{\parallel,2},0))$. From the definition of (\ref{etain}), $\p_{i}\eta_{p} := \frac{\p}{\p x_{\parallel,i}} \eta_{p}$ belongs to tangential plane of the level set for fixed $x_{n} \geq 0$ and hence 
			\[
			(\p_{1}\eta_{p}\times \p_{2}\eta_{p}) \parallel n(x_{\parallel,1}, x_{\parallel,2}).
			\]
			To specify local parametrization of $\xb(t,x,u)$, let us denote
			\[
			\xb(t,x,u) = \eta_p ( x_{\mathbf{b},1} , x_{\mathbf{b},2}, 0  )  .
			\]  
			
			We compute that, with the standard notation
			$g_{p,ij} :=  \p_{i} \eta_{p}\cdot\p_{j} \eta_{p} $,
			\Be\label{dxdtdv}
			\begin{split}
				&\det \bigg(\frac{\p ( x_{\mathbf{b},1}, x_{\mathbf{b},2},\tb )}{\p u } \bigg) \\
				=  &
				\ 
				\det 
				\begin{bmatrix} \tb \text{Id}_{3\times3}  - \int^{t-\tb}_{t}
					\int^{s}_{t} \nabla_{u} \big(E (\tau ,X(\tau ;t,x,u))\big) \dd \tau 
					\dd s
				\end{bmatrix}  
				\\
				&    \times 
				\det \begin{bmatrix}
					\frac{-1}{\sqrt{g_{p,11} (\xb)}}
					\Big[ \frac{\p_{1} \eta_{p}(\xb)}{\sqrt{g_{p,11}(\xb)} } +  \frac{v_{\parallel,1}(\xb) }{v_{n}(\xb)} \frac{\p_{3} \eta_{p}(\xb) }{\sqrt{g_{p,33}(\xb)}}\Big]\\
					\frac{-1}{\sqrt{g_{p,22} (\xb)}}
					\Big[ \frac{\p_{2} \eta_{p}(\xb)}{\sqrt{g_{p,22}(\xb)} } +  \frac{v_{\parallel,2}(\xb) }{v_{n}(\xb)} \frac{\p_{3} \eta_{p}(\xb) }{\sqrt{g_{p,33}(\xb)}}\Big]\\
					\frac{1}{v_{n} (\xb)} \frac{\p_{3} \eta_{p}(\xb) }{\sqrt{g_{p,33}(\xb)}}
				\end{bmatrix},
			\end{split} \Ee
			%
			where we have used direct computations (See \cite{GKTT1,KL1} for details)
			\Be\begin{split}
			\frac{\p x_{\mathbf{b},i}}{\p u_{j}} =&
			- \Big[ \tb e_{j} 
			- \int^{t-\tb}_{t} \int^{s}_{t} 
			\p_{u_j}\big(E(\tau, X(\tau;t,x,u))\big)
			\dd \tau \dd  {s} \Big]
			\notag
			\\
			&   \ \ \ \ \cdot \frac{1}{\sqrt{g_{p,ii} (\xb)}}
			\Big[ \frac{\p_{i} \eta_{p}(\xb)}{\sqrt{g_{p,ii}(\xb)} } +  \frac{v_{\parallel,i}(\xb) }{v_{n}(\xb)} \frac{\p_{3} \eta_{p}(\xb) }{\sqrt{g_{p,33}(\xb)}}\Big],  
			\\
			\nabla_{u} \tb   =& 
			\frac{1}{v_{n} (\xb)} \Big[ \tb I_{3 \times 3}  - \int^{t-\tb}_{t}  \int^{s}_{t}
			\nabla_{u} 
			\big(
			E (\tau, X(\tau;t,x,u))  \big)
			\dd \tau \dd  {s}
			\Big] \\
			&\ \ \     \ \ \ \ \ \ \ \cdot   \frac{\p_{3} \eta_{p}}{\sqrt{g_{p,33}}}
			.
			\end{split}\Ee

			
			Note that we have (\ref{result_X_v}) from  (\ref{decay_phi_2}) and Lemma \ref{est_X_v}. \hide
			
			\Be\label{linear_Xu}
			|\nabla_u X(\tau;t,x,u)| \leq  C e^{C \Lambda_{2} \sqrt{\delta_{2}}} |t-\tau|.
			\Ee\unhide
			Then, from (\ref{result_X_v}) and (\ref{decay_phi_2}),
			\hide
			\Be\begin{split}\label{int_phi_uu}
				&\Big|\int^{t-\tb}_{t} \int^{s}_{t} 
				\nabla_ u \big(E(\tau, X(\tau;t,x,u))\big)
				\dd \tau \dd  {s}\Big| \\
				\leq 	&  \ 	
				\int^t_{t-\tb} \int^t_s
				|\nabla_u X(\tau;t,x,u)| | \nabla_x E (\tau, X(\tau;t,x,u))|
				\dd \tau  \dd s 	.
			\end{split}\Ee
			
			From (\ref{linear_Xu}) and (\ref{decay_phi_2}) \unhide
			\Be\begin{split}\label{double_int_phi_1} 
				&\Big|\int^{t-\tb}_{t} \int^{s}_{t} 
				\nabla_ u \big(E(\tau, X(\tau;t,x,u))\big)
				\dd \tau \dd  {s}\Big| \\
				\leq 	&  \ 	
				\int^t_{t-\tb} \int^t_s
				|\nabla_u X(\tau;t,x,u)| | \nabla_x E (\tau, X(\tau;t,x,u))|
				\dd \tau  \dd s\\\leq  & \ C e^{C  \delta_{2}(\Lambda_{2} )^{-2} } \delta_{2} \int^t_{t-\tb} \int^t_s   |t-\tau | e^{- {\Lambda_{2} } \tau }
				\dd \tau  \dd s \\
				\leq & \  C e^{C \delta_{2}(\Lambda_{2} )^{-2}
				} \delta_{2} \tb  \int^t_{t-\tb} \int^t_s  e^{- {\Lambda_{2} } \tau }
				\dd \tau 
				\dd s\\
				\leq & \ Ce^{C \delta_{2}(\Lambda_{2} )^{-2}} 
				\delta_{2} \tb \times \frac{1}{(\Lambda_2)^2}   ,
			\end{split}\Ee      
			where we have used a direct computation  
			\Be\label{direct_computation}
			\begin{split}
		   \int^t_{t-\tb} \int^t_s  e^{- {\Lambda_{2} } \tau }
				\dd \tau 
				\dd s 
				= \ &\int^{t}_{t-\tb}  
				\frac{e^{- {\Lambda_{2} }  t } - e^{- \Lambda_{2} s}}{- \Lambda_{2} } \dd s\\
				= \  & 
				\frac{e^{- {\Lambda_{2} }  t }}{-  \Lambda_{2}  }\tb - \frac{e^{- {\Lambda_{2} }  t } - e^{- \Lambda_{2}  (t-\tb)}}{\left( {\Lambda_{2} } \right)^2} \\
				= \  & \frac{e^{-\Lambda_2(t-\tb)} }{(-\Lambda_2)^2}
				\Big\{ 
				1- e^{-  {\Lambda_{2} }\tb } (1+  {\Lambda_{2} } \tb ) 
				\Big\}\\
				\leq   \ & \frac{1}{(\Lambda_2)^2}.
			\end{split}
			\Ee
			
			\hide
			
			From (\ref{condition_lambda_infty}) and (\ref{X_v})-(\ref{double_int_phi})
			\Be\label{lower_jacob_l}
			\det\left(\frac{\p X(s^\prime;t^\prime_{l^\prime}, x^\prime_{l^\prime},v^\prime_{l^\prime})}{\p v^\prime_{l^\prime}}\right)\geq \frac{|t^\prime_{l^\prime}-s^\prime|^3}{2}
			.             \Ee

			Then 
			\Bes
			&& \Big| \int^{t-\tb}_{t}
			\int^{s}_{t} \nabla_{u  } \big(\nabla \phi (\tau ,X(\tau ;t,x,u))\big) \dd \tau 
			\dd s\Big|\\
			&=& \Big|\int^{t-\tb}_{t}
			\int^{s}_{t} (\nabla_{u} X(\tau;t,x,u)  \cdot \nabla_x) \nabla_x \phi (\tau ,X(\tau ;t,x,u)) \dd \tau 
			\dd s\Big|\\
			&\leq&  
			\int^{t-\tb}_{t}
			\int^{s}_{t} C_{K, \lambda} |t-\tau|
			K e^{-\lambda_0 \tau }  \dd \tau \dd s\\
			&\leq&  K  C_{K, \lambda}\tb  \times    
			\int_{t-\tb}^t
			\frac{e^{- \lambda_0 s} - e^{- \lambda t}}{\lambda }
			\dd s 
			\\
			&\leq&  
			K  C_{K, \lambda}\tb  \times    
			\frac{e^{-\lambda_0 t}}{(\lambda )^2}
			\Big[
			e^{\lambda_0 \tb} -(1 + \lambda  \tb)
			\Big]\\
			&\leq& \tb \times \frac{K C_{K, \lambda} }{(\lambda )^2} e^{-\lambda  
				(t-\tb)}.
			\Ees 
			%

			\unhide

			Therefore we derive a lower bound of the first determinant of (\ref{dxdtdv}) for $\delta_{2}(\Lambda_{2} )^{-2}   \ll 1$, which is guaranteed by 
			\Be \label{first_determinant}
			(\tb)^3 \det \Big[ I_{3\times 3} +  \frac{C\delta_{2}}{(\Lambda_2)^2} e^{C \delta_{2}(\Lambda_{2} )^{-2}} \Big] 
			\geq  \frac{(\tb)^3}{2},
			\Ee 
			where we have used $\det [I_{n \times n} + O(\delta)] =1 + O_n(\delta)$ for $0< \delta \ll1$.
			
			Clearly the second determinant of (\ref{dxdtdv}) equals
			\Be\label{second_determinant}
			\begin{split}
				&\bigg|\frac{1}{v_{n} (\xb)} \frac{\p_{3} \eta_{p}(\xb) }{\sqrt{g_{p,33}(\xb)}} \\
				&
				\cdot 
				\bigg(
				\frac{-1}{\sqrt{g_{p,11} (\xb)}}
				\Big[ \frac{\p_{1} \eta_{p}(\xb)}{\sqrt{g_{p,11}(\xb)} } +  \frac{v_{\parallel,1}(\xb) }{v_{n}(\xb)} \frac{\p_{3} \eta_{p}(\xb) }{\sqrt{g_{p,33}(\xb)}}\Big]\\
				&   \ \ \ 
				\times 
				\frac{-1}{\sqrt{g_{p,22} (\xb)}}
				\Big[ \frac{\p_{2} \eta_{p}(\xb)}{\sqrt{g_{p,22}(\xb)} } +  \frac{v_{\parallel,2}(\xb) }{v_{n}(\xb)} \frac{\p_{3} \eta_{p}(\xb) }{\sqrt{g_{p,33}(\xb)}}\Big]			
				\bigg)
				\bigg|\\
				&
				=  \frac{1}{ \sqrt{g_{p,11} (\xb)}\sqrt{g_{p,22} (\xb)}} \times \frac{1}{|v_{n}(\xb) |}.
			\end{split}
			\Ee

			Now we check whether the mapping $u \rightarrow (\xb(t,x,u), \tb (t,x,u) )$ is one-to-one for all $t \in \R, \tb(t,x,u)<N,$ and $n(\xb(t,x,u)) \cdot \vb(t,x,u)< - \frac{1}{N}$ with $x \in \O$ or $x \in \p\O$ with $n(x) \cdot u>0$. Assume that there exist $u$ and $\tilde{u}$ with such conditions and satisfy $\xb(t , x , u ) = \xb(t , x , \tilde{u})$ and $\tb(t , x , u ) = \tb(t , x , \tilde{u})$. As (\ref{dxdtdv}) we choose $p \in \p\O$ near $\xb(t , x , u )$ and use the same parametrization. Then by an expansion, for some $\bar{u}  \in \overline{\tilde{u}  u }$, 
			%
			\Be 
			0  = 
			\begin{bmatrix}
				x_{\mathbf{b},1} (t , x , \tilde{u} )\\
				x_{\mathbf{b},2} (t , x , \tilde{u} )\\
				\tb(t , x , \tilde{u} )
			\end{bmatrix} -
			\begin{bmatrix}
				x_{\mathbf{b},1} (t , x , {u} )\\
				x_{\mathbf{b},2} (t , x ,  {u} )\\
				\tb(t , x ,  {u} )
			\end{bmatrix} 
			=  \begin{bmatrix}
				\nabla_u x_{\mathbf{b},1} (t , x , \bar{u} )\\
				\nabla_u x_{\mathbf{b},2} (t , x , \bar{u} )\\
				\nabla_u \tb(t , x , \bar{u} )
			\end{bmatrix}
			(\tilde{u}  - u ).\label{u-tilde_u}
			\Ee 
			This equality can be true only if the determinant of the Jacobian matrix equals zero. And (\ref{first_determinant}) and (\ref{second_determinant}) imply that $\tb(t,x,\bar{u})=0$. But this implies $x \in \p\O$ and hence $\vb(t,x,u) =u$ and $\vb(t,x,\tilde{u}) = \tilde{u}$. Then $n(x) \cdot u<0$ and $n(x) \cdot \tilde{u} <0$ which are out of our domain. 
			%
			
			Now we apply the change of variables to conclude the proof of Lemma \ref{COV_boundary}.
		\end{proof}

		Finally we present the proof of the main result of this section.
		
		\begin{proof}[\textbf{Proof of Proposition \ref{prop_int_alpha}}]
			\textit{Step 1.} 
			%
			%
			%
			\hide
			
			\textit{Step 2.} Choose $N>0$. It is clear that if $\tb(s,x,u)<N$ then $\xb(s,x,u)$ and $\vb(s,x,u)$ are well-defined and $n(\xb(s,x,u)) \cdot \vb(s,x,u) \leq  0$. We claim that if $\tb(s,x,u)<N$ then 
			\Be\label{no_grazing}
			n(\xb(s,x,u)) \cdot \vb(s,x,u) <  0 \ \ \text{almost everywhere in } \bar{\O} \times \R^3.
			\Ee
			Once (\ref{no_grazing}) holds, it is easy to conclude that
			\Be\begin{split}
				\mathbf{1}_{\tb(s,x,u)<N} \mathbf{1}_{n(\xb(s,x,u)) \cdot \vb(s,x,u) <-\frac{1}{N}}
				\rightarrow \mathbf{1}_{\tb(s,x,u)<\infty}\\ \ \ \text{for almost every } u \text{, as} \ 
				N \rightarrow \infty.\end{split}
			\Ee
			By the monotone convergence theorem we only need to prove a uniform bound in $N$ for 
			\Be\label{cutoff_alpha}
			\int_{\R^3} \mathbf{1}_{\tb(s,x,u)<N} \mathbf{1}_{n(\xb(s,x,u)) \cdot \vb(s,x,u) <-\frac{1}{N}}\frac{e^{-C|v-u|^2}}{\alpha(s,x,u)^\sigma} \dd u.
			\Ee
			
			\vspace{4pt}
			
			\unhide
			%
			We apply Lemma \ref{COV_boundary} and deduce that 
			\Be\begin{split}\label{apply_Lemma_COV}
				\int_{|u| \leq N} \frac{\mathbf{1}_{\tb(s,x,u)<t+1}  }{\alpha_{f,\e}(s,x,u)^\sigma} \dd u
				&\leq \int_{\p\O} \int  \frac{|n(\xb(t,x,u)) \cdot \vb(t,x,u)|^{1- \sigma}}{|\tb(t,x,u)|^3} 
				\dd \tb  \dd \xb.
			\end{split}\Ee 
			Here $\tb$ has a lower bound from (\ref{tb_lower})
			\Be\label{lower_tb_lessN}
			\tb  \geq \frac{|\xb  -x|}{N+ \delta_1/\Lambda_1},
			\Ee
			where we have used the fact, from $E$ in (\ref{decay_E}), for $|u|\leq N$ and $0\leq s \leq t$,
			\Be
			\begin{split}\notag
			|V(s;t,x,u)| \leq& \  |u| + \int^t_{0} |E (s,X(s;t,x,u))| \dd s\\ \leq & \ 
			 N + \int^t_{0} \delta_1e^{-\Lambda_1 s} \dd s\\  \leq & \  N+ \delta_1/\Lambda_1.
			\end{split}\Ee
			
			From $V(s;t,x,u) =  \vb(t,x,u) +\int^s_{t-\tb(t,x,u)}E(\tau, X(\tau;t,x,u)) \dd \tau$, we have 
			\Be\begin{split}
				&\xb(t,x,u) - x  \\
				=  & \   \int^{t-\tb(t,x,u)}_t
				\bigg\{
				\vb(t,x,u)+\int^s_{t-\tb(t,x,u)}E(\tau, X(\tau;t,x,u)) \dd \tau
				\bigg\}
				\dd s \\
				=  & \  - \tb(t,x,u) \vb(t,x,u)\\
				&  +  \int^{t-\tb(t,x,u)}_t 
				\int^s_{t-\tb(t,x,u)} 
				E(\tau, X(\tau;t,x,u)) 
				\dd \tau
				\dd s.\end{split}\notag
			\Ee 
			Therefore we can conclude that 
			\Be\begin{split}\notag
				&\vb(t,x,u)  \\
				= & \   \frac{x-\xb(t,x,u)}{\tb(t,x,u)}\\
				&
				- \frac{1}{\tb(t,x,u)} \int_{t-\tb(t,x,u)}^t 
				\int^s_{t-\tb(t,x,u)} 
				E(\tau, X(\tau;t,x,u)) 
				\dd \tau
				\dd s,
			\end{split}\Ee 
			and hence
			\Be 
			\begin{split}\label{n_bv_b}
				&|n(\xb(t,x,u))\cdot\vb(t,x,u)  |\\
				\leq  & \  \frac{\big|(x-\xb(t,x,u)) \cdot n(\xb(t,x,u))\big|}{\tb(t,x,u)}\\
				&  
				+ \frac{\tb(t,x,u)}{2} \max_{t-\tb(t,x,u) \leq \tau \leq t}
				|E(\tau, X(\tau;t,x,u)) 
				\cdot n(\xb(t,x,u))|.
			\end{split}
			\Ee
			Using (\ref{n_bv_b}), we further bound (\ref{apply_Lemma_COV}) by
			\begin{eqnarray}
			 (\ref{apply_Lemma_COV})
			& \leq& \int_{\p\O}  {|(x- \xb(t,x,u)) \cdot n(\xb(t,x,u))|^{1- \sigma}}
			\int \frac{1}{\tb^{4- \sigma}}
			\dd \tb
			\dd S_{\xb}\label{bound1}\\
			&&  + \| E \|_{L^\infty} \int_{\p\O} 
			\int
			\frac{1 }{\tb ^{2+ \sigma}} \dd \tb
			\dd S_{\xb},\label{bound2}
			\end{eqnarray}
			where $\tb$ has a lower bound of (\ref{lower_tb_lessN}).

			By integrating over $\tb$ in (\ref{bound1}) and using the lower bound (\ref{lower_tb_lessN}),
			we deduce that 
			\Be 
			(\ref{bound1})  
			 \lesssim     \int_{\p\O} {|(x- \xb ) \cdot n(\xb )|^{1- \sigma}} \frac{ (N+ \delta_1/\Lambda_1)^{3- \sigma}}{|  x- \xb|^{3- \sigma}} \dd S_{\xb}.\label{bound1_1}
			\Ee
			Recall (\ref{eta}) and choose $p \in \p\O$ such that $\xb$ is close to $p$ locally. Then locally (\ref{bound1_1}) is bounded by 
			\Be\label{bound1_1_eta}
			\begin{split}
	&	(N+ \delta_1/\Lambda_1)^{3- \sigma}	\\
	& \times \int_{|(y_{\parallel,1}, y_{\parallel,2})| \ll 1}\frac{|(x- \eta_p(y_{\parallel,1}, y_{\parallel,2}, 0)) \cdot n(y_{\parallel,1}, y_{\parallel,2})|^{1-\sigma}}{|x- \eta_p(y_{\parallel,1}, y_{\parallel,2}, 0)|^{3- \sigma}} \dd y_{\parallel,1} \dd y_{\parallel,2}.
			\end{split}\Ee
			Without loss of generality we may assume $|x- \eta_p(y_{\parallel,1}, y_{\parallel,2}, 0)| \ll 1$, otherwise the denominator of (\ref{bound1_1_eta}) has a lower bound. Hence $x$ is close to $p$ and we can parametrize $x= \eta_p(x_{\parallel,1}, x_{\parallel,2} , 0 ) - x_{n} n_p(x_{\parallel,1}, x_{\parallel,2} )$ as (\ref{X_local}) or (\ref{etain}). By a re-parametrization we may assume that the coordinate of $p \in \p\O$ is $(y_{\parallel}, y_n) =( y_{\parallel,1}, y_{\parallel,2}, y_n) = (0,0,0)$. Therefore we can assume $(x_{\parallel,1}, x_{\parallel,2},x_n)\sim(0,0,0)\sim (y_{\parallel,1}, y_{\parallel,2},0)$. \hide
			
			 Furthermore, upon a rotation and translation, we may assume that on $p$, therefore at $(0,0,0)$, we have $\eta_p (0,0,0) = 0$, $n_p(0,0) = - e_3$, and $\frac{\p \eta_{p}}{\p x_{\parallel,i}} (0,0,0) \cdot n_p(0,0)= 0$ for $i=1,2$ . Equipped with these assumptions we can expand 
			\Bes
			\eta_p( z_{\parallel,1},z_{\parallel,2} ,0)   =    z_{\parallel }  \cdot \nabla_{x_{\parallel}} \eta_p (0,0,0)+ \frac{1}{2} z_{ \parallel} \cdot \nabla_{x_{\parallel}} ^2 \eta_{p} (0,0,0) \cdot z_{\parallel}  + O_{\| \eta _{p}\|_{C^3}}( |z_{\parallel}|^3),
			\Ees
			for any $|(z_{\parallel,1}, z_{\parallel,2})| \ll1$. Note that $ z_{\parallel}  \cdot \nabla_{x_\parallel} \eta_p (0,0,0) \perp n_p(0,0)$.

			Hence the numerator of (\ref{bound1_1_eta}), for $0 \leq \sigma \leq 1$, \unhide
			
			\hide We expand
			\Bes
			&&  \big(x- \eta_p(x_{\parallel},   0)\big) \cdot \big(- n_p(y_{\parallel})\big)  \\
			&=&
			\big(
			  x_n (- n_p (x_\parallel))+(\eta_p(x_\parallel, 0) - \eta_p (y_\parallel,0)) 
			\big)\cdot \big(- n_p(y_{\parallel})\big)
			\\
			&=&
			x_n + 
			\Big(
			- \int^1_0 
			(x_\parallel - y_\parallel) \cdot \nabla_{x_\parallel} n_p (s x_\parallel +(1-s) y_\parallel)
			\dd s
			\Big)
			\cdot\big( -n_p (y_\parallel)\big) \times  x_n\\
			&&+ 
			\Big(\int^1_0 
			(x_\parallel - y_\parallel) \cdot \nabla_{x_\parallel} \eta_p (s x_\parallel +(1-s) y_\parallel,0)
			\dd s
			\Big)
			\cdot\big( -n_p (y_\parallel)\big) 
			\\
			&=&
			x_n + \| n_p \|_{C^1} |x_\parallel - y_\parallel| O(|x_n|)
\\
&&+  
			\cancel{ (x_\parallel - y_\parallel) \cdot \nabla_{x_\parallel} \eta_p ( y_\parallel,0)\cdot\big( -n_p (y_\parallel)\big)  }
\\
&&+\Big(
\int^1_0(x_\parallel - y_\parallel) \cdot   \int^s_0  (x_\parallel - y_\parallel) \cdot \nabla_{x_\parallel}
 \nabla_{x_\parallel} \eta_p ( \tau x_\parallel +(1- \tau) y_\parallel,0) \dd \tau \Big) \cdot \big(-n_p (y_\parallel)\big)
\dd s 
\\
&=&
x_n + \| \eta_p \|_{C^2} |x_\parallel - y_\parallel| O(|x_n|)
+ \| \eta_p \|_{C^2}|x_\parallel - y_\parallel|^2.
%
			\Ees
			
			\unhide
			
			We expand
			\Be
			\begin{split}\label{(x-eta)n}\notag
			&   \big(x- \eta_p(y_{\parallel},   0)\big) \cdot \big(- n_p(y_{\parallel})\big)  \\
			 =& \ 
			\big(
			  x_n (- n_p (x_\parallel))+(\eta_p(x_\parallel, 0) - \eta_p (y_\parallel,0)) 
			\big)\cdot \big(- n_p(y_{\parallel})\big)
			\\
			 =& \  
			\Big( (- n_p (y_\parallel)) 
			- \int^1_0 
			(x_\parallel - y_\parallel) \cdot \nabla_{x_\parallel} n_p (s x_\parallel +(1-s) y_\parallel)
			\dd s
			\Big)
			\cdot\big( -n_p (y_\parallel)\big) \times  x_n\\
			 &+ 
			\Big(\int^1_0 
			(x_\parallel - y_\parallel) \cdot \nabla_{x_\parallel} \eta_p (s x_\parallel +(1-s) y_\parallel,0)
			\dd s
			\Big)
			\cdot\big( -n_p (y_\parallel)\big) 
			\\
			 =& \ 
			x_n + O_{\| n_p \|_{C^1}}( |x_\parallel - y_\parallel|  |x_n|)
\\
 &+  
			\cancel{ (x_\parallel - y_\parallel) \cdot \nabla_{x_\parallel} \eta_p ( y_\parallel,0)\cdot\big( -n_p (y_\parallel)\big)  }\\
 &-n_p (y_\parallel) \cdot \Big(
\int^1_0   \int^s_0  (x_\parallel - y_\parallel) \cdot \nabla_{x_\parallel}^2
  \eta_p ( \tau x_\parallel +(1- \tau) y_\parallel,0) \cdot (x_\parallel - y_\parallel)  \dd \tau \dd s  \Big) 
\\
 =& \ 
x_n + O_{\| \eta_p \|_{C^2}}( |x_\parallel - y_\parallel|)  |x_n|
+
O_{\| \eta_p \|_{C^2}} (|x_\parallel - y_\parallel|^2),\end{split}\Ee 
%
			where we have used the fact, at the cancellation above, 
			\Be\label{perp_nabla_eta_n}
		 \nabla_{x_\parallel} \eta_p (y_\parallel,0) \perp n_p (y_\parallel).
			\Ee
	Clearly, for $|x_\parallel - y_\parallel| \ll 1$, 
			\Be\label{est_(x-eta)n}
			\begin{split}
	& \frac{1}{2}\big( |x_n| - C_{\| \eta_p \|_{C^2}}  |x_\parallel - y_\parallel|^2\big)\\
	 \leq 	& 	\
			|  \big(x- \eta_p(y_{\parallel},   0)\big) \cdot \big(- n_p(y_{\parallel})\big) | 
		\leq    \  2 \big( |x_n| + C_{\| \eta_p \|_{C^2}}   |x_\parallel - y_\parallel|^2\big).
			\end{split}\Ee\hide
			From the convexity (\ref{convexity_eta}) and (\ref{(x-eta)n}), for $|x_\parallel - y_\parallel| \ll 1$,
			\Be\label{lower_(x-eta)n}
			\begin{split}
		 	|  \big(x- \eta_p(y_{\parallel},   0)\big) \cdot \big(- n_p(y_{\parallel})\big) | 
		\gtrsim    \  |x_n| +  |x_\parallel - y_\parallel|^2.
			\end{split}
			\Ee \unhide

			On the other hand, for $i=1,2$,
			\Be
			\begin{split}\notag
			&   \big(x- \eta_p(y_{\parallel},   0)\big) \cdot  \frac{\p_i \eta_p (y_\parallel,0)}{\sqrt{g_{p,ii} (y_\parallel,0)}}  \\
			 =& \ 
			\big(
			  x_n (- n_p (x_\parallel))+(\eta_p(x_\parallel, 0) - \eta_p (y_\parallel,0)) 
			\big)\cdot  \frac{\p_i \eta_p (y_\parallel,0)}{\sqrt{g_{p,ii} (y_\parallel,0)}} 
			\\
			 =& \   
			\Big( \big(- n_p (y_\parallel)\big)
			- \int^1_0 
			(x_\parallel - y_\parallel) \cdot \nabla_{x_\parallel} n_p (s x_\parallel +(1-s) y_\parallel)
			\dd s
			\Big)
			\cdot 
			 \frac{\p_i \eta_p (y_\parallel,0)}{\sqrt{g_{p,ii} (y_\parallel,0)}} 
			  x_n\\
			 &+ 
			\Big(\int^1_0 
			(x_\parallel - y_\parallel) \cdot \nabla_{x_\parallel} \eta_p (s x_\parallel +(1-s) y_\parallel,0)
			\dd s
			\Big)
			\cdot  \frac{\p_i \eta_p (y_\parallel,0)}{\sqrt{g_{p,ii} (y_\parallel,0)}} 
			\\
			 =&  
		\cancel{\big(- n_p (y_\parallel)\big)	 \cdot 
			 \frac{\p_i \eta_p (y_\parallel,0)}{\sqrt{g_{p,ii} (y_\parallel,0)}} 
			  x_n}
+			 O_{ \| \eta_p \|_{C^2}}(|x_\parallel- y_\parallel|)  x_n, 
			\ \ \ \ \   {from}  \  (\ref{perp_nabla_eta_n})
			  \\
			 &+ \big(
			 (x_\parallel - y_\parallel) \cdot \nabla_{x_\parallel} \eta_p (  y_\parallel,0) 
			\big)
			\cdot  \frac{\p_i \eta_p (y_\parallel,0)}{\sqrt{g_{p,ii} (y_\parallel,0)}} 
			+ O_{\| \eta _p \|_{C^2}} (|x_\parallel - y_\parallel|^2) .
			 %
			\end{split}\Ee
			Hence, for $|x_n|, |x_\parallel - y_\parallel| \ll 1$,  
			\Be\begin{split}\label{(x-eta)eta}
		  &   \bigg(  {	 \sum_{i=1,2}\big| \big(x- \eta_p(y_{\parallel},   0)\big) \cdot  \frac{\p_i \eta_p (y_\parallel,0)}{\sqrt{g_{p,ii} (y_\parallel,0)}}\big|^2}
		 \bigg)^{1/2}
		 \\
		= \  &    |x_\parallel - y_\parallel| +  O_{ \| \eta_p \|_{C^2}}(|x_\parallel- y_\parallel|)  x_n + O_{\| \eta _p \|_{C^2}} (|x_\parallel - y_\parallel|^2) \\
		\geq \ & \frac{1}{2}  |x_\parallel - y_\parallel|.
			\end{split}\Ee
			
			From (\ref{est_(x-eta)n}) and (\ref{(x-eta)eta}), we derive an upper bound of the integrand of (\ref{bound1_1_eta}) and then split it as

		\hide

			\Bes
			&& |x- \eta_p( x_{\parallel} , 0)|^{3 -\sigma}\\
			&=&\Big[ |(x- \eta_p( x_{\parallel} , 0)) \cdot (-n_p (0,0))|^2 + \sum_{i=1,2}  |(x- \eta_p(x_{\parallel} , 0)) \cdot \frac{\p_i \eta_p (0)}{\sqrt{g_{p,ii}(0)}}|^2\Big]^{\frac{3- \sigma}{2}}\\
			&=&\Big[
			\big| x_{n} - \frac{1}{2} (x_{\parallel}
			\cdot \nabla_{\parallel}^2 \eta_p(0) \cdot x_{\parallel}
			) \cdot (-n_p(0,0))  \big|^2
			+ |x_{\parallel}|^2 + O_{\|\eta_{p}\|_{C^3}} ( |x_{\parallel}|^3  )\Big]^{\frac{3- \sigma}{2}}\\
			&\gtrsim &\Big[\big| x_{n} - \frac{1}{2} (x_{\parallel}
			\cdot \nabla_{\parallel}^2 \eta_p(0) \cdot x_{\parallel}
			) \cdot (-n_p(0,0))  \big|^2
			+ |x_{\parallel}|^2\Big]^{\frac{3- \sigma}{2}}.
			\Ees\unhide
			
			\begin{eqnarray} 
			 (\ref{bound1_1_eta}) 
			& \lesssim&     \int_{|x_{\parallel}- y_\parallel| \ll 1} 
			\frac{  |x_n|^{1-\sigma} +   |x_\parallel - y_\parallel|^{2(1-\sigma)}}{
		\Big(\big(  |x_n| -  C_{\| \eta_p \|_{C^2}}  |x_\parallel - y_\parallel|^2\big)^2
		+     |x_\parallel - y_\parallel|
		 ^2\Big)^{\frac{3-\sigma}{2}}
			}
			 \dd y_{\parallel} \notag\\
			&=&   \int_{  | x_{\parallel}- y_\parallel| \leq  |x_{n}| }  + \int_{   |x_{n}|  \leq| x_{ \parallel}-y_\parallel|   }   .\label{bound1_1_eta_split}
			\end{eqnarray}
			
			First we consider the case of $|x_{\parallel}- y_\parallel| \leq   |x_{n}|$. For $|x_\parallel - y_\parallel| \ll1 $, 
			\[ |x_n| -  C_{\| \eta_p \|_{C^2}}  |x_\parallel - y_\parallel|^2
			\geq \frac{|x_n| }{2}.
			\]
			\hide
			$$\big| x_{n} - \frac{1}{2} (x_{\parallel}
			\cdot \nabla_{\parallel}^2 \eta_{p}(0) \cdot x_{\parallel}
			) \cdot (-n_p(0,0))  \big|\geq \frac{|x_{n}|}{2}.$$  
			\unhide
			Hence, by a change of variables $R: = |x_{\parallel}- y_\parallel|^2$,
			\Be\begin{split}\label{bound_11eta_1}
				&\int_{  | x_{\parallel}- y_\parallel| \leq  |x_{n}| }  \ \ \text{in} \ \ (\ref{bound1_1_eta_split}) \\
				\lesssim& \   \int_0^{|x_{n}|^2}  \frac{|x_{n}|^{1-\sigma} + R^{1- \sigma} }{[ {|x_{n}|^2}/ {4 } + R]^{\frac{3-\sigma}{2}}
				}\dd R\\
				\lesssim&  \  \int_0^{|x_{n}|^2} \frac{|x_{n}|^{1-\sigma}}{|x_{n}|^{3- \sigma}} \dd R + \int_0^{|x_{n}|^2} R^
				{- \frac{1}{2} - \frac{\sigma}{2}} \dd R\\
				\lesssim& \  1 + |x_{n}|^{1- \sigma}.
			\end{split} \Ee 
		On the other hand, from the same change of variables,  
			\Be\begin{split}\label{bound_11eta_2}
				&\int_{   |x_{n}|  \leq| x_{\parallel}-y_\parallel|  }    \ \ \text{in} \ \ (\ref{bound1_1_eta_split})   
				\\
				\lesssim& \ \int_{|x_{n}|^2}^{1}  \frac{|x_{n}|^{1- \sigma}}{R^{\frac{3-\sigma}{2}}}\dd R +  \int_{|x_{n}|^2}^{1} 
				R^{- \frac{1}{2} - \frac{\sigma}{2}}
				\dd R  \\
				\lesssim& \  1+ |x_{n}|^{1 - \sigma}.
			\end{split}\Ee
		From (\ref{bound_11eta_1}) to (\ref{bound_11eta_2}), and (\ref{bound1_1}),  (\ref{bound1_1_eta}), we conclude $(\ref{bound1}) \lesssim 1$.

			In a similar way as deriving (\ref{bound1_1}), using (\ref{lower_tb_lessN}), we bound 
			\Be \notag
			(\ref{bound2})
			\lesssim  \frac{\| E\|_{L^\infty} (N + \delta_1/\Lambda_1)^{1+ \sigma}}{1+\sigma} \int_{\p\O}  \frac{
				\dd S_{\xb}
			}{|\xb  -x|^{1+ \sigma}}.
			\Ee 
			Since $\xb$ varies in the 2 dimensional smooth hypersurface $\p \O$, we can easily see that $\frac{
				1
			}{|\xb  -x|^{1+ \sigma}}$ is integrable for $0<\sigma<1$ and hence $(\ref{bound2})$ is bounded. This estimate together with (\ref{bound_11eta_2}) derives a bound of (\ref{apply_Lemma_COV}) and hence (\ref{NLL_split2}).


			\vspace{4pt}
			
			\textit{Step 2.} Now we consider (\ref{NLL_split3}). From $|u| \geq N \gg 1$ and (\ref{decay_E}), we have 
			\Bes
			|X(s;t,x,u) - x| 
			&=&  \Big|\int^t_s V(s^\prime;t,x,u) \dd s^\prime\Big|\\
			&\geq& |u| |t-s| - \int^t_s \int_{s^\prime}^{t} 
			e^{- \Lambda_1 s^{\prime\prime}} \delta_1
			\dd  s^{\prime\prime}\dd s^\prime\\
			&\geq&  
			( |u|  - {\Lambda_1}^{-1} {\delta_1} )|t-s|
			.
			\Ees
			If $N> 2 \Lambda_1^{-1}\delta_1$ then $|X(s;t,x,u) - x|  \geq \frac{|u|}{2} |t-s|$. Since the domain $\O$ is bounded, we have $\sup |X(s;t,x,u) - x|< \infty$. Finally we conclude that 
			\begin{eqnarray}
			|u| \tb(t,x,u) \lesssim_\Omega 1& {for}&|u|\geq N,
			\label{upper_|v|t_b}
			\\
			\tb(t,x,u)\lesssim_{\Omega, N} 1& {for}&|u|\geq N.\label{upper_t_b}
			\end{eqnarray}

			Furthermore from (\ref{upper_|v|t_b}), for $|u| \geq N$
			\Be \begin{split}\label{upper_int_V}
				& \int_{s- \tb(s,x,u)}^s |V(s^\prime;s,x,u)| \dd s^\prime\\
				&\leq  |u|\tb(s,x,u) + 
				\int_{s- \tb(s,x,u)}^s 
				\int_{s^\prime}^s |E (s^{\prime\prime}, X(s^{\prime\prime};s,x,u))| \dd s^{\prime\prime}
				\dd s^\prime\\
				& \lesssim_\O
				1+ \int_{s- \tb(s,x,u)}^s 
				\int_{s^\prime}^s
				\delta_1e^{-\Lambda_1 s^{\prime\prime}}
				\dd s^{\prime\prime}
				\dd s^\prime, \ \ \  {from} \ \ (\ref{decay_E})
				\\
				&\lesssim_\O  1 + \frac{\delta_1}{(\Lambda_1)^2},
			\ \ \  {from} \ \ 	(\ref{direct_computation})
				\\
				& \lesssim_{\O, K , \Lambda_{1}, \delta_1} 1.
			\end{split}\Ee

			\vspace{4pt}
			
			\textit{Step 3.} Now we derive a version of velocity lemma (See several versions of velocity lemma in \cite{Guo_V,Hwang,Guo10,GKTT1}). Without loss of generality we assume that there exists $\xi: \R^3 \rightarrow \R$ such that $\O= \{x \in \R^3: \xi(x)<0\}$ and $\nabla \xi(x) \neq 0$ when $\xi(x)=0$ (See the construction of such $\xi$ in (2.3) of \cite{EGKM2}).
			 We define 
			\Be\label{beta}
			\tilde{\alpha}(t,x,v) : = \sqrt{\xi(x)^2  + |\nabla \xi(x) \cdot u|^2 - 2 (u\cdot \nabla_x^2 \xi(x) \cdot u) \xi(x)
			}.
			\Ee

			For $|u| \geq N$ and $t- \tb(t,x,u)\geq - \e/2$,
			\Be\label{velocity_N}
			 {\alpha}_{f,\e}(t,x,u)^2  \lesssim_{\Omega, N, \delta_1, \Lambda_1} 
				\tilde{\alpha} ( t,x,u)^2 \lesssim_{\Omega, N, \delta_1, \Lambda_1}   {\alpha}_{f,\e}(t,x,u)^2. 
			\Ee 
	 We remark that the equivalent relationship fails in general for small $u$ as $N \downarrow 0$.

			The proof of (\ref{velocity_N}) is now given here: Assume $|u|\geq N$. Then for $(\Lambda_1)^{-1}\delta_1 \leq 1$ and $t-\tb(t,x,u) \leq s \leq t$
			\Be\label{|u|_bound}
			N/2 \leq N- \Lambda_1^{-1}\delta_1\leq |V(s;t,x,u)| \leq N + \Lambda_1^{-1}\delta_1 \leq 2N.
			\Ee

			From a direct computation, 
			\Be\begin{split}\label{velocity_lemma_N}
				& [\p_t + u\cdot \nabla_x +E(t,x) \cdot \nabla_u ]\{ \xi(x)^2  + |\nabla \xi(x) \cdot u|^2 - 2 (u\cdot \nabla_x^2 \xi(x) \cdot u) \xi(x) \}\\
				&=2 \{u \cdot \nabla  \xi \} \xi  \cancel{+ 2 \{u\cdot \nabla ^2 \xi  \cdot u\} \{ u \cdot \nabla_x \xi  \}}-2 u \cdot (u \cdot \nabla  \nabla ^2 \xi  \cdot u) \xi \\
				& \ \ \cancel{ - 2 \{u\cdot \nabla ^2 \xi  \cdot u\}\{ u\cdot \nabla  \xi\}} + 2 \{  E \cdot \nabla \xi  \} \{ \nabla  \xi  \cdot u\}-4 \{ E \cdot \nabla^2 \xi \cdot u \}\xi\\
				&\lesssim |u \cdot \nabla  \xi|^2 + |\xi|^2 + \{|u|+ \frac{1}{|u|}\} (- 2 (u\cdot \nabla_x^2 \xi(x) \cdot u) \xi(x))\\
				& \ \ \ \ \ 
				+ | E  \cdot \nabla \xi  
				|| \nabla  \xi  \cdot u|.
			\end{split}\Ee
			From the Neumann BC ($n(x) \cdot E(t,x) =0$ on $x \in \p\O$), we have 
			\Be\begin{split}\label{phi_xi}
				& |E (t,x) \cdot \nabla \xi(x)| \\
				&\leq |E (t,x_*) \cdot \nabla \xi(x_*)| + \|  E(t) \|_{C^1(\bar{\O})} \| \xi \|_{C^2 (\bar{\O})}
				|x-x_*|\\
				&\lesssim_\O \| E(t) \|_{C^1(\bar{\O})} 
				|\xi(x )|,\end{split}
			\Ee
			where $x_* \in\p\O$ such that $|x-x_*|= \inf_{y \in \p\O}|x-y|$.
			
			By controlling the last term of (\ref{velocity_lemma_N}) by (\ref{phi_xi}) and using (\ref{|u|_bound}), we conclude that 
			\Be\begin{split}\notag
				& \frac{d}{ds} \tilde{\alpha}(s,X(s;t,x,u),V(s;t,x,u))^2 \\
				&\lesssim_\O \Big(1+ |V(s;t,x,u)| + \frac{1}{|V(s;t,x,u)|}\Big)\tilde{\alpha}(s,X(s;t,x,u),V(s;t,x,u))^2\\
				&\lesssim_{\O,R, N} \Big(1+ |V(s;t,x,u)|\Big)\tilde{\alpha}(s,X(s;t,x,u),V(s;t,x,u))^2.
			\end{split}\Ee 
			Then, using (\ref{upper_|v|t_b}), (\ref{upper_t_b}), and (\ref{upper_int_V}), we derive 
			\Be \label{upper_lower_tilde_alpha}
			\begin{split}  &\tilde{\alpha}(t-\tb(t,x,u),X(t-\tb,t,x,u),V(t-\tb,t,x,u))^2\\
				& \lesssim_{\Omega, N, \delta_1, \Lambda_1} 
				\tilde{\alpha} ( t,x,u)^2\\
				& \lesssim_{\Omega, N, \delta_1, \Lambda_1}   \tilde{\alpha}(t-\tb(t,x,u),X(t-\tb,t,x,u),V(t-\tb,t,x,u))^2. 
			\end{split}\Ee
			Note that from $(X(t-\tb,t,x,u),V(t-\tb,t,x,u)) \in \gamma_-$, (\ref{alphaweight}), (\ref{beta}), and the fact  $\nabla \xi (x) \neq 0$ when $\xi(x)=0$ (i.e. $x\in \p\O$), we have 
			\Be\begin{split}\notag
				&\tilde{\alpha}(t-\tb(t,x,u),X(t-\tb,t,x,u),V(t-\tb,t,x,u))\\
				\
				\equiv   \ & |\nabla \xi (X(t-\tb,t,x,u)) \cdot V(t-\tb,t,x,u)|
				\\
				\
				\sim   \ & |n(X(t-\tb,t,x,u)) \cdot V(t-\tb,t,x,u)|
				\\
				\
				\sim   \ & {\alpha}_{f,\e}(t,x,u),
				\ \ \ for \ t- \tb(t,x,u)\geq - \e/2
				.
			\end{split}\Ee
		Here the equivalent relation $``\sim$'' depends on $\Omega$.

			Combining the above estimate with (\ref{upper_lower_tilde_alpha}) we finish the proof of (\ref{velocity_N}).
			
			\vspace{4pt}
			
			\textit{Step 4.} From (\ref{velocity_N}) and (\ref{alphaweight}), we have, for $|u| \geq N$,
\Be
\tilde{\alpha} (t,x,u)\gtrsim \mathbf{1}_{t- \tb(t,x,u)\geq - \e/2} |\nabla \xi (x) \cdot u|  + \mathbf{1}_{t- \tb(t,x,u)\leq - \e/2}.
\Ee
			Therefore we conclude that 
			\Bes
			(\ref{NLL_split3} ) &\lesssim&1+
			\int_{0}^\infty \frac{e^{-C|(v-u) \cdot \nabla \xi|^2}}{|\nabla \xi (x) \cdot u|^\sigma}\dd (\nabla \xi (x)  \cdot u) \\
			&& \ \ \ \ \ \ \times \iint_{\R^2} \frac{1}{|u_{\perp}-v_{\perp}|^{2-\kappa} }e^{-C| v_\perp-u_\perp   |^2} \dd u_\perp \\
			&\lesssim&1+
			\int_{0}^\infty \frac{e^{-C|(v-u) \cdot \nabla \xi|^2}}{|\nabla \xi (x) \cdot u|^\sigma}\dd (\nabla \xi (x)  \cdot u) \times \int_{0}^{\infty} \frac{e^{-C|y|^2}}{|y|^{1-\kappa} } \dd |y|   \\
			&\lesssim& 1,
			\Ees
			where $u_\perp = u - (u \cdot  \frac{\nabla\xi(x)}{|\nabla\xi(x)|})   \frac{\nabla\xi(x)}{|\nabla\xi(x)|}$.\end{proof}

		\hide\Bes
		(\ref{NLL_split3} ) &\lesssim&
		\int_{0}^\infty \frac{e^{-C|(v-u) \cdot \nabla \xi|^2}}{|\nabla \xi (x) \cdot u|^\sigma}\dd (\nabla \xi (x)  \cdot u) \times \iint_{\R^2} e^{-C| v_\perp-u_\perp   |^2} \dd u_\perp \\
		&\lesssim& 1,
		\Ees\unhide

		\section{Nonlinear-Normed Energy Estimates in Weighted ${W}^{1,p}$}
		The main result of this section is the following a priori estimate.
		\begin{proposition}\label{prop_W1p}
			Let us choose $0< \tilde{\vartheta}< \vartheta \ll 1$ and 
			\Be\label{beta_condition}
			\frac{p-2}{p}<\beta< \frac{2}{3} 
			, \quad\text{for}\quad 3< p < 6. 
			\Ee
			Assume $f$ solves (\ref{eqtn_f}), (\ref{phi_f}), (\ref{BC_f}), and $f \in L^\infty((0,T); L^p (\O \times \R^3)) \cap L^1 ((0,T); L^p (\gamma))$ and $w_{\tilde{\vartheta}}\alpha_{f,\e}^\beta \nabla_{x,v} f \in L^\infty((0,T); L^p (\O \times \R^3))  \cap L^1 ((0,T); L^p (\gamma))$ and 
			\Be\label{small_bound_f_infty}
			\sup_{0 \leq t \leq \infty} \| w_\vartheta f(t) \|_\infty \ll 1,
			\Ee
			and 
			\Be\label{integrable_nabla_phi_f}
			\sup_{0 \leq t \leq T} e^{\Lambda_1 t} \|  \nabla_x \phi_f (t) \|_\infty  < \delta_1,
			\Ee
			with 
			\Be\label{delta_1/lamdab_1}
			0< \frac{  \delta_1}{\Lambda_1}  \ll_\O 1.
			\Ee
			Then there exists $C_p>0$ such that 
			\Be\label{est_W}
			\begin{split}
				&\|w_{\tilde{\vartheta}}f(t) \|_p^p
				+\|w_{\tilde{\vartheta}} \alpha_{f,\e}^\beta 
				\p f (t)\|_p^p
				\\
				& 
				\leq  \  C_p e^{ {C_p (1 
						+ \sup_{0 \leq s \leq t} \| \nabla^2 \phi_f(s) \|_\infty) t} }
				\{
				\| |w_{\tilde{\vartheta}}f(0) \|_p^p
				+\| w_{\tilde{\vartheta}} \alpha_{f,\e}^\beta 
				\p f (0)\|_p^p \}
				.\end{split} \Ee
		\end{proposition}

		\begin{lemma}
			For any $0< \delta <1$, we claim that if $(f, \phi_f)$ solves (\ref{phi_f}) then 
			\Be\label{Morrey}
			\|   \phi_f(t) \|_{C^{1, 1- \delta}(\bar{\O})}\leq C_{  \O}  \| w_\vartheta f(t) \|_\infty  \ \ \text{for all} \ \ t\geq 0.
			\Ee\end{lemma}

		\begin{proof} 
			We have, for any $p>1$, 
			\Be\notag
			\begin{split}
				\left\|\int_{\R^3} f(t,x,v) \sqrt{\mu(v)} \dd v\right\|_{L^p (\O)} 
				\leq        |\O|^{1/p} \left( \int_{\R^3} w_\vartheta(v)^{-1} \sqrt{\mu (v)} \dd v\right) \| w_\vartheta f (t) \|_\infty.
			\end{split}\Ee
			Then we apply the standard elliptic estimate to (\ref{phi_f}) and deduce that 
			\Be\label{phi_2p}
			\|\phi_f (t) \|_{W^{2,p} (\O)} \lesssim 
			\| w_\vartheta f(t) \|_\infty
			.\notag
			\Ee
			On the other hand, from the Morrey inequality, we have, for $p>3$ and $\O \subset \R^3$, 
			\[
			\|   \phi_f(t) \|_{C^{1, 1-  {3}/{p}}(\O)}\lesssim_{p, \O} \|   \phi_{f}(t) \|_{W^{2,p} (\O)}.
			\]
			Now we choose $p=3/\delta$ for $0< \delta <1$. Then we can obtain (\ref{Morrey}). Note that the constant $C_\O$ can be chosen independent of $p$ and $\delta$.  
		\end{proof}

		To close the estimate, we use the following lemma crucially.
		\begin{lemma}\label{lemma_apply_Schauder}
			Assume (\ref{beta_condition}). If $\phi_f$ solves (\ref{phi_f}) then there exists $C_{1}
			\geq 0$ such that 
			\Be\label{phi_c,gamma}
			\| \phi_f (t) \|_{C^{2, 1- \frac{3}{p}} (\bar{\O})}    \leq (C_1)^{1/p} \big\{ \|  f(t) \|_p  + \| \alpha_{f,\e}^\beta \nabla_x f(t) \|_{p} \big\}.
			\Ee

		\end{lemma}
		
		\begin{proof}

			Applying the Schauder estimate to (\ref{phi_f}), we deduce 
			\Be
			\begin{split}\label{apply_Schauder}
				\|\phi_{f} (t)\|_{C^{2,1-\frac{3}{p}}(\bar{\O})} &\lesssim_{p,\O} \Big\|\int_{\R^{3}} f(t)\sqrt{\mu} \dd v \Big\|_{C^{0,1-\frac{3}{p}}(\bar{\O})}  \ \ \  {for} \ \ p>3.
			\end{split}\Ee

			By the Morrey inequality, $W^{1,p}   \subset C^{0,1- \frac{3}{p}} $ with $p>3$ for a domain $\O \subset \R^3$ with a smooth boundary $\p\O$, we derive 
			%
			%
			\Be\begin{split}\label{apply_Morrey}
				&\Big\| \int_{\R^{3}} f (t) \sqrt{\mu }\dd v \Big\|_{C^{0,1-\frac{n}{p}} (\bar{\O}) }  \\
				\lesssim &  \ \Big\| \int_{\R^{3}} f  (t)  \sqrt{\mu }\dd v \Big\|_{W^{1,p} (\O) } 
				\\
				\lesssim & \ \left(\int_{\R^3}   {\mu}^{q /2}\dd v\right)^{1/q}   \|  f (t) \|_{L^p (\O \times \R^3)} 
				+\Big\| \int_{\R^{3}} \nabla_{x} f  (t)  \sqrt{\mu }\dd v \Big\|_{L^p (\O) }.
			\end{split}
			\Ee
			By the H\"older inequality 
			\Be\begin{split}\label{nabla_p}
				&\Big|\int_{\R^{3}} \nabla_{x} f(t,x,v) \sqrt{\mu(v)}\dd v\Big|\\
				\leq  & \  \Big\|\frac{\sqrt{\mu(\cdot)}}{ \alpha_{f,\e}(t,x,\cdot) ^{\beta}}  \Big\|_{L^{\frac{p}{p-1} }(\R^{3})}
				\Big\| \alpha_{f,\e}(t,x,\cdot )^{\beta} \nabla_{x} f(t,x,\cdot ) \Big\|_{L^{p} (\R^{3})}\\
				= & \ \underbrace{\left( \int_{\R^{3}}  \frac{ \mu(v)^{\frac{p}{2(p-1)}}}{\alpha_{f,\e} (t,x,v)^{\frac{\beta p}{p-1}}}   \dd v \right)^{\frac{p-1}{p}}}_{(\ref{nabla_p})_1} \|\alpha_{f,\e}(t,x, \cdot)^{\b}\nabla_{x}f(t,x, \cdot)\|_{L^{p} (\R^{3})}.
			\end{split}\Ee
			Note that $\frac{p-2}{p-1}<\frac{\beta p}{p-1}< \frac{2}{3} \frac{p}{p-1}<1$ from (\ref{beta_condition}). We apply Proposition \ref{prop_int_alpha} and conclude that $(\ref{nabla_p})_1\lesssim 1$. Taking $L^p(\O)$-norm on (\ref{nabla_p}) and from (\ref{apply_Morrey}), we conclude (\ref{phi_c,gamma}).\end{proof}

		We need some basic estimates to prove Proposition \ref{prop_W1p}. Recall the decomposition of $L$ in (\ref{L_decomposition}). From (\ref{collision_frequency})
		\Be\label{nabla_nu}
		|\nabla_v \nu(v)| \leq  \int_{\R^3} \int_{\S^2} |\o| \mu(u) \dd \omega \dd u \lesssim 1. 
		\Ee

		We define a notation
		\begin{equation}\label{kzeta}
		\mathbf{k}_{  \varrho}(v,u) := \frac{1}{|v-u| } \exp\left\{- {\varrho} |v-u|^{2}  
		-  {\varrho} \frac{ ||v|^2-|u|^2 |^2}{|v-u|^2}
		\right\}.
		\end{equation}
	It is standard (see \cite{Guo10}) that for $0<\frac{\vartheta}{4}<\varrho$ and $0<\tilde{\varrho}<  \varrho- \frac{\vartheta}{4}$,
		\Be\label{k_vartheta_comparision}
		\mathbf{k}_{  \varrho}(v,u)  {e^{\vartheta |v|^2}}{e^{-\vartheta |u|^2}} \lesssim  \mathbf{k}_{\tilde{\varrho}}(v,u) .
		\Ee 
Moreover, for $0<\frac{\vartheta}{4}<\varrho$, (see the proof of Lemma 7 in \cite{Guo10})
		\Be\label{grad_estimate}
		\int_{\R^3}\mathbf{k}_{  \varrho}(v,u) {e^{\vartheta |v|^2}}{e^{-\vartheta |u|^2}} \dd u  \lesssim \langle v\rangle^{-1}.
		\Ee

		\hide
		\begin{lemma} \label{DKG}
			For $0<\varrho< \frac{1}{8}$,
			\Be\label{vK}\begin{split}
				| \nabla_v Kg(v) | \lesssim    \ \| w g \|_\infty
				+ \int_{\R^3} \mathbf{k}_\varrho (v,u) |\nabla_v g(u)| \dd u,
			\end{split}
			\Ee
			and
			\Be\label{vGamma}
			\begin{split}
				|\nabla_v \Gamma (g,g) (v)| 
				\lesssim  & \ \| w g \|_\infty \int_{\R^3} \mathbf{k}_\varrho (v,u) |\nabla_v g (u)| \dd u \\
				&+ \langle v\rangle  \| w g \|_\infty |\nabla_v g (v)|  + \langle v \rangle w(v)^{-1} \| w g \|_\infty^2.
			\end{split}
			\Ee
			
		\end{lemma}
		\unhide
		From (\ref{k_estimate}) and a direct computation, for $0<\varrho< \frac{1}{8} $,
		\Be \label{nabla_k1}
		\begin{split}
			|\p_{v_{i}} \mathbf{k}_{1}(v,u+v)| &= C_{\mathbf{k}_{1}} \p_{v_{i}}\Big( |u|e^{-\frac{|v|^{2}+|u+v|^{2}}{4}} \Big)  
			\lesssim 
			\mathbf{k}_{\varrho}(v,u+v),
		\end{split} 
		\Ee
		and 
		\Be \label{nabla_k2}
		\begin{split}
			\p_{v_{i}} \mathbf{k}_{2}(v,u+v) &= C_{\mathbf{k}_{2}} \p_{v_{i}}\Big( \frac{1}{|u|}e^{-\frac{|u|^{2}}{8}} e^{-\frac{||v|^{2}-|u+v|^{2}|^{2}}{8|u|^{2}}}  \Big) \\
			&= -\frac{C_{\mathbf{k}_{2}}}{|u|}e^{-\frac{|u|^{2}}{8}} e^{-\frac{||v|^{2}-|u+v|^{2}|^{2}}{8|u|^{2}}} \frac{||v|^{2}-|u+v|^{2}|}{4|u| } \frac{u_{i}}{|u|}  \\
			&\lesssim {|u|^{{-1}}}{e^{-\frac{|u|^{2}}{8}} 
			}
			e^{-\frac{
					||v|^{2}-|u+v|^{2}|^{2}
				}{16|u|^{2}}} 
			\\
			&
			\lesssim \mathbf{k}_{\varrho}(v,u+v) .  \\
		\end{split} 
		\Ee


		We define 
		\Be\label{K_v}
		K_v g (v): = \int_{\R^3}\{ \nabla_v \mathbf{k}_2(v,u+v) - \nabla_v \mathbf{k}_{1} (v, u+v)\} g(u+v) \dd u.
		\Ee
		From (\ref{k_vartheta_comparision}), (\ref{nabla_k1}), and (\ref{nabla_k2}),
		\Be\label{vKsum}\begin{split}
			| w_{\tilde{\vartheta}}  K \nabla_v g(v) |
			&\leq \ \sum_{i}\int_{\R^3} |\mathbf{k}_i  (v,u+v) | \frac{w_{\tilde{\vartheta}}(v)}{w_{\tilde{\vartheta}}(u+v)} |w_{\tilde{\vartheta}} \nabla_v g(u+v)|\}\dd u\\
			&\lesssim \int_{\R^3} \mathbf{k}_{\tilde{\varrho}} (v,u) |w_{\tilde{\vartheta}}\nabla_v g(u)| \dd u, \\
		| w_{\tilde{\vartheta}} K_{v} g(v) |
			&\leq \ \sum_{i}\int_{\R^3} | \nabla_v\mathbf{k}_i  (v,u+v)| \frac{w_{\tilde{\vartheta}}(v)}{w_{\vartheta}(u+v)} |w_{\vartheta} g(u+v) | \dd u  \\
			&\lesssim 
			\int_{\R^3} \mathbf{k}_\varrho (v,u) \frac{w_{\tilde{\vartheta}}(v)}{w_{\vartheta}(u)} |w_{\vartheta} g(u)|  \dd u  \\
			&\lesssim \| w_\vartheta g \|_\infty.  \\
		\end{split}
		\Ee

		The nonlinear Boltzmann operator (\ref{Gamma_def}) equals
		\begin{equation}\label{carleman}
		\begin{split}
		& 
		\int_{\mathbb{R}^{3}}   \int_{\S^2}   | u \cdot \omega|
		g_1(v+ u_\perp) g_2(v + u_\parallel)
		\sqrt{\mu(v+u)} \dd \omega  \mathrm{d}u \\
		& - \int_{\mathbb{R}^{3}}   \int_{\S^2}   | u \cdot \omega|
		g_1(v+u) g_2(v  )
		\sqrt{\mu(v+u)} \dd \omega  \mathrm{d}u,
		\end{split}
		\end{equation}
		where $u_\parallel = (u \cdot \omega)\omega$ and $u_\perp = u - u_\parallel$. Following the derivation of (\ref{k_estimate}) in Chapter 3 of \cite{gl}, by exchanging the role of $\sqrt{\mu}$ and $w^{-1}$, we have 
		\Be\begin{split}\label{bound_Gamma_k}
			|w_{\vartheta}\Gamma  (g, g)|   \lesssim \| w_{\vartheta} g \|_\infty \int_{\R^3} \mathbf{k}_{\tilde{\varrho }} (v,u) |w_{\vartheta}g(u) | \dd u.
		\end{split}\Ee

		By direct computations
		\Be\begin{split}\label{nabla_Gamma}
			&\nabla_v \Gamma(g,g) (v) \\
			=& \  \nabla_{v}\Gamma_{\textrm{gain}}(g, g) - \nabla_{v} \Gamma_{\textrm{loss}}(g, g) \\
			= & \ \Gamma_{\textrm{gain}} (\nabla_v g,g) + \Gamma_{\textrm{gain}} ( g,\nabla_v g)
			\\ & \ - \Gamma_{\textrm{loss}} (\nabla_v g,g) -\Gamma_{\textrm{loss}} ( g,\nabla_vg)+ \Gamma_v (g,g).
		\end{split}\Ee
		Here we have defined
		\Be\begin{split}\label{Gamma_v}
			\Gamma_v (g_{1},g_{2})(v) &:= \Gamma_{v,\textrm{gain}} - \Gamma_{v,\textrm{loss}} \\
			&:=\int_{\mathbb{R}^{3}}   \int_{\S^2}   | u \cdot \omega|
			g_1(v+ u_\perp) g_2(v + u_\parallel)
			\nabla_v\sqrt{\mu(v+u)} \dd \omega  \mathrm{d}u  \\
			&\quad - \int_{\mathbb{R}^{3}}   \int_{\S^2}   | u \cdot \omega|
			g_1(v+u) g_2(v  )
			\nabla_v \sqrt{\mu(v+u)} \dd \omega  \mathrm{d}u.
		\end{split}\Ee
		
		Note that 
		\Be\begin{split}\notag
			&|w_{\tilde{\vartheta}}\Gamma_{\textrm{gain}} (\nabla_v g,g) |+ |w_{\tilde{\vartheta}}\Gamma_{\textrm{gain}} ( g,\nabla_vg)|\\
			\lesssim & \ \| w_{\vartheta} g \|_\infty \big\{|w_{\tilde{\vartheta}}\Gamma_{\textrm{gain}} (|\nabla_v g|, w_{\vartheta}^{-1}) | + |w_{\tilde{\vartheta}}\Gamma_{\textrm{gain}} ( w_{\vartheta}^{-1}, |\nabla_vg|)|\big\} \\
			\lesssim & \ \| w_{\vartheta} g \|_\infty  
			\int_{\R^3} \int_{\S^2} 
			|(v-u) \cdot \omega| \frac{w_{\tilde{\vartheta}}(v)}{w_{\vartheta}(u )}  \Big\{ \frac{|\nabla_v g (u^\prime)|}{w_{\vartheta}(v^\prime) } + \frac{ |\nabla_v g (v^\prime)|}{w_{\vartheta}(u^\prime)}\Big\}
			\dd \omega \dd u.
		\end{split}\Ee

		Then following the derivation of (\ref{k_estimate}) in Chapter 3 of \cite{gl}, by exchanging the role of $\sqrt{\mu}$ and $w_{\vartheta}^{-1}$, we can obtain a bound of 
		\Be\label{bound_Gamma_nabla_vf1}
		\begin{split}
		 &
		|w_{\tilde{\vartheta}} \Gamma_{\textrm{gain}} (\nabla_v g,g) |+ |w_{\tilde{\vartheta}}\Gamma_{\textrm{gain}} ( g,\nabla_v g)|\\
		 \lesssim  & \
		\| w_{\vartheta} g \|_\infty \int_{\R^3} \mathbf{k}_{\varrho} (v,u) \frac{w_{\tilde{\vartheta}}(v)}{w_{\tilde{\vartheta}}(u)} |w_{\tilde{\vartheta}}\nabla_v g(u)| \dd u  \\
		 \lesssim & \ \| w_{\vartheta} g \|_\infty \int_{\R^3} \mathbf{k}_{\tilde{\varrho}} (v,u)  |w_{\tilde{\vartheta}}\nabla_v g(u)| \dd u.  
		\end{split}
		\Ee

		Clearly 
		\Be\label{bound_Gamma_nabla_vf2}
		\begin{split}
			 |w_{\tilde{\vartheta}}\Gamma_{\textrm{loss}}(\nabla_v g, g)|  
			&\lesssim  \| w_{\vartheta} g \|_\infty   \int_{\R^3} 
			\frac{w_{\tilde{\vartheta}}(v)}{w_{\vartheta}(v)w_{\tilde{\vartheta}}(u)}
			|w_{\tilde{\vartheta}}\nabla_v g (u)| \mu(u)^{\frac{1}{2}} 
			\dd u \\
			&\lesssim \|w_{\vartheta} g\|_{\infty} \int_{\R^3} \mathbf{k}_{\tilde{\varrho}} (v,u) |w_{\tilde{\vartheta}}\nabla_v g(u)| \dd u , \\
			|w_{\tilde{\vartheta}}\Gamma_{\textrm{loss}}( g, \nabla_v g)| &\lesssim \langle v\rangle  \| w_{\vartheta} g \|_\infty   
			|w_{\tilde{\vartheta}}\nabla_v g (v)|  . 
		\end{split}
		\Ee
		For $\Gamma_{v,\textrm{loss}}(g,g)$ defined in (\ref{Gamma_v}),
		\Be \label{Gvloss}
		\begin{split}
			&|w_{\tilde{\vartheta}}\Gamma_{v,\textrm{loss}} (g,g)|\\
			&\lesssim \frac{w_{\tilde{\vartheta}}(v)}{w_{\vartheta} (v)} \|w_{\vartheta} g\|_{\infty} \iint_{\R^{3}\times {\S}^{2}} |(u-v) \cdot\omega|  \frac{1}{w_{\vartheta}(u)} |w_{ \vartheta}g(u)| \nabla_{v}\sqrt{\mu(u)} \dd u \dd \omega  \\
			&\lesssim  {\langle v \rangle} \|w_{\vartheta} g\|^{2}_{\infty}. 
		\end{split}
		\Ee	
		For $\Gamma_{v,\textrm{gain}}(g,g)$, following the derivation of (\ref{k_estimate}) in Chapter 3 of \cite{gl}, by exchanging the role of $\sqrt{\mu}$ and $w_{\vartheta}^{-1}$
		\Be \label{Gvgain}
		\begin{split}
			&|w_{\tilde{\vartheta}}\Gamma_{v,\textrm{gain}} (g,g)| \\
			\lesssim & \  \|w_{\vartheta}g\|_{\infty} \iint_{\R^{3}\times {\S}^{2}} |(u-v)\cdot\omega| \frac{w_{\tilde{\vartheta}}(v)}{w_{\vartheta}(v^{\prime})} \frac{w_{\vartheta}g(v^\prime)}{w_{\vartheta}(u^{\prime}) } \nabla_{v}\sqrt{\mu(u)} \dd u \dd\omega  \\
			\lesssim & \   \|w_{\vartheta}g\|_{\infty} \int_{\R^3} \mathbf{k}_{\tilde{\varrho}} (v,u) | w_{\vartheta} g(u)| \dd u.  
		\end{split}
		\Ee
		
		\hide
		we know that $w(v) \lesssim w(v^{\prime})w(u^{\prime})$ from consevation $|u|^{2} + |v|^{2} = |u^{\prime}|^{2} + |v^{\prime}|^{2}$ to derive 
		\Be \label{Gvgain}
		\begin{split}
			|\Gamma_{v,\textrm{gain}} (g,g)| &\lesssim \|w_{\vartheta} g\|_{\infty}  \iint_{\R^{3}\times {\S}^{2}} |(u-v)\cdot\omega| \frac{1}{w_{\vartheta}(u^{\prime})w_{\vartheta}(v^{\prime})} \nabla_{v}\sqrt{\mu}(u) du d\omega  \\
			&\lesssim \|w_{\vartheta}g\|_{\infty} \int_{\R^3} \mathbf{k}_{\varrho} (v,u) |\nabla_v g(u)| \dd u. 
		\end{split}
		\Ee

		We combine (\ref{nabla_Gamma}), (\ref{bound_Gamma_nabla_vf1}), (\ref{bound_Gamma_nabla_vf2}), (\ref{Gvloss}), and (\ref{Gvgain}) to obtain (\ref{vGamma}).
		
		\unhide

		The next result is about estimates of derivatives on the boundary. Assume (\ref{eqtn_f}) and (\ref{BC_f}). We claim that for $(x,v) \in\gamma_-$,
		\Be\label{BC_deriv}\begin{split}
			|\nabla_{x,v} f(t,x,v) |
			\lesssim
			\langle v\rangle\sqrt{\mu(v)}
			\Big(1+ \frac{1}
			{|n(x) \cdot v|}  \Big)
			\times (\ref{BC_deriv_R}) ,
		\end{split}\Ee   
		with
		\Be\begin{split}\label{BC_deriv_R}
			& 
			\int_{n(x) \cdot u>0}  \Big\{
			(
			\langle u\rangle  + |\nabla_x \phi_f|
			)
			|\nabla_{x,v} f(t,x,u) |  \\ 
			&
			\ \ \ \ \    \ \ \  \ \ \  \ \ \   + \langle u \rangle |f|  
			+(1+ \| w_{\vartheta} f \|_\infty) |\int_{\R^3}\mathbf{k}_{\varrho }(u, u^\prime) |f(u^\prime)| \dd u^\prime|
			\\
			&
			\ \ \ \ \    \ \ \  \ \ \  \ \ \   + 
			(\langle u \rangle |f| + \mu(u)^{\frac{1}{4}} ) |\nabla_x \phi_f |
			\Big\}\sqrt{\mu(u)}\{n(x) \cdot u\} \dd u.
		\end{split}
		\Ee
		
		From (\ref{eqtn_f}),
		%
		\Be\label{fn}
		\begin{split}
			&\partial _{n}f(t,x,v)\\
			=&\frac{-1}{n(x)\cdot v}\bigg\{ \partial
			_{t}f+\sum_{i=1}^{2}(v\cdot \tau _{i})\partial _{\tau _{i}}f 
			- \nabla_x \phi_f \cdot \nabla_v f  \\
			& \ \ \ \  \ \ \ \  \   \ \ \ + \frac{v}{2} \cdot \nabla_x \phi_f f + Lf
			-\Gamma 
			(f,f)  + v\cdot \nabla_x \phi_f \sqrt{\mu}\bigg\}. 
		\end{split}\Ee%
		%

		Let $\tau _{1}(x)$ and $\tau _{2}(x)$ be unit tangential vectors to $\partial\Omega$ satisfying $\tau
		_{1}(x)\cdot n(x)=0=\tau _{2}(x)\cdot n(x)$ and $\tau _{1}(x)\times \tau
		_{2}(x)=n(x)$. Define the orthonormal transformation from $\{n,\tau
		_{1},\tau _{2}\}$ to the standard basis $\{\mathbf{e}_{1},\mathbf{e}_{2},%
		\mathbf{e}_{3}\}$, i.e. $\mathcal{T}(x)n(x)=\mathbf{e}_{1},\ \mathcal{T}%
		(x)\tau _{1}(x)=\mathbf{e}_{2},\ \mathcal{T}(x)\tau _{2}(x)=\mathbf{e}_{3},$
		and $\mathcal{T}^{-1}=\mathcal{T}^{T}.$ Upon a change of variable: $%
		u^{  \prime }=\mathcal{T}(x)u,$ we have%
		\begin{equation*}
		n(x)\cdot u=n(x)\cdot \mathcal{T}^{t}(x)u^{\prime }=n(x)^{t}%
		\mathcal{T}^{t}(x)u^{ \prime }=[\mathcal{T}(x)n(x)]^{t}u^{
			\prime }=\mathbf{e}_{1}\cdot u^{  \prime }=u_{1}^{  \prime },
		\end{equation*}%
		then the RHS of the diffuse BC (\ref{BC_f}) equals
		\begin{equation*}
		c_{\mu }%
		\sqrt{\mu (v)}\int_{u_{1}^{  \prime }>0}f(t,x,\mathcal{T}%
		^{t}(x)u^{  \prime })\sqrt{\mu (u^{  \prime })}\{u_{1}^{
			\prime }\}\mathrm{d}u^{  \prime }.
		\end{equation*}%
		Then we can further take tangential derivatives $\partial _{\tau _{i}}$
		as, for $(x,v)\in \gamma _{-},$
		\begin{equation}\label{boundary_tau}
		\begin{split}
		&\partial _{\tau _{i}}f(t,x,v)\\
		& =c_{\mu }\sqrt{\mu (v)}\int_{n(x)\cdot u >0}\partial _{\tau
			_{i}}f(t,x,u)\sqrt{\mu (u)}\{n(x)\cdot u\}%
		\mathrm{d}u\\
		& \ \ +c_{\mu }\sqrt{\mu (v)}\int_{n(x)\cdot u>0}\nabla
		_{v}f(t,x,u)\frac{\partial \mathcal{T}^{t}(x)}{\partial \tau _{i}}%
		\mathcal{T}(x)u\sqrt{\mu (u)}\{n(x)\cdot u\}%
		\mathrm{d}u.
		\end{split}
		\end{equation}%


		
		We can take velocity derivatives directly to (\ref{BC_f}) and obtain that for $%
		(x,v)\in \gamma _{-},$
		\begin{eqnarray}
		\nabla _{v}f(t,x,v) &=&c_{\mu }\nabla _{v}\sqrt{\mu (v)}\int_{n(x)\cdot
			u>0}f(t,x,u)\sqrt{\mu (u)}\{n({x})\cdot
		u\}\mathrm{d}u,\label{boundary_v}\\
		\partial _{t}f(t,x,v) &=&c_{\mu }\sqrt{\mu (v)}\int_{n(x)\cdot u>0}\partial _{t}f(t,x,u)\sqrt{\mu (u)}\{n({x})\cdot
		u\}\mathrm{d}u. \notag
		\end{eqnarray}%
		For the temporal derivative, we use (\ref{eqtn_f}) again to deduce that 
		\Be\label{boundary_t} 
		\begin{split}
			&\partial _{t}f(t,x,v)\\
			=&c_{\mu }\sqrt{\mu (v)}\int_{n(x)\cdot u>0}
			\Big\{
			- u\cdot \nabla_x f + \nabla_x  \phi  \cdot \nabla_v f - \frac{u}{2} \cdot \nabla_x \phi f -   Lf \\
			& \ \ \ \ \ \ \ \  \ \ \ \ \ \ \ \  \ \ \ \ \ \ \ \ \ \ 
			+  \Gamma(f,f) - u\cdot \nabla_x \phi \sqrt{\mu}
			\Big\}
			\sqrt{\mu (u)}\{n({x})\cdot
			u\}\mathrm{d}u.
		\end{split}
		\Ee 
		From (\ref{fn})-(\ref{boundary_t}), (\ref{k_estimate}), and (\ref{bound_Gamma_k}), we conclude (\ref{BC_deriv}).




		\begin{proof}[\textbf{Proof of Proposition \ref{prop_W1p}}]
			\textit{Step 1.} Define 
			\Be\label{def_nu_phi}
			\nu_{\phi_f}(t,x,v):= \nu(v) + \frac{v}{2} \cdot \nabla_x \phi_f. 
			\Ee
			From the assumption (\ref{small_bound_f_infty}), we have that $\nu_{\phi_f}(t,x,v)\gtrsim \frac{\nu(v)}{2}$.
			
			From (\ref{eqtn_f}), (\ref{k_estimate}), and (\ref{bound_Gamma_k}), we can easily obtain that, for $0<\varrho \ll 1$ 
			\Be
			\begin{split}\label{energy_f_p}
				&\|w_{\tilde{\vartheta}} f(t)\|_p^p  + \int^t_0 \| \nu_{\phi_f}^{1/p} w_{\tilde{\vartheta}} f \|_p^p +  \int^t_0 | w_{\tilde{\vartheta}}f|_{p,+}^p\\
				\lesssim & \ \|w_{\tilde{\vartheta}}f( 0)\|_p^p + (1+ \| w_\vartheta f \|_\infty) \int^t_0 \iint_{\O \times \R^3} |w_{\tilde{\vartheta}}f(v)|^{p-1} \\
				&
				\ \ \ \ \ \ \ \ \ \ \ \ \ \ \ \ \ \ \ \ \ \ \ \ \ \ \ \ \ \ \ \ \  \ \ 
				 \times \int_{\R^3} \mathbf{k}_\varrho (v,u)\frac{w_{\tilde{\vartheta}}(v)}{w_{\tilde{\vartheta}}(u)} | w_{\tilde{\vartheta}}f(u)|\dd u\\
				&   +  \int^t_0 \| w_{\tilde{\vartheta}}f \|_p^p + \int_0^t\| w_{\tilde{\vartheta}}\nabla \phi_f\|_p^p +  \int^t_0 | w_{\tilde{\vartheta}}f|_{p,-}^p .
			\end{split}
			\Ee
			Note that by the H\"older inequality, (\ref{grad_estimate}), and (\ref{k_vartheta_comparision}),			\Be\begin{split}\label{k_p_bound}
				&\int_{\R^3}|w_{\tilde{\vartheta}}f(v)|^{p-1} \int_{\R^3} \mathbf{k}_{\tilde{\varrho}} (v,u)  | w_{\tilde{\vartheta}}f(u)|\dd u \dd v  \\
				\lesssim & \ \| w_{\tilde{\vartheta}} f  \|_{L^p_v}^{\frac{1}{p-1}}  \left\| \int_{\R^3} \mathbf{k}_{\tilde{\varrho}} (v,u)^{1/q} \mathbf{k}_{\tilde{\varrho}} (v,u)^{1/p}  | w_{\tilde{\vartheta}} f(u)|\dd u  \right\|_{L^p_v }\\
				\lesssim & \  \| w_{\tilde{\vartheta}} f  \|_{L^p_v}^{\frac{1}{p-1}} \left( \int_{\R^3} \mathbf{k}_{\tilde{\varrho}} (v,u) \dd u\right)^{1/q}\left\| 
				\left( \int_{\R^3} \mathbf{k}_{\tilde{\varrho}}  (v,u) |w_{\tilde{\vartheta}}f(u)|^p \dd u \right)^{1/p} \right\|_{L^p_v }\\
				\lesssim & \ \| w_{\tilde{\vartheta}}f  \|_{L^p_v}^ p  \left( \int_{\R^3} \mathbf{k}_{\tilde{\varrho}} (v,u) \dd u\right)^{1/q} \left( \int_{\R^3} \mathbf{k}_{\tilde{\varrho}}  (v,u)   \dd v \right)^{1/p} \\
				\lesssim &  \ \| w_{\tilde{\vartheta}} f  \|_{L^p_v}^ p.
			\end{split}\Ee
			From a standard elliptic theorem and (\ref{phi_f}), we have 
			\Be\label{phi_p_bound}
			\int_0^t\| \nabla \phi_f\|_p^p\lesssim \int^t_0 \| w_{\tilde{\vartheta}}f \|_p^p. 
			\Ee
			
			Now we focus on $\int^t_0 |w_{\tilde{\vartheta}}f|^p_{p,-}$ in (\ref{energy_f_p}). We plug in (\ref{BC_f}) and then decompose $\gamma_+^\e  \cup \gamma_+   \backslash \gamma_+^\e $ where $\e$ is small but satisfies (\ref{lower_bound_e}). This leads
			\Be\notag
			\begin{split}
				&\int^t_0 |w_{\tilde{\vartheta}}f|_{p,-}^p\\
				\lesssim & \ 
				\int^t_0  \int_{\p\O} \left(
				\int_{\gamma_+^\e(x)} w_{\tilde{\vartheta}}f \sqrt{\mu} \{n \cdot u\} \dd u 
				\right)^p \\
				& \ +  \int^t_0  \int_{\p\O}  \left(
				\int_{\gamma_+(x)\backslash \gamma_+^\e(x)} w_{\tilde{\vartheta}}f \sqrt{\mu} \{n \cdot u\} \dd u 
				\right)^p\\
				\lesssim & \  \Big( \int_{\gamma_+^\e} \sqrt{\mu} \{n \cdot u\} \dd u  \Big)^{p/q}\int^t_0 |w_{\tilde{\vartheta}}f|_{p,+}^p
				+ \int^t_0 \int_{\gamma_+ \backslash \gamma_+^\e} |w_{\tilde{\vartheta}}f\sqrt{\mu}|^p\\
				\lesssim & \  o(1)\int^t_0 |w_{\tilde{\vartheta}}f|_{p,+}^p
				+ \int^t_0 \int_{\gamma_+ \backslash \gamma_+^\e} |w_{\tilde{\vartheta}}f\sqrt{\mu}|^p.
			\end{split}
			\Ee
			From (\ref{eqtn_f}), Lemma \ref{le:ukai}, (\ref{k_estimate}), (\ref{bound_Gamma_k}), and (\ref{k_p_bound})
			\Be\label{bound_p_est}\begin{split}
			\int^t_0 |w_{\tilde{\vartheta}}f|_{p,-}^p\lesssim & \  \| w_{\tilde{\vartheta}}f(0)\|_p^p + o(1) \int^t_0 |w_{\tilde{\vartheta}}f|_{p,+}^p   \\
			&+(1+ \| w_\vartheta f \|_\infty) \int^t_0 \|  w_{\tilde{\vartheta}} f \|_p^p.
			\end{split}\Ee
			
			Collecting terms from (\ref{energy_f_p}), (\ref{k_p_bound}), (\ref{phi_p_bound}), and (\ref{bound_p_est}), we conclude that 
			\Be\label{Lp_estimate_f}\begin{split}
		&	\| w_{\tilde{\vartheta}}f (t) \|_p^p  + \int^t_0 \| \nu_{\phi_f}^{1/p} w_{\tilde{\vartheta}} f \|_p^p + \int^t_0 |w_{\tilde{\vartheta}}f|_{p,+}^p\\
		 \lesssim  & \ \| w_{\tilde{\vartheta}}f(0) \|_p^p + (1+ \| w_\vartheta f \|_\infty) \int^t_0 \| w_{\tilde{\vartheta}} f \|_p^p.
			\end{split}\Ee

			\vspace{4pt}
			
			\textit{Step 2.} By taking derivatives $\p \in\{ \nabla_{x_{i}}, \nabla_{v_{i}} \}$ to (\ref{eqtn_f}),
			\Be \label{eqtn_nabla_f}   [ \p_{t} + v\cdot\nabla_{x} - \nabla_{x}\phi_{f}\cdot\nabla_{v} + \nu_{\phi_{f}, w_{\tilde{\vartheta}}} ]( w_{\tilde{\vartheta}}\p f) =  w_{\tilde{\vartheta}}\mathcal{G}, \Ee
			where
			\Be  
			\begin{split}	\label{mathcal_G}
				\mathcal{G}  =& 
				- \p v \cdot\nabla_{x}f 
				+ \p \nabla \phi_f \cdot \nabla_v f\\&
				+ \p \Gamma(f,f) - \p \big[\nu(v) + \frac{v}{2} \cdot \nabla \phi_f(t,x)\big] f - \p K f - \p(v\cdot\nabla_{x}\phi_{f}\sqrt{\mu}).
			\end{split}\Ee
			Here we have used
			\Be
			\nu_{\phi_f, w_{\tilde{\vartheta}}} = \nu_{\phi_f, w_{\tilde{\vartheta}}} (t,x,v) := \nu(v) + \frac{v}{2} \cdot \nabla \phi_f(t,x)
			+ \frac{\nabla_x \phi_f \cdot \nabla_v w_{\tilde{\vartheta}}}{w_{\tilde{\vartheta}}}
			. \label{nu_phi}
			\Ee
			
			\hide
			Note that we used $K_v$ to denote 
			\begin{equation} \label{Kv}
			K_v f(v) =\int_{\mathbb{R}^3} \left\{ \nabla_u \mathbf{k}_1 + \nabla_ v\mathbf{k}_1 + \nabla_u \mathbf{k}_2 + \nabla_v \mathbf{k}_2\right\}(v,u) f(u) \dd u \lesssim \|wf\|_{\infty},
			\end{equation}
			where the last inequality comes from (\ref{nabla_k1}) and (\ref{nabla_k2}).   \\
			
			Let us choose 
			\Be\label{beta_p_1}
			\frac{p-2}{p}< \beta <\frac{p-1}{p} , \ \ p>3.
			\Ee\unhide

			From (\ref{alpha_invariant}) and (\ref{eqtn_nabla_f}),
			\Be \label{Boltzmann_p} \begin{split}
				&\frac{1}{p} |w_{\tilde{\vartheta}}\alpha_{f,\e}^\beta  \p f |^{p-1} \big[\p_t + v\cdot \nabla_x - \nabla_x \phi_f \cdot \nabla_v
				+ \nu_{\phi_f, w_{\tilde{\vartheta}}} 
				\big] |w_{\tilde{\vartheta}} \alpha_{f,\e}^\beta  \p f | \\
				&=   \alpha_{f,\e}^{\b p} | w_{\tilde{\vartheta}} \p f|^{p-1} \big[\p_t + v\cdot \nabla_x - \nabla_x \phi_f \cdot \nabla_v + \nu_{\phi_f,w_{\tilde{\vartheta}}} \big]|w_{\tilde{\vartheta}} \p f | \\
				&= w_{\tilde{\vartheta}}^p \alpha_{f,\e}^{\b p} |  \p f|^{p-1}    \mathcal{G} .
			\end{split}\Ee
			
			From (\ref{nabla_nu}), (\ref{nu_phi}), (\ref{K_v}), and (\ref{Gamma_v})
			\Be
			\begin{split}\label{G est}
				|\mathcal{G}| &\lesssim |\nabla_{x} f| 
				+  | \nabla^2 \phi_f | |\nabla_v f|
				+ |\Gamma(\p f, f)| + |\Gamma(f, \p f)|+ |K \p f|   \\
				&\quad  + |f|  +  |\Gamma_{v}(f,f)|+ |K_{v}f|\\
				& \quad  + 
				w_\vartheta(v)^{-1/2}(  | \nabla \phi_f |  + | \nabla^2\phi_f | ) (1+\|w_\vartheta f\|_\infty )    .\end{split}
			\Ee
			Now we apply Lemma \ref{lem_Green} to (\ref{Boltzmann_p}) to obtain
			\Be\begin{split} \label{pf_green}
				& \| w_{\tilde{\vartheta}} \alpha_{f,\e}^\beta 
				\p f (t)\|_p^p
				+ \int^t_0 \| \nu_{\phi_f, w_{\tilde{\vartheta}}}^{1/p} w_{\tilde{\vartheta}} \alpha_{f,\e}^\beta 
				\p f   \|_p^p  
				+ 
				\int^t_0 | w_{\tilde{\vartheta}} \alpha_{f,\e}^\beta 
				\p f   |_{p,+}^p  
				\\
				\leq&  \  \| w_{\tilde{\vartheta}} \alpha_{f,\e}^\beta 
				\p f (0)\|_p^p  + 
				\underbrace{\int^t_0 |w_{\tilde{\vartheta}} \alpha_{f,\e}^\beta 
					\p f  |_{p,-}^p  } _{(\ref{pf_green})_{\gamma_-}}\\
					&
				+ \underbrace{\int_{0}^{t} \iint_{\O\times\R^{3}}p \alpha_{f,\e}^{\b p} w_{\tilde{\vartheta}}^p| \p f|^{p-1} |\mathcal{G}| } _{(\ref{pf_green})_{\mathcal{G}}}.
			\end{split}\Ee

			First we consider $(\ref{pf_green})_{\mathcal{G}}$. Directly, the contribution of $|\nabla_{x} f| 
			+  | \nabla^2 \phi_f | |\nabla_v f|$ of (\ref{G est}) in $(\ref{pf_green})_{\mathcal{G}}$ is bounded by
			\Be
			\begin{split}\label{pf_p_nabla}
				(1+  \sup_{0 \leq s \leq t} \| \nabla^2 \phi_f \|_\infty ) \int_0^t \| w_{\tilde{\vartheta}}\alpha_{f,\e}^\beta \p f \|_p^p.
			\end{split}
			\Ee

			From (\ref{vKsum}), (\ref{bound_Gamma_nabla_vf1}), and (\ref{bound_Gamma_nabla_vf2}), the contribution of $|\Gamma(\p f, f)| + |\Gamma_{\mathrm{gain}}(f, \p f)|+ |K \p f|$ of (\ref{G est}) in $(\ref{pf_green})_{\mathcal{G}}$ is bounded by 
			\Be\begin{split}\label{pf_p_K}
				&(1 +\sup_{0 \leq s\leq  t} \|w_\vartheta f(s)\|_{\infty}) \\
				& \times \int^t_0    \iint_{\O\times\R^{3}}  |\alpha_{f,\e}^{\b  } w_{\tilde{\vartheta}}  \p f(v)|^{p-1} \\&\times\int_{\R^3} \alpha_{f,\e}(v)^\beta \mathbf{k}_\varrho(v,u)w_{\tilde{\vartheta}} (v)  |\p f(u)| \dd u  
					\dd v \dd x \dd s.
			\end{split}\Ee
			The estimate of (\ref{pf_p_K}) is carried out in \textit{Step 3}.
			
			From (\ref{bound_Gamma_nabla_vf2}), the contribution of $|\Gamma_{\mathrm{loss}}(f, \p f)|$ of (\ref{G est}) in $(\ref{pf_green})_{\mathcal{G}}$ is bounded by 
			\Be\label{pf_p_0}
			\sup_{0 \leq s\leq t}\| w_\vartheta f(s) \|_\infty
			\int^t_0  \| \nu_{\phi_f, w_{\tilde{\vartheta}}}^{1/p} w_{\tilde{\vartheta}}   \alpha_{f,\e}^\beta \p f \|_p^p
			.
			\Ee

			For the $|f|$ contribution of (\ref{G est}) in $(\ref{pf_green})_{\mathcal{G}}$, we bound
			\Be
			\begin{split}\label{pf_p_1}
				&\int^t_0 
				\iint_{\O \times \R^3}p w_{\tilde{\vartheta}}^p \alpha_{f,\e}^{\beta p } |\p f|^{p-1} 
				|f|  \dd x \dd v \dd s\\
				\lesssim & \ \int^t_0 \iint_{\O \times \R^3}
				| \nu_{\phi_f,w_{\tilde{\vartheta}} }^{1/p}  w_{\tilde{\vartheta}} \alpha_{f,\e}^{\beta} \p f| ^{p-1}
				| w_{\tilde{\vartheta}}f|
				\frac{  
					\alpha_{f,\e}(s,x,v)^{\beta} }{\nu_{\phi_f} (v)^ {{(p-1)}/{p}} }   \dd x \dd v \dd s\\
				\lesssim & \ o(1) \int^t_0 \iint_{\O \times \R^3}
				| \nu_{\phi_f, w_{\tilde{\vartheta}}}^{1/p} w_{\tilde{\vartheta}} \alpha_{f,\e}^{\beta} \p f| ^{p}
				+(1+\delta_1/\Lambda_1 )  \int^t_0 \iint_{\O \times \R^3} |w_{\tilde{\vartheta}} f|^p.
			\end{split}
			\Ee
			Here we have used the fact that, from (\ref{weight}) and (\ref{integrable_nabla_phi_f})
			\Be\label{upper_vb}
			\begin{split}
				&\alpha_{f,\e}(s,x,v)\\  \leq & \ 
				\mathbf{1}_{s +1 \geq \tb(s,x,v)}
				|\vb(s,x,v)| 
				+ \mathbf{1}_{s \leq \tb(s,x,v) +1}\\
				\lesssim &  \ 1+   |v| + \int^{0}_{-1} |\nabla \phi_f (\tau,X(\tau;s,x,v))| \dd \tau\\
				& +  \int^s_{0} |\nabla \phi_f (\tau,X(\tau;s,x,v))| \dd \tau\\
				\lesssim &  \ ( 1+ \| w_{\vartheta} f_0 \|_\infty+  {\delta_1}/{\Lambda_1} ) \langle v\rangle, 
			\end{split}
			\Ee
			and from (\ref{beta_condition}), $
			\frac{  
				\alpha_{f,\e}(s,x,v)^{\beta} }{\nu_{\phi_f} (v)^{{(p-1)}/{p}} }  \lesssim {\delta_1}/{\Lambda_1} \times \frac{\langle v\rangle^\beta }{\langle v\rangle^ {{(p-1)}/{p}} }\lesssim {\delta_1}/{\Lambda_1}.$

			From (\ref{Gvloss}), the contribution of $|\Gamma_{v, \mathrm{loss}} (f,f)|$ of (\ref{G est}) in $(\ref{pf_green})_{\mathcal{G}}$ is bounded by
			\Be
			\begin{split}\label{pf_p_Gamma_v_-}
				& \|w_{\vartheta}f\|_{\infty} \int^t_0 \iint_{\O \times \R^3} p   | \alpha_{f,\e}^\beta \p f|^{p-1} 
				\alpha_{f,\e}(v)^\beta {\langle v \rangle}w_{\vartheta}(v)^{-1} 
				\| f(s,x, \cdot) \|_{L^p }\\
				\lesssim & \ 	
				\|w_{\vartheta}f\|_{\infty} 
				\left\{
				\int^t_0 \iint_{\O \times \R^3}  | \alpha_{f,\e}^\beta \p f|^p
				+   \int^t_0 \iint_{\O \times \R^3}  |  f|^p
				\right\},
			\end{split}
			\Ee
			where we have used, from (\ref{upper_vb}), $\alpha_{f,\e}(v)^\beta {\langle v \rangle}w_{\vartheta}(v)^{-1}  \lesssim w_{\vartheta}(v)^{-1/2}$.
			
			From (\ref{vKsum}) and (\ref{Gvgain}), the contribution of $|\Gamma_{v, \mathrm{gain}}|$ and $|K_vf|$ in $(\ref{pf_green})_{\mathcal{G}}$ is bounded by 
			\Be\begin{split}\label{pf_p_Kv}
				&(1 +\sup_{0 \leq s\leq  t} \|w_\vartheta f \|_{\infty})\\
				&\times  \int^t_0  
				\iint_{\O\times\R^{3}}  \alpha_{f,\e}^{\b p} | w_{\tilde{\vartheta}} \p f (v)|^{p-1} 
				\int_{\R^3} \mathbf{k}_\varrho(v,u) \frac{w_{\tilde{\vartheta}}(v)}{w_{\tilde{\vartheta}}(u)} |f(u)| 
				\\
				\lesssim & \ 
				o(1) \int^t_0 \iint_{\O \times \R^3}
				| \nu_{\phi_f}^{1/p} \alpha_{f,\e}^{\beta} \p f| ^{p}\\
				&
				+(1 +\sup_{0 \leq s\leq  t} \|w_{\vartheta}f(s)\|_{\infty}) \int^t_0 \iint_{\O \times \R^3} |f|^p,
			\end{split}\Ee
			where we have used, for $1/p+1/p^*=1$ and $0< \tilde{\varrho} \ll \varrho$, from (\ref{k_vartheta_comparision}), (\ref{grad_estimate}),
			\Be\begin{split}\notag
				&\int_{\R^3}  \alpha_{f,\e}(v) ^{\b p} |w_{\tilde{\vartheta}} \p f (v)|^{p-1} \int_{\R^3} \mathbf{k}_{\tilde{\varrho}}(v,u) w_{\tilde{\vartheta}}(u) |f(u)| \dd u \dd v \\
				\lesssim & \ \int_{\R^3} \alpha_{f,\e}(v)^{\b p} | w_{\tilde{\vartheta}} \p f (v)|^{p-1} \int_{\R^3}\mathbf{k}_{\tilde{\varrho}}(v,u)^{1/{p^*}} \mathbf{k}_{\tilde{\varrho}}(v,u)^{1/p}  | w_{\tilde{\vartheta}}f(u)| \dd u \dd v\\
				\lesssim & \ \int_{\R^3} \frac{\alpha_{f,\e}(v)^\beta}{\langle v\rangle^{\frac{p-1}{p}}} | \langle v\rangle^{1/p} w_{\tilde{\vartheta}} \alpha_{f,\e}^\beta\p f (v)|^{p-1}\\
				& \ \ \  \times \left(\int_{\R^3} \mathbf{k}_{\tilde{\varrho}}(v,u) \dd u \right)^{1/{p^*}}
				\left(
				\int_{\R^3}  \mathbf{k}_{\tilde{\varrho}}(v,u)   | w_{\tilde{\vartheta}} f(u)|  ^p  \dd u 
				\right)^{1/p} \dd v\\
				\lesssim & \ \left( \int_{\R^3}  |\langle v\rangle^{1/p}w_{\tilde{\vartheta}} \alpha_{f,\e}^\beta\p f (v)|^p  \dd v\right)^{\frac{p-1}{p}}
				\left(
				\int_{\R^3} \int_{\R^3} \mathbf{k}_{\tilde{\varrho}}(v,u)  |w_{\tilde{\vartheta}}f(u)|^p \dd u \dd v
				\right)^{\frac{1}{p}}\\
				\lesssim & \   \left(\int_{\R^3}
				| \nu_{\phi_f}^{1/p} w_{\tilde{\vartheta}}\alpha_{f,\e}^{\beta} \p f| ^{p}\right)^{\frac{p-1}{p}}
				\left( \int_{\R^3} |w_{\tilde{\vartheta}}f|^p\right)^{\frac{1}{p}}.
			\end{split}\Ee

			Note that from the standard elliptic estimate and (\ref{phi_f}),
			\Be\label{nabla_2_phi_p}
			\|   \phi_f (t) \|_{W^{2,p}(\O  ) } \lesssim \left\|\int_{\R^3} f   (t,x,v)\sqrt{\mu(v)} \dd v \right\|_{L^p(\O )}
			\lesssim \|   f(t) \|_{L^p (\O  \times \R^3)}
			.
			\Ee
			Then from (\ref{nabla_2_phi_p}) we bound the contribution of $w_{\tilde{\vartheta}}(v)^{-1/2}(  | \nabla \phi_f |  + | \nabla^2\phi_f | ) (1+\|w_{\vartheta}f\|_\infty )$ of (\ref{G est}) in (\ref{Boltzmann_p}) by 
			\Be \label{pf_p_2}
			\begin{split}
				&(1+\|w_\vartheta f\|_\infty ) \int^t_0 
				\iint_{\O \times \R^3}p | w_{\tilde{\vartheta}} \alpha_{f,\e}^{\beta   }  \p f|^{p-1} 
				\frac{\alpha_{f,\e}(v)^{\beta}}{w_{\tilde{\vartheta}}(v)^{1/2}} (  | \nabla \phi_f |  + | \nabla^2\phi_f | )  
				\\
				\lesssim  & \ (1+\|w_\vartheta f\|_\infty ) \int^t_0 
				\iint_{\O \times \R^3}  | w_{\tilde{\vartheta}} \alpha_{f,\e}^{\beta   }  \p f|^{p-1} 
				w_{\tilde{\vartheta}}^{-1/4}  (  | \nabla \phi_f |  + | \nabla^2\phi_f | )  
				\\
				\lesssim  & \ (1+\|w_\vartheta f\|_\infty ) 
				\bigg\{
				\int^t_0 
				\iint_{\O \times \R^3}  |  w_{\tilde{\vartheta}} \alpha_{f,\e}^{\beta   }  \p f|^{p } 
			\\
			& \ \ \ \ \ \ \ \ \ \ \ \ 	\ \ \ \ \ \ \  \ \ +\int^t_0 
				\int_{\O  }  \| \phi_f \|_{W^{2,p}}^p \int_{\R^3} w_{\tilde{\vartheta}}^{-p/4}  \bigg\}\\
				\lesssim& \ (1+\|w_\vartheta f\|_\infty ) \int^t_0 \iint_{\O \times \R^3} |w_{\tilde{\vartheta}}\alpha_{f,\e}^\beta \p f|^p + \int^t_0 \iint_{\O \times \R^3} |f|^p,
			\end{split}
			\Ee
			where we have used, from (\ref{upper_vb}), $\alpha_{f,\e}(v)^\beta w_\vartheta(v)^{-1/2} \lesssim w_\vartheta(v)^{-1/4}$.

			\hide
			
			Note that $\Gamma_v (f,f)$ is defined in (\ref{Gamma_v}) and bounded by 
			\[
			\Gamma_{v}(f,f) \lesssim \frac{\langle v \rangle}{w(v)}\|wf\|_{\infty}^{2}
			\]
			from (\ref{Gvloss}) and (\ref{Gvgain}). Also we note that 
			\[
			| \Gamma(\p f, f) + \Gamma(f, \p f) | \lesssim \|wf\|_{\infty} \int_{\R^3} \mathbf{k}_{\varrho}(v,u )  |  \p f (u)| \dd u
			\dd v.
			\]

			We further apply Lemma \ref{lem_Green} to (\ref{Boltzmann_p}) and use (\ref{vKsum}), (\ref{bound_Gamma_nabla_vf1}), (\ref{bound_Gamma_nabla_vf2}), (\ref{Gvloss}), (\ref{Gvgain}), and (\ref{nabla_2_phi_p}) to obtain, together with (\ref{pf_p_1}) and (\ref{pf_p_2}),   
			\Be\begin{split}\label{pf_green}
				& \| \alpha^\beta 
				\p f (t)\|_p^p
				+ \int^t_0 \| \nu_{\phi_f}^{1/p} \alpha^\beta 
				\p f   \|_p^p  
				+ \underbrace{\int^t_0 |\alpha^\beta 
					\p f   |_{p,+}^p  }_{(\ref{pf_green})_{1}}
				\\
				\lesssim&  \  \| \alpha^\beta 
				\p f (0)\|_p^p   + 
				\underbrace{
					\int^t_0 |\alpha^\beta 
					\p f  |_{p,-}^p
				}_{(\ref{pf_green})_2}+  (1+ \sup_{0 \leq s\leq  t} \| \nabla^2 \phi_{f} (t) \|_\infty) \int^t_0 \| \alpha^\beta 
				\p f  \|_p^p + \int^t_0 \iint_{\O \times \R^3} |f|^p
				\\
				&+ (1 +\sup_{0 \leq s\leq  t} \|wf(s)\|_{\infty}) \int^t_0  \underbrace{ \iint_{\O\times\R^{3}}
					\alpha(v)^{p\beta}    |  \p f (v)|^{p-1}  
					\int_{\R^3} \mathbf{k}_{\varrho}(v,u )  |  \p f (u)| \dd u
					\dd v
					\dd x }_{(\ref{pf_green})_3} \dd s \\
				&+ (1 +\sup_{0 \leq s\leq  t} \|wf(s)\|_{\infty}) \int^t_0   \underbrace{ \iint_{\O\times\R^{3}} \alpha(v)^{\b p} |\p f(v)|^{p-1} 
					\int_{\R^3} \mathbf{k}_\varrho(v,u)  |f(u)| \dd u  
					\dd v \dd x }_{(\ref{pf_green})_4} \dd s.
			\end{split}\Ee
			\unhide
			
			\hide
			
			Estimates $(\ref{pf_green})_4$ and $(\ref{pf_green})_5$ are easily obtained by H\"older inequality,
			\Be \label{151_4}
			\begin{split}
				(\ref{pf_green})_4 &= \iint_{\O\times\R^{3}} \alpha^{\b p} |\p f|^{p-1} \langle v \rangle |f| \dd v \dd x = \iint_{\O\times\R^{3}} |\langle v \rangle\alpha^{\b}f| |\alpha^{\b}\p f|^{p-1} \dd v \dd x  \\
				&\lesssim \|\alpha^{\b}\p f\|_{p}^{p-1} \|\langle v \rangle\alpha^{\b}f\|_{p} \lesssim \|wf\|_{\infty} \|\alpha^{\b}\p f\|_{p}^{p-1},
			\end{split}
			\Ee  
			and
			\Be \label{151_5}
			\begin{split}
				(\ref{pf_green})_5  &\lesssim \iint_{\O\times\R^{3}} \alpha^{\b p} |\p f|^{p-1} \Big( \sqrt{\mu}^{1-\delta} + \frac{\langle v \rangle}{w(v)} \Big) \dd v \dd x  \\
				&\lesssim \|\alpha^{\b}\p f\|_{p}^{p-1} \Big\|\alpha^{\b}\Big( \sqrt{\mu}^{1-\delta} + \frac{\langle v \rangle}{w(v)} \Big) \Big\|_{p} \lesssim \|\alpha^{\b}\p f\|_{p}^{p-1}.
			\end{split}
			\Ee\unhide
			
			\vspace{4pt}

			\textit{Step 3.} We focus on (\ref{pf_p_K}). 
			With $N>0$, we split the $u$-integration of (\ref{pf_p_K}) into the integrations over $\{|u| \leq N\}$ and $\{|u| \geq N\}$.
			
			For $\{|u| \geq N\}$ and $0< \tilde{\varrho} \ll \varrho$, by Holder inequality with $\frac{1}{p} + \frac{1}{p^*}=1$
			\Be\begin{split}\label{est_k_p_f}
				&\int_{|u| \geq N}\alpha_{f,\e}^\beta(v)  \mathbf{k}_{\tilde{\varrho}} (v,u) |\p f(u)|\\
				\leq & \ 
				\alpha_{f,\e}^\beta (v)
				{\left(\int_{|u| \geq N}  \mathbf{k}_{\tilde{\varrho}}  (v,u) \frac{1}{\alpha_{f,\e}(u)^{\beta {p^*}}}   \right)^{1/{p^*}}}\\
				& \times 
				\left(
				\int_{|u| \geq N} \mathbf{k}_{\tilde{\varrho}}  (v,u)  |  \alpha_{f,\e}^{\beta  }\p f (u)|^p  
				\right)^{1/p}\\
				\lesssim & \ \alpha_{f,\e}^\beta(v) \left(
				\int_{|u| \geq N}  \mathbf{k}_{\tilde{\varrho}}  (v,u)  |\alpha_{f,\e}^{\beta  }\p f (u)|^p \dd u 
				\right)^{1/p},
			\end{split}\Ee
			where we have used Proposition \ref{prop_int_alpha} with $\beta q< \frac{p-1}{p}\frac{p}{p-1}=1$ from (\ref{beta_condition}).
			
			Then the contribution of $\{|u| \geq N\}$ in (\ref{pf_p_K}) is bounded by 
			\Be\begin{split}\label{small result}
				&\int^t_0 \int_{\O} \int_{v \in \R^3}  | \nu_{\phi_f}^{1/p} w_{\tilde{\vartheta}}\alpha_{f,\e}^\beta \p f (v)|^{p-1} \frac{\alpha_{f,\e}(v)^\beta}{ \nu_{\phi_f} (v)^{\frac{p-1}{p}}}\\
				& \ \ \ \   \ \ \ \ \ \ \ \ \times 
				\int_{|u| \geq N} \mathbf{k}_\varrho (v,u)\frac{w_{\tilde{\vartheta}}(v)}{w_{\tilde{\vartheta}}(u)} | w_{\tilde{\vartheta}} \p f(u)| \dd u \dd v \dd x \dd s  \\
				\leq& \ \int^t_0 \int_{\O} \bigg(\int_{v}  |\nu_{\phi_f}^{1/p}w_{\tilde{\vartheta}} \alpha_{f,\e}^\beta \p f (v)|^{p }\bigg)^{1/q}\\
				& \ \ \ \ \ \  \times 
				\bigg(
				\int_{|u| \geq N} |w_{\tilde{\vartheta}}\alpha_{f,\e}^{\beta  }\p f (u)|^p  \int_{v}\mathbf{k}_{\tilde{\varrho}}(v,u)  
				\bigg)^{1/p}\\
				\lesssim & \ 
				o(1)
				\int_0^t \|\nu_{\phi_f}^{1/p} w_{\tilde{\vartheta}} \alpha_{f,\e}^\beta \p f  (s)\|_p^p \dd s
				+\int_0^t \| w_{\tilde{\vartheta}} \alpha_{f,\e}^\beta \p f  (s)\|_p^p \dd s,
			\end{split}\Ee
			where we have used, from (\ref{upper_vb}), $\frac{\alpha_{f,\e}(v)^\beta}{ \nu_\phi (v)^{\frac{p-1}{p}}}\lesssim 1$ for $\beta$ in (\ref{beta_condition}), (\ref{k_vartheta_comparision}), and (\ref{grad_estimate}).

			The contribution of $\{|u| \leq N\}$ in (\ref{pf_p_K}) is bounded by, from the H\"older inequality, 
			\begin{eqnarray}
			&& \int^t_0  \int_\O 
			\int_{\R^3}
			| \nu_{\phi_f}^{1/p} w_{\tilde{\vartheta}}\alpha_{f,\e}^\beta \p f (v)|^{p-1} \notag\\
			&& \ \ \  \times 
			\int_{|u| \leq N} \mathbf{k}_\varrho(v,u) \frac{ w_{\tilde{\vartheta}}(v)}{ w_{\tilde{\vartheta}}(u)}\frac{ \alpha_{f,\e}(v)^{\beta} |  w_{\tilde{\vartheta}} \alpha_{f,\e}^\beta \p f (u)| }{\nu_{\phi_f}(v)^{(p-1)/p}  \alpha_{f,\e}(u)^{\beta}}     \dd u
			\dd v
			\dd x \dd s \notag\\
			&\leq& 
			\int^t_0    \| \nu_{\phi_f}^{1/p} w_{\tilde{\vartheta}}  \alpha_{f,\e}^\beta \p f  (s)\|_p^{p-1}\notag\\
			&& \ \ \   \times 
			\Big[
			\int_\O
			\int_{\R^3}
			\Big(
			\underline{\underline{
					\int_{|u| \leq N} 
					\mathbf{k}_{\tilde{\varrho}} (v,u) 
					\frac{ |
						w_{\tilde{\vartheta}} \alpha_{f,\e}^\beta \p f ( u)|}{ 
						\alpha_{f,\e}(u)^{\beta}}
					\dd u}}
			\Big)^p
			\dd v
			\dd x  \Big]^{1/p}
			\dd s 
			.\label{double_underline}
			\end{eqnarray}
			where we have used (\ref{k_vartheta_comparision}) and the fact $\alpha_{f,\e}^\beta/\nu_{\phi_f}^{\frac{p-1}{p}} \lesssim 1$ from (\ref{upper_vb}) and (\ref{beta_condition}).
			
			By the H\"older inequality, we bound an underlined $u$-integration inside (\ref{double_underline}) as
			\Be
			\|w_{\tilde{\vartheta}} \alpha_{f,\e}^\beta \p f(\cdot) \|_{L^p(\R^3)}
			\times 
			\bigg(
			\int_{\R^3}
			\frac{e^{-p^*\tilde{\varrho}  |v-u|^2}}{|v-u|^{p^*}} \frac{\mathbf{1}_{|u| \leq N}}{\alpha_{f,\e}(u)^{\beta p^*}}
			\dd u
			\bigg)^{1/q}
			\label{double_underline_split},
			\Ee
			where $1/p + 1/p^* =1$.

			\hide
			By the H\"older inequality, we bound it by 
			\begin{equation} \label{whole}
			\int^t_0  \int_\O 
			\int_{\R_{v}^3}
			| \nu_\phi^{1/p}\alpha^\beta \p f (v)|^{p-1} 
			\int_{\R_{u}^3} \mathbf{k}(v,u) \frac{ |  \alpha^\beta \p f (u)| }{ \alpha(u)^{\beta}}     \dd u
			\dd v
			\dd x \dd s.   
			\end{equation}
			We note that $\frac{p-2}{p} < \b < \frac{p-1}{p}$ implies $-1 < (\b-1) p + 1 < 0$ and therefore 
			$
			\frac{\alpha^{\b p}(v)}{\nu_{\phi}(v)^{p-1}} \lesssim \frac{\alpha^{\b p}(v)}{\langle v \rangle^{p-1}} \lesssim 1.   
			$

			Recall a standard estimate of $\mathbf{k}(v,u)\lesssim  \frac{e^{- C|v-u|^2}}{|v-u|}$ in (\ref{estimate_k}). 
			With $M>0$, we split $\{|u| \leq M\} \cup \{|u| \geq M\}$.
			In (\ref{whole}), for $\{|u| \geq M \}$,
			\Be \label{double_underline_split}
			\begin{split}
				&\int^t_0  \int_\O 
				\int_{\R_{v}^3}
				| \nu_\phi^{1/p}\alpha^\beta \p f (v)|^{p-1} 
				\int_{|u| \geq M} \mathbf{k}(v,u) \frac{ |  \alpha^\beta \p f (u)| }{ \alpha(u)^{\beta}}     \dd u
				\dd v
				\dd x \dd s  \\
				&\leq \int^t_0 \int_{\O} \int_{\R^{3}_{v}}  |\nu_{\phi}^{1/p} \alpha^\beta \p f (v)|^{p-1} \\
				&\quad\times 
				\underbrace{ \Big[ \int_{|u| \geq M} \mathbf{k}_\zeta (v,u) \frac{1}{\alpha^{\b q}} \Big] ^{1/q} }
				\Big[ \int_{|u| \geq M} \mathbf{k}_\zeta (v,u) |\alpha^{\b}\p f(u)|^{p} \Big]^{1/p}  dv dx ds \\
				&\leq \ \int^t_0 \int_{\O} \Big[ \int_{\R^{3}_{v}}  |\nu_{\phi}^{1/p}\alpha^\beta \p f (v)|^{p } \Big]^{1/q}
				\Big[
				\int_{|u| \geq M} |\alpha^{\beta  }\p f (u)|^p \int_{\R^{3}_{v}} \mathbf{k}_\zeta(v,u)  
				\Big]^{1/p}  dx  ds  \\
				&\lesssim \ \int_0^t \| \nu_{\phi}^{1/p}\alpha^\beta \p f  (s)\|_{p}^{p-1} \| \alpha^\beta \p f  (s)\|_{p} \dd s  \\
				&\lesssim \ O(\varepsilon)\int_0^t \| \nu_{\phi}^{1/p}\alpha^\beta \p f  (s)\|_{p}^{p} + \int_0^t \| \alpha^\beta \p f  (s)\|_{p}^{p}  \dd s  \\
			\end{split}
			\Ee
			where we used Proposition \ref{prop_int_alpha} for above underbraced term.  \\
			
			For $\{|u| \leq M\}$ part of (\ref{whole}),
			\Be  \label{small}
			\begin{split}
				&\int^t_0  \int_\O 
				\int_{\R_{v}^3}
				| \nu_\phi^{1/p}\alpha_{f,\e}^\beta \p f (v)|^{p-1} 
				\int_{|u| \leq N} \mathbf{k}(v,u) \frac{ |  \alpha_{f,\e}^\beta \p f (u)| }{ \alpha_{f,\e}(u)^{\beta}}     \dd u
				\dd v
				\dd x \dd s  \\
				&\lesssim \int^t_0    \| \nu_\phi^{1/p} \alpha_{f,\e}^\beta \p f  (s)\|_p^{ p/q}  
				\Big[
				\int_\O
				\|  \alpha_{f,\e}^\beta \p f(\cdot) \|^{p}_{L^p(\R^3)}
				\int_{\R_{v}^3}
				|(**)|^{p}
				\dd v
				\dd x  \Big]^{1/p}
				\dd s  \\
			\end{split}
			\Ee
			where  \unhide
			
			It is important to note that 
			\Be\label{convolution}
			\bigg(
			\int_{\R^3}
			\frac{e^{-p^*\tilde{\varrho}  |v-u|^2}}{|v-u|^{p^*}} \frac{\mathbf{1}_{|u| \leq N}}{\alpha_{f,\e}(u)^{\beta p^*}}
			\dd u
			\bigg)^{1/{p^*}}  \leq \bigg|
			\frac{1}{| \cdot |^{p^*}} *  \frac{\mathbf{1}_{|\cdot| \leq N}}{\alpha_{f,\e}(\cdot)^{p^*\beta}}
			\bigg| ^{1/{p^*}}.
			\Ee
			%
			By the Hardy-Littlewood-Sobolev inequality with $$1+ \frac{1}{p/{p^*}} = \frac{1}{3/{p^*}} + 
		\frac{1}{
				\frac{3}{2} \frac{p-1}{p}
			},$$ we have  
			\Be \label{153_1p}
			\begin{split}
				&\left\|
				\bigg|
				\frac{1}{| \cdot |^{p^*}} *  \frac{\mathbf{1}_{|\cdot| \leq N}}{\alpha_{f,\e}(\cdot)^{{p^*}\beta}}
				\bigg| ^{1/{p^*}}\right\|_{L^p(\R^3)} \\
				& = 
				\bigg\|
				\frac{1}{| \cdot |^{p^*}} *  \frac{\mathbf{1}_{|\cdot| \leq N}}{\alpha_{f,\e}(\cdot)^{p^*\beta}}
				\bigg\|_{L^{p/{p^*}} (\R^3)}  ^{1/{p^*}}\\
				&\lesssim 
				\left\| \frac{\mathbf{1}_{|\cdot| \leq N}}{\alpha_{f,\e}(\cdot)^{p^*\beta}}\right\|_{L^{
						\frac{3(p-1)}{2p}
					} (\R^3)}^{1/{p^*}}\\
				&\lesssim
				\left( \int_{\R^3} \frac{\mathbf{1}_{|v| \leq N}}{\alpha_{f,\e}(v) ^{\frac{p}{p-1} \beta 
						\frac{3(p-1)}{2p}
					}
				} \dd v\right)^{
					\frac{2p}{3(p-1)} \frac{p-1}{p}
				}
				\\
				&=
				\left( \int_{\R^3} \frac{\mathbf{1}_{|v| \leq N}}{\alpha_{f,\e}(v) ^{  3\beta/2
					}
				} \dd v\right)^{
					2/3
				}
				.
			\end{split}
			\Ee
			For $3<p<6$, we have $\frac{3}{2}\frac{p-2}{p} < 1$ and $\frac{2}{3} < \frac{p-1}{p}$. 
			Importantly from (\ref{beta_condition}) we have $\frac{3\beta}{2}<1$. Now we apply (\ref{NLL_split2}) in Proposition \ref{prop_int_alpha} to conclude that 
			\Be\notag
			\left( \int_{\R^3} \frac{\mathbf{1}_{|v| \leq M}}{\alpha_{f,\e}(v) ^{ 3\beta/2}
			} \dd v\right)^{2/3}\lesssim_{p, \beta, M,\O} 1 .
			\Ee 
			Finally from (\ref{double_underline}), (\ref{double_underline_split}), (\ref{convolution}), (\ref{153_1p}), and (\ref{small result}) we bound 
			\Be\label{151_3}\begin{split}
			(\ref{pf_p_K}) \lesssim & \ 
			o(1)\int_{0}^{t}  \| \nu_\phi^{1/p} w_{\tilde{\vartheta}}\alpha_{f,\e}^\beta \p f   \|_{p}^{p}\\
			&  + (1+ \sup_{0 \leq s \leq t } \| w_\vartheta f(s) \|_\infty) \int_{0}^{t} \|w_{\tilde{\vartheta}}\alpha_{f,\e}^{\b}\p f \|_{p}^{p}   .
			\end{split}\Ee
			
			Collecting terms from (\ref{pf_p_nabla}), (\ref{pf_p_K}), (\ref{pf_p_0}), (\ref{pf_p_1}), (\ref{pf_p_Gamma_v_-}), (\ref{pf_p_Kv}), (\ref{pf_p_2}), (\ref{small result}), and (\ref{151_3})
			\Be
			\begin{split}\label{final_est_G}
				(\ref{pf_green})_\mathcal{G}
				\lesssim& \ \Big(o(1) + \sup_{0 \leq s\leq t}\| w_\vartheta f(s) \|_\infty\Big)
				\int^t_0  \| \nu_{\phi_f}^{1/p} w_{\tilde{\vartheta}}  \alpha_{f,\e}^\beta \p f \|_p^p\\
				&+  (1 + \sup_{0 \leq s\leq t}\| w_\vartheta f(s) \|_\infty + \sup_{0 \leq s\leq t}\| \nabla^2 \phi_f(s) \|_\infty
				)\\
				&  \ \ \ \times 
				\int^t_0   \| w_{\tilde{\vartheta}}  \alpha_{f,\e}^\beta \p f \|_p^p 
				\\
				& 
				+(1 + \sup_{0 \leq s\leq t}\| w_\vartheta f(s) \|_\infty+ \delta_1/\Lambda_1)
				\int^t_0 \|w_{\tilde{\vartheta}} f \|_p^p  .
			\end{split}
			\Ee
			
			\hide
			Therefore
			\Be \label{small result}
			\begin{split}
				(\ref{small}) &\lesssim \int^t_0    \| \nu_\phi^{1/p} \alpha^\beta \p f  (s)\|_p^{ p/q}
				\|  \alpha^\beta \p f(\cdot) \|_{p}	\dd s  \\
				&\lesssim O(\varepsilon)\int_{0}^{t}  \| \nu_\phi^{1/p} \alpha^\beta \p f  (s)\|_{p}^{p} ds + \int_{0}^{t} \|\alpha^{\b}\p f(s)\|_{p}^{p} ds.
			\end{split}
			\Ee

			From (\ref{small result}) and (\ref{double_underline_split}) we conclude
			\Be \label{151_3}
			\begin{split}
				\int_{0}^{t} (\ref{pf_green})_{3} \dd s &\lesssim  o(1)\int_{0}^{t}  \| \nu_\phi^{1/p} \alpha^\beta \p f  (s)\|_{p}^{p} \dd s + \int_{0}^{t} \|\alpha^{\b}\p f(s)\|_{p}^{p} \dd s .
			\end{split}
			\Ee

			\unhide

			\vspace{4pt}

			\textit{Step 4.}
			We focus on $(\ref{pf_green})_{\gamma_-}$. From (\ref{BC_deriv}) and (\ref{BC_deriv_R}),
			\Be \label{BC_deriv_1}
			\begin{split}
				&
				\int_{n(x) \cdot v<0}
				|n(x) \cdot v|^{\beta p}
				|w_{\tilde{\vartheta}}  \nabla_{x,v} f(t,x,v) |^p
				|n(x) \cdot v| \dd v
				\\
				&\lesssim  \int_{n(x) \cdot v<0}
				\langle v\rangle^p \mu(v)^{\frac{p}{2}}w_{\tilde{\vartheta}} ^p
				 \\
				& \  \ \ \ \ \ \ \ \ \ \ \   \   \times\Big(|n(x) \cdot v|^{\beta p+1}+ |n(x) \cdot v|^{(\beta-1) p+1} 
				\Big) |(\ref{BC_deriv_R})|^p  \dd v.
			\end{split}\Ee
			Note that for $0< \tilde{\vartheta} \ll_p 1$ we have $\mu(v)^{\frac{p}{2}}w_{\tilde{\vartheta}} ^p\lesssim e^{C|v|^{2}}$ for some $C>0$ when $|v|\gg1$. 
			
			On the other hand, from (\ref{beta_condition}), we have   
			\Be\label{L1_loc}
			(\beta-1) p + 1> \frac{p-2}{p} p - p+1 = -1, \ \ \  |n(x) \cdot v|^{(\beta-1) p+1} \in L^1_{loc}(\R^3).
			\Ee
			
			\hide
			\Be\notag
			\begin{split}
				&|\nabla_{x,v} f(t,x,v) |\\
				\lesssim& \ 
				\langle v\rangle\sqrt{\mu(v)}
				\Big(1+ \frac{1}
				{|n(x) \cdot v|}  \Big)\\
				& \times
				\int_{n(x) \cdot u>0}  \Big\{
				(
				\langle u\rangle  + |\nabla_x \phi_f|
				)
				|\nabla_{x,v} f(t,x,u) | \\ &  \ \ \ \ \   \ \ \  \ \ \  \ \ \  \ \ \   + 
				\langle u \rangle (1+ |\nabla_x \phi_f |) |f|+ |\nabla_x \phi_f| \mu(u)^{\frac{1}{4}}
				\\
				&  \ \ \ \ \   \ \ \  \ \ \  \ \ \  \ \ \  
				+(1+ \| w f \|_\infty) |\int_{\R^3}\mathbf{k}_{\varrho }(u, u^\prime)|f(u^\prime)| \dd u^\prime|
				\Big\}\sqrt{\mu(u)}\{n(x) \cdot u\} \dd u.
			\end{split}\Ee   
			\unhide

			Now we bound $|(\ref{BC_deriv_R})|^p$. For the first line of (\ref{BC_deriv_R}), 
			%
			we split the $u$-integration into $\gamma_+^\e(x) \cup \gamma_+ (x) \backslash \gamma_+^\e(x)$ where $\e$ is small but satisfies (\ref{lower_bound_e}). By the H\"older inequality 
			\Be
			\begin{split}
				& \bigg\{\int_{n(x) \cdot u>0}  |w_{\tilde{\vartheta}}\alpha_{f,\e}^\beta\nabla_{x,v} f(s,x,u)| \{w_{\tilde{\vartheta}}\alpha_{f,\e}( u)\}^{- \beta} \\
				& \ \ \ \ \ \ \ \ \ \ \ \times \langle u\rangle \sqrt{\mu(u)} \{n(x) \cdot u\}  \dd u\bigg\}^p
				\label{W1p_bdry}
				\\
				&\lesssim 
				\bigg\{\int_{\gamma_+^\e (x)}  |w_{\tilde{\vartheta}}\alpha_{f,\e}^\beta\nabla_{x,v} f(s,x,u)|^p\{n(x) \cdot u\}  \dd u\bigg\}  \\
				& \ \ \ \  \times 
				 { \bigg\{ \int_{\gamma_+^\e (x)}
					\{w_{\tilde{\vartheta}}\alpha_{f,\e}( u)\}^{- \beta p^*}
					|n(x) \cdot u|  \mu ^{\frac{q}{4}}\dd u 
					\bigg\}^{p/{p^{*} }} }  \\
				& +\bigg\{\int_{\gamma_+(x) \backslash\gamma_+^\e (x)}  | w_{\tilde{\vartheta}}\alpha_{f,\e}^\beta\nabla_{x,v} f(s,x,u)|^p
				\mu ^{\frac{p}{8}}
				\{n(x) \cdot u\}  \dd u\bigg\}  \\
				& \ \ \ \   \times  { \bigg\{ \int_{\gamma_+(x) \backslash\gamma_+^\e (x)}
					\{w_{\tilde{\vartheta}}\alpha_{f,\e}(s,x,u)\}^{- \beta p^* }
					|n(x) \cdot u|\mu ^{\frac{p^{*} }{8}}\dd u 
					\bigg\}^{p/{p^{*} }} },  
			\end{split}
			\Ee
			for $p^{*} := \frac{p}{p-1}$. Note that $\alpha_{f,\e}(s,x,u)\neq |n(x) \cdot u|$ for $(x,u) \in \gamma_+$ in general. From (\ref{beta_condition}), $\beta p^*<1$. From (\ref{NLL_split2}) and (\ref{NLL_split3}) with $v=0$, we have $\alpha_{f,\e}^{- \beta p^*} |n(x) \cdot u |   \lesssim \alpha_{f,\e}^{- \beta p^*}   \in L^1_{loc} (\{ u \in \R^3\})$. Since $\mathbf{1}_{\gamma^\e_+ (x)}(v) \downarrow 0$ almost everywhere in $\R^3$ as $\e \downarrow 0$, by the dominant convergence theorem, for (\ref{delta_1/lamdab_1}), we choose $\e : = \frac{2 \delta_1}{\Lambda_1}\ll_\O 1$  
			\Be\begin{split}\label{W1p_bdry_1}
				(\ref{W1p_bdry}) \lesssim& \  o(1)  \int_{\gamma_+^\e (x)}  |w_{\tilde{\vartheta}}\alpha_{f,\e}^\beta\nabla_{x,v} f(s,x,u)|^p\{n(x) \cdot u\}  \dd u\\
				&+\int_{\gamma_+(x) \backslash\gamma_+^\e (x)}  |w_{\tilde{\vartheta}}\alpha_{f,\e}^\beta\nabla_{x,v} f(s,x,u)|^p \mu(u)^{{p}/{8}}\{n(x) \cdot u\}  \dd u.
			\end{split} 
			\Ee

			From Lemma \ref{le:ukai} and (\ref{Boltzmann_p}), the last term of (\ref{W1p_bdry_1}) has a bound as 
			\Be\begin{split}\label{nongrazing_nabla_f}
				&\int^t_0 \int_{\p\O} \int_{\gamma_+(x) \backslash\gamma_+^\e (x)}  | w_{\tilde{\vartheta}}\alpha_{f,\e}^\beta\nabla_{x,v} f(s,x,u)|^p \mu(u)^{{p}/{8}}\{n(x) \cdot u\}  \dd u \dd S_x \dd s\\
				 &\lesssim \ \| w_{\tilde{\vartheta}} \alpha_{f,\e}^\beta\nabla_{x,v} f(0)  \mu^{{1}/{8}}\|^p_p + \int^t_0 \| w_{\tilde{\vartheta}} \alpha_{f,\e}^\beta \nabla_{x,v} f \|_p^p + (\ref{8 Last}),
			\end{split} 
			\Ee
			where, from (\ref{eqtn_nabla_f}), (\ref{mathcal_G}),
			%
			%
			\begin{eqnarray} 
			&&
			\int^t_0\iint_{ \O \times \R^3}
			\big| \p_t + v\cdot \nabla_x - \nabla_x \phi_f \cdot \nabla_v + \nu_{\phi_f,w_{\tilde{\vartheta}}}  |\nonumber\\
			&& \ \ \ \ \ \ \ \ \ \ \ \ \  \times 
			w_{\tilde{\vartheta}}\mu^{1/8}
			\alpha_{f,\e}^\beta|\nabla_{x,v} f 
			|^p \label{8 Last} \\
			&\leq& \int^t_0
			\iint_{ \O \times \R^3} p\alpha_{f,\e}^{\beta p}  |\nabla_{x,v} f|^{p-1}|w_{\tilde{\vartheta}}  \mu^{1/8}|^p
			\big|
			\mathcal{G} \big|\label{8 Last_1} \\
			&& + \int^t_0\iint_{ \O \times \R^3} 
			|\nabla_x \phi_f| \mu^{0+} 
			|\alpha_{f,\e}^\beta\nabla_{x,v} f |^p .  \label{8 Last_2}
			\end{eqnarray}
			\hide
			&\lesssim 
			\int^t_0\iint_{ \O \times \R^3} p\alpha^{\beta p}  |\nabla_{x,v} f|^{p-1}  \mu^{p/8} |\mathcal{G}|  
			+ ( 1 + \| \nabla \phi_{f}\|_{\infty} ) \|\alpha^{\b}\p f\|_{p}^{p}  \\
			&\lesssim 
			(o(1) + \sup_{0 \leq s\leq t}\| w_\vartheta f(s) \|_\infty)
			\int^t_0  \| \nu_{\phi_f}^{1/p} w_{\tilde{\vartheta}}\alpha^\beta \p f \|_p^p\\
			& \  \ +  (1 + \sup_{0 \leq s\leq t}\| w_{\vartheta} f(s) \|_\infty + \sup_{0 \leq s\leq t}\| \nabla^2 \phi_f(s) \|_\infty
			)
			\int^t_0   \| w_{\tilde{\vartheta}} \alpha_{f,\e}^\beta \p f \|_p^p 
			\\
			&  \ \ \  
			+(1 + \sup_{0 \leq s\leq t}\| w f(s) \|_\infty+ \delta_1/\Lambda_1) \int^t_0  \|  w_{\tilde{\vartheta}}  f \|_p^p.
		\end{eqnarray}\unhide
		Clearly $(\ref{8 Last_1})\lesssim (\ref{final_est_G})$. And, from (\ref{integrable_nabla_phi_f}),  
		\Be\notag
		(\ref{8 Last_2}) \lesssim \delta_1 \int^t_0 \| w_{\tilde{\vartheta}} \alpha_{f, \e}^\beta \nabla_{x,v} f \|_p^p.
		\Ee

		Now we consider the third term of (\ref{BC_deriv_R}). From the trace theorem $W^{1,p} (\O) \rightarrow W^{1- \frac{1}{p}, p} (\p\O)$ and (\ref{nabla_2_phi_p})
		\Be \label{phi_bdry_estimate}
		\| \nabla \phi_f \|_{L^p (\p\O)}\lesssim \| \nabla \phi _f \|_{W^{1- \frac{1}{p}, p} (\p\O)}\lesssim \| \nabla \phi_f \|_{W^{1,p} (\O)}\lesssim \| w_{\tilde{\vartheta}} f \|_{L^p (\O \times \R^3)}.
		\Ee
		Then 
		\Be
		\begin{split}\label{W1p_bdry_2}
			&\Big\{
			\int_{\p\O}\int_{n(x) \cdot u>0} 
			( \langle u\rangle \| w_\vartheta f\|_\infty  + \mu(u)^{\frac{1}{4}} )|\nabla \phi_f|  \sqrt{\mu(u)} \{n(x) \cdot u\} \dd u \dd S_{x}
			\Big\}^p\\
		&	\lesssim    (1+ \| w_\vartheta f \|_\infty) \| \nabla \phi_f \|_{L^p(\p\O)}^p\\
		&	\lesssim   (1+ \| w_\vartheta f \|_\infty)  \| w_{\tilde{\vartheta}} f\|^{p}_{L^p (\O \times \R^3)}.
		\end{split}
		\Ee
		
		For the second term of (\ref{BC_deriv_R}), by the H\"older inequality with $\frac{1}{p}+ \frac{1}{q}=1$ for $3<p<6$
		\Be
		\begin{split}\label{W1p_bdry_3}
			&\bigg\{
			\int_{n \cdot u>0} 
			\langle u\rangle |f| + (1+ \| w_{\vartheta} f \|_\infty) \int_{\R^3}\frac{ \mathbf{k}_\varrho (u,u^\prime)^{1/q}}{|n\cdot u^\prime|^{1/p} }
			 \mathbf{k}_\varrho (u,u^\prime)^{1/p}|f(u^\prime)|\\
		& \ \ \ \ \ \ \ \ \ \ \ \   \ \ \ \ \ \ \ \ \ \ \ \   \ \ \ \ \ \ \ \ \ \ \ \  \ \ \ \ 	 \  \times  |n\cdot u^\prime|^{1/p} \dd u^\prime
			\sqrt{\mu } \{n  \cdot u\} \dd u 
			\bigg\}^p\\
			\lesssim & \ \int_{n\cdot u>0} |f|^p \{n\cdot u\} \dd u\\
			&  + (1+ \| w_{\vartheta} f \|_\infty) \Big(\int_{\R^3} \mathbf{k}_{\varrho} (u,u^\prime) |n \cdot u^\prime|^{-q/p} \dd u^\prime\Big)^{p/q}\\
			& \ \ \ \ \ \ \ \ \ \ \ \  \ \ \ \ \ \ \ \ \  \times  \int_{\R^3}\int_{\R^3} \mathbf{k}_\varrho (u,u^\prime) |f(u^\prime)|^p |n \cdot u^\prime| \dd u^\prime \dd u\\
			\lesssim & \ (1+ \| w_{\vartheta} f \|_\infty) \int_{n\cdot u>0} |f|^p \{n\cdot u\} \dd u.
		\end{split}
		\Ee
		
		Collecting terms from (\ref{BC_deriv_1}), (\ref{W1p_bdry_1}), (\ref{8 Last}), (\ref{W1p_bdry_2}), and (\ref{W1p_bdry_3}) we derive that 
		\Be \label{mid} \begin{split} 
			&(\ref{pf_green})_{\gamma_-}\\
			\lesssim & \
			\|  w_{\tilde{\vartheta}}\alpha_{f,\e}^\beta\nabla_{x,v} f(0)  \mu(u)^{{1}/{8}}\|^p_p
			\\
			& + o(1) \int^t_0 | w_{\tilde{\vartheta}}\alpha_{f,\e}^\beta \p f |_{p,+}^p\\ 
			&+ (o(1) + \sup_{0 \leq s\leq t}\| w_{\vartheta} f(s) \|_\infty)
			\int^t_0  \| \nu_{\phi_f}^{1/p} w_{\tilde{\vartheta}}\alpha_{f,\e}^\beta \p f \|_p^p\\
			& + (1 + \sup_{0 \leq s\leq t}\| w_{\vartheta} f(s) \|_\infty+\sup_{0 \leq s\leq t}\| \nabla^2 \phi_f(s) \|_\infty   ) \int^t_0  \|   \alpha_{f,\e}^\beta w_{\tilde{\vartheta}} \p f \|_p^p 
			\\
			&    
			+(1 + \sup_{0 \leq s\leq t}\| w_{\vartheta} f(s) \|_\infty) \int^t_0  
			\big( \| w_{\tilde{\vartheta}} f \|_p^p+   | w_{\tilde{\vartheta}}f|_{p,+}^p \big) .
		\end{split} \Ee

		\hide

		Using Lemma \ref{le:ukai} and (\ref{8 Last}) to (\ref{mid}),
		\Be \label{151_2}
		\begin{split}
			& (\ref{pf_green})_2 = \int_{0}^{t} |\alpha^\beta 
			\p f(s)  |_{p,-}^p ds   \\
			&\lesssim O(\e) (1 + \|\phi_{f}\|_{C^{2}})^{p} \int^t_0 |\alpha^\beta \p f (s) |_{p,+}^p ds  + O(\varepsilon) (1+\|wf\|_{\infty}) \int_{0}^{t} \| \nu_\phi^{1/p} \alpha^\beta \p f  (s)\|_{p}^{p}  \\
			&\quad + \|\alpha^{\b}\p f (0) \|_{p}^{p} + \int_{0}^{t} (1 + \|wf\|_{\infty})^{2p}(1 + \|\phi_{f}\|_{C^{2}})^{p} ( 1 + \|\alpha^{\b}\p f\|_{p}) \|\alpha^{\b}\p f\|_{p}^{p-1}  ds \\
		\end{split}
		\Ee
		\unhide
		
		\vspace{4pt}

		\textit{Step 5.} From (\ref{pf_green}), (\ref{final_est_G}), (\ref{mid}) we have 
		\Be\begin{split}\label{W1p_pf}
			&\| w_{\tilde{\vartheta}} \alpha_{f,\e}^\beta \p f(t) \|_p^p  +  \int^t_0\|  \nu_{\phi_f, w_{\tilde{\vartheta}}}^{1/p} w_{\tilde{\vartheta}}\alpha_{f,\e}^\beta \p f(t) \|_{p,+}^p
			+ \int^t_0 | w_{\tilde{\vartheta}}\alpha_{f,\e}^\beta \p f|_{p}^p\\
			\lesssim & \  \| w_{\tilde{\vartheta}} \alpha_{f,\e}^\beta \p f(0) \|_p^p+(1 + \sup_{0 \leq s \leq t} \| w f(s) \|_\infty)\int^t_0 \big( \|w_{\tilde{\vartheta}} f \|_p^p + | w_{\tilde{\vartheta}}f|_{p,+}^p \big)   \\
			&+ (1 + \sup_{0 \leq s \leq t} \| w f(s) \|_\infty
			+ \sup_{0 \leq s \leq t} \| \nabla^2 \phi_f(s) \|_\infty
			)\int^t_0\| w_{\tilde{\vartheta}} \alpha_{f,\e}^\beta \p f(t) \|_p^p.
		\end{split}\Ee
		
		Multiplying small number to (\ref{W1p_pf}) and adding to (\ref{Lp_estimate_f}) we derive that 
		\Be\label{final_est_W1p}
		\begin{split}
			&\| w_{\tilde{\vartheta}} f(t) \|_p^p+\| w_{\tilde{\vartheta}}\alpha_{f,\e}^\beta \p f(t) \|_p^p\\
			&  + \int^t_0 \big( \| \nu_{\phi_f,w_{\tilde{\vartheta}}}^{1/p}  w_{\tilde{\vartheta}} f \|_p^p  + \|  \nu_{\phi_f,w_{\tilde{\vartheta}}}^{1/p} w_{\tilde{\vartheta}}\alpha_{f,\e}^\beta \p f  \|_p^p  \big) \\
			&+ \int^t_0 | w_{\tilde{\vartheta}}  f|_{p,+}^p + \int^t_0 |w_{\tilde{\vartheta}} \alpha_{f,\e}^\beta \p f|_{p,+}^p \\
			\lesssim & \  \| w_{\tilde{\vartheta}} f(0) \|_p^p+\|  w_{\tilde{\vartheta}}\alpha_{f,\e}^\beta \p f(0) \|_p^p+ 
			\int^t_0
			\| w_{\tilde{\vartheta}} f \|_p^p
			\\& 
			+ (1 
			+ \sup_{0 \leq s \leq t} \| \nabla^2 \phi_f(s) \|_\infty) \int^t_0 
			\| w_{\tilde{\vartheta}}  \alpha_{f,\e}^\beta \p f \|_p^p,
		\end{split}\Ee
		where we have used (\ref{small_bound_f_infty}).  Then by the Gronwall's inequality, we deduce (\ref{est_W}).  \hide

		\vspace{4pt}

		\textit{Step 5.} By simple modification we can prove 
		\Be\label{tilde_w_W1p}
		\| \tilde{w} \alpha^\beta \nabla_{x,v} f(t) \|_{L^p (\O \times \R^3)}^p + \int^t_0 \| \langle v \rangle^{1/p}\tilde{w}  \alpha^\beta \nabla_{x,v} f(t) \|_{L^p (\O \times \R^3)}^p  \lesssim_t 1.
		\Ee
		The equation for $\tilde{w} \p f$ equals 
		\Be\label{eqtn_nabla_f_tilde_w}
		[\p_t + v\cdot \nabla_x - \nabla_x \phi_f \cdot \nabla_v + \tilde{\nu}_{\phi_f,\tilde{w}}] \p f 
		=   \tilde{w} \mathcal{G},
		\Ee
		where ${\nu}_{\phi_g,\tilde{w}}$ in (\ref{nu_w}) and $\mathcal{G}$ in (\ref{mathcal_G}). Note that, from {\color{red}(??)}, for some $0<\tilde{\varrho}\ll \varrho$
		\Be \label{K_tilde_w}
		\tilde{w} (v)| K \p f(v)| 
		\lesssim  \int_{\R^3} \mathbf{k}_\varrho (v,u) \frac{\tilde{w}(v)}{\tilde{w}(u)} | \tilde{w}\p f (u)| \dd u 
		\lesssim  \int_{\R^3} \mathbf{k}_{\tilde{\varrho}} (v,u)   | \tilde{w}\p f (u)| \dd u.
		\Ee
		Similarly we have 
		\Be\begin{split}\label{Gamma_tilde_w}
			&\tilde{w} (v) |\Gamma( \tilde{w}^{-1} \|w f\|_\infty,\p f)|\lesssim  \tilde{w} (v) |\Gamma(w^{-1}\|w f\|_\infty, \tilde{w}^{-1} \tilde{w}\p f)|\\
			\lesssim & \ 
			\|w f\|_\infty \int_{\R^3} \mathbf{k}_{\tilde{\varrho}} (v,u)   | \tilde{w}\p f (u)| \dd u.
		\end{split}\Ee
		The rest of proof is same as \textit{Step 1}-\textit{Step 4}.\hide
		
		We define
		\Be \label{def W}
		\mathcal{W}(t) := 
		\| f(t) \|_p^p
		+\| \alpha^\beta 
		\p f (t)\|_p^p  .
		\Ee
		From (\ref{final_est_W1p}),
		\Be\label{eqtn_W}
		\mathcal{W} (t) \lesssim \mathcal{W} (0)  + (1 
		+ \sup_{0 \leq s \leq t} \| \nabla^2 \phi_f(s) \|_\infty)
		\int^t_0 \mathcal{W} (s) \dd s.
		\Ee

		\vspace{4pt}
		
		\textit{Step 5. }

		and assume $\|\alpha^{\b}\p f(t)\|_{p} \geq 1$ WLOG, since we apply growth estimate. We use (\ref{151_2}), (\ref{151_3}), (\ref{151_4}), and (\ref{151_5}) to (\ref{pf_green}) with sufficiently small $\varepsilon \ll 1$ to obtain
		\Be
		\begin{split}
			\mathcal{W}(t) &\lesssim \| \alpha^\beta \p f (0)\|_{p}^{p} + \int_{0}^{t}\mathcal{P}\big( \|wf\|_{\infty},  \|\phi_{f}\|_{C^{2}} \big) \mathcal{W}(s) ds, 
		\end{split}	
		\Ee
		and Gronwall's inequality yield
		\Be \label{e_growth}
		\mathcal{W}(t) := \| \alpha^\beta 
		\p f (t)\|_p^p
		+ \int^t_0 \| \nu_{\phi_f}^{1/p} \alpha^\beta 
		\p f (s) \|_p^p  ds \lesssim \| \alpha^\beta \p f (0)\|_{p}^{p} \ e^{ t C_{p}(f,\phi_{f}) },
		\Ee
		where 
		\[
		C_{p}(f,\phi_{f}) = \mathcal{P}\big( \|wf\|_{\infty},  \|\phi_{f}\|_{C^{2}} \big) = (1 + \|wf\|_{\infty})^{2p}(1 + \|\phi_{f}\|_{C^{2}})^{p}.
		\]
		
		From Schauder estimate,
		\Be
		\begin{split}\label{apply_schauder}
			\|\phi_{f}\|_{C^{2,1-\frac{3}{p}}} &\lesssim_{p,\O} \Big\|\int_{\R^{3}} f\sqrt{\mu} dv \Big\|_{C^{0,1-\frac{n}{p}}(\O)}  \\
			&\lesssim \|wf\|_{\infty} + \Big\|\int_{\R^{3}} \nabla_{x} f\sqrt{\mu} dv \Big\|_{L^{p}(\O)}  \\
			&\lesssim  \|wf\|_{\infty} + \Big( \int_{\R^{3}} \Big| \frac{\sqrt{\mu}}{\alpha^{\b}(x,v)} \Big|^{\frac{p}{p-1}} dv \Big)^{\frac{p-1}{p}} \|\alpha^{\b}\nabla_{x}f\|_{p}  \\
			&\lesssim e^{t C_{p}(f,\phi_{f})},
		\end{split}
		\Ee
		since we can apply Proposition \ref{prop_int_alpha} with $\b \frac{p}{p-1} < 1$. Now, let us choose 
		\[
		\delta_{2} = 1-\frac{3}{p}, \quad \lambda_{0} > \frac{1}{\delta_{2}} C_{p}(f,\phi_{f}).
		\]
		From (\ref{morrey}), $\|\nabla_{x}\phi_{f}\|_{C^{0,1-\frac{3}{p_{1}}}} \lesssim \|wf\|_{\infty} \lesssim e^{-\lambda t}$ holds for any $p_{1} > 3$. We choose $p_{1} \gg 3$ so that $\frac{3}{p_{1}}\lambda_{0} < \lambda$. Finally we apply (\ref{phi_interpolation}) to conclude
		\Be\begin{split}
			\|\nabla^2_x \phi(t )\|_{L^\infty_x}
			&\leq
			e^{\frac{3}{p_{1}}\lambda_0t}[\nabla_x \phi(t)]_{C^{0,1-\frac{3}{p_{1}}}_x}
			+ e^{-(1-\frac{3}{p}) \lambda_0t}[\nabla^2 \phi(t)]_{C_x^{0, 1-\frac{3}{p}}}  \\
			&\lesssim e^{ -\lambda_{1} t } ,
		\end{split}\Ee
		where
		\[
		\lambda_{1} := \min\{ \lambda-\frac{3}{p}\lambda_{0}, (1-\frac{3}{p})\lambda_{0} - C_{p}(f,\phi_{f}) \}.	
		\]
		\unhide\unhide
	\end{proof}

	\section{$L^3_xL^{1+}_v$-Estimate of $\nabla_v f$ and $L^{1+}$-Stability}
	\hide In this section we prove the uniqueness statement assuming the extra condition (\ref{extra_con_uniqueness}). We need the following lemma.
	\begin{lemma}
		Now we claim that for $0<\kappa\leq2$ and $0< \sigma<1$
		\Be\label{nonlocal_finite}
		\sup_{t,x}\int_{\R^3}
		\frac{e^{-C|v-u|}}{|v-u|^{2-\kappa}} \frac{\mathbf{1}_{\tb(t,x,v) \leq 2}}{|n(\xb(t,x,v)) \cdot \vb(t,x,v)|^\sigma} \dd u \lesssim 1.
		\Ee
	\end{lemma}\unhide

	\hide In this section we prove the uniqueness statement assuming the extra condition (\ref{extra_con_uniqueness}). We need the following lemma.
	\begin{lemma}
		Now we claim that for $0<\kappa\leq2$ and $0< \sigma<1$
		\Be\label{nonlocal_finite}
		\sup_{t,x}\int_{\R^3}
		\frac{e^{-C|v-u|}}{|v-u|^{2-\kappa}} \frac{\mathbf{1}_{\tb(t,x,v) \leq 2}}{|n(\xb(t,x,v)) \cdot \vb(t,x,v)|^\sigma} \dd u \lesssim 1.
		\Ee
	\end{lemma}\unhide
	
	\begin{proposition}\label{prop_better_f_v}Assume $f$ and $\phi_{f}$ solve (\ref{eqtn_f}), (\ref{phi_f}), (\ref{BC_f}), and satisfy estimates (\ref{main_Linfty}), (\ref{W1p_main}), (\ref{integrable_nabla_phi_f}) with the condition (\ref{delta_1/lamdab_1}). We also assume extra initial condition
		\Be \label{Extra_uniq}
		\| w_{\tilde{\vartheta}} \nabla_v f_{0} \|_{L^3_{x,v}} < \infty.
		\Ee 
		Then 
		\Be\label{bound_nabla_v_g_global}
		\| \nabla_v f(t) \|_{L^3_x (\O)L^{1+\delta}_v(\R^3)} \lesssim_t 1 \ \ \text{for all } \ t\geq 0.
		\Ee
		%
		%
		%
	\end{proposition}
	
	Once we have Proposition \ref{prop_better_f_v}, we can prove the following stability result.
	\begin{proposition}\label{L1+stability}Suppose $f$ and $g$ solve (\ref{eqtn_f}), (\ref{phi_f}), (\ref{BC_f}), and satisfy (\ref{main_Linfty}). Also we assume $f,g$ satisfy (\ref{bound_nabla_v_g_global}). Then
		\Be\label{1+delta_stability}\begin{split}
	&	\| f(t) - g(t) \|_{L^{1+\delta}(\O \times \R^3)}^{1+\delta}
		+ \int^t_0 \| f(s) - g(s) \|_{L^{1+\delta}(\gamma)} ^{1+\delta}\dd s\\
&		\lesssim_t \| f_0 - g_0 \|_{L^{1+\delta}(\O \times \R^3)}^{1+\delta}.
		\end{split}\Ee
	\end{proposition}
	\begin{proof}
		Assume that $f$ and $g$ solve (\ref{eqtn_nabla_f}). %
		Then 
		\begin{equation}\label{eqtn_f-g}
		\begin{split}
		& \p_t [f-g] + v\cdot \nabla_x [f-g] - \nabla_x  \phi_f  \cdot \nabla_v [f-g]\\
		&+ \frac{v}{2} \cdot \nabla_x \phi_f[f-g]+  \nu [f-g] \\
		= & \nabla_x \phi_{f-g} \cdot \nabla_v g\\
		& + K[f-g]  - \frac{v}{2} \cdot \nabla_x \phi_{f-g} g +
		\Gamma(f,f) -\Gamma(g,g) - v\cdot \nabla_x \phi_{f-g} \sqrt{\mu} .
		\end{split}\end{equation} 
		By Lemma \ref{lem_Green} for $L^{1+ \delta}$-space with $0 < \delta \ll 1$, we obtain 
		\Be\begin{split}\label{f-g_energy}
			&\| [f- g](t)\|_{1+ \delta}^{1+ \delta} + \int^t_0 \| \nu_{\phi_f}^{ {1}/{1+ \delta}} [f-g] \|_{1+ \delta}^{1+ \delta}  + \int_0^t |[f-g]|_{1+ \delta, + }^{1 + \delta} \\
			\leq & \ \| [f- g](0)\|_{1+ \delta}^{1+ \delta}+
			\int^t_0 \iint_{\O \times \R^3}
			|\text{RHS of } (\ref{eqtn_f-g}) | |f-g|^{\delta}\\
			&
			+\int_0^t |[f-g]|_{1+ \delta, - }^{1 + \delta} .
		\end{split}
		\Ee

		
		\textit{Step 1. } For $0 < \delta \ll 1$, by the H\"older inequality with $1=\frac{1}{\frac{3(1+ \delta)}{2- \delta}} + \frac{1}{3} +  \frac{1}{\frac{1+ \delta}{\delta}}$ and the Sobolev embedding $W^{1, 1+ \delta} (\O)\subset L^{\frac{3(1+ \delta)}{2- \delta}}(\O)$ when $\O \subset \R^3$,
		\Be\label{gronwall_f-g}
		\begin{split}
			&\int^t_0 \iint_{\O \times \R^3} |\nabla_x \phi_{f-g} \cdot \nabla_v g|  |f-g|^{\delta}  \\
			& \lesssim 
			\int^t_0 \| \nabla_x \phi_{f-g}\|_{L^{ \frac{3(1+ \delta)}{2- \delta}}_x} \| \nabla_v g \|_{L^{3}_xL^{1+ \delta}_v} \left\| |f-g|^\delta\right\|_{L_{x,v}^{\frac{1+ \delta}{\delta}}}\\
			& \lesssim \sup_{0 \leq s \leq t} \| \nabla_v g (s) \|_{L^{3}_xL^{1+ \delta}_v}\times \int^t_0 \| [f-g](s) \|_{1+ \delta}^{1+ \delta} \dd s.
		\end{split} \Ee

		A simple modification of (\ref{k_p_bound}) and (\ref{phi_p_bound}) as 
		\Be
		\begin{split}\notag
			&\int^t_0 \int_x \int_v \int_u \mathbf{k}_\varrho(v,u) |f(u)-g(u)| |f(v)-g(v)|^\delta\\
			\lesssim & \ \int^t_0 \int_x \int_v \int_u\mathbf{k}_\varrho(v,u)^{\frac{1}{1+\delta}} |f(u)-g(u)|  \mathbf{k}_\varrho(v,u)^{\frac{\delta}{1+\delta}}|f(v)-g(v)|^\delta\\
			\lesssim & \ \int^t_0 \int_x \int_{v }  |f(v)- g(v)|^{1+ \delta}\int_u\mathbf{k}_\varrho (v,u)\\
			\lesssim & \ \int^t_0   \|f - g \|^{1+ \delta}_{1+\delta}  ,
		\end{split}
		\Ee
		leads to 
		\Be\label{gronwall_f-g_last}\begin{split}
			&\int^t_0 \iint_{\O \times \R^3} | \text{the }  2^{\text{nd} }   \text{ line of RHS of } (\ref{eqtn_f-g})|  |f-g|^{\delta}  \\
			\lesssim& \sup_{0 \leq s \leq t}\big\{ 1+ \| w_{\vartheta}f(s) \|_\infty + \|w_{\vartheta} g(s) \|_\infty
			\big\}  \int^t_0 \| f-g \|_{1+ \delta}^{1+ \delta}.
		\end{split}\Ee

		Then following the proof of (\ref{bound_p_est}) and applying (\ref{gronwall_f-g}) to (\ref{gronwall_f-g_last}), we can obtain
		\Be\label{gamma_+,1+delta}
		\begin{split}
			&\int_0^t |[f-g]|_{1+ \delta, - }^{1 + \delta} \\
			&\lesssim   o(1) \int_0^t |[f-g]|_{1+ \delta, + }^{1 + \delta}+\| [f- g](0)\|_{1+ \delta}^{1+ \delta} \\
			&  +
			\sup_{0 \leq s \leq t}\big\{ 1+ \| \nabla_v g (s) \|_{L^{3}_xL^{1+ \delta}_v} + \| w_{\vartheta} f(s) \|_\infty + \|w_{\vartheta} g(s) \|_\infty
			\big\}\\
			& \ \ \ \ \times   \int^t_0 \| f-g \|_{1+ \delta}^{1+ \delta}.
		\end{split}
		\Ee
		Using (\ref{bound_nabla_v_g_global}), (\ref{f-g_energy}), (\ref{gronwall_f-g}), (\ref{gronwall_f-g_last}), (\ref{gamma_+,1+delta}) and applying the Gronwall inequality, we prove $L^{1+\delta}$-stability (\ref{1+delta_stability}) \hide in a small time interval on $[0,t]$ for $0< t \ll 1$. Now we split the time interval as $[0,t] \cup [t,2t] \cup [2t,3t] \cup \cdots$ and derive (\ref{bound_nabla_v_g}) for each interval. Then we prove (\ref{bound_nabla_v_g_global}) and then (\ref{1+delta_stability})\unhide for all time $t \geq 0$.\end{proof}
	\begin{proof}[\textbf{Proof of Proposition \ref{prop_better_f_v}}] 
		\textit{Step 1. } Note that from (\ref{eqtn_f}) and (\ref{boundary_v}), we have
		\Be\begin{split}\label{eqtn_g_v}
			&[\p_t + v\cdot \nabla_x - \nabla_x \phi_f \cdot \nabla_v + \nu(v) + \frac{v}{2} \cdot \nabla_x \phi_f  ] \p_v f \\
			=&  -\p_x f- \frac{1}{2} \p_x \phi_f f - \p_v \nu f + \p_v(Kf) + \p_v (\Gamma (f,f)) + |\nabla_x \phi_f| \langle v\rangle^2\sqrt{\mu} , 
		\end{split}\Ee
		with the boundary bound for $(x,v) \in\gamma_-$
		\Be\label{bdry_g_v}
		\big|\p_v f  \big| \lesssim   |v| \sqrt{\mu} \int_{n \cdot u>0} |f| \sqrt{\mu} \{n \cdot u \} \dd u \ \  {on } \ \gamma_-.
		\Ee
		From (\ref{nabla_nu}), (\ref{vKsum}), (\ref{bound_Gamma_nabla_vf1}), (\ref{bound_Gamma_nabla_vf2}), (\ref{Gvloss}), and (\ref{Gvgain}), we obtain the following bound along the characteristics
		\begin{align}
		& |\p_v f(t,x,v)|\nonumber
		\\
		 \leq &   \mathbf{1}_{\tb(t,x,v)> t}  
		|\p_v f(0,X(0;t,x,v), V(0;t,x,v))|\label{g_initial}\\
		& +    \ \mathbf{1}_{\tb(t,x,v)<t}
		\mu(\vb)^{\frac{1}{4}}  \int_{n(\xb) \cdot u>0} 
		| f(t-\tb, \xb, u) |\sqrt{\mu} \{n(\xb) \cdot u\} \dd u\label{g_bdry}\\
		&  +    \int^t_{\max\{t-\tb, 0\}} 
		|\p_x f(s, X(s;t,x,v),V(s;t,x,v))|
		\dd s\label{g_x}\\
		&   +  \int^t_{\max\{t-\tb, 0\}} 
		(1+ \| w_{\vartheta} f \|_\infty)
		\int_{\R^3} \mathbf{k}_\varrho (V(s),u) |\p_v f(s,X(s),u)| \dd u 
		\dd s\label{g_K}\\
		& +   
		\int^t_{\max\{t-\tb, 0\}} 
		|\nabla_x \phi_f (s, X(s;t,x,v))| \mu ^{1/4}
		\dd s \label{g_phi}\\
		&+  
		\int^t_{\max\{t-\tb, 0\}}  
		(1+ \delta_1)  \delta_1  |{w}_{\tilde{\vartheta}}( V(s;t,x,v))|^{-1} 
		\dd s 
		\label{g_infty}
		,
		\end{align}
		where $\delta_1$ is in (\ref{integrable_nabla_phi_f}).
		
		Note that if $|v| > 2 \frac{\delta_1}{\Lambda_1}$, then from (\ref{integrable_nabla_phi_f}) and (\ref{delta_1/lamdab_1}), for $0 \leq s \leq t$,
		\Be\label{V_lower_bound_v}
		\begin{split}
			|V(s;t,x,v)| &\geq |v| - \int^t_0 |\nabla_x \phi (\tau;t,x,v)| \dd \tau \\
			&\geq |v| 
			- \delta_1/\Lambda_1
			\\
			&    \geq \frac{|v|}{2}.
		\end{split}
		\Ee
		Therefore 
		\Be\label{tilde_w_integrable}
		\sup_{s,t,x}\left\|    \frac{1}{{w}_{\tilde{\vartheta}} (V(s;t,x,v)) }    \right\|_{ L^{r}_v} \lesssim 1 \  \   { for \ any } \     1 \leq r \leq \infty .
		\Ee

		We derive 
		\Be\begin{split}\label{est_g_initial}
			&\| (\ref{g_initial})\|_{L^3_x L^{1+ \delta}_v}\\
			\lesssim & \ \bigg(
			\int_{\O}
			\left(\int_{\R^3} |w_{\tilde{\vartheta}} \p_v f(0,X(0 ), V(0 ))|^3 
			\right)\\
			&  \ \ \ \times 
			\left(
			\int_{\R^3} \frac{\dd v}{|w_{\tilde{\vartheta}} (V(0 ))|^{(1+ \delta) \frac{3}{2-\delta}}}
			\right)^{\frac{2-\delta}{1+ \delta}}
			\bigg)^{1/3} \\
			\lesssim & \
			\left(\iint_{\O \times \R^3} |w_{\tilde{\vartheta}}(V(0 )) \p_v f(0,X(0 ), V(0 ))|^3 \dd v \dd x\right)^{1/3}\\
			\lesssim & \ \| w_{\tilde{\vartheta}} \p_v f (0) \|_{L^3_{x,v}},
		\end{split}
		\Ee
		where we have used a change of variables $(x,v) \mapsto (X(0;t,x,v), V(0;t,x,v))$ and (\ref{tilde_w_integrable}).

		\hide
		Also we use $|V(0;t,x,v)| \gtrsim |v|$ for $|v|\gg 1$, from (\ref{decay_phi}), and hence $\tilde{w}(V(0;t,x,v))^{- (1+ \delta) \frac{3}{2-\delta}} \in L_v^1 (\R^3)$.\unhide
		
		Clearly 
		\Be\label{est_g_bdry}
		\|(\ref{g_bdry})\|_{L^3_x L^{1+ \delta}_v} + \|(\ref{g_infty})\|_{L^3_x L^{1+ \delta}_v} \lesssim \sup_{0 \leq s \leq t} \| w_{\vartheta} f(s) \|_\infty.
		\Ee
		
		From $W^{1,2}(\O)\subset L^6(\O)\subset L^2(\O)$ for a bounded $\O \subset \R^3$, and the change of variables $(x,v) \mapsto (X(s;t,x,v), V(s;t,x,v))$ for fixed $s\in(\max\{t-\tb,0\},t)$,\Be
		\begin{split}\label{est_g_phi}
			\|(\ref{g_phi})\|_{L^3_x L^{1+ \delta}_v} 
			\lesssim & \ \int^t_0 \| \mu^{1/8}  \nabla_x \phi_f (s,X(s;t,x,v))  \|_{L^3_{x,v} }\| 
			\mu^{1/8}
			\|_{L^{
					\frac{3(1+ \delta)}{2- \delta}
				}_v}\\
			\lesssim & \ \int^t_0 \| \nabla_x \phi_f (s)   \|_{L^3_{x} }
			\lesssim    \int^t_{\max\{t-\tb,0\}} \| \phi_f (s) \|_{W^{2,2}_{x} } 
			\\
			\lesssim & \ \int^t_0 \|w_{\tilde{\vartheta}} f(s) \|_{2}.
		\end{split}
		\Ee

		\vspace{4pt}
		
		\textit{Step 2. }  We claim 
		\Be\label{est_g_x}
		\|(\ref{g_x})\|_{L^3_x L^{1+ \delta}_v} \lesssim \int^t_0 \| w_{\tilde{\vartheta}} \alpha_{f,\e}^\beta \p_x f (s) \|_{L^p_{x,v}}. 
		\Ee
		Now we have for $3<p<6$, by the H\"older inequality $\frac{1}{1+ \delta}= \frac{1}{ \frac{p+p \delta}{p-1 - \delta}}+ \frac{1}{p}$,
		\Be\label{init_p_xf}
		\begin{split}
			&\left\|\left\| \int^t_{\max\{t-\tb, 0\}}
			\p_x f(s,X(s;t,x,v),V(s;t,x,v)) \dd s
			\right\|_{L_{v}^{1+ \delta}(  \R^3)}\right\|_{L^{3}_x}\\
			\lesssim & \ \left\|\left\| \int^t_{\max\{t-\tb, 0\}} \frac{ w_{\tilde{\vartheta}} \alpha_{f,\e}^{\beta} \p_x f(s,X(s ),V(s ))}{w_{\tilde{\vartheta}}\alpha_{f,\e}(s,X(s ),V(s ))^{\beta} }
			\dd s
			\right\|_{L_{v}^{1+ \delta}(  \R^3)}\right\|_{L^{3}_x}
			\\
			\lesssim & \ \left\| \frac{w_{\tilde{\vartheta}}(v)^{-1}}{\alpha_{f,\e}(t,x,v)^\beta}\right\|_{L_v^{\frac{p+p \delta}{p-1 - \delta}} (\R^3) }\\
			& \times \left\|
			\left\| \int^t_0 w_{\tilde{\vartheta}}\alpha_{f,\e}^\beta \p_x f(s,X(s;t,x,v), V(s;t,x,v)) \dd s\right\|_{L_{v}^p( \R^3)}\right\|_{L^{3}_x}\\
			\lesssim & \ \left\| \frac{ w_{\tilde{\vartheta}}(v)^{-1}}{\alpha_{f,\e}(t,x,v)^\beta}\right\|_{L_v^{\frac{p+p \delta}{p-1 - \delta}} (\R^3) }
			\times \int^t_0 \| w_{\tilde{\vartheta}} \alpha_{f,\e}^\beta \p_x f (s) \|_{L^p_{x,v}} \dd s 
			,
		\end{split}\Ee
		where we have used $\alpha_{f,\e}(t,x,v)=\alpha_{f,\e}(s,X(s;t,x,v),V(s;t,x,v))$ for $t-\tb(t,x,v)\leq s \leq t$ and the change of variables $(x,v) \mapsto (X(s;t,x,v), V(s;t,x,v))$ and the Minkowski inequality.
		
		For $\beta$ in (\ref{W1p_initial}), we have $
		\beta \frac{p}{p-1} < 1$ since $\frac{2}{3} < \frac{p-1}{p}$ for $3 < p$. Therefore, we can choose $0 < \delta \ll 1$ so that $\beta$ in (\ref{W1p_initial}) satisfies 
		\Be\label{less_than_1}
		\beta \times \frac{p+p \delta}{p-1 - \delta} < 1.
		\Ee
		\hide
		For $\beta$ in (\ref{W1p_initial}) we have $
		\beta \frac{p+p \delta}{p-1 - \delta}> \frac{p+p\delta}{p} \frac{p-2}{p-1-\delta}.$ Therefore for $0<\delta\ll1$ we can choose $\beta$ in (\ref{W1p_initial}) and satisfies 
		\Be\label{less_than_1}
		\beta \times \frac{p+p \delta}{p-1 - \delta} < 1.
		\Ee\unhide
		We apply Proposition \ref{prop_int_alpha} to conclude that 
		\Be\label{bound_g_x}
		\sup_{t,x} \left\| \frac{w_{\tilde{\vartheta}} (v)^{-1}}{\alpha_{f,\e}(t,x,v)^\beta}\right\|_{L_v^{\frac{p+p \delta}{p-1 - \delta}} (\R^3) }^{\frac{p+p \delta}{p-1 - \delta}} 
		= \sup_{t,x} \int_{\R^3}  \frac{
			e^{- \tilde{\vartheta} \frac{p+p \delta}{p-1 - \delta} |v|^2 }
		}{\alpha_{f,\e}(t,x,v)^{\beta \frac{p+p \delta}{p-1 - \delta}   }}   \dd v \lesssim 1.
		\Ee
		\hide
		We have 
		\Be
		\begin{split}
			| \tilde{w} \mathcal{G}|\lesssim&  \tilde{w} (\ref{G est})\\
			\lesssim & \tilde{w} |\nabla_{x,v} g| +
		\end{split}
		\Ee

		{\color{red}DETAILS}\unhide
		Finally, from (\ref{init_p_xf}), (\ref{bound_g_x}), and (\ref{W1p_main}), we conclude the claim (\ref{est_g_x}).
		
		\hide
		Now for $6>p>3$ and $\beta$ in (\ref{W1p_initial}) and a bounded domain $\O$, by the Minkowski inequality and a change of variables $(x,v) \mapsto (X(s;t,x,v), V(s;t,x,v))$ for fixed $s$, whose Jacobian equals 1,
		\Be\notag
		\begin{split}
			&\left\| (\ref{init_p_xf}) \right\|_{L^{3}_x}\\
			\lesssim& \  \left\| \int^t_0 \tilde{w} \alpha^\beta \p_x f(s,X(s;t,x,v), V(s;t,x,v)) \dd s\right\|_{L ^p(\O \times \R^3)}\\
			\lesssim & \ \int^t_0 \| \tilde{w} \alpha^\beta \p_x f(s) \|_p^p \dd s.
		\end{split}
		\Ee
		From (\ref{tilde_w_W1p}) we conclude the claim (\ref{est_g_x}).
		\unhide
		
		\vspace{4pt}
		
		\textit{Step 3. } We consider (\ref{g_K}). We split the $u$-integration of (\ref{g_K}) into two parts with $N\gg 1$ as 
		\begin{eqnarray}
		&&\int_{|u| \leq N} \mathbf{k}_\varrho (V(s), u) |\nabla_v f(s,X(s ), u) | \dd u \label{g_K_split1}\\
		&+& \int_{|u| \geq N} \mathbf{k}_\varrho (V(s), u) |\nabla_v f(s,X(s ), u) | \dd u .\label{g_K_split2}
		\end{eqnarray}

		First we bound (\ref{g_K_split1}). From the change of variables $(x,v) \mapsto (X(s;t,x,v), V(s;t,x,v))$ for $t-\tb(t,x,v)\leq s \leq t$
		\Be\label{g_K_COV}
		\begin{split}
			&\left\|\int_{|u| \leq N} \mathbf{k}_\varrho (V(s;t,x,v), u) |\nabla_v f(s,X(s;t,x,v ), u) | \dd u  \right\|_{L^3_x L^3_v}\\
			= & \ \left\|\int_{|u| \leq N} \mathbf{k}_\varrho (v, u) |\nabla_v f(s,x, u) | \dd u  \right\|_{L^3_x L^3_v}.
		\end{split}
		\Ee
		If $|v|\geq 2N$ then $|v-u|^2\gtrsim |v|^2$ and $\mathbf{k}_\varrho (v,u) \lesssim \frac{e^{-C|v|^2}}{|v-u|^2}$ for $|v|\geq 2N$ and $|u| \leq N$. For $0 < \delta \ll 1$ with $\frac{3(1+ \delta)}{1-2\delta}>3$,  
		\Be
		\begin{split}\label{g_K_COV1}
			&(\ref{g_K_COV})\\
			\lesssim & \ C_N \left\| \left\| \int_{|u| \leq N} \mathbf{k}_\varrho (v, u) |\nabla_v f(s,x, u) | \dd u \right\|_{L^{\frac{3(1+ \delta)}{1-2\delta}}_v
				(\{|v| \leq 2N\})
			} \right\|_{L^3_x  }\\
			+& \left\| \left\|  e^{-C|v|^2}\right\|_{L^{3/2}_v} \left\| \int_{|u| \leq N} \frac{1}{|v-u|} |\nabla_v f(s,x, u) | \dd u \right\|_{L^{\frac{3(1+ \delta)}{1-2\delta}}_v
				(\{|v| \geq 2N\})
			} \right\|_{L^3_x  }\\
			\lesssim & \ \left\|   \left\|   \frac{1}{|v- \cdot |} * |\nabla_v f(s,x, \cdot ) |   \right\|_{L^{\frac{3(1+ \delta)}{1-2\delta}}_v 
			} \right\|_{L^3_x  }.
		\end{split}
		\Ee
		Then by the Hardy-Littlewood-Sobolev inequality with $1+ \frac{1}{\frac{3(1+ \delta)}{1- 2\delta}} = \frac{1}{3} + \frac{1}{1+ \delta}$, we derive that   
		\Be\begin{split}\notag
			(\ref{g_K_COV1}) 
			\lesssim & \ \left\| \| \nabla_v f(s, x,v)  \|_{L^{1+ \delta}_v  }\right\|_{L^3_x}
			= \| \nabla_v f(s) \|_{L^3_x L^{1+ \delta}_v}
			.
		\end{split}\Ee 
		Combining the last estimate with (\ref{g_K_COV}), (\ref{g_K_COV1}), we prove that 
		\Be\label{est_g_K_split1}
		\|(\ref{g_K_split1})\|_{L^3_x L^{1+ \delta}_v} \lesssim    \| \nabla_v f(s) \|_{L^3_x L^{1+ \delta}_v}.
		\Ee
		
		Now we consider (\ref{g_K_split2}). Choose $0< \iota \ll  1$. We write as 
		\Be\notag
		\begin{split}
			&(\ref{g_K_split2}) \\
			= \  &\int_{|u| \geq N} \frac{1}{w_{\tilde{\vartheta}} (V(s;t,x,v))^{1- \iota}} \frac{w_{\tilde{\vartheta}}(V(s;t,x,v))  }{w_{\tilde{\vartheta}}(u)} \frac{\mathbf{k}_\varrho (V(s;t,x,v), u)}{\alpha_{f,\e}(s,X(s;t,x,v), u)^\beta} \\
			& \ \  \times \frac{ w_{\tilde{\vartheta}}(u)}{w_{\tilde{\vartheta}}(V(s;t,x,v))^{ \iota}} \alpha_{f,\e}(s,X(s;t,x,v), u)^\beta | \nabla_v f(s,X(s;t,x,v), u)| \dd u .
		\end{split} 
		\Ee

		By the H\"older inequality with $\frac{1}{p} + \frac{1}{p^*}=1$ with $3<p<6$
		\Be
		\begin{split}\label{g_K_split2_Holder}
			&|(\ref{g_K_split2})| \\
			&\lesssim  \ \frac{1}{w_{\tilde{\vartheta}} (V(s;t,x,v))^{1- \iota}}\\
			& \ \ \ \ \times  \left\|  \frac{w_{\tilde{\vartheta}}(V(s;t,x,v))  }{w_{\tilde{\vartheta}}(u)} \frac{\mathbf{k}_\varrho (V(s;t,x,v), u)}{\alpha_{f,\e}(s,X(s;t,x,v), u)^\beta} \right\|_{L^{p^*} (\{ |u|\geq N \})} \\
			&  \ \ \ \  \times  \bigg\|  \frac{w_{\tilde{\vartheta}}(u) }{w_{\tilde{\vartheta}}(V(s;t,x,v))^{ \iota}} \alpha_{f,\e}(s,X(s;t,x,v), u)^\beta\\
		&  \ \ \ \ \ \   \ \ \ \  \ \ \ \  \ \ \ \  \ \ \ \  \ \ \ \  \times 	 | \nabla_v f(s,X(s;t,x,v), u)| \bigg\|_{L^p_u (\R^3)}.
		\end{split}
		\Ee
		
		Then by the H\"older inequality with $\frac{1}{1+ \delta} = \frac{1}{p} + \frac{1}{\frac{(1+ \delta)p}{p - (1+ \delta)}}$,
		\Be
		\begin{split}\notag
			&\|(\ref{g_K_split2})\|_{L^{1+ \delta}_v}\\
			\lesssim & \ \left\| \frac{1}{w_{\tilde{\vartheta}} (V(s;t,x,v))^{1- \iota}} \right\|_{L_v^{\frac{(1+ \delta) p}{ p- (1+ \delta)}}} \\
			& \times \sup_v \left\|  \frac{w_{\tilde{\vartheta}}(V(s;t,x,v))  }{w_{\tilde{\vartheta}}(u)} \frac{\mathbf{k}_\varrho (V(s;t,x,v), u)}{\alpha_{f,\e}(s,X(s;t,x,v), u)^\beta} \right\|_{L^{p^*} (\{ |u|\geq N \})} \\
			& \times \left\| \left\|  \frac{w_{\tilde{\vartheta}}(u) }{w_{\tilde{\vartheta}}(V(s;t,x,v))^{ \iota}} \alpha_{f,\e}(s,X(s;t,x,v), u)^\beta | \nabla_v f(s,X(s;t,x,v), u)| \right\|_{L^p_u  }\right\|_{L^{p}_v}.
		\end{split}
		\Ee
		
		Note that, from (\ref{k_vartheta_comparision}), $\mathbf{k}_\varrho (v,u) \frac{e^{\tilde{\vartheta} |v|^2}}{e^{\tilde{\vartheta} |u|^2}} \lesssim \mathbf{k}_{\tilde{\varrho}} (v,u)$ for some $0<\tilde{\varrho}< \varrho$. Hence we derive, using (\ref{tilde_w_integrable}) 
		\Be
		\begin{split}\notag
			&\left\|\|(\ref{g_K_split2})\|_{L^{1+ \delta}_v}\right\|_{L^3_x}\\
			\lesssim_\O &  \ \sup_{X,V} \left\| \frac{e^{- \frac{\tilde{\vartheta}}{10}|V-u|^2  }}{|V-u|}  \frac{1}{\alpha_{f,\e}(s,X , u)^\beta} \right\|_{L^{p^*} (\{ |u|\geq N \})} \\
			&\times  \left\|  \frac{\tilde{w}(u) }{\tilde{w}(V(s;t,x,v))^{ \iota}} \alpha_{f,\e}(s,X(s;t,x,v), u)^\beta | \nabla_v f(s,X(s;t,x,v), u)| \right\|_{L^p_{u,v,x}  } .
		\end{split}
		\Ee
		Finally using (\ref{NLL_split3}) in Proposition \ref{prop_int_alpha} with $ \frac{p-2}{p-1}<\beta p^*< 1$ from (\ref{W1p_initial}) and applying the change of variables $(x,v) \mapsto (X(s;t,x,v), V(s;t,x,v))$, we derive that 
		\Be
		\begin{split}\label{g_K_split2_Holder2}
			&\left\|\|(\ref{g_K_split2})\|_{L^{1+ \delta}_v}\right\|_{L^3_x}\\
			\lesssim_\O &    
			\left\| \left\|  \frac{1 }{w_{\tilde{\vartheta}}(v)^{ \iota}}w_{\tilde{\vartheta}}(u) \alpha_{f,\e}(s,x, u)^\beta | \nabla_v f(s,x, u)| \right\|_{L^p_v} \right\|_{L^p_{u ,x}  } \\
			\lesssim \ & \left\|  \frac{1 }{w_{\tilde{\vartheta}}(v)^{ \iota}} \right\|_{L^p_v}    \left\| w_{\tilde{\vartheta}}(u) \alpha_{f,\e}(s,x, u)^\beta | \nabla_v f(s,x, u)|   \right\|_{L^p_{u ,x}  }\\
			\lesssim \ &   \left\| w_{\tilde{\vartheta}} \alpha_{f,\e}^\beta | \nabla_v f(s )|   \right\|_{L^p  }.
		\end{split}
		\Ee
		Combining (\ref{g_K_split2_Holder}) and (\ref{g_K_split2_Holder2}) we conclude that 
		\Be
		\|(\ref{g_K_split2})\|_{L^3_x L^{1+ \delta}_v} \lesssim   \|w_{\tilde{\vartheta}} \alpha_{f,\e}^\beta   \nabla_v f(s )  \|_{L^p_{x,v}}.\label{est_g_K_split2}
		\Ee
		
		Finally from (\ref{est_g_K_split1}) and (\ref{est_g_K_split2}), and using the Minkowski inequality, we conclude that 
		\Be\label{est_g_K}\begin{split}
			& \| (\ref{g_K}) \|_{L^3_v L^{1+ \delta}_x} \\
			\lesssim & \  (1+ \| w_{\vartheta} f\|_\infty) \int^t_0  \big[
			\| \nabla_v f (s) \|_{L^3_x L^{1+ \delta}_v}
			+
			\| w_{\tilde{\vartheta}} \alpha_{f,\e}^\beta | \nabla_v f(s )| \|_{L^p_{x,v}}\big] \dd s. 
		\end{split}\Ee
		
		Collecting terms from (\ref{g_initial})-(\ref{g_phi}), and (\ref{est_g_initial}), (\ref{g_bdry}), (\ref{est_g_phi}), (\ref{est_g_x}), (\ref{est_g_K}), we derive
		\Be\begin{split}\label{bound_nabla_v_g}
			&\sup_{0 \leq s \leq t}\| \nabla_vf(s) \|_{L^3_xL^{1+ \delta}_v} \\
			\lesssim & \ 1+
			\| w_{\tilde{\vartheta}} \nabla_v f(0) \|_{L^3_{x,v}} + 
			\sup_{0 \leq s \leq t}
			\| w_{\vartheta} f(s) \|_\infty +\int^t_0 \| w_{\tilde{\vartheta}} f(s)\|_2\\
			&  
			+
			t (1+ \| w_{\vartheta} f \|_\infty) \sup_{0 \leq s \leq t}
			\|w_{\tilde{\vartheta}} \alpha_{f,\e}^\beta \nabla_{x,v} f(s) \|_{p}\\
			&
			+(1+ \| w_{\vartheta} f\|_\infty) \int^t_0  
			\| \nabla_v f (s) \|_{L^3_x L^{1+ \delta}_v}
			,
		\end{split} \Ee
		where the first two lines of RHS is bounded due to the assumptions (\ref{main_Linfty}), (\ref{W1p_main}). \hide

		Then by the Holder inequality $\frac{1}{3} = \frac{1}{\frac{3p}{p-3}} + \frac{1}{p}$
		\Be\label{g_K_split2_Holder_Holder}
		\begin{split}
			&\|(\ref{g_K_split2_Holder})\|_{L^3_x L^{1+ \delta}_v }\\
			\lesssim  & \  \left\| \left\|   \frac{\mathbf{k}_\varrho (V(s;t,x,v), u)}{\tilde{w} (u)\alpha(s, X(s;t,x,v), u )^\beta} \right\|_{L^{p^*} (\{ |u|\geq N \})  } \right\|_{L_x^{\frac{3p}{p-3} } L^{1+ \delta}_v } \\
			& \times \| \tilde{w} \alpha^\beta \p_v g(s, X(s;t,x,v), u) \|_{L^p_{x,u }}  .
		\end{split}
		\Ee
		
		\Be\begin{split}\label{g_K_split2_Holder_Holder1}
			&\left\|  \left\|   \frac{\mathbf{k}_\varrho (V(s;t,x,v), u)}{\tilde{w} (u)\alpha(s, X(s;t,x,v), u )^\beta} \right\|_{L^{p^*}(\{ |u|\geq N \})   } \right\|_{L^{1+ \delta}_v }
			\\
			\lesssim&  \ 
			\left\| \frac{1}{\tilde{w} (V(s;t,x,v))} \left\|   \frac{\mathbf{k}_{\tilde{\varrho}} (V(s;t,x,v), u)}{ \alpha(s, X(s;t,x,v), u )^\beta} \right\|_{L^{p^*}(\{ |u|\geq N \})  }
			\right\|_{ L^{1+ \delta}_v  }\\
			\lesssim & \  
			\left\|    \tilde{w} (V(s;t,x,v))^{-1}   \right\|_{ L^{1+ \delta}_v  }
			\times \sup_v  \left\|   \frac{\mathbf{k}_{\tilde{\varrho}} (V(s;t,x,v), u)}{ \alpha(s, X(s;t,x,v), u )^\beta} \right\|_{L^{p^*}(\{ |u|\geq N \})  }.
		\end{split}\Ee

		On the other hand from (\ref{NLL_split3}) in Proposition \ref{prop_int_alpha} the last term of (\ref{g_K_split2_Holder_Holder1}) is bounded.

		From (\ref{g_K_split2_Holder}), (\ref{g_K_split2_Holder_Holder}), (\ref{g_K_split2_Holder_Holder1})
		\Be
		\|(\ref{g_K_split2})\|_{L^3_x L^{1+ \delta}_x}\lesssim  \| \tilde{w} \alpha^\beta \nabla_v g (s) \|_{L^p_{x,v}}
		\Ee

		Now we check that a change of variables $v \mapsto V(s;t,x,v)$ for fixed $s,t,x$ with $|t-s|\ll \lambda_\infty$ in (\ref{main_Linfty}) is well-defined and has a Jacobian 
		\Be\label{V_v_short}
		\begin{split}
			&\det \left(\frac{\p V(s;t,x,v)}{\p v}\right)\\
			= & \  \det \left(\text{Id}_{3 \times 3}- \int^s_t \frac{\p X(\tau;t,x,v)}{\p v} \cdot \nabla_x \nabla_x \phi (\tau, X(\tau;t,x,v)) \dd \tau \right)\\
			\gtrsim & \ 1  - \int^t_s |t- \tau| e^{- \frac{\lambda_\infty}{2} \tau } \dd \tau\\
			\gtrsim & \ 1- \frac{|t-s|}{\lambda_\infty}\\
			\gtrsim & \ 1 \ \ \  \text{for} \ \   |t-s| \ll \lambda_\infty,
		\end{split}
		\Ee
		where we have used (\ref{result_X_v}).

		We need the following estimate
		\Be\label{X_x_V_v_local}
		|\nabla_x X(s;t,x,v)| + |\nabla_v V(s;t,x,v)| \lesssim e^{C|t-s|},
		\Ee
		where $C>0$ depends on $\sup_{s \leq \tau \leq t }\|\nabla_x^2 \phi_g(\tau) \|_\infty$. This is a direct consequence of Gronwall inequality.

		Note that for fixed $t,x,v$ with $|t-s| \ll 1$ in (\ref{main_Linfty}), the change of variables $x \mapsto X(s;t,x,v)$ is well-defined and has a Jacobian
		\Be\label{X_x_jacobian}
		\begin{split}
			&\det \left(\frac{\p X(s;t,x,v)}{\p x}\right)\\
			= & \  \det \left( \text{Id}_{3 \times 3} - \int^s_t \int^{\tau }_t \frac{\p X(\tau^\prime;t,x,v)}{\p x} \cdot \nabla_x \nabla_x \phi (\tau^\prime, X(\tau^\prime;t,x,v)) \dd \tau^\prime \dd \tau \right)\\
			\gtrsim &  \ 1 - O\left(|t-s|^2\right)\\
			\gtrsim& \ 1 .
		\end{split}
		\Ee
		Here we have used (\ref{main_Linfty}) and (\ref{X_x_V_v_local}).
		
		Then from the above change of variables with (\ref{X_x_jacobian}) and the Minkowski inequality we derive
		\Be
		\begin{split}
			&\left\|(\ref{g_K_split1})
			\right\|_{L^3_x (\O)}\\
			\lesssim & \ \int_{|u| \leq N} \mathbf{k}_\varrho (V(s), u) \|\nabla_v g(s,x, u) \|_{L^3_x  } \dd u
		\end{split}
		\Ee

		If $|v|\geq 10N\gg  \frac{\| w f_0 \|_\infty}{\lambda_\infty} $ then from (\ref{Morrey})
		\Be\notag
		\begin{split}
			|V(s;t,x,v)| &\geq |v| - \int^t_s |\nabla_x \phi (\tau;t,x,v)| \dd \tau \\
			&\geq |v| - C \| w f_0 \|_\infty   \int^t_s e^{- \lambda_\infty \tau } \dd \tau \\
			&\geq |v| - \frac{C}{\lambda_\infty}\| w f_0 \|_\infty \\
			& \geq \frac{|v|}{2}.
		\end{split}
		\Ee
		Therefore if $|u| \leq N$ and $|v|\geq 10N$ then 
		\[
		|V(s;t,x,v)- u|^2 \gtrsim  |v|^2,
		\]
		which implies that 
		\Be \notag
		\mathbf{k}_\varrho (V(s;t,x,v), u) 
		\lesssim    \frac{e^{-C|v  |^2}}{|V(s;t,x,v) - u|} .
		\Ee
		Therefore we have 
		\Be
		\begin{split}\label{g_K_small}
			&\left\|\int_{|u | \leq N} \mathbf{k}_\varrho (V(s;t,x,v), u) |\p_v  g(s,X(s;t,x,v), u)| \dd u\right\|_{L^{1+ \delta}_v (\R^3)}\\
			\lesssim & \ \left\{1+\| e^{-C|v|^2} \|_{L^{3/2}_v (\R^3)}\right\}\\
			& \times \left\| 
			\frac{1}{|V(s;t,x,v)-\cdot |} *|\nabla_v g (s,X(s;t,x,v), \cdot)|
			\right\|_{L^{\frac{3(1+ \delta)}{1-2\delta}}_v (\R^3)}.
		\end{split}
		\Ee
		Using (\ref{V_v_short}) and the Hardy-Littlewood-Sobolev inequality with $1+ \frac{1}{\frac{3(1+ \delta)}{1- 2\delta}} = \frac{1}{3} + \frac{1}{1+ \delta}$, we derive that   
		\Be\begin{split}
			(\ref{g_K_small}) \lesssim& \  \left\| 
			\frac{1}{|v-\cdot |} *|\nabla_v g (s,X(s;t,x,v), \cdot)|
			\right\|_{L^{\frac{3(1+ \delta)}{1-2\delta}}_v (\R^3)}\\
			\lesssim & \  \| \nabla_v g(s,X(s;t,x,V(s;t,x,v)), \cdot)\|_{L^{1+ \delta}_v (\R^3)}
			.
		\end{split}\Ee

		Therefore we have 
		\Be
		\begin{split}
			&\left\|\int_{\R^3} \mathbf{k}_\varrho (V(s;t,x,v), u) |\p_v  g(s,X(s;t,x,v), u)| \dd u\right\|_{L^{1+ \delta}_v (\R^3)}\\
			\lesssim & \
			\int_{\R^3} \mathbf{k}_{\tilde{\varrho}} (V(s;t,x,v), u) |\nabla_v  g(s,X(s;t,x,v), u)|  \dd u
			\\
			\lesssim & \ \left\| 
			\frac{1}{|V(s;t,x,v)-\cdot |} *|\nabla_v g (s,X(s;t,x,v), \cdot)|
			\right\|_{L^{\frac{3(1+ \delta)}{1-2\delta}}_v (\R^3)}\\
			\lesssim & \  \| \nabla_v g(s,X(s;t,x,v), \cdot)\|_{L^{1+ \delta}_v (\R^3)}
		\end{split}
		\Ee

		\vspace{4pt}
		
		\textit{Step .}
		
		\unhide
		By the Gronwall inequality we prove (\ref{bound_nabla_v_g_global}).\hide
		we derive that 
		\Be\label{improve_v}
		\sup_{0 \leq s \leq   t} \| \nabla_v g (s)\|_{L^3_x L^{1+ \delta}_v} \lesssim 1.
		\Ee
		From (\ref{improve_v}), we conclude the uniqueness by combining (\ref{f-g_energy}), (\ref{gronwall_f-g}), (\ref{gronwall_f-g_last}), (\ref{gamma_+,1+delta}), (\ref{improve_v}) and applying the Gronwall inequality. 
		
		\unhide\end{proof}

	\section{Local Existence}

	\begin{theorem}\label{local_existence}
		Let $0< \tilde{\vartheta}< \vartheta \ll1$. Assume that for sufficiently small $M>0$, $F_0= \mu+ \sqrt{\mu}f_0\geq 0$, and 
		\[
		\| w_\vartheta f_0 \|_\infty 
		\leq \frac{M}{2},
		\]
		and $ \| w_{\tilde{\vartheta}} \alpha_{f_0 ,\e}^\beta \nabla_{x,v} f _0 \|_{p} <\infty$ for $0< \e \ll1$ and (\ref{compatibility_condition}), where $(p, \beta)$ satisfies (\ref{beta_condition}), and $\| w_{\tilde{\vartheta}} \nabla_v f_0 \|_{L^{3}_{x,v}}<+\infty$. 
		
	Then there exists $T^*(M)>0$ and a unique solution $F(t,x,v) =\mu+ \sqrt{\mu} f(t,x,v)\geq 0$ to (\ref{Boltzmann_E}), (\ref{diffuse_BC}), (\ref{Field}), (\ref{Poisson}), and (\ref{phi_BC}) in $[0, T^*(M)) \times \O \times \R^3$ such that 
		\Be\begin{split}\label{infty_local_bound} \sup_{0 \leq t \leq T*} 
			\| w_\vartheta f  (t) \|_{\infty}
			\leq M.
		\end{split}\Ee 
		
		Moreover
		\Be\label{31_local_bound}
		\sup_{0 \leq t \leq T^*}\| \nabla_v f (t) \|_{L^3_xL^{1+ \delta}_v}< \infty \ \ \text{for} \ 0< \delta \ll1,
		\Ee
		and 
		\Be\begin{split}\label{W1p_local_bound}
			\sup_{0 \leq t \leq T*} 
			\Big\{
			\| w_{\tilde{\vartheta}}\alpha_{f,\e }^\beta \nabla_{x,v} f (t) \|_{p} ^p
			+ 
			\int^t_0  |w_{\tilde{\vartheta}} \alpha_{f,\e}^\beta \nabla_{x,v} f (t) |_{p,+}^p
			\Big\}
			< \infty.
		\end{split}\Ee
		
		Furthermore, $\|  w_\vartheta f  (t) \|_{\infty}$, $\| \nabla_v f (t) \|_{L^3_xL^{1+ \delta}_v}$ and $\| w_{\tilde{\vartheta}}\alpha_{f,\e }^\beta \nabla_{x,v} f (t) \|_{p} ^p$,
		\\
		 and $\int^t_0  |w_{\tilde{\vartheta}} \alpha_{f,\e }^\beta \nabla_{x,v} f (t) |_{p,+}^p$ are continuous in $t$.
		
		\hide
		
		Furthermore if we add one more extra condition $\| \nabla_v f_0 \|_{L^3_x L^{1+ \delta} _v}<+\infty$ for any $0<\delta \ll1$ then we have a unique solution and  
		\Be\label{31_local_bound}
		\sup_{0 \leq t \leq T^*}\| \nabla_v f (t) \|_{L^3_xL^{1+ \delta}_v}< \infty.
		\Ee

		\unhide\end{theorem}

	\begin{proof}

	Within the whole proof of Theorem \ref{local_existence} we consider a sequence for $\ell\geq 0$
		\Be
		\begin{split}\label{Fell}
			\p_t F^{\ell+1} + v\cdot \nabla_x F^{\ell+1} - \nabla {\phi^\ell} \cdot \nabla_v F^{\ell+1}\\
			= Q_{\text{gain}} (F^\ell, F^\ell) - Q_{\text{loss}} (F^\ell, F^{\ell+1}),\\
			- \Delta \phi^\ell= \int_{\R^3} F^\ell \dd v - \rho_0, \ \ \int_\O \phi^\ell \dd x =0, \ \ \frac{\p \phi^\ell}{\p n}\Big|_{\p\O }=0,
		\end{split}\Ee
		and, on $(x,v) \in \gamma_-$,
		\Be\label{F_ell_BC}
		F^{\ell+1}(t,x,v) = c_\mu \mu \int_{n(x) \cdot v>0} F^\ell (t,x,u) \{n(x) \cdot u\} \dd u,
		\Ee
		and $F^{\ell+1} (0,x,v) = F_0 (x,v)$ for $\ell \geq0$. We set $F^{0}(t,x,v) \equiv \mu$ and $\phi^0\equiv 0$.
		 
	\vspace{2pt}
	
		\textit{Step 1. } We claim that 
		\Be\label{positive_preserving}
\textit{If  }  \ 0\leq F^\ell ,\textit{  then } \ 0\leq F^{\ell+1}.
		\Ee
		  We define 
		\Be
		\label{nu_f}
		\nu(F) := \iint_{\R^3 \times \S^2}
		|(v-u) \cdot \o| F(u)
		\dd \o \dd u. 
		\Ee
		Clearly we have 
		\Be\label{positive_nu_Q_bdry}
		\nu(F^\ell)\geq 0 \ \ \ \text{and} \ \ \ Q_{\mathrm{gain}}(F^\ell, F^\ell) \geq 0 \ \ \ \text{and} \ \ \ 
		F^{\ell+1} |_{\gamma_-} \geq 0.
		\Ee

		Denote the characteristics $(X^\ell, V^\ell)$ which solves
		\Be\label{XV_ell}
		\begin{split}
			\frac{d}{ds}X^\ell(s;t,x,v) &= V^\ell(s;t,x,v),\\
			\frac{d}{ds} V^\ell(s;t,x,v) &= - \nabla \phi^\ell (s, X^\ell(s;t,x,v)).\end{split}
		\Ee
		Note that it is well-known that $\nabla\phi^\ell$ in (\ref{Fell}) is quasi-Lipschitz continuous and (\ref{XV_ell}) induces a H\"older continuous characteristics $(X^\ell(s;t,x,v), V^\ell(s;t,x,v))$. (see Chapter 8 of \cite{MB} for example) 
		
		From (\ref{Fell}), for $t- \tb (t,x,v)\leq s\leq t$,
		\Be \label{trajec_seq}
		\begin{split}
			&\frac{d}{ds}\Big\{e^{-\int^t_s \nu(F^\ell)(\tau, X^\ell (\tau), V^\ell (\tau))  \dd \tau } F^{\ell+1} (s,X^\ell(s  ), V^\ell (s  ))\Big\}\\
			= & \ e^{-\int^t_s \nu(F^\ell)(\tau, X^\ell (\tau), V^\ell (\tau))  \dd \tau }  Q_{\text{gain}} (F^\ell, F^\ell) (s,X^\ell(s ), V^\ell (s )),
		\end{split}
		\Ee
		where $(X^\ell(s ), V^\ell (s )): = (X^\ell(s;t,x,v), V^\ell (s;t,x,v))$. From the representation along the trajectory and (\ref{positive_nu_Q_bdry}), we prove that $F^{\ell+1}\geq 0$, and hence (\ref{positive_preserving}).

		\vspace{2pt}

		\textit{Step 2-1. } We prove, from \textit{Step 2-1} to \textit{Step 2-3}, that by choosing $M \ll1$ and $T^*=T^*(M)\ll 1$ and 
		\Be\label{uniform_h_ell}
		\sup_{0 \leq t \leq T^*} \max_{\ell} \| w_{\vartheta}  f^\ell (t) \|_\infty \leq M.
		\Ee
		 We define 
		\Be
		\label{h}
		h^{\ell} (t,x,v) : = w_{\vartheta}(v) f^\ell (t,x,v).
		\Ee
		By an induction hypothesis we assume 
		\Be \label{hypAssump}
			\sup_{0 \leq t \leq T^*}\| h^{\ell}(t) \|_\infty\leq M.
		\Ee
		Then $h^{\ell+1}$ solves 
		\Be\label{fell_local}
		\begin{split}
			&[\p_t + v\cdot\nabla_x - \nabla_x \phi^{\ell} \cdot \nabla_v + \nu + \frac{v}{2} \cdot\nabla \phi^\ell  
			- \frac{\nabla_x \phi^\ell \cdot \nabla_v w_{\vartheta}}{w_{\vartheta}}
			] h^{\ell+1}\\
			=&  \  K_{w_{\vartheta}} h^{\ell}- v\cdot \nabla \phi^\ell w_{\vartheta} \sqrt{\mu} + w_{\vartheta}\Gamma_{\text{gain}} (\frac{h^\ell}{w_{\vartheta}}, \frac{h^\ell}{w_{\vartheta}}) -   w_{\vartheta}\Gamma_{\text{loss}} (\frac{h^\ell}{w_{\vartheta}}, \frac{h^{\ell+1}}{w_{\vartheta}}),\\
			& - \Delta \phi^\ell = \int_{\R^3} f^\ell \sqrt{\mu} \dd v , \ \ \int_\O \phi^\ell \dd x =0, \ \ \frac{\p \phi^\ell}{\p n}\Big|_{\p\O }=0,
		\end{split}
		\Ee 
		where $K_{w_{\vartheta} }( \ \cdot \ )=w_{\vartheta}  K( \frac {1}{w_{\vartheta} } \ \cdot)$. The boundary condition is
		\Be\label{bdry_local}
		h^{\ell+1} |_{\gamma_-} = c_\mu w_{\vartheta} \sqrt{\mu} \int_{n \cdot u>0}h^\ell w_{\vartheta}^{-1}\sqrt{\mu}   \{n \cdot u\} \dd u. 
		\Ee

		\hide
		The perturbation of $F^{\ell+1} = \mu + \sqrt{\mu} f^{\ell+1}$ solves
		\Be\label{fell_local}
		\begin{split}
			&[\p_t + v\cdot\nabla_x - \nabla_x \phi^{\ell} \cdot \nabla_v + \nu + \frac{v}{2} \cdot\nabla \phi^\ell  ] f^{\ell+1}\\
			=&  \ Kf^{\ell}- v\cdot \nabla \phi^\ell \sqrt{\mu} + \Gamma_{\text{gain}} (f^\ell, f^\ell) -   \Gamma_{\text{loss}} (f^\ell, f^{\ell+1}),\\
			& - \Delta \phi^\ell = \int_{\R^3} f^\ell \sqrt{\mu} \dd v , \ \ \int_\O \phi^\ell \dd x =0, \ \ \frac{\p \phi^\ell}{\p n}\Big|_{\p\O }=0,
		\end{split}
		\Ee 
		and, on $(x,v) \in \gamma_-$,
		\Be\label{bdry_local}
		f^{\ell+1}(t,x,v) = c_\mu \sqrt{\mu} \int_{n(x) \cdot v>0} f^\ell (t,x,u)\sqrt{\mu(u)} \{n(x) \cdot u\} \dd u,
		\Ee
		and $f^{\ell+1} (0,x,v) = f_0 (x,v)$.

		\unhide
		
		We define
		\Be\label{nu_ell}
		\nu^\ell (t,x,v) : = \nu(v) + \frac{v}{2} \cdot \nabla \phi^\ell - \frac{\nabla_x \phi^\ell \cdot \nabla_v w_{\vartheta}}{w_{\vartheta}}+ \nu(\sqrt{\mu} \frac{h^\ell}{w_{\vartheta}}).
		\Ee
		From (\ref{hypAssump}), for $M\ll 1$, $\| \nabla \phi^\ell \|_\infty \ll1$ and hence
		   \\
		\Be\label{lower_nu_l}
		\nu^\ell (t,x,v)\geq \frac{\nu_0}{2} \langle v\rangle  .
		\Ee
		Let 
		\Be\label{g_ell}
		g^\ell : = -   v\cdot \nabla \phi^\ell \sqrt{\mu} +  \Gamma_{\mathrm{gain}} (\frac{h^\ell}{w_{\vartheta}}, \frac{h^\ell}{w_{\vartheta}}).
		\Ee
		Note that 
		\Be\label{bound_g_ell}
		| w_{\vartheta} g^\ell |  \lesssim  \| h^\ell \|_\infty +  \langle v\rangle  \| h^\ell \|_\infty^ 2,
		\Ee
		where we have used
		\Be\label{est_Gamma}
		|w_{\vartheta}\Gamma(\frac{h}{w_{\vartheta}}, \frac{h}{w_{\vartheta}})| \lesssim \langle v\rangle \| h \|_\infty^2.
		\Ee

		Along the trajectory from (\ref{fell_local})
		\Be\begin{split}\label{duhamel_local}
			&\frac{d}{ds} \Big\{ e^{- \int^t_s \nu^\ell (\tau, X^\ell (\tau), V^\ell (\tau)) \dd \tau}
			h^{\ell+1} (s, X^\ell(s;t,x,v), V^\ell (s;t,x,v)) 
			\Big\}\\
			=& \ e^{- \int^t_s \nu^\ell (\tau, X^\ell (\tau), V^\ell (\tau)) \dd \tau}\\
			& \ \times 
			\Big\{K_{w_{\vartheta}} h^\ell (s, X^\ell (s), V^\ell(s))+ w_{\vartheta} g^\ell (s, X^\ell (s), V^\ell(s))\Big\}.
		\end{split}\Ee
		By (\ref{bdry_local}) and (\ref{bound_g_ell}) we represent $h^{\ell+1}$ along the stochastic cycles:
		\Be
		\begin{split}\label{cycle}
			&t^{\ell}_1 (t,x,v):= 
			\sup\{ s<t:
			X^\ell(s;t,x,v) \in \p\O
			\}
			,\\
			&x^\ell_1 (t,x,v ) := X^\ell (t^{\ell}_1 (t,x,v);t,x,v)
			,\\
			&t^{\ell-1}_2 (t,x,v, v_1)\\
			 & \ \ \ \ \ := \sup\{ s<t^\ell_1:
			X^{\ell-1}(s;t^{\ell}_1 (t,x,v),x^{\ell}_1 (t,x,v),v_1) \in \p\O
			\}
			,\\
			&x^{\ell-1}_2 (t,x,v, v_1 ) := X^{\ell-1} (t^{\ell-1}_2 (t,x,v,v_1);t^\ell_1(t,x,v),x^\ell_1(t,x,v),v_1)
			,\\
		\end{split}
		\Ee
		and inductively 
		\Be
		\begin{split}\label{cycle_ell}
			& t^{\ell-(k-1)}_k (t,x,v, v_1, \cdots, v_{k-1})  \\
			&:=     \sup\big\{ s<t^{\ell-(k-2)}_{k-1} 
			:
			X^{\ell-1}(s;t_{k-1}^{\ell - (k-2)}  , x_{k-1}^{\ell - (k-2)} ,v_{k-1}) \in \p\O
			\big\},\\
			& x_k^{\ell - (k-1)} (t,x,v, v_1, \cdots, v_{k-1})\\
			&:= X^{\ell- (k-2)} (t_k^{\ell- (k-1)}; t_{k-1}^{\ell- (k-2)},x_{k-1}^{\ell- (k-2)} , v_{k-1})
			.
		\end{split}
		\Ee
		Here,
		\Be\begin{split}\notag
			t^{\ell-(i-1)}_{i } &:= t^{\ell-(i-1)}_{i }
			(t,x,v,v_1, \cdots, v_{i-1}),\\
			x^{\ell-(i-1)}_{i } &:= x^{\ell-(i-1)}_{i }
			(t,x,v,v_1, \cdots, v_{i-1}).\end{split}\Ee
		
		From (\ref{duhamel_local}) and (\ref{bdry_local}), we have 
		\Be
		\begin{split}\label{h_ell_local}
			&h^{\ell+1} (t,x,v)\\
			= & \ \mathbf{1}_{t_{1}^\ell   \leq 0}e^{-      \int^t_0     \nu^\ell 
			} h^{\ell+1} (0,X^\ell (0),V^\ell (0) )   \\
			&+ \int_{\max\{t_{1}^\ell ,0\}}^{t}e^{-  \int^t_s     \nu^\ell }|[K_{w_{\vartheta} }h^\ell  +w_{\vartheta} g^\ell](s,X^\ell(s;t,x,v), V^\ell(s;t,x,v))|\dd s \\
			&+ \mathbf{1}_{t_{1}^\ell  \geq  0}e^{-      \int^t_{t_1^\ell}     \nu^\ell 
			} h^{\ell+1}(t_1^\ell,X^\ell (t_{1}^\ell  ;t,x,v  ), V^\ell (t_{1}^\ell  ;t,x,v  )).
		\end{split}
		\Ee
		We define
		\begin{equation}
		\tilde{w}_{\vartheta} (v)\equiv \frac{1}{ w_{\vartheta}  (v)\sqrt{\mu (v)}}
		.\label{tweight}
		\end{equation}
		From (\ref{bdry_local}),
		\Be\begin{split}\notag
			&\textit{the last line of} \ (\ref{h_ell_local})\\
			 =  & \ \mathbf{1}_{t_{1}^\ell  \geq  0}e^{-      \int^t_{t_1^\ell}     \nu^\ell 
			}  \frac{1}{\tilde{w_{\vartheta}} (V^\ell (t^\ell_1))} \int_{n(x_1^\ell)\cdot v_1>0} h^{\ell}    (t^\ell_1, x^\ell_1, v_1)
			\tilde{w}_{\vartheta} (v_1) 
			c_\mu \mu  \{n(x_1^\ell) \cdot v_1\} \dd v_1.
		\end{split}\Ee
		
		We define $\mathcal{V}(x)=\{v \in\mathbb{R}^3 : n(x)\cdot v >0\}$
		with a probability measure $\dd\sigma=\dd\sigma(x)$ on $\mathcal{V}(x)$ which is given by
		\begin{equation}
		\dd\sigma \equiv
		c_\mu
		\mu
		(v)\{n(x)\cdot v\}\dd v.\label{smeasure}
		\end{equation}
		Let
		\Be\label{mathcal_V}
		\mathcal{V}_j: = \{v_j \in \R^3: n(x^{\ell - (j-1)}_j) \cdot v_j >0\}.
		\Ee
		
		Then inductively we obtain from (\ref{h_ell_local}), (\ref{duhamel_local}) and (\ref{bdry_local}), 
		\Be\begin{split}\label{h_iteration}
			&|h^{\ell+1} (t,x,v)| \\
			\leq & \ \mathbf{1}_{t_{1}^\ell   \leq 0}e^{-      \int^t_0     \nu^\ell 
			}|h^{\ell+1} (0,X^\ell (0),V^\ell (0))|    \\
			&  +
			\int_{\max\{t_{1}^\ell ,0\}}^{t}e^{-  \int^t_s     \nu^\ell }|[K_{w_{\vartheta} }h^\ell  +w_{\vartheta} g^\ell](s,X^\ell(s;t,x,v), V^\ell(s;t,x,v))|\dd s   \\
			&    +\mathbf{1}_{t_{1}^\ell >0}   
			\frac{
				e^{-  \int^t_{t_1^\ell}     \nu^\ell }
			}{\tilde{w_{\vartheta}}_{\varrho}(V^\ell (t_1^\ell))}\int_{\prod_{j=1}^{k-1}\mathcal{V}_{j}}|H| , 
		\end{split}\Ee
		where $|H|$ is bounded by
		\begin{align}
		 &\sum_{l=1}^{k-1}\mathbf{1}_{\{t^{\ell-l}_{l+1}\leq
			0<t_{l}^{\ell - (l-1)}\}}   |h^{\ell-l} (0,  
		X^{\ell-l}(0
		; v_l
		)
		,V^{\ell-l}(0; v_l))|\dd\Sigma _{l}(0) \label{h1} \\
		 +& \sum_{l=1}^{k-1}\int_{\max\{ t_{l+1}^{\ell-l}, 0 \}}^{t_l^{\ell - (l-1)}}\mathbf{1}_{\{t_{l+1}^{\ell-l}\leq
			0<t_{l}^{\ell - (l-1)}\}} \nonumber
		\\
		 & \ \ \   \times
		|[K_{w_{\vartheta} }h^{\ell-l} +w_{\vartheta}   g^{\ell-l}](s,
		X^{\ell-l}(s; v_l), V^{\ell-l}(s; v_l)
		|\dd \Sigma
		_{l}(s)\dd s  \label{h2} \\
		 +& \mathbf{1}_{\{0<t_{k}^{\ell - (k-1)}\}}|h^{\ell-(k-1)} (t_{k}^{\ell - (k-1)},x_{k}^{\ell - (k-1)},v_{k-1})|\dd\Sigma
		_{k-1}(t_{k}^{\ell - (k-1)}),  \label{h5}
		\end{align}%
		with
		\Be\begin{split}
			\label{measure} 
			d\Sigma _{l}^{k-1}(s) = \{\Pi _{j=l+1}^{k-1}\dd\sigma _{j}\}& \times\{
			e^{
				-\int^{t_l}_s
				\nu_\phi 
			}
			\tilde{w_{\vartheta}}(v_{l})\dd\sigma _{l}\}\\
			&\times \Pi
			_{j=1}^{l-1}\{{{
					e^{
						-\int^{t_j}_{t_{j+1}}
						\nu_\phi 
					}
					\dd\sigma _{j}}}\},
		\end{split} \Ee 
		and 
		\Be\begin{split}\label{X_ell_l}
			X^{\ell-l} (s; v_l)&:= X^{\ell-l}(s;t_l^{\ell- (l-1)},x_l^{\ell- (l-1)},v_l),\\
			V^{\ell-l} (s;v_l)&:= V^{\ell-l}(s;t_l^{\ell- (l-1)},x_l^{\ell- (l-1)},v_l).
		\end{split}\Ee

		\vspace{2pt}
		
		\textit{Step 2-2. } We claim that there exists $T>0$ and $k_{0} >0$ such that for all $k\geq  k_{0}$ and for all $(t,x,v) \in [0,T] \times \bar{\O} \times \R^{3}$, we have 
		\Be\label{small_k}
		\int_{\prod_{j=1}^{k-1} \mathcal{V}_{j}} \mathbf{1}_{\{ t^{k} (t,x,v,v^{1}, \cdots , v^{k-1}) >0 \}} \dd \Sigma_{k-1}^{k-1} \lesssim_{\O} \Big\{\frac{1}{2}\Big\}^{ k/5}.
		\Ee
		The proof of the claim is a modification of a proof of Lemma 14 of \cite{GKTT1}.
		
		For $0<\delta\ll 1$ we define
		\Be\label{V_j^delta}
		\mathcal{V}_{j}^{\delta} := \{ v^{j } \in \mathcal{V}^{j} : |v^{j} \cdot n(x^{j})| > \delta, \ |v^{j}| \leq \delta^{-1} \}.
		\Ee 
		
		Choose   
		\Be\label{large_T}
		T= \frac{2}{\delta^{2/3} (1+ \| \nabla \phi \|_\infty)^{2/3}}.
		\Ee
		We claim that 
		\Be\label{t-t_lowerbound}
		|t^j -t^{j+1}|\gtrsim   \delta^3, \ \ \text{for} \ v^j \in \mathcal{V}^\delta_j,  \ 0 \leq t\leq T, \ 0 \leq t^j.
		\Ee
		
		For $j \geq 1$
		\Bes
		&&\Big| \int^{t^{j+1}}_{t^{j}} V(s;t^{j}, x^{j}, v^{j}) \dd s   \Big|^{2}\\
		&=& |x^{j+1} -x^{j}|^{2}\\
		&\gtrsim& |(x^{j+1} -x^{j}) \cdot n(x^{j})|\\
		&=&\Big| \int^{t^{j+1}}_{t^{j }}  V(s;t^{j},x^{j},v^{j}) \cdot  n(x^{j})  \dd s 
		\Big|\\
		&=&\Big| \int_{t^{j}}^{t^{j+1}}  
		\Big(
		v^{j} - \int^{s}_{t^{j}} \nabla \phi (\tau, X(\tau;t^j,x^j,v^j))    \dd \tau
		\Big)\cdot   n(x^{j}) 
		\dd s 
		\Big|\\
		&\geq& |v^{j} \cdot n(x^{j})| |t^{j}-t^{j+1}|
		- \Big|
		\int^{t^{j+1}}_{t^{j }} \int^{s}_{t^{j }}  \nabla \phi (\tau, X(\tau; t^{j},x^{j},v^{j})) \cdot n(x^{j})\dd \tau \dd s 
		\Big|.
		\Ees
		Here we have used the fact if $x,y \in \p\O$ and $\p\O$ is $C^2$ and $\O$ is bounded then $|x-y|^2\gtrsim_\O |(x-y) \cdot n(x)|$. Hence
		\Be\label{lower_tb}
		\begin{split}
			&|v^{j} \cdot n(x^{j})| \\
			& \lesssim \frac{1}{|t^{j} - t^{j+1}| } \Big| \int^{t^{j+1}}_{t^{j}}
			V(s;t^{j},x^{j},v^{j}) \dd s
			\Big|^{2}\\
			& \ \  + 
			\frac{1}{|t^{j} - t^{j+1}| } \Big| \int^{t^{j+1}}_{t^{j}} 
			\int^{s}_{t^{j}} 
			\nabla \phi (\tau, X(\tau;t^{j},x^{j},v^{j}) )\cdot n(x^{j})\dd \tau  \dd s
			\Big|\\
			& \lesssim 
			|t^{j} - t^{j+1}|   \big\{ |v^j|^2 + |t^j - t^{j+1}|^3 \|\nabla \phi\|^2_\infty
			\\
			& \ \ \ \ \ \ \ \  \ \ \ \ \ \  \ \ \ \     +    \frac{1}{2}\sup_{t^{j+1} \leq \tau \leq t^{j}} 
			|  \nabla \phi (\tau, X(\tau;t^{j},x^{j},v^{j}) )\cdot n(x^{j})|
			\big\}.
		\end{split}
		\Ee
		For $v^j \in \mathcal{V}^\delta_j$, $0 \leq t\leq T$, and $t^j\geq 0$,
		\Be\notag
		|v^j \cdot n(x^j)|\lesssim |t^j -t^{j+1}|
		\{
		\delta^{-2} + T^3 \| \nabla \phi \|_\infty^2 + \| \nabla \phi \|_\infty
		\}.
		\Ee
		We choose $T$ as (\ref{large_T}) then prove (\ref{t-t_lowerbound}).

		Therefore if $t^{k}  \geq 0$ then there can be at most $\left\{\left[\frac{C_\Omega}{%
			\delta^3}\right]+1\right\}$ numbers of $v^m \in \mathcal{V}_m^\delta$ for $1\leq m
		\leq k-1$. Equivalently there are at least $k-2- \left[\frac{C_\Omega}{%
			\delta^3}\right]$ numbers of $v^{i} \in \mathcal{V}_{i} \backslash \mathcal{V}_{i}^\delta$ for $0 \leq i \leq m$. 
		
		Let us choose $k=N \times \left(\left[\frac{C_\Omega}{\delta^3}\right]%
		+1\right)$ and $N=\left(\left[\frac{C_\Omega}{\delta^3}\right]%
		+1\right)   \gg C>1$. Then we have 
		\begin{eqnarray*}
			&&\int_{\prod_{j=1}^{k-1}\mathcal{V}_j} 1_{\{t^{k}(t,x,v,v^{1},\cdots, v^{k-1})
				> 0\}} d\Sigma_{k-1}^{k-1} \\
			&\leq& \sum_{m=1}^{\left[\frac{C_\Omega}{\delta^3}\right]+1} \int_{\Big\{\begin{array}{ccc}{\text{there are exactly } m \text{ of } v_{i} \in \mathcal{V}_{i}^\delta} \\ {\small \text{ and } \ k-1-m \
						of \ v_{i} \in \mathcal{V}_{i} \backslash \mathcal{V}_{i}^{\delta}}\end{array}\Big\}}
			\prod_{j=1}^{k-1} C_0 \mu(v^j)^{1/4} \mathrm{d} v^j \\
			&\leq& \sum_{m=1}^{\left[\frac{C_\Omega}{\delta^3}\right]+1} \left(%
			\begin{array}{ccc}
				k-1  \\
				m
			\end{array}%
			\right) \left\{ \int_{\mathcal{V}}C_0 \mu(v)^{1/4}\mathrm{d} v\right\}^{m} \left\{ \int_{%
				\mathcal{V}\backslash \mathcal{V}^\delta}C_0 \mu(v)^{1/4}\mathrm{d} v\right\}^{k-1-m} \\
			&\leq& \left(\left[\frac{C_\Omega}{\delta^3}\right]+1\right) \{k-1\}^{\left[%
				\frac{C_\Omega}{\delta^3}\right]+1} \{ \delta\}^{k-2-\left[\frac{%
					C_\Omega}{\delta^3}\right]} \left\{ \int_{\mathcal{V}} C_0\mu(v)^{1/4}
			dv\right\}^{\left[\frac{C_\Omega}{\delta^3}\right]+1} \\
			&\leq& \{CN\}^{\frac{k}{N}} \left\{\frac{k}{N}\right\}^{\frac{k}{N}} \left\{\frac{k}{N}\right\}^{-\frac{k}{N} \frac{N^2}{20}}\\
			&\leq&  \left\{\frac{k}{N}\right\}^{\frac{k}{N} \left(- \frac{N^2}{20} + 3\right) } \\ 
		&	\leq&
			\left\{\frac{1}{2}\right\}^{ k}
			,
		\end{eqnarray*}
		where we have chosen $k=N \times \left(\left[\frac{C_\Omega}{\delta^3}\right]%
		+1\right)$ and $N=\left(\left[\frac{C_\Omega}{\delta^3}\right]%
		+1\right)   \gg C>1$.
		
		\vspace{2pt}
		
		\textit{Step 2-3. }  From (\ref{bound_g_ell}), (\ref{lower_nu_l}), and (\ref{h_iteration})-(\ref{measure}), and (\ref{small_k}), if we choose $\ell \geq k_0$ and $0 \leq t \leq T$ where $k_0$ and $T$ in (\ref{small_k})
		\Be
		\begin{split}\label{stochastic_h^ell}
			& | h^{\ell+1} (t,x,v )|\\
			\leq & \  C_k\| e^{- \frac{\nu_0 }{2}  t} h_0 \|_\infty \\
			& + \ \int^t_{\max\{t^\ell_1, 0 \}}
			e^{ - \frac{\nu_0}{2} (t-s)}\int_{\R^3}
			\mathbf{k}_\varrho (V^\ell (s ),u) |h^\ell (s, X^\ell (s ), u)|
			\dd u
			\dd s\\
			& +C_k \sup_l \int^{t_l^{\ell- (l-1)}}_{\max\{ t^{\ell-1}_{l+1}, 0  \}}   e^{ - \frac{\nu_0}{2} (t-s)}\int_{\R^3} \int_{\R^3} 
			\mathbf{k}_\varrho (V^{\ell-l} (s; v_l), u) 
			\\
			& \ \ \ \ \ \ \ \ \ \ \  \times
			|h^{\ell-l} (s, X^{\ell-l}(s; v_l)  ,u )|
			\{n(x_l)\cdot v_l\}\frac{\sqrt{\mu (v_l)}}{w_{\vartheta}(v_l)}
			\dd v_l
			\dd u
			\dd s \\
			& +  \int^t_{\max\{ t^\ell_{1}, 0  \}} 
			\langle V^\ell (s;t,x,v) \rangle     e^{- \int^t_s \frac{ \nu^\ell (\tau)}{2} \dd \tau}
			\| 
			e^{- \frac{\nu_0}{2} (t-s)}  h^\ell (s) \|_\infty^2
			\dd s \\
			&+ C_k \sup_l \int^{t_l^{\ell - (l-1)}}_{\max\{ t^{\ell-l}_{l+1} ,0 \}} \langle V^{\ell-l} (s;v_l) \rangle  \\
			& \ \ \ \ \ \ \ \ \ \ \  \times
			e^{
				- \int^{t_l^{\ell - (l-1)}}_s \frac{ \nu^{\ell-l } (\tau)}{2} \dd \tau 
			}
			\| e^{-\frac{\nu_0}{2} (t-s)}   h^{\ell-l} (s) \|_\infty^2 \dd s\\
			& +  \int^t_{\max\{ t^\ell_{1}, 0  \}} 
			\| 
			e^{- \frac{\nu_0}{2} (t-s)}  \nabla  \phi^\ell (s) \|_\infty 
			\dd s \\
			&+ C_k \sup_l \int^{t_l^{\ell - (l-1)}}_{\max\{ t^{\ell-l}_{l+1} ,0 \}}  
			\| e^{-\frac{\nu_0}{2} (t-s)}  \nabla   \phi ^{\ell-l} (s) \|_\infty  \dd s\\
			& +\Big\{\frac{1}{2}\Big\}^{ k/5} \| e^{- \frac{\nu_0 }{2}  (t - t_k^{\ell - (k-1)}) }  h(t_k^{\ell - (k-1)})  \|_\infty,
		\end{split}
		\Ee
		where we used the abbreviation of (\ref{X_ell_l}).
		
		From $\int^t_{0}\langle V^{\ell-l} (s;v_l) \rangle 
		e^{
			- \int^{t_l^{\ell - (l-1)}}_s \frac{ \nu^{\ell-l } (\tau)}{2} \dd \tau 
		} \dd s \lesssim 1$ and (\ref{grad_estimate}), we derive that 
		\Be\begin{split}\label{uniform_h^ell}
			\| h^{\ell+1} (t) \|_\infty \lesssim_k& \  \| h(0) \|_\infty + o(1) \| h (t_{k}^{\ell - (k-1)}) \|_\infty\\
			&
			+ t  \max_{l \geq 0 } \sup_{0 \leq s \leq t}   \| h^{\ell-l} (s) \|_\infty
			+ \max_{l \geq 0 } \sup_{0 \leq s \leq t}  \| h^{\ell-l} (s) \|_\infty^2 .
		\end{split}\Ee
		By taking supremum in $\ell$ and choosing $M \ll1$ and $0\leq t\leq T^* \leq T$ with $T^*\ll 1$, we conclude (\ref{uniform_h_ell}).
		
		{\color{red}
		}
		\vspace{4pt}
		
		\textit{Step 3. } We claim that there exist $T^{**} \ll 1$ (and $T^{**} < T^*(M)$) and $C>0$ such that the sequence $F^\ell = \mu+ \sqrt{\mu} f^\ell$ in (\ref{Fell}) satisfies 
		\Be \label{unif_Em}
		   \max_{ \ell\geq0}\sup_{0 \leq t \leq T^{**}} \mathcal{E}^\ell(t)
			\leq C \{ \| w_{\tilde{\vartheta}} f _0 \|_p^p
			+
			\| w_{\tilde{\vartheta}} \alpha_{f_{0 },\e}^\beta \nabla_{x,v} f  _0 \|_{p}^p \}< \infty.
		 \Ee 
		 where we define, for $0 < \e \ll 1$,
		\Be\label{mathcal_E}\begin{split}
		\mathcal{E}^{\ell+1}(t):=& \ 
		\| w_{\tilde{\vartheta}} f^{\ell+1} (t) \|_p^p
		+
		\| w_{\tilde{\vartheta}}\alpha_{f^{\ell},\e}^\beta \nabla_{x,v} f^{\ell+1}(t) \|_{p}^p\\&
		+ 
		\int^t_0  | w_{\tilde{\vartheta}}\alpha_{f^{\ell},\e}^\beta \nabla_{x,v} f^{\ell+1}(t) |_{p,+}^p.
		\end{split}\Ee
		\hide
		For proving (\ref{W1p_local_bound}), we recall the sequence defined in (\ref{Fell})
		\Be \begin{split}\label{Fell_1}
			\p_t F^{\ell+1} + v\cdot \nabla_x F^{\ell+1} - \nabla {\phi^{\ell }} \cdot \nabla_v F^{\ell+1}= Q_{\text{gain}} (F^{\ell}, F^{\ell}) - Q_{\text{loss}} (F^{\ell}, F^{\ell+1}),\\
			- \Delta \phi^{\ell }= \int_{\R^3} F^{\ell } \dd v - \rho_0, \ \ \int_\O \phi^{\ell } \dd x =0, \ \ \frac{\p \phi^{\ell }}{\p n}\Big|_{\p\O }=0,
		\end{split}\Ee
		with (\ref{F_ell_BC}) and $F^{\ell+1} (0,x,v) = F_0 (x,v)$.  \\
		\unhide
		\hide
  	{\color{red} 
		We can solve (\ref{Fell}) and (\ref{F_ell_BC})  from \textit{Step 1}-\textit{Step 4} by considering (\ref{F_ell_BC}) as the in-flow boundary condition with a fixed external potential and using another sequences as 
		\Be \begin{split}\label{Fell_2}
			\p_t F^{\ell+1,m+1} + v\cdot \nabla_x F^{\ell+1,m+1} - \nabla {\phi^{\ell  }} \cdot \nabla_v F^{\ell+1,m+1}\\
			= Q_{\text{gain}} (F^{\ell+1,m}, F^{\ell+1,m}) - Q_{\text{loss}} (F^{\ell+1,m}, F^{\ell+1,m+1}),\\
			- \Delta \phi^{\ell }= \int_{\R^3} F^{\ell } \dd v - \rho_0, \ \ \int_\O \phi^{\ell } \dd x =0, \ \ \frac{\p \phi^{\ell }}{\p n}\Big|_{\p\O }=0,\\
			F^{\ell+1, m+1}|_{\gamma_-} = c_\mu \mu \int_{n \cdot u>0} F^{\ell} \{n \cdot u\} \dd u.
		\end{split}\Ee

		Following \textit{Step 1}-\textit{Step 5}, it is straightforward to prove $L^\infty$ uniform bound of $f^{\ell+1, m+1} = \{ F^{\ell+1, m+1} - \mu\}/ \sqrt{\mu}$ as (\ref{uniform_h_ell}). Then we obtain $(F^{\ell+1}, \phi^\ell)$ solving (\ref{Fell}) as a limit of the sequence.    \\
			} 
		
		\unhide
		It suffices to prove the following induction statement: there exist $T^{**} \ll 1$ (and $T^{**} < T^*(M)$) and $C>0$ such that  
		\Be \label{induc_hypo}
		\begin{split}
			\text{ if }   \   \max_{0\leq m\leq\ell}\sup_{0 \leq t \leq T^{**}} \mathcal{E}^m(t)
			\leq C \{ \| w_{\tilde{\vartheta}} f _0 \|_p^p
			+
			\| w_{\tilde{\vartheta}} \alpha_{f_{0 },\e}^\beta \nabla_{x,v} f  _0 \|_{p}^p \}< \infty\\
			\text{ then } \  \sup_{0 \leq t \leq T^{**}}
			\mathcal{E}^{\ell+1} (t)\leq C \{ \| w_{\tilde{\vartheta}} f _0 \|_p^p
			+
			\| w_{\tilde{\vartheta}} \alpha_{f_{0 },\e}^\beta \nabla_{x,v} f_0 \|_{p}^p \}. 
		\end{split} \Ee 
		
		The proof is similar to the proof of Proposition \ref{prop_W1p}. By taking derivatives $\p \in\{\nabla_x,\nabla_v \}$ to $f^{\ell+1} = \{F^{\ell+1} - \mu\} / \sqrt{\mu}$ to (\ref{Fell}) 
		\Be\label{eqtn_nabla_f_ell}
		[\p_t + v\cdot \nabla_x - \nabla \phi^{\ell } \cdot \nabla_v+ \nu 
		+ \frac{v}{2} \cdot \nabla \phi^{\ell } 
		] \p f^{\ell+1} = \mathcal{G}^{\ell+1},
		\Ee
		where $\mathcal{G}^{\ell+1}$ equals 
		\Be
		\begin{split}	\label{mathcal_G_ell}
			\mathcal{G}^{\ell+1}  =& 
			- \p v \cdot\nabla_{x}f^{\ell+1} 
			+ \p \nabla \phi^{\ell}\cdot \nabla_v f^{\ell+1}\\&
			+ \p \Big\{  \Gamma(f^{\ell},f^{\ell})  +  \Gamma_{\text{loss}}(f^{\ell}, f^{\ell}) - \Gamma_{\text{loss}}(f^{\ell}, f^{\ell+1})  \Big\}   \\
			&- \p \big[\nu(v) + \frac{v}{2} \cdot \nabla \phi^{\ell}(t,x)\big] f^{\ell+1} - \p K f^{\ell} - \p(v\cdot\nabla_{x}\phi^{\ell}\sqrt{\mu}).
		\end{split}\Ee  
	\vspace{1cm}
	The boundary condition is given by (\ref{BC_f}) but replacing $f$ in the LHS by $f^{\ell+1 }$ and $f$ in the RHS by $f^\ell$. 
	
	Recall (\ref{nu_phi}). From (\ref{uniform_h_ell}) we have $\nu_{\phi^{\ell}}\geq \frac{\nu(v)}{2}$. We can follow the proof of (\ref{pf_green}), (\ref{pf_p_1}), (\ref{pf_p_2}) from Lemma \ref{lem_Green}, and (\ref{Gvgain}), (\ref{Gvloss}), (\ref{bound_Gamma_nabla_vf1}), (\ref{bound_Gamma_nabla_vf2}) to obtain 
		\Be \label{mathcalE}
		\begin{split}
			\mathcal{E}^{\ell+1}(s)
			\lesssim&  \   \| w_{\tilde{\vartheta}}\alpha_{f^{\ell },\e}^{\b} 
			\nabla_{x,v} f_0\|_p^p+ 
			\underbrace{
				\int^t_0 |w_{\tilde{\vartheta}}\alpha_{f^{\ell },\e}^{\b} 
				\nabla_{x,v} f^{\ell+1}|_{p,-}^p
			}_{(\ref{mathcalE})_{\text{boundary}}}  
			\\
			& +  t \big(1+ \|w_{\vartheta}f^{\ell}\|_{\infty} +  \|w_{\vartheta}f^{\ell+1}\|_{\infty} +\| \nabla^2 \phi^{\ell } \|_\infty
			\big) \\
			& \ \ \ \ \ \ \times \max_{m=\ell, \ell+1} \sup_{0 \leq s \leq t}\mathcal{E}^{m}(s).
		\end{split}
		\Ee
		There is only one extra term to (\ref{mathcal_G}) which is $ \p \big\{ \Gamma_{\text{loss}}(f^{\ell}, f^{\ell}) - \Gamma_{\text{loss}}(f^{\ell}, f^{\ell+1}) \big\} $ in (\ref{mathcal_G_ell}). Its contribution in $(\ref{pf_green})_{\mathcal{G}}$ can be controlled by   
		\Be\begin{split}\notag
			&(1 + \sup_{0 \leq s\leq  t} \|w_\vartheta f^{\ell}(s)\|_{\infty} + \sup_{0 \leq s\leq  t} \|w_\vartheta f^{\ell+1}(s)\|_{\infty} ) \\
			& \times  \int^t_0      \iint_{\O\times\R^{3}}  |\alpha_{f^{\ell},\e}^{\b  } w_{\tilde{\vartheta}}  \p f^{\ell+1}(v)|^{p-1} \\
			& \times  
				\int_{\R^3} \alpha_{f^{\ell},\e}(v)^\beta \mathbf{k}_\varrho(v,u)w_{\tilde{\vartheta}} (v) \big\{  |\p f^{\ell}(u)| + |\p f^{\ell+1}(u)| \big\} \dd u  
				\dd v \dd x   \dd s,
		\end{split}\Ee
		and then further bounded as, by the way to treat (\ref{pf_p_K}) and (\ref{pf_p_0}),
		\Be\notag
		\sup_{0 \leq s\leq t} \max_{m=\ell,\ell+1}\| w_\vartheta f^{m}(s) \|_\infty \int^t_0  \max_{m=\ell,\ell+1}\| \nu_{\phi_{f^{\ell}}, w_{\tilde{\vartheta}}}^{1/p} w_{\tilde{\vartheta}}   \alpha_{f^{\ell},\e}^\beta \p f^{\ell+1} \|_p^p  .
		\Ee
		Clearly it is bounded the RHS of (\ref{mathcalE}).

	Now we focus on $(\ref{mathcalE})_{\text{boundary}}$. Note that $\alpha_{f^{\ell }, \e}= |n\cdot v|$ on $\gamma_-$. On $\gamma_-$, $|\p f^{\ell+1}|$ is bounded by (\ref{BC_deriv}) replacing $|f|, |\nabla_{x,v} f|, |\nabla \phi_f|$ by $|f^{\ell}|+ |f^{\ell+1}|, |\nabla_{x,v} f^\ell | + |\nabla_{x,v} f^{\ell+1} |, |\nabla \phi^\ell| + |\nabla \phi^{\ell-1}|$. Now we follow \textit{Step 4} in the proof of Proposition \ref{prop_W1p}. We have the same decomposition (\ref{W1p_bdry}) replacing all $\nabla_{x,v} f$ by $\nabla_{x,v} f^\ell$. Then we obtain (\ref{W1p_bdry_1}) but replacing all $\nabla_{x,v} f$ by $\nabla_{x,v} f^\ell$. Then we apply (\ref{case:nondecay}) of Lemma \ref{le:ukai} to $\gamma_+ (x) \backslash \gamma^\e _+ (x)$ contribution in (\ref{W1p_bdry}). And this leads (\ref{nongrazing_nabla_f}) but replacing all $\nabla_{x,v}f,f ,\phi_f$ by $\nabla_{x,v} f^\ell, f^\ell,\phi^{\ell-1}$. Finally we end up with 
		\Be
		\begin{split}\label{control_bdry_ell}
			(\ref{mathcalE})_{\text{boundary}}  
			\lesssim     &  \ 
			\| w_{\tilde{\vartheta}}\alpha^\beta_{f_0, \e} \nabla_{x,v} f_0 \|_p^p  +o(1) \max_{0 \leq m \leq \ell}\sup_{0 \leq s \leq t}\mathcal{E}^{m}(s)\\
			&+t (1+  \|w_{\vartheta}f^{\ell }\|_{\infty} +\| \nabla^2 \phi^{\ell-1 } \|_\infty) \max_{0 \leq m \leq \ell} \sup_{0 \leq s \leq t}
			\mathcal{E}^{m}(s).
		\end{split}\Ee
		
		On the other hand, from Lemma \ref{lemma_apply_Schauder}, 
		\Be\label{phi_ell<E}
		\| \nabla^2 \phi^\ell (t)\|_\infty + \| \nabla^2 \phi^{\ell-1}(t) \|_\infty 
		\lesssim [\mathcal{E}^{\ell}(t) + \mathcal{E}^{\ell-1}(t)]^{1/p}.
		\Ee
		
		From (\ref{mathcalE}), (\ref{control_bdry_ell}), (\ref{phi_ell<E}), and the induction hypothesis in (\ref{induc_hypo}), we derive 
		\Be\notag
		\sup_{m\leq \ell+1}\sup_{0 \leq s\leq t}\mathcal{E}^{m} (s) \lesssim \| w_{\tilde{\vartheta}} \alpha_{f_0,\e}^{\b} 
		\nabla_{x,v} f_0\|_p^p + t \sup_{0 \leq s\leq t}\sup_{m\leq \ell}\mathcal{E}^{m} (s).
		\Ee
		Finally we choose small $T^{**}>0$ and conclude (\ref{induc_hypo}).
		
		{\color{red}
	} 
		
		\hide

		Following \textit{Step 3} in the proof of Proposition \ref{prop_W1p} we obtain (\ref{mid}) replacing $\alpha_f,f, \p f$ by $\alpha_{f^\ell}, f^\ell, \p f^\ell$.
		
		Following \textit{Step 1} in the proof of Proposition we obtain (\ref{Lp_estimate_f})
		\Be
		\begin{split}
			&\| f^{\ell+1} (t)\|_p^p+ \int^t_0 \| \nu_{\phi^{\ell+1}}^{1/p} f^{\ell+1} \|_p^p
			+ \int^t_0 |f^{\ell+1}|_{p,+}^p\\
			\lesssim &
			\| f(0) \|_p^p + 
			o(1) \int^t_0 |f^{\ell }|^{p}_{p,+}  
			+ (1+ \sup_\ell\| w_{\vartheta} f^{\ell} \|_\infty)\{ \int^t_0 \| f ^\ell\|_p^p +\int^t_0 \| f ^{\ell+1}\|_p^p \}.
		\end{split}
		\Ee
		\hide
		
		First we bound $(\ref{mathcalE})_\mathbf{k}$. Following (\ref{est_k_p_f}) and (\ref{small result}), we bound $\{|u|\geq N\}$ contribution of $(\ref{mathcalE})_\mathbf{k}$ by 
		\Be
		\int^t_0 \| \alpha^{\beta}_{f^\ell,\e}   \p f^{\ell+1} \|_p^p
		\Ee

		Since $\alpha$ depend on $f$, from (\ref{eqtn_nabla_f_ell}) and (\ref{alpha_invariant}),
		\Be \notag
		\begin{split}
			& |\alpha^\beta_{f^{\ell}}  \p f^{\ell+1} |^{p-1} \big[\p_t + v\cdot \nabla_x - \nabla_x \phi_f \cdot \nabla_v 
			\big] |\alpha^\beta_{f^{\ell}}  \p f^{\ell+1} | \\
			&= \alpha^{\b p}_{f^{\ell}} |\p f^{\ell+1}|^{p-1} \big( \nu(\sqrt{\mu} f^{\ell})\p f^{\ell+1} - \nu_{\phi^{\ell}}\p f^{\ell+1} + \mathcal{G}^{\ell+1} \big) ,\quad \nu_{\phi^{\ell}} :=   \nu(v) + \frac{v}{2} \cdot \nabla \phi^\ell.
		\end{split}
		\Ee
		
		If we define
		\Be\begin{split}
			\mathcal{E}_{m}(t) &:= \| w_{\vartheta}f^m (t) \|^{p}_{L^\infty (\bar{\O} \times \R^3)} + \| \alpha_{f^{\ell}}^{\b} 
			\p f^{\ell+1} (t)\|^{p}_p
			+ \int^t_0 \| \nu_{\phi^{\ell}}^{1/p} \alpha_{f^{\ell}}^\beta 
			\p f^{\ell+1}   \|^{p}_p  
			+ \int^t_0 |\alpha_{f^{\ell}}^{\b} 
			\p f^{\ell+1}|^{p}_{p,+}  ,  \\
		\end{split}\Ee
		from (\ref{Gvgain}), (\ref{Gvloss}), (\ref{bound_Gamma_nabla_vf1}), (\ref{bound_Gamma_nabla_vf2}), (\ref{BDD h}), and Lemma \ref{lem_Green}, we obtain 
		
		Now we consider supremum in time interval $t\in[0,T^{*}]$ and define
		\Be \label{E def}
		E_{\ell+1} := \sup_{t\in[0,T^{*}]}\mathcal{E}_{\ell+1}(t).
		\Ee
		For underbraced $(A)$, $(B)$, $(C)$, and $(D)$ in (\ref{mathcalE}), we can use estimates of $(\ref{pf_green})_{2}$, $(\ref{pf_green})_{3}$, $(\ref{pf_green})_{4}$, and $(\ref{pf_green})_{5}$ in section 4. From (\ref{151_2}), (\ref{151_3}), (\ref{151_4}), and (\ref{151_5}),  
		\Be \label{AtoD}
		\begin{split}
			(C) &\lesssim E_{\ell+1},  \\
			(D) &\lesssim E_{\ell+1}^{\frac{p-1}{p}},  \\
			(B) &\lesssim E_{\ell} + E_{\ell+1},  \\
			(A) &\lesssim \|\alpha^{\b}\p f (0) \|_{p}^{p} + o(1)E_{\ell} + o(1)E_{\ell+1} + o(1)(1 + E_{\ell}^{1/p} + E_{\ell+1}^{1/p} )(E_{\ell} + E_{\ell+1})  \\
			&\quad + (1 + \|\phi_{f}\|_{C^{2}})^{p} T^{*} E_{\ell}^{\frac{p-1}{p}}
		\end{split}
		\Ee
		
		\hide
		\Be  
		\begin{split}
			& (\ref{pf_green})_2 = \int_{0}^{t} |\alpha^\beta 
			\p f(s)  |_{p,-}^p ds   \\
			&\lesssim o(1) (1 + \|\phi_{f}\|_{C^{2}})^{p} \int^t_0 |\alpha^\beta \p f (s) |_{p,+}^p ds  + o(1) (1+\|w_{\vartheta}f\|_{\infty}) \int_{0}^{t} \| \nu_\phi^{1/p} \alpha^\beta \p f  (s)\|_{p}^{p}  \\
			&\quad + \int_{0}^{t} (1 + \|w_{\vartheta}f\|_{\infty})^{2p}(1 + \|\phi_{f}\|_{C^{2}})^{p} ( 1 + \|\alpha^{\b}\p f\|_{p}) \|\alpha^{\b}\p f\|_{p}^{p-1}  ds \\
		\end{split}
		\Ee
		\unhide
		
		From (\ref{AtoD}), (\ref{mathcalE}) gives
		\Be
		\begin{split}
			E_{\ell+1} 
			&\lesssim \|h_{0}\|^{p}_{\infty} +  \| \alpha_{f^{\ell}}^{\b} \nabla_{x,v} f^{\ell+1} (0)\|_p^p  \\
			&\quad + T^{*}E_{\ell}^{\frac{1}{p}} E_{\ell+1} + (1 + \|\nabla_{x}^{2}\phi^{\ell}\|_{\infty}) T^{*} E_{\ell+1} \\
			&\quad + (1 + \|\nabla_{x}^{2}\phi^{\ell}\|_{\infty}) T^{*} E_{\ell+1} + (1 + E_{\ell}^{\frac{1}{p}})(1 + \|\nabla_{x}^{2}\phi^{\ell}\|_{\infty}) T^{*} E_{\ell+1}^{\frac{p-1}{p}}  \\
			&\quad + (1 + E_{\ell}^{\frac{1}{p}} + E_{\ell+1}^{\frac{1}{p}}) T^{*} (E_{\ell} + E_{\ell+1})  \\
			&\quad + \|\alpha_{f^{\ell}}^{\b} \nabla_{x,v} f (0) \|_{p}^{p} + o(1)E_{\ell} + o(1)E_{\ell+1} + o(1)(1 + E_{\ell}^{1/p} + E_{\ell+1}^{1/p} )(E_{\ell} + E_{\ell+1})  \\
			&\quad + (1 + \|\phi_{f}\|_{C^{2}})^{p} T^{*} E_{\ell}^{\frac{p-1}{p}} . \\
		\end{split}
		\Ee
		
		\hide
		\Be
		\begin{split}
			E_{\ell+1} 
			&\lesssim \|h_{0}\|^{p}_{\infty} +  \| \alpha_{f^{\ell}}^{\b} \nabla_{x,v} f^{\ell+1} (0)\|_p^p  \\
			&\quad + T^{*}E_{\ell}^{\frac{1}{p}} E_{\ell+1} + (1 + \|\nabla_{x}^{2}\phi^{\ell}\|_{\infty}) T^{*} E_{\ell+1} \\
			&\quad + (1 + \|\nabla_{x}^{2}\phi^{\ell}\|_{\infty}) T^{*} E_{\ell+1} + (1 + E_{\ell}^{\frac{1}{p}})(1 + \|\nabla_{x}^{2}\phi^{\ell}\|_{\infty}) T^{*} E_{\ell+1}^{\frac{p-1}{p}}  \\
			&\quad + (1 + E_{\ell}^{\frac{1}{p}} + E_{\ell+1}^{\frac{1}{p}}) T^{*} (E_{\ell} + E_{\ell+1})  \\
			&\quad + \|\alpha_{f^{\ell}}^{\b} \nabla_{x,v} f (0) \|_{p}^{p} + o(1)E_{\ell} + o(1)E_{\ell+1} + o(1)(1 + E_{\ell}^{1/p} + E_{\ell+1}^{1/p} )(E_{\ell} + E_{\ell+1})  \\
			&\quad + (1 + \|\phi_{f}\|_{C^{2}})^{p} T^{*} E_{\ell}^{\frac{p-1}{p}} . \\
		\end{split}
		\Ee
		\unhide
		
		\noindent From Young's inequality and assumption $E_{\ell} \leq M \ll 1$,
		\Be \label{ineq}
		\begin{split}
			E_{\ell+1} 
			&\lesssim E(0) + (1 + \|\phi^{\ell}\|_{C^{2}}) ( T^{*} + \Lambda(E_{\ell}^{1/p}) + o(1)  ) E_{\ell+1}  \\
			&\quad +  (1 + \|\phi^{\ell}\|_{C^{2}})^{p} \big[ (T^{*}+o(1)) E_{\ell} + T^{*} \big], 
		\end{split}
		\Ee
		where $\Lambda(E_{\ell}^{1/p})$ is a polynomial of $E_{\ell}^{1/p}$ which goes to zero as $E_{\ell}^{1/p} \rightarrow 0$. From (\ref{ineq}), we can choose sufficiently small $T^{*} \ll 1$, $M \ll 1$, and $E(0) \ll 1$ so that $E_{\ell+1} \leq M$.

		\hide
		\Be\label{fell_local}
		\begin{split}
			&[\p_t + v\cdot\nabla_x - \nabla_x \phi^{\ell} \cdot \nabla_v + \nu + \frac{v}{2} \cdot\nabla \phi^\ell  ] f^{\ell+1}\\
			=&  \ Kf^{\ell}- v\cdot \nabla \phi^\ell \sqrt{\mu} + \Gamma_{\text{gain}} (f^\ell, f^\ell) -   \Gamma_{\text{loss}} (f^\ell, f^{\ell+1}),\\
			& - \Delta \phi^\ell = \int_{\R^3} f^\ell \sqrt{\mu} \dd v , \ \ \int_\O \phi^\ell \dd x =0, \ \ \frac{\p \phi^\ell}{\p n}\Big|_{\p\O }=0,
		\end{split}
		\Ee

		\unhide

		\hide

		We claim that if $\max_{0\leq m\leq\ell+1}\sup_{0 \leq t \leq T*}\| w_{\vartheta}f^m (t) \|_\infty \leq M\ll 1$ then there exits $T^{**}  \ll1$ such that  
		\Be\begin{split}
			\text{if}   \ \ 
			\max_{0\leq m\leq\ell}\sup_{0 \leq t \leq T^{**}}
			\Big\{
			\| \alpha^\beta \nabla_{x,v} f^m(t) \|_{p}
			+ 
			\int^t_0 \| \alpha^\beta \nabla_{x,v} f^m(t)\|_{\gamma, p}
			\Big\} \leq L,
		\end{split}\Ee
		then 
		\Be\begin{split}
			\sup_{0 \leq t \leq T^{**}}
			\Big\{
			\| \alpha^\beta \nabla_{x,v} f^{\ell+1}(t) \|_{p}
			+ 
			\int^t_0 \| \alpha^\beta \nabla_{x,v} f^{\ell+1} (t)\|_{\gamma, p}
			\Big\} \leq L,
		\end{split}\Ee
		
		\unhide\unhide

		\unhide
		
		\vspace{4pt}

		\textit{Step 4. } Now we claim that the sequence $F^\ell = \mu+ \sqrt{\mu} f^\ell$ in (\ref{Fell}) satisfies 
		\Be\label{31_local_bound_m}
		\sup_{\ell} \sup_{0 \leq t \leq T^{**}}\| \nabla_v f^\ell (t) \|_{L^3_x L^{1+ \delta}_v}\lesssim 1. 
		\Ee 
		To prove this, similar as above step, we claim that 
		\Be\label{induction_31}\begin{split}
		 {if} \ \ \sup_{m \leq \ell} \sup_{0 \leq t \leq T^{**}}\| \nabla_v f^m (t) \|_{L^3_x L^{1+ \delta}_v}\lesssim 1,\\
		 \quad  {then}\quad \sup_{0 \leq t \leq T^{**}}\| \nabla_v f^{\ell+1} (t) \|_{L^3_x L^{1+ \delta}_v}\lesssim 1.
	\end{split}	\Ee
		 Recall that $\nabla_v f^{\ell+1}$ solves (\ref{eqtn_nabla_f_ell}) and (\ref{mathcal_G_ell}) with $\p= \nabla_v$.
		\hide\Be\begin{split}\label{eqtn_g_v_ell}
			&[\p_t + v\cdot \nabla_x - \nabla_x \phi_{f^{\ell}} \cdot \nabla_v + \nu(v) + \frac{v}{2} \cdot \nabla_x \phi_{f^{\ell}}  ] \p_v f^{\ell+1} \\
			=&  -\p_x f^{\ell+1}- \frac{1}{2} \p_x \phi_{f^{\ell}} f^{\ell+1} - \p_v \nu f^{\ell+1} + \p_v(Kf^{\ell})   \\
			&\quad + \p_v \big\{ \Gamma (f^{\ell},f^{\ell}) +  \Gamma_{\text{loss}} (f^{\ell},f^{\ell}) - \Gamma_{\text{loss}} (f^{\ell},f^{\ell+1})   \big\} + |\nabla_x \phi_f| \langle v\rangle^2\sqrt{\mu} , 
		\end{split}\Ee
	similar as	(\ref{eqtn_g_v}). 
	Then similar as (\ref{g_initial})-(\ref{g_infty}), we obtain 
	\begin{equation}
	\begin{split}
	&|\p_v f^{\ell+1}(t,x,v)|\nonumber
	\\
	&\leq   \mathbf{1}_{\tb(t,x,v)> t}  
	|\p_v f^{\ell+1}(0,X(0;t,x,v), V(0;t,x,v))| \\
	& +   \ \mathbf{1}_{\tb(t,x,v)<t}
	\mu(\vb)^{\frac{1}{4}}  \int_{n(\xb) \cdot u>0} 
	| f^{\ell+1}(t-\tb, \xb, u) |\sqrt{\mu} \{n(\xb) \cdot u\} \dd u  \\
	&  + \int^t_{\max\{t-\tb, 0\}} 
	|\p_x f^{\ell+1}(s, X(s;t,x,v),V(s;t,x,v))|
	\dd s  \\
	&   + \int^t_{\max\{t-\tb, 0\}} 
	(1+ \| w_{\vartheta} f^{\ell} \|_\infty + \| w_{\vartheta} f^{\ell+1} \|_\infty)
	\int_{\R^3} \mathbf{k}_\varrho (V(s),u) \{ |\p_v f^{\ell}(s,X(s),u)| + |\p_v f^{\ell+1}(s,X(s),u)| \}  \dd u 
	\dd s  \\
	& + 
	\int^t_{\max\{t-\tb, 0\}} 
	|\nabla_x \phi_{f^{\ell}} (s, X(s;t,x,v))| \mu ^{1/4}
	\dd s   \\
	&+
	\int^t_{\max\{t-\tb, 0\}}  
	(1+ \delta_1)  \delta_1  |{w}_{\tilde{\vartheta}}( V(s;t,x,v))|^{-1} 
	\dd s 	,
		\end{split}
	\end{equation}\unhide
	Then we follow exactly the same proof of (\ref{tilde_w_integrable})-(\ref{est_g_K}) and conclude that 
		\Be\begin{split}\notag
			&\| \nabla_v f^{\ell+1} (t) \|_{L^3_x L^{1+ \delta}_v}\\
			\lesssim& \ \| w_{\tilde{\vartheta}} \p_v f^{\ell+1}(0)\|_{L^3_{x,v}} + \sup_{0\leq s \leq t}\{ \| w_\vartheta f^\ell(s)\|_\infty+  \| w_\vartheta f^{\ell+1}(s)\|_\infty \}\\
			& + t \sup_{0\leq s \leq t} \big\{\| f^{\ell}(s)\|_2  +\| w_{\tilde{\vartheta}} \alpha^\beta_{f^{\ell-1}, \e} \nabla_{x,v} f^{\ell} \|_p  +\| w_{\tilde{\vartheta}} \alpha^\beta_{f^{\ell}, \e} \nabla_{x,v} f^{\ell+1} \|_p   \\
			&\quad\quad \quad\quad\quad  + \| \nabla_v f^{\ell} (s) \|_{L^3_x L^{1+ \delta}_v} + \| \nabla_v f^{\ell+1} (s) \|_{L^3_x L^{1+ \delta}_v}  \big\}.
		\end{split}\Ee
		From (\ref{induc_hypo}) and (\ref{uniform_h_ell}), we can choose $T^{**}\ll 1$ and conclude (\ref{induction_31}) and hence (\ref{31_local_bound_m}).

		{\color{red}
	} 
		\vspace{4pt}
		
		\textit{Step 5. } Now we claim that the convergence of whole sequence 
		\Be\label{Cauchy_L1+}
		f^{\ell}\rightarrow f \ \text{strongly in } \ L^{1+}(\O \times \R^3) \cap L^{1+}(\gamma).  
		\Ee
		
		Note that $f^{\ell+1} - f^\ell$ satisfies $(f^{\ell+1} - f^\ell)|_{t=0} \equiv0$ and (\ref{eqtn_f-g}) with $f=f^{\ell+1}$ and $g= f^\ell$ to get
		\begin{equation}\notag
		\begin{split}
		& \p_t [f^{\ell+1}-f^{\ell}] + v\cdot \nabla_x [f^{\ell+1}-f^{\ell}] - \nabla_x  \phi_{f^{\ell}}  \cdot \nabla_v [f^{\ell+1}-f^{\ell}]\\
		&+ \frac{v}{2} \cdot \nabla_x \phi_{f^{\ell}}[f^{\ell+1}-f^{\ell}]+  \nu [f^{\ell+1}-f^{\ell}] \\
		&=  \nabla_x \phi_{f^{\ell}-f^{\ell-1}} \cdot \nabla_v f^{\ell-1} + K[f^{\ell}-f^{\ell-1}]  - \frac{v}{2} \cdot \nabla_x \phi_{f^{\ell}-f^{\ell-1}} f^{\ell-1}    \\ & - v\cdot \nabla_x \phi_{f^{\ell}-f^{\ell-1}} \sqrt{\mu} 
		 +  \{ \Gamma(f^{\ell},f^{\ell}) - \Gamma(f^{\ell-1},f^{\ell-1}) \}\\
		 &  + \{ \Gamma_{\text{loss}}(f^{\ell},f^{\ell}) - \Gamma_{\text{loss}}(f^{\ell-1},f^{\ell-1}) \}   - \{ \Gamma_{\text{loss}}(f^{\ell},f^{\ell+1}) - \Gamma_{\text{loss}}(f^{\ell-1},f^{\ell})  \} .
		\end{split}\end{equation} 
		\hide
		\Be\begin{split}\label{f_ell-ell}
			&[\p_t + v\cdot \nabla_x - \nabla_x \phi^\ell \cdot \nabla_v+ \nu_{\phi^\ell}] (f^{\ell+1} - f^\ell)\\
			= & \ 
			(\nabla \phi^{\ell }-  \nabla \phi^{\ell-1}) \cdot \nabla_v f^\ell
			+K(f^{\ell } - f^{\ell-1}) - v\cdot \nabla (\phi^{\ell} - \phi^{\ell-1})\sqrt{\mu}\\
			& + w \Gamma_{\text{gain}} (\frac{f^\ell}{w}, \frac{f^\ell - f^{\ell-1}}{w})+ w \Gamma_{\text{gain}} (\frac{f^\ell- f^{\ell-1}}{w}, \frac{f^{\ell-1}}{w})\\
			&- w \Gamma_{\text{loss}} ( \frac{f^\ell}{w}, \frac{f^{\ell+1} - f^\ell}{w})-w \Gamma_{\text{loss}} ( \frac{f^\ell - f^{\ell-1}}{w}, \frac{ f^\ell}{w}).
		\end{split}\Ee\unhide
	\hide	From (\ref{31_local_bound_m}), 
			\Be\label{bound_nabla_v_g_global_ell}
			\| \nabla_v (f^{\ell+1} - f^{\ell})(t) \|_{L^3_x (\O)L^{1+\delta}_v(\R^3)} \lesssim_t 1 \ \ \text{for all } \ t\geq 0,
			\Ee
			which corresponds to (\ref{bound_nabla_v_g_global}). \unhide
	With help of uniform estimatjes (\ref{uniform_h_ell}) and (\ref{31_local_bound_m}), we can modify RHSs of (\ref{f-g_energy}), (\ref{gronwall_f-g}), (\ref{gronwall_f-g_last}), and (\ref{gamma_+,1+delta}) in the proof of Proposition \ref{L1+stability} to get
		\Be \notag
		\begin{split}
		&\sup_{0 \leq s \leq t}\| f^{\ell+1} (s)- f^\ell (s) \|_{L^{1+ \delta}(\O \times \R^3)} + 
		\Big(\int^t_0 | f^{\ell+1} (s)- f^\ell (s)  |_{L^{1+ \delta}(\gamma)} \Big)^{\frac{1}{1+ \delta}}
		 \\
		&\lesssim  t \big( 1 + \sup_{0\leq s \leq t }\sup_{m}  \| h^{m}(s)\|_{\infty} +   \sup_{0 \leq s \leq t} \sup_{m} \| \nabla_v f^m (t) \|_{L^3_x L^{1+ \delta}_v}  \big)\\
		&   \ \ \ \times  \sup_{0 \leq s \leq t}\| f^{\ell } (s)- f^{\ell-1} (s) \|_{L^{1+ \delta}(\O \times \R^3)}  \\
		&\lesssim  O(t)\sup_{0 \leq s \leq t}\| f^{\ell } (s)- f^{\ell-1} (s) \|^{1+ \delta}_{L^{1+ \delta}(\O \times \R^3)}.
		\end{split}
		\Ee
		Inductively we derive stability
		\Be \label{stab}\begin{split}
		&\sup_{0 \leq s \leq t}\| f^{\ell } (s)- f^m (s) \|_{L^{1+ \delta}(\O \times \R^3)}^{1+ \delta}
		+ \int^t_0 |f^{\ell } (s)- f^m (s) |_{L^{1+\delta}  (\gamma)}\\
		 \leq & \  O(t^{1+ \delta})^{\min \{\ell, m\}},
	\end{split}	\Ee
		and hence (\ref{Cauchy_L1+}). 
		
		\vspace{4pt}

		\textit{Step 6. } We combine (\ref{uniform_h_ell}) and (\ref{Cauchy_L1+}) to get unique weak-$\ast$ convergence (up to subsequence if necessary),
		 $(w_{\vartheta} f^\ell, w_{\vartheta} f^{\ell+1}) \overset{\ast}{\rightharpoonup} (w_{\vartheta} f, w_{\vartheta} f)$ weakly$-*$ in $L^\infty (\R \times \O \times \R^3) \cap L^\infty (\R \times \gamma)$. For $\varphi \in C^\infty_c (\R \times \bar{\O} \times \R^3)$, 
		\Be\begin{split}\label{weak_form_ell}
			&\int_0^T\iint_{\O \times \R^3} f^{\ell+1} [-\p_t  - v\cdot \nabla_x  + \nu   ] \varphi\\
			& \ \ \ \ \ \  \ \ \ \ \ 
			+\underbrace{ f^{\ell+1} \{\nabla_x \phi^\ell \cdot \nabla_v\varphi
				+ \frac{v}{2} \cdot \nabla_x \phi^\ell \varphi\}}_{(\ref{weak_form_ell})_\phi}
			\\
			=& \int_0^T \iint_{\O \times \R^3} K f^\ell \varphi - v\cdot \nabla_x \phi^\ell \sqrt{\mu}  \varphi 
			\\& +   \underbrace{\Gamma_{\text{gain}} ( {f^\ell} , {f^\ell} ) \varphi}_{(\ref{weak_form_ell})_\text{gain}} -   \underbrace{\Gamma_{\text{loss}} ( {f^\ell} , {f^{\ell+1}}  ) \varphi }_{(\ref{weak_form_ell})_{\text{loss}}}\\
			&+ \int^T_0 \int_{\gamma_+} f^{\ell+1} \varphi - \int^T_0 \int_{\gamma_-} 
			c_\mu \sqrt{\mu} \int_{n \cdot u>0} f^\ell \sqrt{\mu}\{n \cdot u\} \dd u 
			\varphi.
		\end{split}\Ee
Except the underbraced terms in (\ref{weak_form_ell}) all terms converges to limits with $f$ instead of $f^{\ell+1}$ or $f^\ell$.

		We define, for $(t,x,v) \in \mathbb{R} \times  \bar{\Omega} \times \mathbb{R}^{3}$ and for $0 < \delta \ll 1$, 
		\begin{equation}\label{Z_dyn}
		\begin{split}
		f_{\delta}^\ell(t,x,v) 
		:=  &  \ \kappa_\delta (x,v) f^\ell(t,x,v)
		\\
		: =  &  \  \chi\Big(\frac{|n(x) \cdot v|}{\delta}-1\Big) 
		\Big[ 1- \chi(\delta|v|) \Big]  \chi\Big(\frac{|v|}{\delta}-1\Big)  f^\ell(t,x,v )
		.
		\end{split}
		\end{equation}
		Note that $f_\delta (t,x,v)=0$ if either $|n(x) \cdot v| \leq \delta$, $|v| \geq \frac{1}{\delta}$, or $|v| \leq \delta$.

		%
		
		From (\ref{uniform_h_ell})
		\Be\notag
		\begin{split}
			&\left|\int^T_0 \iint (\ref{weak_form_ell})_\text{loss}
			- \int^T_0 \iint \Gamma_{\text{loss}} (f,f) \varphi\right|
			\\
			\leq &\left| \int^T_0 \iint_{\O \times \R^3}
			\int_{\R^3}   |v-u| \{f^\ell ( u) - f(u) \}\sqrt{\mu(u)}  \dd u f^{\ell+1} (v) \varphi(t,x,v
			) \dd v \dd x \dd t\right|\\
			 +&\left|\int^T_0 \iint_{\O \times \R^3}
			\int_{\R^3}   |v-u|  f(u) \sqrt{\mu(u)}  \dd u \{ f^{\ell+1} (v)- f(v) \} \varphi(t,x,v
			) \dd v \dd x \dd t\right|.
		\end{split}\Ee
		The second term converges to zero from the weak$-*$ convergence in $L^\infty$ and (\ref{uniform_h_ell}). The first term is bounded by, from (\ref{uniform_h_ell}),
		\Be \label{difference_f^ell-f}
		\begin{split}
			&\left[\int^T_0 \left\|  \int_{\R^3} \kappa_\delta (x,u) (f^\ell(t,x,u) - f(t,x,u)) \langle u\rangle \sqrt{\mu(u)} \dd u \right\|^2_{L^2(\O \times \R^3)}\right]^{1/2}\\
			&\times \sup_{0 \leq t \leq T} \|w_{\vartheta} f^{\ell+1}(t)\|_\infty 
			+O(\delta).
		\end{split} \Ee
		
		On the other hand, from Lemma \ref{extension_dyn}, we have an extension $\bar{f}^\ell(t,x,v)$ of $\kappa_\delta (x,u)  f^\ell(t,x,u)$. We apply the average lemma (see Theorem 7.2.1 in page 187 of \cite{gl}, for example) to $\bar{f}^\ell(t,x,v)$. From (\ref{fell_local}) and (\ref{uniform_h_ell})
		
		\Be\label{H_1/4}
		\sup_\ell\left\| \int_{\R^3}   \bar{f}^\ell(t,x,u)  
		\langle u\rangle \sqrt{\mu(u)} \dd u
		\right\|_{H^{1/4}_{t,x} (\R \times \R^3)}< \infty.
		\Ee\hide
		\Be
		\sup_\ell\left\| \int_{\R^3}   f^\ell(t,x,u)  
		\langle u\rangle \sqrt{\mu(u)} \dd u
		\right\|_{H^{1/4}_{t,x} (\R \times \R^3)}< \infty.
		\Ee\unhide
		Then by $H^{1/4} \subset\subset L^2$, up to subsequence, we conclude that 
		\Be\begin{split}
		&\int_{\R^3} \kappa_\delta (x,u)  f^\ell(t,x,u)  \langle u\rangle \sqrt{\mu(u)} \dd u \\
		\rightarrow  & \ \int_{\R^3} \kappa_\delta (x,u)  f(t,x,u) \langle u\rangle \sqrt{\mu(u)} \dd u \  \  {strongly \ in } \  L^2_{t,x}.
		\end{split}\Ee
		So we conclude that $(\ref{difference_f^ell-f})\rightarrow 0$ as $\ell \rightarrow \infty$.
		
		For $(\ref{weak_form_ell})_{\text{gain}}$ let us use a test function $\varphi_1(v) \varphi_2 (t,x)$. From the density argument, it suffices to prove a limit by testing with $\varphi(t,x,v)$.
		
		We use a standard change of variables $(v,u) \mapsto(v^\prime,u^\prime)$ and $(v,u) \mapsto(u^\prime,v^\prime)$ (for example see page 10 of \cite{gl}) to get    
		\begin{align} 
		 &\int^T_0 \iint  (\ref{weak_form_ell})_{\text{gain}} - \int^T_0 \iint \Gamma_{\text{gain}}(f,f) \varphi \nonumber\\
		  =& \ \int^T_0 \iint\Gamma_{\text{gain}}(f^\ell - f,f^\ell) \varphi + \int^T_0 \iint\Gamma_{\text{gain}}( f,f^\ell-f) \varphi \nonumber \\
		 = &   \int^T_0 \iint_{\O \times \R^3}\nonumber\\
	 & \   	\times  \left( \int_{\R^3}  \int_{\S^2} 
		(   f^\ell (t,x,u) -  f(t,x,u)  )   \sqrt{\mu(u^\prime)}     |(v-u) \cdot \o|  \varphi_1(v^\prime)            \dd \o\dd u\right)\nonumber\\
		 & \   \times 
		f^\ell (t,x,v)   
		\varphi_2 (t,x)
		\dd v \dd x\dd t  \label{weak_fell-f1}\\
		   +&  \int^T_0 \iint_{\O \times \R^3} \bigg( \int_{\R^3}  \int_{\S^2} \nonumber \\
	 & \    	\times	(   f^\ell (t,x,u) -  f(t,x,u)  )   \sqrt{\mu(v^\prime)}     |(v-u) \cdot \o|  \varphi_1(u^\prime)            \dd \o\dd u\bigg)\nonumber\\
		 & \  \times 
		f (t,x,v)   
		\varphi_2 (t,x)
		\dd v \dd x\dd t  . \label{weak_fell-f2}
		\end{align} 
		For $N\gg 1$ we decompose the integration of (\ref{weak_fell-f1}) and (\ref{weak_fell-f2}) using 
		\Be
		\begin{split}\label{decomposition_N}
			1=& \{1- \chi (|u|-N)\}\{1- \chi (|v|-N)\}\\
			& +  \chi (|u|-N) + \chi (|v|-N) - \chi (|u|-N)   \chi (|v|-N).
		\end{split} \Ee
		Note that $\{1- \chi (|u|-N)\}\{1- \chi (|v|-N)\} \neq 0$ if $|v| \leq N+1$ and $|u| \leq N+1$, and if $ \chi (|u|-N) + \chi (|v|-N) - \chi (|u|-N)   \chi (|v|-N) \neq 0$ then either $|v|\geq N$ or $|u| \geq N$. From (\ref{uniform_h_ell}), the second part of (\ref{weak_fell-f1}) and (\ref{weak_fell-f2}) from (\ref{decomposition_N}) are bounded by 
		\Be
		\begin{split}
			&\int^T_0 \iint_{\O \times \R^3} \int_{\R^3} \int_{\S^2} [\cdots]  \times \{\chi (|u|-N) + \chi (|v|-N) - \chi (|u|-N)   \chi (|v|-N)\}\\
			\leq &  \ \sup_\ell \| w_{\vartheta}f^\ell \|_\infty \| w_{\vartheta} f\|_\infty \times 
			\{
			e^{-\frac{\vartheta}{2} |v|^2} e^{-\frac{\vartheta}{2} |u|^2} 
			\}\{ \mathbf{1}_{|v|\geq N} +  \mathbf{1}_{|u|\geq N}  \}\\
			\leq& \  O(\frac{1}{N}). \notag
		\end{split}
		\Ee
		
		Now we only need to consider the parts with $ \{1- \chi (|u|-N)\}\{1- \chi (|v|-N)\}$. Then 
		\Be\begin{split}\label{bound_weak_fell-f1}
			&(\ref{weak_fell-f1})\\
			= & \int^T_0 \iint_{\O \times \R^3}    \int_{\R^3}  
			(   f^\ell (t,x,u) -  f(t,x,u)  ) \\
			& \ \ \ \ \ \ \ \ \ \ \ \  \times  \{1- \chi (|u|-N)\}\left(\int_{\S^2}  \sqrt{\mu(u^\prime)}     |(v-u) \cdot \o|  \varphi_1(v^\prime)            \dd \o \right)\dd u \\
			& \ \ \ \ \ \ \ \ \ \ \ \  \times 
			\{1- \chi (|v|-N)\}  f^\ell (t,x,v)   
			\varphi_2 (t,x)
			\dd v \dd x\dd t.
		\end{split}
		\Ee
		Let us define 
		\Be\label{Phi_v}
		\Phi _v(u ) :=  \{1- \chi (|u|-N)\}\int_{\S^2}  \sqrt{\mu(u^\prime)}     |(v-u) \cdot \o|  \varphi_1(v^\prime)            \dd \o \ \ \text{for} \ |v| \leq N+1.
		\Ee          
		
		For $0<\delta\ll1$ we have $O(\frac{N^3}{\delta^3})$ number of $v_i \in \R^3$ such that $\{v \in \R^3: |v| \leq N+1\}\subset\bigcup_{i=1}^{O(\frac{N^3}{\delta^3})} B(v_i, \delta)$. Since (\ref{Phi_v}) is smooth in $u$ and $v$ and compactly supported, for $0<\e\ll1$ we can always choose $\delta>0$ such that
		\Be\label{Phi_v_continuous}
		|\Phi_v(u) - \Phi_{v_i} (u)|< \e  \ \ \text{if} \ v \in B(v_i, \delta).
		\Ee
		
		Now we replace $\Phi_v(u)$ in the second line of (\ref{bound_weak_fell-f1}) by $\Phi_{v_i}(u)$ whenever $v \in B(v_i, \delta)$. Moreover we use $\kappa_\delta$-cut off in (\ref{Z_dyn}). If $v$ is included in several balls then we choose the smallest $i$.  From (\ref{Phi_v_continuous}) and (\ref{uniform_h_ell}) the difference of (\ref{bound_weak_fell-f1}) and the one with $\Phi_{v_i}(u)$ can be controlled and we conclude that
		\Be
		\begin{split}\label{bound_weak_fell-f1_ell}
			(\ref{bound_weak_fell-f1}) & =\{O( \e)+ O(\delta)\} \sup_{\ell} \| w_{\vartheta} f^\ell \|_\infty  ^2\\
			&+  \int^T_0 \int_{\O  } \sum_{i} \int_{\R^3} \mathbf{1}_{v \in B(v_i,\delta) } \\
			& \ \ \ \ \ \ \ \ \ \  \times  \int_{\R^3}  
			\kappa_\delta(x,u)  (   f^\ell (t,x,u) -  f(t,x,u)  ) \Phi_{v_i} (u)  \dd u \\
			&  \ \ \ \ \ \ \ \ \ \  \times 
			\{1- \chi (|v|-N)\}  f^\ell (t,x,v)   
			\varphi_2 (t,x)
			\dd v \dd x\dd t.
		\end{split}\Ee
		From Lemma \ref{extension_dyn} and the average lemma
		\Be\label{H_1/4_gain}
		\max_{ 1\leq i \leq O(\frac{N^3}{\delta^3})}\sup_\ell\left\| \int_{\R^3} \kappa_\delta (x,u) f^\ell(t,x,u)  
		\Phi_{v_i} (u) \dd u
		\right\|_{H^{1/4}_{t,x} (\R \times \R^3)}< \infty.
		\Ee
		For $i=1$ we extract a subsequence $\ell_1 \subset\mathcal{I}_1$ such that 
		\Be\label{strong_converge_extract}
		\begin{split}
	&	\int_{\R^3} \kappa_\delta (x,u) f^{\ell_1}(t,x,u)  
		\Phi_{v_i} (u) \dd u \\
		\rightarrow & \int_{\R^3} \kappa_\delta (x,u) f (t,x,u)  
		\Phi_{v_i} (u) \dd u \   \ { strongly \  in }  \ L^2_{t,x}.
		\end{split}\Ee
		Successively we extract subsequences $\mathcal{I}_{O(\frac{N^3}{\delta^3})} \subset \cdots\subset \mathcal{I}_2 \subset \mathcal{I}_1$. Now we use the last subsequence $\ell \in\mathcal{I}_{O(\frac{N^3}{\delta^3})}$ and redefine $f^\ell$ with it. Clearly we have (\ref{strong_converge_extract}) for all $i$. Finally we bound the last term of (\ref{bound_weak_fell-f1_ell}) by
		\Be\notag
		\begin{split}
			&C_{\varphi_2,N} \max_i \int^T_0\left\|
			\int_{\R^3} \kappa_\delta (x,u) (f^\ell(t,x,u)- f(t,x,u)  )
			\Phi_{v_i} (u) \dd u
			\right\|_{L^2_{t,x}}\\
			& \ \times  \sup_\ell \| w_{\vartheta} f^\ell \|_\infty\\
			&\rightarrow 0 \ \ \text{as} \ \ \ell \rightarrow \infty.
		\end{split}
		\Ee
		Together with (\ref{bound_weak_fell-f1_ell}) we prove $(\ref{weak_fell-f1})\rightarrow 0$. Similarly we can prove $(\ref{weak_fell-f2})      \rightarrow 0$.  
		
		Now we consider $(\ref{weak_form_ell})_{\phi}$. From 
		\Be\begin{split}\notag
			- ( \Delta \phi^\ell-  \Delta \phi) 
			= \int  \kappa_\delta (f^\ell- f)\sqrt{\mu} +   \int (1-\kappa_\delta) (f^\ell -f)\sqrt{\mu} ,
		\end{split} \Ee
		we have
		\Be\begin{split}
			\| \nabla_x \phi^\ell - \nabla_x \phi \|_{L^2_{t,x}} 
			\leq  \left\| \int  \kappa_\delta (f^\ell- f)\sqrt{\mu}\right\|_{L^2_{t,x}}
			+ O(\delta) \sup_\ell \| w_{\vartheta} f ^\ell \|_\infty.
		\end{split}
		\Ee
		Then following the previous argument, we prove $  \nabla_x \phi^\ell \rightarrow  \nabla_x \phi$ strongly in $L^2_{t,x}$ as $\ell \rightarrow \infty$. Combining with $w_{\vartheta} f^\ell \overset{\ast}{\rightharpoonup} w_{\vartheta} f$ in $L^\infty$, we prove $\int^T_0 \iint_{\O \times \R^3}(\ref{weak_form_ell})_{\phi}$ converges to $\int^T_0 \iint_{\O \times \R^3}f  \{\nabla_x \phi \cdot \nabla_v\varphi
		+ \frac{v}{2} \cdot \nabla_x \phi \varphi\}$. This proves the existence of a (weak) solution $f \in L^\infty$. 
		
		\vspace{4pt}

		\textit{Step 7. }  We claim (\ref{W1p_local_bound}).
	\hide	\Be\begin{split}  \notag
			&\sup_{0 \leq t \leq T^{**}}  \|w_{\tilde{\vartheta}} f  (t) \|_p^p
			+
			\sup_{0 \leq t \leq T^{**}}  \| w_{\tilde{\vartheta}}\alpha_{f ,\e}^\beta \nabla_{x,v} f (t) \|_{p}^p
			+ 
			\int^{T^{**}}_0  | w_{\tilde{\vartheta}}\alpha_{f ,\e}^\beta \nabla_{x,v} f (t) |_{p,+}^p \\
			&\lesssim  \   \| w_{\tilde{\vartheta}}f _0 \|_p^p
			+
			\| w_{\tilde{\vartheta}}\alpha_{f_{0 },\e}^\beta \nabla_{x,v} f  _0 \|_{p}^p
			.
		\end{split}\Ee
		This implies .\unhide
		By the weak lower-semicontinuity of $L^p$ we know that (if necessary we further extract a subsequence out of  the subsequence of \textit{Step 6})
		\[
			w_{\tilde{\vartheta}}\alpha^\beta_{f^\ell, \e} \nabla_{x,v} f^{\ell+1} \rightharpoonup \mathcal{F},  \]
			and 
			\[
			\sup_{0 \leq t \leq T^{**}}\|\mathcal{F}(t)\|_p^p \leq \liminf   \sup_{0 \leq t \leq T^{**}}  \| w_{\tilde{\vartheta}}\alpha_{f^{\ell} ,\e}^\beta \nabla_{x,v} f^{\ell+1} (t) \|_{p}^p,
			\]
		and 
		\[ 
		\int^{T^{**}}_0  | \mathcal{F}  |_{p,+}^p \leq \liminf   
		\int^{T^{**}}_0  | w_{\tilde{\vartheta}}\alpha_{f^{\ell} ,\e}^\beta \nabla_{x,v} f^{\ell+1} (t) |_{p,+}^p.
		\]
		We need to prove that 
		\Be\label{a_nabla_f_seq}
			\mathcal{F}
			=w_{\tilde{\vartheta}}\alpha^{\beta}_{f,\e} \nabla_{x,v}f  \ \  {almost \  everywhere  \ except } \  \gamma_0.
		\Ee

		We claim that, up to some subsequence, for any given smooth test function $\psi \in C_c^\infty(\bar{\O} \times \R^3\backslash \gamma_0)$ 
\Be
\begin{split}
\label{lim_Fell=F}
&\lim_{\ell  \rightarrow \infty}\int^t_0\iint_{\O \times \R^3} w_{\tilde{\vartheta}} \alpha^\beta_{f^{\ell },\e} \nabla_{x,v} f^{\ell +1} \psi \dd x\dd v\\
= & \  \int^t_0\iint_{\O \times \R^3} w_{\tilde{\vartheta}} \alpha^\beta_{f ,\e} \nabla_{x,v} f  \psi \dd x\dd v.
\end{split}\Ee 
We note that we need to extract a single subsequence, let say $\{\ell_*\} \subset \{\ell\}$, satisfying (\ref{lim_Fell=F}) for all test functions in $C_c^\infty(\bar{\O} \times \R^3\backslash \gamma_0)$. Of course the convergent rate needs not to be uniform and it could vary with test functions. 

For each $N\in\mathbb{N}$ we define a set
\Be\label{set_N}\begin{split}
\mathcal{S}_N:=& \Big\{ (x,v) \in \bar{\O} \times \R^3: \text{dist}(x, \p\O) \leq \frac{1}{N}  \ \text{and} \ |n(x) \cdot v| \leq \frac{1}{N}  \Big\}
	\\	&  \cup \{|v|>N\}.
\end{split}\Ee
For a given test function we can always find $N\gg1$ such that 
		\Be\label{test_function_condition}
		supp(\psi) \subset (\mathcal{S}_N)^c:=\bar{\O} \times \R^3 \backslash \mathcal{S}_N.
		 \Ee 
		 
We will exam (\ref{lim_Fell=F}) by the identity obtained from the integration by parts
		\begin{eqnarray}
			&&\int^t_0\iint_{\O \times \R^3} w_{\tilde{\vartheta}} \alpha^\beta_{f^{\ell },\e} \nabla_{x,v} f^{{\ell }+1} \psi \dd x\dd v
			\notag
			\\
			&= &  
			-  \int^t_0\iint_{\O \times \R^3} \alpha^\beta_{f^{\ell },\e}  f^{{\ell }+1} \nabla_{x,v}(w_{\tilde{\vartheta}}\psi) \dd x\dd v\label{af_l_1}
			  \\
			&& + \int^t_0 \iint_{\gamma} n  \alpha^\beta_{f^{\ell },\e}   f^{{\ell }+1} (w_{\tilde{\vartheta}}\psi)
			\label{af_l_2}\\
			&&-  {\int^t_0\iint_{\O \times \R^3} \nabla_{x,v}\alpha^\beta_{f^{\ell },\e}  f^{{\ell }+1} (w_{\tilde{\vartheta}}\psi) \dd x\dd v}. \label{af_l_3}
		\end{eqnarray}
	We finish this step by proving the convergence of (\ref{af_l_1}) and (\ref{af_l_2}). From (\ref{alphaweight}) and (\ref{uniform_h_ell}), if $(x,v) \in (\mathcal{S}_N)^c$ then
		\Be
		\begin{split}\notag
	&	\sup_{\ell \geq 0}	| \alpha^\beta_{f^\ell,\e}(t,x,v)|  \\
	\lesssim  &  |v|^\beta + (t+\e)^\beta\sup_{\ell \geq 0} \| \nabla \phi_{f^\ell} \|_\infty^\beta \\
	\lesssim  & \ N^\beta  +(T^{**}+\e)^\beta \sup_{\ell \geq 0} \| w_{\vartheta} f^\ell \|_\infty^\beta
	 \\
\leq   &
	C_N <+\infty  .
		\end{split}
		\Ee
	Hence we extract a subsequence (let say $\{\ell_N\}$) out of subsequence in \textit{Step 6} such that $
	 \alpha_{f^{\ell_N},\e}^\beta  \overset{\ast}{\rightharpoonup} A \in L^\infty  \textit{ weakly}-* \textit{ in }  L^\infty ((0,T^{**}) \times (\mathcal{S}_N)^c) \cap L^\infty ((0,T^{**})  \times (\gamma \cap  (\mathcal{S}_N)^c)).$ Note that $\alpha_{f^{\ell_N},\e}^\beta$ satisfies $[\p_t + v\cdot \nabla_x - \nabla_x \phi^{\ell_N}\cdot \nabla_v]  \alpha_{f^{\ell_N},\e}^\beta=0$ and $\alpha_{f^{\ell_N},\e}^\beta|_{\gamma_-} =|n \cdot v |^\beta$. By passing a limit in the weak formulation we conclude that $[\p_t + v\cdot \nabla_x - \nabla_x \phi_f \cdot \nabla_v]  A =0$ and $A |_{\gamma_-} =|n \cdot v |^\beta$. By the uniqueness of the Vlasov equation ($\nabla \phi_f \in W^{1,p}$ for any $p<\infty$) we derive $A=  \alpha_{f ,\e}^\beta$ almost everywhere and hence conclude that 
	\Be\label{converge_alpha_ell}\begin{split}
	& \alpha_{f^{\ell_N},\e}^\beta  \overset{\ast}{\rightharpoonup}   \alpha_{f ,\e}^\beta  \textit{ weakly}-*\\
	& \textit{ in }  L^\infty ((0,T^{**})  \times (\mathcal{S}_N)^c) \cap L^\infty ((0,T^{**})  \times (\gamma \cap  (\mathcal{S}_N)^c)).
	\end{split}\Ee	
Now the convergence of (\ref{af_l_1}) and (\ref{af_l_2}) is a direct consequence of strong convergence of (\ref{Cauchy_L1+}) and the weak$-*$ convergence of (\ref{converge_alpha_ell}).

		\vspace{4pt}

		\textit{Step 8. } We devote the entire \textit{Step 8} to prove the convergence of (\ref{af_l_3}).

		\textit{Step 8-a. } Let us choose $(x,v) \in (\mathcal{S}_N)^c$. From (\ref{weight}) 
	\Be\label{alpha=1}
If \  \ \tb^{f^\ell} \geq t+ \e \  \ then \  \ 	\alpha_{f^\ell, \e}(t,x,v)=1.
	\Ee
	 From now we only consider that case 
	\Be\label{tb_upperT^**}
	\tb^{f^\ell} (t,x,v) \leq \e +  t.
	\Ee
	
	If $|v|\geq 2 (\e + T^{**})\sup_{\ell} \| \nabla \phi^\ell \|_\infty$ then 
	\Be\notag
	\begin{split}
	|V^{f^\ell}(s;t,x,v)| &\geq |v| - \int^t_s \| \nabla \phi^\ell(\tau) \|_\infty \dd \tau \\
	&\geq  (\e + T^{**})\sup_{\ell} \| \nabla \phi^\ell \|_\infty \ \ \ for \ all \  \ell  \ and \ s \in [-\e ,T^{**}].
	\end{split}
	\Ee
	Then following the argument to get (\ref{velocity_N}), which are based on the estimates (\ref{velocity_lemma_N}) and (\ref{phi_xi}), we derive 
	\Be\begin{split}\notag
	& |\tilde{\alpha}(s,X^{f^\ell}(s;t,x,v),V^{f^\ell}(s;t,x,v))| \\
	\geq & \   \frac{1}{C_\O}\tilde{\alpha} (t,x,v) e^{-C_\O \frac{|t-s|}{(\e + T^{**})\sup_{\ell} \| \nabla \phi^\ell  \|_\infty}}\\
	\geq & \ 
	\frac{e^{- \frac{C_\O }{\sup_{\ell} \| \nabla \phi^\ell  \|_\infty}}}{C_\O } \times \frac{1}{N}
	   \ \ \ for \ all \  \ell  \ and \ s \in [-\e ,T^{**}].
	\end{split}\Ee
Especially at $s= t-\tb^{f^\ell}(t,x,v)$, from (\ref{beta}),
\Be\label{alpha_|v|>}
|n(\xb^{f^\ell}) \cdot \vb^{f^\ell}|
 \geq \frac{e^{- \frac{C_\O }{\sup_{\ell} \| \nabla \phi^\ell  \|_\infty}}}{C_\O } \times \frac{1}{N} \ \ \ for \ all \ \ell. 
\Ee

	 	\vspace{4pt}
	 
	 \textit{Step 8-b. }  From now on we assume (\ref{tb_upperT^**}). From$|v|\leq 2 (\e + T^{**})\sup_{\ell} \| \nabla \phi^\ell \|_\infty$ and (\ref{hamilton_ODE}),
	\Be\label{upper_|V|}
	\begin{split}
	 |V^{f^\ell}(s;t,x,v)|\leq 3 (\e + T^{**})\sup_{\ell} \| \nabla \phi^\ell  \|_\infty \ \    for \  \ s \in [-\e, T^{**}].
	\end{split}\Ee
	 Let $(X_n^{f^\ell}, X_\parallel^{f^\ell},V_n^{f^\ell}, V_\parallel^{f^\ell})$ satisfy (\ref{dot_Xn_Vn}), (\ref{def_V_parallel}), and (\ref{hamilton_ODE_perp}) with $E=- \nabla \phi^\ell$. 
	 
	 	 Let us define 
		\Be \label{tau_1}
		\tau_1  
		: = \sup  \big\{ \tau \geq  0: V_n^{f^\ell} (s;t,x,v)\geq 0 \ for \ all \ s \in [t-\tb^{f^\ell}(t,x,v),  \tau ]
		\big\}
		. 
		\Ee
		Since $(X^{f^\ell} (s;t,x,v),V^{f^\ell} (s;t,x,v))$ is $C^1$ (note that $\nabla\phi^\ell \in C^1_{t,x}$) in $s$ we have $V_n^{f^\ell} (\tau_1 ;t,x,v)=0$.

	 We claim that, there exists some constant $\delta_{**} = O_{\e, T^{**}, \sup_{\ell} \| \nabla \phi^\ell  \|_{C^1}}(\frac{1}{N})$ in (\ref{choice_delta**}) which does not depend on $\ell$ such that 
	 \Be\label{claim_Vperp_growth}
	 \begin{split}
	 If& \ \ 0\leq V_n^{f^\ell} (t-\tb^{f^\ell}(t,x,v);t,x,v) < \delta_{**} \  and \  (\ref{upper_|V|}),
	 \\
	&then \ \ V_n^{f^\ell} (s;t,x,v) \leq e^{C|s-(t-\tb^{f^\ell}(t,x,v))|^2}V_n^{f^\ell} (t-\tb^{f^\ell}(t,x,v);t,x,v)\\ &  for  \ s \in [ t-\tb^{f^\ell},\tau_1].
	 \end{split}\Ee
For the proof we regard the equations (\ref{dot_Xn_Vn}), (\ref{def_V_parallel}), and (\ref{hamilton_ODE_perp}) as the forward-in-time problem with an initial datum at $s=t-\tb^{f^\ell} (t,x,v)$. 	Clearly we have $X_n^{f^\ell}(t-\tb^{f^\ell}(t,x,v);t,x,v)=0$ and $V_n^{f^\ell}(t-\tb^{f^\ell}(t,x,v);t,x,v)\geq0$ from Lemma \ref{cannot_graze}. Again from Lemma \ref{cannot_graze}, if $V_n^{f^\ell}(t-\tb^{f^\ell}(t,x,v);t,x,v)=0$ then $X_n^{f^\ell}(s;t,x,v)=0$ for all $s\geq t-\tb^{f^\ell} (t,x,v)$. From now on we assume $V_n^{f^\ell} (t-\tb^{f^\ell}(t,x,v);t,x,v)] >0$. From (\ref{hamilton_ODE_perp}), as long as $t - \tb^{f^\ell} (t,x,v) \leq s \leq T^{**}$ and 
	\Be\label{small_s_V_n} 
	V_n^{f^\ell}(s;t,x,v) \geq 0 \ \ and \ \ 
	X_\perp^{f^\ell} (s;t,x,v) \leq \frac{1}{N} \ll 1,
	\Ee  
	then we have
	\Be\label{hamilton_ODE_perp_bound}\begin{split}
	\dot{V}_n^{f^\ell} (s) 
		= & \  \underbrace{[V^{f^\ell}_\parallel (s)\cdot \nabla^2 \eta (X^{f^\ell}_\parallel(s)) \cdot V^{f^\ell}_\parallel(s) ] \cdot n(X^{f^\ell}_\parallel(s)) }_{\leq 0  \    from \ (\ref{convexity_eta})} \\
	- 	& \underbrace{\nabla \phi^\ell  (s , X^{f^\ell} (s ) ) \cdot [-n(X^{f^\ell}_\parallel(s)) ]}_{
		=O(1) \sup_\ell  \|  \nabla \phi^\ell \|_{C^1} \times X_n^{f^\ell} (s)
		 \   from \  (\ref{expansion_E})} \\
		  - &\underbrace{X_n^{f^\ell} (s) [V^{f^\ell}_\parallel(s) \cdot \nabla^2 n (X^{f^\ell}_\parallel(s)) \cdot V^{f^\ell}_\parallel(s)]  \cdot n(X^{f^\ell}_\parallel(s)) }_{
		=O(1)  \{3 (\e + T^{**})\sup_{\ell} \| \nabla \phi^\ell  \|_\infty \}^2 \times X_n^{f^\ell} (s) \ from \ 
		(\ref{upper_|V|})
		}\\
		 \leq & \ 
		C (1+ \e  + T^{**})^2 ( \sup_{\ell} \| \nabla \phi^\ell  \|_{C^1}  \sup_{\ell} \| \nabla \phi^\ell  \|_{\infty} )
		  \times X_n^{f^\ell} (s).
	\end{split}\Ee	
	
Let us consider (\ref{hamilton_ODE_perp_bound}) together with $\dot{X}^{f^\ell}_{n}(s;t,x,v) = V^{f^\ell}_{n}(s;t,x,v)$. Then, as long as $s$ satisfies (\ref{small_s_V_n}), 
\Be
\begin{split}\notag
V_n^{f^\ell} (s)& =    V_n^{f^\ell} (t-\tb^{f^\ell} ) + \int^s_{t-\tb^{f^\ell}} \dot{V}_n^{f^\ell}(\tau) \dd \tau \\
&\leq V_n^{f^\ell} (t-\tb^{f^\ell} )\\
& \ \ +\int^s_{t-\tb^{f^\ell}} 
C (1+ \e  + T^{**})^2  ( \sup_{\ell} \| \nabla \phi^\ell  \|_{C^1}  \sup_{\ell} \| \nabla \phi^\ell  \|_{\infty} )
		  \times X_n^{f^\ell} (\tau)
\dd \tau\\
&= V_n^{f^\ell} (t-\tb^{f^\ell} ) \\
& \ \ + C (1+ \e  + T^{**})^2  ( \sup_{\ell} \| \nabla \phi^\ell  \|_{C^1}  \sup_{\ell} \| \nabla \phi^\ell  \|_{\infty} ) \\
& \ \ \times \int^s_{t-\tb^{f^\ell}} 
	\int^\tau_{t-\tb^{f^\ell}} V_n^{f^\ell}(\tau^\prime) \dd \tau^\prime
\dd \tau.
\end{split}
\Ee
Following the same argument of the proof of Lemma \ref{est_X_v}, we derive that 
	\Be
\begin{split}\notag
V_n^{f^\ell} (s) \leq&  \  V_n^{f^\ell} (t-\tb^{f^\ell} )\\
& +
C (1+ \e  + T^{**})^2  ( \sup_{\ell} \| \nabla \phi^\ell  \|_{C^1}  \sup_{\ell} \| \nabla \phi^\ell  \|_{\infty} )\\
& \ \ \times 
\int^s_{t-\tb^{f^\ell}} 
|s -  (t-\tb^{f^\ell}) |  V_n^{f^\ell}(\tau^\prime)  \dd \tau^\prime.
\end{split}
\Ee
From the Gronwall's inequality, we derive that, as long as (\ref{small_s_V_n}) holds, 	
\Be\label{upper_bound_Vn}\begin{split}
&V_n^{f^\ell} (s;t,x,v)\\
 \leq&  \  V_n^{f^\ell} (t-\tb^{f^\ell}(t,x,v)) \\
& \times e^{C (1+ \e  + T^{**})^2 ( \sup_{\ell} \| \nabla \phi^\ell  \|_{C^1}  \sup_{\ell} \| \nabla \phi^\ell  \|_{\infty} ) \times |s -  (t-\tb^{f^\ell}(t,x,v)) | ^2}.
\end{split}
\Ee	

Now we verify the conditions of (\ref{small_s_V_n}) for all $- \e \leq t - \tb^{f^\ell} (t,x,v) \leq s \leq T^{**}$. Note that we are only interested in the case of $V_n^{f^\ell} (t-\tb^{f^\ell}(t,x,v);t,x,v)< \delta_{**}$. From the argument of (\ref{hamilton_ODE_perp_bound}), ignoring negative curvature term,
		 \Be
		 \begin{split}\notag
		&| X_n^{f^\ell} (s;t,x,v)| \\
		\leq& 
	\ 	(\e+ T^{**}) | V_n^{f^\ell} (\tb^{f^\ell};t,x,v)  | \\
	&
		+   C[1 + (\e + T^{**})^2 \sup_\ell \| \nabla \phi^\ell \|_\infty]  \sup_\ell\| \nabla \phi^\ell \|_{C^1} \\
		& \ \ \ \  \times 
		 \int^s_{t-\tb^{f^\ell}} \int^\tau_{t-\tb^{f^\ell}} |X_n^{f^\ell}(\tau;t,x,v)| \dd \tau  \dd s\\
		  \leq & \  (\e+ T^{**}) | V_n^{f^\ell} (\tb^{f^\ell};t,x,v)  | \\
		  &  +  C  \int_{t-\tb^{f^\ell}} ^s |\tau- (t-\tb^{f^\ell})| |X_n^{f^\ell}(\tau;t,x,v)|  \dd \tau.
		 \end{split}
		 \Ee
		 Then by the Gronwall's inequality we derive that, in case of (\ref{tb_upperT^**}),
		 \Be\label{upper_X_n}\begin{split}
		& | X_n^{f^\ell} (s;t,x,v)| \\
		\leq & \  C_{\e+ T^{**}} | V_n^{f^\ell} (t-\tb^{f^\ell};t,x,v)  | \ \ for  \ all \  -\e \leq t-\tb^{f^\ell} \leq s \leq t  \leq T^{**}.
		\end{split} \Ee 
If we choose 
\Be\label{choice_delta**}
\delta_{**} = \frac{o(1)}{ |T^{**} + \e |}\times \frac{1}{N}, 
\Ee
then (\ref{upper_bound_Vn}) holds for $- \e \leq t - \tb^{f^\ell} (t,x,v) \leq s \leq T^{**}$. Hence we complete the proof of (\ref{claim_Vperp_growth}).

	\hide

	\Be\notag
	\dot{X}_{n}(s;t,x,v) = V_{n}(s;t,x,v) , \ \  \dot{X}_{\parallel}(s;t,x,v) = V_{\parallel}(s;t,x,v).
	\Ee

	Then the proof of (\ref{no_graze}) asserts that  
		\Be\label{alpha_lower_bound}
		|V_n^f(s;t,x,v)| \gtrsim 1 \ \ \text{for all } \ (t,x,v) \in [0,\frac{\e}{2}] \times  \text{supp} (\psi) .
		\Ee
		Now we consider $\alpha_{f^\ell, \e}(t,x,v)$.

		 On the other hand if $\tb^{f^\ell}(t,x,v)< 2\e$, since (\ref{uniform_h_ell}), from (\ref{hamilton_ODE_perp}),
		\Be
		\begin{split}\notag
			&\sup_{\ell \geq 0}|V^f_n(s;t,x,v) - V^{f^\ell }_n(s;t,x,v)|\\
			\lesssim & \ \sup_{\ell \geq 0}\max\{ \tb^{f}, \tb^{f^\ell}\} \times \left\{N^2  + \| \nabla \phi_f \|_\infty + \| \nabla \phi_{f^\ell} \|_\infty\right\}\\
			\lesssim & \ \e (1+ N^2).
		\end{split}\Ee
	Therefore, from (\ref{alpha_lower_bound})	 for small $\e>0$, we prove the lower bound 
		\Be
	\inf_{\ell }	|n(\xb^{f^\ell})  \cdot \vb^{f^\ell}| \gtrsim 1 \ \ \text{for all } \ (t,x,v) \in [0,\frac{\e}{2}] \times  \text{supp} (\psi) .
		\Ee
		
		\unhide

		\vspace{4pt}

		\textit{Step 8-c. } Suppose that (\ref{upper_|V|}) holds and $0 \leq V_n^{f^\ell} (t-\tb^{f^\ell}(t,x,v);t,x,v) < \delta_{**}$ with $\delta_{**}$ of  (\ref{choice_delta**}). Recall the definition of $\tau_1$ in (\ref{tau_1}). Inductively we define $\tau_2  
		: = \sup  \big\{ \tau \geq  0: V_n^{f^\ell} (s;t,x,v)\leq 0 \ for \ all \ s \in [\tau_1,  \tau ]
		\big\}$ and $\tau_3, \tau_4, \cdots$. Clearly such points can be countably many at most in an interval of $[t-\tb^{f^\ell},t]$. Suppose $ \lim_{k\rightarrow \infty}\tau_k=t   $. Then choose $k_0\gg1$ such that $|\tau_{k_0} -t| \ll  |V_n^{f^\ell} (t-\tb^{f^\ell} ;t,x,v)|$. Then, for $s \in [\tau_{k_0}, t] $, from (\ref{hamilton_ODE_perp_bound}) and (\ref{upper_|V|}),
		\Be\begin{split}\label{upper_V_n_0}
		|V_n^{f^\ell} (t;t,x,v)| \lesssim   {|V_n^{f^\ell} (t-\tb^{f^\ell} ;t,x,v)|} .
		\end{split}\Ee 
		
		Now we assume that $\tau_{k_0} < t \leq \tau_{k_0+1}$. From the definition of $\tau_i$ in  (\ref{tau_1}) we split the case in two.
		
		\textit{\underline{Case 1:} Suppose $V_n^{f^\ell} (s;t,x,v)>0$ for $s \in (\tau_{k_0}, t)$. }

		From (\ref{hamilton_ODE_perp_bound}) and (\ref{upper_X_n})
		\Be
		\begin{split}\label{upper_V_n_1}
		 V_n^{f^\ell} (t;t,x,v)  \lesssim \int^{T^{**}}_{\tau_{k_0}} X_n^{f^\ell}(s) \lesssim |V_n^{f^\ell} (t-\tb^{f^\ell};t,x,v)|.
		\end{split}
		\Ee

		\textit{\underline{Case 2:} Suppose $V_n^{f^\ell} (s;t,x,v)<0$ for $s \in (\tau_{k_0}, t)$. }
		
		Suppose 
		\Be\label{assumption_lowerbound_V_n}\begin{split}
		-V_n^{f^\ell} (t;t,x,v)&= |V_n^{f^\ell} (t;t,x,v)| \\
		&\geq |V_n^{f^\ell} (t-\tb^{f^\ell};t,x,v)|^A \ \ for \ any \  0<A<\frac{1}{2}. 
		\end{split}\Ee \hide	
In this step we claim that 
\Be\label{max_Xn}
X_n^{f^\ell} (
t-\tb^{f^\ell}(t,x,v)+
\tau_*^\ell (t,x,v);t,x,v)\geq \frac{1}{2} \frac{| V_n^{f^\ell} (t- \tb^{f^\ell} (t,x,v);t,x,v)|^2}{1+\{3 (\e + T^{**})\sup_{\ell} \| \nabla \phi^\ell  \|_\infty\}^2 }.
\Ee		
\unhide
From (\ref{hamilton_ODE_perp_bound}), now taking account of the curvature term this time, we derive that 	\Be
	\begin{split}\notag
	-V_n^{f^\ell} (t;t,x,v)  \leq & \   \int_{\tau_{k_0}}^{t} (-1) [V^{f^\ell}_\parallel (s)\cdot \nabla^2 \eta (X^{f^\ell}_\parallel(s)) \cdot V^{f^\ell}_\parallel(s) ] \cdot n(X^{f^\ell}_\parallel(s))  \dd s\\
	& +C | V_n^{f^\ell} (t- \tb^{f^\ell} (t,x,v);t,x,v)|,
	\end{split}
	\Ee		
	where we have used (\ref{upper_|V|}) and (\ref{upper_X_n}). From (\ref{assumption_lowerbound_V_n}) the above inequality implies that, for $| V_n^{f^\ell} (t- \tb^{f^\ell} (t,x,v);t,x,v)|\ll 1$, 
	\Be\notag
	|V_n^{f^\ell} (t-\tb^{f^\ell};t,x,v)|^A \lesssim \int_{\tau_{k_0}}^{t} (-1) [V^{f^\ell}_\parallel (s)\cdot \nabla^2 \eta (X^{f^\ell}_\parallel(s)) \cdot V^{f^\ell}_\parallel(s) ] \cdot n(X^{f^\ell}_\parallel(s))  \dd s. 
	\Ee
Note that $|\frac{d}{ds} V^{f^\ell}_\parallel (s)|$ and $|\frac{d}{ds}X^{f^\ell}_\parallel (s)|$ are all bound from $\nabla \phi^\ell \in C^1$, (\ref{upper_|V|}), and (\ref{upper_X_n}). Hence we have 
\Be\label{lower_bound_V_||^2}\begin{split}
	&\frac{1}{2}|V_n^{f^\ell} (t-\tb^{f^\ell};t,x,v)|^A \\
	\lesssim  & \ \int_{\tau_{k_0}}^{t- |V_n^{f^\ell} (t-\tb^{f^\ell};t,x,v)|^A} (-1) [V^{f^\ell}_\parallel (s)\cdot \nabla^2 \eta (X^{f^\ell}_\parallel(s)) \cdot V^{f^\ell}_\parallel(s) ] \\
	& \  \ \ \ \   \  \ \ \ \  \  \ \ \ \   \ \  \  \ \ \ \    \  \ \ \ \  \cdot n(X^{f^\ell}_\parallel(s))  \dd s. 
	\end{split}\Ee
	On the other hand, if $t- |V_n^{f^\ell} (t-\tb^{f^\ell};t,x,v)|^A\leq \tau_{k_0}$ then $|t-\tau_{k_0}| \leq |V_n^{f^\ell} (t-\tb^{f^\ell};t,x,v)|^A$, which implies that, from (\ref{hamilton_ODE_perp_bound}), (\ref{upper_|V|}), and (\ref{upper_X_n}),
	\Be\label{upper_V_n_2}
	|V_n^{f^\ell}(t;t,x,v)| \lesssim  |V_n^{f^\ell} (t-\tb^{f^\ell};t,x,v)|^A  \ \ for \ any \  0<A<\frac{1}{2}. 
	\Ee
	
	Now we consider $X_n^{f^\ell} (t;t,x,v)$. From (\ref{hamilton_ODE_perp_bound}) and $\dot{X}^{f^\ell}_{n}(s;t,x,v) = V^{f^\ell}_{n}(s;t,x,v)$ together with (\ref{upper_X_n}) and (\ref{upper_|V|})
	\Be\label{X_n<0}
	\begin{split}
	&X_n^{f^\ell} (t;t,x,v)\\
	 \lesssim & \   | V_n^{f^\ell} (t-\tb^{f^\ell};t,x,v)  | \\
	 &
	 + \int_{\tau_{k_0}}^t \int_{\tau_{k_0}}^\tau\underbrace{[V^{f^\ell}_\parallel (s)\cdot \nabla^2 \eta (X^{f^\ell}_\parallel(s)) \cdot V^{f^\ell}_\parallel(s) ] \cdot n(X^{f^\ell}_\parallel(s)) }_{\leq 0}  \dd s  \dd \tau \\
	\lesssim  & \   | V_n^{f^\ell} (t-\tb^{f^\ell};t,x,v)  | \\
	&
	 +  |V_n^{f^\ell} (t-\tb^{f^\ell};t,x,v)|^A\\
	 \times &  \int_{\tau_{k_0}}^{t- |V_n^{f^\ell} (t-\tb^{f^\ell};t,x,v)|^A}[V^{f^\ell}_\parallel (s)\cdot \nabla^2 \eta (X^{f^\ell}_\parallel(s)) \cdot V^{f^\ell}_\parallel(s) ] \cdot n(X^{f^\ell}_\parallel(s))   \dd s   \\
	 \lesssim & \  | V_n^{f^\ell} (t-\tb^{f^\ell};t,x,v)  | -  |V_n^{f^\ell} (t-\tb^{f^\ell};t,x,v)|^{2A}
	 \ \ \ from \ (\ref{lower_bound_V_||^2})
	 \\
	 \lesssim & \ | V_n^{f^\ell} (t-\tb^{f^\ell};t,x,v)  | -  |V_n^{f^\ell} (t-\tb^{f^\ell};t,x,v)|^{1-}\\
	 < & \ 0,
	\end{split}
	\Ee
	for $ | V_n^{f^\ell} (t-\tb^{f^\ell};t,x,v)  | \ll 1$. Clearly this cannot happen since $x \in \bar{\O}$ and $x_n\geq 0$. Therefore our assumption (\ref{assumption_lowerbound_V_n}) was wrong and we conclude (\ref{upper_V_n_2}).

		\vspace{4pt}

		\textit{Step 8-d. } From (\ref{claim_Vperp_growth}), (\ref{upper_V_n_0}), (\ref{upper_V_n_1}), and (\ref{upper_V_n_2}) in \textit{Step 8-a} and \textit{Step 8-b}, we conclude that the same estimate (\ref{upper_V_n_2}) for $|V_n^{f^\ell} (t-\tb^{f^\ell};t,x,v)|\ll 1$ in the case of (\ref{tb_upperT^**}) and (\ref{upper_|V|}). Finally from (\ref{alpha=1}), (\ref{alpha_|v|>}), (\ref{claim_Vperp_growth}), and (\ref{upper_V_n_2}) 
		Therefore we conclude that 
		\Be\label{lower_bound_V_n_final}
		|V_n^{f^\ell} (t-\tb^{f^\ell} (t,x,v);t,x,v )|\geq \left(\frac{1}{N}\right)^{1/A} \ \ (t,x,v) \in [0,T^{**}] \times (\mathcal{S}_N)^c.
		\Ee

	\hide

 From (\ref{hamilton_ODE_perp}) and (\ref{upper_|V|})
		\Be\begin{split}\label{lower_t_*1}
		&\{  3 (\e + T^{**})\sup_{\ell} \| \nabla \phi^\ell  \|_\infty\}^2 \times  \tau_*^\ell (t,x,v)
		\\
		\geq & \ \int^{t-\tb^{f^\ell} +\tau_*^\ell }_{t-\tb^{f^\ell}} (-1) [V^{f^\ell}_\parallel (s)\cdot \nabla^2 \eta (X^{f^\ell}_\parallel(s)) \cdot V^{f^\ell}_\parallel(s) ] \cdot n(X^{f^\ell}_\parallel(s))  \dd s 
		\\
\geq & \ 	V_n^{f^\ell} (\tau_*^\ell (t,x,v);t,x,v)	
+  \int^{t-\tb^{f^\ell} +\tau_*^\ell}_{t-\tb^{f^\ell}}  \nabla \phi^\ell  (s , X^{f^\ell} (s ) ) \cdot [-n(X^{f^\ell}_\parallel(s)) ] \dd s \\
& +   \int^{t-\tb^{f^\ell} +\tau_*^\ell }_{t-\tb^{f^\ell}}  X_n^{f^\ell} (s) [V^{f^\ell}_\parallel(s) \cdot \nabla^2 n (X^{f^\ell}_\parallel(s)) \cdot V^{f^\ell}_\parallel(s)]  \cdot n(X^{f^\ell}_\parallel(s)).
		\end{split}
		\Ee
	Note that, from (\ref{upper_bound_Vn}), the last two terms of the above estimates are bounded above by 
	\Be\begin{split}\label{lower_t_*2}
&
\{ \| \nabla \phi^{\ell} \|_{C^1}
+9 (\e + T^{**})^2\sup_{\ell} \| \nabla \phi^\ell  \|_\infty^2 \}
\int^{t-\tb^{f^\ell} +\tau_*^\ell }_{t-\tb^{f^\ell}}  X_n^{f^\ell} (s)\\
\leq & \ \{ \| \nabla \phi^{\ell} \|_{C^1}
+9 (\e + T^{**})^2\sup_{\ell} \| \nabla \phi^\ell  \|_\infty^2 \}
 \int^{t-\tb^{f^\ell} +\tau_*^\ell }_{t-\tb^{f^\ell}}  \int^\tau_{t-\tb^{f^\ell}}
V_n^{f^\ell} (\tau;t,x,v)
 \dd \tau   \dd s \\
 \leq & \  C_* \times  |\tau^\ell_*(t,x,v)|^2 \times V_n^{f^\ell} (t-\tb^{f^\ell}(t,x,v)) ,
	\end{split}\Ee	
		where 
		\Be\notag
		C_* := \{ \| \nabla \phi^{\ell} \|_{C^1}
+9 (\e + T^{**})^2\sup_{\ell} \| \nabla \phi^\ell  \|_\infty^2 \}
		e^{C (1+ \e  + T^{**})^2 (1+ \sup_{\ell} \| \nabla \phi^\ell  \|_{C^1} )^2 \times | \e + T^{**}| ^2}.
		\Ee
		By considering $\tau_*^\ell (t,x,v)\leq \e + T^{**} \ll1 $ and $\tau_*^\ell (t,x,v)\geq \e + T^{**}$ separately, from (\ref{lower_t_*1}) and (\ref{lower_t_*2}), we derive
		\Be\begin{split} \label{lower_t_*3}
		  \tau_*^\ell (t,x,v) 
		 \geq  & \  \tau_{**}^\ell (t,x,v)\\& \ :=  \min \left\{ \e + T^{**} ,  \frac{(1- C_* ( \e + T^{**})^2)}{1+\{  3 (\e + T^{**})\sup_{\ell} \| \nabla \phi^\ell  \|_\infty\}^2} \times V_n^{f^\ell} (t- \tb^{f^\ell} (t,x,v);t,x,v) \right\}	. 
	\end{split}	\Ee

		Now we study the lower bound of $X_n^{f^\ell} ( t- \tb^{f^\ell}+ \tau_{**}^\ell (t,x,v);t,x,v)$. From (\ref{upper_bound_Vn}) and (\ref{hamilton_ODE_perp_bound}), 
		\Be\label{lower_bound_Xn}
		\begin{split}
		&X_n^{f^\ell} ( t- \tb^{f^\ell}(t,x,v) + \tau_{*}^\ell (t,x,v);t,x,v)\\
		\geq & \ X_n^{f^\ell} ( t- \tb^{f^\ell}(t,x,v) + \tau_{**}^\ell (t,x,v);t,x,v)\\
		\geq & \    \int^{ t- \tb^{f^\ell}(t,x,v) + \tau_{**}^\ell  (t,x,v) }_{t- \tb^{f^\ell}(t,x,v) } V_n ^{f^\ell} (s;t,x,v) \dd s\\
		\geq & \  V_n^{f^\ell} (t- \tb^{f^\ell} (t,x,v);t,x,v) \times   \tau_{**}^\ell (t,x,v) \\
		& \ - 
		C (1+ \e  + T^{**})^2  ( \sup_{\ell} \| \nabla \phi^\ell  \|_{C^1}  \sup_{\ell} \| \nabla \phi^\ell  \|_{\infty} )\times  |\tau_{**}^\ell (t,x,v)|^2
		\\
		\geq & \ 
		(1- C_* ( \e + T^{**})^2) \big\{
		1- 
		C (1- C_* ( \e + T^{**})^2) (1+ \e  + T^{**})^2  ( \sup_{\ell} \| \nabla \phi^\ell  \|_{C^1}  \sup_{\ell} \| \nabla \phi^\ell  \|_{\infty} ) 
		\big\}\\
		& 
		\times
		\frac{| V_n^{f^\ell} (t- \tb^{f^\ell} (t,x,v);t,x,v)|^2}{1+\{3 (\e + T^{**})\sup_{\ell} \| \nabla \phi^\ell  \|_\infty\}^2 }\\
		\geq & \ \frac{1}{2} \frac{| V_n^{f^\ell} (t- \tb^{f^\ell} (t,x,v);t,x,v)|^2}{1+\{3 (\e + T^{**})\sup_{\ell} \| \nabla \phi^\ell  \|_\infty\}^2 },
		\end{split}
		\Ee
		as long as $\sup_{\ell} \| \nabla \phi^\ell  \|_{C^1}  \sup_{\ell} \| \nabla \phi^\ell  \|_{\infty}  \ll 1$ and $\e + T^{**} \ll 1.$

		\vspace{4pt}

		\textit{Step 8-c. } Suppose that (\ref{upper_|V|}) holds and $0 \leq V_n^{f^\ell} (t-\tb^{f^\ell}(t,x,v);t,x,v) < \delta_{**}$ with $\delta_{**}$ of  (\ref{choice_delta**}). Let us define 
		\Be\label{tau^ell}
		\tau_*^\ell (t,x,v)
		: = \sup  \big\{ \tau \geq  0: V_n^{f^\ell} (s;t,x,v)\geq 0 \ for \ all \ s \in [t-\tb^{f^\ell}(t,x,v), 
		\tau ]
		\big\}
		. 
		\Ee
		With two parameters $0< B < A \leq 1$,
		\Be\label{AB_choice}
		A, B
		\Ee
		 we split the case in three.  
		
			\vspace{2pt}
			
			\textit{\underline{Case 1.} } Suppose $|V_\parallel^{f^\ell} (\tau_*^\ell (t,x,v);t,x,v)|\geq |V_n^{f^\ell} (t-\tb^{f^\ell}(t,x,v);t,x,v)|^A$. Set 
			\Be\label{triangle_1}
			 \vartriangle =  (\e + T^{**} )^{1/2} \times |V_n^{f^\ell} (t-\tb^{f^\ell}(t,x,v);t,x,v)|^{\frac{1- 2A}{2}}.
			\Ee

			Note that 
			\Be\notag
			|V_\parallel^{f^\ell} (\tau_*^\ell + \vartriangle ;t,x,v)|\geq  |V_n^{f^\ell} (t-\tb^{f^\ell}(t,x,v);t,x,v)|^A - \| \nabla \phi^\ell \|_\infty \times  \vartriangle,
			\Ee
			and, from (\ref{convexity_eta}), 
			\Be
			\begin{split}
			&X^{f^\ell}_n(\tau_*^\ell + \vartriangle ;t,x,v)\\
			= & \ X^{f^\ell}_n(\tau_*^\ell   ;t,x,v) + \int^{\vartriangle}_0  \int_0^s \dot{V}_n^{f^\ell}(\tau_*^\ell + \tau;t,x,v) \dd \tau  \dd s \\
			\leq & \ (\e + T^{**})|V_n^{f^\ell} (  t-\tb^{f^\ell} ;t,x,v)| - C_\O \int^{\vartriangle}_0  \int_0^s  |V_n^{f^\ell}(t-\tb^{f^\ell};t,x,v)|^{2A} \dd \tau  \dd s\\
			\leq & \  0. 
			\end{split}
			\Ee 
			
			\vspace{2pt}
			
			\textit{\underline{Case 2.} } Suppose 
			\Be
			\begin{split}\label{case_2}
		 	|V_\parallel^{f^\ell} (\tau_*^\ell (t,x,v);t,x,v)|&\leq |V_n^{f^\ell} (t-\tb^{f^\ell}(t,x,v);t,x,v)|^A\\
		 	 and    \  |\nabla \phi^\ell(\tau_*^\ell (t,x,v);t,x,v)|&\leq |V_n^{f^\ell} (t-\tb^{f^\ell}(t,x,v);t,x,v)|^B.
			\end{split}\Ee

			Note that, from (\ref{XV_ell}), for $\tau \geq 0$
			\Be\notag
			|V^{f^\ell}(\tau^\ell_*-\tau;t,x,v)| \leq |V^{f^\ell}(\tau^\ell_* ;t,x,v)|+ \int_0^{\tau}  \Big|
			\frac{V^{f^\ell}(\tau^\ell_*-\tau^\prime )}{|V^{f^\ell}(\tau^\ell_*-\tau^\prime )|} \cdot 
			\nabla \phi (\tau^\ell_*-\tau^\prime, X^{f^\ell} (\tau^\ell_*-\tau^\prime ))
			\Big|\dd \tau^{\prime}.
			\Ee 
			This, together with (\ref{case_2}), deduces that
			\Be
			\begin{split}
			&\Big|
			\frac{V^{f^\ell}(\tau^\ell_*-s )}{|V^{f^\ell}(\tau^\ell_*-s )|} \cdot 
			\nabla \phi (\tau^\ell_*-s, X^{f^\ell} (\tau^\ell_*- s ))
			\Big|\\
			 \leq & \ |
			\nabla \phi (\tau^\ell_*-s, X^{f^\ell} (\tau^\ell_*- s ))
			 |\\
			 \leq&  \  |V_n^{f^\ell} (t-\tb^{f^\ell} ;t,x,v)|^B\\
			 & + \int
			\end{split}
			\Ee

			\Be
			\begin{split}
			&|\nabla \phi^\ell (X^{f^\ell}(\tau^\ell_*-s;t,x,v))| \\
			\leq&  \  |V_n^{f^\ell} (t-\tb^{f^\ell}(t,x,v);t,x,v)|^B
			\end{split}
			\Ee
			
			\vspace{2pt}
			
			\textit{\underline{Case 3.} }

		\vspace{4pt}

		\textit{Step 8-d. }

		Now we are ready to prove the claim (\ref{alpha_lower_bound_ell}). Choose $(x,v)$ in $(\mathcal{S}_N)^c$ in (\ref{test_function_condition}) and $- \e    \leq t \leq  T^{**}.$ If $|v|\leq 2 (\e + T^{**})\sup_{\ell} \| \nabla \phi^\ell \|_\infty$ and $0\leq V_n^{f^\ell} (t-\tb^{f^\ell}(t,x,v);t,x,v)< \delta_{**}$ then from (\ref{claim_Vperp_growth}) we have 
		\Be\label{lower_Vn1}
  V_n^{f^\ell} (t-\tb^{f^\ell}(t,x,v);t,x,v)
\geq  e^{-C|T^{**} + \e |^2}	V_n^{f^\ell} (s;t,x,v) 
\geq \frac{1}{e^{ C|T^{**} + \e |^2}} \times \frac{1}{N}.	
		\Ee

		Finally from (\ref{alpha=1}), (\ref{alpha_|v|>}), and (\ref{lower_Vn1}), we conclude that \Be\label{alpha_lower_bound_ell}
\inf_{\substack{ \ell \in \mathbb{N}, \ 0 \leq t \leq T^{**}
\\
 (x,v) \in (\mathcal{S}_N)^c  }}	\alpha_{f^\ell, \e} (t,x,v) \geq \frac{1}{N} \times \min \left\{    \frac{e^{- \frac{C_\O }{\sup_{\ell} \| \nabla \phi^\ell  \|_\infty}}}{C_\O }     ,\frac{e^{-C (1+ \e  + T^{**})^4 (1+ \sup_{\ell} \| \nabla \phi^\ell  \|_{C^1} )^2}}{ |T^{**} + \e |},  e^{ -C|T^{**} + \e |^2}\right\}.
	\Ee

			\unhide

		 From (2.36), (2.37), (2.40), and (2.41) in Lemma 2.4 in \cite{KL1},
		\Be\notag\begin{split}
&	\sup_{ \substack{ \ell \in \mathbb{N}
, \ 
 (x,v) \in (\mathcal{S}_N)^c , \\  - \e \leq t - \tb^{f^\ell} (t,x,v)   \leq t \leq  T^{**}}}		|\nabla_{x,v} \alpha^\beta_{f^\ell,\e}(t,x,v)|\\
  \lesssim & \  \frac{1}{|V_n^{f^\ell} (t-\tb^{f^\ell};t,x,v)|^{2- \beta}}  \lesssim_{\e,N,T^{**}  } 1.
	\end{split}	\Ee
	Hence we extract another subsequence out of all previous steps (and redefine this as $\{\ell_N\}$) such that 
	\Be\label{converge_D_alpha_ell}
	\nabla_{x,v} \alpha_{f^{\ell_N},\e}^\beta  \overset{\ast}{\rightharpoonup}  	\nabla_{x,v} \alpha_{f ,\e}^\beta  \textit{ weakly}-* \textit{ in } L^\infty
			( (-\e, T^{**}) \times (\mathcal{S}_N)^c
			).
	\Ee	
Note that the limiting function is identified from (\ref{converge_alpha_ell}).  Clearly the convergence of (\ref{af_l_3}) is an easy consequence of strong convergence of (\ref{Cauchy_L1+}) and the weak$-*$ convergence of (\ref{converge_D_alpha_ell}).

		Now we extract the final subsequence $\{\ell_*\}$ from the previous subsequence: By the Cantor's diagonal argument we define
		\Be\label{ell_*}
		\ell_*= \ell_{ \ell} .
		\Ee
		Then clearly (\ref{lim_Fell=F}) holds with this subsequence for any test function $\psi$. For any $\psi \in C^\infty_c (\bar{\O} \times \R^3 \backslash \gamma_0)$ there exists $N_{\psi} \in \mathbb{N}$ such that $supp (\psi) \subset (\mathcal{S}_{N_\psi})^c$. Then all the proofs work. 
		
		This implies (\ref{a_nabla_f_seq}) from (\ref{af_l_1}), (\ref{af_l_2}), (\ref{af_l_3}). Positivity $F= \mu+ \sqrt{\mu}f\geq 0$ comes from \textit{Step 1} and \textit{Step 6}.

		\hide
		Hence we have 
		\Be\label{uniform_bound}
		\begin{split}
			\sup_{\ell\gg 1}\sup_{0 \leq t \leq \frac{\e}{2}}\left\{
			\|(\alpha_{f^\ell,\e}(t))^\beta\|_{L^\infty (\text{supp} (\psi))}+ 
			\| \nabla_{x,v} (\alpha_{f^\ell,\e}(t))^\beta \|_{L^\infty (\text{supp} (\psi))}\right\}< \infty,\end{split}
		\Ee
		and therefore up to subsequence
		\hide
		On the other hand from (\ref{alpha_lower_bound}) and $\nabla_{x,v} \alpha_{f^\ell} \sim \frac{1}{\alpha_{f^\ell}}$ we conclude the second bound of (\ref{uniform_bound}). \hide
		
		From the change of variables $\mathcal{A}_-$ in Lemma ??
		\Be\notag
		\begin{split}
			&\| \nabla_{x,v} (\alpha_{f^\ell})^\beta \|_{L^{1+\delta} (\O \times \{ v \in \R^3: |v| \leq N\})}^{1+\delta}\\
			\lesssim & \ \iint (\alpha_{f^\ell})^{(\beta-2)(1+ \delta)} \dd x \dd v\\
			\lesssim & \ T\iint_{\gamma_-} |n(x) \cdot v|^{(\beta-2)(1+ \delta)+1} \dd S_x \dd v\\
			\lesssim & \ T 
		\end{split}
		\Ee
		where we have used $(\beta-2)(1+ \delta)+1>-1$.\unhide
		Since $(\alpha_{f^\ell})^\beta$ is bounded in $L^\infty$ and $\nabla_{x,v}(\alpha_{f^\ell})^\beta$ is bounded in $L^{1+\delta}$, we deduce that, up to subsequence,\unhide
		\Be\begin{split}\label{weak_conv}
			(\alpha_{f^\ell,\e})^\beta   \overset{\ast}{\rightharpoonup} (\alpha_{f,\e})^\beta  \ \  &\text{ in } L^\infty
			(\text{supp} (\psi)) \ \text{and} \  L^\infty
			( \gamma\cap\text{supp} (\psi))
			,\\
			\nabla_{x,v}(\alpha_{f^\ell,\e})^\beta \overset{\ast}{\rightharpoonup} \nabla_{x,v}(\alpha_{f,\e})^\beta  \ \  &\text{ in } L^\infty
			(\text{supp} (\psi)).
		\end{split}\Ee

			Now we apply (\ref{Cauchy_L1+}), (\ref{converge_alpha_ell}), and (\ref{converge_D_alpha_ell}), and conclude that the RHS of (\ref{af_l}) converges to the same terms (\ref{af_l}) with replacing $f^{\ell+1}, \alpha_{f^\ell,\e}, \nabla \alpha_{f^\ell,\e}$ by $f , \alpha_{f ,\e}, \nabla \alpha_{f ,\e}$. This implies (\ref{a_nabla_f_seq}).
		
		\unhide

		\vspace{4pt}

		\textit{Step 9. }  Choose $t>t^\prime\geq0$. Instead of expanding $h(t,x,v)$ at $t=0$ as (\ref{h_ell_local}), we expand at $t^\prime$. Then by the iteration we have (\ref{h_iteration}) replacing $0$ by $t^\prime$. Collecting all terms at time $t^\prime$, we have
		\Be\begin{split}\label{expansion_t_prime}
			&\| h(t)\|_\infty \\
			\leq &  \| h(t^\prime)\|_\infty \bigg\{
			\mathbf{1}_{t_1\leq t^\prime} e^{- \int^t_{t^\prime} \nu}  \\
			& \ \ \ \ \ \ \  
			+  \mathbf{1}_{t_1 \geq  t^\prime} \frac{e^{- \int^t_{t_1 }\nu}}{\widetilde{w}_\vartheta (V(t_1))}
			\int_{\prod_{j=1}^{k-1} \mathcal{V}_j} \sum_{l=1}^{k-1}\mathbf{1}_{\{t^{\ell-l}_{l+1}\leq
				t^\prime<t_{l}^{\ell - (l-1)}\}} \dd\Sigma _{l}(t^\prime) \bigg\}. 
		\end{split} \Ee
		Since (\ref{smeasure}) is a probability measure and $|e^{-\int^t_{t^\prime} \nu} - 1| \ll |t-t^\prime|$ for $|t-t^\prime| \ll 1$, 
		\Be 
		|(\ref{expansion_t_prime})-  \| h(t^\prime)\|_\infty| \leq  O(|t-t^\prime|)+ \int_{\prod_{j=1}^{k-1} \mathcal{V}_{j}} \mathbf{1}_{\{ t^{k} (t,x,v,u^{1}, \cdots , u^{k-1}) >0 \}} \dd \Sigma_{k-1}^{k-1} 
		.\notag
		\Ee
		Then by (\ref{h_iteration}) we have $
		\| h(t)\|_\infty - \| h (t^\prime) \|_\infty < \frac{1}{2^k}
		+O_k(|t-t^\prime|).$ For large $k$, choosing $|t-t^\prime|\ll1$, we can prove $
		\| h(t)\|_\infty - \| h (t^\prime) \|_\infty\ll 1$ as $|t-t^\prime| \ll 1$. 
		
		Now we can expand $h(t^\prime,x,v)$ at $t$ by (\ref{h_ell_local}). Following the same argument we have $\| h(t^\prime) \|_\infty - \| h(t ) \|_\infty \ll 1$ as $|t-t^\prime| \ll 1$. Hence $\| w_{\vartheta}f(t) \|_\infty$ is continuous in $t$.

		The continuity of $\| \nabla_v f (t) \|_{L^3_xL^{1+ \delta}_v}$ and $\| w_{\tilde{\vartheta}}\alpha_{f,\e }^\beta \nabla_{x,v} f (t) \|_{p} ^p
		+ 
		\int^t_0  |w_{\tilde{\vartheta}} \alpha_{f,\e }^\beta \nabla_{x,v} f (t) |_{p,+}^p$ is an easy consequence of (\ref{g_initial})-(\ref{g_infty}), and (\ref{W1p_pf}) as well.
		\hide

		Then 
		\Be
		\begin{split}\notag
			& (\ref{af_l})_1 - \int^t_0\iint_{\O \times \R^3} \nabla_{x,v} \alpha^\beta_{f} f \psi  \\
			= & \ \int^t_0\iint_{\O \times \R^3}  \nabla_{x,v}\alpha^\beta_{f^\ell} \{ f^{\ell+1}
			-f\}
			\psi + \int^t_0\iint_{\O \times \R^3}\{\nabla_{x,v}\alpha^\beta_{f^\ell}
			- \nabla_{x,v}\alpha^\beta_{f }\}  f
			\psi  .
		\end{split}
		\Ee
		The first term of RHS is bounded by $\| \nabla_{x,v} \alpha^\beta_{f^\ell} \|_{L^{\infty}} \| f^{\ell+1} - f \|_{L^{1+} } \|\psi\|_{L^{(1+ \delta)^\prime}}\lesssim \| f^{\ell+1} - f \|_\infty \rightarrow 0$ as $\ell \rightarrow \infty$. The second term of RHS goes to zero from (\ref{weak_conv}).
		
		On the other hand,
		\Be\begin{split}\notag
			&(\ref{af_l})_2 - \int^t_0\iint_{\O \times \R^3} \alpha^\beta_{f} f \nabla_{x,v} \psi \dd x \dd v\\
			=& \ \int^t_0\iint_{\O \times \R^3} \alpha^\beta_{f^\ell}  \{ f^{\ell+1} - f\}\nabla_{x,v} \psi  + \int^t_0\iint_{\O \times \R^3} \{\alpha_{f^\ell}^\beta - \alpha^\beta_f \} f\nabla_{x,v} \psi.
		\end{split}\Ee
		The first term of RHS is bounded by $\| \alpha^\beta_{f^\ell}\|_{L^\infty} \| f^{\ell+1}-f\|_{L^\infty}\rightarrow 0$ as $\ell \rightarrow \infty$. From (\ref{weak_conv}) and $f \nabla_{x,v}\psi \in L^1$ the second term of RHS converges to zero as $\ell \rightarrow \infty$.
		
		For $(\ref{af_l})_3$ we utilize the trace theorem and deduce $f^{\ell} \rightarrow f$ strongly in $L^\infty (\p\O \times \R^3)$. Then 
		\Be
		\begin{split} 
			&(\ref{af_l})_3 - \int^t_0\iint_\gamma \alpha^\beta_f f \psi \\
			=& \  \int^t_0\iint_\gamma \alpha^\beta_{f^\ell} \{f^{\ell+1} -f \} \psi 
			+\int^t_0 \iint_\gamma \{\alpha^\beta_{f^\ell} - \alpha^\beta_f \} f \psi.
		\end{split}
		\Ee
		We split with some $0<\delta \ll1$ as $\iint_\gamma \alpha^\beta_{f^\ell} \{f^{\ell+1} -f \} \psi = \iint_{\gamma^\delta} + \iint_{\gamma \backslash\gamma^\delta}$. By (\ref{uniform_bound}), we derive $\iint_{\gamma^\delta} \lesssim O(\delta) \| \alpha^\beta_{f^\ell}\|_\infty \{ \| f^{\ell+1}\|_\infty + \| f\|_\infty\}< O(\delta)$. For $\gamma\backslash \gamma^\delta$, we use the trace theorem to bound as $t \| \alpha^\beta_{f^\ell}\|_\infty\| f^{\ell+1} - f \|_{L^{1+ }}\rightarrow 0$ as $\ell \uparrow \infty$.

		The first term of RHS is bounded by $\| \alpha^\beta_{f^\ell} \|_{L^{\infty } (\p\O \times \R^3)}\| f^{\ell+1} - f \|_{L^{\infty } (\p\O \times \R^3)} \rightarrow 0$ as $\ell \rightarrow 0$. The second term of RHS converges to zero due to (\ref{weak_conv}).\unhide\hide

		\textit{Step 4. } We claim that 
		\Be
		\sup_\ell\|f^\ell\|_{L^3_x L^{1+}_v} < \infty,
		\Ee
		and 
		\Be
		f^{\ell} \ \text{ is Cauchy in } \ L^{1+}(\O \times \R^3).  
		\Ee
		
		Note that $f^{\ell+1} - f^\ell$ satisfies $(f^{\ell+1} - f^\ell)|_{t=0} \equiv0$ and 
		\Be\begin{split}\label{f_ell-ell}
			&[\p_t + v\cdot \nabla_x - \nabla_x \phi^\ell \cdot \nabla_v+ \nu_{\phi^\ell}] (f^{\ell+1} - f^\ell)\\
			= & \ 
			(\nabla \phi^{\ell }-  \nabla \phi^{\ell-1}) \cdot \nabla_v f^\ell
			+K(f^{\ell } - f^{\ell-1}) - v\cdot \nabla (\phi^{\ell} - \phi^{\ell-1})\sqrt{\mu}\\
			& + w \Gamma_{\text{gain}} (\frac{f^\ell}{w}, \frac{f^\ell - f^{\ell-1}}{w})+ w \Gamma_{\text{gain}} (\frac{f^\ell- f^{\ell-1}}{w}, \frac{f^{\ell-1}}{w})\\
			&- w \Gamma_{\text{loss}} ( \frac{f^\ell}{w}, \frac{f^{\ell+1} - f^\ell}{w})-w \Gamma_{\text{loss}} ( \frac{f^\ell - f^{\ell-1}}{w}, \frac{ f^\ell}{w}).
		\end{split}\Ee
		Then we follow the proof of Proposition \ref{L1+stability}
		\Be
		\sup_{0 \leq s \leq t}\| f^{\ell+1} (s)- f^\ell (s) \|_{L^{1+ \delta}(\O \times \R^3)}\leq O(t)\sup_{0 \leq s \leq t}\| f^{\ell } (s)- f^{\ell-1} (s) \|_{L^{1+ \delta}(\O \times \R^3)}
		\Ee
		Then we deduce that 
		\Be
		\sup_{0 \leq s \leq t}\| f^{\ell } (s)- f^m (s) \|_{L^{1+ \delta}(\O \times \R^3)} \leq O(t)^{\min \{\ell, m\}}.
		\Ee

		\unhide
	\end{proof}

	\hide
	We have 
	\Be\label{Gamma_infty_2}
	\|\nu^{-1/2} \Gamma(f,f)\|_2\lesssim \| w f \|_\infty \| f \|_2.
	\Ee
	From $|(v-u) \cdot \o|\leq |v| + |u|$,
	\Be
	\begin{split}
		&  \nu(v)^{-1} \Gamma(f,f)(v) \\
		\lesssim& \ \int_{\R^3}\int_{\S^2} |(v-u) \cdot \o| \langle v\rangle^{-1} \sqrt{\mu(u)}\{ |f(v^\prime)| |f(u^\prime)| + |f(u)| |f(v)| \} \dd \o \dd u\\
		\lesssim & \ \int_{\R^3}\int_{\S^2} \mu(u)^{1/4} \{ |f(v^\prime)| |f(u^\prime)| + |f(u)| |f(v)| \} \dd \o \dd u
	\end{split}\Ee
	Then by the standard change of variables $(v^\prime,u^\prime) \rightarrow (v,u)$ 
	\Be
	\begin{split}
		& \left\| 
		\iint_{\R^3 \times \S^2}  \mu(u)^{1/4} |f(u^\prime)| |f(v^\prime)|\dd \o \dd u
		\right\|_{L^2_v} \\
		\lesssim& \ \| w f\|_\infty \left[\iint_{\R^3 \times \S^2} \mu(u) w(u^\prime)^{-1} \dd \o \dd u\right]^{1/2}\\
		& \times 
		\left[
		\iiint_{\R^3 \times \R^3 \times \S^2} |f(v^\prime)|^2 \mu(u) w(u^\prime)^{-1} \dd \o \dd u^\prime \dd v^\prime
		\right]^{1/2}\\
		\lesssim & \ \| w f\|_\infty \| f \|_2
	\end{split}
	\Ee
	
	Clearly $ 
	\left\| 
	\iint_{\R^3 \times \S^2}  \mu(u)^{1/4} |f(u)| |f(v)|\dd \o \dd u
	\right\|_{L^2} 
	\lesssim   \ \| w f \|_\infty   \|f \|_{L^2}.
	$\unhide
	
	\hide

	\section{$L^2$ coercivity}

	The main purpose of this section is to prove the following:
	
	\begin{proposition}
		\label{dlinearl2}Suppose $(f, \phi_f)$ solves (\ref{eqtn_f}), (\ref{phi_f}), (\ref{BC_f}), and (\ref{phi_BC}). Then 
		%
		%
		there is $0<\lambda_{L^2} \ll 1$ such that for $0 \leq s \leq t$, 
		\Be\begin{split}\label{completes_dyn}
			& \| e^{\lambda_{L^2} t}f(t)\|_2^2
			+ \| e^{\lambda_{L^2} t} \nabla \phi (t) \|_2^2\\
			&
			+  \int_s^t \| e^{\lambda_{L^2} \tau}  f (\tau)\|_\nu^2 
			+ \| e^{\lambda_{L^2} \tau} \nabla \phi_f(\tau)\|_2^2 
			\mathrm{d} \tau 
			+  \int_s^t | e^{\lambda_{L^2} \tau} f |^2_{\gamma, 2 }    \\
			\lesssim & \ \|e^{\lambda_{L^2} s} f(s)\|_2^2 +   \|e^{\lambda_{L^2} s} \nabla \phi_f(s)\|_2^2  \\&   + \sup_{s \leq \tau \leq t} \| w f (\tau) \|_\infty \int_s^t \| e^{\lambda_{L^2} \tau}  f (\tau)\|_\nu^2  .  
		\end{split}\Ee
	\end{proposition}
	
	In order to prove the proposition we need the following:
	
	\begin{lemma}
		\label{dabc}
		There exists
		a function $G(t)$ such that, for all $0\le s\leq t$, $G(s)\lesssim
		\|f(s)\|_{2}^{2}$ and 
		\begin{equation}\begin{split}\label{estimate_dabc}
		&\int_{s}^{t}\|\mathbf{P}f(\tau)\|_{\nu }^{2}  + \int^t_s \| \nabla \phi_f \|_2^2\\
		\lesssim & \ 
		G(t)-G(s)
		+  \int_{s}^{t}\|(\mathbf{I}-\mathbf{P}%
		)f(\tau)\|_{\nu }^{2} +\int_{s}^{t}|(1-P_{\gamma })f(\tau)|_{2,+}^{2} \\
		& +\int_{s}^{t}\| \nu^{-1/2} { {\Gamma(f,f)} } \|_{2}^{2}.
		\end{split}\end{equation}
	\end{lemma}
	
	\begin{proof}[Proof of Proposition \ref{dlinearl2}] 
		\textit{Step 1. }  Without loss of generality we prove the result with $s=0$. We have an $L^2$-estimate from $e^{\lambda_{L^2} t} \times (\ref{eqtn_f})$
		\Be\notag
		\begin{split}
			&\|e^{\lambda t}  f(t) \|_2^2 - \| f(0) \|_2^2 +  \int^{t}_{0} | e^{\lambda \tau }
			(1-P_{\gamma} ) f^{j}|_{2,+}^{2}
			\\
			& + \int_0^t  \iint_{\O \times \R^3} v \cdot \nabla \phi_f e^{2\lambda \tau} |f|^2+
			2\int^t_0 \iint_{\O \times \R^3}e^{\lambda \tau}2  f Lf\\
			= &\ 2\int^t_0 \iint_{\O \times \R^3} e^{2\lambda \tau}f \Gamma(f,f) - 2\int^t_0 e^{2\lambda \tau}\int_{\O  } \nabla \phi_f \cdot \int_{\R^3} v \sqrt{\mu} f\\
			& 
			+2 \lambda  \int_0^t \|e^{\lambda \tau}  f(\tau) \|_2^2
			.
		\end{split}
		\Ee 
		
		On the other hand multiplying $\sqrt{\mu(v)} \phi_f(t,x)$ with a test function $\psi(t,x)$ to (\ref{eqtn_f}) and applying the Green's identity, we obtain
		\Be\begin{split}\notag
			&
			\int_{\O} \nabla \phi_f(t,x)  \cdot \int_{\R^3}  v \sqrt{\mu} f \dd v\dd x
			\\
			= & \ \int_{\O} \phi_f(t,x) \p_\tau\left(\int_{\R^3} f (\tau) \sqrt{\mu} \dd v\right)  \dd x
			+ \iint_{\p\O \times \R^3}  \phi_f(t,x) f \sqrt{\mu} \{n \cdot v\} \dd v\dd S_x .
		\end{split}\Ee
		From (\ref{null_flux}), the last boundary contribution equals zero. Now we use (\ref{phi_f}) and (\ref{phi_BC}) and deduce that 
		\Be
		\begin{split}\notag
			& \int^t_0e^{2 \lambda_{L^2} \tau}\int_{\O} \phi_f(t,x) \p_\tau\left(\int_{\R^3} f(\tau)  \sqrt{\mu} \dd v\right)  \dd x \dd \tau\\
			= & \ -  \int^t_0e^{2 \lambda_{L^2} \tau}\int_{\O} \phi_f(t,x) \p_\tau  \Delta_x  \phi_f(\tau,x) \dd x \dd \tau\\
			= &  \  \frac{1}{2}\int^t_0e^{2 \lambda_{L^2} \tau}\int_{\O}  \p_\tau  |\nabla_x  \phi_f(\tau,x)|^2 \dd x \dd \tau
			\\
			= & \ \frac{1}{2}   \left(\int_{\O} e^{2 \lambda_{L^2} t} |\nabla_x \phi_f (t,x)|^2  \dd x\right)- \frac{1}{2}   \left(\int_{\O} |\nabla_x \phi_f (0,x)|^2  \dd x\right)\\
			& \ 
			- \lambda \int^t_0 e^{2 \lambda_{L^2} \tau} \int_\O |\nabla_x  \phi_f(\tau,x)|^2 \dd x \dd \tau
			.
		\end{split}
		\Ee
		
		Hence we derive 
		\Be\notag
		\begin{split}
			&\| e^{\lambda_{L^2} t} f(t) \|_2^2 + \| e^{\lambda_{L^2} t}\phi_f (t) \|_2^2  + \int_0^t  \iint_{\O \times \R^3}
			e^{2\lambda_{L^2} \tau}
			v \cdot \nabla \phi_f |f|^2\\
			& +
			2C\int^t_0 \iint_{\O \times \R^3}\| e^{\lambda_{L^2} \tau} (\mathbf{I} - \mathbf{P}) f \|_\nu^2
			+ \int^{t}_{0} | e^{\lambda_{L^2} \tau }
			(1-P_{\gamma} ) f^{j}|_{2,+}^{2}
			\\
			\lesssim &\ 
			\| f(0) \|_2^2+ \| \phi_f (0) \|_2^2 +
			\int^t_0 \| e^{\lambda_{L^2} \tau}  \nu^{-1/2} \Gamma(f,f) \|_2^2\\
			& 
			+\{ \lambda_{L^2} + o(1)\} \int^t_0 \| e^{\lambda_{L^2} \tau} f \|_\nu^2 + \lambda_{L^2} \int^t_0 \| e^{\lambda_{L^2} \tau} \nabla_x \phi_f \|_2^2
			.
		\end{split}
		\Ee
		Now we apply Lemma \ref{dabc} and add $o(1) \times (\ref{estimate_dabc})$ to the above inequality and choose $0< \lambda_{L^2} \ll 1$ to conclude (\ref{completes_dyn}) except the full boundary control. 
		
		\vspace{4pt}
		
		\textit{Step 2. } Note that $P_\gamma f = z(t,x)\sqrt{\mu(v)}$ for a suitable function $z(t,x)$ on the boundary. Then for $0 < \e \ll 1$
		\Be
		\begin{split}\notag
			| P_\gamma f |_{\gamma,2}^2 
			= & \ \int_{\p\O } |z(t,x)|^2 \dd x \times \int_{\R^3} \mu (v) |n(x) \cdot v| \dd v\\
			\lesssim & \  \int_{\p\O } |z(t,x)|^2 \dd x \times \int_{ \gamma_+(x) \backslash \gamma_+^\e(x)} \mu (v)^{3/2} |n(x) \cdot v| \dd v\\
			= & \  | \mathbf{1}_{ \gamma_+  \backslash \gamma_+^\e }\mu^{1/4} P_\gamma f   |_{ +,2}^2.
		\end{split}
		\Ee
		Since $P_\gamma f = f - (1- P_\gamma) f$ on $\gamma_+$ we have $
		| \mathbf{1}_{ \gamma_+  \backslash \gamma_+^\e } \mu^{1/4} P_\gamma f|_{ +,2}^2
		\lesssim |\mathbf{1}_{ \gamma_+  \backslash \gamma_+^\e }\mu^{1/4}  f |_{ +,2}^2 + |(1- P_\gamma) f|_{ +,2}^2.$ Therefore
		\Be\label{Pgamma_bound}
		\int^t_0 | P_\gamma f|_{\gamma,2}^2 \lesssim \int^t_0 |\mathbf{1}_{ \gamma_+  \backslash \gamma_+^\e }\mu^{1/4}  f |_{ +,2}^2 +\int^t_0 |(1- P_\gamma) f|_{ +,2}^2.
		\Ee
		Note that 
		\Be\begin{split}\notag
			&\big|[\p_t + v\cdot \nabla_x - \nabla \phi \cdot \nabla_v] (\mu^{1/4}  f )\big|\\
			\lesssim&\ \mu^{1/4} \{ |v|| \nabla_x \phi|  f  +  |v|| \nabla_x \phi|  +  |Lf| + |\Gamma (f,f) | \}
		\end{split}\Ee
		By the trace theorem Lemma \ref{le:ukai},  
		\Be\label{est_trans_f}
		\begin{split}
			& \int^t_0 |\mathbf{1}_{ \gamma_+  \backslash \gamma_+^\e }\mu^{1/4}  f |_{ +,2}^2\\
			\lesssim & \ \| f_0 \|_2 + (1 + \| wf \|_\infty) \int^t_0 \| f \|_2^2 + \int^t_0\| \nabla \phi \|_2^2.
		\end{split}
		\Ee
		Adding $o(1) \times (\ref{Pgamma_bound})$ to the result of \text{Step 1} and using (\ref{est_trans_f}) we conclude (\ref{completes_dyn}).\end{proof}


	\begin{proof}[Proof of Lemma \ref{dabc}]
		From the Green's identity (\ref{}), a solution $(f,\phi_f)$ satisfies
		\begin{equation}\begin{split}
		&\iint_{\O \times \R^3} f(t) \psi(t) - \iint_{\O \times \R^3} f(s) \psi(s)
		- \int_s^t \iint_{\O \times \R^3} f \p_t\psi \\
		&+
		\int^t_s\int_{\gamma }\psi f -\int^t_s\iint_{\Omega \times  \R^{3}} {v}\cdot \nabla
		_{x}\psi f + \int_s^t \iint_{\O \times \R^3} \sqrt{\mu} f \nabla_x \phi_f \cdot \nabla_v \Big[\frac{1}{\sqrt{\mu}} \psi\Big]
		\\
		&=\int^t_s\iint_{\Omega \times \R^{3}}\psi \{ -L  (\mathbf{I}-\mathbf{P%
		})f + \Gamma(f,f)\}- \underbrace{\int^t_s\iint_{\Omega \times  \R^{3}} \psi v\cdot \nabla_x \phi_f \sqrt{\mu}}_{(\ref{weakformulation})_{\phi_f}}  .  \label{weakformulation}
		\end{split}\end{equation}%
		
		As \cite{EGKM, EGKM2} we use a set of test functions:
		\begin{eqnarray}
		&&\psi_a  \equiv   (|v|^{2}-\beta_{a} )\sqrt{\mu }v\cdot\nabla_x\phi _{a},\label{phia}\\
		&&\psi^{i,j}_{b,1}\equiv (v_{i}^{2}-\beta_ b)\sqrt{\mu }\partial _{j}\phi _{b}^{j}, \quad i,j=1,2,3,\label{phibj} \\
		&&\psi^{i,j}_{b,2}\equiv |v|^{2}v_{i}v_{j}\sqrt{\mu }\partial _{j}\phi _{b}^{i}(x),\quad i\neq j,\label{phibij}\\
		&&\psi_c \equiv(|v|^{2}-\beta_c )\sqrt{\mu }v \cdot \nabla_x \phi _{c}(x),\label{phic}
		\end{eqnarray}
		where
		\Be\begin{split}\label{phi_abc}
			- \Delta \phi_a &= a,  \ \ \frac{\p \phi_a}{\p n} \Big|_{\p\O} =0,\\
			- \Delta \phi_b^j &= b_j,   \ \ \phi^j_b|_{\p\O} =0,\  \text{and}  \ - \Delta \phi_c  = c,   \ \ \phi_c|_{\p\O} =0,
		\end{split}
		\Ee
		and $\beta_a=10$, $\beta_b=1$, and $\beta_c= 5$ such that for all $i=1,2,3,$
		\Be \label{beta}
		\begin{split}
			&\int_{{\R}^{3}}(|v|^{2}-\beta_a )(\frac{|v|^{2}}{2}-\frac{3}{2}%
			) v_{i} ^{2}\mu(v)\dd v =0 ,\\
			&\int_{\R^3} (v_i^2 - \beta_b) \mu(v) \dd v=0 ,\\
			&\int_{\R^3} (|v|^{2}-\beta_c )v_{i}^{2}\mu(v) \dd v=0. 
		\end{split}
		\Ee
		
		\vspace{4pt}

		\textit{Step 1. } Estimate of $(\ref{weakformulation})_{\phi_f}$: From (\ref{phia})-(\ref{phic}) and (\ref{beta}), we have $(\ref{weakformulation})_{\phi_f}\equiv0$ for $\psi_{b,1}^{i,j}$, $\psi_{b,2}^{i,j}$, and $\psi_c$. For $\psi_a$, $(\ref{weakformulation})_{\phi_f}$ equals
		\Be\begin{split}\notag
			\int_{\R^3} (|v|^2- \beta_a) (v_1)^2\mu \dd v\times\int^t_s\int_{\O  } \nabla \phi_a \cdot \nabla \phi_f.
		\end{split}\Ee
		Note that $\phi_a = \phi_f$ from the definitions of (\ref{phi_f}) and (\ref{phia}). Therefore 
		\Be
		(\ref{weakformulation})_{\phi_f} = C \int^t_s \| \nabla \phi_f \|_2^2 \ \ \text{for} \ \ \psi= \psi_a.
		\Ee
		
		{\color{red} 
			Complete the rest of proof. I copy \cite{EGKM} for your convenience. Make it shorter.}
		
		{\it Step 2}. {Estimate of} \ ${c}$
		
		We have
		\begin{equation*}
		{v}\cdot \nabla _{x}\psi_c =\sum_{i,j=1}^{d}(|v|^{2}-\beta_c )v_{i}\sqrt{\mu }%
		v_{j}\partial _{i j}\phi _{c}(x),
		\end{equation*}%
		so that the left hand side of (\ref{weakformulation}) takes the form, for $i=1,\cdots,d,$
		\begin{eqnarray}
		&&\int_{\partial \Omega \times \mathbf{R}^{3}}(n(x)\cdot v)(|v|^{2}-\beta_c )\sqrt{\mu }%
		\sum_{i=1}^d v_{i}\partial _{i}\phi _{c}f\notag\\
		&&-\iint_{\Omega \times \mathbf{R}%
			^{3}}(|v|^{2}-\beta_c )\sqrt{\mu }\Big\{\sum_{i,j=1}^d v_{i}v_{k}\partial _{ij}\phi
		_{c} \Big\}f .  \label{lweak}
		\end{eqnarray}%
		We decompose
		\begin{eqnarray}
		f_{\gamma } &=&P_{\gamma }f +\mathbf{1}_{\gamma
			_{+}}(1-P_{\gamma })f+\mathbf{1}_{\gamma _{-}}r,\text{ \ \  \ \ \ \ \ \ \ \ \ \ \ \ \ \ \ \ \ on
		}\gamma,  \label{bsplit} \\
		f &=&\Big\{a+v\cdot b+c[\frac{|v|^{2}}{2}-\frac{3}{2}]\Big\}\sqrt{\mu }+(\mathbf{I}-%
		\mathbf{P})f,\text{ \ \ \ \ \  on }\Omega\times\mathbf{R}^3. \label{insidesplit}
		\end{eqnarray}%
		We shall choose $\beta_c $ so
		that, for all $i$
		\begin{equation}
		\int (|v|^{2}-\beta_c )v_{i}^{2}\mu(v) dv=0.  \label{beta}
		\end{equation}%
		Since $\mu= \frac{1}{2\pi}e^{-\frac{|v|^2}{2}}$ the desired value of $\beta_c$  is $\beta_c=5$.
		Because of the choice of $\beta_c$, there is no $a$ contribution in the bulk and no $P_{\gamma }f$ contribution at
		the boundary in (\ref{lweak}).
		
		Therefore, substituting (\ref{bsplit}) and (\ref{insidesplit}) into (\ref%
		{lweak}), since the $b$ terms and the off-diagonal $c$ terms also vanish by oddness in $v$ in the bulk, the left hand side of (\ref{weakformulation}) becomes%
		\begin{eqnarray*}
			&& \sum_{i=1}^d \int_{\gamma }(|v|^{2}-\beta_c ){v_{i}}^{2}n_{i}\sqrt{\mu }\partial _{i}\phi
			_{c}[(1-P_{\gamma })f \mathbf{1}_{\gamma _{+}}+r\mathbf{1}_{\gamma _{-}}]
			\\
			&&-\sum_{i=1}^d \int_{\mathbf{R}^3}(|v|^{2}-\beta_c )v_{i}^{2} (\frac{|v|^{2}}{2}-\frac{3}{2})\mu(v)dv \int_{\Omega } \partial
			_{ii}\phi _{c}(t,x)c(t,x) dx \\
			&&-\sum_{i=1}^d \iint_{\Omega \times \mathbf{R}^{3}}(|v|^{2}-\beta_c )v_{i}\sqrt{\mu }%
			( {v}\cdot \nabla_x )\partial _{i}\phi _{c}(\mathbf{I}-\mathbf{P})f.
		\end{eqnarray*}%
		For $\beta_c =5$, $\int_{\mathbf{R}^3}(|v|^{2}-\beta_c )v_{i}^{2} (\frac{|v|^{2}}{2}-\frac{3}{2})\mu(v)dv= 10\pi\sqrt{2\pi}$.
		Therefore, we obtain from (\ref{rweak})
		\begin{multline*}
		- 10\pi\sqrt{2\pi} \int_{\Omega }\Delta_x \phi _{c}(x)c(x)
		\lesssim \ \|c\|_{2 }\{|(1-P_{\gamma })f |_{2,+}+\|(\mathbf{I}-\mathbf{P}%
		)f\|_{2}+\|g\|_{2}+|r|_{2}\},
		\end{multline*}%
		where we have used the elliptic estimate and the trace estimate:
		\begin{equation*}
		|\nabla_x \phi _{c}|_{2 }\lesssim \|\phi _{c}\|_{H^{2} }\lesssim \|c\|_{2}.
		\end{equation*}%
		Since $-\Delta_x \phi_c = c$, from (\ref{phic}) we obtain
		\begin{equation*}
		\|c\|_{2 }^{2} \ \lesssim \ \big\{|(1-P_{\gamma })f |_{2,+}+\|(\mathbf{I}-\mathbf{P}%
		)f\|_{2}+\|g\|_{2}+|r|_{2}\big\}\|c\|_{2 },\end{equation*}%
		and hence
		\begin{equation}
		\|c\|_{2 }^2 \ \lesssim \ |(1-P_{\gamma })f |_{2,+}^{2}+\|(\mathbf{I}-\mathbf{P}%
		)f\|_{2}^{2}+\|g\|_{2}^{2}+|r|_{2}^{2}.  \label{cestimate}
		\end{equation}
		
		\noindent{\it Step 2.} {Estimate of} \ $ {b}$
		
		We shall establish the estimate of $b$ by estimating   $(\partial_i \partial_j \Delta^{-1}b_j) b_i$ for all $i,j=1,\dots, d$, and  $(\partial_j\partial_j \Delta^{-1} b_i) b_i$ for $i\neq j$.
		
		We fix $i,j$.  To estimate $\partial_i \partial_j \Delta^{-1}b_j b_i$ we choose as test function in (\ref{weakformulation})
		\begin{equation}
		\psi = \psi^{i,j}_b\equiv (v_{i}^{2}-\beta_ b)\sqrt{\mu }\partial _{j}\phi _{b}^{j}, \quad i,j=1,\dots, d,  \label{phibj}
		\end{equation}%
		where $\beta_b$ is a constant to be determined, and
		\begin{equation}
		-\Delta_x \phi _{b}^{j}(x)=b_{j}(x)  , \ \ \ \phi_b^{j}|_{\partial\Omega}=0.\label{jb}
		\end{equation}%
		From the standard elliptic estimate%
		\begin{equation*}
		\|\phi _{b}^{j}\|_{H^{2} }\lesssim \|b\|_{2 }.
		\end{equation*}%
		Hence the right hand side of (\ref{weakformulation}) is now bounded by%
		\begin{equation}
		\|b\|_{2 }\big\{\|(\mathbf{I}-\mathbf{P})f\|_{2}+\|g\|_{2}\big\}.
		\label{rweakbj}
		\end{equation}%
		Now substitute (\ref{bsplit}) and (\ref{insidesplit}) into the left hand side of (\ref{weakformulation}). Note that $(v_{i}^{2}-\beta_b  )\{n(x)\cdot  {v}\}\mu $ is odd in $v$, therefore $P_{\gamma
		}f $ contributions to (\ref{weakformulation}) vanishes. Moreover, by (\ref{insidesplit}), the $a, c$ contributions to (\ref{weakformulation}) also vanish by oddness. Therefore the left hand side of (\ref{weakformulation}) takes the form %
		\begin{eqnarray}
		&&\int_{\partial \Omega \times \mathbf{R}^{3}}(n(x)\cdot  {v})(v_{i}^{2}-\beta_b  )\sqrt{\mu }%
		\partial _{j}\phi _{b}^{j}f-\iint_{\Omega \times \mathbf{R}%
			^{3}}(v_{i}^{2}-\beta_b  )\sqrt{\mu }\{\sum_{l}v_{l} \partial_{l j}\phi
		_{b}^{j}\}f  \notag \\
		&=&\int_{\partial \Omega \times \mathbf{R}^{3}}(n(x)\cdot  {v})  (v_{i}^{2}-\beta_b  )
		\sqrt{\mu }\partial _{j}\phi _{b}^{j}[(1-P_{\gamma })f +r]\mathbf{1}%
		_{\gamma _{+}}  \label{sti_6} \\
		&&-\int_{\Omega } { {  \int_{\mathbf{R}%
					^{3}}\sum_{l}(v_{i}^{2}-\beta_b )v_{l}^{2}\mu \partial _{lj}\phi
				_{b}^{j}(x)b_{l}dv } } dx \label{sti_8} \\
		&&-\iint_{\Omega \times \mathbf{R}^{3}}\sum_{l}(v_{i}^{2}-\beta_b  )v_{l}\sqrt{\mu }%
		\partial _{l j} \phi _{b}^{j}(x)(\mathbf{I}-\mathbf{P})f.  \notag
		\end{eqnarray}%
		Furthermore, since $\mu (v)=\frac{1}{2\pi}\prod_{i=1}^3  e^{-\frac{|v_i|^2}{2}}$ we can
		choose $\beta_b >0$ such that for all $i,$
		\begin{equation}
		\int_{\mathbf{R}^{3}}[(v_{i})^{2}-\beta_b ]\mu (v) dv=  \int_{%
			\mathbf{R}}[v_{1}^{2}-\beta_b]  e^{-\frac{|v_1|^2}{2}} dv_1=0.  \label{alpha}
		\end{equation}%
		We remark that the choice (\ref{alpha}) also plays a crucial rule in the dynamical estimate (\ref{i_7}). Since $\mu (v)=\frac{1}{2\pi}e^{-\frac{|v|^{2}}{2}}$, the desired value is $\beta_b=1$.
		
		Note that for such chosen $\beta_b
		, $ and for $i\neq k$, by an explicit computation
		\begin{eqnarray*}
			\int (v_{i}^{2}-\beta_b )v_{k}^{2}\mu dv  &=& \int (v_{1}^{2}-\beta_b
			)v_{2}^{2}  \frac{1}{2\pi}e^{-\frac{|v_1|^2}{2}} e^{-\frac{|v_2|^2}{2}} e^{-\frac{|v_3|^2}{2}} dv\\
			&=&  \int_{\mathbf{R}}
			(v_{1}^{2}-\beta_b ) e^{-\frac{|v_1|^2}{2}} dv_1=0, \\
			\int (v_{i}^{2}-\beta_b )v_{i}^{2}\mu dv & = &  \int_{\mathbf{R}}
			[v_{1}^{4}-\beta_b v_{1}^{2}]e^{-\frac{|v_1|^2}{2}}dv_{1}= 2\sqrt{2\pi}\neq 0.
		\end{eqnarray*}%
		Therefore, (\ref{sti_8}) becomes, by (\ref{phibj}),
		\begin{eqnarray*}
			&&-\iint_{\Omega\times\mathbf{R}^3} (v_{i}^{2}-\beta_b )v_{i}^{2}\mu dv \partial _{ij}\phi^j
			_{b}(x)b_{i}+\sum_{k(\neq i)}\underbrace{[\int_{\mathbf{R}^3} (v_{i}^{2}-\beta_b )v_{k}^{2}\mu
				]}_{=0}\int_{\Omega}\partial _{kj}\phi _{b}^j (x)b_{k}  \\
			&&=2\sqrt{2\pi}\int_{\Omega }(\partial
			_{i}\partial _{j}\Delta ^{-1}b_{j})b_{i}.
		\end{eqnarray*}%
		Hence we have the following estimate for all $i,j$, by (\ref{rweakbj}):
		\begin{eqnarray}
		\left|\int_{\Omega }\partial _{i}\partial _{j}\Delta ^{-1}b_{j}b_{i}\right|
		\lesssim  |(1-P_{\gamma })f |_{2,+}^{2}+\|(\mathbf{I}-\mathbf{P}%
		)f\|^{2}_{2}+ \|g\|_{2}^{2}+|r|_{2}^{2} +\varepsilon \|b\|_{2}^{2}.  \label{stiijjkiki}
		\end{eqnarray}%
		In order to estimate $ \partial _{j}(\partial
		_{j}\Delta ^{-1}b_{i})b_{i}$ for $i\neq j$, we choose as test function in (\ref{weakformulation})
		\begin{equation}
		\psi =
		|v|^{2}v_{i}v_{j}\sqrt{\mu }\partial _{j}\phi _{b}^{i}(x),\quad i\neq j,
		\label{phibij}
		\end{equation}%
		where $\phi_b^i$ is given by (\ref{jb}). Clearly, the right hand side of (\ref{weakformulation}) is again bounded by (\ref{rweakbj}%
		). We substitute again (\ref{bsplit}) and (\ref{insidesplit}) into the left hand side of (\ref{weakformulation}). The $P_{\gamma}f$ contribution and $a,c$ contributions vanish again due to oddness. Then the left hand side of (\ref{weakformulation}) becomes
		\begin{eqnarray}
		&&\int_{\partial \Omega \times \mathbf{R}^{3}}\{n\cdot  {v}\}|v|^{2}v_{i}v_{j}\sqrt{\mu }%
		\partial _{j}\phi _{b}^{i}f-\iint_{\Omega \times \mathbf{R}%
			^{3}}|v|^{2}v_{i}v_{j}\sqrt{\mu }\{\sum_{k}v_{k}\partial _{kj}\phi
		_{b}^{i}\}f  \notag \\
		&=&\int_{\partial \Omega \times \mathbf{R}^{3}} \{n\cdot  {v}\} |v|^{2} v_{i}v_{j}
		\sqrt{\mu }\partial_j\phi ^{i}_b[(1-P_{\gamma })f +r]\mathbf{1}_{\gamma _{+}}
		\label{stij_6} \\
		&&- {\iint_{\Omega \times \mathbf{R}^{3}}|v|^{2}v_{i}^{2}v_{j}^{2}%
			\mu \lbrack \partial _{ij}\phi _{b}^{i}b_{j}+\partial _{jj}\phi
			_{b}^{i}(x)b_{i}]}  \label{stij_8} \\
		&&-\iint_{\Omega \times \mathbf{R}^{3}}|v|^{2}v_{i}v_{j}v_{k}\sqrt{\mu }%
		\partial _{kj}\phi _{b}^{i}(x)[\mathbf{I}-\mathbf{P}]f.  \label{stij_9}
		\end{eqnarray}%
		Note that (\ref{stij_8}) is evaluated as
		\begin{equation*}
		7\sqrt{2\pi} \int_{\Omega }\{(\partial _{i}\partial _{j}\Delta
		^{-1}b_{i})b_{j}+(\partial _{j}\partial _{j}\Delta ^{-1}b_{i})b_{i}\}.
		\end{equation*}%
		Furthermore, by $(\ref{jb})$, $|\partial _{j}\phi _{b}^{i}|_{2}\lesssim
		\|\phi _{b}^{i}\|_{H^{2}}\lesssim \|b\|_{2}$, so that
		\begin{equation*}
		(\ref{stij_6})+(\ref{stij_9})\lesssim \|b\|_{2}\big\{|(1-P_{\gamma
		})f |_{2,+} +|r|_2 +\|(\mathbf{I}-\mathbf{P})f\|_{2} \big\}.
		\end{equation*}
		Combining (\ref{stiijjkiki}), we have the following estimate for $i\neq j,$
		\begin{eqnarray}
		&&\left|\int_{\Omega} \partial _{j}\partial _{j}\Delta ^{-1}b_{i}b_{i}\right|\   \notag \\
		&\lesssim &\ \left|\int_{\Omega} \partial _{i}\partial _{j}\Delta
		^{-1}b_{i}b_{j}\right|+ |(1-P_{\gamma })f |_{2,+}^{2}+ \|(\mathbf{I}-\mathbf{P})f\|^{2}_{2}+\|g\|_{2}^{2}+ |r|_2^2 +\varepsilon \|b\|_2^2  \notag \\
		&\lesssim &  |(1-P_{\gamma })f |_{2,+}^{2} +\|(\mathbf{I}-\mathbf{P})f\|^{2}_{2}+\|g\|_{2}^{2}+ |r|_2^2 +\varepsilon \|b\|_2^2.
		\label{stijijiijj}
		\end{eqnarray}%
		Moreover, by (\ref{stiijjkiki}), for $i=j=1,2,\dots,d,$
		\begin{eqnarray}
		&& \ \ \left|\int_{\Omega} \partial _{j}\partial _{j}\Delta ^{-1}b_{j}b_{j}\right| \   \notag \\
		&&  \lesssim \ |(1-P_{\gamma })f |_{2,+}^{2} +\|(\mathbf{I}-\mathbf{P})f\|^{2}_{2}+\|g\|_{2}^{2}+ |r|_2^2 +\varepsilon \|b\|_2^2.\label{stjjjj}
		\end{eqnarray}
		Combining (\ref{stijijiijj}) and (\ref{stjjjj}), we sum over $j=1,2,\dots,d$, to obtain, for all $i=1,2,\dots,d,$
		\begin{eqnarray}
		\| b_i \|_2 \lesssim \ |(1-P_{\gamma })f |_{2,+}^{2} +\|(\mathbf{I}-\mathbf{P})f\|^{2}_{2}+\|g\|_{2}^{2}+ |r|_2^2.\label{b_i}
		\end{eqnarray}
		
		\noindent{\it Step 3.} {Estimate of} \ ${a}$
		
		The estimate for $a$ is more delicate because it requires the zero mass condition
		$$\iint_{\Omega\times\mathbf{R}^3} f \sqrt{\mu}dx dv= 0=\int_{\Omega} a dx.$$  We choose as test function%
		\Be
		\psi =\psi_a  \equiv  (|v|^{2}-\beta_{a} )v\cdot\nabla_x\phi _{a}\sqrt{\mu }
		= \sum_{i=1}^d(|v|^{2}-\beta_{a} )v_{i}\partial _{i}\phi _{a}\sqrt{\mu } , \label{phia}
		\Ee
		where
		\begin{equation*}
		-\Delta_x \phi _{a}(x)=   a(x) ,\text{ \ \ \ \   }\frac{%
			\partial }{\partial n}\phi_a|_{\partial\Omega} =0.
		\end{equation*}%
		It follows from the elliptic estimate with $\int_{\Omega}a =0$ that we have
		\begin{equation*}
		\|\phi \|_{H^{2} }\lesssim \|a\|_{2 }.
		\end{equation*}%
		Since $\int_{\mathbf{R}^3}(\frac{|v|^2}{2}-\frac{3}{2})(v_i)^2 \mu(v) dv\neq 0$, we can choose $\beta_a>0 $ so that, for all $i,$
		\begin{equation}
		\int_{\mathbf{R}^{3}}(|v|^{2}-\beta_a )(\frac{|v|^{2}}{2}-\frac{3}{2}%
		)(v_{i})^{2}\mu(v) =0.\label{betaalpha}
		\end{equation}%
		Since $\mu(v)=\frac{1}{2\pi}e^{-\frac{|v|^2}{2}}$, the desired value is $\beta_a =10$.
		Plugging $\psi_a $ into (\ref{weakformulation}) and its right hand side is again bounded by%
		\begin{equation*}
		\|a\|_{2}\big\{\|(\mathbf{I}-\mathbf{P})f\|_{2}+\|g\|_{2}\big\}.
		\end{equation*}%
		By (\ref{bsplit}) and (\ref{insidesplit}), since the $c$ contribution vanishes in (\ref{weakformulation}) due to our choice of $\beta_a $ and the $b$ contribution vanishes in (\ref{weakformulation}) due to the oddness, the right hand side of (\ref{weakformulation}%
		) takes the form of
		\begin{eqnarray}
		&&\sum_{i=1}^d\int_{\gamma } \{ n \cdot  {v}\} (|v|^{2}-\beta_{a} )v_{i}\sqrt{\mu }\partial _{i}\phi
		_{a}(x) [P_{\gamma }f+(I-P_{\gamma })f\mathbf{1}_{\gamma _{+}}+r%
		\mathbf{1}_{\gamma _{+}}] \ \ \ \ \ \ \label{aboundary} \\
		&&-\sum_{i,k=1}^d \iint_{\Omega\times\mathbf{R}^3} (|v|^{2}-\beta_a )v_{i}v_{k}\partial _{ik}\phi _{a}(x)a(x)\mu(v) \label{abulk}\\
		&&-\sum_{i,k=1}^d\iint_{\Omega\times\mathbf{R}^3} (|v|^{2}-\beta_a )v_{i}v_{k}\partial _{ik}\phi _{a}(x)(\mathbf{I}-%
		\mathbf{P})f.\label{a2bulk}
		\end{eqnarray}%
		We make an orthogonal decomposition at the boundary,
		\begin{equation*}
		{v}_i=(  {v}\cdot n)n_i +  ({v}_{\perp })_i = v_n n_i + ({v}_{\perp })_i .
		\end{equation*}%
		The contribution of $P_{\gamma}f=z_{\gamma}(x)\sqrt{\mu}$ in (\ref{aboundary}) is
		\begin{eqnarray*}
			&&\int_{\gamma}(|v|^{2}-\beta_a ) {v}\cdot \nabla_x \phi _{a}(x) {v}_{n}\mu(v) z_{\gamma}(x)  \\
			&=&\int_{\gamma }(|v|^{2}-\beta_a ) {v}_{n} \frac{\partial \phi _{a}}{%
				\partial n} {v}_{n}\mu(v) z_{\gamma}(x) \\
			&&+ \int_{\gamma }(|v|^{2}-\beta_a ) v_{\bot}\cdot \nabla_x \phi_a  {v}_{n}\mu(v) z_{\gamma}(x) .
		\end{eqnarray*}%
		The crucial choice of (\ref{betaalpha}) makes the first term vanish due to the Neumann boundary condition, while the second term also
		vanishes due to the oddness of $ ({v}_{\bot})_i  {v}_n$ for all $i$.
		Therefore, (\ref{aboundary}) and (\ref{a2bulk}) are bounded by%
		\begin{equation*}
		\|a\|_{2 }\big\{\|(\mathbf{I}-%
		\mathbf{P})f\|_{2}+ |(1-P_{\gamma })f |_{+ } +|r|_{2} \big\}.
		\end{equation*}%
		The second term (\ref{abulk}), for $k\neq i$ vanishes due to the oddness. Hence we only have the $k=i$ contribution:
		\begin{equation*}
		\sum_{i=1}^d\iint_{\Omega \times \mathbf{R}^{3}}(|v|^{2}-\beta_a )(v_{i})^{2}\mu \partial
		_{ii}\phi _{a}a.
		\end{equation*}%
		Using $-\Delta_x \phi_{a}=a$ we obtain
		\begin{equation}
		\|a\|_{2 }^{2}\lesssim \|(\mathbf{I}-\mathbf{P})f\|_{2}^{2}+|(1-P_{\gamma
		})f |_{2,+}^{2}+|r|_{2}^{2}+\|g\|_{2}^{2}.\label{aestimate}
		\end{equation}
		
		\noindent {\it Step 1.} Estimate of  ${ \nabla_x  \Delta_N^{-1} \partial_{t}a   =   \nabla_x \partial_t \phi_a}$ in (\ref{timeweak})

		In the weak formulation(with time integration over $[t,t+\varepsilon]$), if we choose the test function $\psi=\phi \sqrt{\mu }$ with $\phi(x)$ dependent only of $x$, then we get
		(note that $L f$ and $g$ against $\phi(x) \sqrt{\mu }$
		are zero)
		\begin{eqnarray*}
			\sqrt{2\pi}\int_{\Omega }[a(t+\varepsilon )-a(t)]\phi (x)
			= 2\pi\sqrt{2\pi}\int_{t}^{t+\varepsilon }\int_{\Omega }(b\cdot \nabla _{x})\phi (x)
			+\int_{t}^{t+\varepsilon }\int_{\gamma_-}  r \phi\sqrt{\mu},
		\end{eqnarray*}%
		where $\int_{\mathbf{R}^3}\mu(v)dv=\sqrt{2\pi}, \   \int_{\mathbf{R}^{3}}(v_{1})^{2}\mu (v)dv=2\pi\sqrt{2\pi}$ and we have used the splitting (\ref{bsplit}) and (\ref{insidesplit}). Taking difference
		quotient, we obtain, for almost $t$,
		\begin{equation*}
		\int_{\Omega }\phi \partial _{t}a=\sqrt{2\pi}\int_{\Omega }(b\cdot \nabla_x )\phi
		+\frac{1}{{2\pi}}\int_{\gamma_-}  r  \phi \sqrt{\mu}.
		\end{equation*}
		Notice that, for $\phi =1$, from (\ref{dlinearcondition}) the right hand side of the above equation is zero so that, for all $t> 0 $,%
		\begin{equation*}
		\int_{\Omega }\partial _{t}a(t )dx =0.
		\end{equation*}
		On the other hand, for all $\phi(x)
		\in H^{1}\equiv H^{1}(\Omega )$, we have, by the trace theorem $|\phi|_2 \lesssim \|\phi \|_{H^1}$,
		\begin{eqnarray*}
			\left|\int_{\Omega }\phi (x)\partial _{t}adx\right| &\lesssim & |r|_{2} |\phi  |_{ 2 }+\|b\|_{2 }\|\phi
			\|_{H^{1} } \\
			&\lesssim & \{ \ \|b(t)\|_{2}  +|r|_{2} \} \ \|\phi \|_{H^{1}}.
		\end{eqnarray*}%
		Therefore we conclude that, for all $t>0$,
		\begin{equation*}
		\|\partial _{t}a(t)\|_{(H^{1})^{\ast }}\ \lesssim \
		\|b(t)\|_{2} +|r|_{2},
		\end{equation*}%
		where $(H^{1})^{\ast }\equiv(H^{1}(\Omega ))^{\ast }$ is the dual space of $H^{1}(\Omega )$
		with respect to the dual pair $$\langle A,B \rangle=\int_{\Omega }A(x)B(x)dx,$$ for $A\in H^{1}$
		and $B\in (H^1)^{\ast}$.
		
		By the standard elliptic theory, we can solve the Poisson equation with the Neumann
		boundary condition
		\begin{equation}
		-\Delta \Phi_a =\partial _{t}a(t)\ ,\ \ \ \ \ \frac{\partial \Phi_a }{\partial n}%
		\Big{|}_{\partial \Omega }=0, \notag 
		\end{equation}%
		with the crucial condition $\int_{\Omega }\partial _{t}a(t,x)dx=0$ for all $t>0$.   Notice that $\Phi_a =-\Delta _{N}^{-1}\partial _{t}a=\partial_t \phi_a$ where $\phi_a$ is defined in (\ref{phia}).
		Moreover we have
		\begin{eqnarray*}
			\|\nabla_x \partial_t \phi_a \|_{2} &=& \|\Delta _{N}^{-1}\partial
			_{t}a(t)\|_{H^{1}}=\|\Phi_a \|_{H^1}\\
			&\lesssim &\|\partial _{t}a(t)\|_{(H^{1})^{\ast }}
			\lesssim  \|b(t)\|_{2}+ |r|_{2}.
		\end{eqnarray*}%
		Therefore we conclude, for almost all $t> 0,$
		\begin{equation}
		\|\nabla_x \partial _{t}\phi _{a}(t)\|_{2 }\lesssim
		\|b(t)\|_{2 } +|r|_{2}.  \label{-1a}
		\end{equation}%

		\vskip .3cm
		\noindent{\it Step 2.} {Estimate of} $\nabla_x  \Delta^{-1}\partial _{t}{b}^{j} = \nabla_x \partial_t \phi_b^i$ in (\ref{timeweak})
		
		In (\ref{timeweak}) we choose a test function $\psi =\phi (x)v_{i}%
		\sqrt{\mu }$.  Since $\mu(v)=\frac{1}{2\pi}e^{-\frac{|v|^2}{2}}$,  $\int v_{i}v_{j}\mu(v)dv=\int
		v_{i}v_{j}(\frac{|v|^{2}}{2}-\frac{3}{2})\mu(v)dv=2\pi\sqrt{2\pi}\delta _{i,j}$ and we get
		\begin{eqnarray*}
			&&2\pi\sqrt{2\pi}\int_{\Omega }[b_i(t+\varepsilon )-b_i(t)]\phi\\
			&=&-\int_{t}^{t+\varepsilon }\int_{\gamma }f\phi v_{i}\sqrt{\mu }  +2\pi\sqrt{2\pi} \int_{t}^{t+\varepsilon
			}\int_{\Omega}\partial _{i}\phi [a+ c] \\
			&&+\int_{t}^{t+\varepsilon }\iint_{\Omega \times \mathbf{R}%
				^{3}}\sum_{j=1}^{d}v_{j}v_{i}\sqrt{\mu }\partial _{j}\phi (\mathbf{I}-\mathbf{P})f + \int_{t}^{t+\varepsilon} \iint_{\Omega\times\mathbf{R}^3} \phi v_i g \sqrt{\mu}.
		\end{eqnarray*}%
		Taking difference quotient, we obtain
		\begin{eqnarray*}
			&&\int_{\Omega }\partial _{t}b_{i}(t)\phi\notag\\
			&=& \frac{-1}{2\pi\sqrt{2\pi}}\int_{\gamma}f(t)v_{i}\phi \sqrt{\mu } + \int_{\Omega }\partial _{i}\phi [a(t)+ c(t)]\\
			&&+\frac{1}{2\pi\sqrt{2\pi}}\Big\{\iint_{\Omega
				\times \mathbf{R}^{3}}\sum_{j=1}^{d}v_{j}v_{i}\sqrt{\mu }\partial _{j}\phi(\mathbf{I}-\mathbf{P})f(t)+\iint_{\Omega\times\mathbf{R}^3} \phi v_i g(t) \sqrt{\mu}\Big\}.
		\end{eqnarray*}%
		For fixed $t> 0$, we choose $\phi=\Phi_b^i $ solving
		\begin{eqnarray*}
			-\Delta \Phi^{i}_{b}   = \partial _{t}b_{i}(t) , \ \ \ \ \
			\Phi^{i}_{b}|_{\partial \Omega }  = 0.
		\end{eqnarray*}%
		Notice that $\Phi_b^i = -\Delta^{-1}\partial_t b_i = \partial_t \phi_b^i$ where $\phi_b^i$ is defined in (\ref{jb}). The boundary terms vanish because of the Dirichlet
		boundary condition on $\Phi^{i}_{b}$.  Then we have, for $t\geq 0$,
		\begin{eqnarray*}
			\int_{\Omega }|\nabla_x \Delta ^{-1}\partial _{t}b_{i}(t)|^{2}dx
			&=&\int_{\Omega }|\nabla_x \Phi^{i}_{b} |^{2}dx=-\int_{\Omega }\Delta  \Phi^{i}_{b} \Phi^{i}_{b} dx \\
			&\lesssim &\varepsilon \{\|\nabla_x \Phi
			^{i}_{b}\|_{2}^{2} + \|  \Phi
			^{i}_{b}\|_{2}^{2} \}+\|a(t)\|_{2}^{2}+\|c(t)\|_{2}^{2}\\
			&& \ \   +\|(\mathbf{I}-\mathbf{P}%
			)f(t)\|_{2}^{2} +\|g(t)\|_2^2 \\
			&\lesssim &\varepsilon  \|\nabla_x \Phi
			^{i}_{b}\|_{2}^{2}  +\|a(t)\|_{2}^{2}+\|c(t)\|_{2}^{2}\\
			&& \ \   +\|(\mathbf{I}-\mathbf{P}%
			)f(t)\|_{2}^{2} +\|g(t)\|_2^2,
		\end{eqnarray*}%
		where we have used the Poincar$\acute{e}$ inequality.
		Hence, for all $t>0$
		\begin{eqnarray}
		&&  \ \  \|\nabla_x \partial_t \phi_b^i (t)\|_{2}=\|\nabla_x \Delta ^{-1}\partial _{t}b_{i}(t)\|_{2}\notag\\
		&&  \lesssim \ \
		\|a(t)\|_{2}+\|c(t)\|_{2}+\|(\mathbf{I}-%
		\mathbf{P})f(t)\|_{2}+\|g(t)\|_2.\label{-1b}
		\end{eqnarray}
		
		\vskip .3cm
		\noindent{\it Step 3.} {Estimate of  }$\nabla_x \Delta^{-1} \partial_t c= \nabla_x \partial _{t}\phi _{c} $ in (\ref{timeweak})
		
		In the weak formulation, we choose a test function $\phi (x)(\frac{|v|^{2}}{2}-%
		\frac{3}{2})\sqrt{\mu }$.  Since $\int_{\mathbf{R}^3}\mu (v)(\frac{|v|^{2}}{2}-\frac{3}{2})dv=0, \ \int_{\mathbf{R}^3} \mu(v) v_i v_j (\frac{|v|^2-3}{2})=2\pi\sqrt{2\pi}\delta_{i,j}$ and $%
		\int_{\mathbf{R}^3}\mu (v)(\frac{|v|^{2}}{2}-\frac{3}{2})^{2}dv=3\pi\sqrt{2\pi}$ we get
		\begin{eqnarray*}
			&&3\pi\sqrt{2\pi} \int_{\Omega} \phi(x)  c(t+\varepsilon,x) dx -3\pi\sqrt{2\pi} \int_{\Omega} \phi(x)  c(t,x) dx  \\
			& &=  2\pi\sqrt{2\pi}\int_{t}^{t+\varepsilon}\int_{\Omega }b\cdot \nabla_x \phi
			-\int_{t}^{t+\varepsilon}\int_{\gamma }(\frac{|v|^{2}}{2}-\frac{3}{2})\sqrt{\mu }%
			\phi f \\
			&& \ \ +\int_{t}^{t+\varepsilon}\iint_{\Omega \times \mathbf{R}^{3}}(\mathbf{I}-\mathbf{P})f(%
			\frac{|v|^{2}}{2}-\frac{3}{2})\sqrt{\mu }(v\cdot \nabla_x )\phi \\
			&& \ \ +\int_{t}^{t+\varepsilon}\iint_{\Omega \times \mathbf{R}^{3}} \phi g (%
			\frac{|v|^{2}}{2}-\frac{3}{2})\sqrt{\mu }
			,
		\end{eqnarray*}%
		and taking difference quotient, we obtain
		\begin{eqnarray*}
			&&\int_{\Omega} \phi(x)\partial_t c(t,x) dx\\
			&& =\frac{2}{3} \int_{\Omega} b(t)\cdot \nabla_x \phi -\frac{1}{3\pi\sqrt{2\pi}}\int_{\gamma} (\frac{|v|^2}{2}-\frac{3}{2})\sqrt{\mu} \phi f(t)\\
			&& \ \ + \frac{1}{3\pi\sqrt{2\pi}}\iint_{\Omega \times \mathbf{R}^{3}}(\mathbf{I}-\mathbf{P})f(t)(%
			\frac{|v|^{2}}{2}-\frac{3}{2})\sqrt{\mu }(v\cdot \nabla_x )\phi \\
			&& \ \ + \frac{1}{3\pi\sqrt{2\pi}}\iint_{\Omega \times \mathbf{R}^{3}} \phi g(t) (%
			\frac{|v|^{2}}{2}-\frac{3}{2})\sqrt{\mu }.
		\end{eqnarray*}
		For fixed $t> 0$, we define a test function $\phi=\Phi_c$ where $\phi_c$ is defined in (\ref{phic}). The boundary terms vanish because of the Dirichlet
		boundary condition on $\Phi_{c}$.  Then we have, for $t>0$,
		\begin{eqnarray*}
			-\Delta \Phi _{c}=\partial _{t}c(t) , \ \ \ \ \
			\Phi_{c}|_{\partial \Omega }=0.
		\end{eqnarray*}%
		Notice that $\Phi_c = -\Delta^{-1}\partial_t c(t) = \partial_t \phi_c(t)$ in (\ref{phic}).
		We follow the same procedure of estimates $\nabla_x \Delta ^{-1}\partial _{t}a$
		and $\nabla_x \Delta ^{-1}\partial _{t}b$ to have
		\begin{eqnarray*}
			&&\|\nabla_x \Delta ^{-1}\partial _{t}c(t)\|_{2}^{2}  =
			\int_{\Omega} |\nabla_x \Phi_c(x)|^2 dx\\
			&&= \int_{\Omega} \Phi_c(x) \partial_t c(t,x) dx\\
			&&\lesssim \
			\varepsilon \{ \|\nabla_x \Phi_c \|_2^2 + \|  \Phi_c \|_2^2 \}+ \|b(t)\|_{2}^2
			+ \|(\mathbf{I}-\mathbf{P})f(t)\|_2^2 + \|g(t)\|_2^2\\
			&&\lesssim \ \varepsilon   \|\nabla_x \Phi_c \|_2^2  + \|b(t)\|_{2}^2
			+ \|(\mathbf{I}-\mathbf{P})f(t)\|_2^2 + \|g(t)\|_2^2,
		\end{eqnarray*}%
		where we have used the Poincar$\acute{e}$ inequality. Finally we have, for all $t>0$,
		\begin{eqnarray}
		&&\|\nabla_x \partial _{t}\phi _{c}\|_{2}\lesssim \| \nabla_x \Delta^{-1} \partial_t c(t) \|_2\notag\\
		&& \lesssim \ \|b(t)\|_{2}+\|(%
		\mathbf{I}-\mathbf{P})f(t)\|_{2} + \|g(t)\|_2.  \label{-1c}
		\end{eqnarray}
		
		\vskip.3cm
		\noindent{\it Step 4.} {Estimate of } ${a,b,c}$ contributions in (\ref{timeweak})
		
		To estimate $c$ contribution in (\ref{timeweak}), we plug (\ref{phic}) into (\ref{timeweak}) to have from (\ref{insidesplit})
		\begin{eqnarray*}
			&& \int_{0}^{t}\iint_{\Omega \times \mathbf{R}^{3}}(|v|^{2}-\beta_c )v_{i}\sqrt{%
				\mu }\partial _{t}\partial _{i}\phi _{c}f  \\
			&=& \sum_{j=1}^{d}\int_{0}^{t}\iint_{\Omega \times \mathbf{R}^{3}}(|v|^{2}-\beta_c
			)v_{i}v_{j}\mu(v) \partial _{t}\partial _{i}\phi _{c}b_{j}\\
			& &+\int_{0}^{t}\iint_{\Omega \times \mathbf{R}^{3}}(|v|^{2}-\beta_c )v_{i} %
			\sqrt{\mu }\partial _{t}\partial _{i}\phi _{c} (\mathbf{I}-\mathbf{P})f .
		\end{eqnarray*}%
		The second line has non-zero contribution only for $j=i$ which leads to zero by the
		definition of $\beta_c $ in (\ref{beta}). We thus have from (\ref{-1c}) and Lemma \ref{steadyabc}
		\begin{eqnarray}
		&&\left|\int_{0}^{t}\iint_{\Omega \times \mathbf{R}^{3}}(|v|^{2}-\beta_c )v_{i}\sqrt{%
			\mu }\partial _{t}\partial _{i}\phi _{c}f \right| \notag\\
		&\lesssim&
		\int_{0}^{t}\Big\{\|b\|_{2}+\|(\mathbf{I}-\mathbf{P})f\|_{2} +\|g\|_2\Big\}\|(\mathbf{I}-%
		\mathbf{P})f\|_{2} \notag \\
		&\lesssim &\varepsilon \int_{0}^{t}\|b \|_{2}^2+ \int_{0}^{t}\Big\{ \|(\mathbf{I}%
		-\mathbf{P})f \|_{2}^{2} +\|g\|_2^2\Big\}\notag\\
		&\lesssim& \int_0^t \Big\{ \|(\mathbf{I}-\mathbf{P})f\|_{\nu}^2 + \|g\|_2^2 + |(1-P_{\gamma})f|_{2,+}^2 + |r|_2^2
		\Big\}
		.\label{ct}
		\end{eqnarray}
		Combined with Lemma \ref{steadyabc} we conclude
		\begin{eqnarray}
		\int_0^t \|c(s)\|_2^2 ds &\lesssim& G(t)-G(0)\label{timec}\\
		&+&\int_s^t \Big\{ \|(\mathbf{I}-\mathbf{P})f(s)\|_{\nu}^2 + \|g(s)\|_2^2 + |(1-P_{\gamma})f(s)|_{2,+}^2 + |r(s)|_2^2
		\Big\}ds.\notag
		\end{eqnarray}
		To estimate $b$ contribution in (\ref{timeweak}), we plug (\ref{phibj}) into (\ref{timeweak}) to have from (\ref{insidesplit}) %
		\begin{eqnarray}
		&& \int_{0}^{t}\iint_{\Omega \times \mathbf{R}^{3}}(v_{i}^{2}-\beta_b )\sqrt{%
			\mu }\partial _{t}\partial _{j}\phi _{b}^j f   \notag \\
		&=& {{\int_{0}^{t}\iint_{\Omega \times \mathbf{R}^{3}}(v_{i}^{2}- \beta_b ){\mu
				}\partial _{t}\partial _{j}\phi _{b}^{j} \{\frac{|v|^{2}}{2}-\frac{3}{2}\}c }}
		\label{i_7} \\
		&&+ \int_{0}^{t}\iint_{\Omega \times \mathbf{R}^{3}}(v_{i}^{2}- \beta_b ) %
		\sqrt{\mu }\partial_t\partial _{j}\phi _{b}^{j} (\mathbf{I}-\mathbf{P})f, \notag
		\end{eqnarray}%
		where we used (\ref{alpha}) to remove the $a$ contribution. We thus have from (\ref{-1b}) and Lemma \ref{steadyabc}
		\begin{eqnarray}
		&&\left|\int_{0}^{t}\iint_{\Omega \times \mathbf{R}^{3}}(v_{i}^{2}-\beta_b )\sqrt{\mu
		}\partial _{t}\partial _{j}\phi _{b}^jf \right|\notag\\
		&\lesssim
		&\int_{0}^{t}\Big\{\|a\|_{2}+\|c\|_{2}+\|(\mathbf{I}-\mathbf{P}%
		)f\|_{2} +\|g\|_2\Big\}\Big\{\|c\|_{2}+\|(\mathbf{I}-\mathbf{P})f\|_{2}\Big\}  \notag\\
		&\lesssim &
		\int_{0}^{t}\|\mathbf{P}f\|_{2}^{2}
		+ \|(\mathbf{I}-\mathbf{P})f\|_{2}^{2}
		+\|c\|_{2}^{2}\notag\\
		&\lesssim& \int_0^t \Big\{ \|(\mathbf{I}-\mathbf{P})f\|_{\nu}^2 + \|g\|_2^2 + |(1-P_{\gamma})f|_{2,+}^2 + |r|_2^2
		\Big\}
		.\label{bt1}
		\end{eqnarray}%
		Next we plug (\ref{phibij}) into (\ref{timeweak}) and we have from (\ref{-1b}) and Lemma \ref{timeweak}
		\begin{eqnarray}
		&&  \int_{0}^{t}\int_{\Omega \times \mathbf{R}^{3}} |v|^{2} v_{i}v_{j}%
		\sqrt{\mu }\partial _{t}\partial _{j}\phi _{b}^i f  \notag \\
		&=&  \int_{0}^{t}\int_{\Omega \times \mathbf{R}^{3}} |v|^{2} v_{i}v_{j}\sqrt{\mu }\partial_{t}\partial_{j}\phi _{b}^i(\mathbf{I}-\mathbf{%
			P})f  \notag \\
		&\lesssim &\int_{0}^{t}\{\|a\|_{2}+\|c\|_{2}+\|(\mathbf{I}-\mathbf{P}%
		)f\|_{2} +\|g\|_2 \}\|(\mathbf{I}-\mathbf{P})f\|_{2}  \notag \\
		&\lesssim &\int_0^t \Big\{ \|(\mathbf{I}-\mathbf{P})f\|_{\nu}^2 + \|g\|_2^2 + |(1-P_{\gamma})f|_{2,+}^2 + |r|_2^2
		\Big\}
		.  \label{bt2}
		\end{eqnarray}
		Combining this with Lemma \ref{steadyabc} we conclude
		\begin{eqnarray}
		\int_0^t \|b(s)\|_2^2 ds &\lesssim& G(t)-G(0)\label{timeb}\\
		&+&\int_0^t \Big\{ \|(\mathbf{I}-\mathbf{P})f(s)\|_{\nu}^2 + \|g(s)\|_2^2 + |(1-P_{\gamma})f(s)|_{2,+}^2 + |r(s)|_2^2
		\Big\}ds.\notag
		\end{eqnarray}
		
		Finally in order to estimate $a$ contribution in (\ref{timeweak}) we plug (\ref{phia}) for into (\ref{timeweak}). We estimate%
		\begin{eqnarray}
		&& \int_{0}^{t}\int_{\Omega\times\mathbf{R}^3} (|v|^{2}-\beta_a )v_{i}\mu \partial _{t}\partial _{i}\phi
		_{a}f  \label{at} \\
		&=& \int_{0}^{t}\int_{\Omega\times\mathbf{R}^3}  (|v|^{2}- \beta_a )(v_{i})^{2}\mu \partial _{t}\partial
		_{i}\phi _{a} b_{i}  \notag \\
		&& +\int_{0}^{t}\int_{\Omega\times\mathbf{R}^3}  (|v|^{2}- \beta_a )v_{i}\mu \partial _{t}\partial
		_{i}\phi _{a} (\mathbf{I}-\mathbf{P})f  \notag \\
		&\lesssim &\int_{0}^{t}\{\|b(t)\|_{2,\Omega}+ |r|_{2}\}\{\|b\|_{2,\Omega}+\|(\mathbf{I}-\mathbf{P})f\|_{2}\}  \notag
		\\
		&\lesssim &\int_0^t \Big\{ \|(\mathbf{I}-\mathbf{P})f\|_{\nu}^2 + \|g\|_2^2 + |(1-P_{\gamma})f|_{2,+}^2 + |r|_2^2
		\Big\} .
		\notag
		\end{eqnarray}%
		Combining this with Lemma \ref{steadyabc} we conclude
		\begin{eqnarray}
		\int_0^t \|a(s)\|_2^2 ds &\lesssim& G(t)-G(0)\label{timea}\\
		&+&\int_0^t \Big\{ \|(\mathbf{I}-\mathbf{P})f(s)\|_{\nu}^2 + \|g(s)\|_2^2 + |(1-P_{\gamma})f(s)|_{2,+}^2 + |r(s)|_2^2
		\Big\}ds.\notag
		\end{eqnarray}
		From (\ref{timec}), (\ref{timeb}) and (\ref{timea}), we proved the lemma.
	\end{proof}\unhide

	\hide
	
	\section{$L^2$ coercivity}

	The main purpose of this section is to prove the following:
	
	\begin{proposition}
		\label{dlinearl2}Suppose $(f, \phi)$ solves (\ref{eqtn_f}), (\ref{phi_f}), (\ref{BC_f}), and (\ref{phi_BC}). Then 
		%
		%
		there is $0<\lambda \ll 1$ such that for $0 \leq s \leq t$, 
		\Be\begin{split}\label{completes_dyn}
			& \| e^{\lambda t}f(t)\|_2^2
			+ \| e^{\lambda t} \nabla \phi (t) \|_2^2\\
			&
			+  \int_s^t \| e^{\lambda \tau}  f (\tau)\|_\nu^2 
			+ \| e^{\lambda \tau} \nabla \phi_f(\tau)\|_2^2 
			\mathrm{d} \tau 
			+  \int_s^t | e^{\lambda \tau} f |^2_{\gamma, 2 }    \\
			\lesssim & \ \|e^{\lambda s} f(s)\|_2^2 +   \|e^{\lambda s} \nabla \phi_f(s)\|_2^2  \\&   + \sup_{s \leq \tau \leq t} \| w f (\tau) \|_\infty \int_s^t \| e^{\lambda \tau}  f (\tau)\|_\nu^2  .  
		\end{split}\Ee
	\end{proposition}
	
	In order to prove the proposition we need the following:
	
	\begin{lemma}
		\label{dabc}
		There exists
		a function $G(t)$ such that, for all $0\le s\leq t$, $G(s)\lesssim
		\|f(s)\|_{2}^{2}$ and 
		\begin{equation}\begin{split}\label{estimate_dabc}
		&\int_{s}^{t}\|\mathbf{P}f(\tau)\|_{\nu }^{2}  + \int^t_s \| \nabla \phi_f \|_2^2\\
		\lesssim & \ 
		G(t)-G(s)
		+  \int_{s}^{t}\|(\mathbf{I}-\mathbf{P}%
		)f(\tau)\|_{\nu }^{2} +\int_{s}^{t}|(1-P_{\gamma })f(\tau)|_{2,+}^{2} \\
		& +\int_{s}^{t}\| \nu^{-1/2} { {\Gamma(f,f)} } \|_{2}^{2} + \int_{s}^{t} \|wf(\tau)\|^{2}_{\infty} \|\mathbf{P}f(\tau)\|_{2}^{2} .
		\end{split}\end{equation}
	\end{lemma}
	
	\begin{proof}[Proof of Proposition \ref{dlinearl2}] 
		\textit{Step 1. }  Without loss of generality we prove the result with $s=0$. We have an $L^2$-estimate from $e^{\lambda t} \times (\ref{eqtn_f})$
		\Be\notag
		\begin{split}
			&\|e^{\lambda t}  f(t) \|_2^2 - \| f(0) \|_2^2 +  \int^{t}_{0} | e^{\lambda \tau }
			(1-P_{\gamma} ) f^{j}|_{2,+}^{2}
			\\
			& + \int_0^t  \iint_{\O \times \R^3} v \cdot \nabla \phi_f e^{2\lambda \tau} |f|^2+
			2\int^t_0 \iint_{\O \times \R^3}e^{\lambda \tau}2  f Lf\\
			= &\ 2\int^t_0 \iint_{\O \times \R^3} e^{2\lambda \tau}f \Gamma(f,f) - 2\int^t_0 e^{2\lambda \tau}\int_{\O  } \nabla \phi_f \cdot \int_{\R^3} v \sqrt{\mu} f\\
			& 
			+2 \lambda  \int_0^t \|e^{\lambda \tau}  f(\tau) \|_2^2
			.
		\end{split}
		\Ee 
		
		On the other hand multiplying $\sqrt{\mu(v)} \phi_f(t,x)$ with a test function $\psi(t,x)$ to (\ref{eqtn_f}) and applying the Green's identity, we obtain
		\Be\begin{split}\notag
			&
			\int_{\O} \nabla \phi_f(t,x)  \cdot \int_{\R^3}  v \sqrt{\mu} f \dd v\dd x
			\\
			= & \ \int_{\O} \phi_f(t,x) \p_\tau\left(\int_{\R^3} f (\tau) \sqrt{\mu} \dd v\right)  \dd x
			+ \iint_{\p\O \times \R^3}  \phi_f(t,x) f \sqrt{\mu} \{n \cdot v\} \dd v\dd S_x .
		\end{split}\Ee
		From (\ref{null_flux}), the last boundary contribution equals zero. Now we use (\ref{phi_f}) and (\ref{phi_BC}) and deduce that 
		\Be
		\begin{split}\notag
			& \int^t_0e^{2 \lambda \tau}\int_{\O} \phi_f(t,x) \p_\tau\left(\int_{\R^3} f(\tau)  \sqrt{\mu} \dd v\right)  \dd x \dd \tau\\
			= & \ -  \int^t_0e^{2 \lambda \tau}\int_{\O} \phi_f(t,x) \p_\tau  \Delta_x  \phi_f(\tau,x) \dd x \dd \tau\\
			= &  \  \frac{1}{2}\int^t_0e^{2 \lambda \tau}\int_{\O}  \p_\tau  |\nabla_x  \phi_f(\tau,x)|^2 \dd x \dd \tau
			\\
			= & \ \frac{1}{2}   \left(\int_{\O} e^{2 \lambda t} |\nabla_x \phi_f (t,x)|^2  \dd x\right)- \frac{1}{2}   \left(\int_{\O} |\nabla_x \phi_f (0,x)|^2  \dd x\right)\\
			& \ 
			- \lambda \int^t_0 e^{2 \lambda \tau} \int_\O |\nabla_x  \phi_f(\tau,x)|^2 \dd x \dd \tau
			.
		\end{split}
		\Ee
		
		Hence we derive 
		\Be\notag
		\begin{split}
			&\| e^{\lambda t} f(t) \|_2^2 + \| e^{\lambda t}\nabla\phi_f (t) \|_2^2  + \int_0^t  \iint_{\O \times \R^3}
			e^{2\lambda \tau}
			v \cdot \nabla \phi_f |f|^2\\
			& +
			2C\int^t_0 \iint_{\O \times \R^3}\| e^{\lambda \tau} (\mathbf{I} - \mathbf{P}) f \|_\nu^2
			+ \int^{t}_{0} | e^{\lambda \tau }
			(1-P_{\gamma} ) f^{j}|_{2,+}^{2}
			\\
			\lesssim &\ 
			\| f(0) \|_2^2+ \| \phi_f (0) \|_2^2 +
			\int^t_0 \| e^{\lambda \tau}  \nu^{-1/2} \Gamma(f,f) \|_2^2\\
			& 
			+\{ \lambda + o(1)\} \int^t_0 \| e^{\lambda \tau} f \|_\nu^2 + \lambda \int^t_0 \| e^{\lambda \tau} \nabla_x \phi_f \|_2^2
			.
		\end{split}
		\Ee
		Now we apply Lemma \ref{dabc} and add $o(1) \times (\ref{estimate_dabc})$ to the above inequality and choose $0< \lambda \ll 1$ to conclude (\ref{completes_dyn}) except the full boundary control. 
		
		\vspace{4pt}
		
		\textit{Step 2. } Note that $P_\gamma f = z(t,x)\sqrt{\mu(v)}$ for a suitable function $z(t,x)$ on the boundary. Then for $0 < \e \ll 1$
		\Be
		\begin{split}\notag
			| P_\gamma f |_{\gamma,2}^2 
			= & \ \int_{\p\O } |z(t,x)|^2 \dd x \times \int_{\R^3} \mu (v) |n(x) \cdot v| \dd v\\
			\lesssim & \  \int_{\p\O } |z(t,x)|^2 \dd x \times \int_{ \gamma_+(x) \backslash \gamma_+^\e(x)} \mu (v)^{3/2} |n(x) \cdot v| \dd v\\
			= & \  | \mathbf{1}_{ \gamma_+  \backslash \gamma_+^\e }\mu^{1/4} P_\gamma f   |_{ +,2}^2.
		\end{split}
		\Ee
		Since $P_\gamma f = f - (1- P_\gamma) f$ on $\gamma_+$ we have $
		| \mathbf{1}_{ \gamma_+  \backslash \gamma_+^\e } \mu^{1/4} P_\gamma f|_{ +,2}^2
		\lesssim |\mathbf{1}_{ \gamma_+  \backslash \gamma_+^\e }\mu^{1/4}  f |_{ +,2}^2 + |(1- P_\gamma) f|_{ +,2}^2.$ Therefore
		\Be\label{Pgamma_bound}
		\int^t_0 | P_\gamma f|_{\gamma,2}^2 \lesssim \int^t_0 |\mathbf{1}_{ \gamma_+  \backslash \gamma_+^\e }\mu^{1/4}  f |_{ +,2}^2 +\int^t_0 |(1- P_\gamma) f|_{ +,2}^2.
		\Ee
		Note that 
		\Be\begin{split}\notag
			&\big|[\p_t + v\cdot \nabla_x - \nabla \phi \cdot \nabla_v] (\mu^{1/4}  f )\big|\\
			\lesssim&\ \mu^{1/4} \{ |v|| \nabla_x \phi|  f  +  |v|| \nabla_x \phi|  +  |Lf| + |\Gamma (f,f) | \}
		\end{split}\Ee
		By the trace theorem Lemma \ref{le:ukai},  
		\Be\label{est_trans_f}
		\begin{split}
			& \int^t_0 |\mathbf{1}_{ \gamma_+  \backslash \gamma_+^\e }\mu^{1/4}  f |_{ +,2}^2\\
			\lesssim & \ \| f_0 \|_2 + (1 + \| wf \|_\infty) \int^t_0 \| f \|_2^2 + \int^t_0\| \nabla \phi \|_2^2.
		\end{split}
		\Ee
		Adding $o(1) \times (\ref{Pgamma_bound})$ to the result of \text{Step 1} and using (\ref{est_trans_f}) we conclude (\ref{completes_dyn}).\end{proof}


	\begin{proof}[Proof of Lemma \ref{dabc}]
		From the Green's identity (\ref{}), a solution $(f,\phi_f)$ satisfies
		\begin{equation} \label{weakformulation}
		\begin{split}
		&\iint_{\O \times \R^3} f(t) \psi(t) - \iint_{\O \times \R^3} f(s) \psi(s)
		- \underbrace{ \int_s^t \iint_{\O \times \R^3} f \p_t\psi }_{(T)} +
		\underbrace{ \int^t_s\int_{\gamma }\psi f (v\cdot n(x)) }_{(B)}  \\
		& - \underbrace{ \int^t_s\iint_{\Omega \times  \R^{3}} \mathbf{P}f {v}\cdot \nabla _{x}\psi  }_{(C)} - \int^t_s\iint_{\Omega \times  \R^{3}} (\mathbf{I-P})f {v}\cdot \nabla _{x}\psi + \underbrace{ \int_s^t \iint_{\O \times \R^3} \sqrt{\mu} f \nabla_x \phi_f \cdot \nabla_v \Big[\frac{1}{\sqrt{\mu}} \psi\Big] }_{(\ref{weakformulation})_{P}}
		\\
		&=\int^t_s\iint_{\Omega \times \R^{3}}\psi \{ -L  (\mathbf{I}-\mathbf{P%
		})f + \Gamma(f,f)\}- \underbrace{\int^t_s\iint_{\Omega \times  \R^{3}} \psi v\cdot \nabla_x \phi_f \sqrt{\mu}}_{(\ref{weakformulation})_{\phi_f}}  .  
		\end{split}
		\end{equation}%
		
		As \cite{EGKM, EGKM2} we use a set of test functions:
		\Be \label{tests}
		\begin{split}
			\psi_a  &\equiv   (|v|^{2}-\beta_{a} )\sqrt{\mu }v\cdot\nabla_x\phi _{a},  \\
			\psi^{i,j}_{b,1} &\equiv (v_{i}^{2}-\beta_ b)\sqrt{\mu }\partial _{j}\phi _{b}^{j}, \quad i,j=1,2,3,   \\
			\psi^{i,j}_{b,2} &\equiv |v|^{2}v_{i}v_{j}\sqrt{\mu }\partial _{j}\phi _{b}^{i}(x),\quad i\neq j,  \\
			\psi_c &\equiv (|v|^{2}-\beta_c )\sqrt{\mu }v \cdot \nabla_x \phi _{c},  \\
		\end{split}
		\Ee
		where $\phi_{a}(t,x)$, $\phi_{b}(t,x)$, and $\phi_{c}(t,x)$ solve
		\Be\begin{split}\label{phi_abc}
			- \Delta \phi_a &= a(t,x),  \quad \p_{n}\phi_{a} \vert_{\p\O} = 0,  \\
			- \Delta \phi_b^j &= b_j(t,x),   \ \ \phi^j_b|_{\p\O} =0,\  \text{and}  \ - \Delta \phi_c  = c(t,x),   \ \ \phi_c|_{\p\O} =0,
		\end{split}
		\Ee
		and $\beta_a=10$, $\beta_b=1$, and $\beta_c= 5$ such that for all $i=1,\notag2,3,$
		\Be \label{defbeta}
		\begin{split}
			&\int_{{\R}^{3}}(|v|^{2}-\beta_a )(\frac{|v|^{2}}{2}-\frac{3}{2}%
			) v_{i} ^{2}\mu(v)\dd v =0 ,\\
			&\int_{\R^3} (v_i^2 - \beta_b) \mu(v) \dd v=0 ,\\
			&\int_{\R^3} (|v|^{2}-\beta_c )v_{i}^{2}\mu(v) \dd v=0. 
		\end{split}
		\Ee
		
		\vspace{4pt}

		\textit{Step 1. } Estimate of $(\ref{weakformulation})_{\phi_f}$: From (\ref{tests}) and (\ref{defbeta}), we have $(\ref{weakformulation})_{\phi_f}\equiv0$ for $\psi_{b,1}^{i,j}$, $\psi_{b,2}^{i,j}$, and $\psi_c$. For $\psi = \psi_a$, $(\ref{weakformulation})_{\phi_f}$ equals
		\Be \label{extra}
		\begin{split}
			(\ref{weakformulation})_{\phi_f} \big\vert_{\psi=\psi_{a}} &= \int_{\R^3} (|v|^2- \beta_a) (v_1)^2\mu \dd v \int^t_s\int_{\O  } \nabla \phi_a \cdot \nabla \phi_f \\
			&= C \int^t_s \| \nabla \phi_f \|_2^2 ,  \\
			(\ref{weakformulation})_{\phi_f} \big\vert_{\psi=\psi_{b} \ \psi_{c}} &= 0,  \\
		\end{split}\Ee
		because $\phi_a = \phi_f$ from the definitions of (\ref{phi_f}) and (\ref{phi_abc}). \\
		
		Estimate of $(\ref{weakformulation})_{P}$: From (\ref{tests}),
		\Be \label{extraP}
		\begin{split}
			(\ref{weakformulation})_{P} &= \int_s^t \iint_{\O \times \R^3} \sqrt{\mu} f \nabla_x \phi_f \cdot \nabla_v \Big[\frac{1}{\sqrt{\mu}} \psi\Big],\quad  \psi = \psi_{a,b,c},  \\
			&\lesssim \int_{s}^{t} \|wf\|_{\infty} \int_{\O}  \nabla_{x}\phi_{f}\cdot\nabla_{x}\phi_{a,b,c}  \lesssim \int_s^t \|wf(\tau)\|_{\infty} \|\mathbf{P}f(\tau)\|_{2}^{2}, 
		\end{split}
		\Ee
		by elliptic estimate $\|\nabla_{x}\phi_{a,b,c}\|_{2} \lesssim \|\phi_{a,b,c}\|_{H^{2}} \lesssim \|\mathbf{P}f\|_{2} $.  \\
		
		{\it Step 2}. {Estimate of} \ ${c}$ : We apply $\phi_{c}$ to (\ref{weakformulation}). From oddness in velocity integration and (\ref{defbeta}), $(C)$ becomes,
		\Be \label{cC}
		\begin{split}
			\int^t_s\iint_{\Omega \times  \R^{3}} \mathbf{P}f {v}\cdot \nabla
			_{x}\psi_{c} &= C_{1} \int^t_s \|c(\tau)\|_{2}^{2},
		\end{split}	
		\Ee
		where $\int_{\mathbf{R}^3}(|v|^{2}-\beta_c )v_{i}^{2} (\frac{|v|^{2}}{2}-\frac{3}{2})\mu(v)dv= 10\pi\sqrt{2\pi}$. For boundary integral $(B)$, we decompose $f_{\gamma } = P_{\gamma }f +\mathbf{1}_{\gamma_{+}}(1-P_{\gamma })f$. Then from (\ref{defbeta}) and trace theorem $|\nabla\phi_{c}|_{2} \lesssim \|\phi_{c}\|_{H^{2}} \lesssim \|c\|_{2}$,
		\Be \label{cB}
		\begin{split}
			&\int^t_s\int_{\gamma }\psi_{c} f (v\cdot n(x)) = \int^t_s\int_{\gamma }\psi_{c} \mathbf{1}_{\gamma_{+}} (1-P_{\gamma})f d\gamma  \\
			&\lesssim \varepsilon\int^t_s \|c(\tau)\|_{2}^{2} + C_{\varepsilon}\int^t_s |(1-P_{\gamma})f(\tau)|^{2}_{2,+} ,\quad \varepsilon \ll 1. \\
		\end{split}
		\Ee
		If we define
		\Be \label{Re}
		\begin{split}
			Re &:= (\ref{weakformulation})_{P} + (\ref{weakformulation})_{\phi_{f}} +  \int^t_s\iint_{\Omega \times \R^{3}}\psi \{ L  (\mathbf{I}-\mathbf{P%
			})f - \Gamma(f,f)\}   \\
			&\quad - \int^t_s\iint_{\Omega \times  \R^{3}} (\mathbf{I-P})f {v}\cdot \nabla _{x}\psi,
		\end{split}
		\Ee
		\Be \label{c Re}
		\int_{s}^{t} \iint_{\O\times \R^{3}} \psi_{c} Re\vert_{\psi_{c}} \lesssim \varepsilon \int_{s}^{t} \|c\|_{2}^{2} + \int_{s}^{t} \|(\mathbf{I-P})f(\tau)\|_{\nu}^{2} + \int_{s}^{t} \|\nu^{-1/2}\Gamma(f,f)(\tau)\|_{2}^{2},
		\Ee
		from elliptic estimate and Young's inequality. We also use even/oddness in velocity integration , (\ref{defbeta}), and Young's inequality to estimate,
		\Be \label{cT}
		\begin{split}
			(T)\vert_{\psi_{c}} &= \int_s^t \iint_{\O \times \R^3} f \p_t\psi_{c} = \int_s^t \iint_{\O \times \R^3} (\mathbf{I-P})f \p_t\psi_{c}  \\
			&\lesssim \varepsilon \int_s^t \|\nabla\Delta^{-1}\p_{t}c(\tau)\|_{2}^{2} + \int_{s}^{t} \|(\mathbf{I-P})f(\tau)\|_{\nu}^{2}.
		\end{split}
		\Ee
		
		Now, we choose a new test function $\psi_{c}^{t} := (\frac{|v|^{2}}{2} - \frac{3}{2})\sqrt{\mu } \p_{t}\phi_{c} (t,x)$. Note that $\p_{t}\phi_{c}$ solves $-\Delta \p_{t}\phi_{c} = \p_{t}c(t,x)$ with $\p_{t}\phi_{c}(t,x)|_{\p\O} = 0$. We taking difference quotient for $\p_{t}f$ in (\ref{eqtn_f}) and it replace first three terms in the LHS of (\ref{weakformulation}). With help of Poincar\'e inequality $\|\p_{t}\phi_{c}\|_{2} \lesssim \|\nabla\p_{t}\phi_{c}\|_{2}$, we can also compute
		\Be \label{ct_extra}
		\begin{split}
			(\ref{weakformulation})_{\phi_f} \big|_{\psi=\psi_{c}^{t}} &= \int^t_s\iint_{\Omega \times  \R^{3}} (\frac{|v|^{2}}{2} - \frac{3}{2})\sqrt{\mu } \p_{t}\phi_{c} v\cdot \nabla_x \phi_f \sqrt{\mu}  \\
			&\lesssim \varepsilon \int^t_s \|\nabla\Delta^{-1}\p_{t}c(\tau)\|_{2}^{2} +\int^t_s \|a(\tau)\|_{2}^{2},	\\
			(\ref{weakformulation})_{P} \big|_{\psi=\psi_{c}^{t}} &= \int_s^t \iint_{\O \times \R^3} \sqrt{\mu} f \nabla_x \phi_f \cdot 2v \p_{t}\phi_{c}    \\
			&\lesssim \varepsilon \int^t_s \|\nabla\Delta^{-1}\p_{t}c(\tau)\|_{2}^{2} + \int^t_s \big(\|a(\tau)\|_{2}^{2} + \|b(\tau)\|_{2}^{2}\big) + \int_{s}^{t} \|(\mathbf{I-P})f\|_{\nu}^{2}	\\
		\end{split}
		\Ee
		Since $\psi_{c}^{t}$ vanishes when it acts with $Lf$ and $\Gamma(f,f)$, and boundary integral $(B)$ vanishes by Dirichlet boundary condition of $\phi_{c}$ , from (\ref{ct_extra}) and (\ref{weakformulation}), we obtain
		\begin{equation} \label{c time}
		\begin{split}
		&\int_s^t \int_{\Omega} \p_{t}\phi_{c}(\tau,x)\partial_t c(\tau,x) dx = \int_s^t \|\nabla\Delta^{-1} \p_{t} c(\tau)\|_{2}^{2}    \\
		&\lesssim \varepsilon \int^t_s \|\nabla\Delta^{-1}\p_{t}c(\tau)\|_{2}^{2} +\int^t_s \big( \|a(\tau)\|_{2}^{2} + \|b(\tau)\|_{2}^{2} \big)  + \int_{s}^{t} \|(\mathbf{I-P})f(\tau)\|_{\nu}^{2} . \\	
		\end{split}
		\end{equation}
		
		We combine (\ref{weakformulation}), (\ref{extra}), (\ref{extraP}), (\ref{cB}), (\ref{c Re}), (\ref{cT}), and (\ref{c time}) with $\varepsilon \ll 1$ to obtain
		\Be \label{c est}
		\begin{split}
			\int^t_s \|c(\tau)\|_{2}^{2} &\lesssim G_{c}(t) - G_{c}(s) + \int_{s}^{t} \|(\mathbf{I-P})f(\tau)\|_{\nu}^{2} + \int^t_s |(1-P_{\gamma})f(\tau)|_{2,+} \\
			&\quad + \int_{s}^{t} \|\nu^{-1/2}\Gamma(f,f)(\tau)\|_{2}^{2} + \int_{s}^{t} \|wf(\tau)\|^{2}_{\infty} \|\mathbf{P}f(\tau)\|_{2}^{2}  \\
			&\quad + \varepsilon\int^t_s \big(\|a(\tau)\|_{2}^{2} + \|b(\tau)\|_{2}^{2}\big),
		\end{split}
		\Ee
		for $\varepsilon \ll 1$ where $G_{c}(t) := \iint_{\O\times \R^{3}} f(t)\psi_{c}(t) \lesssim \|f(t)\|_{2}^{2}$.  \\
		
		{\it Step 3.} {Estimate of} \ ${a}$ : From mass conservation $\int_{\O} a(t,x) dv = 0$, $\phi_{a}$ in (\ref{phi_abc}) is well-defined. Moreover, we choose $\phi_{a}$ so that has mean zero, $\int_{\O} \phi_{a}(t,x) dx = 0$. Therefore, Poincar\'e inequality $\|\phi_{a}\|_{2} \lesssim \|\nabla \phi_{a}\|_{2}$ holds and these are also true for $\p_{t}\phi_{a}$ which solves same elliptic equation with Neumann boundary condition.  
		
		By even/oddness in velocity integral and $\b_{a}$ defined in (\ref{defbeta}), we can replace $c$ into $a$ in estimates (\ref{cC}) and (\ref{c Re}). For boundary integral $(B)$, we decompose $f_{\gamma } = P_{\gamma }f +\mathbf{1}_{\gamma_{+}}(1-P_{\gamma })f$. From Neumann boundary condition $\p_{n}\phi_{a} = 0$ and oddness in velocity integral, $\int_{\gamma} \psi_{a} P_{\gamma}f (v\cdot n(s)) = 0$ and we obtain similar esimate as (\ref{cB}),
		\Be \label{aB}
		\begin{split}
			&\int^t_s\int_{\gamma }\psi_{a} f (v\cdot n(x)) = \int^t_s\int_{\gamma }\psi_{a} \mathbf{1}_{\gamma_{+}} (1-P_{\gamma})f d\gamma  \\
			&\lesssim \varepsilon\int^t_s \|a(\tau)\|_{2}^{2} + C_{\varepsilon}\int^t_s |(1-P_{\gamma})f(\tau)|^{2}_{2,+} ,\quad \varepsilon \ll 1. \\
		\end{split}
		\Ee
		For $(T)$, from oddness,
		\Be \label{aT}
		\begin{split}
			(T)\vert_{\psi_{a}} &= \int_s^t \iint_{\O \times \R^3} f \p_t\psi_{a} = \int_s^t \iint_{\O \times \R^3} [ \sqrt{\mu} (b\cdot v) + (\mathbf{I-P})f] \p_t\psi_{a}  \\
			&\lesssim \varepsilon \int_s^t \|\nabla\Delta^{-1}\p_{t}a(\tau)\|_{2}^{2} + \int_{s}^{t} \|b(\tau)\|_{2}^{2}+ \int_{s}^{t} \|(\mathbf{I-P})f(\tau)\|_{\nu}^{2}.
		\end{split}
		\Ee
		Now let us estimate $\int_{s}^{t} \|\nabla\Delta^{-1}\p_{t}a(\tau)\|_{2}^{2}$ which appear in (\ref{cT}) type estimate. We use new test function $\psi_{a}^{t} = \varphi(x)\sqrt{\mu}$. It easy to check 
		\Be \label{at_extra}
		\begin{split}
			(\ref{weakformulation})_{\phi_f} \big|_{\psi=\sqrt{\mu}} &= \int^t_s\iint_{\Omega \times  \R^{3}} \sqrt{\mu } v\cdot \nabla_x \phi_f \varphi(x) \sqrt{\mu} =0, 	\\
			(\ref{weakformulation})_{P} \big|_{\psi=\sqrt{\mu}} &= \int_s^t \iint_{\O \times \R^3} \sqrt{\mu} f \nabla_x \phi_f \cdot \nabla_v \varphi(x) = 0,  \\
		\end{split}
		\Ee
		and from taking difference quotient and null condition (\ref{null_flux}), we obtain, for almost $t$,
		\begin{equation*}
		\int_{\Omega } \varphi \partial _{t}a=\sqrt{2\pi}\int_{\Omega }(b\cdot \nabla_x )\varphi.
		\end{equation*}
		In particular, if we choose $\varphi = 1$, we directly get $\int_{\O} \p_{t} a = 0$, so Neumann problem 
		\begin{equation}
		-\Delta \Phi_a =\partial _{t}a(t)\ ,\ \ \ \ \ \frac{\partial \Phi_a }{\partial n}%
		\Big{|}_{\partial \Omega }=0, \notag 
		\end{equation}
		is well-defined. For dual pair $(H^{1})^{\ast }\equiv(H^{1}(\Omega ))^{\ast }$ with respect to $\langle A,B \rangle=\int_{\Omega } A \cdot B dx$ for $A\in H^{1}$ and $B\in (H^1)^{\ast}$,
		\begin{equation} \label{a time}
		\|\nabla_x \partial_t \phi_a \|_{2} = \|\Delta^{-1}\partial
		_{t}a(t)\|_{H^{1}}=\|\Phi_a \|_{H^1} \lesssim \|\partial _{t}a(t)\|_{(H^{1})^{\ast }}
		\lesssim  \|b(t)\|_{2}.  \\
		\end{equation}
		
		We change $c$ into $a$ in (\ref{cC}) and (\ref{c Re}) and combine with (\ref{weakformulation}), (\ref{extra}), (\ref{extraP}), (\ref{aB}), (\ref{aT}), and (\ref{a time}) with $\varepsilon \ll 1$ to obtain
		\Be \label{a est}
		\begin{split}
			&\int^t_s \|a(\tau)\|_{2}^{2} + \int^t_s \|\nabla \phi_{f}(\tau)\|_{2}^{2} \\
			&\lesssim G_{a}(t) - G_{a}(s) + \int_{s}^{t} \|(\mathbf{I-P})f(\tau)\|_{\nu}^{2} + \int^t_s |(1-P_{\gamma})f(\tau)|_{2,+} \\
			&\quad + \int_{s}^{t} \|\nu^{-1/2}\Gamma(f,f)(\tau)\|_{2}^{2} + \int_{s}^{t} \|wf(\tau)\|^{2}_{\infty} \|\mathbf{P}f(\tau)\|_{2}^{2}  + \int^t_s \|b(\tau)\|_{2}^{2} ,  \\
		\end{split}
		\Ee
		for $\varepsilon \ll 1$ where $G_{a}(t) := \iint_{\O\times \R^{3}} f(t)\psi_{c}(t) \lesssim \|f(t)\|_{2}^{2}$.  \\
		
		{\it Step 4.} {Estimate of} \ $ {b}$ : For fixed $i,j$, we choose test function $\psi = \psi_{b,1}^{i,j}$ in (\ref{tests}) where $\b_{b}$ and $\phi_{b}$ are defined in (\ref{defbeta}) and (\ref{phi_abc}). From oddness in velocity integration and definition of $\b_{b}$, $(C)$ in (\ref{weakformulation}) yields
		\Be \label{bC}
		\begin{split}
			(C)\vert_{\psi_{b,1}^{i,j}} &:= \int_s^t \iint_{\O\times \R^{3}} \mathbf{P}f v\cdot\nabla \psi_{b,1}^{i,j} = C_{3} \int_s^t \int_{\O} b_{i} (\p_{ij} \Delta^{-1} b_{j}), 
		\end{split}
		\Ee
		where $C_{3} := \int_{\R^{3}} (v_{i}^{2} - \b_{b}) v_{i}^{2} \mu dv = 2\sqrt{2\pi}$. For boundary integration, contribution of $P_{\gamma}f$ vanishes by oddness. 
		\Be \label{bB}
		\begin{split}
			(B)\vert_{\psi_{b,1}^{i,j}} &:= \int_s^t \int_{\gamma} \psi_{b,1}^{i,j} \mathbf{1}_{\gamma_{+}} (1-P_{\gamma})f (v\cdot n(x)) \lesssim \varepsilon \int_s^t \|b(\tau)\|^{2}_{2} + \int_s^t |(1-P_{\gamma})f|_{2,+}^{2},
		\end{split}
		\Ee
		and similar as (\ref{cT}) and (\ref{c Re}), we use oddness and definition of $\b_{b}$ to vanish contribution of $a$ and $b$. We obtain
		\Be \label{bT}
		\begin{split}
			(T)\vert_{\psi_{b,1}^{i,j}} &\lesssim  \varepsilon \int_s^t \|\nabla\Delta^{-1}\p_{t} b_{j}(\tau)\|_{2}^{2} + \int_{s}^{t} \|c(\tau)\|_{2}^{2} + \int_{s}^{t} \|(\mathbf{I-P})f(\tau)\|_{\nu}^{2},
		\end{split}
		\Ee
		\Be \label{b Re}
		\int_{s}^{t} \iint_{\O\times \R^{3}} \psi_{b,1}^{i,j} Re\vert_{\psi_{b,1}^{i,j}} \lesssim \varepsilon \int_{s}^{t} \|b(\tau)\|_{2}^{2} + \int_{s}^{t} \|(\mathbf{I-P})f(\tau)\|_{\nu}^{2} + \int_{s}^{t} \|\nu^{-1/2}\Gamma(f,f)(\tau)\|_{2}^{2},
		\Ee
		
		Next, we try test function $\psi_{b,2}^{i,j}$ with $i \neq j$ to obtain 
		\Be \label{bC2}
		\begin{split}
			(C)\vert_{\psi_{b,2}^{i,j}} &:= \int_s^t \iint_{\O\times \R^{3}} (b\cdot v)\sqrt{\mu} v\cdot\nabla \psi_{b,2}^{i,j} = C_{4} \int_s^t \int_{\O} \big(  b_{j} (\p_{ij} \Delta^{-1} b_{i}) + b_{i} (\p_{jj} \Delta^{-1} b_{i}) \big) 
		\end{split}
		\Ee
		by oddness in velocity integral where $C_{4} := 7\sqrt{2\pi}$. We also have the following three estimates using oddness of velocity integral,
		\Be \label{bB2T2Re2}
		\begin{split}
			(B)\vert_{\psi_{b,2}^{i,j}} &:= \int_s^t \int_{\gamma} \psi_{b,1}^{i,j} \mathbf{1}_{\gamma_{+}} (1-P_{\gamma})f (v\cdot n(x)) \lesssim \varepsilon \int_s^t \|b(\tau)\|^{2}_{2} + \int_s^t |(1-P_{\gamma})f|_{2,+}^{2},  \\
			(T)\vert_{\psi_{b,2}^{i,j}} &\lesssim  \varepsilon \int_s^t \|\nabla\Delta^{-1}\p_{t} b_{i}(\tau)\|_{2}^{2} + \int_{s}^{t} \|(\mathbf{I-P})f(\tau)\|_{\nu}^{2},  \\
			\int_{s}^{t} \iint_{\O\times \R^{3}} \psi_{b,2}^{i,j} Re\vert_{\psi_{b,2}^{i,j}} &\lesssim \varepsilon \int_{s}^{t} \|b(\tau)\|_{2}^{2} + \int_{s}^{t} \|(\mathbf{I-P})f(\tau)\|_{\nu}^{2} + \int_{s}^{t} \|\nu^{-1/2}\Gamma(f,f)(\tau)\|_{2}^{2}.
		\end{split}
		\Ee
		
		To obtain estimate for $\|\nabla\Delta^{-1} \p_{t} b_{j}\|_{2}$, we use a test function $\psi_{b,j}^{t} := v_{j}\sqrt{\mu } \p_{t}\phi^{j}_{b} (t,x)$. Note that $\p_{t}\phi_{b}^{j}$ solves $-\Delta \p_{t}\phi_{b}^{j} = \p_{t}b_{j}(t,x)$ with $\p_{t}\phi_{b}^{j}(t,x)|_{\p\O} = 0$. We taking difference quotient for $\p_{t}f$ in (\ref{eqtn_f}) and with help of Poincar\'e inequality, we get
		\Be \label{bt_extra}
		\begin{split}
			(\ref{weakformulation})_{\phi_f} \big|_{\psi=\psi_{b,j}^{t}} &= \int^t_s\iint_{\Omega \times  \R^{3}} v_{j} \sqrt{\mu } \p_{t}\phi_{b}^{j} v\cdot \nabla_x \phi_f \sqrt{\mu}  \\
			&\lesssim \varepsilon \int^t_s \|\nabla\Delta^{-1}\p_{t}b_{j}(\tau)\|_{2}^{2} +\int^t_s \|a(\tau)\|_{2}^{2},	\\
			(\ref{weakformulation})_{P} \big|_{\psi=\psi_{b,j}^{t}} &= \int_s^t \iint_{\O \times \R^3} \sqrt{\mu} f \p_{j} \phi_f \cdot \p_{t}\phi_{b}^{j}    \\
			&\lesssim \varepsilon \int^t_s \|\nabla\Delta^{-1}\p_{t}b_{j}(\tau)\|_{2}^{2} + \int^t_s \|a(\tau)\|_{2}^{2} + \int_{s}^{t} \|(\mathbf{I-P})f\|_{\nu}^{2}	\\
		\end{split}
		\Ee
		Since $\psi_{b,j}^{t}$ vanishes when it acts with $Lf$ and $\Gamma(f,f)$, and boundary integral $(B)$ vanishes by Dirichlet boundary condition of $\phi_{c}$ , from (\ref{bt_extra}) and (\ref{weakformulation}), we obtain
		\begin{equation} \label{b time}
		\begin{split}
		&\int_s^t \int_{\Omega} \p_{t}\phi_{b}^{j}(\tau,x)\partial_t b_{j}(\tau,x) dx = \int_s^t \|\nabla\Delta^{-1} \p_{t} b_{j}(\tau)\|_{2}^{2}    \\
		&\lesssim \varepsilon \int^t_s \|\nabla\Delta^{-1}\p_{t}b_{j}(\tau)\|_{2}^{2} +\int^t_s  \|a(\tau)\|_{2}^{2}  + \int_{s}^{t} \|(\mathbf{I-P})f(\tau)\|_{\nu}^{2} . \\	
		\end{split}
		\end{equation}
		
		Now we combine (\ref{weakformulation}), (\ref{extra}), (\ref{extraP}), (\ref{bC}), (\ref{bB}), (\ref{bT}), (\ref{b Re}), (\ref{bC2}), and (\ref{bB2T2Re2}) for all $i,j$ with proper constant weights. In particular we note that RHS of (\ref{bC}) is cancelled by the first term on the RHS of (\ref{bC2}). Therefore,
		\Be \label{b est}
		\begin{split}
			&\int^t_s \|b(\tau)\|_{2}^{2}  = - \sum_{i,j} \int_s^t \int_{\O} b_{i} (\p_{jj} \Delta^{-1} b_{i})   \\
			&\lesssim G_{b}(t) - G_{b}(s) + \int_{s}^{t} \|(\mathbf{I-P})f(\tau)\|_{\nu}^{2} + \int^t_s |(1-P_{\gamma})f(\tau)|_{2,+} \\
			&\quad + \int_{s}^{t} \|\nu^{-1/2}\Gamma(f,f)(\tau)\|_{2}^{2} + \int_{s}^{t} \|wf(\tau)\|^{2}_{\infty} \|\mathbf{P}f(\tau)\|_{2}^{2}  \\
			&\quad + \int^t_s \|c(\tau)\|_{2}^{2} + \varepsilon\int^t_s \|a(\tau)\|_{2}^{2} , \quad G_{b}(t) \lesssim \|f(t)\|_{2}^{2},\quad \varepsilon \ll 1. \\
		\end{split}
		\Ee
		
		Finally we combine (\ref{c est}), (\ref{a est}), and (\ref{b est}) with $\varepsilon \ll 1$ to conclude (\ref{estimate_dabc}).\end{proof}\unhide

	\section{$L^2$ coercivity}

	The main purpose of this section is to prove the following:
	
	\begin{proposition}
		\label{dlinearl2}Suppose $(f, \phi)$ solves (\ref{eqtn_f}), (\ref{phi_f}), and (\ref{BC_f}). Then 
		%
		%
		there is $0<\lambda_{2} \ll 1$ such that for $0 \leq s \leq t$, 
		\Be\begin{split}\label{completes_dyn}
			& \| e^{\lambda_{2} t}f(t)\|_2^2
			+ \| e^{\lambda_{2} t} \nabla \phi (t) \|_2^2\\
			&
			+  \int_s^t \| e^{\lambda_{2} \tau}  f (\tau)\|_\nu^2 
			+ \| e^{\lambda_{2} \tau} \nabla \phi_f(\tau)\|_2^2 
			\mathrm{d} \tau 
			+  \int_s^t | e^{\lambda_{2} \tau} f |^2_{2,+}    \\
			\lesssim & \ \|e^{\lambda_{2} s} f(s)\|_2^2 +   \|e^{\lambda_{2} s} \nabla \phi_f(s)\|_2^2  \\&   + \sup_{s \leq \tau \leq t} \| w_{\vartheta} f (\tau) \|_\infty \int_s^t \| e^{\lambda_{2} \tau}  f (\tau)\|_\nu^2  .  
		\end{split}\Ee
	\end{proposition}
	
	The null space of linear operator $L$ is a five-dimensional subspace of $L_{v}^{2} (\R^{3})$ spanned by orthonormal vectors $\big\{ \sqrt{\mu}, v\sqrt{\mu},  \frac{|v|^{2}-3}{2} \sqrt{\mu}\big\}$ and the projection of $f$ onto the null space $N(L)$ is denoted by  
	\begin{equation} \label{Pabc}
	\mathbf{P}f (t,x,v) \ := \ \Big\{ a(t,x) + v\cdot  b(t,x) + \frac{|v|^{2}-3}{2} c(t,x)\Big\} \sqrt{\mu }.
	\end{equation} 
	In order to prove the proposition we need the following:
	
	\begin{lemma}
		\label{dabc}
		There exists
		a function $G(t)$ such that, for all $0\le s\leq t$, $G(s)\lesssim
		\|f(s)\|_{2}^{2}$ and 
		\begin{equation}\begin{split}\label{estimate_dabc}
		&\int_{s}^{t}\|\mathbf{P}f(\tau)\|_{\nu }^{2}  + \int^t_s \| \nabla \phi_f \|_2^2\\
		\lesssim & \ 
		G(t)-G(s)
		+  \int_{s}^{t}\|(\mathbf{I}-\mathbf{P}%
		)f(\tau)\|_{\nu }^{2} +\int_{s}^{t}|(1-P_{\gamma })f(\tau)|_{2,+}^{2} \\
		& +\int_{s}^{t}\| \nu^{-1/2} { {\Gamma(f,f)} } \|_{2}^{2} + \int_{s}^{t} \|w_{\vartheta}f(\tau)\|^{2}_{\infty} \|\mathbf{P}f(\tau)\|_{2}^{2} .
		\end{split}\end{equation}
	\end{lemma}
	
	\begin{proof}[\textbf{Proof of Proposition \ref{dlinearl2}}] 
		\textit{Step 1. }  Without loss of generality we prove the result with $s=0$. We have an $L^2$-estimate from $e^{\lambda_{2} t} \times (\ref{eqtn_f})$
		\Be\notag
		\begin{split}
			&\|e^{\lambda_{2} t}  f(t) \|_2^2 - \| f(0) \|_2^2 +  \int^{t}_{0} | e^{\lambda_{2} \tau }
			(1-P_{\gamma} ) f^{j}|_{2,+}^{2}
			\\
			& + \int_0^t  \iint_{\O \times \R^3} v \cdot \nabla \phi_f e^{2\lambda_{2} \tau} |f|^2+
			2\int^t_0 \iint_{\O \times \R^3}e^{\lambda_{2} \tau}2  f Lf\\
			= &\ 2\int^t_0 \iint_{\O \times \R^3} e^{2\lambda_{2} \tau}f \Gamma(f,f) - 2\int^t_0 e^{2\lambda_{2} \tau}\int_{\O  } \nabla \phi_f \cdot \int_{\R^3} v \sqrt{\mu} f\\
			& 
			+2 \lambda_{2}  \int_0^t \|e^{\lambda_{2} \tau}  f(\tau) \|_2^2,
		\end{split}
		\Ee 
		where 
		\Be \label{Pgamma}
		P_\g f := c_{\mu}\sqrt{\mu(v)} \int_{n(x)\cdot u>0} f(u) \sqrt{%
			\mu(u)} \{n(x) \cdot u\} \mathrm{d} u.
		\Ee 
		On the other hand multiplying $\sqrt{\mu(v)} \phi_f(t,x)$ with a test function $\psi(t,x)$ to (\ref{eqtn_f}) and applying the Green's identity, we obtain
		\Be\begin{split}\notag
			&
			\int_{\O} \nabla \phi_f(t,x)  \cdot \int_{\R^3}  v \sqrt{\mu} f \dd v\dd x
			\\
			= & \ \int_{\O} \phi_f(t,x) \p_\tau\left(\int_{\R^3} f (\tau) \sqrt{\mu} \dd v\right)  \dd x
			+ \iint_{\p\O \times \R^3}  \phi_f(t,x) f \sqrt{\mu} \{n \cdot v\} \dd v\dd S_x .
		\end{split}\Ee
		From (\ref{null_flux}), the last boundary contribution equals zero. Now we use (\ref{phi_f}) and deduce that 
		\Be
		\begin{split}\notag
			& \int^t_0e^{2 \lambda_{2} \tau}\int_{\O} \phi_f(t,x) \p_\tau\left(\int_{\R^3} f(\tau)  \sqrt{\mu} \dd v\right)  \dd x \dd \tau\\
			= & \ -  \int^t_0e^{2 \lambda_{2} \tau}\int_{\O} \phi_f(t,x) \p_\tau  \Delta_x  \phi_f(\tau,x) \dd x \dd \tau\\
			= &  \  \frac{1}{2}\int^t_0e^{2 \lambda_{2} \tau}\int_{\O}  \p_\tau  |\nabla_x  \phi_f(\tau,x)|^2 \dd x \dd \tau
			\\
			= & \ \frac{1}{2}   \left(\int_{\O} e^{2 \lambda_{2} t} |\nabla_x \phi_f (t,x)|^2  \dd x\right)- \frac{1}{2}   \left(\int_{\O} |\nabla_x \phi_f (0,x)|^2  \dd x\right)\\
			& \ 
			- \lambda_{2} \int^t_0 e^{2 \lambda_{2} \tau} \int_\O |\nabla_x  \phi_f(\tau,x)|^2 \dd x \dd \tau
			.
		\end{split}
		\Ee
		
		Hence we derive 
		\Be\notag
		\begin{split}
			&\| e^{\lambda_{2} t} f(t) \|_2^2 + \| e^{\lambda_{2} t}\nabla\phi_f (t) \|_2^2  + \int_0^t  \iint_{\O \times \R^3}
			e^{2\lambda_{2} \tau}
			v \cdot \nabla \phi_f |f|^2\\
			& +
			2C\int^t_0 \iint_{\O \times \R^3}\| e^{\lambda_{2} \tau} (\mathbf{I} - \mathbf{P}) f \|_\nu^2
			+ \int^{t}_{0} | e^{\lambda_{2} \tau }
			(1-P_{\gamma} ) f^{j}|_{2,+}^{2}
			\\
			\lesssim &\ 
			\| f(0) \|_2^2+ \| \phi_f (0) \|_2^2 +
			\int^t_0 \| e^{\lambda_{2} \tau}  \nu^{-1/2} \Gamma(f,f) \|_2^2\\
			& 
			+\{ \lambda_{2} + o(1)\} \int^t_0 \| e^{\lambda_{2} \tau} f \|_\nu^2 + \lambda_{2} \int^t_0 \| e^{\lambda_{2} \tau} \nabla_x \phi_f \|_2^2
			.
		\end{split}
		\Ee
		Now we apply Lemma \ref{dabc} and add $o(1) \times (\ref{estimate_dabc})$ to the above inequality and choose $0< \lambda_{2} \ll 1$ to conclude (\ref{completes_dyn}) except the full boundary control. 
		
		\vspace{4pt}
		
		\textit{Step 2. } Note that from (\ref{Pgamma}), $P_\gamma f = z(t,x)\sqrt{\mu(v)}$ for a suitable function $z(t,x)$ on the boundary. Then for $0 < \e \ll 1$
		\Be
		\begin{split}\notag
			| P_\gamma f |_{\gamma,2}^2 
			= & \ \int_{\p\O } |z(t,x)|^2 \dd x \times \int_{\R^3} \mu (v) |n(x) \cdot v| \dd v\\
			\lesssim & \  \int_{\p\O } |z(t,x)|^2 \dd x \times \int_{ \gamma_+(x) \backslash \gamma_+^\e(x)} \mu (v)^{3/2} |n(x) \cdot v| \dd v\\
			= & \  | \mathbf{1}_{ \gamma_+  \backslash \gamma_+^\e }\mu^{1/4} P_\gamma f   |_{2,+}^2.
		\end{split}
		\Ee
		Since $P_\gamma f = f - (1- P_\gamma) f$ on $\gamma_+$ we have $$
		| \mathbf{1}_{ \gamma_+  \backslash \gamma_+^\e } \mu^{1/4} P_\gamma f|_{2, +}^2
		\lesssim |\mathbf{1}_{ \gamma_+  \backslash \gamma_+^\e }\mu^{1/4}  f |_{ 2,+}^2 + |(1- P_\gamma) f|_{ 2,+}^2.$$ Therefore
		\Be\label{Pgamma_bound}
		\int^t_0 | P_\gamma f|_{\gamma,2}^2 \lesssim \int^t_0 |\mathbf{1}_{ \gamma_+  \backslash \gamma_+^\e }\mu^{1/4}  f |_{ 2,+}^2 +\int^t_0 |(1- P_\gamma) f|_{ 2,+}^2.
		\Ee
		Note that 
		\Be\begin{split}\notag
			&\big|[\p_t + v\cdot \nabla_x - \nabla \phi \cdot \nabla_v] (\mu^{1/4}  f )\big|\\
			\lesssim&\ \mu^{1/4} \{ |v|| \nabla_x \phi|  f  +  |v|| \nabla_x \phi|  +  |Lf| + |\Gamma (f,f) | \}.
		\end{split}\Ee
		By the trace theorem Lemma \ref{le:ukai},  
		\Be\label{est_trans_f}
		\begin{split}
			& \int^t_0 |\mathbf{1}_{ \gamma_+  \backslash \gamma_+^\e }\mu^{1/4}  f |_{ 2,+}^2\\
			\lesssim & \ \| f_0 \|_2 + (1 + \| wf \|_\infty) \int^t_0 \| f \|_2^2 + \int^t_0\| \nabla \phi \|_2^2.
		\end{split}
		\Ee
		Adding $o(1) \times (\ref{Pgamma_bound})$ to the result of \text{Step 1} and using (\ref{est_trans_f}) we conclude (\ref{completes_dyn}).\end{proof}


	\begin{proof}[\textbf{Proof of Lemma \ref{dabc}}]
		From the Green's identity, a solution $(f,\phi_f)$ satisfies
		\begin{equation} \label{weakformulation}
		\begin{split}
		&\iint_{\O \times \R^3} f(t) \psi(t) - \iint_{\O \times \R^3} f(s) \psi(s)\\
		&
		- \underbrace{ \int_s^t \iint_{\O \times \R^3} f \p_t\psi }_{(\ref{weakformulation})_{T}} +
		\underbrace{ \int^t_s\int_{\gamma }\psi f (v\cdot n(x)) }_{(\ref{weakformulation})_{B}}  \\
		& - \underbrace{ \int^t_s\iint_{\Omega \times  \R^{3}} \mathbf{P}f {v}\cdot \nabla _{x}\psi  }_{(\ref{weakformulation})_{C}} - \int^t_s\iint_{\Omega \times  \R^{3}} (\mathbf{I-P})f {v}\cdot \nabla _{x}\psi\\
		& + \underbrace{ \int_s^t \iint_{\O \times \R^3} \sqrt{\mu} f \nabla_x \phi_f \cdot \nabla_v \Big[\frac{1}{\sqrt{\mu}} \psi\Big] }_{(\ref{weakformulation})_{P}}
		\\
		=& \ \int^t_s\iint_{\Omega \times \R^{3}}\psi \{ -L  (\mathbf{I}-\mathbf{P%
		})f + \Gamma(f,f)\}\\
		&- \underbrace{\int^t_s\iint_{\Omega \times  \R^{3}} \psi v\cdot \nabla_x \phi_f \sqrt{\mu}}_{(\ref{weakformulation})_{\phi_f}}  .  
		\end{split}
		\end{equation}%
		
		As \cite{EGKM,EGKM2} we use a set of test functions:
		\Be \label{tests}
		\begin{split}
			\psi_a  &\equiv   (|v|^{2}-\beta_{a} )\sqrt{\mu }v\cdot\nabla_x\varphi _{a},  \\
			\psi^{i,j}_{b,1} &\equiv (v_{i}^{2}-\beta_ b)\sqrt{\mu }\partial _{j}\varphi _{b}^{j}, \quad i,j=1,2,3,   \\
			\psi^{i,j}_{b,2} &\equiv |v|^{2}v_{i}v_{j}\sqrt{\mu }\partial _{j}\varphi _{b}^{i}(x),\quad i\neq j,  \\
			\psi_c &\equiv (|v|^{2}-\beta_c )\sqrt{\mu }v \cdot \nabla_x \varphi_{c},  \\
		\end{split}
		\Ee
		where $\varphi_{a}(t,x)$, $\varphi_{b}(t,x)$, and $\varphi_{c}(t,x)$ solve
		\Be\begin{split}\label{phi_abc}
			- \Delta \varphi_a &= a(t,x),  \quad \p_{n}\varphi_{a} \vert_{\p\O} = 0,  \\
			- \Delta \varphi_b^j &= b_j(t,x),   \ \ \varphi^j_b|_{\p\O} =0,\  \text{and}  \ - \Delta \varphi_c  = c(t,x),   \ \ \varphi_c|_{\p\O} =0,
		\end{split}
		\Ee
		and $\beta_a=10$, $\beta_b=1$, and $\beta_c= 5$ such that for all $i=1,\notag2,3,$
		\Be \label{defbeta}
		\begin{split}
			&\int_{{\R}^{3}}(|v|^{2}-\beta_a )(\frac{|v|^{2}}{2}-\frac{3}{2}%
			) v_{i} ^{2}\mu(v)\dd v =0 ,\\
			&\int_{\R^3} (v_i^2 - \beta_b) \mu(v) \dd v=0 ,\\
			&\int_{\R^3} (|v|^{2}-\beta_c )v_{i}^{2}\mu(v) \dd v=0. 
		\end{split}
		\Ee
		
		\vspace{4pt}

		\textit{Step 1. } Estimate of $(\ref{weakformulation})_{\phi_f}$: From (\ref{tests}) and (\ref{defbeta}), we have $(\ref{weakformulation})_{\phi_f}\equiv0$ for $\psi_{b,1}^{i,j}$, $\psi_{b,2}^{i,j}$, and $\psi_c$. For $\psi = \psi_a$, $(\ref{weakformulation})_{\phi_f}$ equals
		\Be \label{extra}
		\begin{split}
			(\ref{weakformulation})_{\phi_f} \big\vert_{\psi=\psi_{a}} &= \int_{\R^3} (|v|^2- \beta_a) (v_1)^2\mu \dd v \int^t_s\int_{\O  } \nabla \varphi_a \cdot \nabla \phi_f \\
			&= C \int^t_s \| \nabla \phi_f \|_2^2 ,  \\
			(\ref{weakformulation})_{\phi_f} \big\vert_{\psi=\psi_{b},\psi_{c}} &= 0,  \\
		\end{split}\Ee
		because $\varphi_a = \phi_f$ from the definitions of (\ref{phi_f}) and (\ref{phi_abc}). \\
		
		Estimate of $(\ref{weakformulation})_{P}$: From (\ref{tests}),
		\Be \label{extraP}
		\begin{split}
			(\ref{weakformulation})_{P} &= \int_s^t \iint_{\O \times \R^3} \sqrt{\mu} f \nabla_x \phi_f \cdot \nabla_v \Big[\frac{1}{\sqrt{\mu}} \psi\Big],\quad  \psi = \psi_{a,b,c},  \\
			&\lesssim \int_{s}^{t} \|w_{\vartheta}f\|_{\infty} \int_{\O}  \nabla_{x}\phi_{f}\cdot\nabla_{x}\varphi_{a,b,c}  \lesssim \int_s^t \|w_{\vartheta}f(\tau)\|_{\infty} \|\mathbf{P}f(\tau)\|_{2}^{2}, 
		\end{split}
		\Ee
		by elliptic estimate $\|\nabla_{x}\varphi_{a,b,c}\|_{2} \lesssim \|\varphi_{a,b,c}\|_{H^{2}} \lesssim \|\mathbf{P}f\|_{2} $.  \\
		
		{\it Step 2}. {Estimate of} \ ${c}$ : We apply $\varphi_{c}$ to (\ref{weakformulation}). From oddness in velocity integration and (\ref{defbeta}), $(\ref{weakformulation})_{C}$ becomes,
		\Be \label{cC}
		\begin{split}
			\int^t_s\iint_{\Omega \times  \R^{3}} \mathbf{P}f {v}\cdot \nabla
			_{x}\psi_{c} &= C_{1} \int^t_s \|c(\tau)\|_{2}^{2},
		\end{split}	
		\Ee
		where $\int_{\mathbf{R}^3}(|v|^{2}-\beta_c )v_{i}^{2} (\frac{|v|^{2}}{2}-\frac{3}{2})\mu(v)dv= 10\pi\sqrt{2\pi}$. For boundary integral $(\ref{weakformulation})_{B}$, we decompose $f_{\gamma } = P_{\gamma }f +\mathbf{1}_{\gamma_{+}}(1-P_{\gamma })f$. Then from (\ref{defbeta}) and trace theorem $|\nabla\varphi_{c}|_{2} \lesssim \|\varphi_{c}\|_{H^{2}} \lesssim \|c\|_{2}$,
		\Be \label{cB}
		\begin{split}
			&\int^t_s\int_{\gamma }\psi_{c} f (v\cdot n(x)) = \int^t_s\int_{\gamma }\psi_{c} \mathbf{1}_{\gamma_{+}} (1-P_{\gamma})f d\gamma  \\
			&\lesssim \varepsilon\int^t_s \|c(\tau)\|_{2}^{2} + C_{\varepsilon}\int^t_s |(1-P_{\gamma})f(\tau)|^{2}_{2,+} ,\quad \varepsilon \ll 1. \\
		\end{split}
		\Ee
		If we define
		\Be \label{Re}
		\begin{split}
			Re &:= (\ref{weakformulation})_{P} + (\ref{weakformulation})_{\phi_{f}} +  \int^t_s\iint_{\Omega \times \R^{3}}\psi \{ L  (\mathbf{I}-\mathbf{P%
			})f - \Gamma(f,f)\}   \\
			&\quad - \int^t_s\iint_{\Omega \times  \R^{3}} (\mathbf{I-P})f {v}\cdot \nabla _{x}\psi,
		\end{split}
		\Ee
		then 
		\Be\begin{split} \label{c Re}
		&\int_{s}^{t} \iint_{\O\times \R^{3}} \psi_{c} Re\vert_{\psi_{c}}\\
		 \lesssim & \  \varepsilon \int_{s}^{t} \|c\|_{2}^{2} + \int_{s}^{t} \|(\mathbf{I-P})f(\tau)\|_{\nu}^{2} + \int_{s}^{t} \|\nu^{-1/2}\Gamma(f,f)(\tau)\|_{2}^{2},
	\end{split}	\Ee
		from elliptic estimate and Young's inequality. We also use even/oddness in velocity integration , (\ref{defbeta}), and Young's inequality to estimate,
		\Be \label{cT}
		\begin{split}
			(\ref{weakformulation})_{T}\vert_{\psi_{c}} &= \int_s^t \iint_{\O \times \R^3} f \p_t\psi_{c} = \int_s^t \iint_{\O \times \R^3} (\mathbf{I-P})f \p_t\psi_{c}  \\
			&\lesssim \varepsilon \int_s^t \|\nabla\Delta^{-1}\p_{t}c(\tau)\|_{2}^{2} + \int_{s}^{t} \|(\mathbf{I-P})f(\tau)\|_{\nu}^{2}.
		\end{split}
		\Ee
		
		Now, we choose a new test function $\psi_{c}^{t} := (\frac{|v|^{2}}{2} - \frac{3}{2})\sqrt{\mu } \p_{t}\varphi_{c} (t,x)$. Note that $\p_{t}\varphi_{c}$ solves $-\Delta \p_{t}\varphi_{c} = \p_{t}c(t,x)$ with $\p_{t}\varphi_{c}(t,x)|_{\p\O} = 0$. We taking difference quotient for $\p_{t}f$ in (\ref{eqtn_f}) and it replace first three terms in the LHS of (\ref{weakformulation}). With help of Poincar\'e inequality $\|\p_{t}\varphi_{c}\|_{2} \lesssim \|\nabla\p_{t}\varphi_{c}\|_{2}$, we can also compute
		\Be \label{ct_extra}
		\begin{split}
			(\ref{weakformulation})_{\phi_f} \big|_{\psi=\psi_{c}^{t}} &= \int^t_s\iint_{\Omega \times  \R^{3}} (\frac{|v|^{2}}{2} - \frac{3}{2})\sqrt{\mu } \p_{t}\varphi_{c} v\cdot \nabla_x \phi_f \sqrt{\mu}  \\
			&\lesssim \varepsilon \int^t_s \|\nabla\Delta^{-1}\p_{t}c(\tau)\|_{2}^{2} +\int^t_s \|a(\tau)\|_{2}^{2},	\\
			(\ref{weakformulation})_{P} \big|_{\psi=\psi_{c}^{t}} &= \int_s^t \iint_{\O \times \R^3} \sqrt{\mu} f \nabla_x \phi_f \cdot 2v \p_{t}\varphi_{c}    \\
			&\lesssim \varepsilon \int^t_s \|\nabla\Delta^{-1}\p_{t}c(\tau)\|_{2}^{2} + \int^t_s \big(\|a(\tau)\|_{2}^{2} + \|b(\tau)\|_{2}^{2}\big) \\
			& \ \ \ + \int_{s}^{t} \|(\mathbf{I-P})f\|_{\nu}^{2}.	\\
		\end{split}
		\Ee
		Since $\psi_{c}^{t}$ vanishes when it acts with $Lf$ and $\Gamma(f,f)$, and boundary integral $(\ref{weakformulation})_{B}$ vanishes by Dirichlet boundary condition of $\varphi_{c}$ , from (\ref{ct_extra}) and (\ref{weakformulation}), we obtain
		\begin{equation} \label{c time}
		\begin{split}
		&\int_s^t \int_{\Omega} \p_{t}\varphi_{c}(\tau,x)\partial_t c(\tau,x) dx = \int_s^t \|\nabla\Delta^{-1} \p_{t} c(\tau)\|_{2}^{2}    \\
		&\lesssim \varepsilon \int^t_s \|\nabla\Delta^{-1}\p_{t}c(\tau)\|_{2}^{2} +\int^t_s \big( \|a(\tau)\|_{2}^{2} + \|b(\tau)\|_{2}^{2} \big)  \\
		& \ \ \ + \int_{s}^{t} \|(\mathbf{I-P})f(\tau)\|_{\nu}^{2} . \\	
		\end{split}
		\end{equation}
		
		We combine (\ref{weakformulation}), (\ref{extra}), (\ref{extraP}), (\ref{cB}), (\ref{c Re}), (\ref{cT}), and (\ref{c time}) with $\varepsilon \ll 1$ to obtain
		\Be \label{c est}
		\begin{split}
			&\int^t_s \|c(\tau)\|_{2}^{2} \\
			&\lesssim G_{c}(t) - G_{c}(s) + \int_{s}^{t} \|(\mathbf{I-P})f(\tau)\|_{\nu}^{2} + \int^t_s |(1-P_{\gamma})f(\tau)|^{2}_{2,+} \\
			&\quad + \int_{s}^{t} \|\nu^{-1/2}\Gamma(f,f)(\tau)\|_{2}^{2} + \int_{s}^{t} \|w_{\vartheta}f(\tau)\|^{2}_{\infty} \|\mathbf{P}f(\tau)\|_{2}^{2}  \\
			&\quad + \varepsilon\int^t_s \big(\|a(\tau)\|_{2}^{2} + \|b(\tau)\|_{2}^{2}\big),
		\end{split}
		\Ee
		for $\varepsilon \ll 1$ where $G_{c}(t) := \iint_{\O\times \R^{3}} f(t)\psi_{c}(t) \lesssim \|f(t)\|_{2}^{2}$.  \\
		
		{\it Step 3.} {Estimate of} \ ${a}$ : From mass conservation $\int_{\O} a(t,x) dv = 0$, $\varphi_{a}$ in (\ref{phi_abc}) is well-defined. Moreover, we choose $\varphi_{a}$ so that has mean zero, $\int_{\O} \varphi_{a}(t,x) dx = 0$. Therefore, Poincar\'e inequality $\|\varphi_{a}\|_{2} \lesssim \|\nabla \varphi_{a}\|_{2}$ holds and these are also true for $\p_{t}\varphi_{a}$ which solves same elliptic equation with Neumann boundary condition.  
		
		By even/oddness in velocity integral and $\b_{a}$ defined in (\ref{defbeta}), we can replace $c$ into $a$ in estimates (\ref{cC}) and (\ref{c Re}). For boundary integral $(\ref{weakformulation})_{B}$, we decompose $f_{\gamma } = P_{\gamma }f +\mathbf{1}_{\gamma_{+}}(1-P_{\gamma })f$. From Neumann boundary condition $\p_{n}\varphi_{a} = 0$ and oddness in velocity integral, $\int_{\gamma} \psi_{a} P_{\gamma}f (v\cdot n(s)) = 0$ and we obtain similar esimate as (\ref{cB}),
		\Be \label{aB}
		\begin{split}
			&\int^t_s\int_{\gamma }\psi_{a} f (v\cdot n(x)) = \int^t_s\int_{\gamma }\psi_{a} \mathbf{1}_{\gamma_{+}} (1-P_{\gamma})f d\gamma  \\
			&\lesssim \varepsilon\int^t_s \|a(\tau)\|_{2}^{2} + C_{\varepsilon}\int^t_s |(1-P_{\gamma})f(\tau)|^{2}_{2,+} ,\quad \varepsilon \ll 1. \\
		\end{split}
		\Ee
		For $(\ref{weakformulation})_{T}$, from oddness,
		\Be \label{aT}
		\begin{split}
			(\ref{weakformulation})_{T}\vert_{\psi_{a}} &= \int_s^t \iint_{\O \times \R^3} f \p_t\psi_{a} = \int_s^t \iint_{\O \times \R^3} [ \sqrt{\mu} (b\cdot v) + (\mathbf{I-P})f] \p_t\psi_{a}  \\
			&\lesssim \varepsilon \int_s^t \|\nabla\Delta^{-1}\p_{t}a(\tau)\|_{2}^{2} + \int_{s}^{t} \|b(\tau)\|_{2}^{2}+ \int_{s}^{t} \|(\mathbf{I-P})f(\tau)\|_{\nu}^{2}.
		\end{split}
		\Ee
		Now let us estimate $\int_{s}^{t} \|\nabla\Delta^{-1}\p_{t}a(\tau)\|_{2}^{2}$ which appear in (\ref{cT}) type estimate. We use new test function $\psi_{a}^{t} = \bar{\varphi}(x)\sqrt{\mu}$. It easy to check 
		\Be \label{at_extra}
		\begin{split}
			(\ref{weakformulation})_{\phi_f} \big|_{\psi=\sqrt{\mu}} &= \int^t_s\iint_{\Omega \times  \R^{3}} \sqrt{\mu } v\cdot \nabla_x \phi_f \bar{\varphi}(x) \sqrt{\mu} =0, 	\\
			(\ref{weakformulation})_{P} \big|_{\psi=\sqrt{\mu}} &= \int_s^t \iint_{\O \times \R^3} \sqrt{\mu} f \nabla_x \phi_f \cdot \nabla_v \bar{\varphi}(x) = 0,  \\
		\end{split}
		\Ee
		and from taking difference quotient and null condition (\ref{null_flux}), we obtain, for almost $t$,
		\begin{equation*}
		\int_{\Omega } \bar{\varphi} \partial _{t}a=\sqrt{2\pi}\int_{\Omega }(b\cdot \nabla_x )\bar{\varphi}.
		\end{equation*}
		In particular, if we choose $\bar{\varphi} = 1$, we directly get $\int_{\O} \p_{t} a = 0$, so Neumann problem 
		\begin{equation}
		-\Delta \Phi_a =\partial _{t}a(t)\ ,\ \ \ \ \ \frac{\partial \Phi_a }{\partial n}%
		\Big{|}_{\partial \Omega }=0  \notag 
		\end{equation}
		is well-defined. For dual pair $(H^{1})^{\ast }\equiv(H^{1}(\Omega ))^{\ast }$ with respect to $\langle A,B \rangle=\int_{\Omega } A \cdot B dx$ for $A\in H^{1}$ and $B\in (H^1)^{\ast}$,
		\begin{equation} \label{a time}
		\|\nabla_x \partial_t \varphi_a \|_{2} = \|\Delta^{-1}\partial
		_{t}a(t)\|_{H^{1}}=\|\Phi_a \|_{H^1} \lesssim \|\partial _{t}a(t)\|_{(H^{1})^{\ast }}
		\lesssim  \|b(t)\|_{2}.  \\
		\end{equation}
		
		We change $c$ into $a$ in (\ref{cC}) and (\ref{c Re}) and combine with (\ref{weakformulation}), (\ref{extra}), (\ref{extraP}), (\ref{aB}), (\ref{aT}), and (\ref{a time}) with $\varepsilon \ll 1$ to obtain
		\Be \label{a est}
		\begin{split}
			&\int^t_s \|a(\tau)\|_{2}^{2} + \int^t_s \|\nabla \phi_{f}(\tau)\|_{2}^{2} \\
			&\lesssim G_{a}(t) - G_{a}(s) + \int_{s}^{t} \|(\mathbf{I-P})f(\tau)\|_{\nu}^{2} + \int^t_s |(1-P_{\gamma})f(\tau)|^{2}_{2,+} \\
			&\quad + \int_{s}^{t} \|\nu^{-1/2}\Gamma(f,f)(\tau)\|_{2}^{2} + \int_{s}^{t} \|wf(\tau)\|^{2}_{\infty} \|\mathbf{P}f(\tau)\|_{2}^{2}  + \int^t_s \|b(\tau)\|_{2}^{2} ,  \\
		\end{split}
		\Ee
		for $\varepsilon \ll 1$ where $G_{a}(t) := \iint_{\O\times \R^{3}} f(t)\psi_{c}(t) \lesssim \|f(t)\|_{2}^{2}$.  \\
		
		{\it Step 4.} {Estimate of} \ $ {b}$ : For fixed $i,j$, we choose test function $\psi = \psi_{b,1}^{i,j}$ in (\ref{tests}) where $\b_{b}$ and $\varphi_{b}$ are defined in (\ref{defbeta}) and (\ref{phi_abc}). From oddness in velocity integration and definition of $\b_{b}$, $(\ref{weakformulation})_{C}$ in (\ref{weakformulation}) yields
		\Be \label{bC}
		\begin{split}
			(\ref{weakformulation})_{C}\vert_{\psi_{b,1}^{i,j}} &:= \int_s^t \iint_{\O\times \R^{3}} \mathbf{P}f v\cdot\nabla \psi_{b,1}^{i,j} = C_{3} \int_s^t \int_{\O} b_{i} (\p_{ij} \Delta^{-1} b_{j}), 
		\end{split}
		\Ee
		where $C_{3} := \int_{\R^{3}} (v_{i}^{2} - \b_{b}) v_{i}^{2} \mu dv = 2\sqrt{2\pi}$. For boundary integration, contribution of $P_{\gamma}f$ vanishes by oddness. 
		\Be \label{bB}
		\begin{split}
			(\ref{weakformulation})_{B}\vert_{\psi_{b,1}^{i,j}} &:= \int_s^t \int_{\gamma} \psi_{b,1}^{i,j} \mathbf{1}_{\gamma_{+}} (1-P_{\gamma})f (v\cdot n(x)) \\ & \lesssim \varepsilon \int_s^t \|b(\tau)\|^{2}_{2} + \int_s^t |(1-P_{\gamma})f|_{2,+}^{2},
		\end{split}
		\Ee
		and similar as (\ref{cT}) and (\ref{c Re}), we use oddness and definition of $\b_{b}$ to vanish contribution of $a$ and $b$. We obtain
		\Be \label{bT}
		\begin{split}
			&(\ref{weakformulation})_{T}\vert_{\psi_{b,1}^{i,j}}\\
			 &\lesssim  \varepsilon \int_s^t \|\nabla\Delta^{-1}\p_{t} b_{j}(\tau)\|_{2}^{2} + \int_{s}^{t} \|c(\tau)\|_{2}^{2} + \int_{s}^{t} \|(\mathbf{I-P})f(\tau)\|_{\nu}^{2},
		\end{split}
		\Ee
		\Be \label{b Re}\begin{split}
		\int_{s}^{t} \iint_{\O\times \R^{3}} \psi_{b,1}^{i,j} Re\vert_{\psi_{b,1}^{i,j}} \lesssim& \varepsilon \int_{s}^{t} \|b(\tau)\|_{2}^{2} + \int_{s}^{t} \|(\mathbf{I-P})f(\tau)\|_{\nu}^{2} \\&+ \int_{s}^{t} \|\nu^{-1/2}\Gamma(f,f)(\tau)\|_{2}^{2}.
		\end{split}\Ee
		
		Next, we try test function $\psi_{b,2}^{i,j}$ with $i \neq j$ to obtain 
		\Be \label{bC2}
		\begin{split}
			(\ref{weakformulation})_{C}\vert_{\psi_{b,2}^{i,j}} & := \int_s^t \iint_{\O\times \R^{3}} (b\cdot v)\sqrt{\mu} v\cdot\nabla \psi_{b,2}^{i,j}\\
			& = C_{4} \int_s^t \int_{\O} \big(  b_{j} (\p_{ij} \Delta^{-1} b_{i}) + b_{i} (\p_{jj} \Delta^{-1} b_{i}) \big) .
		\end{split}
		\Ee
		by oddness in velocity integral where $C_{4} := 7\sqrt{2\pi}$. We also have the following three estimates using oddness of velocity integral,
		\Be \label{bB2T2Re2}
		\begin{split}
			(\ref{weakformulation})_{B}\vert_{\psi_{b,2}^{i,j}} &:= \int_s^t \int_{\gamma} \psi_{b,1}^{i,j} \mathbf{1}_{\gamma_{+}} (1-P_{\gamma})f (v\cdot n(x)) \\
			& \lesssim \varepsilon \int_s^t \|b(\tau)\|^{2}_{2} + \int_s^t |(1-P_{\gamma})f|_{2,+}^{2},  \\
			(\ref{weakformulation})_{T}\vert_{\psi_{b,2}^{i,j}} &\lesssim  \varepsilon \int_s^t \|\nabla\Delta^{-1}\p_{t} b_{i}(\tau)\|_{2}^{2} + \int_{s}^{t} \|(\mathbf{I-P})f(\tau)\|_{\nu}^{2},  \\
			\int_{s}^{t} \iint_{\O\times \R^{3}} \psi_{b,2}^{i,j} Re\vert_{\psi_{b,2}^{i,j}} &\lesssim \varepsilon \int_{s}^{t} \|b(\tau)\|_{2}^{2} + \int_{s}^{t} \|(\mathbf{I-P})f(\tau)\|_{\nu}^{2} \\
			& + \int_{s}^{t} \|\nu^{-1/2}\Gamma(f,f)(\tau)\|_{2}^{2}.
		\end{split}
		\Ee
		
		To obtain estimate for $\|\nabla\Delta^{-1} \p_{t} b_{j}\|_{2}$, we use a test function $\psi_{b,j}^{t} := v_{j}\sqrt{\mu } \p_{t}\varphi^{j}_{b} (t,x)$. Note that $\p_{t}\varphi_{b}^{j}$ solves $-\Delta \p_{t}\varphi_{b}^{j} = \p_{t}b_{j}(t,x)$ with $\p_{t}\varphi_{b}^{j}(t,x)|_{\p\O} = 0$. We taking difference quotient for $\p_{t}f$ in (\ref{eqtn_f}) and with help of Poincar\'e inequality, we get
		\Be \label{bt_extra}
		\begin{split}
			(\ref{weakformulation})_{\phi_f} \big|_{\psi=\psi_{b,j}^{t}} &= \int^t_s\iint_{\Omega \times  \R^{3}} v_{j} \sqrt{\mu } \p_{t}\varphi_{b}^{j} v\cdot \nabla_x \phi_f \sqrt{\mu}  \\
			&\lesssim \varepsilon \int^t_s \|\nabla\Delta^{-1}\p_{t}b_{j}(\tau)\|_{2}^{2} +\int^t_s \|a(\tau)\|_{2}^{2},	\\
			(\ref{weakformulation})_{P} \big|_{\psi=\psi_{b,j}^{t}} &= \int_s^t \iint_{\O \times \R^3} \sqrt{\mu} f \p_{j} \phi_f \cdot \p_{t}\varphi_{b}^{j}    \\
			&\lesssim \varepsilon \int^t_s \|\nabla\Delta^{-1}\p_{t}b_{j}(\tau)\|_{2}^{2} + \int^t_s \|a(\tau)\|_{2}^{2}\\& \ \ \ \  + \int_{s}^{t} \|(\mathbf{I-P})f\|_{\nu}^{2}.	\\
		\end{split}
		\Ee
		Since $\psi_{b,j}^{t}$ vanishes when it acts with $Lf$ and $\Gamma(f,f)$, and boundary integral $(\ref{weakformulation})_{B}$ vanishes by Dirichlet boundary condition of $\varphi_{c}$ , from (\ref{bt_extra}) and (\ref{weakformulation}), we obtain
		\begin{equation} \label{b time}
		\begin{split}
		&\int_s^t \int_{\Omega} \p_{t}\varphi_{b}^{j}(\tau,x)\partial_t b_{j}(\tau,x) dx = \int_s^t \|\nabla\Delta^{-1} \p_{t} b_{j}(\tau)\|_{2}^{2}    \\
		&\lesssim \varepsilon \int^t_s \|\nabla\Delta^{-1}\p_{t}b_{j}(\tau)\|_{2}^{2} +\int^t_s  \|a(\tau)\|_{2}^{2}  + \int_{s}^{t} \|(\mathbf{I-P})f(\tau)\|_{\nu}^{2} . \\	
		\end{split}
		\end{equation}
		
		Now we combine (\ref{weakformulation}), (\ref{extra}), (\ref{extraP}), (\ref{bC}), (\ref{bB}), (\ref{bT}), (\ref{b Re}), (\ref{bC2}), and (\ref{bB2T2Re2}) for all $i,j$ with proper constant weights. In particular, we note that RHS of (\ref{bC}) is cancelled by the first term on the RHS of (\ref{bC2}). Therefore,
		\Be \label{b est}
		\begin{split}
			&\int^t_s \|b(\tau)\|_{2}^{2}  = - \sum_{i,j} \int_s^t \int_{\O} b_{i} (\p_{jj} \Delta^{-1} b_{i})   \\
			&\lesssim G_{b}(t) - G_{b}(s) + \int_{s}^{t} \|(\mathbf{I-P})f(\tau)\|_{\nu}^{2} + \int^t_s |(1-P_{\gamma})f(\tau)|^{2}_{2,+} \\
			&\quad + \int_{s}^{t} \|\nu^{-1/2}\Gamma(f,f)(\tau)\|_{2}^{2} + \int_{s}^{t} \|w_{\vartheta}f(\tau)\|^{2}_{\infty} \|\mathbf{P}f(\tau)\|_{2}^{2}  \\
			&\quad + \int^t_s \|c(\tau)\|_{2}^{2} + \varepsilon\int^t_s \|a(\tau)\|_{2}^{2} , \quad G_{b}(t) \lesssim \|f(t)\|_{2}^{2},\quad \varepsilon \ll 1. \\
		\end{split}
		\Ee
		
		Finally we combine (\ref{c est}), (\ref{a est}), and (\ref{b est}) with $\varepsilon \ll 1$ to conclude (\ref{estimate_dabc}).\end{proof}

	\section{Global Existence and Exponential decay}\label{sec_global}
	
	We start the section by proving a crucial interpolation in H\"older spaces.

	\begin{proof}[\textbf{Proof of Lemma \ref{lemma_interpolation}}]
	Let ${\O}_1$ be an open bounded subset of $\R^3$ containing the closure $\bar{\O}$. Suppose ${\phi (t)} \in C^{2,D_2}(\O)$. From a standard extension theorem (e.g. see Lemma 6.37 of \cite{GT} in page 136) there exists a function $\bar{\phi}(t) \in C^{2,D_2} (\O_1)$ and $\bar{\phi}(t)\equiv 0$ in $\R^3 \backslash \O_1$ such that $\phi(t) \equiv \bar{\phi}(t)$ in $\O$ and 
	\Be\label{Holder_extension}
	\begin{split}
	\| \bar{\phi}(t)\|_{C^{1,1-D_1} (\O_1)}  
	&\leq   \  C_{\O,\O_1,D_1, D_2} \| \phi(t)\|_{C^{1,1-D_1} (\O)},\\
	 and \ \ 
	\| \bar{\phi}(t)\|_{C^{2,D_2} (\O_1)} &\leq C_{\O,\O_1,D_1, D_2} \| \phi(t)\|_{C^{2,D_2} (\O)},
	\end{split}\Ee	
	where $C_{\O,\O_1,D_1, D_2}$ does not depend on $\phi(t)$ and $t$.

		Choose arbitrary points $x,y$ in $\R^3$. For $0\leq s\leq 1$, we have $(1-s)x + sy \in \overline{xy}$. Then we derive that 
		\Be\begin{split}\notag
			& [(y-x) \cdot \nabla]\nabla \bar{\phi} (t,(1-s)x + sy)  \\
			=& \ \frac{ 
				[(y-x) \cdot \nabla]\nabla \bar{\phi} (t,(1-s)x + sy) -  [(y-x) \cdot \nabla]\nabla \bar{\phi} (t,x)
			}{|(1-s) x + s y - x|^{D_2}}
			\\
			& \times 
			|(1-s) x + s y - x|^{D_2}\\
			+&   \Big(\frac{y-x}{|y-x|} \cdot \nabla\Big)\nabla \bar{\phi} (t,x) |y-x|\\
			=& \ O( |x-y|^{ 1+ D_2} ) s^{D_2} [\nabla^2 \bar{\phi} (t)]_{C^{0, D_2}}
			+  \Big(\frac{y-x}{|y-x|} \cdot \nabla\Big)  \nabla \bar{\phi} (t,x) |y-x|.
		\end{split}\Ee
		Taking an integration on $s \in [0,1]$, we obtain that 
		\Be\begin{split}\label{C_2_alpha}
			&\left| \Big(\frac{y-x}{|y-x|} \cdot \nabla\Big)  \nabla \bar{\phi} (t,x)\right|\\
			\leq & \  \frac{1}{|y-x|}
			\left| \int_0^1
			[(y-x) \cdot \nabla]\nabla \bar{\phi} (t,(1-s)x + sy) \dd s \right|\\
			&
			+\frac{1}{1+ D_2} |x-y|^{D_2} [\nabla^2 \bar{\phi}(t)]_{C^{0, D_2}}.
		\end{split}\Ee
		
		On the other hand, from an expansion along $s$,
		\[
		\nabla \bar{\phi}(t,y)- \nabla \bar{\phi}(t,x)  =  \int^1_0 [(y-x) \cdot \nabla]\nabla \bar{\phi} (t,(1-s)x + sy) \dd s.
		\]
		We plug this identity into (\ref{C_2_alpha}) and deduce that for $0<D_1<1$
		\Be\begin{split}\label{C2_interpolation}
			&\left|\left(\frac{x-y}{|x-y|} \cdot \nabla \right) \nabla \bar{\phi} (t,x)\right|\\ 
			\leq&  \ \frac{|\nabla \bar{\phi} (t,x) - \nabla \bar{\phi}(t,y)|}{|x-y|}
			+\frac{1}{1+ D_2}  |x-y|^{  {D_2}}    [\nabla^2 \bar{\phi}(t)]_{C_x^{0, {D_2}}}\\
			\leq& \  \frac{1}{|x-y|^{D_1}} [\nabla \bar{\phi} (t)]_{C^{0, 1-D_1}}+ \frac{1}{1 + D_2}  |x-y|^{  {D_2}}    [\nabla^2 \bar{\phi}(t)]_{C_x^{0, {D_2}}}.
		\end{split}\Ee
		Now let us choose 
		\[
		|x-y| = e^{- \Lambda_0 t}, \ \ \ \hat{\omega}:=\frac{x-y}{|x-y|} \in \mathbb{S}^2.
		\]  
		From (\ref{C2_interpolation})
		\Be\notag
		|\big(\hat{\omega} \cdot \nabla\big) \nabla \bar{\phi} (t,x)|
		\leq e^{D_1 \Lambda_0 t} [\nabla \bar{\phi} (t)]_{C^{0, 1- D_1}} + \frac{1}{1+ D_2}
		e^{- D_2 \Lambda_0 t}   [\nabla^2 \bar{\phi}(t)]_{C_x^{0, {D_2}}}.
		\Ee
		Taking supremum in $x$ and $\hat{\o}$ to the above inequality and using $$
		\|\nabla^2_x \bar{\phi}(t )\|_{L^\infty_x } 
		=  \sup_x\sup_{\hat{\omega} \in \mathbb{S}^2}|\big(\hat{\omega} \cdot \nabla\big) \nabla \bar{\phi} (t,x)|,$$ we get
		\Be\begin{split}\notag
			\|\nabla^2_x \bar{\phi}(t )\|_{L^\infty (\O_1)}
			\leq
		&	\ e^{D_1 \Lambda_0t}[\nabla_x \bar{\phi}(t)]_{C^{0,1-D_1}(\O_1)}\\
&			+ e^{- D_2 \Lambda_0t}[\nabla^2 \bar{\phi}(t)]_{C^{0, D_2}(\O_1)}.
		\end{split}\Ee
		Finally from (\ref{Holder_extension}) and the above estimate we conclude (\ref{phi_interpolation}).
	\end{proof}

	Now we are ready to prove the global-in-time result. 
	
	
	\begin{proof}[\textbf{Proof of Theorem \ref{main_existence}}]

	\textit{Step 1. }	For $0< M \ll 1$ and $0<\delta_*\ll 1$, we first assume that an initial datum satisfies
		\Be
		\begin{split}
		\label{initial_M}
		\| w_\vartheta f_0 \|_\infty  + \| w_{\tilde{\vartheta}}  f_0 \|_p + \|w_{\tilde{\vartheta}} \alpha_{f_0,\e}^\beta \nabla_{x,v} f_0 \|_p   \leq \delta_*  M,    \\
		    \| w_{\tilde{\vartheta}} \nabla_v f_0\|_{L^3(\O \times \R^3)}+ \| \nabla_x^2 \phi_f(0) \|_{\infty}< \infty.
		\end{split}\Ee
We will choose $M, \delta_*$ later. For the sake of convenience we choose a large constant $L\gg \max\left(M, \| \nabla_x^2 \phi_f (0) \|_{\infty}\right)$. In order to use the continuation argument along the lines of the local existence theorem, Theorem \ref{local_existence}, we set 
		\Be\label{global_T}
		\begin{split}
			T 
			=& \sup_{t } \Big\{ t \geq 0:     \| e^{\lambda_\infty s} w_\vartheta f(t) \|_\infty  + \| w_{\tilde{\vartheta}} f (t) \|_p    \leq  M ,
			\\
			& \ \ \ \ \ \ \ \     \ \ \     {and}   \ \  \|w_{\tilde{\vartheta}} \alpha_{f,\e}^\beta \nabla_{x,v} f (t) \|_p^p
			+ \int^t_0  | w_{\tilde{\vartheta}}\alpha_{f,\e}^\beta \nabla_{x,v} f(t) |_{p,+}^p  
			<\infty, \\
			& \ \ \ \ \ \ \ \      \ \ \    {and} \   \ \| \nabla_v f(t) \|_{L^3_x(\O) L^{1+ \delta}_v (\R^3)}< \infty, \\
			& \ \ \ \ \ \ \ \    \ \ \    {and}  \  \  \| \nabla_x^2 \phi_f (t) \|_{\infty}\leq L
			\Big\}.
		\end{split}
		\Ee
		Here for fixed $\delta \ll 1$, we choose $\lambda_{\infty}$ such that 
		\Be\label{condition_lambda_infty}
		\begin{split}
			 &  
			 20\sqrt{   C C_2 M }
			  \leq \lambda_\infty \leq  \min \left( \frac{\lambda_2}{2}, \frac{\nu_0}{4}\right), 
 \ \  for \   M \ll 1,
		\end{split}
		\Ee
		where $\lambda_{2}$ is obtained in Proposition \ref{dlinearl2}. Note that from (\ref{Morrey}) the condition (\ref{delta_1/lamdab_1}) holds for $M \ll 1$.
		
		\hide

		\vspace{4pt}

		\textit{Step 2. } From Proposition \ref{prop_W1p}, (\ref{initial_M}), and (\ref{phi_c,gamma_M}), we have for $t \leq T$ 
		%
		\hide\[
		e^{ {C (1 
				+ \sup_{0 \leq s \leq t} \| \nabla^2 \phi_f(s) \|_\infty) t} } \leq e^{C(1+ (C_1 M)^{1/p})t}\delta_* M.
		\]

		\Be\label{est_W}
		\begin{split}
			&\| f(t) \|_p^p
			+\| \alpha_f^\beta 
			\p f (t)\|_p^p \\
			& 
			\lesssim  \   e^{ {C (1 
					+ \sup_{0 \leq s \leq t} \| \nabla^2 \phi_f(s) \|_\infty) t} }
			\{
			\| f(0) \|_p^p
			+\| \alpha_f^\beta 
			\p f (0)\|_p^p \}
			.\end{split} \Ee

		\unhide
		\Be\label{growth_W1p}
		\| \alpha_f^\beta \nabla_{x,v} f(t) \|_p ^p
		+ \| f(t) \|_p^p
		\leq 
		C_p e^{ C_p(1+ (C_1 M)^{1/p}) t} 
		\times \delta_* M  .
		\Ee

		\vspace{4pt}
		
		\unhide
		
		\textit{Step 2. } We claim that 
		\Be\label{decay_C2}
	\sup_{0 \leq t \leq T}	e^{  \frac{\lambda_\infty}{2}t} \| \nabla_x^2 \phi_f(t) \|_\infty \leq
		C_2 M , \ \   with   \ \ C_2 := C_{  \O }+ (C_1 C_p )^{1/p}\delta_*. 
		\Ee
		Here $C_{\O }$ appears in (\ref{Morrey}), and $C_1$ in (\ref{phi_c,gamma}), and $C_p$ in Proposition \ref{prop_W1p}.

		\hide
		We claim that, for $\delta_1>0$, $\lambda_0>0$, and $3<p<6$ 
		\Be\label{decay_phi_C2}
		\begin{split}
			&\| \nabla_x^2 \phi_f(t) \|_\infty\\
			\leq & \  e^{\delta_1 \lambda_0 t}  \| w f (t)\|_\infty   \\
			& \  +   e^{- (1- \frac{3}{p}) \lambda_0 t} 
			e^{ C\int^t_0 1+ \| w f(s) \|_\infty  + \| \phi_f (s) \|_{C^2}  \dd s} \times  \| \alpha^\beta \nabla_{x,v} f_0 \|_p
		\end{split}
		\Ee

		From (\ref{global_T}), (\ref{phi_c,gamma}) and (\ref{growth_W1p}), for all $p>3$
		\Be\label{growth_C2gamma}
		\| \phi_f (t) \|_{C^{2, 1- \frac{3}{p}} (\bar{\O})} \lesssim e^{CMt} M.
		\Ee\unhide
		
		From (\ref{Morrey}) and (\ref{global_T}), for $0 \leq t \leq T$, for all $D_1>0$
		\Be\begin{split}\label{decay_C0,alpha}
			\|  \phi_f (t)\|_{C ^{1,1- D_1}(\bar{\O})}  
			\leq  C_{ \O} \| w_{\vartheta} f (t) \|_\infty 
			\leq  C_{  \O}M e^{-  {\lambda_\infty}  t}.
		\end{split}\Ee
		
		On the other hand, from Proposition \ref{prop_W1p}, (\ref{est_W}), and (\ref{initial_M}), 
		 we derive that for $0\leq t \leq T$ 
		%
		\hide\[
		e^{ {C (1 
				+ \sup_{0 \leq s \leq t} \| \nabla^2 \phi_f(s) \|_\infty) t} } \leq e^{C(1+ (C_1 M)^{1/p})t}\delta_* M.
		\]

		\Be\label{est_W}
		\begin{split}
			&\| w_{\tilde{\vartheta}} f(t) \|_p^p
			+\| \alpha_f^\beta 
			\p f (t)\|_p^p \\
			& 
			\lesssim  \   e^{ {C (1 
					+ \sup_{0 \leq s \leq t} \| \nabla^2 \phi_f(s) \|_\infty) t} }
			\{
			\| f(0) \|_p^p
			+\| \alpha_f^\beta 
			\p f (0)\|_p^p \}
			.\end{split} \Ee

		\unhide
		\Be\label{growth_W1p}
		\begin{split}
			& 
			\| f(t) \|_p^p+\| w_{\tilde{\vartheta}}  \alpha_{f,\e}^\beta \nabla_{x,v} f(t) \|_p ^p 
			+ \int^t_0 | w_{\tilde{\vartheta}}  \alpha_{f,\e}^\beta \nabla_{x,v} f(s)  |_{p,+} ^p
			\\
			\leq  &    \ 
			C_p e^{ C_p(1+ L ) t} 
			\times( \delta_* M)^p  .
		\end{split}
		\Ee
		Now we use Lemma \ref{lemma_apply_Schauder}, from (\ref{phi_c,gamma}), for $p>3$ and $0\leq t \leq T$,
		\Be\label{phi_c,gamma_M}
		 \| \phi_f (t) \|_{C^{2, 1-\frac{3}{p}} (\bar{\O})} \leq 
		(C_1   C_p )^{1/p}  e^{ \frac{1}{p}C_p(1+ L ) t} 
			\times \delta_*  M  .
		\Ee

	\hide	
		On the other hand, from (\ref{phi_c,gamma}) and (\ref{growth_W1p}), 
		\Be\label{phi_c,gamma_M}
		\| \phi_f (t) \| _{C^{2, 1- \frac{3}{p}} (\bar{\O})} \leq  (C_1 C_p \delta_*  M )^{1/p}   e^{\frac{C_p(1+ (C_1 M)^{1/p})}{p}t}.
		\Ee 
	\unhide

	Finally we use an interpolation between $C^{1,1-D_1}(\bar{\O})$ and $C^{2, 1- \frac{3}{p}}(\bar{\O})$ and derive an estimate of $C^2(\bar{\O})$: Applying Lemma \ref{lemma_interpolation} and (\ref{phi_interpolation}) with $D_2 = 1-\frac{3}{p}$, from (\ref{growth_W1p}) and (\ref{decay_C0,alpha}), we derive that for all $0< D_1<1$, $3<p<6$, $\Lambda_0>0$, and $0 \leq t \leq T$,
		\Be 
		\begin{split}\label{apply_Lemma_interpolation}
			  \| \nabla_x^2 \phi_f(t) \|_\infty 
			\leq & \  e^{- [\lambda_\infty - D_1 \Lambda_0] t} C_{ \O} M  \\
			& + e^{- [(1- \frac{3}{p}) \Lambda_0
			-\frac{1}{p}C_p(1+ L ) 
			] t} 
			(C_1   C_p )^{1/p}  
		 \delta_*  M  .
			%
			%
		\end{split}
		\Ee 
		Then we choose 
		\Be\label{choice_lambda_0_delta_1}
		\Lambda_0 =  
			\frac{\frac{\lambda_\infty}{2} + \frac{C_p}{p}
			(1+ L)}{1- \frac{3}{p}}  \ \  and  \ then    \ D_1= \frac{\lambda_\infty}{2\Lambda_0}.
		\Ee
		In conclusion we have, for all $0 \leq t \leq T$, 
		\[  \| \nabla_x^2 \phi_f(t) \|_\infty 
			\leq
			e^{- \frac{\lambda_\infty}{2} t} [C_\O+(C_1   C_p )^{1/p}   
		 \delta_*  ]M .
		\]
		As long as $M \ll L$ then $ \| \nabla_x^2 \phi_f(t) \|_\infty\leq L$ for all $0 \leq t \leq T$ and hence the claim (\ref{decay_C2}) holds.


		\hide

		Let us choose 
		\Be
		\lambda_0:= \frac{\frac{\lambda}{2} + C(1+ 2M)}{1- 3/p} \ \ \ \text{and} \ \ \ \delta_1:= \frac{\lambda}{2 \lambda_0}.
		\Ee
		Then 
		\Be
		\begin{split}
			\| \nabla_x^2 \phi_f(t) \|_\infty&\leq e^{\delta_1 \lambda_0 t} e^{-\lambda t}M + 
			e^{- (1- \frac{3}{p}) \lambda_0 t}  e^{C(1+ 2M) }\frac{M}{2}\\
			&\leq e^{- \frac{\lambda}{2}}\frac{3M}{2}.
		\end{split}\Ee
		
		\unhide
		\vspace{4pt}
		
		\hide
		
		\textit{Step 4. } We claim that there exist $T>0$ and $k_{0} >0$ such that for all $k\geq  k_{0}$ and for all $(t,x,v) \in [0,T] \times \bar{\O} \times \R^{3}$, we have 
		\Be
		\int_{\prod_{j=1}^{k-1} \mathcal{V}_{j}} \mathbf{1}_{\{ t^{k} (t,x,v,u^{1}, \cdots , u^{k-1}) >0 \}} \dd \Sigma_{k-1}^{k-1} \lesssim_{\O} \Big\{\frac{1}{2}\Big\}^{-k/5}.
		\Ee
		Here we define 
		\Be\label{mathcal_V}
		\mathcal{V}_j: = \{v^j \in \R^3: n(x^j) \cdot v^j >0\},
		\Ee
		and
		\Be\begin{split}
			\label{measure} 
			d\Sigma _{l}^{k-1}(s) &= \{\Pi _{j=l+1}^{k-1}\dd\sigma _{j}\}\times\{
			e^{
				-\int^{t_l}_s
				\nu_\phi 
			}
			\tilde{w}(v_{l})\dd\sigma _{l}\}\times \Pi
			_{j=1}^{l-1}\{{{
					e^{
						-\int^{t_j}_{t_{j+1}}
						\nu_\phi 
					}
					\dd\sigma _{j}}}\}.
		\end{split} \Ee 
		The proof of the claim is a modification of a proof of Lemma 14 of \cite{GKTT1}.
		
		For $0<\delta\ll 1$ we define
		\Be\label{V_j^delta}
		\mathcal{V}_{j}^{\delta} := \{ v^{j } \in \mathcal{V}^{j} : |v^{j} \cdot n(x^{j})| > \delta, \ |v^{j}| \leq \delta^{-1} \}.
		\Ee 
		
		Choose   
		\Be\label{large_T}
		T= \frac{2}{\delta^{2/3} (1+ \| \nabla \phi \|_\infty)^{2/3}}.
		\Ee
		We claim that 
		\Be\label{t-t_lowerbound}
		|t^j -t^{j+1}|\gtrsim   \delta^3, \ \ \text{for} \ v^j \in \mathcal{V}^\delta_j,  \ 0 \leq t\leq T, \ 0 \leq t^j.
		\Ee
		
		For $j \geq 1$
		\Bes
		&&\Big| \int^{t^{j+1}}_{t^{j}} V(s;t^{j}, x^{j}, v^{j}) \dd s   \Big|^{2}\\
		&=& |x^{j+1} -x^{j}|^{2}\\
		&\gtrsim& |(x^{j+1} -x^{j}) \cdot n(x^{j})|\\
		&=&\Big| \int^{t^{j+1}}_{t^{j }}  V(s;t^{j},x^{j},v^{j}) \cdot  n(x^{j})  \dd s 
		\Big|\\
		&=&\Big| \int_{t^{j}}^{t^{j+1}}  
		\Big(
		v^{j} - \int^{s}_{t^{j}} \nabla \phi (\tau, X(\tau;t^j,x^j,v^j))    \dd \tau
		\Big)\cdot   n(x^{j}) 
		\dd s 
		\Big|\\
		&\geq& |v^{j} \cdot n(x^{j})| |t^{j}-t^{j+1}|
		- \Big|
		\int^{t^{j+}}_{t^{j }} \int^{s}_{t^{j }}  \nabla \phi (\tau, X(\tau; t^{j},x^{j},v^{j})) \cdot n(x^{j})\dd \tau \dd s 
		\Big|.
		\Ees
		Here we have used the fact if $x,y \in \p\O$ and $\p\O$ is $C^2$ and $\O$ is bounded then $|x-y|^2\gtrsim_\O |(x-y) \cdot n(x)|$.

		Hence
		\Be\label{lower_tb}
		\begin{split}
			&|v^{j} \cdot n(x^{j})| \\
			& \lesssim \frac{1}{|t^{j} - t^{j+1}| } \Big| \int^{t^{j+1}}_{t^{j}}
			V(s;t^{j},x^{j},v^{j}) \dd s
			\Big|^{2}\\
			& \ \  + 
			\frac{1}{|t^{j} - t^{j+1}| } \Big| \int^{t^{j+1}}_{t^{j}} 
			\int^{s}_{t^{j}} 
			\nabla \phi (\tau, X(\tau;t^{j},x^{j},v^{j}) )\cdot n(x^{j})\dd \tau  \dd s
			\Big|\\
			& \lesssim 
			|t^{j} - t^{j+1}|   \big\{ |v^j|^2 + |t^j - t^{j+1}|^3 \|\nabla \phi\|^2_\infty
			\\
			& \ \ \ \ \ \ \ \  \ \ \ \ \ \  \ \ \ \     +    \frac{1}{2}\sup_{t^{j+1} \leq \tau \leq t^{j}} 
			|  \nabla \phi (\tau, X(\tau;t^{j},x^{j},v^{j}) )\cdot n(x^{j})|
			\big\}.
		\end{split}
		\Ee
		For $v^j \in \mathcal{V}^\delta_j$, $0 \leq t\leq T$, and $t^j\geq 0$,
		\Be\notag
		|v^j \cdot n(x^j)|\lesssim |t^j -t^{j+1}|
		\{
		\delta^{-2} + T^3 \| \nabla \phi \|_\infty^2 + \| \nabla \phi \|_\infty
		\}.
		\Ee
		If we choose $T$ as (\ref{large_T}) then we can prove (\ref{t-t_lowerbound}).

		{\color{red}Complete the proof by modifying the proof of Lemma 14 of \cite{GKTT1}}. 
		
		\vspace{4pt}
		
		\unhide
		

		\vspace{4pt}
		
		\textit{Step 3. } We claim that there exists $T_\infty\gg1$ such that, for $N\in \mathbb{N}$, $t \in [NT_\infty, (N+1)T_\infty]$, and $(N+1)T_\infty \leq T$, 
		\Be\begin{split}
			&\| w_\vartheta f (t) \|_\infty \\
			\leq&  \   e^{-\frac{\nu_0}{2} (t-NT)} \| w_\vartheta f (NT_\infty) \|_\infty \\
			&
			+o(1) \sup_{NT_\infty \leq s \leq t}e^{- \frac{\nu_0}{2}  (t-s)} \| w_\vartheta f(s) \|_\infty
			\\
			&+ C_{T_\infty} \int^t_{NT_\infty} e^{- \frac{\nu_0}{2}  (t-s)} \| f(s) \|_{L^2_{x,v}} \dd s \\
			&+C_{T_\infty}  \int^t_{NT_\infty} e^{- \frac{\nu_0}{2}  (t-s)} \| \nabla \phi_f (s) \|_\infty \dd s.
		\end{split}\label{decay_time_interval}
		\Ee
		For the sake of simplicity we present a proof of (\ref{decay_time_interval}) for $N=0$. The proof for $N >0$ can be easily obtained by considering $f(NT_\infty)$ as an initial datum.  
		
		As (\ref{h}) we define $h(t,x,v): = w_\vartheta f(t,x,v)$. Then $h$ solves (\ref{fell_local}) and (\ref{bdry_local}) with exchanging all $(h^{\ell}, h^{\ell+1}, \phi^\ell)$ to $(h,h, \phi_f)$. We define 
		\Be\label{nu_w}
		\nu_{\phi_f, w_\vartheta} : = \nu(v) + \frac{v}{2} \cdot \nabla_x \phi_f - \frac{\nabla_x \phi_{f} \cdot \nabla_v w_\vartheta}{w_\vartheta}.
		\Ee
		From (\ref{global_T}) and (\ref{Morrey}), for $0 \leq t \leq T$ 
		\Be\begin{split}\label{lower_nu_phi_f,w}
			\nu_{\phi_f,w_\vartheta}\geq&  \ \big\{\nu_0- \frac{\| \nabla \phi_f \|_\infty}{2}  - 2 \vartheta \| \nabla \phi_f \|_\infty  \big\}\langle v\rangle \\
			\geq & \ \big\{\nu_0 - (\frac{1}{2} -2 \vartheta)M \big\}   \langle v\rangle\\
			\geq & \ \frac{\nu_0}{2} \langle v\rangle. 
		\end{split}\Ee
		Then $h$ solves (\ref{duhamel_local}) along the trajectory with deleting all superscriptions of $\ell$ and $\ell+1$ and exchanging $\nu^\ell$ to $\nu_{\phi_f, w_\vartheta}$ and with new $g$
		\Be\label{g_global}
		g: = - v\cdot \nabla \phi_f \sqrt{\mu} + \Gamma (\frac{h}{w_\vartheta}, \frac{h}{w_\vartheta}).
		\Ee
We define a stochastic cycles  
		\Be\notag
		(t_l(t,x,v,v_1, \cdots, v_{l-1}),x_l(t,x,v,v_1, \cdots, v_{l-1})),
		\Ee
		by deleting all superscriptions in (\ref{cycle}) and (\ref{cycle_ell}). Then $h$ has a bound (\ref{stochastic_h^ell}) with deleting all superscriptions of $\ell$ and $\ell+1$.

		For any large $m\gg 1$ we define
		\Be\label{k_m}
		\mathbf{k}_{\varrho,m} (v,u) = \mathbf{1}_{|v-u| \geq \frac{1}{m}, |v| \leq m} \mathbf{k}_\varrho (v,u),
		\Ee
		such that $\sup_v \int_{\R^3} | \mathbf{k}_{\varrho,m} (v,u) - \mathbf{k}_\varrho (v,u) | \dd u \lesssim \frac{1}{m}$, 
		and $|\mathbf{k}_{\varrho,m} (v,u)| \lesssim_m 1$.
		
		Furthermore we split the time interval as, for each $\ell, l$
		\Be\label{time_splitting_l_ell}
		\begin{split}
			&\{\max\{  t _{l+1}, 0  \} \leq s\leq  t  _l\}  \\
			= & \ 
			\{\max\{  t _{l+1}, 0  \} \leq s\leq  t _l- \delta\} 
			\cup \{  t _l - \delta \leq s\leq  t _l\} ,
		\end{split}\Ee
		where we choose a small constant $0<\delta \ll_k 1$ later in (\ref{choice_m_delta}). 
		
		Following (\ref{stochastic_h^ell}) with an extra decomposition of (\ref{k_m}), we derive that
		\Be\label{double_teration}
		\begin{split}
			& |h  (t,x,v)|\\
			\leq & \  O(k) \| e^{- \frac{\nu_0}{2} t} h_0 \|_\infty\\
			& + O(k) \sup_{0 \leq s \leq t} \| e^{- \frac{\nu_0}{2} (t-s)} h  (s) \|_\infty^2\\
			& + O(k)  \int^t_0  \| e^{-\frac{\nu_0}{2} (t-s)}   \nabla   \phi_{\frac{h}{w_{\vartheta}}} (s) \|_\infty \dd s \\
			&+ \left\{O(k)  (\delta + \frac{1}{m})   + \Big\{ \frac{1}{2} \Big\}^{k/5}\right\}  \sup_{0 \leq s \leq t} \| e^{- \frac{\nu_0}{2} (t-s)} h  (s) \|_\infty \\
			& + \ \int^{t-\delta}_{\max\{t _1, 0 \}}
			e^{ - \frac{\nu_0}{2} (t-s)}\\
			&\ \ \ \ \  \ \ \ \   \times \int_{
				\R^3} \mathbf{k}_{\varrho,m} (V(s;t,x,v),u)
			|h  (s, X (s;t,x,v), u)|
			\dd u
			\dd s\\
			& +O(k) \sup_l  \int^{t_l -\delta}_{\max\{ t_{l+1} , 0  \}}   e^{ - \frac{\nu_0}{2} (t-s)}\\
			& \ \ \ \ \  \ \ \ \   \times
			\int_{|u| \leq m} \int_{|v_l| \leq m} 
			|h (s, X (s;t_l,x_l, v_l)  ,u )| \dd v_l
			\dd u
			\dd s .
		\end{split}
		\Ee
		Note that $|h  (s, X (s;t,x,v), u)|$ has the similar bound as
		\hide
		\Be
		\begin{split}
			&|h  (s, X (s;t,x,v), u)|\\
			\leq & \  O(k) \| e^{- \frac{\nu_0}{2} t} h_0 \|_\infty\\
			& + O(k) \sup_{0 \leq s \leq t} \| e^{- \frac{\nu_0}{2} (t-s)} h  (s) \|_\infty^2\\
			& + O(k) \int^s_0  \| e^{-\frac{\nu_0}{2} (t-s)}  \nabla   \phi_{\frac{h}{w}} (s) \|_\infty   \\
			&+ \left\{O(k) \left(\delta + \frac{1}{m}\right) + \Big\{ \frac{1}{2} \Big\}^{k/5}\right\}   \sup_{0 \leq s \leq t} \| e^{- \frac{\nu_0}{2} (t-s)} h  (s) \|_\infty \\
			& + \ \int^{t-\delta}_{\max\{t _1, 0 \}}
			e^{ - \frac{\nu_0}{2} (t-s)}\\
			\ \ \ \ \  \ \ \ \   \times
			\int_{
				|u| \leq m}
			|h  (s, X (s;t,x,v), u)|
			\dd u
			\dd s\\
			& +O(k) \int^{t_l -\delta}_{\max\{ t^{\ell-1} , 0  \}}   e^{ - \frac{\nu_0}{2} (t-s)}\\
			& \ \ \ \ \  \ \ \ \ \  \  \times
			\int_{|u| \leq m} \int_{|v_l| \leq m} 
			|h (s, X (s;t_l,x_l, v_l)  ,u )| \dd v_l
			\dd u
			\dd s .
		\end{split}
		\Ee\unhide
		\Be\label{double_teration_s}
		\begin{split}
			& |h  (s, X(s;t,x,v), u)|\\
			\leq &  \ O(k) \| e^{- \frac{\nu_0}{2} s} h_0 \|_\infty\\
			& + O(k) \sup_{0 \leq s^\prime \leq s} \| e^{- \frac{\nu_0}{2} (s-s^\prime)} h (s^\prime) \|_\infty^2\\
			& + O(k)  \int^s_0 \| e^{-\frac{\nu_0}{2} (s-s^\prime)}  \nabla   \phi_{\frac{h}{w_{\vartheta}}}  (s^\prime) \|_\infty \dd s^\prime \\
			&+ \left\{O(k) \left(\delta + \frac{1}{m}\right) + \Big\{ \frac{1}{2} \Big\}^{k/5}\right\}   \sup_{0 \leq s^\prime \leq s} \| e^{- \frac{\nu_0}{2} (s-s^\prime)} h (s^\prime) \|_\infty \\
			& + \ \int^{s-\delta}_{\max\{t^\prime_1, 0 \}}
			e^{ - \frac{\nu_0}{2} (s-s^\prime)}\\
			&\ \ \ \ \  \ \ \ \   \times
			\int_{
				|u^\prime| \leq m}
			|h  (s^\prime, X (s^\prime;s,X(s;t,x,v),u), u^\prime)|
			\dd u^\prime
			\dd s^\prime\\
			& +O(k) \sup_{l,l^\prime}   \int^{t^\prime_{l^\prime} -\delta}_{\max\{ t^{\prime}_{l^\prime+1}, 0  \}}   e^{ - \frac{\nu_0}{2} (s-s^\prime)}
			\\
			& \ \ \ \ \  \ \ \ \    \times\int_{|u^\prime| \leq m} \int_{|v^\prime_{l^\prime}| \leq m} 
			|h (s^\prime, X(s^\prime;t^\prime_{l^\prime}, x^\prime_{l^\prime}, v^\prime_{l^\prime})  ,u^\prime )| \dd v^\prime_{l^\prime}
			\dd u^\prime
			\dd s^\prime ,
		\end{split}
		\Ee  
		where 
		\Be\notag
		\begin{split} t^\prime_{l^\prime} &= t_{l^\prime} (s, X(s;t,x,v), u,v^\prime_1, \cdots, v^\prime_{l^\prime-1}),\\
			x^\prime_{l^\prime} &= x_{l^\prime} (s, X(s;t,x,v), u,v^\prime_1, \cdots, v^\prime_{l^\prime-1}).
		\end{split}
		\Ee
		
		\hide

		From (\ref{small_k}) for $k\gg 1$
		\Be
		\mathbf{1}_{t_{1}^\ell >0}   
		\frac{
			e^{-  \int^t_{t_1^\ell}     \nu^\ell }
		}{\tilde{w}_{\varrho}(V^\ell (t_1^\ell))}\int_{\prod_{j=1}^{k-1}\mathcal{V}_{j}}  (\ref{h5}) \leq o(1)\max_{ l \geq 0 } \| h^{\ell-l} \|_\infty
		\Ee

		From Lemma 3 in \cite{Guo10} for $\varrho>0$ and $-2\varrho<  \vartheta <2\varrho$ and $\zeta\in\mathbb{R}$, we have \begin{equation} 
		\int_{\mathbb{R}^{3}}
		\mathbf{k} (v,u)
		\frac{\langle
			v\rangle ^{\zeta }e^{\vartheta |v|^{2}}}{\langle u\rangle ^{\zeta }e^{\vartheta
				|u|^{2}}}\mathrm{d}u \ \lesssim \ \langle v\rangle ^{-1}.   \label{int_k}  
		\end{equation}
		Then $|K_w h^{\ell-l}| \leq \int_{\R^3} \mathbf{k} (v,u) \frac{w(v)}{w(u)} \dd u \times \| h^{\ell-l} \|_\infty \lesssim \| h^{\ell-l} \|_\infty$. Together with (\ref{bound_g_ell}) 
		%

		Then $|h(t,x,v)|$ have the same bound of (\ref{h_iteration}), (\ref{h1}), (\ref{h2}), (\ref{h5}).

		. 
		
		We split $\mathbf{k}_w = \mathbf{k}_m + ( \mathbf{k}_w - \mathbf{k}_m)$
		\Be\label{double_teration_ell}
		\begin{split}
			& |h^{\ell+1} (t,x,v)|\\
			\leq & O(k) \| e^{- \frac{\nu_0}{4} t} h_0 \|_\infty\\
			& + O(k)\max_l\sup_{0 \leq s \leq t} \| e^{- \frac{\nu_0}{2} (t-s)} h^{\ell-l} (s) \|_\infty^2\\
			& + O(k) \max_l  \| e^{-\frac{\nu_0}{2} (t-s)}  \nabla   \phi ^{\ell-l} (s) \|_\infty \\
			&+ o(1) \max_l\sup_{0 \leq s \leq t} \| e^{- \frac{\nu_0}{2} (t-s)} h^{\ell-l} (s) \|_\infty \\
			& + \ \int^t_{\max\{t^\ell_1, 0 \}}
			e^{ - \frac{\nu_0}{2} (t-s)}\int_{
				|u| \leq m}
			|h^\ell (s, X^\ell (s;t,x,v), u)|
			\dd u
			\dd s\\
			& +O(k) \sup_l \int^{t_l^{\ell- (l-1)}}_{\max\{ t^{\ell-1}_{l+1}, 0  \}}   e^{ - \frac{\nu_0}{2} (t-s)}
			\int_{|u| \leq m} \int_{|v_l| \leq m} 
			|h^{\ell-l} (s, X^{\ell-l}(s;t_l,x_l, v_l)  ,u )| \dd v_l
			\dd u
			\dd s .
		\end{split}
		\Ee

		We plug this again to achieve 
		\Be\label{double_teration_twice}
		\begin{split}
			& |h (t,x,v)|\\
			\leq & O(k ) \| e^{- \frac{\nu_0}{4} t} h_0 \|_\infty\\
			& + O(k )\max_l\sup_{0 \leq s \leq t} \| e^{- \frac{\nu_0}{2} (t-s)} h  (s) \|_\infty^2\\
			& + O(k) \max_l  \| e^{-\frac{\nu_0}{2} (t-s)}  \nabla   \phi  _f (s) \|_\infty \\
			&+ \{ \frac{1}{N} + \{\frac{1}{2}\}^{-k/5}  \}  \max_l\sup_{0 \leq s \leq t} \| e^{- \frac{\nu_0}{2} (t-s)} h  (s) \|_\infty \\
			& + \ \int^t_{\max\{t_1, 0 \}}
			e^{ - \frac{\nu_0}{2} (t-s)}\int_{
				|u| \leq m}
			\int^s_{\max\{ t_1^\prime, 0 \}} 
			e^{- \frac{\nu_0}{2} (s-s^\prime)} \\
			& \ \ \ \ \ \  \times 
			\int_{|u^\prime| \leq m} 
			|h(s^\prime, X(s^\prime; s, X(s;t,x,v), u), u^\prime)| \dd u^\prime
			\dd s^\prime
			\dd u
			\dd s\\
			& +O(k) \sup_l \int^{t_l^{\ell- (l-1)}}_{\max\{ t^{\ell-1}_{l+1}, 0  \}}   e^{ - \frac{\nu_0}{2} (t-s)}
			\int_{|u| \leq m} \int_{|v_l| \leq m} 
			|h^{\ell-l} (s, X^{\ell-l}(s;t_l,x_l, v_l)  ,u )| \dd v_l
			\dd u
			\dd s 
			%
			%
		\end{split}
		\Ee


		\hide
		
		We define
		\begin{equation}
		\tilde{w} (v)\equiv \frac{1}{ w  (v)\sqrt{\mu (v)}}
		,\label{tweight}
		\end{equation}
		and $\mathcal{V}(x)=\{v \in\mathbb{R}^3 : n(x)\cdot v >0\}$
		with a probability measure $\dd\sigma=\dd\sigma(x)$ on $\mathcal{V}(x)$ which is given by
		\begin{equation}
		\dd\sigma \equiv\mu
		(v)\{n(x)\cdot v\}\dd v.\label{smeasure}
		\end{equation}

		Denote $h=w  f$ and $K_{w }( \ \cdot \ )=w  K( \frac {1}{w } \ \cdot)$.   
		\Be\begin{split}
			&|h (t,x,v)| \\
			\leq & \ \mathbf{1}_{t_{1}\leq 0}e^{-      \int^t_0     \nu_\phi 
			}|h (0,x-t{v},v)|  \notag \\
			&  +
			\int_{\max\{t_1,0\}}^{t}e^{-  \int^t_s     \nu_\phi }|[K_{w }h  +w g](s,X(s;t,x,v), V(s;t,x,v))|\dd s  \notag \\
			&    +\mathbf{1}_{t_{1}>0}   
			\frac{
				e^{-  \int^t_{t_1}     \nu_\phi }
			}{\tilde{w_{\varrho}}(v)}\int_{\prod_{j=1}^{k-1}\mathcal{V}_{j}}|H|  \notag
		\end{split}\Ee
		where $|H|$ is bounded by
		\begin{eqnarray}
		&&\sum_{l=1}^{k-1}\mathbf{1}_{\{t_{l+1}\leq
			0<t_{l}\}}|h (0,
		X(0;t_l,x_l,v_l)
		,V(0;t_l,x_l,v_l))|\dd\Sigma _{l}(0)  \label{h1} \\
		&&+\sum_{l=1}^{k-1}\int_{\max\{ t_{l+1}, 0 \}}^{t_l}\mathbf{1}_{\{t_{l+1}\leq
			0<t_{l}\}} 
		\\
		&& \ \ \ \ \  \ \ \ \ \   \times|[K_{w }h +w   g](s,
		X(s;t_l,x_l,v_l), V(s;t_l,x_l,v_l)
		|\dd \Sigma
		_{l}(s)\dd s  \ \ \ \  \ \ \ \ \ \ \label{h2} \\
		&&+\mathbf{1}_{\{0<t_{k}\}}|h (t_{k},x_{k},v_{k-1})|\dd\Sigma
		_{k-1}(t_{k}),  \label{h5}
		\end{eqnarray}%
		where $\dd\Sigma _{l}^{k-1}(s)$ is defined in (\ref{measure}).
		
		\unhide
		
		\unhide

		Plugging (\ref{double_teration_s}) into (\ref{double_teration}) we conclude that 
		\Be\label{double_teration_double}
		\begin{split}
			& |h  (t,x,v)|\\
			\leq & \  O_m(k)\{ \| e^{- \frac{\nu_0}{2} t} h_0 \|_\infty  +   \sup_{0 \leq s \leq t} \| e^{- \frac{\nu_0}{2} (t-s)} h  (s) \|_\infty^2
			\}\\
			& + O_m(k)  \int^t_0  \| e^{-\frac{\nu_0}{2} (t-s)}   \nabla   \phi_{\frac{h}{w_{\vartheta}}} (s) \|_\infty \dd s \\
			&+ \left\{O_m(k) \delta +  O(k)\frac{1}{m}  + O_m(1)\Big\{ \frac{1}{2} \Big\}^{k/5}\right\}  \sup_{0 \leq s \leq t} \| e^{- \frac{\nu_0}{2} (t-s)} h  (s) \|_\infty \\
			& + \ \int^{t }_{0}
			\int_{
				|u^\prime| \leq m}
			%
			%
			\int^{s-\delta}_{0} e^{ - \frac{\nu_0}{2} (t-s^\prime)} \int_{
				|u| \leq m}\\
				& \ \ \ \ \  \ \ \ \   \times
			|h  (s^\prime, X (s^\prime;s,X(s;t,x,v),u), u^\prime)| \dd u
			\dd s^\prime \dd u^\prime
			\dd s\\
			& +O(k) \sup_l  \int^{t_l -\delta}_{0}   e^{ - \frac{\nu_0}{2} (t-s)} \\
			 & \ \ \ \ \  \ \ \ \
			\times
			\int_{|u| \leq m} \int_{|v_l| \leq m} 
			|h (s, X (s;t_l,x_l, v_l)  ,u )| \dd v_l
			\dd u
			\dd s \\
			& +O (k) \sup_{l,l^\prime}    \int^{t }_{0}
			\int_{
				|u| \leq m} \int^{t^\prime_{l^\prime} -\delta}_{0}   e^{ - \frac{\nu_0}{2} (t-s^\prime)}  \\
			& \ \ \ \ \  \ \ \ \    \times\int_{|u^\prime| \leq m} \int_{|v^\prime_{l^\prime}| \leq m} 
			|h (s^\prime, X(s^\prime;t^\prime_{l^\prime}, x^\prime_{l^\prime}, v^\prime_{l^\prime})  ,u^\prime )| \dd v^\prime_{l^\prime}
			\dd u^\prime
			\dd s^\prime.
		\end{split}\Ee
		
		Choose $T\gg 1$ and $k\gg 1$ in (\ref{small_k}) and (\ref{large_T}). Then we choose 
		\Be
		m =  {k^2} \ \ \text{and} \ \ \delta = \frac{1}{m^3 k},\label{choice_m_delta}
		\Ee
		so that $O_m(k) \delta +  O(k)\frac{1}{m}  + O_m(1) \{ \frac{1}{2} \}^{k/5}  \ll1$.
		
		Note that 
		\Be\label{X_v}
		\begin{split}
			&\frac{\p  X(s ;t _{l },x_{l },v _{l }) }{\p v _{l }}\\
			= &   - ( t _{l  }-s ) \mathrm{Id}_{3 \times 3}\\
			& -
			\int^{s }_{t _{l }} \int^\tau_{t _{l }}  \frac{\p X(\tau^\prime;t _{l }, x _{l }, v _{l })}{\p v _{l }} \cdot \nabla_x 
			\left(          
			\nabla_x \phi_{\frac{h}{w_{\vartheta}}} (\tau^\prime, X(\tau^\prime;t _{l },x _{l },v _{l }))
			\right) \dd \tau^\prime\dd \tau.
		\end{split}
		\Ee

		\hide
		\Be\notag
		\| \nabla_x^2 \phi_f(t) \|_\infty \leq
		C_2 M^{1/p} 
		e^{- \frac{\lambda_\infty}{2}t}   \ \ \text{with} \ C_2 := C_{\frac{\lambda_\infty}{2 \lambda_0}, \O }+ (C_1 C_p \delta_*)^{1/p},
		\Ee\unhide

		Now we use Lemma \ref{est_X_v}. Note that from (\ref{decay_C2}), the condition (\ref{decay_phi_2}) of Lemma \ref{est_X_v} is satisfied with $\Lambda_{2} = \frac{\lambda_\infty}{2}$ and $\delta_{2} = C_2 M$. From Lemma \ref{est_X_v} and (\ref{result_X_v}) we have 
		\Be\label{X_v_l}
		\left|\frac{\p X(\tau^\prime;t _{l },x   _{l  },v _{l })}{\p v _{l }}\right| \leq Ce^{ \frac{4C C_2 M}{(\lambda_\infty)^2}   }|t _{l } -\tau^\prime|.
		\Ee
		From (\ref{X_v_l}) and (\ref{decay_C2}), the second term of RHS in (\ref{X_v}) is bounded by 
		\Be\begin{split}\label{double_int_phi}
			&C C_2 M e^{ \frac{4CC_2 M}{(\lambda_\infty)^2}    }  \int^{t _{l }}_{s }   \int^{t _{l }}_\tau (t _{l } - \tau^\prime) e^{-\frac{\lambda_\infty}{2} \tau^\prime}
			\dd \tau^\prime \dd \tau \\ 
			\leq & \ 
	 \frac{4	C C_2 M }{(\lambda_\infty)^2}	 e^{ \frac{4	C C_2 M }{(\lambda_\infty)^2}}
			|t _{l } -s | ,
		\end{split}\Ee  
		where we have used (\ref{direct_computation}). 
		From our choice of $\lambda_\infty$ in (\ref{condition_lambda_infty}), we have 
		\Be\notag
		 \frac{4	C C_2 M }{(\lambda_\infty)^2}	 e^{ \frac{4	C C_2 M }{(\lambda_\infty)^2}} < \frac{1}{10}. 
		\Ee
		Therefore from (\ref{X_v}), for $0 \leq s  \leq t_l - \delta$
\hide		\Be\label{X_v_estimate}
		\left|\frac{\p  X(s ;t _{l },x _{l },v _{l }) }{\p v _{l }} 
		+ (t _{l }-s^\prime) \mathrm{Id}_{3 \times 3}\right|
		\leq  4Ce^{C \frac{\lambda_\infty}{2} \sqrt{C_2 M^{1/p}}} C_2 \frac{M^{1/p}}{(\lambda_\infty)^2} |t _{l }-s^\prime|.
		\Ee       
		From (\ref{condition_lambda_infty}) and (\ref{choice_m_delta}) we have  
		%
		%
\unhide		\Be\label{lower_jacob_l}
		\begin{split}
		&\det\left(\frac{\p X(s ;t _{l }, x _{l},v _{l })}{\p v _{l }}\right)
		\\
	= 	&  \  \det\left( - (t_l - s ) \text{Id}_{3 \times 3}
	+  o(1)
	\right) \\ 
		\gtrsim & \   |t_l - s |^3\\
		\gtrsim & \   \delta^3.
	\end{split}	\Ee
		We can obtain the exactly same lower bound for both
		\Be
		\begin{split}\notag
		&\det \left( \frac{\p X(s;t^\prime_{l^\prime}, x^\prime_{l^\prime},v^\prime_{l^\prime})}{\p v^\prime_{l^\prime}}\right) \  for \   0 \leq s^\prime \leq s-\delta, \\
		 & \det \left( \frac{\p X(s^\prime;s,X(s;t,x,v), u)}{\p u} \right)   \ for \  0 \leq s^\prime \leq t^\prime_{l^\prime} - \delta.
		\end{split}
		\Ee
		
		Now we apply the change of variables 
		\Bes
		v_l &\mapsto& X(s;t_l, x_l,v_l),\\
		v^\prime_{l^\prime} &\mapsto& X(s;t^\prime_{l^\prime}, x^\prime_{l^\prime},v^\prime_{l^\prime}),\\
		u &\mapsto& X(s^\prime;s,X(s;t,x,v), u), 
		\Ees 
		and conclude (\ref{decay_time_interval}) from (\ref{double_teration_double}) and (\ref{choice_m_delta}).

		Applying (\ref{decay_time_interval}) successively, we achieve that 
		\Be\begin{split}\label{global_decay_N}
			&\| w_{\vartheta} f(t) \|_\infty\\
			\leq & \  e^{\frac{\nu_0}{2}} e^{- \frac{\nu_0}{2} t} \| w_{\vartheta} f(0) \|_\infty 
			+o(1) \frac{e^{\nu_0 T_\infty}}{1- e^{- \nu_0 T_\infty}}
			\sup_{0 \leq s \leq t} e^{- \nu_0 (t-s)} \| w_{\vartheta} f(s) \|_\infty
			\\
			&
			+\underbrace{ C_{T_\infty} e^{\frac{\nu_0}{2}} \int^t_0 e^{- \frac{\nu_0}{2} (t-s)} \| f(s) \|_2 \dd s}_{(\ref{global_decay_N})_{L^2}}\\
			&
    			+  \underbrace{C_{T_\infty} e^{\frac{\nu_0}{2}}\int^t_0 e^{- \frac{\nu_0}{2} (t-s)}  \| \nabla \phi_f (s) \|_\infty \dd s
			}_{(\ref{global_decay_N})_{\phi_f}}
			,
		\end{split}\Ee
		where we have used
		\[
	 e^{\nu_0 T_\infty} \{ 1+ e^{- \nu_0 T_\infty} +  \cdots +  e^{- \nu_0 NT_\infty}\}= \frac{e^{\nu_0 T_\infty}}{1- e^{- \nu_0 T_\infty}}.
		\]
		
		\vspace{4pt}
		
		\textit{Step 4. } 	From Proposition \ref{dlinearl2} and (\ref{completes_dyn}) we have 
		\Be\begin{split}\label{completes_dyn_final}
			& \| e^{\lambda_2 t}f(t)\|_2^2
			+ \| e^{\lambda_2 t} \nabla \phi (t) \|_2^2\\
			&
			+  \int_0^t \| e^{\lambda_2 \tau}  f (\tau)\|_\nu^2 
			+ \| e^{\lambda_2 \tau} \nabla \phi_f(\tau)\|_2^2 
			\mathrm{d} \tau 
			+  \int_0^t | e^{\lambda_2 \tau} f |^2_{2,+ }    \\
			\lesssim & \ \| f_0\|_2^2 +   \| \nabla \phi_{f_0}\|_2^2   .  
		\end{split}\Ee
		Hence 
		\Be\label{global_decay_N_L2}
		\begin{split}
		(\ref{global_decay_N})_{L^2}\lesssim& \   t e^{- \min ( \frac{\nu_0}{2}, \lambda_2) \times  t}
		\{\| f_0\|_2  +   \| \nabla \phi_{f_0}\|_2   \}\\
		\lesssim& \  e^{- \min ( \frac{\nu_0}{4}, \frac{\lambda_2}{2}) \times  t}
		\{\| f_0\|_2  +   \| \nabla \phi_{f_0}\|_2   \}.
		\end{split}\Ee

	Now we consider $(\ref{global_decay_N})_{\phi_f}$. In order to close the estimate in (\ref{global_decay_N}) we need to improve the decay rate of $ \| \nabla \phi_f (s) \|_\infty$. We claim that, for $\theta_{2,r,p}>0$ (which is specified in (\ref{theta_2rp})),
		\Be\label{better_decay_phi_C1}
		\| \nabla_x \phi_f (s) \|_\infty \lesssim e^{- (1+ \theta_{2,r,p}) \lambda_\infty s}
		\{\sup_{t\geq 0} \|e^{\lambda_2 t} f(s) \|_2  + 
		\sup_{t\geq 0} \|e^{\lambda_\infty t} f(s) \|_\infty\}.
		\Ee

		By Morrey's inequality for $\O \subset \R^3$ and $r>3$
		\Be
		\begin{split}\label{Morry_phi_infty}
			\| \nabla_x \phi_f \|_\infty 
			\lesssim   \| \nabla_x \phi_f \|_{C^{0, 1-3/r} (\O)} 
			\lesssim  \| \nabla_x \phi_f \|_{W^{1,r} (\O)}.
		\end{split}\Ee
		Then applying the standard elliptic estimate to (\ref{phi_f}), we get  
		\begin{eqnarray}
		\| \nabla_x \phi_f(t) \|_{W^{1,2} (\O)}
		\lesssim 
		\left\| 
			\int_{\R^3} f(t) \sqrt{\mu}  \dd v
			\right\|_{L^2 (\O)}
		 \lesssim e^{- \lambda_2 t}  \sup_{t\geq 0} \|e^{\lambda_2 t} f(t) \|_2,\label{decay_phi_1r} \\
			\| \nabla_x \phi_f(t) \|_{W^{1,p} (\O)} \lesssim
			\left\| 
			\int_{\R^3} f(t) \sqrt{\mu}  \dd v
			\right\|_{L^p (\O)}
			\lesssim
			e^{- \lambda_\infty t} \sup_{t\geq 0} \|e^{\lambda_\infty t} f(t) \|_\infty
			.\label{decay_phi_1p}
		 \end{eqnarray}
		 
		 Now we use the standard interpolation: For $p>r>3$,
		 \Be\begin{split}\notag
			\| \nabla_x \phi_f \|_{W^{1,r} (\O)}
			\lesssim 
			\| \nabla_x \phi_f(t) \|_{W^{1,2} (\O)}^{\theta_{2,r,p} } \| \nabla_x \phi_f(t)\|_{W^{1,p} (\O)}^{1-\theta_{2,r,p} },
		\end{split}
		\Ee
		for 
		\Be\label{theta_2rp}
		\theta_{2,r,p} : = \frac{\frac{1}{r}- \frac{1}{p}}{\frac{1}{2} - \frac{1}{p}}> \frac{2}{3}\cdot \frac{p-3}{p-2}.
		\Ee
		Then we derive 
		\Be\begin{split}\label{decay_interpolation}
			& \sup_{t\geq 0 }\|e^{  [\theta_{2,r,p} \lambda_2  + (1- \theta_{2,r,p}) \lambda_\infty  ] t } \nabla_x \phi_f (t) \|_\infty \\
			\lesssim & \ \left(\sup_{t\geq 0} \|e^{\lambda_2 t} f(t) \|_2\right)^{\theta_{2,r,p}} \left(\sup_{t\geq 0} \|e^{\lambda_\infty t} f(t) \|_\infty\right)^{1-\theta_{2,r,p}}\\
			\lesssim & \  \sup_{t\geq 0} \|e^{\lambda_2 t} f(t) \|_2  + o(1)
			\sup_{t\geq 0} \|e^{\lambda_\infty t} f(t) \|_\infty .
		\end{split}\Ee
		From our choice (\ref{condition_lambda_infty}) and $0<p-3 \ll 1$,
		\Be\label{decay_phiC1}
		\theta_{2,r,p} \lambda_2  + (1- \theta_{2,r,p}) \lambda_\infty  \geq (1+ \theta_{2,r,p}) \lambda_\infty.
		\Ee
		
		From (\ref{decay_interpolation})
		\Be
		\begin{split}\label{est_gdN_phif}
		(\ref{global_decay_N})_{\phi_f}
			\lesssim \ &\int^t_0 e^{- \frac{\nu_0}{2} (t-s)}  e^{-(1+ \theta_{2,r,p}) \lambda_\infty s}
			 \|e^{\lambda_2 s} f(s) \|_2 \dd s  \\
			 &+o(1) \int^t_0e^{- \frac{\nu_0}{2} (t-s)}    e^{-(1+ \theta_{2,r,p}) \lambda_\infty s}
		  \| e^{\lambda_ \infty s}w_{\vartheta} f(s) \|_\infty 
			\dd s\\
			\lesssim \ & e^{- \min (\frac{\nu_0}{4}, \lambda_\infty) \times t}
		\{	 \| f_0\|_2  +   \| \nabla \phi_{f_0}\|_2 \}
			 \ \  from \ (\ref{completes_dyn_final})\\
			 & + o(1)e^{- \min (\frac{\nu_0}{4}, \lambda_\infty) \times t}\sup_{0 \leq s\leq t} \| e^{\lambda_ \infty s}w_{\vartheta} f(s) \|_\infty .
		\end{split}
		\Ee
		
		Multiplying $e^{\lambda_\infty t}$ and taking $\sup_{t \geq 0}$ to (\ref{global_decay_N}) with $\lambda_\infty \leq \min \left(\frac{\nu_0}{4}, \frac{\lambda_2}{2}\right)$, and from (\ref{global_decay_N_L2}) and (\ref{est_gdN_phif}), we obtain that 
		\Be\label{f_infty_1}
		 \begin{split}
			 \sup_{t \geq 0}e^{\lambda_\infty t}\| w_{\vartheta} f(t) \|_\infty 
			\lesssim & \    \|w_{\vartheta} f(0) \|_\infty + \| f_0\|_2  +   \| \nabla \phi_{f_0}\|_2  \\
			& + o(1) \sup_{0 \leq s\leq t} e^{\lambda_\infty s} \| w_{\vartheta} f(s) \|_\infty.
 		\end{split}
		\Ee
		By absorbing the last (small) term, we conclude that \hide
		
		From (\ref{choice_L_infty})
		\Be\notag
		\begin{split}
			\lambda_\infty- \min \{ {\nu_0} , \lambda_2- \delta \}
			\geq &   -\frac{\lambda_2 - \delta}{2},\\
			\lambda_\infty- 
			\min \big\{
			{\nu_0} 
			,
			\lambda_\infty + \delta \frac{\theta}{1-\theta}
			\big\}\geq &  - \delta \frac{\theta}{1-\theta}.
		\end{split}
		\Ee
		Note that 
		\Be\begin{split}\notag
			\int^t_0 e^{  -\frac{\lambda_2 - \delta}{2} s} \dd s &\lesssim  \frac{1}{\lambda_2- \delta},\\
			\int^t_0 e^{  -    \delta \frac{\theta}{1-\theta} s} \dd s &\lesssim  \frac{1-\theta}{\delta \theta}.
		\end{split}
		\Ee
		
		Using (\ref{f_infty_1}) and choosing $o(1) \ll_{\delta, \theta, \lambda_2} 1$ we conclude that 
		\Be\label{Linfty_L2}
		\sup_{t \geq 0}e^{\lambda_\infty t}\| w_\vartheta f(t) \|_\infty \lesssim \| w_\vartheta f(0) \|_\infty + 
		\sup_{t \geq 0} e^{\lambda_2 t} \| f(t) \|_2.
		\Ee

		From (\ref{Linfty_L2}) and (\ref{completes_dyn_final}) we conclude that, for some $\mathfrak{C}\gg1$,
		 \unhide
		\Be\label{Linfty_final}
		\sup_{0 \leq t \leq T}e^{\lambda_\infty t}\| w_\vartheta f(t) \|_\infty  
		\leq \mathfrak{C}\delta_* M.
		\Ee
		If we choose $\delta_* \ll  1/ \mathfrak{C}$ then by the local existence theorem (Theorem \ref{local_existence}) and continuity of $\| w_\vartheta f(t) \|_\infty$, $\| w_{\tilde{\vartheta}} f (t) \|_p^p+ \|w_{\tilde{\vartheta}} \alpha_{f,\e}^\beta \nabla_{x,v} f (s) \|_p^p
			+ \int^t_0  | w_{\tilde{\vartheta}}\alpha_{f,\e}^\beta \nabla_{x,v} f(s) |_{p,+}^p$, and $\| \nabla_v f(t) \|_{L^3_x(\O) L^{1+ \delta}_v (\R^3)}$, we conclude that $T= \infty$.

		Then the estimates of (\ref{W1p_main}), (\ref{nabla_v f_31}), and (\ref{stability_1+}) are direct consequence of Proposition \ref{prop_W1p}, Lemma \ref{lemma_apply_Schauder}, Proposition \ref{prop_better_f_v}, and Proposition \ref{L1+stability}.\end{proof}

	\appendix

	\section{Auxiliary Results and Proofs}
	\begin{proof}[\textbf{Proof of (\ref{alpha_invariant})}]
		From (\ref{hamilton_ODE}), for $t-\tb(t,x,v)<s\leq t$,
		\Be\begin{split}\notag
			\xb(s,X(s;t,x,v),V(s;t,x,v))  =& \  \xb(t,x,v),\\
			\vb(s,X(s;t,x,v),V(s;t,x,v))  =& \  \vb(t,x,v).
		\end{split}\Ee
		Therefore 
		\Be\notag
		\begin{split}
			& [\p_t + v\cdot \nabla_x - \nabla_x \phi_f \cdot \nabla_v] \alpha_{f,\e}(t,x,v)  \\
			=   & \    \frac{d}{ds} \alpha_{f,\e}(s,X(s;t,x,v),V(s;t,x,v)) \big|_{s=t}\\
			= & \ \frac{d}{ds} \alpha_{f,\e}(t,x,v) \\
			= & \    0.
		\end{split}
		\Ee
		
		From (\ref{tb}) and (\ref{hamilton_ODE}),
		\Be\begin{split}\notag
			\tb(s,X(s;t,x,v), V(s;t,x,v)) 
			= \tb(t,x,v)- (t-s).
		\end{split}\Ee
		Therefore 
		\Be\notag
		\begin{split}
			& [\p_t + v\cdot \nabla_x - \nabla_x \phi_f  \cdot \nabla_v] (t-\tb(t,x,v)) \\
			=   & \    \frac{d}{ds} [s- \tb(s,X(s;t,x,v),V(s;t,x,v))] \big|_{s=t}\\
			= & \ \frac{d}{ds} [t-\tb(t,x,v)]
			=    0.
		\end{split}
		\Ee
		These prove (\ref{alpha_invariant}).
		\hide
		\Be \notag
		\begin{split}
			& [\p_{t} + v\cdot\nabla_{x} - \nabla_{x}\phi\cdot\nabla_{v}] \alpha(t,x,v)  \\
			&= \frac{d}{ds} \alpha(s,X(s;t,x,v), V(s;t,x,v)) \Big\vert_{s=t}  \\
			&= \mathbf{1}_{ t+1\geq \tb(t,x,v)  } (t,x,v) \frac{d}{ds} |n(\xb(s,X(s), V(s)))\cdot \vb(s,X(s), V(s)) | \Big\vert_{s=t}  \\
			&\quad + |n(\xb(t,x,v))\cdot \vb(t,x,v)| \frac{d}{ds} \mathbf{1}_{ s+1 \geq \tb(s,X(s),V(s))  } \Big\vert_{s=t} \\
			&= \mathbf{1}_{ t+1\geq \tb(t,x,v)  } (t,x,v) \frac{d}{ds} |n(t,x,v)\cdot \vb(t,x,v) | \Big\vert_{s=t}  \\
			&\quad + |n(\xb(t,x,v))\cdot \vb(t,x,v)| \frac{d}{ds} \mathbf{1}_{ t+1 \geq \tb(t,x,v) } \Big\vert_{s=t}  \\
			&= 0,
		\end{split}
		\Ee
		where $X(s)=X(s;t,x,v)$ and $V(s) = V(s;t,x,v)$. We have used deterministic properties:
		\Be
		\begin{split}
			n(\xb(s,X(s), V(s)))\cdot \vb(s,X(s), V(s)) &= n(\xb(t,x,v))\cdot \vb(t,x,v) ,  \\
			t - \tb(t,x,v) &= s - \tb(s,X(s), V(s)),
		\end{split}
		\Ee
		in the last step. \unhide
	\end{proof}

	\begin{proof}[\textbf{Proof of (\ref{k_vartheta_comparision})}] The proof follows the argument of Lemma 7 in \cite{Guo10}. Note
		\Be\notag
		\begin{split}
			\mathbf{k}_{  \varrho}(v,u) \frac{e^{\vartheta |v|^2}}{e^{\vartheta |u|^2}} 
			=  \frac{1}{|v-u| } \exp\left\{- {\varrho} |v-u|^{2}  
			-  {\varrho} \frac{ ||v|^2-|u|^2 |^2}{|v-u|^2} + \vartheta |v|^2 - \vartheta |u|^2
			\right\}.
		\end{split}\Ee
		%
		Let $v-u=\eta $ and $u=v-\eta $. Then the exponent equals
		\begin{eqnarray*}
			&&- \varrho|\eta |^{2}-\varrho\frac{||\eta |^{2}-2v\cdot \eta |^{2}}{%
				|\eta |^{2}}-\vartheta \{|v-\eta |^{2}-|v|^{2}\} \\
			&=&-2 \varrho |\eta |^{2}+ 4 \varrho v\cdot \eta - 4 \varrho\frac{|v\cdot
				\eta |^{2}}{|\eta |^{2}}-\vartheta \{|\eta |^{2}-2v\cdot \eta \} \\
			&=&(-2 \varrho-\vartheta  )|\eta |^{2}+(4 \varrho+2\vartheta )v\cdot \eta -%
			4 \varrho\frac{\{v\cdot \eta \}^{2}}{|\eta |^{2}}.
		\end{eqnarray*}%
		If $0<\vartheta <4 \varrho$ then the discriminant of the above quadratic form of 
		$|\eta |$ and $\frac{v\cdot \eta }{|\eta |}$ is 
		\begin{equation*}
		(4 \varrho+2\vartheta )^{2}-4
		(-2 \varrho-\vartheta  )(-%
		4 \varrho)
		=4\vartheta ^{2}- 16 \varrho \vartheta<0.
		\end{equation*}%
		Hence, the quadratic form is negative definite. We thus have, for $%
		0<\tilde{\varrho}< \varrho - \frac{\vartheta}{4}  $, the following perturbed quadratic form is still negative definite 
		\[
		-(\varrho - \tilde{\varrho})|\eta |^{2}-(\varrho - \tilde{\varrho})\frac{||\eta
			|^{2}-2v\cdot \eta |^{2}}{|\eta |^{2}}-\vartheta \{|\eta |^{2}-2v\cdot \eta \}  \leq 0.
		\]
		Therefore we conclude (\ref{k_vartheta_comparision}).
		\hide, for given $|v|\geq 1,$ we make another change of variable $\eta
		_{\shortparallel }=\{\eta \cdot \frac{v}{|v|}\}\frac{v}{|v|},$ and $\eta
		_{\perp }=\eta -\eta _{||}$ so that $|v\cdot \eta |=|v||\eta
		_{\shortparallel }|$ and $|v-v^{\prime }|\geq |\eta _{\perp }|.$ We can
		absorb $\{1+|\eta |^{2}\}^{|\beta |}$, $|\eta |\{1+|\eta |^{2}\}^{|\beta |}$
		by $e^{\frac{C_{\theta }}{4}|\eta |^{2}}$, and bound the integral in (\ref%
		{wk}) by (\ref{exponent}): 
		\begin{eqnarray*}
			&&C_{\beta }\int_{\mathbf{R}^{2}}(\frac{1}{|\eta _{_{\perp }}|}+1)e^{-\frac{%
					C_{\theta }}{4}|\eta |^{2}}\left\{ \int_{-\infty }^{\infty }e^{-C_{\theta
				}|v|\times |\eta _{||}|}d|\eta _{||}|\right\} d\eta _{\perp } \\
			&\leq &\frac{C_{\beta }}{|v|}\int_{\mathbf{R}^{2}}(\frac{1}{|\eta _{_{\perp
				}}|}+1)e^{-\frac{C_{\theta }}{4}|\eta _{\perp }|^{2}}\left\{ \int_{-\infty
			}^{\infty }e^{-C_{\theta }|y|}dy\right\} d\eta _{\perp }\text{ \ \ }%
			(y=|v|\times |\eta _{||}|).
		\end{eqnarray*}%
		We thus deduce our lemma since both integrals are finite.

		\unhide
	\end{proof}

	Recall $\kappa_\delta(x,v)$ in (\ref{Z_dyn}). Let us denote $f_{\delta}(t,x,v) 
	:=      \kappa_\delta (x,v) f(t,x,v)$. We assume that $f(s,x,v)=e^s f_0(x,v)$ for $s<0$. Then 
	\Be\begin{split}\notag
	\| f_{\delta} \|_{L^{2} (\mathbb{R} \times \Omega \times \mathbb{R}^{3})} 
	&\lesssim   \| f \|_{L^{2} (\mathbb{R}_{+} \times \Omega \times \mathbb{R}^{3})}
	+ \| f_{0}\|_{L^{2} (\Omega \times \mathbb{R}^{3})} , \\
	 \| f_{\delta} \|_{L^{2} ( \mathbb{R} \times \gamma)} & \lesssim  \| f_{\gamma} \|_{L^{2} ( \mathbb{R}_{+} \times \gamma)} + \| f_{0} \|_{L^{2} (\gamma)}.
	\end{split}\Ee
	
	\hide
	Note that, at the boundary $(x,v) \in \gamma:=\partial\Omega \times \mathbb{R}^{3}$,  
	\begin{equation}\label{Z_support_dyn}
	f_{\delta}(t,x,v)|_{\gamma}\equiv 0, \ \  \text{for}   \ |n(x) \cdot v| \leq \delta \  \text{ or } \  |v|\geq \frac{1}{\delta}  .
	\end{equation}
	\unhide
	
	\hide

	The main goal of this section is the following:
	\begin{proposition}\label{prop_3} Assume $g \in L^{2} (\mathbb{R}_{+} \times \Omega \times \mathbb{R}^{3})$, $f_{0} \in L^{2} (\Omega \times \mathbb{R}^{3})$, and $f_{\gamma}\in L^{2} (  \mathbb{R}_{+}\times \gamma)$. Let $f \in L^{\infty}(  \mathbb{R}_{+}; L^{2} (\Omega \times \mathbb{R}^{3}))$ solves (\ref{linear_dyn}) in the sense of distribution and satisfies $f(t,x,v)   =   f_{\gamma}(t,x,v) $ on $ \mathbb{R}_{+} \times \gamma$ and $f(0,x,v)  = f_{0} (x,v)$ on $\Omega \times \mathbb{R}^{3}.$
		
		Then
		\begin{equation}\label{S1}
		\begin{split}
		|a  (t,x)| + |b(t,x)| + |c(t,x)|  
		\leq  \ \mathbf{S}_{1}f(t,x) + \mathbf{S}_{2} f(t,x),\\
		\mathbf{S}_{1} f(t,x) : =   \ 4 \int_{\mathbb{R}^{3}} 
		| f_{\delta} (t,x,v)|
		\langle v\rangle^{2} \sqrt{\mu(v)}\dd v,  \\
		\mathbf{S}_{s} f(t,x) : =    4 \int_{\mathbb{R}^{3}}
		| (\mathbf{I} - \mathbf{P}) f (t,x,v)| \langle v\rangle ^{2}\sqrt{\mu(v)} \dd v
		+ 2\chi(t) \int_{\mathbb{R}^{3}} |f_{0} (x,v)| \langle v \rangle^{2} \sqrt{\mu(v)} \dd v
		,
		\end{split}
		\end{equation} 
		where $f_{\delta}$ is defined in (\ref{Z_dyn}). 
		
		Moreover
		\begin{equation}\label{decom_Pf}
		\begin{split}
		\| \mathbf{S}_{1}f \|_{L^{3}_{x} L^{2}_{t}}   & \ \lesssim   \ \| w^{-1} f\|_{L^{2}_{t,x,v}}+  \| w^{-1} g\|_{L^{2}_{t,x,v} }+ \|   f\|_{L^{2}(\mathbb{R}_{+} \times\gamma)} ,\\
		\| \mathbf{S}_{2}f \|_{L^{2}_{t,x} }  & \  \lesssim   \  \| (\mathbf{I}- \mathbf{P}) f\|_{L^{2}_{t,x,v} }
		+ \| f_{0} \|_{L^{2}_{x,v}}
		,
		\end{split}
		\end{equation}
		for $w= e^{\beta |v|^{2}}$ with $0< \beta\ll 1$.
		
	\end{proposition}
	
	\bigskip
	
	We need several lemmas to prove Proposition \ref{prop_3}.

	\unhide

	\begin{lemma} \label{extension_dyn}Assume $\O$ is convex in (\ref{convexity_eta}) and $\sup_{0 \leq t \leq T}\|E(t)\|_{L^\infty (\O)} < \infty$. Let $\bar{E}(t,x) = \mathbf{1}_{\O}(x) E(t,x)$ for $x \in \R^3$. There exists $\bar{f}(t,x,v) \in L^{2}( \mathbb{R} \times   \mathbb{R}^{3} \times \mathbb{R}^{3})$, an extension of $f_{\delta}$, such that 
		\begin{equation}\notag
		\bar{f}  |_{\Omega \times \mathbb{R}^{3}}\equiv f_{\delta}    \  \text{ and } \  \bar{f}  |_{\gamma}\equiv f_{ \delta} |_{\gamma}   \  \text{ and } \ \bar{f } |_{t=0} \equiv f_{\delta} |_{t=0}.
		\end{equation} Moreover, in the sense of distributions on $\mathbb{R} \times \mathbb{R}^{3} \times \mathbb{R}^{3}$,
		\begin{equation}\label{eq_barf_dyn}
		[\partial_{t} +  v\cdot \nabla_{x} + \bar{E} \cdot \nabla_{v}]\bar{f} =  h_{1} + h_{2}  ,
		\end{equation}
		where
		\Be\begin{split}\label{barf_h}
			h_{1} (t,x,v)=& \ 
			\kappa_\delta(x,v)    \mathbf{1}_{t \in [0,\infty)} 
			[ \partial_{t} + v\cdot \nabla_{x}   + E
			\cdot \nabla_{v}  ] f\\
			&
			+ \kappa_\delta(x,v)  \mathbf{1}%
			_{t \in ( - \infty, 0 ]} e^t
			[1  + v\cdot \nabla_{x}   + E
			\cdot \nabla_{v}] f_0 \kappa_\delta(x,v) \\
			& 
			+ f (t,x,v) [v\cdot \nabla_{x} + E
			\cdot \nabla_{v}] \kappa_\delta(x,v),\\
			h_2 (t,x,v) = & \  \mathbf{1}_{t-\tb^{EX}(t,x,v)\leq 0}  
			\mathbf{1}_{\xb^{EX} (t,x,v) \in \p \O}   
			e^{t-\tb^{EX} (x,v)} f_0 (\xb^{EX} (x,v),v) \\
			&+    \mathbf{1}_{t-\tb^{EX}(t,x,v)\leq 0}  
			\mathbf{1}_{\xf^{EX} (t,x,v) \in \p \O}   
			e^{t-\tf^{EX} (x,v)} f_0 (\xf^{EX} (x,v),v),
		\end{split}\Ee
		where $t_{\mathbf{b}}^{EX}, x_{\mathbf{b}}^{EX}, t_{\mathbf{f}}^{EX}, x_{\mathbf{f}}^{EX}$ are defined in (\ref{def_tb_EX}).
		
		Moreover,
		\Be\begin{split}\label{estimate_h1_h2}
			\| h_{1} \|_{  L^{2}( \mathbb{R} \times   \mathbb{R}^{3} \times \mathbb{R}^{3})}  
			\lesssim  & \  \|  [ \partial_{t} + v\cdot \nabla_{x}   + E
			\cdot \nabla_{v}  ] f\|_{L^{2}(
				\mathbb{R}_{+} \times 
				\Omega \times \mathbb{R}^{3})} \\
			&  
			+  \| f \|_{L^{2} (\R \times \Omega \times \mathbb{R}^{3})} 
			+ \| [ v\cdot \nabla_{x} + E\cdot \nabla_{v} ] f_{0} \|_{L^{2 } (\Omega \times \mathbb{R}^{3})}
			,\\
			\| h_{2} \|_{L^{2} ( \mathbb{R} \times \mathbb{R}^{3} \times \mathbb{R}^{3})}\lesssim & \ 
			\| f_0 \|_{L^2(\gamma)}.
		\end{split}\Ee
		
	\end{lemma}
	
	\begin{proof} In the sense of distributions  
		\begin{equation} \label{eqtn_f_delta}
		\partial_{t} f_{\delta}+ v\cdot \nabla_{x} f_{\delta} + E
		\cdot \nabla_{v} f_{\delta} = h_1 \  \  in  \  (\ref{barf_h}).
		\end{equation}
		Clearly $| [v\cdot \nabla_{x} + E
		\cdot \nabla_{v}] \kappa_\delta(x,v)| \lesssim_\delta 1$.\hide
		\begin{eqnarray}
		&&\Big|\{v\cdot \nabla_{x} + \e^{2} \Phi \cdot \nabla_{v}\} [1-\chi(\frac{%
			n(x) \cdot v}{\delta}) \chi \big( \frac{ \xi(x )}{\delta}\big) ]
		\chi(\delta|v|)\Big|  \label{der_chi} \\
		&=&\Big| - \frac{1}{\delta} \{v\cdot \nabla_{x} n(x) \cdot v + \e^{2} \Phi
		\cdot n(x) \} \chi^{\prime} \big(\frac{n(x) \cdot v}{\delta} \big) \chi %
		\big( \frac{ \xi(x )}{\delta}\big) \chi(\delta|v|)  \notag \\
		&& - \ \frac{1}{\delta} v\cdot \nabla_{x} \xi(x) \chi^{\prime} \big( \frac{%
			\xi(x)}{\delta}\big) \chi (\frac{n(x) \cdot v}{\delta}) \chi(\delta|v|) + \e%
		^{2}\delta \Phi \cdot \frac{v}{|v|} \chi^{\prime} (\delta|v|)[1-\chi(\frac{%
			n(x) \cdot v}{\delta}) \chi \big( \frac{\xi(x)}{\delta}\big) ] \Big|  \notag
		\\
		&\leq& \frac{4}{\delta}( |v|^{2}\|\xi\|_{C^2} + \e^{2}\|\Phi\|_\infty )
		\chi(\delta|v|) + \frac{C_{\Omega}}{\delta} |v|\chi(\delta|v|) + \e^{2}
		\delta\|\Phi\|_\infty \mathbf{1}_{|v| \leq {2}{\delta}^{-1}}  \notag \\
		&\lesssim & {\delta^{-3}} \mathbf{1}_{|v| \leq 2 \delta^{-1}}.  \notag
		\end{eqnarray}\unhide

		For $x \in \R^3 \backslash \bar{\O}$ we define
		\Be\begin{split}\label{def_tb_EX}
			\tb^{EX}(t,x,v) &:= \sup\{s \geq 0: x-(t-\tau) v \in \R^3 \backslash \bar{\O}
			\  \  for \ all   \ \tau \in (t-s,t)
			\},\\
			\tf^{EX}(t,x,v) &:= \sup\{s \geq 0: x-(t-\tau) v \in \R^3 \backslash \bar{\O}
			\ \  for \  all   \ \tau \in (t ,t+s)
			\},
		\end{split}\Ee
		and $\xb^{EX}(t,x,v) = x- (t-\tb^{EX}(t,x,v))v$, $\xf^{EX}(t,x,v) = x- (t-\tf^{EX}(t,x,v))v$.
		
		We define, for $x \in \R^3 \backslash \bar{\O}$,
		\Be\label{def_f_E}
		\begin{split}
			f_E (t,x,v) =& \mathbf{1}_{\xb^{EX} (t,x,v) \in \p \O} f_\delta(t-\tb^{EX}(t,x,v), \xb^{EX} (t,x,v),v)\\
			+& \mathbf{1}_{\xf^{EX} (t,x,v) \in \p \O}f_\delta(t-\tf^{EX}(t,x,v), \xf^{EX} (t,x,v),v).
		\end{split}\Ee
		Recall that, from (\ref{Z_dyn}), $f_\delta\equiv0$ when $v=0$ and hence $f_E\equiv0$ for $v=0$. Since $\O$ is convex if $v\neq 0$ then $\{\xb^{EX} (t,x,v) \in \p \O\} \cap \{\xf^{EX} (t,x,v) \in \p \O\}= \emptyset$. Note that 
		\Be\label{no_jump_bdry}
		f_E(t,x,v) = f_\gamma(t,x,v) = f_\delta (t,x,v) \ \  for \   x \in \p\O.
		\Ee
		In the sense of distribution, in $\R \times [\R^3 \backslash \bar{\O}] \times \R^3$ 
		\Be
		\label{eqtn_f_E}
		\p_t f_E + v\cdot \nabla_x f_E = h_2   \text{ in }  (\ref{barf_h}).
		\Ee

		We define 
		\Be\label{def_bar_f}
		\bar{f}(t,x,v) : = \mathbf{1}_{\O} (x)  f_\delta (t,x,v)
		+ \mathbf{1}_{\R^3 \backslash \bar{\O}} (x) f_E (t,x,v).
		\Ee
		From (\ref{eqtn_f_delta}), (\ref{no_jump_bdry}), and (\ref{eqtn_f_E}) we prove (\ref{eq_barf_dyn}). The estimates of (\ref{estimate_h1_h2}) are direct consequence of Lemma \ref{lem COV} and Lemma \ref{lem COV2}.
		\hide

		This proves the second line of (\ref{force_Z_dyn}). Since $\big[1-\chi(\frac{%
			n(x) \cdot v}{\delta}) \chi ( \frac{\xi(x)}{\delta}) \big] \chi(\delta|v|)
		\leq 1$, we prove the first line of (\ref{force_Z_dyn}) directly. 
		
		\vspace{4pt}
		
		\noindent\textit{Step 2. } We claim that if $0 \leq \xi (x) \leq \tilde{C}
		\delta^{4}, \ |n(x) \cdot v|> \delta$ and $|v| \leq \frac{1}{\delta}$ then
		either $\xi(\tilde{x}_{\mathbf{f}} (x,v))=\tilde{C}\delta^{4}$ or $\xi(%
		\tilde{x}_{\mathbf{b}}(x,v)) =\tilde{C}\delta^{4}$.
		
		To show this, if $v \cdot n(x) \geq \delta$, we take $s>0$, while if $v\cdot
		n(x) \leq - \delta$ then we take $s<0$. From (\ref{intfor}), 
		\begin{eqnarray}
		&&\xi(X(s; 0, x, v) ) = \xi(x) + \int_{0}^{s} V( \tau; 0, x, v) \cdot
		\nabla_{x} \xi( X( \tau; 0, x, v) ) \mathrm{d}\tau  \notag \\
		&=&\xi(x) + \int^{s}_{0} \{ v+ O(1)\e^{2} \|\Phi\|_{\infty} \tau \} \cdot \{
		\nabla_{x} \xi(x) + O(1)\|\xi\|_{C^{2}} (|v| + \e^{2} \|\Phi\|_{\infty}
		\tau) \tau \}  \notag \\
		&=&\xi(x) + v \cdot \nabla_{x} \xi(x) s +O(1) \|\xi\|_{C^{2}} \big\{ |v|^{2}
		s^{2} + \e^{2} \|\Phi \|_{\infty} s ^{2} + \e^{2} \|\Phi\|_{\infty}|v|s^{3}
		+ \e^{4} \| \Phi \|_{\infty}^{2} s^{4} \big\}.  \notag
		\end{eqnarray}
		From $\xi(x)\ge 0$, 
		\begin{eqnarray}
		\xi(X(s; 0, x, v) ) &\geq & \delta |s| \Big\{ 1 - \frac{ \|\xi\|_{C^{2}} }{%
			\delta}\Big[ |v|^{2} |s| + \e^{2} \|\Phi \|_{\infty} |s| + \e^{2}
		\|\Phi\|_{\infty}|v\|s|^{2} + \e^{4} \| \Phi \|_{\infty}^{2} |s|^{3} \Big] %
		\Big\}  \notag \\
		&\ge & \delta |s| \Big\{ 1 - \frac{ \|\xi\|_{C^{2}} }{\delta}\Big[\frac{1}{%
			\delta^{2}} |s| + \e^{2} \|\Phi \|_{\infty} | s| + \e^{2} \|\Phi\|_{\infty} 
		\frac{1}{\delta}|s|^{2} + \e^{4} \| \Phi \|_{\infty}^{2} |s|^{3} \Big] \Big\}
		\notag \\
		&\geq& \delta |s| \Big\{ 1- \Big[ \frac{1}{4} + \frac{ \e^{2} \delta^{2} \|
			\Phi \|_{\infty}}{4 } + \frac{\e^{2} \delta^{4} \| \Phi \|_{\infty}}{16} + 
		\frac{\e^{4}\delta^{8} \| \Phi \|_{\infty}^{2}}{64} \Big] \Big\} \ \geq \ 
		\frac{\delta |s| }{2},  \label{xi_est}
		\end{eqnarray}
		for $0 \leq |s| \leq \frac{\delta^{3}}{ 4(1+ \|\xi \|_{C^{2}})}$ and $0< \e %
		\ll 1$. Therefore 
		\begin{eqnarray*}
			\xi(Y(s;0,x,v)) = \xi(X(\frac{s}{\e}; 0,x,v )) \geq \delta \frac{|s|}{2\e} ,
		\end{eqnarray*}
		for all $0\leq |s| \leq \frac{ \e\delta^{3}}{ 4(1+ \|\xi \|_{C^{2}})}$ with $%
		0< \varepsilon \ll1$. Especially with $\e s_{*} = +\frac{ \varepsilon%
			\delta^{3}}{4(1+ \| \xi \|_{C^{2}})}$ for $n(x) \cdot v >\delta$ and $\e %
		s_{*} = -\frac{\varepsilon\delta^{3} }{4(1+ \| \xi \|_{C^{2}})}$ for $n(x)
		\cdot v <\delta$, 
		\begin{equation*}
		\xi (Y(s_{*};0,x,v)) > \tilde{C} \delta^{4}.
		\end{equation*}
		Therefore, by the intermediate value theorem, we prove our claim.
		
		\vspace{10pt}
		
		\noindent\textit{Step 3. } We define $f_{E}(t,x,v)$ for $(x,v) \in [ \mathbb{%
			R}^{3} \backslash \bar{\Omega}] \times \mathbb{R}^{3}$: 
		\begin{equation}  \label{ZE_dyn}
		\begin{split}
		f_{E}(t,x,v) & \ := \ f_{\delta}(t- \e t_{\mathbf{b}}^{*}(x,v),x_{\mathbf{b}%
		}^{*}(x,v), v_{\mathbf{b}}^{*}(x,v)) \chi\big( \frac{\xi(x)}{\tilde{C}
			\delta^{4}} \big) \chi \big( {t_{\mathbf{b}}^{*} (x,v)} \big) \ \ \ \text{if}
		\ \ x_{\mathbf{b}}^{*}(x,v) \in \partial\Omega, \\
		& \ := \ f_{\delta}(t+ \e t_{\mathbf{f}}^{*}(x,v),x_{\mathbf{f}}^{*}(x,v),
		v_{\mathbf{f}}^{*}(x,v)) \chi\big( \frac{\xi(x)}{\tilde{C} \delta^{4}} \big) %
		\chi \big( {t_{\mathbf{f}}^{*} (x,v)} \big) \ \ \ \ \text{if} \ \ x_{\mathbf{%
				f}}^{*}(x,v) \in \partial\Omega, \\
		& \ := \ 0 \ \ \ \ \ \ \ \ \ \ \ \ \ \ \ \ \ \ \ \ \ \ \ \ \ \ \ \ \ \ \ 
		\text{if} \ \ x_{\mathbf{b}}^{*}(x,v) \notin \partial\Omega \ \text{and} \
		x_{\mathbf{f}}^{*}(x,v) \notin \partial\Omega.
		\end{split}%
		\end{equation}
		
		We check that $f_{E}$ is well-defined. It suffices to prove the following: 
		\begin{equation*}
		\begin{split}
		& \text{If} \ \ x_{\mathbf{b}}^{*}(x,v) \in \partial\Omega \text{ and } x_{%
			\mathbf{f}}^{*}(x,v) \in \partial\Omega \\
		& \ \text{then} \ \ f_{\delta}( t- \e t_{\mathbf{b}}^{*} (x,v), x_{\mathbf{b}%
		}^{*}(x,v), v_{\mathbf{b}}^{*}(x,v)) \chi\big( \frac{\xi(x)}{\tilde{C}
			\delta^{4}} \big ) = 0 = f_{\delta}( t+ \e t_{\mathbf{f}}^{*} (x,v), x_{%
			\mathbf{f}}^{*}(x,v), v_{\mathbf{f}}^{*}(x,v)) \chi\big( \frac{\xi(x)}{%
			\tilde{C} \delta^{4}} \big) .
		\end{split}%
		\end{equation*}
		If $|n(x_{\mathbf{b}}^{*}(x,v)) \cdot v_{\mathbf{b}}^{*}(x,v)| \leq \delta$
		or $|v_{\mathbf{b}}^{*}(x,v)| \geq \frac{1}{\delta}$ then $f_{\delta}( t- \e %
		t^{*}_{\mathbf{b}}(x,v), x_{\mathbf{b}}^{*}(x,v), v_{\mathbf{b}}^{*}(x,v))=0$
		due to (\ref{Z_support_dyn}). If $n(x_{\mathbf{b}}^{*}(x,v)) \cdot v_{%
			\mathbf{b}}^{*}(x,v) > \delta$ and $|v_{\mathbf{b}}^{*}(x,v)|\leq \frac{1}{%
			\delta}$ then, due to \textit{Step 2}, $\xi(x_{\mathbf{f}}^{*}(x,v)) = \xi(
		x_{\mathbf{f}}^{*}( x_{\mathbf{b}}^{*}(x,v), v_{\mathbf{b}}^{*}(x,v) ) ) = 
		\tilde{C}\delta^{4}$ so that $x_{\mathbf{f}}^{*}(x,v) \notin \partial\Omega$.
		
		On the other hand, if $|n(x_{\mathbf{f}}^{*}(x,v)) \cdot v_{\mathbf{f}%
		}^{*}(x,v)| \leq \delta$ or $|v_{\mathbf{f}}^{*}(x,v)|\geq \frac{1}{\delta}$
		then $f_{\delta}( t+ \e t^{*}_{\mathbf{f}}(x,v) ,x_{\mathbf{f}}^{*}(x,v), v_{%
			\mathbf{f}}^{*}(x,v))=0$ due to (\ref{Z_support_dyn}). If $n(x_{\mathbf{f}%
		}^{*}(x,v)) \cdot v_{\mathbf{f}}^{*}(x,v) <-\delta$ and $|v_{\mathbf{f}%
		}^{*}(x,v)|\leq \frac{1}{\delta}$ then, due to \textit{Step 2}, $\xi(x_{%
			\mathbf{b}}^{*}(x,v)) = \xi( x_{\mathbf{b}}^{*}( x_{\mathbf{f}}^{*}(x,v), v_{%
			\mathbf{f}}^{*}(x,v) ) ) = \tilde{C}\delta^{4}$ so that $x_{\mathbf{b}%
		}^{*}(x,v) \notin \partial\Omega$.
		
		Note that 
		\begin{equation}  \label{fE=fd_dyn}
		f_{E}(t,x,v) = f_{\delta} (t,x,v) \ \ \ \text{for all } x \in \partial\Omega.
		\end{equation}
		If $x \in \partial\Omega$ and $n(x) \cdot v >\delta$ then $(x_{\mathbf{b}%
		}^{*}(x,v), v_{\mathbf{b}}^{*}(x,v)) = (x ,v)$. From the definition (\ref%
		{ZE_dyn}), for those $(x,v)$, we have $f_{E}(t,x,v) = f_{\delta} (t,x,v)$.
		If $x \in \partial\Omega$ and $n(x) \cdot v< -\delta$ then $(x_{\mathbf{f}%
		}^{*}(x,v), v_{\mathbf{f}}^{*}(x.v)) = (x,v)$. From the definition (\ref%
		{ZE_dyn}), we conclude (\ref{fE=fd_dyn}) again. Otherwise, if $-\delta< n(x)
		\cdot v< \delta$ then $f_{E}|_{\partial\Omega} \equiv 0 \equiv
		f_{\delta}|_{\partial\Omega}$.
		
		\vspace{10pt}
		
		\noindent\textit{Step 4. } We claim that $f_{E}(x,v) \in L^{2}([\mathbb{R}%
		^{3} \backslash \bar{\Omega}] \times \mathbb{R}^{3})$.
		
		From the definition of (\ref{ZE_dyn}), we have $f_{E}(x,v) \equiv 0$ if $x_{%
			\mathbf{b}}^{*} (x,v) \notin \partial\Omega$ and $x_{\mathbf{f}}^{*} (x,v)
		\notin \partial\Omega$. Therefore, from (\ref{ZE_dyn}), 
		\begin{eqnarray}
		&&\int_{-\infty}^{\infty} \iint_{[ \mathbb{R}^{3 } \backslash \Omega ]\times 
			\mathbb{R}^{3}} |f_{E} (t,x,v)|^{2} \mathrm{d} x \mathrm{d} v \mathrm{d} t 
		\notag \\
		&=&\int_{-\infty}^{\infty} \iint_{ [ \mathbb{R}^{3 } \backslash \Omega ]
			\times \mathbb{R}^{3}} \mathbf{1}_{x_{\mathbf{b}}^{*} (x,v) \in
			\partial\Omega} |f_{E} |^{2} + \int_{-\infty}^{\infty} \iint _{ [ \mathbb{R}%
			^{3 } \backslash \Omega ] \times \mathbb{R}^{3}} \mathbf{1}_{x_{\mathbf{f}%
			}^{*} \in \partial\Omega} |f_{E} |^{2}  \notag \\
		&=& \int_{-\infty}^{\infty} \iint_{ [ \mathbb{R}^{3 } \backslash \Omega
			]\times \mathbb{R}^{3}} \mathbf{1}_{x_{\mathbf{b}}^{*} (x,v) \in
			\partial\Omega} |f_{\delta} ( t- \e t^{*}_{\mathbf{b}} , x_{\mathbf{b}}^{*}
		, v_{\mathbf{b}}^{*} )|^{2} \big| \chi \big( \frac{\xi(x)}{\tilde{C}
			\delta^{4}} \big) \big|^{2} |\chi ( {t_{\mathbf{b}}^{*} } ) |^{2}\mathrm{d}
		x \mathrm{d} v \mathrm{d} t  \label{fE1_dyn} \\
		&&+ \int_{-\infty}^{\infty} \iint_{ [ \mathbb{R}^{3 } \backslash \Omega ]
			\times \mathbb{R}^{3}} \mathbf{1}_{x_{\mathbf{f}}^{*} (x,v) \in
			\partial\Omega} |f_{\delta} ( t+ \e t^{*}_{\mathbf{f}} , x_{\mathbf{f}}^{*}
		, v_{\mathbf{f}}^{*} )|^{2} \big| \chi \big( \frac{\xi(x)}{\tilde{C}
			\delta^{4}} \big) \big|^{2} |\chi ( {t_{\mathbf{f}}^{*} } ) |^{2} \mathrm{d}
		x \mathrm{d} v \mathrm{d} t ,  \label{fE2_dyn}
		\end{eqnarray}
		where $(t^{*}_{\mathbf{b}},x^{*}_{\mathbf{b}},v^{*}_{\mathbf{b}} )$ and $%
		(t^{*}_{\mathbf{f}},x^{*}_{\mathbf{f}},v^{*}_{\mathbf{f}} )$ are evaluated
		at $(x,v)$.
		
		By (\ref{bdry_int2}), 
		\begin{eqnarray*}
			(\ref{fE1_dyn}) &\leq& \int^{\infty}_{- \infty} \mathrm{d} t
			\int_{\partial\Omega} \int_{n(x) \cdot v >0} \int_{0}^{ \min\{ t_{\mathbf{f}%
				}^{*}(x,v), 1 \}}\mathrm{d} S_{x} \mathrm{d} v \mathrm{d} s \{|n(x) \cdot v|
			+ O(\e)(1+ |v|)s \} \\
			\\
			&& \ \ \ \ \ \times \big| f_{\delta} \big( t- \e s , x_{\mathbf{b}}^{*}
			(X(s;0,x,v), V(s;0,x,v)), v_{\mathbf{b}}^{*} (X(s;0,x,v), V(s;0,x,v)) \big) %
			\big|^{2} \\
			&\leq&\int^{\infty}_{-\infty}\mathrm{d} t \int_{\partial\Omega} \int_{n(x)
				\cdot v >0} \int_{0}^{ 1} \big| f_{\delta} (t,x,v ) \big|^{2} \{|n(x) \cdot
			v| + O(\e)(1+ |v|)s \} \mathrm{d} s \mathrm{d} v \mathrm{d} S_{x} \\
			&\lesssim& \int_{-\infty}^{\infty} \mathrm{d} t \int_{\partial\Omega}
			\int_{n(x) \cdot v >0} \big| f_{\delta} (t,x,v ) \big|^{2} |n(x) \cdot v| 
			\mathrm{d} v \mathrm{d} S_{x} \lesssim \| f_{\delta} \|_{L^{2} (\mathbb{R}
				\times \partial\Omega \times \mathbb{R}^{3})}^{2},
		\end{eqnarray*}
		where we have used the fact, from (\ref{Z_dyn}), $O(\e)(1+|v|) |s| \leq O(\e%
		) (1+ \frac{1}{\delta})\lesssim \delta \lesssim |n(x) \cdot v|$ for $(x,v)
		\in \text{supp} ( f_{\delta} )$, and, for $n(x) \cdot v>0$, $%
		x\in\partial\Omega$, and $0 \leq s \leq t_{\mathbf{f}}^{*}(x,v)$, 
		\begin{equation*}
		(x_{\mathbf{b}}^{*} (X(s;0,x,v), V(s;0,x,v)), v_{\mathbf{b}}^{*}
		(X(s;0,x,v), V(s;0,x,v)) ) = (x_{\mathbf{b}}^{*} (x,v), v_{\mathbf{b}}^{*}
		(x,v) ) = (x,v),
		\end{equation*}
		and $t_{\mathbf{b}}^{*} (X(s;0,x,v), V(s;0,x,v)) =s$ and the change of
		variables $t- \e s \mapsto t$. Similarly we can show $(\ref{fE2_dyn})
		\lesssim \| f_{\delta} \|_{L^{2} (\mathbb{R} \times \partial\Omega \times 
			\mathbb{R}^{3})}^{2}$.
		
		\vspace{10pt}
		
		\noindent\textit{Step 5. } We claim that, in the sense of distributions on $%
		\mathbb{R} \times [\Omega_{\tilde{C} \delta^{4}}\backslash \bar{\Omega}]
		\times \mathbb{R}^{3}$, 
		\begin{eqnarray}
		&& \e \partial_{t } f_{E}+ v\cdot \nabla_{x} f_{E} + \e^{2} \Phi \cdot
		\nabla_{v} f_{E}  \label{eq_ZE_dyn} \\
		& = & \ \frac{1}{\tilde{C} \delta^{4}}v \cdot \nabla_{x} \xi(x)
		\chi^{\prime} \big( \frac{\xi(x)}{\tilde{C} \delta^{4}} \big) \Big[ %
		f_{\delta}( t- \e t_{\mathbf{b}}^{*} (x,v) , x_{\mathbf{b}}^{*}(x,v), v_{%
			\mathbf{b}}^{*}(x,v) ) \chi(t_{\mathbf{b}}^{*} (x,v))\mathbf{1}_{x_{\mathbf{b%
			}}^{*}(x,v) \in \partial\Omega}  \notag \\
		&& \ \ \ \ \ \ \ \ \ \ \ \ \ \ \ \ \ \ \ \ \ \ \ \ \ \ \ \ \ \ \ \ \ +
		f_{\delta}( t+ \e t_{\mathbf{f}}^{*} (x,v) , x_{\mathbf{f}}^{*}(x,v), v_{%
			\mathbf{f}}^{*}(x,v) ) \chi(t_{\mathbf{f}}^{*} (x,v)) \mathbf{1}_{x_{\mathbf{%
					f}}^{*}(x,v) \in \partial\Omega} \Big]  \notag \\
		&& + f_{\delta}( t- \e t_{\mathbf{b}}^{*}(x,v), x_{\mathbf{b}}^{*} (x,v), v_{%
			\mathbf{b}}^{*}(x,v)) \chi \big( \frac{\xi(x)}{\tilde{C} \delta^{4}} \big) %
		\chi^{\prime} (t_{\mathbf{b}}^{*}(x,v))\mathbf{1}_{x_{\mathbf{b}}^{*}(x,v)
			\in \partial\Omega}  \notag \\
		&& - f_{\delta}( t+ \e t_{\mathbf{f}}^{*} (x,v), x_{\mathbf{b}}^{*} (x,v),
		v_{\mathbf{b}}^{*}(x,v)) \chi \big( \frac{\xi(x)}{\tilde{C} \delta^{4}} %
		\big) \chi^{\prime} (t_{\mathbf{f}}^{*}(x,v))\mathbf{1}_{x_{\mathbf{f}%
			}^{*}(x,v) \in \partial\Omega} .  \notag
		\end{eqnarray}
		
		Note that 
		\begin{eqnarray*}
			&&[ \e \partial_{t} + v\cdot\nabla_{x} + \e^{2} \Phi \cdot \nabla_{v} ]f(t- %
			\e t_{\mathbf{b}}^{*} (x,v), x_{\mathbf{b}}^{*} (x,v) , v_{\mathbf{b}}^{*}
			(x,v)) \\
			& &= \underbrace{ [ \e \partial_{t} + v\cdot\nabla_{x} + \e^{2} \Phi \cdot
				\nabla_{v} ] (t-\e t_{\mathbf{b}}^{*} (x,v)) } \times \partial_{t} f(t- \e %
			t_{\mathbf{b}}^{*} (x,v), x_{\mathbf{b}}^{*} (x,v) , v_{\mathbf{b}}^{*}
			(x,v)) \\
			& & \ \ \ + [ v\cdot\nabla_{x} + \e^{2} \Phi \cdot \nabla_{v} ] f(s, x_{%
				\mathbf{b}}^{*} (x,v) , v_{\mathbf{b}}^{*} (x,v)) |_{s= t- \e t_{\mathbf{b}%
				}^{*} (x,v)}, \\
			&&[ \e \partial_{t} + v\cdot\nabla_{x} + \e^{2} \Phi \cdot \nabla_{v} ]f(t+ %
			\e t_{\mathbf{f}}^{*} (x,v), x_{\mathbf{f}}^{*} (x,v) , v_{\mathbf{f}}^{*}
			(x,v)) \\
			&&= \underbrace{ [ \e \partial_{t} + v\cdot\nabla_{x} + \e^{2} \Phi \cdot
				\nabla_{v} ] (t+\e t_{\mathbf{f}}^{*} (x,v)) } \times \partial_{t} f(t+ \e %
			t_{\mathbf{f}}^{*} (x,v), x_{\mathbf{f}}^{*} (x,v) , v_{\mathbf{f}}^{*}
			(x,v)) \\
			& & \ \ \ + [ v\cdot\nabla_{x} + \e^{2} \Phi \cdot \nabla_{v} ] f(s, x_{%
				\mathbf{f}}^{*} (x,v) , v_{\mathbf{f}}^{*} (x,v)) |_{s= t + \e t_{\mathbf{f}%
				}^{*} (x,v)}.
		\end{eqnarray*}
		The underbraced terms vanish 
		because $[ v\cdot \nabla_{x} + \e^{2} \Phi \cdot \nabla_{v} ] (t- \e t_{%
			\mathbf{b}}^{*} (x.v)) = \frac{d}{ds} \Big|_{s=0} (t- \e t_{\mathbf{b}}^{*}
		(X(s;0,x,v) , V(s;0,x,v) )) = \frac{d}{ds} \Big|_{s=0} (t- \e s ) = - \e  $,
		and $[ v\cdot \nabla_{x} + \e^{2} \Phi \cdot \nabla_{v} ] (t+ \e t_{\mathbf{f%
		}}^{*} (x.v)) = \frac{d}{ds} \Big|_{s=0} (t+ \e t_{\mathbf{f}}^{*}
		(X(s;0,x,v) , V(s;0,x,v) )) = \frac{d}{ds} \Big|_{s=0} (t- \e s + \e t_{%
			\mathbf{f}}^{*} (x,v)) = - \e  $. Moreover, in the sense of distributions on 
		$[\Omega_{\tilde{C} \delta^{4}}\backslash \bar{\Omega}] \times \mathbb{R}^{3}
		$, 
		\begin{eqnarray}
		&&v\cdot \nabla_{x} f_{E} + \e^{2} \Phi \cdot \nabla_{v} f_{E}  \notag \\
		& = & \ \frac{1}{\tilde{C} \delta^{4}}v \cdot \nabla_{x} \xi(x)
		\chi^{\prime} \big( \frac{\xi(x)}{\tilde{C} \delta^{4}} \big) \Big[ %
		f_{\delta}( x_{\mathbf{b}}^{*}(x,v), v_{\mathbf{b}}^{*}(x,v) ) \chi(t_{%
			\mathbf{b}}^{*} (x,v))\mathbf{1}_{x_{\mathbf{b}}^{*}(x,v) \in \partial\Omega}
		\notag \\
		& & \ \ \ \ \ \ \ \ \ \ \ \ \ \ \ \ \ \ \ \ \ \ \ \ \ \ \ \ \ \ \ \ \ +
		f_{\delta}( x_{\mathbf{f}}^{*}(x,v), v_{\mathbf{f}}^{*}(x,v) ) \chi(t_{%
			\mathbf{f}}^{*} (x,v)) \mathbf{1}_{x_{\mathbf{f}}^{*}(x,v) \in
			\partial\Omega} \Big]  \label{eq_ZE} \\
		&& + f_{\delta}(x_{\mathbf{b}}^{*} (x,v), v_{\mathbf{b}}^{*}(x,v)) \chi %
		\big( \frac{\xi(x)}{\tilde{C} \delta^{4}} \big) \chi^{\prime} (t_{\mathbf{b}%
		}^{*}(x,v))\mathbf{1}_{x_{\mathbf{b}}^{*}(x,v) \in \partial\Omega}  \notag \\
		&& - f_{\delta}(x_{\mathbf{f}}^{*} (x,v), v_{\mathbf{f}}^{*}(x,v)) \chi %
		\big( \frac{\xi(x)}{\tilde{C} \delta^{4}} \big) \chi^{\prime} (t_{\mathbf{f}%
		}^{*}(x,v))\mathbf{1}_{x_{\mathbf{f}}^{*}(x,v) \in \partial\Omega} .  \notag
		\end{eqnarray}
		For $\phi \in C_{c}^{\infty}( [\Omega_{\tilde{C} \delta^{4}}\backslash \bar{%
			\Omega}] \times \mathbb{R}^{3} )$, we choose small $t>0$ such that $X(s; 0,
		x,v) \in \Omega_{\tilde{C} \delta^{4}}\backslash \bar{\Omega}$ for all $|s|
		\leq t$ and all $(x,v) \in \text{supp}(\phi)$. Then, from (\ref{ZE_dyn}),
		for $(X(s),V(s))= (X(s;0,x,v), V(s;0,x,v))$, 
		\begin{eqnarray*}
			&&\frac{d}{ds}f_{E}(X(s), V(s)) \\
			&=& \frac{d}{ds} \Big[ f_{\delta}( x_{\mathbf{b}}^{*}(X(s),V(s)), v_{\mathbf{%
					b}}^{*} (X(s),V(s)) ) \chi (t_{\mathbf{b}}^{*} (X(s),V(s))) \mathbf{1}_{x_{%
					\mathbf{b}}^{*} (X(s),V(s)) \in \partial\Omega} \\
			&& \ \ \ \ + f_{\delta}( x_{\mathbf{f}}^{*} (X(s),V(s)), v_{\mathbf{f}}^{*}
			(X(s),V(s)) ) \chi (t_{\mathbf{f}}^{*} (X(s),V(s))) \mathbf{1}_{x_{\mathbf{f}%
				}^{*} (X(s),V(s)) \in \partial\Omega} \Big] \times \chi \big( \frac{\xi(X(s
				))}{\tilde{C} \delta^{4}} \big) \\
			&&+ \Big[ f_{\delta}( x_{\mathbf{b}}^{*}(X(s),V(s)), v_{\mathbf{b}}^{*}
			(X(s),V(s)) ) \chi (t_{\mathbf{b}}^{*} (X(s),V(s))) \mathbf{1}_{x_{\mathbf{b}%
				}^{*} (X(s),V(s)) \in \partial\Omega} \\
			&& \ \ \ \ + f_{\delta}( x_{\mathbf{f}}^{*} (X(s),V(s)), v_{\mathbf{f}}^{*}
			(X(s),V(s)) ) \chi (t_{\mathbf{f}}^{*} (X(s),V(s))) \mathbf{1}_{x_{\mathbf{f}%
				}^{*} (X(s),V(s)) \in \partial\Omega} \Big] \times \frac{d}{ds}\chi \big( 
			\frac{\xi(X(s ))}{\tilde{C} \delta^{4}} \big).
		\end{eqnarray*}
		From $(x_{\mathbf{b}}^{*}(X(s;0,x,v),V(s;0,x,v)), v_{\mathbf{b}%
		}^{*}(X(s;0,x,v),V(s;0,x,v)))= (x_{\mathbf{b}}^{*}(x,v), v_{\mathbf{b}%
		}^{*}(x,v))$ and \newline
		$(x_{\mathbf{f}}^{*}(X(s;0,x,v),V(s;0,x,v)), v_{\mathbf{f}%
		}^{*}(X(s;0,x,v),V(s;0,x,v)))= (x_{\mathbf{f}}^{*}(x,v), v_{\mathbf{f}%
		}^{*}(x,v))$ and \newline
		$t_{\mathbf{f}}^{*}(X(s;0,x,v),V(s;0,x,v) )=t_{\mathbf{f}}^{*}(x,v ) -s$ and 
		$t_{\mathbf{b}}^{*}(X(s;0,x,v),V(s;0,x,v) )=t_{\mathbf{b}}^{*}(x,v )+s$, 
		\begin{eqnarray}
		&& \frac{d}{ds}f_{E}(X(s), V(s))  \notag \\
		&=&\Big[ f_{\delta}( x_{\mathbf{b}}^{*}(x,v), v_{\mathbf{b}}^{*}(x,v) )
		\chi^{\prime}( t_{\mathbf{b}}^{*} (X(s),V(s))) \mathbf{1}_{x_{\mathbf{b}%
			}^{*}(x,v) \in \partial\Omega}  \notag \\
		&& \ \ \ - f_{\delta}( x_{\mathbf{f}}^{*}(x,v), v_{\mathbf{f}}^{*}(x,v) )
		\chi^{\prime}( t_{\mathbf{f}}^{*} (X(s),V(s))) \mathbf{1}_{x_{\mathbf{f}%
			}^{*}(x,v) \in \partial\Omega} \Big] \chi \big( \frac{\xi(X(s ))}{\tilde{C}
			\delta^{4}} \big)  \notag \\
		& & + \Big[ f_{\delta}( x_{\mathbf{b}}^{*}(x,v), v_{\mathbf{b}}^{*}(x,v) )
		\chi (t_{\mathbf{b}}^{*}(X(s),V(s)) ) \mathbf{1}_{x_{\mathbf{b}}^{*}(x,v)
			\in \partial\Omega}  \label{derfE} \\
		&& \ \ \ + f_{\delta}( x_{\mathbf{f}}^{*}(x,v), v_{\mathbf{f}}^{*}(x,v) )
		\chi (t_{\mathbf{f}}^{*}(X(s),V(s)) )\mathbf{1}_{x_{\mathbf{f}}^{*}(x,v) \in
			\partial\Omega} \Big] \frac{1}{\tilde{C} \delta^{4}}V(s ) \cdot \nabla_{x}
		\xi(X(s )) \chi^{\prime} \big( \frac{\xi(X(s ))}{\tilde{C} \delta^{4}} \big).
		\notag
		\end{eqnarray}
		
		By the change of variables $(x,v) \mapsto (X(s; 0, x,v), V(s;0,x,v))$, for
		sufficiently small $s$, 
		\begin{eqnarray}
		&&-\iint_{[\Omega_{\tilde{C} \delta^{4}} \backslash \bar{\Omega}] \times 
			\mathbb{R}^{3} } f_{E} (x,v) \{ v\cdot \nabla_{x} + \e^{2} \Phi \cdot
		\nabla_{v} \}\phi(x,v) \mathrm{d} x \mathrm{d} v  \label{iniz} \\
		&=&- \iint_{[\Omega_{\tilde{C} \delta^{4}} \backslash \bar{\Omega}] \times 
			\mathbb{R}^{3} } f_{E} (X(s), V(s)) \{ V(s)\cdot \nabla_{X} + \e^{2} \Phi
		\cdot \nabla_{V} \} \phi(X(s), V(s)) \mathrm{d} x \mathrm{d} v  \notag \\
		&=& - \iint_{[\Omega_{\tilde{C} \delta^{4}} \backslash \bar{\Omega}] \times 
			\mathbb{R}^{3} } f_{E} (X(s), V(s)) \frac{d}{ds} \phi(X(s), V(s)) \mathrm{d}
		x \mathrm{d} v.  \notag
		\end{eqnarray}
		
		Since the change of variables $(x,v) \mapsto (X(s; 0, x,v), V(s;0,x,v))$ has
		unit Jacobian, it follows that, for $s$ sufficiently small, 
		\begin{equation*}
		\iint_{[\Omega_{\tilde{C} \delta^{4}} \backslash \bar{\Omega}] \times 
			\mathbb{R}^{3} } f_{E} ((X(s), V(s)) \phi(X(s), V(s)))=\iint_{[\Omega_{%
				\tilde{C} \delta^{4}} \backslash \bar{\Omega}] \times \mathbb{R}^{3} } f_{E}
		(x,v) \phi(x,v) ,
		\end{equation*}
		and hence 
		\begin{equation*}
		\frac{d}{ds} \iint_{[\Omega_{\tilde{C} \delta^{4}} \backslash \bar{\Omega}]
			\times \mathbb{R}^{3} } f_{E} ((X(s), V(s)) \phi(X(s), V(s)) =0.
		\end{equation*}
		Therefore we can move the $s$-derivative on $f_E$: By (\ref{derfE}), 
		\begin{eqnarray*}
			&&(\ref{iniz}) \\
			&=& \iint_{[\Omega_{\tilde{C} \delta^{4}} \backslash \bar{\Omega}] \times 
				\mathbb{R}^{3} } \frac{d}{ds} f_{E} (X(s), V(s)) \phi(X(s), V(s)) \mathrm{d}
			x \mathrm{d} v \\
			&=& \iint_{[\Omega_{\tilde{C} \delta^{4}} \backslash \bar{\Omega}] \times 
				\mathbb{R}^{3} } \Big[ f_{\delta}( x_{\mathbf{b}}^{*}(x,v), v_{\mathbf{b}%
			}^{*}(x,v) ) \chi^{\prime}( t_{\mathbf{b}}^{*} (X(s),V(s))) \mathbf{1}_{x_{%
					\mathbf{b}}^{*}(x,v) \in \partial\Omega} \\
			&& \ \ \ \ \ \ \ \ \ \ \ \ \ \ - f_{\delta}( x_{\mathbf{f}}^{*}(x,v), v_{%
				\mathbf{f}}^{*}(x,v) ) \chi^{\prime}( t_{\mathbf{f}}^{*} (X(s),V(s))) 
			\mathbf{1}_{x_{\mathbf{f}}^{*}(x,v) \in \partial\Omega} \Big] \chi \big( 
			\frac{\xi(X(s ))}{\tilde{C} \delta^{4}} \big) \phi(X(s ), V(s )) \\
			& +& \iint_{[\Omega_{\tilde{C} \delta^{4}} \backslash \bar{\Omega}] \times 
				\mathbb{R}^{3} } \Big[ f_{\delta}( x_{\mathbf{b}}^{*}(x,v), v_{\mathbf{b}%
			}^{*}(x,v) ) \chi (t_{\mathbf{b}}^{*}(X(s),V(s)) ) \mathbf{1}_{x_{\mathbf{b}%
				}^{*}(x,v) \in \partial\Omega} \\
			&& + f_{\delta}( x_{\mathbf{f}}^{*}(x,v), v_{\mathbf{f}}^{*}(x,v) ) \chi (t_{%
				\mathbf{f}}^{*}(X(s),V(s)) )\mathbf{1}_{x_{\mathbf{f}}^{*}(x,v) \in
				\partial\Omega} \Big] \\
			&&\hskip 5cm \times \frac{1}{\tilde{C} \delta^{4}}V(s ) \cdot \nabla_{x}
			\xi(X(s ))\chi^{\prime} \big( \frac{\xi(X(s ))}{\tilde{C} \delta^{4}} \big)%
			\phi(X(s ), V(s )) .
		\end{eqnarray*}
		From the change of variable $(X(s;0,x,v), V(s;0,x,v)) \mapsto (x,v)$, 
		\begin{eqnarray*}
			(\ref{iniz}) &=& \iint_{[\Omega_{\tilde{C} \delta^{4}} \backslash \bar{\Omega%
				}] \times \mathbb{R}^{3} } \big[ f_{\delta}(x_{\mathbf{b}}^{*} (x,v), v_{%
				\mathbf{b}}^{*}(x,v)) \chi^{\prime} (t_{\mathbf{b}}^{*}(x,v)) \mathbf{1}_{x_{%
					\mathbf{b}}^{*}(x,v) \in \partial\Omega} \\
			&&- f_{\delta}(x_{\mathbf{f}}^{*} (x,v), v_{\mathbf{f}}^{*}(x,v))
			\chi^{\prime} (t_{\mathbf{f}}^{*}(x,v)) \mathbf{1}_{x_{\mathbf{f}}^{*}(x,v)
				\in \partial\Omega} \big] \chi \big(\frac{\xi(x)}{\tilde{C} \delta^{4}} %
			\big) \phi(x,v) \\
			& +& \iint_{[\Omega_{\tilde{C} \delta^{4}} \backslash \bar{\Omega}] \times 
				\mathbb{R}^{3} } \Big[ f_{\delta}( x_{\mathbf{b}}^{*}(x,v), v_{\mathbf{b}%
			}^{*}(x,v) ) \chi (t_{\mathbf{b}}^{*}(x,v) ) \mathbf{1}_{x_{\mathbf{b}%
				}^{*}(x,v) \in \partial\Omega} \\
			&& \ \ \ \ \ \ + f_{\delta}( x_{\mathbf{f}}^{*}(x,v), v_{\mathbf{f}%
			}^{*}(x,v) ) \chi (t_{\mathbf{f}}^{*}(x,v) )\mathbf{1}_{x_{\mathbf{f}%
				}^{*}(x,v) \in \partial\Omega} \Big] \frac{1}{\tilde{C} \delta^{4}}v \cdot
			\nabla_{x} \xi(x) \chi^{\prime} \big( \frac{\xi(x)}{\tilde{C} \delta^{4}} %
			\big)\phi(x,v).
		\end{eqnarray*}
		Hence (\ref{eq_ZE_dyn}) is proved.
		
		On the other hand, following the bounds of (\ref{fE1_dyn}) and (\ref{fE2_dyn}%
		) in \textit{Step 4} we prove the third line of (\ref{force_Z_dyn}).
		
		\vspace{8pt}
		
		\noindent\textit{Step 6. } We define $\bar{f}(t,x,v)$ for $(t,x,v) \in 
		\mathbb{R} \times \mathbb{R}^{3} \times\mathbb{R}^{3}$: 
		\begin{equation}  \label{bar_Z_dyn}
		\begin{split}
		\bar{f}(t,x,v) & \ : = \ f_{\delta}(t,x,v) \mathbf{1}_{(x,v) \in \bar{\Omega}
			\times \mathbb{R}^{3}} + f_{E}(t,x,v) \mathbf{1}_{(x,v) \in [\mathbb{R}^{3}
			\backslash \bar{\Omega}] \times \mathbb{R}^{3}}.
		\end{split}%
		\end{equation}
		For $\phi \in C^{\infty}_{c}(\mathbb{R} \times \mathbb{R}^{3} \times \mathbb{%
			R}^{3})$, by Lemma \ref{timedepgreen}, 
		\begin{eqnarray*}
			&&- \int_{- \infty}^{\infty} \mathrm{d} t \iint_{\mathbb{R}^{3} \times 
				\mathbb{R}^{3}} \bar{f}(t,x,v) \{\e \partial_{t} + v\cdot \nabla_{x} + \e%
			^{2} \Phi \cdot \nabla_{v} \} \phi(t,x,v) \mathrm{d} x \mathrm{d} v \\
			&=&- \int_{- \infty}^{\infty} \mathrm{d} t\iint_{\Omega \times \mathbb{R}%
				^{3}} f_{\delta}(t,x,v) \{ \e \partial_{t} + v\cdot \nabla_{x} + \e^{2} \Phi
			\cdot \nabla_{v} \} \phi(t,x,v) \mathrm{d} x \mathrm{d} v \\
			&&- \int_{- \infty}^{\infty} \mathrm{d} t \iint_{ [\mathbb{R}^{3} \backslash 
				\bar{\Omega}] \times \mathbb{R}^{3}} f_{E}(t,x,v) \{\e \partial_{t} + v\cdot
			\nabla_{x} + \e^{2} \Phi \cdot \nabla_{v} \} \phi(t,x,v) \mathrm{d} x 
			\mathrm{d} v \\
			&=& \int_{- \infty}^{\infty} \mathrm{d} t \int_{\gamma} f_{\delta} (t,x,v)
			\phi(t,x,v) \{n(x) \cdot v\} \mathrm{d} S_{x} \mathrm{d} v + \int_{-
				\infty}^{\infty} \mathrm{d} t \int_{\gamma} f_{E}(t,x,v) \phi(t,x,v) \{-n(x)
			\cdot v\} \mathrm{d} S_{x} \mathrm{d} v \\
			&&+ \int_{- \infty}^{\infty} \mathrm{d} t \iint_{\Omega\times \mathbb{R}%
				^{3}} \{\e \partial_{t} + v\cdot \nabla_{x} + \e^{2} \Phi \cdot \nabla_{v}
			\} f_{\delta}(t,x,v) \phi(t,x,v) \mathrm{d} x \mathrm{d} v \\
			&&+ \int_{- \infty}^{\infty} \mathrm{d} t \iint_{ [\Omega_{\tilde{C}
					\delta^{4}} \backslash \bar{\Omega}]\times \mathbb{R}^{3}} \{ \e %
			\partial_{t} + v\cdot \nabla_{x} + \e^{2} \Phi \cdot \nabla_{v} \} {f}%
			_{E}(x,v) \phi(x,v) \mathrm{d} x \mathrm{d} v,
		\end{eqnarray*}
		where the contributions of $\{t= \infty\}$ and $\{t= -\infty\}$ vanish since 
		$\phi(t) \in C_{c}^{\infty} (\mathbb{R})$.
		
		From (\ref{fE=fd_dyn}), the boundary contributions are cancelled: 
		\begin{equation*}
		\int_{- \infty}^{\infty} \int_{\gamma} f_{\delta} (t,x,v) \phi(t,x,v) 
		\mathrm{d} \gamma \mathrm{d} t - \int_{- \infty}^{\infty} \int_{\gamma}
		f_{E}(t,x,v) \phi(t,x,v) \mathrm{d} \gamma \mathrm{d} t =0.
		\end{equation*}
		
		Further from (\ref{eq_Z_dyn}) and (\ref{eq_ZE_dyn}), we prove that $\bar{f}$
		solves (\ref{eq_barf_dyn}) in the sense of distributions on $\mathbb{R}
		\times \mathbb{R}^{3} \times \mathbb{R}^{3}$.\unhide
	\end{proof}
	
	\hide
	\begin{proof}
		
		\noindent\textit{Step 1. }  In the sense of distributions on $[0, \infty)\times \Omega \times \mathbb{R}^{3}$,
		\begin{equation}\label{eq_Z_dyn}
		\begin{split}
		&\e \partial_{t} f_{\delta}+ v\cdot \nabla_{x} f_{\delta} + \e^{2} \Phi \cdot \nabla_{v}   f_{\delta} \\
		& \ =   \Big[1-\chi(\frac{n(x) \cdot v}{\delta})  \chi \big( \frac{ \xi(x) }{\delta}\big)\Big] \chi(\delta|v|)  \\
		& \ \ \ \ \ \ \ \ \times
		\big[ \mathbf{1}_{t \in [0,\infty)}
		g   + \mathbf{1}_{t \in ( - \infty, 0 ]}
		\chi (t)
		\{ \e \frac{\chi^{\prime} (t)}{ \chi(t)} + v \cdot \nabla_{x} + \e^{2} \Phi \cdot \nabla_{v}  \} f_{0} (x,v)\big] \\
		& \ \ \  + \big[\mathbf{1}_{t \in [0,\infty)} f + \mathbf{1}_{t \in (- \infty, 0 ]} \chi(t) f_{0} (x,v)
		\big]
		\{v\cdot \nabla_{x}  + \e^{2} \Phi \cdot \nabla_{v}\} \Big(  [1-\chi(\frac{n(x) \cdot v}{\delta}) \chi \big( \frac{\xi(x)}{\delta}\big) ] \chi(\delta|v|)\Big).
		\end{split}
		\end{equation}
		From \textit{Step 1} of the proof of Lemma \ref{extension} and $(\ref{der_chi})$, we prove the first and the second line of (\ref{force_Z_dyn}).

		\vspace{4pt}
		
		\noindent\textit{Step 2. } We claim that if $ 0 \leq \xi (x) \leq  \tilde{C} \delta^{4}, \ |n(x) \cdot v|>  \delta$ and $ |v| \leq \frac{1}{\delta}$ then either $\xi(\xft (x,v))=\tilde{C}\delta^{4}$ or $\xi(\xbt(x,v)) =\tilde{C}\delta^{4}$. 
		
		From \textit{Step 2} of the proof of Lemma \ref{extension},
		\begin{eqnarray*}
			\xi(Y(s;0,x,v)) 
			=  \xi(X(\frac{s}{\e}; 0,x,v  )) \geq  \delta \frac{|s|}{2\e} ,\end{eqnarray*}
		for all $0\leq |s| \leq \frac{ \e\delta^{3}}{ 4(1+ \|\xi \|_{C^{2}})}$ with $0< \varepsilon \ll1$. Especially with $\e s_{*} = +\frac{  \varepsilon\delta^{3}}{4(1+ \| \xi \|_{C^{2}})}$ for $n(x) \cdot v >\delta$ and $\e s_{*} = -\frac{\varepsilon\delta^{3} }{4(1+ \| \xi \|_{C^{2}})}$ for $n(x) \cdot v <\delta$,
		\[
		\xi (Y(s_{*};0,x,v)) > \tilde{C} \delta^{4}.
		\]
		Therefore, by the intermediate value theorem, we prove our claim.

		\vspace{10pt}
		
		\noindent\textit{Step 3. } We define $f_{E}(t,x,v)$ for $(x,v) \in [ \mathbb{R}^{3} \backslash \bar{\Omega}] \times \mathbb{R}^{3}$:
		\begin{equation}\label{ZE_dyn}
		\begin{split}
		f_{E}(t,x,v) & \ :=  \  f_{\delta}(t- \e \tb^{*}(x,v),\xb^{*}(x,v), \vb^{*}(x,v)) \chi\big(   \frac{\xi(x)}{\tilde{C} \delta^{4}}  \big) \chi \big( {\tb^{*} (x,v)}  \big)   \ \ \ \text{if} \ \ \xb^{*}(x,v) \in \partial\Omega,\\
		& \ :=  \    f_{\delta}(t+ \e \tf^{*}(x,v),\xf^{*}(x,v), \vf^{*}(x,v))  \chi\big(   \frac{\xi(x)}{\tilde{C} \delta^{4}}  \big)  \chi \big( {\tf^{*} (x,v)}  \big)  \ \ \ \ \text{if} \ \ \xf^{*}(x,v) \in \partial\Omega,\\
		& \ :=  \  0  \ \ \    \ \  \ \ \ \  \ \   \ \ \ \  \ \ \ \   \ \ \ \   \ \ \ \  \ \ \ \ \text{if} \ \ \xb^{*}(x,v) \notin \partial\Omega \   \text{and} \  \xf^{*}(x,v) \notin \partial\Omega.
		\end{split}
		\end{equation}
		
		We check that $f_{E}$ is well-defined. It suffices to prove the following:
		\begin{equation*}\begin{split}
		& \text{If}  \ \ \xb^{*}(x,v) \in \partial\Omega  \text{ and } \xf^{*}(x,v) \in \partial\Omega\\
		& \  \text{then} \  \   f_{\delta}( t- \e t_{\mathbf{b}}^{*} (x,v),  \xb^{*}(x,v), \vb^{*}(x,v))  \chi\big(   \frac{\xi(x)}{\tilde{C} \delta^{4}}  \big )   =  0 =  f_{\delta}(
		t+ \e t_{\mathbf{f}}^{*} (x,v),
		\xf^{*}(x,v), \vf^{*}(x,v)) \chi\big(   \frac{\xi(x)}{\tilde{C} \delta^{4}}  \big)  .
		\end{split}\end{equation*}
		If $|n(\xb^{*}(x,v)) \cdot \vb^{*}(x,v)| \leq \delta$ or $|\vb^{*}(x,v)| \geq \frac{1}{\delta}$ then $ f_{\delta}(
		t- \e t^{*}_{\mathbf{b}}(x,v),
		\xb^{*}(x,v), \vb^{*}(x,v))=0$ due to (\ref{Z_support_dyn}). If $n(\xb^{*}(x,v)) \cdot \vb^{*}(x,v) > \delta$ and $|\vb^{*}(x,v)|\leq \frac{1}{\delta}$ then, due to \textit{Step 2}, $\xi(\xf^{*}(x,v))  = \xi( \xf^{*}( \xb^{*}(x,v), \vb^{*}(x,v)  ) ) = \tilde{C}\delta^{4}$ so that $\xf^{*}(x,v) \notin \partial\Omega$.
		
		On the other hand, if $|n(\xf^{*}(x,v)) \cdot \vf^{*}(x,v)| \leq \delta$ or $|\vf^{*}(x,v)|\geq \frac{1}{\delta}$ then $  f_{\delta}( t+ \e t^{*}_{\mathbf{f}}(x,v) ,\xf^{*}(x,v), \vf^{*}(x,v))=0$ due to (\ref{Z_support_dyn}). If $n(\xf^{*}(x,v)) \cdot \vf^{*}(x,v) <-\delta$ and $|\vf^{*}(x,v)|\leq \frac{1}{\delta}$ then, due to \textit{Step 2}, $\xi(\xb^{*}(x,v))  = \xi( \xb^{*}( \xf^{*}(x,v), \vf^{*}(x,v)  ) ) = \tilde{C}\delta^{4}$ so that $\xb^{*}(x,v) \notin \partial\Omega$.
		
		Note that 
		\begin{equation}\label{fE=fd_dyn}
		f_{E}(t,x,v) = f_{\delta} (t,x,v) \ \ \ \text{for all } x \in \partial\Omega.
		\end{equation}
		If $x \in \partial\Omega$ and $n(x) \cdot v >\delta$ then $(\xb^{*}(x,v), \vb^{*}(x,v)) = (x ,v)$. From the definition (\ref{ZE_dyn}), for those $(x,v)$, we have $f_{E}(t,x,v) = f_{\delta} (t,x,v)$. If $x \in \partial\Omega$ and $n(x) \cdot v< -\delta$ then $(\xf^{*}(x,v), \vf^{*}(x.v)) = (x,v)$. From the definition (\ref{ZE_dyn}), we conclude (\ref{fE=fd_dyn}) again. Otherwise, if $-\delta< n(x) \cdot v< \delta$ then $f_{E}|_{\partial\Omega} \equiv 0 \equiv f_{\delta}|_{\partial\Omega}$.

		\vspace{10pt}
		
		\noindent\textit{Step 4. } We claim that $f_{E}(x,v) \in L^{2}([\mathbb{R}^{3} \backslash \bar{\Omega}] \times \mathbb{R}^{3})$. 
		
		From the definition of (\ref{ZE_dyn}), we have $f_{E}(x,v) \equiv 0$ if $\xb^{*} (x,v) \notin \partial\Omega$ and $\xf^{*} (x,v) \notin \partial\Omega$. Therefore, from (\ref{ZE_dyn}),
		\begin{eqnarray}
		&&\int_{-\infty}^{\infty} \iint_{[ \mathbb{R}^{3 } \backslash \Omega ]\times \mathbb{R}^{3}} |f_{E}  (t,x,v)|^{2} \dd x \dd v \dd t\nonumber\\
		&=&\int_{-\infty}^{\infty}  \iint_{ [ \mathbb{R}^{3 } \backslash \Omega ] \times \mathbb{R}^{3}}  \mathbf{1}_{\xb^{*} (x,v) \in \partial\Omega} |f_{E}  |^{2}  + \int_{-\infty}^{\infty}  \iint _{ [ \mathbb{R}^{3 } \backslash \Omega ] \times \mathbb{R}^{3}}  \mathbf{1}_{\xf^{*}   \in \partial\Omega}  |f_{E}  |^{2} \nonumber \\
		&=& \int_{-\infty}^{\infty} \iint_{ [ \mathbb{R}^{3 } \backslash \Omega ]\times \mathbb{R}^{3}} \mathbf{1}_{\xb^{*} (x,v) \in \partial\Omega}  |f_{\delta} (
		t- \e t^{*}_{\mathbf{b}} ,
		\xb^{*}  , \vb^{*}  )|^{2} \big| \chi \big( \frac{\xi(x)}{\tilde{C} \delta^{4}}  \big)  \big|^{2}
		|\chi ( {\tb^{*}  }  ) |^{2}\dd x \dd v \dd t
		\label{fE1_dyn} \\
		&&+   \int_{-\infty}^{\infty} 
		\iint_{ [ \mathbb{R}^{3 } \backslash \Omega ] \times \mathbb{R}^{3}} \mathbf{1}_{\xf^{*} (x,v) \in \partial\Omega}  |f_{\delta} (
		t+ \e t^{*}_{\mathbf{f}} ,
		\xf^{*} , \vf^{*} )|^{2} \big| \chi \big( \frac{\xi(x)}{\tilde{C} \delta^{4}}  \big)  \big|^{2}
		|\chi ( {\tf^{*}  }  ) |^{2} \dd x \dd v \dd t
		,    \label{fE2_dyn}
		\end{eqnarray}
		where $(t^{*}_{\mathbf{b}},x^{*}_{\mathbf{b}},v^{*}_{\mathbf{b}} )$ and $(t^{*}_{\mathbf{f}},x^{*}_{\mathbf{f}},v^{*}_{\mathbf{f}} )$ are evaluated at $(x,v)$.
		
		By (\ref{bdry_int2}),
		\begin{eqnarray*}
			(\ref{fE1_dyn}) &\leq&
			\int^{\infty}_{- \infty} \dd t
			\int_{\partial\Omega} \int_{n(x) \cdot v >0} \int_{0}^{ \min\{ \tf^{*}(x,v), 1  \}}\dd S_{x}  \dd v  \dd s
			\{|n(x) \cdot v| + O(\e)(1+ |v|)s  \}\\
			\\
			&& \  \ \ \ \   \times
			\big|
			f_{\delta} \big( t-  \e s   ,  \xb^{*} (X(s;0,x,v), V(s;0,x,v)), \vb^{*} (X(s;0,x,v), V(s;0,x,v))  \big)  \big|^{2}   \\
			&\leq&\int^{\infty}_{-\infty}\dd t \int_{\partial\Omega} \int_{n(x) \cdot v >0} \int_{0}^{ 1}  \big| f_{\delta} (t,x,v )  \big|^{2}   \{|n(x) \cdot v| + O(\e)(1+ |v|)s  \} \dd s \dd v \dd S_{x}\\
			&\lesssim&
			\int_{-\infty}^{\infty} \dd t
			\int_{\partial\Omega} \int_{n(x) \cdot v >0}    \big| f_{\delta} (t,x,v )  \big|^{2}    |n(x) \cdot v|   \dd v \dd S_{x} \lesssim  \| f_{\delta} \|_{L^{2} (\mathbb{R} \times   \partial\Omega \times \mathbb{R}^{3})}^{2},
		\end{eqnarray*}
		where we have used the fact, from (\ref{Z_dyn}), $O(\e)(1+|v|) |s| \leq O(\e) (1+ \frac{1}{\delta})\lesssim \delta \lesssim |n(x) \cdot v|$ for $(x,v) \in \text{supp} ( f_{\delta}  )$, and, for $n(x) \cdot v>0$, $x\in\partial\Omega$, and $0 \leq s \leq \tf^{*}(x,v)$,  
		\[
		(\xb^{*} (X(s;0,x,v), V(s;0,x,v)), \vb^{*} (X(s;0,x,v), V(s;0,x,v))  ) 
		= (\xb^{*} (x,v), \vb^{*} (x,v)  ) 
		= (x,v),
		\]
		and $\tb^{*} (X(s;0,x,v), V(s;0,x,v)) =s$ and the change of variables $t- \e s \mapsto t$. Similarly we can show $(\ref{fE2_dyn}) \lesssim  \| f_{\delta} \|_{L^{2} (\mathbb{R} \times   \partial\Omega \times \mathbb{R}^{3})}^{2}$.

		\vspace{10pt}
		
		\noindent\textit{Step 5. }   We claim that, in the sense of distributions on $\mathbb{R} \times 
		[\Omega_{\tilde{C} \delta^{4}}\backslash \bar{\Omega}] \times \mathbb{R}^{3}$,
		\begin{eqnarray} 
		&& \e \partial_{t } f_{E}+ v\cdot \nabla_{x} f_{E} + \e^{2} \Phi \cdot \nabla_{v} f_{E} \label{eq_ZE_dyn}\\
		& = &  \ \frac{1}{\tilde{C} \delta^{4}}v \cdot \nabla_{x} \xi(x) \chi^{\prime}  \big( \frac{\xi(x)}{\tilde{C} \delta^{4}} \big) \Big[  f_{\delta}( t- \e \tb^{*} (x,v) , \xb^{*}(x,v), \vb^{*}(x,v) )  \chi(\tb^{*} (x,v))\mathbf{1}_{\xb^{*}(x,v) \in \partial\Omega} \nonumber\\
		&& \ \ \ \ \ \ \ \ \ \ \  \ \ \ \ \ \ \ \ \ \ \  \ \ \ \ \ \ \ \ \ \ \   + 
		f_{\delta}( t+ \e \tf^{*} (x,v) , \xf^{*}(x,v), \vf^{*}(x,v) ) \chi(\tf^{*} (x,v)) \mathbf{1}_{\xf^{*}(x,v) \in \partial\Omega}
		\Big]  \notag \\
		&&   +  f_{\delta}( t- \e \tb^{*}(x,v),   \xb^{*} (x,v), \vb^{*}(x,v)) \chi \big( \frac{\xi(x)}{\tilde{C} \delta^{4}} \big)  \chi^{\prime} (\tb^{*}(x,v))\mathbf{1}_{\xb^{*}(x,v) \in \partial\Omega}\notag\\
		&& -  f_{\delta}( t+ \e \tf^{*} (x,v),  \xb^{*} (x,v), \vb^{*}(x,v)) \chi \big( \frac{\xi(x)}{\tilde{C} \delta^{4}} \big)  \chi^{\prime} (\tf^{*}(x,v))\mathbf{1}_{\xf^{*}(x,v) \in \partial\Omega} 
		. 
		\notag
		\end{eqnarray}
		
		Note that 
		\begin{eqnarray*}
			&&[ \e \partial_{t} + v\cdot\nabla_{x} + \e^{2} \Phi \cdot \nabla_{v} ]f(t- \e t_{\mathbf{b}}^{*} (x,v), x_{\mathbf{b}}^{*} (x,v) , v_{\mathbf{b}}^{*} (x,v))\\
			& &= \underbrace{ [ \e \partial_{t} + v\cdot\nabla_{x} + \e^{2} \Phi \cdot \nabla_{v} ] (t-\e t_{\mathbf{b}}^{*} (x,v)) } \times \partial_{t} f(t- \e t_{\mathbf{b}}^{*} (x,v), x_{\mathbf{b}}^{*} (x,v) , v_{\mathbf{b}}^{*} (x,v))\\
			& & \ \ \ + [   v\cdot\nabla_{x} + \e^{2} \Phi \cdot \nabla_{v} ] f(s, x_{\mathbf{b}}^{*} (x,v) , v_{\mathbf{b}}^{*} (x,v)) |_{s= t- \e t_{\mathbf{b}}^{*} (x,v)},\\
			&&[ \e \partial_{t} + v\cdot\nabla_{x} + \e^{2} \Phi \cdot \nabla_{v} ]f(t+ \e t_{\mathbf{f}}^{*} (x,v), x_{\mathbf{f}}^{*} (x,v) , v_{\mathbf{f}}^{*} (x,v))\\
			&&= \underbrace{ [ \e \partial_{t} + v\cdot\nabla_{x} + \e^{2} \Phi \cdot \nabla_{v} ] (t+\e t_{\mathbf{f}}^{*} (x,v)) } \times \partial_{t} f(t+ \e t_{\mathbf{f}}^{*} (x,v), x_{\mathbf{f}}^{*} (x,v) , v_{\mathbf{f}}^{*} (x,v))\\
			& & \ \ \ + [   v\cdot\nabla_{x} + \e^{2} \Phi \cdot \nabla_{v} ] f(s, x_{\mathbf{f}}^{*} (x,v) , v_{\mathbf{f}}^{*} (x,v)) |_{s= t + \e t_{\mathbf{f}}^{*} (x,v)}.
		\end{eqnarray*}
		If the underbraced term vanishes then we can apply the \textit{Step 5} of the proof of Lemma \ref{extension} to conclude (\ref{eq_ZE_dyn}). 
		
		This is true because $[ v\cdot \nabla_{x}  + \e^{2} \Phi \cdot \nabla_{v}  ]   (t- \e \tb^{*} (x.v)) 
		=   \frac{d}{ds} \Big|_{s=0}  (t- \e \tb^{*} (X(s;0,x,v) , V(s;0,x,v) )) 
		=   \frac{d}{ds} \Big|_{s=0}  (t- \e s ) 
		=  - \e  $,
		and $[ v\cdot \nabla_{x}  + \e^{2} \Phi \cdot \nabla_{v}  ]   (t+ \e \tf^{*} (x.v)) 
		= \frac{d}{ds} \Big|_{s=0}  (t+ \e \tf^{*} (X(s;0,x,v) , V(s;0,x,v) )) 
		=   \frac{d}{ds} \Big|_{s=0}  (t- \e s +  \e \tf^{*} (x,v)) 
		=  - \e  $. 
		
		On the other hand, following the bounds of (\ref{fE1_dyn}) and (\ref{fE2_dyn}) in \textit{Step 4} we prove the third line of (\ref{force_Z_dyn}).
		
		\vspace{8pt}
		
		\noindent\textit{Step 6. } We define $\bar{f}(t,x,v)$ for $(t,x,v) \in \mathbb{R} \times  \mathbb{R}^{3} \times\mathbb{R}^{3}$:
		\begin{equation}\label{bar_Z_dyn}
		\begin{split}
		\bar{f}(t,x,v) &  \ : = \ f_{\delta}(t,x,v) \mathbf{1}_{(x,v) \in \bar{\Omega} \times \mathbb{R}^{3}}   +   f_{E}(t,x,v)   \mathbf{1}_{(x,v) \in   [\mathbb{R}^{3} \backslash \bar{\Omega}] \times \mathbb{R}^{3}}.
		\end{split}
		\end{equation} 
		For $\phi \in C^{\infty}_{c}(\mathbb{R} \times \mathbb{R}^{3} \times \mathbb{R}^{3})$, by Lemma \ref{timedepgreen},
		\begin{eqnarray*}
			&&- \int_{- \infty}^{\infty} \dd t \iint_{\mathbb{R}^{3} \times \mathbb{R}^{3}} \bar{f}(t,x,v)
			\{\e \partial_{t} + v\cdot \nabla_{x} + \e^{2} \Phi \cdot \nabla_{v} \} \phi(t,x,v) \dd x  \dd v\\
			&=&-  \int_{- \infty}^{\infty} \dd t\iint_{\Omega \times \mathbb{R}^{3}} f_{\delta}(t,x,v)   \{ \e \partial_{t} + v\cdot \nabla_{x} + \e^{2} \Phi \cdot \nabla_{v} \} \phi(t,x,v)   \dd x   \dd v \\
			&&- \int_{- \infty}^{\infty} \dd t \iint_{   [\mathbb{R}^{3} \backslash \bar{\Omega}] \times \mathbb{R}^{3}} f_{E}(t,x,v)   \{\e \partial_{t} + v\cdot \nabla_{x} + \e^{2} \Phi \cdot \nabla_{v} \} \phi(t,x,v)   \dd x  \dd v \\
			&=&  \int_{- \infty}^{\infty} \dd t \int_{\gamma}  f_{\delta} (t,x,v) \phi(t,x,v) \{n(x) \cdot v\}    \dd S_{x}   \dd v +
			\int_{- \infty}^{\infty} \dd t \int_{\gamma}   f_{E}(t,x,v) \phi(t,x,v)   \{-n(x) \cdot v\}  \dd S_{x}  \dd v\\
			&&+   \int_{- \infty}^{\infty} \dd t  \iint_{\Omega\times \mathbb{R}^{3}}  \{\e \partial_{t} +  v\cdot \nabla_{x} + \e^{2} \Phi \cdot \nabla_{v} \} f_{\delta}(t,x,v)   \phi(t,x,v)  \dd x \dd v\\
			&&+ \int_{- \infty}^{\infty} \dd t \iint_{    [\Omega_{\tilde{C} \delta^{4}} \backslash \bar{\Omega}]\times \mathbb{R}^{3}}  \{ \e \partial_{t} + v\cdot \nabla_{x} + \e^{2} \Phi \cdot \nabla_{v} \} {f}_{E}(x,v)   \phi(x,v)  \dd x  \dd v,
		\end{eqnarray*}
		where the contributions of $\{t= \infty\}$ and $\{t= -\infty\}$ vanish since $\phi(t)  \in C_{c}^{\infty} (\mathbb{R})$.
		
		From (\ref{fE=fd_dyn}), the boundary contributions are cancelled: $$
		\int_{- \infty}^{\infty}  
		\int_{\gamma}  f_{\delta} (t,x,v) \phi(t,x,v)   \dd \gamma \dd t - 
		\int_{- \infty}^{\infty}
		\int_{\gamma}   f_{E}(t,x,v) \phi(t,x,v)    \dd \gamma   \dd t  =0.$$
		
		Further from (\ref{eq_Z_dyn}) and (\ref{eq_ZE_dyn}), we prove that $\bar{f}$ solves (\ref{eq_barf_dyn}) in the sense of distributions on $\mathbb{R} \times \mathbb{R}^{3} \times \mathbb{R}^{3}$.\end{proof}
	
	Recall $T>0$ in (\ref{large_time}). With such choice $T>0$, 
	\begin{eqnarray*}
		\bar{f} (t,x,v) &=&  - \int^{\e T}_{0}
		\frac{1}{\e}h( t+s  , Y(t+s; t,x,v) ,W(t+s; t,x,v)  ) \dd s\\
		&=&  - \int^{\e T}_{0}\frac{1}{\e}
		h( t+s  , X( t + \frac{s}{\e}; t,x,v) ,V(t + \frac{s}{\e} ; t,x,v)  ) \dd s.
	\end{eqnarray*}
	Note that, from (\ref{ZE_dyn}),  
	\begin{equation}\label{supp_fE_dyn}
	\bar{f}  (x,v) \equiv 0, \ \  \ \text{for} \ \ \xi(x) > 2\tilde{C} \delta^{4}  \ \ \text{or} \ \ |v|>  {2}{  \delta^{-1}}  \ \ \text{or} \ \ \ |v| <  {\delta}/{2} .
	\end{equation}
	Therefore
	\begin{equation}\label{duhamel_f_dyn}
	|\bar{f}(t,x,v)|  \ \leq \   \frac{1}{\e} \int_{0}^{\e T}  \mathbf{1}_{\frac{\delta}{2} \leq |v| \leq  {  \frac{2}{\delta}}}
	|h(  t+s,X(t + \frac{s}{\e};t,x,v), V(t + \frac{s}{\e};t,x,v))|\dd s.
	\end{equation} 
	
	\begin{definition}
		For fixed $T$ in (\ref{large_time}) and $\delta>0$ and a smooth function $\phi \in   L^{1}(\mathbb{R}^{3})$, we define the average operator $S$ as
		\begin{equation} \label{SSSS_dyn} 
		\begin{split}S(h)(t,x) &= \frac{1}{\e} \int_{0}^{\e T}\int_{ \frac{\delta}{2}\leq |v|  \leq \frac{2}{ \delta}
		}  h(t+s,X(t+ \frac{s}{\e};t,x,v),V( t+ \frac{s}{\e};t,x,v)) \phi (v) 
		\dd v \dd s\\ 
		&= \frac{1}{\e}\int_{t}^{t+\e T}\int_{ \frac{\delta}{2}\leq |v|  \leq \frac{2}{ \delta}
		}  h( s,X(  \frac{s-t}{\e};0,x,v),V(   \frac{s-t}{\e};0,x,v)) \phi (v) 
		\dd v \dd s.
		\end{split}
		\end{equation}
	\end{definition}

	\unhide

	\vspace{8pt}
	
	\textbf{Acknowledgements.} The authors thank Hongxu Chen to find errors in the original manuscript and fix them. They also thank referee(s) for useful comments which help us to improve the clarity of presentation. The authors thank Yan Guo for his interest and discussions. They also thank Cl\'ement Mouhot, Lello Esposito, Rossana Marra, Misha Feldman, Hyung Ju Hwang, and St\'ephane Mischler for their interest to this project. C.K. especially thanks James Callen(Center for Plasma Theory and Computation) for discussions on the several relevant kinetic models. The authors also thank kind hospitality of MFO at Oberwolfach, ICERM, KAIST-CMC, math/applied math departments of Brown, Cambridge, Princeton, USC(during a summer school organized by Juhi Jang), UMN, UIC, UT-Austin, POSTECH, NTU, Lyon 1, and Paris-Dauphine during this research. The research is supported in part by National Science Foundation under Grant No.1501031, Wisconsin Alumni Research Foundation, and the Herchel Smith foundation of the University of Cambridge.

	
	
	
\hide

\unhide

	\hide

	\unhide
	\end{document}